%% file: twisted_v2.tex
\documentclass[11pt,twoside,a4paper]{report}
\usepackage{color,amscd,amsmath,amssymb,amsthm,latexsym, stmaryrd,authblk,mathabx, shuffle,mathrsfs,hyperref,tikz,tkz-tab} 
\usepackage[latin1]{inputenc}
\usepackage[T1]{fontenc} 
\usepackage[english]{babel}
\providecommand{\AMSclass}[1]{\textbf{\textit{AMS classification.}} #1}

\input{xy}
\xyoption{all}

\theoremstyle{plain}
\newtheorem{theo}{Theorem}[section]
\newtheorem{lemma}[theo]{Lemma}
\newtheorem{cor}[theo]{Corollary}
\newtheorem{prop}[theo]{Proposition}
\newtheorem{defi}[theo]{Definition}

\theoremstyle{remark}
\newtheorem{remark}{Remark}[section]
\newtheorem{notation}{Notations}[section]
\newtheorem{example}{Example}[section]

\input{arbres.tex}
\input{posets.tex}

\input{graphes.tex}

\newcommand{\bfC}{\mathbf{C}}
\newcommand{\bfD}{\mathbf{D}}

\newcommand{\calA}{\mathcal{A}}
\newcommand{\calB}{\mathcal{B}}
\newcommand{\calC}{\mathcal{C}}
\newcommand{\calD}{\mathcal{D}}
\newcommand{\calE}{\mathcal{E}}
\newcommand{\calF}{\mathcal{F}}
\newcommand{\calG}{\mathcal{G}}
\newcommand{\calH}{\mathcal{H}}
\newcommand{\calI}{\mathcal{I}}

\newcommand{\calM}{\mathcal{M}}
\newcommand{\calO}{\mathcal{O}}
\newcommand{\calP}{\mathcal{P}}
\newcommand{\calQ}{\mathcal{Q}}
\newcommand{\calR}{\mathcal{R}}
\newcommand{\calS}{\mathcal{S}}
\newcommand{\calT}{\mathcal{T}}

\newcommand{\calU}{\mathcal{U}}

\newcommand{\K}{\mathbb{K}}

\newcommand{\N}{\mathbb{N}}
\newcommand{\Z}{\mathbb{Z}}
\newcommand{\Q}{\mathbb{Q}}

\newcommand{\sym}{\mathfrak{S}}
\newcommand{\cl}{\mathrm{cl}}

\newcommand{\id}{\mathrm{Id}}
\newcommand{\coinv}{\mathrm{coInv}}

\newcommand{\Com}{\mathbf{Com}}

\newcommand{\QSym}{\mathbf{QSym}}
\newcommand{\NSym}{\mathbf{NSym}}
\newcommand{\WQSym}{\mathbf{WQSym}}
\newcommand{\FQSym}{\mathbf{FQSym}}
\newcommand{\Sym}{\mathbf{Sym}}

\newcommand{\Char}{\mathrm{Char}}
\newcommand{\InfChar}{\mathrm{InfChar}}

\newcommand{\prim}{\mathrm{Prim}}

\newcommand{\im}{\mathrm{Im}}
\newcommand{\vect}{\mathrm{Vect}}
\renewcommand{\ker}{\mathrm{Ker}}
\newcommand{\Hom}{\mathrm{Hom}}

\newcommand{\val}{\mathrm{val}}

\newcommand{\ic}{\mathrm{ic}}
\newcommand{\comp}{\mathcal{C}\hspace{-0.5mm}\mathit{omp}}
\newcommand{\exch}{\mathbf{Exch}}
\newcommand{\elim}{\mathrm{elim}_\emptyset}
\newcommand{\nd}{\mathrm{nd}}
\newcommand{\reg}{\mathrm{reg}}
\newcommand{\HE}{\mathit{HE}}
\newcommand{\co}{\mathrm{co}}
\newcommand{\ind}{\mathrm{ind}}
\newcommand{\qsh}{\mathrm{qsh}}
\newcommand{\CS}{\mathrm{CS}}
\newcommand{\AC}{\mathrm{AC}}
\newcommand{\cocom}{\mathrm{cocom}}
\newcommand{\conc}{\mathrm{conc}}
\newcommand{\pp}{+\hspace{-1mm}+}
\newcommand{\sh}{\mathrm{sh}}

\newcommand{\CB}{\mathrm{CB}}
\newcommand{\lcom}{\ell\Com}
\newcommand{\double}{\mathrm{double}}
\newcommand{\qc}{\mathrm{qcom}}
\newcommand{\PI}{\mathrm{PI}}
\newcommand{\pack}{\mathrm{pack}}
\newcommand{\red}{\mathrm{Red}}

\begin{document}

\title{Combinatorial twisted bialgebras and combinatorial twisted double bialgebras}
\date{}
\author{Lo\"ic Foissy}
\affil{\small{Univ. Littoral Côte d'Opale, UR 2597
LMPA, Laboratoire de Mathématiques Pures et Appliquées Joseph Liouville
F-62100 Calais, France}.\\ Email: \texttt{foissy@univ-littoral.fr}}

\maketitle

\begin{abstract}

We provide in this text a new framework for a family of combinatorial Hopf algebras. We adopt for this the language of species and twisted bialgebras (also known as Hopf monoids). 
We define the category of combinatorial twisted bialgebras (briefly, CTBs), a class of twisted bialgebras satisfying strong combinatorial conditions of compatibilities of the product and the coproduct with the underlying combinatorial structure. 
The considered morphisms also have to respect the combinatorial structure, in some sense.
Many classical examples fit into this framework, including Hopf algebras of graphs, rooted trees, posets, hypergraphs, symmetric functions, boolean maps, etc. Most of known combinatorial Hopf algebras can be obtained from a CTB by application of a Fock functor, or by duality. 

A central result is the introduction of a universal CTB denoted by $\mathbf{Exch}$, built from collections of compositions satisfying a certain exchange condition. 
It has a special character $\epsilon$, such that the pair $(\mathbf{Exch},\epsilon)$ plays a role analogous to that of quasi-symmetric functions for graded and connected Hopf algebras: it is the terminal object in the category of CTBs equipped with combinatorial characters. 
This leads to universal morphisms that encode important information about the combinatorial structure of coproducts.
The natural order on the set of compositions on a fixed set is used to define three important sub-objects $\mathbf{Exch}_-$, $\mathbf{Exch}_+$ and $\mathbf{Exch}_{+\hspace{-2mm}+}$ of $\mathbf{Exch}$. 

We also studies polynomial invariants associated with CTBs. It establishes a correspondence between characters and polynomial invariants. 
Well-known invariants such as the chromatic polynomial of graphs and the Ehrhart polynomials of posets arise naturally within this framework. 
We show that, for a given degree, only finitely many polynomial invariants can occur. We show that these universal polynomial invariants can be used to give cancellation-free formulas for the antipode or the Eulerian idempotent, under conditions on the CTB,
and use this for an example based on rooted trees, coming from a Hopf algebra used in linguistics. Results are also obtained for families of polynomial invariants with the alternating sign condition, or satisfying a reciprocity principle.

We also extend the theory to combinatorial twisted double bialgebras (briefly, CTDBs), which involve two interacting coproduct structures. These generalize previously known examples arising from graphs, posets, hypergraphs, and rooted trees. 
The twisted approach preserves the combinatorial richness of these structures while providing a unified algebraic framework.
We show that not any CTB can be made a CTDB, leading to the introduction of the notion of quasi-cocommutativity. In particular, we show that $\mathbf{Exch}_+$ is not a CTDB, and construct a sub-object which is, in a unique way, and maximal for the inclusion.

In the cocommutative case, simplifications naturally occur, and we show that in this case, a CTB of hypergraphs plays a central role, especially two sub-objects based on connectedness systems and antichain covers.
We also show that, in the cocommutative case, the universal polynomial invariants are always chromatic polynomials of hypergraphs, no matter which CTB is considered.

Finally, forgetting the combinatorial underlying structure, we study the algebraic structure of $\mathrm{Exch}$, with the help of a twisted version of Loday and Ronco's rigidity theorem for infinitesimal bialgebras. \\

\AMSclass{16T05 16T30 05E16 05E05 05C065 06A11}
\end{abstract}

\tableofcontents

\chapter*{Introduction}

\section*{Combinatorial Hopf algebras}

In the past decades, numerous Hopf algebras or bialgebras based on combinatorial objects appear in the literature. An important example occurred in the nineties with Connes and Kreimer's work on the Renormalization in Quantum Field Theory
\cite{Connes1998,Connes2000,Connes2001}, where this mathematical object, an inductive procedure of extractions of poles, was described as a Birkhoff decomposition in the group of characters of certain Hopf algebras based on Feynman graphs.
A toy-model for these objects is the celebrated Butcher-Connes-Kreimer Hopf algebra of rooted trees, which was quickly related by duality to another Hopf algebra, introduced by Grossman and Larson in a context of differential calculus \cite{Grossman2005,Grossman89,Grossman90,Panaite2000,Hoffman2003}. 
Similar Hopf algebras of rooted trees also appeared in Hairer's theory of regularity structures, introduced to solve stochastics PDEs \cite{Hairer2014,Bruned2015,Bruned2019}, with a heavy complementary structure on the vertices and on the edges. 
Other ones also appeared recently in linguistics, in Marcolli, Chomsky and Berwick's work \cite{Marcolli2025}. 
To stay with trees, we can also cite other Hopf algebras of planar trees or binary trees, with a similar coproduct of admissible cuts \cite{Foissy1,Foissy2,Holtkamp2003,Loday1998}, some related to the theory of operads. 

Other Hopf algebras with a strong combinatorial flavor appear in the theory of representations of family of groups (especially symmetric groups), such as the Hopf algebra of symmetric functions $\Sym$,
 its generalization the Hopf algebra of quasi-symmetric functions $\QSym$ and its dual the Hopf algebra of quasi-symmetric functions $\NSym$ \cite{Gelfand1995,Hazewinkel2005,Malvenuto1995}.
Permutations themselves can be organized as a Hopf algebra into several ways, some commutative \cite{Hivert2008}, the most interesting being perhaps the one known as the Malvenuto-Reutenauer Hopf algebra, also called the Hopf algebra of free quasi-symmetric functions $\FQSym$ \cite{Duchamp2002,Malvenuto1995,Foissy5}.
In a similar spirit, Hopf-algebraic structures are defined on packed words (or equivalently, set compositions) \cite{Novelli2006-2}, parking functions \cite{Novelli2007}, uniform block permutations \cite{Aguiar2008}, etc.

Many results of graph theory have been lifted to the level of Hopf algebras, thanks to bialgebraic structures on graphs \cite{Foissy36,Foissy58}. In particular, the chromatic polynomial can be seen as a Hopf algebra morphism from graphs to polynomials, as well as the Fortuin and Kasteleyn polynomial (a variant of Tutte's polynomial).
These results are extended to hypergraphs, giving several notions of chromatic polynomials or symmetric functions \cite{Foissy44,Ebrahimi-Fard2022,Aval2019}, mixed graphs \cite{Foissy45} or to posets and quasi-posets \cite{Malvenuto2011,Foissy37,Foissy33}.
One can also cite noncrossing partitions, organized as a Hopf algebra which plays a major role in the theory of free probability \cite{Foissy43,Ebrahimi-Fard2020,Nica2006}.

Despite all these examples and applications, there is no consensual definition of a combinatorial Hopf algebra. A tentative has been given by Aguiar, Bergeron and Sottile \cite{Aguiar2006-2}. 
In their setting, a combinatorial Hopf algebra (CHA) is a pair formed by a graded and connected Hopf algebra $H$ and a character $f$ of $H$. The main interest is that $\QSym$, with a particular character $\epsilon$, is a terminal object in the category of CHAs as they define it.
This allows to define morphisms from combinatorial objects to quasi-symmetric functions in a purely algebraic way, such as the chromatic polynomial symmetric function of a graph \cite{Stanley1995}.
We won't follow this definition here, but we will use the formalism of twisted bialgebras (with the terminology of \cite{Patras2004}), also called Hopf monoids, developed in \cite{Aguiar2010}.

\section*{Species and twisted objects}

Let us detail what these objects are. The category of species was introduced by Joyal \cite{Joyal1981,Joyal1986}, in order to give a better understanding of classical manipulations of formal series in combinatorics.
Formally, a species is a functor from the category of finite sets, with bijections, to the category of sets, or of $\K$-vector spaces for linear species. In other words, for a species $\calQ$:
\begin{itemize}
\item For any finite set $X$, $\calQ[X]$ is a set, or is a $\K$-vector space if $\calQ$ is a linear species.
\item For any bijection $\sigma:X\longrightarrow Y$, $\calQ[\sigma]:\calQ[X]\longrightarrow \calQ[Y]$ is a map, linear if $\calQ$ is a linear species.
\end{itemize}
These correspondence are of course compatible with the composition of bijections, see Definition \ref{defi1.1.1} for the full details. Here are classical examples of species, which will be used all along this text.
\begin{itemize}
\item Graphs. For any finite set $X$, $\calG[X]$ is the set of graphs whose set of vertices is $X$. For example,
\begin{align*}
\calG[\{x\}]&=\{\grdun{$x$}\},\\
\calG[\{x,y\}]&=\left\{\begin{array}{c}
\grdun{$x$}\grdun{$y$}, \grddeux{$x$}{$y$}
\end{array}\right\},\\
\calG[\{x,y,z\}]&=\left\{\begin{array}{c}
\grdun{$x$}\grdun{$y$}\grdun{$z$},\grdun{$x$}\grddeux{$y$}{$z$},\grdun{$y$}\grddeux{$x$}{$z$},\grdun{$z$}\grddeux{$x$}{$y$},\grdtroisdeux{$x$}{$z$}{$y$},\grdtroisdeux{$y$}{$z$}{$x$},\grdtroisdeux{$z$}{$y$}{$x$},\grdtroisun{$x$}{$y$}{$z$}
\end{array}\right\}.
\end{align*}
\item Posets. For any finite set $X$, $\calP[X]$ is the set of posets whose set of vertices is $X$, that is to say partial orders on $X$. Representing posets by their Hasse graphs, the edges being oriented upwards,
\begin{align*}
\calP[\{x\}]&=\{\tdun{$x$}\},\\
\calP[\{x,y\}]&=\left\{\begin{array}{cl}
\tdun{$x$}\tdun{$y$},\tddeux{$x$}{$y$},\tddeux{$y$}{$x$}
\end{array}\right\},\\
\calP[\{x,y,z\}]&=\left\{\begin{array}{cl}
\tdun{$x$}\tdun{$y$}\tdun{$z$},\tdun{$x$}\tddeux{$y$}{$z$},\tdun{$x$}\tddeux{$z$}{$y$},\tdun{$y$}\tddeux{$x$}{$z$},\tdun{$y$}\tddeux{$z$}{$x$},\tdun{$z$}\tddeux{$x$}{$y$},\tdun{$z$}\tddeux{$y$}{$x$},\\
\tdtroisun{$x$}{$z$}{$y$},\tdtroisun{$y$}{$z$}{$x$},\tdtroisun{$z$}{$y$}{$x$},
\pdtroisun{$x$}{$y$}{$z$},\pdtroisun{$y$}{$x$}{$z$},\pdtroisun{$z$}{$x$}{$y$},
\tdtroisdeux{$x$}{$y$}{$z$},\tdtroisdeux{$x$}{$z$}{$y$},\tdtroisdeux{$y$}{$x$}{$z$},
\tdtroisdeux{$y$}{$z$}{$x$},\tdtroisdeux{$z$}{$x$}{$y$},\tdtroisdeux{$z$}{$y$}{$x$}
\end{array}\right\}.
\end{align*}
\item Subsets. For any finite set $X$, $\calS[X]$ is the set of subsets of $X$. For example,
\begin{align*}
\calS[\{x\}]&=\{\emptyset,\:\{x\}\},\\
\calS[\{x,y\}]&=\{\emptyset,\: \{x\},\: \{y\},\: \{x,y\}\},\\
\calS[\{x,y,z\}]&=\{\emptyset,\:\{x\},\:\{y\},\:\{z\},\:\{x,y\},\:\{x,z\},\:\{y,z\},\:\{x,y,z\}\}.
\end{align*}
\item Compositions. For any finite set $X$, $\comp[X]$ is the set of compositions (or ordered partitions) on $X$, that is to say families $\pi=(X_1,\ldots,X_k)$ of nonempty subsets of $X$, such that $X=X_1\sqcup\cdots\sqcup X_k$. We denote by $\comp$ the species of compositions. For example,
\begin{align*}
\comp[\{x\}]&=\{(\{x\})\},\\
\comp[\{x,y\}]&=\{(\{x,y\}),\:(\{x\},\{y\}),\: (\{y\},\{x\})\},\\
\comp[\{x,y,z\}]&=\left\{\begin{array}{c}
(\{x,y,z\}),\:(\{x,y\},\{z\}),\:(\{x,z\},\{y\}),\:(\{y,z\},\{x\}),\\
(\{x\},\{y,z\}),\:(\{y\},\{x,z\}),\:(\{z\},\{x,z\}),\\
(\{x\},\{y\},\{z\}),\:(\{x\},\{z\},\{y\}),\:(\{y\},\{x\},\{z\}),\\
(\{y\},\{z\},\{x\}),\:(\{z\},\{x\},\{y\}),\:(\{z\},\{y\},\{x\})
\end{array}\right\}.\end{align*}
We will lighten the notations by deleting some brackets and comas, as follows:
\begin{align*}
\comp[\{x\}]&=\{(x)\},\\
\comp[\{x,y\}]&=\{(xy),\:(x,y),\: (y,x)\},\\
\comp[\{x,y,z\}]&=\left\{\begin{array}{c}
(xyz),\:(xy,z),\:(xz,y),\:(yz,x),\\
(x,yz),\:(y,xz),\:(z,xz),\\
(x,y,z),\:(x,z,y),\:(y,x,z),\\
(y,z,x),\:(z,x,y),\:(z,y,x)
\end{array}\right\}.
\end{align*}\end{itemize}

Linear species form a symmetric monoidal category for the Cauchy tensor product, which allows to define algebras, coalgebras, bialgebras, etc, in the category of species. We now give a short review of these objects, see Definitions \ref{defi1.1.2}, \ref{defi1.1.3} and \ref{defi1.1.4} for the full details.
In a twisted algebra $\calQ$, if $X$ and $Y$ are disjoint finite sets, there exists a map $m_{X,Y}:\calQ[X]\otimes \calQ[Y]\longrightarrow \calQ[X\sqcup Y]$, with an associativity constraint. 
In a twisted coalgebra $\calQ$, if $X$ and $Y$ are disjoint finite sets, there exists a map $\Delta_{X,Y}:\calQ[X\sqcup Y]\longrightarrow\calQ[X]\otimes \calQ[Y]$, with a coassociativity constraint. 
A twisted bialgebra is both a twisted algebra and a twisted coalgebra, with a compatibility between the product and the coproduct. Here are some examples:
\begin{itemize}
\item Graphs. The linearization of $\calG$ is a twisted bialgebra. The product of two graphs with disjoint sets of vertices is their disjoint union. For example,
\[m_{\{x,y,z\},\{u,v\}}\left(\hspace{-2mm}\begin{array}{c}
\tdtroisun{$x$}{$z$}{$y$}
\end{array}\hspace{-1mm}\otimes\hspace{-2mm} \begin{array}{c}
\tddeux{$u$}{$v$}
\end{array}\hspace{-2mm}\right)
=\tdtroisun{$x$}{$z$}{$y$}\tddeux{$u$}{$v$}.\]
The coproduct is given by restrictions. For example, writing the coproduct as the sum over all possible decompositions of $\{x,y,z\}$ as a disjoint union of two subsets,
\begin{align*}
\Delta\left(\hspace{-2mm}\begin{array}{c}
\tdtroisun{$x$}{$z$}{$y$}
\end{array}\hspace{-2mm}\right)&=\tdtroisun{$x$}{$z$}{$y$}\otimes 1+\tdun{$x$}\otimes\tdun{$y$}\tdun{$z$}+\tdun{$y$}\otimes\tddeux{$x$}{$z$}+\tdun{$z$}\otimes\tddeux{$x$}{$y$}+1\otimes \tdtroisun{$x$}{$z$}{$y$}+\tdun{$y$}\tdun{$z$}\otimes\tdun{$x$}+\tddeux{$x$}{$z$}\otimes\tdun{$y$}+\tddeux{$x$}{$y$}\otimes\tdun{$z$},\\
\Delta\left(\hspace{-2mm}\begin{array}{c}
\grdtroisun{$x$}{$z$}{$y$}
\end{array}\hspace{-2mm}\right)&=\grdtroisun{$x$}{$z$}{$y$}\otimes 1+\tdun{$x$}\otimes\tddeux{$y$}{$z$}+\tdun{$y$}\otimes\tddeux{$x$}{$z$}+\tdun{$z$}\otimes\tddeux{$x$}{$y$}+1\otimes \grdtroisun{$x$}{$z$}{$y$}+\tddeux{$y$}{$z$}\otimes\tdun{$x$}+\tddeux{$x$}{$z$}\otimes\tdun{$y$}+\tddeux{$x$}{$y$}\otimes\tdun{$z$}.
\end{align*} 
The twisted bialgebra $\calG$ is commutative and cocommutative. 
\item Posets. The linearization of $\calP$ is a twisted bialgebra. The product of two posets with disjoint sets of vertices is their disjoint union. For example,
\[m_{\{x,y,z\},\{u,v\}}\left(\hspace{-2mm}\begin{array}{c}
\tdtroisun{$x$}{$z$}{$y$}
\end{array}\hspace{-1mm}\otimes\hspace{-2mm} \begin{array}{c}
\tddeux{$u$}{$v$}
\end{array}\hspace{-2mm}\right)
=\tdtroisun{$x$}{$z$}{$y$}\tddeux{$u$}{$v$}.\]
The coproduct $\Delta_{X,Y}(P)$ is given by restrictions of $P$ to $X$ and $Y$, with the complementary condition that $Y$ should be an upper part of $P$; otherwise, the coproduct is zero. For example, 
\begin{align*}
\Delta\left(\hspace{-2mm}\begin{array}{c}
\tdtroisun{$x$}{$z$}{$y$}
\end{array}\hspace{-2mm}\right)&=\tdtroisun{$x$}{$z$}{$y$}\otimes 1+\tdun{$x$}\otimes\tdun{$y$}\tdun{$z$}+\tddeux{$x$}{$y$}\otimes\tdun{$z$}+\tddeux{$x$}{$z$}\otimes\tdun{$y$}+1\otimes \tdtroisun{$x$}{$z$}{$y$},\\
\Delta\left(\hspace{-2mm}\begin{array}{c}
\tdtroisdeux{$x$}{$y$}{$z$}
\end{array}\hspace{-2mm}\right)&=\tdtroisdeux{$x$}{$y$}{$z$}\otimes 1+\tddeux{$x$}{$y$}\otimes\tdun{$z$}+\tdun{$x$}\otimes\tddeux{$y$}{$z$}+1\otimes \tdtroisdeux{$x$}{$z$}{$y$}.
\end{align*} 
The twisted bialgebra $\calP$ is commutative but not cocommutative. 
\item Compositions. The linearization $\ell\comp$ of $\comp$ is a twisted bialgebra, with the quasi-shuffle product $\squplus$ and the deconcatenation coproduct. For example,
\begin{align*}
\squplus_{\{x,y,z\},\{u,v\}}((xy,z)\otimes (u,v))&=(xy,z,u,v)+(xy,zu,v)+(xy,u,z,v)+(xyu,z,v)\\
&+(u,xy,z,v)+(xy,u,zv)+(xyu,zv)+(u,xy,zv)\\
&+(xy,u,v,z)+(xyu,v,z)+(u,xy,v,z)+(u,xyv,z)\\
&+(u,v,xy,z),
\end{align*}
and
\begin{align*}
\Delta((xy,z,u,v))&=(xy,z,u,v)\otimes 1+(xy,z,u)\otimes (v)+(xy,z)\otimes (u,v)\\
&+(xy)\otimes (z,u,v)+1\otimes (xy,z,u,v).
\end{align*}
This twisted bialgebra is commutative but not cocommutative. It can be seen as a twisted version of $\QSym$, in the sense that it satisfies also a terminal property, similar to the one of \cite{Aguiar2006-2}. Let us explicit this. We define a character $\epsilon$ of $\ell\comp$ by
\[\epsilon((X_1,\ldots,X_k))=\begin{cases}
1\mbox{ if }k\leq 1,\\
0\mbox{ otherwise}.
\end{cases}\]
If $\calQ$ is a twisted bialgebra, such that $\calQ[\emptyset]$ is one-dimensional, and if $f$ is a character of $\calQ$,
there exists a unique morphism of twisted bialgebras $\phi$ from $\calQ$ to $\ell\comp$ such that $\epsilon\circ \phi=f$,
see Theorem \ref{theo1.5.6}.
\end{itemize}

Several functors from species to vector spaces are introduced in \cite{Aguiar2010}, which send a twisted bialgebra to a (classical) bialgebra. We will focus in this paper on the bosonic Fock functor $\calF$, which sends a species $\calQ$ to the vector space
\[\calF[\calQ]=\bigoplus_{n=0}^\infty \coinv(\calQ[\{1,\ldots,n\}]),\]
where $\coinv(\calQ[\{1,\ldots,n\}])$ is the quotient of $\calQ[\{1,\ldots,n\}]$ by the subspace generated by elements $\calQ[\sigma](b)-b$, with $b\in \calQ[\{1,\ldots,n\}]$ and $\sigma$ in the $n$-th symmetric group $\sym_n$.
Roughly speaking, this functor consists to "forget" the elements of the underlying sets, and keeping only the structure of the elements of $\calQ$. For example:
\begin{itemize}
\item Graphs. The bialgebra $\calF[\calG]$ has for basis the set of isomorphism classes (isoclasses) of graphs:
\[\left\{1,\:\tun,\:\tun\tun,\:\tdeux,\:\tun\tun\tun,\:\tun\tdeux,\:\ttroisun,\:\grtroisun,\:\ldots\right\}\]
 The product is the disjoint union:
\[m\left(\hspace{-2mm}\begin{array}{c}
\ttroisun
\end{array}\hspace{-1mm}\otimes\hspace{-2mm} \begin{array}{c}
\tdeux
\end{array}\hspace{-2mm}\right)
=\ttroisun\tdeux.\]
The coproduct is given by restrictions:
\begin{align*}
\Delta\left(\hspace{-2mm}\begin{array}{c}
\ttroisun
\end{array}\hspace{-2mm}\right)&=\ttroisun\otimes 1+\tun\otimes\tun\tun+\tun\otimes\tdeux+\tun\otimes\tdeux+1\otimes \ttroisun+\tun\tun\otimes\tun+\tdeux\otimes\tun+\tdeux\otimes\tun,\\
\Delta\left(\hspace{-2mm}\begin{array}{c}
\grtroisun
\end{array}\hspace{-2mm}\right)&=\grtroisun\otimes 1+\tun\otimes\tdeux+\tun\otimes\tdeux+\tun\otimes\tdeux+1\otimes \grtroisun+\tdeux\otimes\tun+\tdeux\otimes\tun+\tdeux\otimes\tun.
\end{align*}
This is the bialgebra of graphs studied in \cite{Foissy36,Foissy58}.
\item Posets. The bialgebra $\calF[\calP]$ has for basis the set of isoclasses of posets:
\[\left\{\begin{array}{c}
1,\:\:\tun,\:\:\tun\tun,\:\:\tdeux,\:\:\tun\tun\tun,\:\tun\tdeux,\:\ttroisun,\:\ptroisun,\:\ttroisdeux,\:\ldots
\end{array}\right\}.\]
The product is the disjoint union:
\[m\left(\hspace{-2mm}\begin{array}{c}
\ttroisun
\end{array}\hspace{-1mm}\otimes\hspace{-2mm} \begin{array}{c}
\tdeux
\end{array}\hspace{-2mm}\right)
=\ttroisun\tdeux.\]
The coproduct is given by restriction to lower and upper parts:
\begin{align*}
\Delta\left(\hspace{-2mm}\begin{array}{c}
\ttroisun
\end{array}\hspace{-2mm}\right)&=\ttroisun\otimes 1+\tun\otimes\tun\tun+\tdeux\otimes\tun+\tdeux\otimes\tun+1\otimes \ttroisun,\\
\Delta\left(\hspace{-2mm}\begin{array}{c}
\ttroisdeux
\end{array}\hspace{-2mm}\right)&=\ttroisdeux\otimes 1+\tdeux\otimes\tun+\tun\otimes\tdeux+1\otimes \ttroisdeux.
\end{align*} 
This is the Hopf algebra of posets of \cite{Foissy37}, see also \cite{Malvenuto2011,Foissy33} for generalizations to quasi-posets or double posets. Seeing rooted trees as Hasse graphs, $\calF[\calP]$ contains the Butcher-Connes-Kreimer Hopf algebra.
\item Compositions. The isoclass of a composition $(X_1,\ldots,X_k)$ is identified with the sequence of positive integers $(|X_1|,\dots, |X_k|)$. Therefore, $\calF[\comp]$ is a Hopf algebra on these sequences. For example, if $a,b,c,d\geq 1$,
\begin{align*}
(a,b)\squplus (c,d)&=(a,b,c,d)+(a,b+c,d)+(a,c,b,d)+(a+c,b,d)+(c,a,b,d)\\
&+(a,c,b+d)+(a+c,b+d)+(c,a,b+d)\\
&+(a,c,d,b)+(a+c,d,b)+(c,a,d,b)+(c,a+d,b).
\end{align*}
The coproduct is the deconcatenation:
\begin{align*}
\Delta((a,b,c,d))&=(a,b,c,d)\otimes 1+(a,b,c)\otimes (d)+(a,b)\otimes (c,d)+(a)\otimes (b,c,d)+1\otimes (a,b,c,d).
\end{align*}
This is the Hopf algebra $\QSym$ of quasi-symmetric functions, mentioned before. It is possible to obtain in a similar way the quasi-shuffle Hopf algebras associated to a semigroup of \cite{Hoffman2000,Hoffman2020}.
\end{itemize}
Mostly all combinatorial Hopf algebras (whatever we mean by this) can be applied by the application of a Fock functor to a twisted bialgebra, see \cite{Aguiar2010} for more examples.
We mention here also the noncommutative Butcher-Connes-Kreimer Hopf algebra \cite{Holtkamp2003,Foissy1,Foissy2} and its dual, the Grossman-Larson Hopf algebra \cite{Grossman2005,Grossman89}, the Hopf algebra of free quasi-symmetric functions \cite{Malvenuto1995,Duchamp2002}, 
the Hopf algebra of noncrossing partitions used in the theory of free probability \cite{Nica2006,Ebrahimi-Fard2020}, and various bialgebras on hypergraphs \cite{Ebrahimi-Fard2022,Aval2019,Foissy44}.

\section*{Combinatorial twisted bialgebras}

We focus in this text on Hopf algebras coming from a twisted bialgebra by the application of a Fock functor (mostly, the bosonic functor $\calF$), assuming that this twisted bialgebra $\ell\calQ$ is the linearization of a species $\calQ$ 
(that is to say, for any finite set $X$, $\calQ[X]$ is a basis of $\ell\calQ[X]$), with the following properties:
\begin{itemize}
\item $\calQ[\emptyset]$ is reduced to the unit of $\ell\calQ$.
\item The product $m_{X,Y}(a\otimes b)$ of two elements $a\in \calQ[X]$ and $b\in \calQ[Y]$ belongs to $\calQ[X\sqcup Y]$, or is zero.
\item The coproduct $\Delta_{X,Y}(b)$ of an element of $\calQ[X\sqcup Y]$ is a tensor of an element of $\calQ[X]$ and an element of $\calQ[Y]$, or is zero.
\end{itemize} 
This is what we shall call a combinatorial twisted bialgebra (briefly, CTB), see Definition \ref{defi2.1.1}. This includes the Butcher-Connes-Kreimer Hopf algebra and its noncommutative version, as well as its decorated ones,
the Hopf algebras of hypergraphs, of posets, of quasi-posets, of noncrossing partitions, the Hopf algebra $\Sym$ of symmetric functions, and many more. 
This does not include the Hopf algebra of quasi-symmetric functions $\QSym$, or its generalizations on uniform block partitions, packed words or parking functions, or the Grossman-Larson Hopf algebras, but includes their duals.
However, all the known CHAs do not fit here anyway: one may think to the Hopf algebra of quasi-symmetric functions $\FQSym$, or its dual (as it is self-dual), at least with its usual basis of permutations; a change of basis may insert it in this setting, but not in a trivial way.

In the sequel, with an abuse of notation, we shall denote the CTB by $\calQ$ (which is a species) instead of $\ell\calQ$ (which is a linear species).
We organize CTBs into a category: if $\calQ$ and $\calQ'$ are two CTBs, a morphism of CTBs from $\calQ$ to $\calQ'$ is a morphism of twisted bialgebras, sending any element of $\calQ$ to an element of $\calQ'$, or to $0$ (Definition \ref{defi2.1.2}).
We also define a notion of quotients (Definition \ref{defi2.1.6}). However, if any morphism $\phi:\calQ\longrightarrow \calQ'$
of CTBs induces a morphism $\overline{\phi}:\calQ/\ker(\phi)\longrightarrow \im(\phi)$, this morphism $\overline{\phi}$ has no reason to be injective, to the difference with Noether's theorem. Nevertheless, it does not take the value 0 on $\calQ/\ker(\phi)$.

The strong combinatorial conditions imposed to a CTB $\calB$ makes its algebraic structure rather rigid. We introduce indecomposable elements (Definition \ref{defi2.2.1}), namely elements of $\calB$ which cannot be written as a product of two elements of $\calB$ different from the unit. 
An easy induction shows that any element of $\calB$ can be decomposed into a product of indecomposable elements, and the existence of the coproduct imposes that this decomposition is unique, up to the order of the factors (Theorem \ref{theo2.2.4}).
In other words, as a twisted algebra, $\calB$ is generated by the subspecies of indecomposable elements, with certain commutation relations between these elements.
For graphs, hypergraphs, posets, rooted trees, and so on, this notion of indecomposability coincide with the notion of connectedness. For other cases (the example $\exch$ to come), it can be more subtle to describe.
This decomposition implies that to any element $b\in \calB[X]$, we can associate an equivalence $\sim_b^{(\ic)}$ on $X$, whose classes are the indecomposable factors of the decomposition of $b$ into indecomposable elements.
The number of classes of $\sim_b^{(\ic)}$ is denoted by $\ic(b)$. This decomposition will play an important role in the sequel. 
If $X$ is a finite set, then for any $b\in \calB[X]$, obviously, $\ic(b)\leq |X|$. Elements such that $\ic(b)=|X|$ are called the discrete elements of $\calB$, and they form a sub-CTB of $\calB$, denoted by $\bfD\calB$ (Definition \ref{defi2.2.7} and Proposition \ref{prop2.2.8}).
Discrete graphs are the graphs without edges, also called totally disconnected. Discrete posets are posets whose order is the equality of the underlying set.

\section*{The CTB $\exch$}

We already mentioned Aguiar, Bergeron and Sottile's characterization of $\QSym$ (more precisely, of the pair $(\QSym,\epsilon))$ as the terminal object in their category of CHAs.
We give a generalization of this result in our framework of CTBs, with $\exch$, in the third chapter.
For any finite set $X$, we consider nonempty subsets $I\subseteq \comp[X]$, with the following condition: if $X=X_1\sqcup X_2$, $\pi_1,\pi'_1\in \comp[X_1]$ and $\pi_2,\pi'_2\in \comp[X_2]$ are such that the concatenations 
$(\pi_1,\pi_2)$ and $(\pi'_1,\pi'_2)$ both belong to $I$, then $(\pi'_1,\pi_2)\in I$. This is what we call the exchange condition (Definition \ref{defi3.1.1}), and this defines a species $\exch$.
This species is rather big, and not so easy to handle with. Let us mention, as an illustration, that if $|X|=3$, then $|\exch(X)|=8191$, and that we do not know the cardinality of $\exch[X]$ when $|X|\geq 4$. 
Nevertheless, we define a CTB structure on $\exch$. The product is an extension of the quasi-shuffle product $\squplus$ on $\comp$. For example, if $X_1$, $X_2$, $Y$ are pairwise disjoint nonempty finite sets, then
\begin{align*}
\{(X_1,X_2),\: (X)\}\squplus \{(Y)\}&=\left\{\begin{array}{c}
(X_1,X_2,Y),\:(X_1,Y,X_2),\:(Y,X_1,X_2),\\
(X_1\sqcup Y,X_2),\:(X_1,X_2\sqcup Y),\:(X$, $Y),\:(Y,X),\:(X\sqcup Y)
\end{array}\right\}.\end{align*}
The coproduct is a sort of deconcatenation coproduct. For example, if $X$, $Y$, $Z$ are pairwise disjoint nonempty finite sets, 
if $I=\{(X\sqcup Y,Z),\:(X,Y,Z),\:(Y,X,Z),\:(X\sqcup Z,Y)\}$,
\begin{align*}
\Delta(I)&=1\otimes I+\{(X)\}\otimes \{(Y,Z)\}+\{(Y)\}\otimes \{(X,Z)\}+\{(X\sqcup Y),\:(X,Y),\:(Y,X)\}\otimes \{(Z)\}\\
&+\{(X\sqcup Z)\}\otimes \{(Y)\}+I\otimes 1.
\end{align*}
This example show that the coproduct is not cocommutative. We obtain indeed a CTB (Theorem \ref{theo3.1.16}). By the way, We can notice that the definition of $\squplus$ can be extended to the whole $\calS(\comp)$ (forgetting the exchange condition), and remains associative and commutative. 
However, the extension of $\Delta$ to the whole $\calS(\comp)$ is not coassociative.

A character of a CTB $\calB$ sends any element $b$ of $\calB$ to a scalar and is multiplicative. A character of $\calB$ is combinatorial if it only takes the values $0$ and $1$ on $\calB$ (Definition \ref{defi3.2.1}). 
If $f$ is a combinatorial character of a CTB $\calB$, it defines a twisted subalgebra $\calB_f$ of $\calB$, considering only the elements of $\calB$ sent to $1$. 
Such a twisted subalgebra is said to be special (Definition \ref{defi3.2.2} and Proposition \ref{prop3.2.5}). Not all the twisted subalgebras of $\calB$ are special, but sub-CTBs of $\calB$ are. An example of combinatorial character of $\exch$, defined in Proposition \ref{prop3.2.7} is given by
\begin{align}
\label{EQ1}
&\forall I\in \exch[X],&\epsilon(I)&=\begin{cases}
1\mbox{ if }(X)\in I,\\
0\mbox{ otherwise}.
\end{cases}\end{align}
This pair $(\exch,\epsilon)$ generalizes $\QSym$ and its terminal property: Theorem \ref{theo3.2.9} states that if $\calB$ is a CTB and $f$ is a combinatorial character of $\calB$, there exists a unique morphism $\Gamma_f^\calB:\calB\longrightarrow \exch$ of CTBs, such that $\epsilon\circ \Gamma_f=f$.
Let us describe this morphism. For any $b\in \calB[X]$, $\calC_f(b)$ is the set of compositions $(X_1,\ldots,X_k)\in \comp[X]$ such that
\[f^{\otimes k}\circ \Delta_{X_1,\ldots,X_k}(b)=1.\]
It may be empty. If not, it satisfies the exchange condition. Then 
\begin{align*}
\Gamma_f^\calB(b)=\begin{cases}
\calC_f(b)\mbox{ if }\calC_f(b)\neq \emptyset,\\
0\mbox{ otherwise}.
\end{cases}\end{align*}
 Any CTB $\calB$ has at least two special twisted subalgebras, $\calB$ and $\bfD\calB$, corresponding to two characters $\epsilon_-^\calB$ and $\epsilon_+^\calB$: for any finite set $X$, and any $b\in \calB[X]$,
\begin{align*}
\epsilon^\calB_-(b)&=1,&\epsilon^\calB_+(b)&=\begin{cases}
1\mbox{ if }\ic(b)=|X|,\\
0\mbox{ otherwise}.
\end{cases}\end{align*}
This corresponds to two morphisms of CTB $\Gamma_+^\calB$ and $\Gamma_-^\calB$ from $\calB$ to $\exch$. For any finite set $X$, and any $b\in \calB[X]$,
\begin{align*}
\calC_-(b)&=\{(X_1,\ldots,X_k)\in \comp[X]\mid \Delta_{X_1,\ldots,X_k}(b)\neq 0\},\\
\calC_+(b)&=\{(X_1,\ldots,X_k)\in \calC_-(b)\mid \Delta_{X_1,\ldots,X_k}(b)\in \bfD\calB[X_1]\otimes \cdots\otimes \bfD\calB[X_k]\},
\end{align*}
and 
\begin{align*}
\Gamma^\calB_-(b)&=\calC_-(b),&
\Gamma^\calB_+(b)&=\begin{cases}
\calC_+(b)\mbox{ if }\calC_+(b)\neq \emptyset,\\
0\mbox{ otherwise},
\end{cases}\end{align*}
see Corollary \ref{cor3.2.11} and Theorem \ref{theo3.2.13} for more details on these morphisms. For example:
\begin{itemize}
\item Graphs. For the CTB $\calG$ of graphs, if $G\in \calG[X]$, then $\Gamma_-^\calG(G)=\comp[X]$, and $\Gamma_+^\calG(G)$ is the set of compositions $(X_1,\ldots,X_k)$ of $X$ into totally disconnected subgraphs of $G$.
\item Posets. For the CTB $\calP$ of posets, if $P\in \calP[X]$, then $\Gamma_-^\calG(G)$ is the set of compositions $(X_1,\ldots,X_k)$ of $X$ such that if $x\in X_i$ and $y\in X_j$, with $x\leq_P y$, then $i\leq j$,
whereas $\Gamma_+^\calG(G)$ is the set of compositions $(X_1,\ldots,X_k)$ of $X$ such that if $x\in X_i$ and $y\in X_j$, with $x<_\calP y$, then $i<j$.
\end{itemize}
These two morphisms are universal, in the sense that they can be defined for any CTB, and give precious pieces of information on the coproduct and its combinatorial structure.
The morphism $\Gamma_+^\calB$, which can take the value $0$ on $\calB$, a major difference with $\Gamma_-^\calB$, leads to the notion of non-degenerate CTB (Definition \ref{defi3.2.14}), for which $\calC_+(b)$ is never empty.
Graphs, hypergraphs, posets are non-degenerate, but quasi-posets are. 
We show in Theorem \ref{theo3.2.17} how to restrict ourselves to the non-degenerate case by taking a quotient $\nd\calB$
which holds most of the information of $\calB$.

For any finite set $X$, the set $\comp[X]$ is partially ordered by the refinement.
For example, here is the Hasse graph of of $\comp[\{x,y,z\}]$.
\begin{align*}
\xymatrix{(x,z,y)&(x,y,z)&(y,x,z)&(y,z,x)&(z,y,x)&(z,x,y)\\
(x,yz)\ar@{-}[u]\ar@{-}[ru]&(xy,z)\ar@{-}[u]\ar@{-}[ru]&(y,xz)\ar@{-}[u]\ar@{-}[ru]&(yz,x)\ar@{-}[u]\ar@{-}[ru]&
(z,xy)\ar@{-}[u]\ar@{-}[ru]&(xz,y)\ar@{-}[u]\ar@{-}[lllllu]&\\
&&&(xyz)\ar@{-}[lllu]\ar@{-}[llu]\ar@{-}[lu]\ar@{-}[u]\ar@{-}[ru]\ar@{-}[rru]&&}
\end{align*}
In general, $\comp[X]$ has a unique minimal element, $(X)$, and $|X|!$ maximal elements, given by compositions of length $|X|$. 
This allows to define three subspecies of $\exch$ (Definition \ref{defi3.1.1}):
\begin{itemize}
\item $\exch_+$, made by elements of $\exch$ that are upper parts for the refinement order.
\item $\exch_-$, made by elements of $\exch$ that are lower parts for the refinement order.
\item $\exch_{\pp}$, made by elements of $\exch$ uniquely containing maximal elements for this order.
\end{itemize}
The subspecies $\exch_+$ and $\exch_-$ are sub-CTB of $(\exch,\squplus,\Delta)$, see Theorem \ref{theo3.1.16}, 
but $\exch_{\pp}$ is not: it is stable under $\Delta$, but not $\squplus$. 
But we can give it a shuffle product $\shuffle$, which makes it a CTB (Theorem \ref{theo3.4.1}).
We prove that for any CTB $\calB$, $\Gamma_+^\calB$ takes its values in $\exch_+$ (Theorem \ref{theo3.2.13}),
and $\Gamma_-^\calB$ takes its values in $\exch_-$ (Corollary \ref{cor3.2.11}).
To complete these universal morphisms, we also define a morphism $\Gamma_{\pp}^\calB$ from $\calB$ to $\exch_{\pp}$ for any CTB $\calB$ in Theorem \ref{theo3.4.1}.
These universal morphisms can be applied to $\exch$ itself. We shortly denote $\Gamma_\pm$ instead of $\Gamma_\pm^{\exch}$. Theorem \ref{theo3.2.15} shows that $\Gamma_-$ is a projection from $\exch$ to $\exch_-$,
sending any $I\in \exch[X]$ to the lower part of $\comp[X]$ generated by $I$, 
whereas Theorem \ref{theo3.2.16} shows that $\Gamma_+$ is a projection from $\exch$ to $\exch_+$.

it is important to notice that $\ell\exch$ and $\ell\comp$ are not isomorphic, although they both satisfy a terminal property: a key point is that it is not in the same category. It is in the category of CTBs for $\exch$, and in the category of twisted bialgebras $\calB$ such that $\calB[\emptyset]$ is one-dimensional for $\comp$. However, $\ell\comp$ is a quotient of $\ell\exch$, through the morphism $\Theta$, sending $I\in \exch[X]$ to the sum of its elements (Proposition \ref{prop3.3.1}).
In particular, $\comp$ is not a CTB, as the product of two nonempty compositions is a sum of several compositions, and not a single composition. 

Although the Malvenuto-Reutenauer does not come from a CTB (at least, in a simple way), it is closely related to $\exch_{\pp}$. For $n\in \N$, the elements of $\exch_{\pp}[n]$ are in one-to-one correspondence with the elements of $\sym_n$. 
Through this identification, we prove (Proposition \ref{propMalvenutoReutenauer}) that the following is a surjective map of Hopf algebras:
\[\left\{\begin{array}{rcl}
\calF_F[\ell\exch_{\pp}]&\longrightarrow&\FQSym\\
I\in \exch_{\pp}[n]&\longmapsto&\displaystyle \sum_{\sigma \in I}\sigma.
\end{array}\right.\]
A similar result holds for $\exch$, if one replaces $\FQSym$ by the Hopf algebra of packed words $\WQSym$ \cite{Novelli2006-2}, see Proposition \ref{prop3.3.1}

\section*{Polynomial invariants}

Before the description of the next results of this text, let us give some details on the notion of polynomial invariants.
The polynomial algebra $\K[T]$ is a bialgebra, with the multiplicative coproduct $\Delta$ defined by
\[\Delta(T)=T\otimes 1+1\otimes T.\]
If $(B,m,\Delta)$ is a bialgebra, a polynomial invariant is a bialgebra map from $B$ to $\K[X]$.
Here are examples of polynomial invariants.
\begin{itemize}
\item Graphs. Let us take $B=\calF[\calG]$, the bialgebra of graphs of \cite{Foissy36,Foissy58}.
A first polynomial invariant is given by
\[\phi_{\pp}(G)=T^{|V(G)|},\]
where $G$ is a graph and $V(G)$ the set of its vertices. A more interesting one is given by the chromatic polynomial.
Recall that if $G$ is a graph, a valid $N$-coloration of $G$ is a map $c:V(G)\longrightarrow \{1,\ldots,N\}$ such that if two vertices $x$, $y$ of $G$ are related by an edge, then $c(x)\neq c(y)$. There exists a unique polynomial $\phi_{\mathrm{chr}}(G)$ such that for any $N\geq 1$, $\phi_{\mathrm{chr}}(G)(N)$ is the number of valid $N$-colorations of $G$.
Then $\phi_{\mathrm{chr}}$ is a polynomial invariant of $\calF[\calG]$. 
Generalizations on hypergraphs can be found in \cite{Ebrahimi-Fard2022,Aval2019,Foissy44}.
\item Posets. If $P$ is a poset, we denote by $V(P)$ the set of its vertices, and we denote by $n_P$ the number of linear extensions of $P$, that is to say the number of total orders on $V(P)$ that contains the partial order of $P$. 
A first polynomial invariant is given by
\[\phi_{\pp}(P)=\frac{n_P}{|V(P)|!}T^{|V(P)|}.\]
If $N\geq 1$, we consider the sets
\begin{align*}
E_P(N)&=\{c:V(P)\longrightarrow \{1,\ldots,N\}\mid \forall x,y\in V(P),\: x\leq_P y\Longrightarrow c(x)\leq c(y)\},\\
E^{\mathrm{str}}_P(N)&=\{c:V(P)\longrightarrow \{1,\ldots,N\}\mid \forall x,y\in V(P),\: x<_P y\Longrightarrow c(x)<c(y)\}.
\end{align*}
There exists two polynomials $\phi_-(P)$ and $\phi_+(P)$ such that for any $N\geq 1$,
\begin{align*}
\phi_-(P)(N)&=|E_P(N)|,&\phi_+(P)(N)&=|E_P^{\mathrm{str}}(N)|.
\end{align*}
Then $\phi_-$ and $\phi_+$ are both polynomial invariants of $\calF[\calP]$, called the Ehrhart and the strict Ehrhart polynomials \cite{Foissy37}.
\end{itemize}
We proved in \cite{Foissy40} that, under a condition of connectedness, there exists a one-to-one correspondence between characters and polynomial invariants of a given bialgebra $B$.
We generalize this notion to the twisted frame (Definition \ref{defi1.2.20}). If $\calQ$ is a twisted bialgebra, a polynomial invariant of $\calQ$ is a polynomial of the (classical) bialgebra $\calF[\calQ]$.
 In other terms, a polynomial invariant $\phi$ is given by a map $\phi_X:\calQ[X]\longrightarrow \K[T]$, for any finite set $X$, in a way compatible with bijections, and satisfies the following conditions:
\begin{itemize}
\item If $X$ and $Y$ are disjoint finite sets, $a\in \calQ[X]$ and $b\in \calQ[Y]$, then $\phi_{X\sqcup Y}(ab)=\phi_X(a)\phi_Y(b)$.
\item If $X$ is a finite set, and $b\in \calQ[X]$,
\begin{align*}
\phi_X(b)(T+T')&=\sum_{X=X'\sqcup X''}(\phi_{X'}\otimes \phi_{X''})\circ \Delta_{X',X''}(b)(T,T').
\end{align*}
\end{itemize}
We will abuse the notations and we will write $\phi(b)$ instead of $\phi_X(b)$ in the sequel.
We give in the twisted frame a one-to-one correspondence between characters and polynomial invariants of a twisted bialgebra $\calQ$ in Theorem \ref{theo1.2.23}.
With the condition that $\calQ[\emptyset]$ is one-dimensional, for any character $f$ of $\calQ$, there exists a unique polynomial invariant $\phi_f$ of $\calQ$ such that for any finite set $X$, for any $b\in \calQ[X]$, $\phi_f(b)(1)=f(b)$. Moreover, for any $x\in \calQ[X]$,
\[\phi(b)=\sum_{k=1}^{|X|}\sum_{(X_1,\ldots,X_k)\in \comp[X]} f^{\otimes k}\circ \Delta_{X_1,\ldots,X_k}(x) \frac{T(T-1)\cdots (T-k+1)}{k!}.\]
If $\calB$ is a CTB, we can apply this to any character of $\calB$, and especially to $\epsilon_\pm^\calB$. We obtain in this way two polynomial invariants $\phi_+^\calB$ and $\phi_-^\calB$, with a general combinatorial interpretation (proposition \ref{prop3.5.6}).
For graphs, $\phi^\calG_+$ is the chromatic polynomial, and $\phi^\calG_-$ is the polynomial invariant described earlier. For posets, $\phi^\calP_+$ is the strict Ehrhart polynomial, and $\phi^\calP_-$ is the non-strict Ehrhart polynomial.
A more specific example is given by the character $\epsilon$ of $\exch$, see (\ref{EQ1}). The associated polynomial invariant is denoted by $\Phi$ (Proposition \ref{prop3.4.3}): for any any $I\in \exch[X]$,
\[\Phi(I)=\sum_{k=0}^{|X|} |I\cap \comp_k[X]|\frac{T(T-1)\cdots (T-k+1)}{k!}.\]
For any CTB $\calB$, Proposition \ref{prop3.5.4} gives that
\begin{align*}
\phi_\pm^\calB&=\Phi\circ \Gamma_\pm^\calB.
\end{align*}
This implies the unexpected result that there only a finite number of possible polynomials $\phi_+^\calB(b)$ of a given degree (Remark \ref{rk3.5.3}). For example, if $|X|=3$, the set of obtainable $\phi^\calB_+(b)$ is
\[\left\{\begin{array}{c}
T^3,\:(T-1)T^2,\:(T-1)^2T,\:(T+1)(T-1)T,\:(T-2)(T-1)T,\:\\[3mm]
\dfrac{(T-1)^2T}{2},\:\dfrac{(T-1)T^2}{2},\:\dfrac{(T-2)(T-1)T}{2},\\[3mm]
\:\dfrac{(2T-1)(T-1)T}{2},\:\dfrac{(2T-3)(T-1)T}{2},\:\dfrac{(2T+1)(T-1)T}{2},\:\\[3mm]
\dfrac{(T-2)(T-1)T}{3},\:\dfrac{2(T-2)(T-1)T}{3},\:\dfrac{(2T-1)(T-1)T}{3},\:\\[3mm]
\dfrac{(T-2)(T-1)T}{6},\:\dfrac{(2T-1)(T-1)T}{6},\:\dfrac{(4T+1)(T-1)T}{6},\:\\[3mm]
\dfrac{(4T-5)(T-1)T}{6},\:\dfrac{5(T-2)(T-1)T}{6},\:\dfrac{(5T-1)(T-1)T}{6},\:\\[3mm]
\dfrac{(5T+2)(T-1)T}{6},\:\dfrac{(5T-7)(T-1)T}{6},\:\dfrac{(5T-4)(T-1)T}{6}
\end{array}\right\}.\]
A similar result holds for $\phi_-$. For example, the set of obtainable $\phi^\calB_-(b)$ on a set of cardinality $3$ is
\[\left\{\begin{array}{c}
T,\:T^2,\:(2T-1)T,\:(3T-2)T,\: \dfrac{(T+1)T}{2},\:\dfrac{(3T-1)T}{2},\:\dfrac{(5T-3)T}{2},\\[3mm]
T^3,\:(T^2+2T-1)T,\:\dfrac{(T+1)T^2}{2},\:\dfrac{(T^2+3T-2)T}{2},\\[3mm]
\dfrac{(T^2+6T-4)T}{3},\:\dfrac{(2T-1)(T+2)T}{3},\:\dfrac{(T^2+3T-1)T}{3},\\[3mm]
\dfrac{(T+2)(T+1)T}{6},\:\dfrac{(T^2+12T-7)T}{6},\:\dfrac{(T^2+15T-10)T}{6},\\[3mm]
\dfrac{(T^2+6T-1)T}{6},\:\dfrac{(T^2+9T-4)T}{6},\:\dfrac{(2T+1)(T+1)T}{6},\\[3mm]
\dfrac{(2T-1)(T+5)T}{6},\:\dfrac{(4T-1)(T+1)T}{6},\:\dfrac{(5T-2)(T+1)}{6}
\end{array}\right\}.\]
If $|X|=n$ and $b\in \calB[X]$, then $\phi^\calB_+(b)$ and $\phi^\calB_-(b)$ are of degree $\leq n$ and share the same monomial $\phi^\calB_{\pp}(b)$ of degree $n$, which itself is a polynomial invariant, related to the CTB $\exch_{\pp}$ (Propositions \ref{prop3.5.1} and \ref{prop3.5.8}).
The leading coefficient $f^\calB_{\pp}(b)$ of $\phi^\calB_+(b)$ or $\phi^\calB_-(b)$ is studied in Proposition \ref{prop3.5.17}, where we details its forbidden values. For example:
\begin{itemize}
\item For $|X|=4$, $4!f^\calB(b)\in \{1,\ldots,22,\:24\}$.
\item For $|X|=5$, $5!f^\calB(b)\in \{1,\ldots,114,\:116,\:118,\:120\}$.
\item For $|X|=6$, $\displaystyle 6!f_{\pp}^\calB(b)\in \{1,\ldots,720\}\setminus \left\{\begin{array}{c}
695, 697, 699, 701, 703, 705, 706, 707, 709, \\
710, 711, 712, 713, 715, 716, 717, 718, 719
\end{array}\right\}$.
\end{itemize}

\section*{Twisted double bialgebras}

Before introducing combinatorial twisted double bialgebras, let us give some reminders on classical double bialgebras.
Let $(B,m,\delta)$ be a bialgebra. The category of right comodules is a symmetric monoidal category:
if $(M,\rho)$ and $(M',\rho')$ are two right comodules over $B$, then $M\otimes M'$ is also a right comodule over $B$, with the coaction defined by
\[m_{1,3,24}\circ (\rho \otimes \rho'),\]
where
\[m_{1,3,24}:\left\{\begin{array}{rcl}
M\otimes B\otimes M'\otimes B&\longrightarrow&M\otimes M'\otimes B\\
x\otimes a\otimes x'\otimes b&\longmapsto&x\otimes x'\otimes ab.
\end{array}\right.\]
Therefore, it is possible to define bialgebras in the category of right comodules over $B$, and this is precisely what a pair of cointeracting bialgebras is. More precisely, it is given by two bialgebras $(A,m,\Delta)$ and $(B,m,\delta)$, and a map $\rho:A\longrightarrow A\otimes B$, such that:
\begin{itemize}
\item $\rho$ is a right coaction of $B$ over $A$.
\item The product $m$ of $A$ is a right comodule morphism:
\[\rho\circ m=m_{12,34}\circ (\rho\otimes \rho).\]
The unit $1_A$ of $A$ is a right comodule morphism, that is to say $\rho(1_A)=1_A\otimes 1_B$.
This item is equivalent to the fact that $\rho$ is an algebra morphism from $A$ to $A\otimes B$.
\item The coproduct $\Delta$ of $A$ is a right comodule morphism:
\[(\Delta\otimes \id_B)\circ \rho=m_{1,3,24}\circ (\rho\otimes \rho)\circ \Delta.\]
\item The counit $\varepsilon_\Delta$ of $A$ is a right comodule morphism: 
\begin{align*}
&\forall a\in A,&(\varepsilon_\Delta\otimes \id_B)\circ \rho(a)&=\varepsilon_\Delta(a)1_B.
\end{align*}\end{itemize}
A double bialgebra is a pair of cointeracting bialgebras $(A,m,\Delta)$, and $(B,m,\delta)$, with the coaction $\rho$, such that $(A,m)=(B,m)$, and $\rho=\delta$. Except from the base field $\K$, the simplest case of double bialgebra is the polynomial algebra $\K[T]$, with its two multiplicative coproducts defined by
\begin{align*}
\Delta(T)&=T\otimes 1+1\otimes T,&\delta(T)&=T\otimes T.
\end{align*}
Various combinatorial examples are known. It seems that the first one is due to Calaque, Ebrahimi-Fard and Manchon \cite{Calaque2011}, where a second coproduct $\delta$ on the Butcher-Connes-Kreimer Hopf algebra with the help of contraction of edges, making it double.
Other examples can be found on graphs \cite{Foissy36}, oriented graphs \cite{Manchon2012}, hypergraphs \cite{Ebrahimi-Fard2022,Foissy44}, quasi-posets or finite topologies \cite{Foissy34,Foissy37}, noncrossing partitions \cite{Foissy38,Foissy43}, etc.

Let us now give a review on theoretical results on double bialgebras, before turning to their twisted version. Most of these results come from \cite{Foissy40}.
Let $(B,m,\Delta,\delta)$ be a double bialgebra. The set $\Char(B)$ of characters of $B$ inherits two products $\star$ and $*$, dual respectively to $\Delta$ and $\delta$. The cointeraction implies that for any $f$, $g$, $h\in \Char(B)$,
\[(f\star g)*h=(f*h)\star(g*h).\]
Let $(B',m',\Delta')$ be another bialgebra. We denote by $\Hom(B,B')$ the set of bialgebra maps from $(B,m,\Delta)$ to $(B',m',\Delta')$. Then the monoid $(\Char(B),\star)$ acts (on the right) on $\Hom(B,B')$ as follows:
\begin{align*}
\curvearrowleft&:\left\{\begin{array}{rcl}
\Hom(B,B')\times\Char(B)&\longrightarrow&\Hom(B,B')\\
(\phi,f)&\longmapsto&\phi\curvearrowleft f=(\phi\otimes f)\circ \delta.
\end{array}\right.
\end{align*}
In the case where $B'=(\K[X],m,\Delta)$ (case of polynomial invariants), under a condition of connectedness, $\Hom(B,\K[X])$ is a one-to-one correspondence with $\Char(B)$, through the map
\begin{align*}
\theta_B&:\left\{\begin{array}{rcl}
\Hom(B,B')&\longrightarrow&\Char(B)\\
\phi&\longmapsto&\left\{\begin{array}{rcl}
B&\longrightarrow&\K\\
x&\longmapsto&\phi(x)(1).
\end{array}\right.
\end{array}\right.
\end{align*}
The polynomial invariant $\phi_+$ associated to the counit $\epsilon_\delta$ of the coproduct $\delta$ is called the fundamental polynomial invariant. It is the unique one which is also a bialgebra morphism from $(B,m,\delta)$ to $(\K[X],m,\delta)$. 
Moreover, the inverse of $\theta_B$ is given by
\begin{align*}
\theta_B^{-1}&:\left\{\begin{array}{rcl}
\Char(B)&\longrightarrow&\Hom(B,B')\\
f&\longmapsto&\phi_+\curvearrowleft f.
\end{array}\right.
\end{align*}
Still under the condition of connectedness, $(B,m,\Delta)$ is a Hopf algebra, whose antipode is denoted by $S_B$. Considering the character of $B$ defined by
\begin{align*}
&\forall b\in B,&f(b)&=\phi_+(x)(-1),
\end{align*}
then for any $x\in B$,
\[S_B(x)=(f\otimes \id_B)\circ \delta(x).\]

Let us now describe the twisted framework we shall use. Let $(\calQ,m,\Delta)$ a a twisted bialgebra. For any finite set $X$, we denote by $\calE[X]$ the set of equivalences on $X$. 
A double structure on $\calQ$ is given by a map $\delta_\sim:\calQ[X]\longrightarrow \calQ[X/{\sim}]\otimes \calQ[X]$, satisfying a compatibility with the species structure, a coassociativity and a counit axiom, and compatibilities with $m$ and $\Delta$, 
which insures that the application of the Fock functor $\calF$ to $\calQ$ gives a double bialgebra. We won't detail these compatibilities here, see Definition \ref{defi1.3.1}. For example:
\begin{itemize}
\item Graphs. A coproduct $\delta$ is defined on $\calG$, by contraction of connected subgraphs, see Example \ref{ex4.1.1}. Summing over all equivalences, we have
\begin{align*}
\delta\left(\hspace{-2mm} \begin{array}{c}
\tdun{$x$}
\end{array}\hspace{-2mm}\right)&=\tdun{$x$}\otimes \tdun{$x$},\\
\delta\left(\hspace{-2mm} \begin{array}{c}
\tddeux{$x$}{$y$}
\end{array}\hspace{-2mm}\right)&=\tddeux{$x$}{$y$}\otimes \tdun{$x$}\tdun{$y$}+\tddun{$xy$}\otimes \tddeux{$x$}{$y$},\\
\delta\left(\hspace{-2mm} \begin{array}{c}
\tdtroisun{$x$}{$z$}{$y$}
\end{array}\hspace{-2mm}\right)&=\tdtroisun{$x$}{$z$}{$y$}\otimes \tdun{$x$}\tdun{$y$}\tdun{$z$}+\tdddeux{$xy$}{$z$}\otimes\tddeux{$x$}{$y$}\tdun{$z$}+\tdddeux{$xz$}{$y$}\otimes\tddeux{$x$}{$z$}\tdun{$y$}+\tdddun{$xyz$}\otimes\tdtroisun{$x$}{$z$}{$y$},\\
\delta\left(\hspace{-2mm} \begin{array}{c}
\grdtroisun{$x$}{$z$}{$y$}
\end{array}\hspace{-2mm}\right)&=\grdtroisun{$x$}{$z$}{$y$}\otimes \tdun{$x$}\tdun{$y$}\tdun{$z$}+\tdddeux{$xy$}{$z$}\otimes\tddeux{$x$}{$y$}\tdun{$z$}+\tdddeux{$xz$}{$y$}\otimes\tddeux{$x$}{$z$}\tdun{$y$}
+\tdddeux{$yz$}{$x$}\otimes\tddeux{$y$}{$z$}\tdun{$x$}+\tdddun{$xyz$}\otimes\grdtroisun{$x$}{$z$}{$y$}.
\end{align*}
\item Posets. A coproduct $\delta$ is defined on $\calP$, by contraction of connected and convex sub-posets, see Example \ref{ex4.1.2}. For example,
\begin{align*}
\delta\left(\hspace{-2mm} \begin{array}{c}
\tdun{$x$}
\end{array}\hspace{-2mm}\right)&=\tdun{$x$}\otimes \tdun{$x$},\\
\delta\left(\hspace{-2mm} \begin{array}{c}
\tddeux{$x$}{$y$}
\end{array}\hspace{-2mm}\right)&=\tddeux{$x$}{$y$}\otimes \tdun{$x$}\tdun{$y$}+\tddun{$xy$}\otimes \tddeux{$x$}{$y$},\\
\delta\left(\hspace{-2mm} \begin{array}{c}
\tdtroisun{$x$}{$z$}{$y$}
\end{array}\hspace{-2mm}\right)&=\tdtroisun{$x$}{$z$}{$y$}\otimes \tdun{$x$}\tdun{$y$}\tdun{$z$}+\tdddeux{$xy$}{$z$}\otimes\tddeux{$x$}{$y$}\tdun{$z$}+\tdddeux{$xz$}{$y$}\otimes\tddeux{$x$}{$z$}\tdun{$y$}+\tdddun{$xyz$}\otimes\tdtroisun{$x$}{$z$}{$y$},\\
\delta\left(\hspace{-2mm} \begin{array}{c}
\tdtroisdeux{$x$}{$y$}{$z$}
\end{array}\hspace{-2mm}\right)&=\tdtroisdeux{$x$}{$y$}{$z$}\otimes \tdun{$x$}\tdun{$y$}\tdun{$z$}+\tdddeux{$xy$}{$z$}\otimes\tddeux{$x$}{$y$}\tdun{$z$}+\tddddeux{$x$}{$yz$}\otimes\tddeux{$y$}{$z$}\tdun{$x$}+\tdddun{$xyz$}\otimes\tdtroisdeux{$x$}{$y$}{$z$}.
\end{align*}\end{itemize}

Let $\calQ$ be a double bialgebra. We show in Theorem \ref{theo1.3.6} that there exists a unique polynomial invariant of $\calQ$, compatible with the coproduct $\delta$. 
It will be called the fundamental polynomial invariant. For graphs, this is the chromatic polynomial, and for posets the strict Ehrhart polynomial.
Again, the action of the monoid of characters of $\calQ$ (dual to the coproduct $\delta)$ acts freely on the set of polynomial invariants of $\calB$, with as a single generator the universal polynomial invariant (Proposition \ref{prop1.3.7}). 
As in the classical setting, the antipode of $(\calQ,m,\Delta)$ is obtained from the universal polynomial invariant and $\delta$, through an evaluation in $-1$ (Theorem \ref{theo1.3.10}).

Let us now turn to our combinatorial frame. A double structure on a CTB $\calB$ is a coproduct $\delta$ making a twisted double bialgebra, with the combinatorial extra condition that for any finite set $X$, for any $b\in \calB[X]$, and any $\sim\in \calE[X]$,
$\delta_\sim(b)$ is a tensor of an element of $\calB[X/{\sim}]$ with an element of $\calB[X]$, or is zero (Definition \ref{defi4.1.1}). Such an object will be called a combinatorial twisted double bialgebra (briefly, CTDB). For example, $\calG$ and $\calP$ are CTDBs.
An important question is to know if a given CTB can be made into a CTDB by the addition of a coproduct $\delta$, and if such a $\delta$ is unique. We give partial answers to these questions. 
Firstly, if $\calB$ is a CTDB, we associate to any $b\in \calB$ a set of equivalences $\calE[b]$, whose elements are the equivalences $\sim$ such that $\delta_\sim(b)\neq 0$. 
These sets are partially ordered by the inclusion, and we show that they have a unique minimal element $\sim_b^{(m)}$ (always the equality of the underlying set), and a unique maximal element $\sim_b^{(M)}$, included in $\sim_b^{(\ic)}$ (Proposition \ref{prop4.1.7}).
This leads to a sub-CTDB $\bfC\calB$ of $\calB$, made of elements $b$ such that $\sim_b^{(m)}=\sim_b^{(M)}$. We show that this twisted subalgebra of $\calB$ is generated by elements which are primitive for $\Delta$, and group-like for $\delta$ (Proposition \ref{prop4.1.14}).
In particular, it contains $\bfD\calB$. A particular class of CTDBs is given by the ones for which $\bfC\calB=\bfD\calB$: these are called regular, and we give in Proposition \ref{prop4.1.15} several characterizations of regularity. 
The CTDB of graphs $\calG$ and of posets $\calP$ are regular, and more generally, if $(\calB,m,\Delta,\delta)$ is a CTDB such that $(\calB,m,\Delta)$ is non-degenerate, then $\calB$ is regular (Proposition \ref{prop4.1.16}).
The irregularity of a CTDB has to be understood as a defect of graduation, in some sense, and we explain how to reduce ourselves to the regular case. For a regular CTDB $\calB$, the fundamental polynomial invariant is $\phi_+^\calB$, defined earlier.
In fact, regular CTDBs are closed to operads (see \cite{Dotsenko2016,Loday2012} for details on these objects): the coproduct $\delta$ of a regular CTDB induces an operadic structure on the indecomposable elements, see Theorem \ref{theo4.1.22}.
However, we won't explore this further in this text.

\section*{Quasi-cocommutativity}

In order to answer the question if a given CTB can be made double, we introduce the technical notion of quasi-cocommutativity (Definition \ref{defi4.2.1}). In rough words, a CTB $\calB$ is quasi-cocommutative if it is commutative,
and if the coproduct is given by a restriction operation on the underlying combinatorial objects, at least when the factors in the coproduct are indecomposable. We prove that non-degenerate CTDBs are quasi-cocommutative, which implies that some CTBs cannot be made double.
This is the case for $\exch$, $\exch_+$ and $\exch_-$, and also for one of the CTB on hypergraphs, whose fundamental polynomial invariant is the chromatic polynomial of \cite{Aval2019}, see Example \ref{ex4.2.3}.
The question is now to determine if a quasi-cocommutative CTB can be made double. We firstly study, for a quasi-cocommutative CTB $\calB$, the structure of the equivalences associated to any element $b$ of $\calB$ in the case it is made double.
We show that these sets $\calE[b]$ are entirely defined by the decomposition into indecomposable elements and the coproduct $\Delta$ (Definition \ref{defi4.2.5}), as well as the right factor $b\mid\sim$ of $\delta_\sim(b)$ if $\sim\in \calE[b]$.
There is a little bit more latitude for the left factor $b/{\sim}$, but not so much, as $\calC_+(b,\sim)=\calC_+(b/{\sim})$ is also entirely determined by the decomposition into indecomposable elements and $\Delta$ (Definition \ref{defi4.2.6}).
These characterizations are enough to construct a coaction-like map $\Xi:\calB\longrightarrow \exch_+\otimes \calB$ (Theorem \ref{theo4.2.10}) for any quasi-cocommutative CTB, which is $(\Gamma_+^\calB\otimes \id_{\calB})\circ \delta$ if $\delta$ makes $\calB$ double.
It turns out this map $\Xi$ is enough to generalize the formula for the antipode for CTDBs to quasi-cocommutative CTBs (Corollary \ref{cor4.2.13}), and to show that the polynomial invariant $\phi_+^\calB$ is compatible with $\Xi$ (Proposition \ref{prop4.2.14}),
as well as it is compatible with $\delta$ in the case of regular CTDBs. Moreover, a formula for the coefficients of the polynomial invariant $\phi^\calB_+$ can also be given with the help of the map $\Xi$ (Theorem \ref{theo4.2.17}).

Of course, there exist non quasi-cocommutative CTBs, but we show that any commutative CTB $\calB$ contains a quasi-cocommutative sub-CTB $\calB^{\qc}$, maximal for the inclusion (Theorem \ref{theo4.2.18}), 
and that this construction is functorial in the category of commutative CTBs (Proposition \ref{prop4.2.19}). We already mentioned that $\exch_+$ is not quasi-cocommutative, so cannot be made double. 
A good candidate to be made double would be its strict sub-CTB $\exch_+^{\qc}$. However, this is still too big, and we show that there exists a sub-CTB $\exch_+^{\double}$ of $\exch_+^{\qc}$, maximal for the inclusion, which can be made double.
Moreover, the coproduct $\delta$ making it double is unique (Theorem \ref{theo4.3.1}). 
A consequence is that if $\calB$ is a CTB such that $\Gamma_+^\calB$ is injective, then it can be made double in at most one way (Corollary \ref{cor4.3.2}): this is the case for the CTB of graphs $\calG$ and of posets $\calP$, whose coproduct $\delta$ is consequently unique (Example \ref{ex4.3.1}). 
For a given CTDB $\calB$, $\Gamma_+^\calB$ is the unique candidate for a morphism of CTDBs from $\calB$ to $\exch_+^{\qc}$, but it is not necessarily the case.
We give in Theorem \ref{theo4.3.5} a necessary and sufficient condition for this, in terms of conservation of the number of indecomposable components. This condition is always satisfied when $\Gamma_+^\calB$ is injective (Corollary \ref{cor4.3.6}).

An example of quasi-cocommutative CTB based on trees is studied in the sixth chapter. 
It is related to the Butcher-Connes-Kreimer Hopf algebra \cite{Connes1998}, and is used in \cite{Bruned2015} for the study of stochastics PDEs and in \cite{Marcolli2025} for linguistics, see \cite{Xiao2026}.
At the level of classical Hopf algebras, the Butcher-Connes-Kreimer coproduct is given by admissible cuts. For example,
\begin{align*}
\Delta^{CK}\left(\hspace{-2mm}\begin{array}{c}
\ttroisun
\end{array}\hspace{-2mm}\right)&=\ttroisun\otimes 1+\tun\otimes\tun\tun+2\tdeux\otimes\tun+1\otimes \ttroisun,\\
\Delta^{CK}\left(\hspace{-2mm}\begin{array}{c}
\ttroisdeux
\end{array}\hspace{-2mm}\right)&=\ttroisdeux\otimes 1+\tdeux\otimes\tun+\tun\otimes\tdeux+1\otimes \ttroisdeux.
\end{align*} 
It can be easily be lifted to the level of CTB, and gives a sub-CTB of $\calP$. For the coproduct that we study in Chapter six, we add a condition of conservation of the number of leaves, which leads us to deleted some terms in the coproduct:
\begin{align*}
\Delta^{L}\left(\hspace{-2mm}\begin{array}{c}
\ttroisun
\end{array}\hspace{-2mm}\right)&=\ttroisun\otimes 1+2\tdeux\otimes\tun+1\otimes \ttroisun,\\
\Delta^{L}\left(\hspace{-2mm}\begin{array}{c}
\ttroisdeux
\end{array}\hspace{-2mm}\right)&=\ttroisdeux\otimes 1+1\otimes \ttroisdeux.
\end{align*} 
This coproduct $\Delta^L$ conserves the number of leaves and the number of vertices. It induces a quotient to reduced rooted trees, that is to say to trees such that no vertex has exactly one child:
\begin{align*}
\Delta^{\ell}\left(\hspace{-2mm}\begin{array}{c}
\ttroisun
\end{array}\hspace{-2mm}\right)&=\ttroisun\otimes 1+2\tun\otimes\tun+1\otimes \ttroisun,\\
\Delta^{\ell}\left(\hspace{-2mm}\begin{array}{c}
\tun
\end{array}\hspace{-2mm}\right)&=\tun\otimes 1+1\otimes \tun.
\end{align*} 
This coproduct $\Delta^\ell$ conserves the number of leaves, but not of vertices.

The coproduct $\Delta^L$ can be lifted to the CTB level into two ways, each one with its qualities and its defects.
The first one is based on vertices: for any finite set $X$, we consider the set of rooted forests $\calR\calF[X]$ whose set of vertices is $X$, and this defines a species $\calR\calF$ (Definition \ref{defi6.1.1}).
The coproduct $\Delta^L$ is lifted to $\calR\calF$, making it a quasi-cocommutative CTB (Theorem \ref{theo6.1.3}).
This allows us to describe easily the antipode of $\calR\calF$ in Corollary \ref{cor6.1.7}. 
However, because $\Delta^\ell$ does not conserve the number of vertices, it is not possible to obtain a twisted version of $\Delta^\ell$ from it, unfortunately. For this, we introduce a second CTB, based on leaves: 
for any finite set $X$, $\calR\calF_\ell[X]$ is the set of isoclasses of rooted forests whose set of leaves is $X$ (Definition \ref{defi6.2.1}).
This defines a CTB $\calR\calF_\ell$ (Theorem \ref{theo6.2.2}). However, it is not quasi-cocommutative. 
The advantage of this construction is the existence of a morphism of CTBs from $\calR\calF_\ell$ to one of non quasi-cocommutative CTB of hypergraphs (Theorem \ref{theo6.2.3}). This morphism is non injective, and induces consequently a quotient of $\calR\calF_\ell$ embedded in hypergraphs:
this quotient is the species of reduced rooted forests $\calR\calR\calF_\ell$, and the induced coproduct is a CTB version of $\Delta^\ell$ (Theorem \ref{theo6.2.8}). 
We show that the Fock functor applied to $\calR\calF$ and to $\calR\calF_\ell$ both gives the same object, precisely the bialgebra of rooted trees with the coproduct $\Delta^L$ (Theorem \ref{theo6.3.1}), 
and the quasi-cocommutativity of $\Delta^L$ gives a cancellation-free formula for its antipode.
Combining both constructions, we obtain a cancellation-free formula for the (more complicated) antipode of $\Delta^\ell$ on reduced forests (Proposition \ref{prop6.3.4}), with a combinatorial interpretation of its coefficients in terms of partitions in trees, which answers a question in \cite{Xiao2026}.

 \section*{Structures related to the Hadamard product}

The category of species has a second monoidal structure, the Hadamard product $\boxtimes$: 
if $\calP$ and $\calQ$ are two species, and $X$ is a finite set,
\[\calP\boxtimes \calQ[X]=\calP[X]\otimes \calQ[X].\]
We explore in the fifth chapter some applications of this tensor product. A first, immediate application is that the Hadamard product of two CTBs is also a CTB (Proposition \ref{prop5.1.2}).
This allows to construct new CTBs from old ones: for example, the Hopf algebra of mixed graphs of \cite{Foissy45} can be obtained by the Hadamard product of the CTB of graphs and the CTB of oriented graphs. Unfortunately, this does not work for CTDBs.

We define an internal product on a CTB $\calB$ as an associative map from $\calB\boxtimes\calB$ to $\calB$, compatible with the coproduct $\Delta$ and the product $m$ of $\calB$ (Definition \ref{defi5.1.4}).
For example, the intersection is an internal product on $\exch_-$ and $\exch_+$ (Proposition \ref{prop5.1.5}); The union of sets of edges or hyperedges is an internal product on graphs or hypergraphs (Proposition \ref{prop5.1.5}); the transitive closure of the union is an internal product on quasi-posets (Proposition \ref{prop5.1.8}). 
It is noticeable that none of these internal products subsist at the application of the bosonic Fock functor. Nevertheless, they are compatible with the full Fock functor.

Another application, probably less anecdotal than the internal product, is given by functors of colorations. If $\Omega$ is a set, we define a CTB $\calD_\Omega$: for any finite set $X$, $\calD_\Omega[X]$ is the set of maps from $X$ to $\Omega$, what we will call $\Omega$-colorations. 
Its product is given by disjoint union, and the coproduct by restriction. 
For any CTB $\calB$, $\calB\boxtimes \calD_\Omega$ is the $\Omega$-colored version of $\calB$, and $\calF_\Omega[\calB]=\calF[\calB\boxtimes \calD_\Omega]$ has for basis the set of isoclasses of $\Omega$-colored elements of $\calB$ (Proposition \ref{prop5.2.2}).
For example, when applied to the CTB of rooted trees, we obtain the colored versions of the Butcher-Connes-Kreimer Hopf algebra \cite{Connes1998}, where vertices are decorated by elements of the set $\Omega$.
When $\calB$ is a CTDB and $\Omega$ is a commutative semigroup, we prove that $\calF_\Omega[\calB]$ is also a CTDB (Proposition \ref{prop5.2.5}).
This is used to define the formal flow associated to a CTB $\calB$, following what is done for Novikov algebras in \cite{Bandiera2025}. A (finite) totally ordered set $\Omega=\{\omega_1,\ldots,\omega_n\}$ being fixed, consider in $\calF_\Omega[\calB]$ a family of infinitesimal characters $(\iota_\omega)_{\omega\in \Omega}$ indexed by $\Omega$, which acts only on elements of degree $1$, according to their coloration.
Then consider in the group of characters of $\calF_\Omega[\calB]$, the element 
\[\exp(\iota_{\omega_1})*\cdots*\exp(\iota_{\omega_n}).\]
This is also a character of $\calF_\Omega[\calB]$, so can be uniquely written as the exponentiation of an infinitesimal character 
\[f_\Omega=\ln\left(\exp(\iota_{\omega_1})*\cdots*\exp(\iota_{\omega_n})\right).\]
As we see, this is closely related to the Baker-Campbell-Hausdorff formula in the exponentiation of the Lie algebra of infinitesimal characters of $\calF_\Omega[\calB]$.
The polynomial invariant associated to the character $\exp(\iota_{\omega_1})*\cdots*\exp(\iota_{\omega_n})$ is the formal flow $\phi_{f_\Omega}$ of $\calF_\Omega[\calB]$ (Definition \ref{defi5.2.7}), and $f$ is obtained as
\[f_\Omega=\dfrac{\mathrm{d}\phi_{f_\Omega}}{\mathrm{d}T}(0).\]
We give a combinatorial description of the formal flow in Theorem \ref{theo5.2.9}. For any pair $(b,\ell)$ of $\Omega$-decorated element of $\calB$, $\phi_{f_\Omega}(b,\ell)$ counts (up to a factorial) the number of $\calB$-partitions of $(b,\ell)$, 
as given in Definition \ref{defi5.2.6}. In the case of posets, we recover the notion of $P$-partitions \cite{Malvenuto1993}.

\section*{Functorial sub-objects of a CTB $\calB$}

We introduce, for a given CTB $\calB$, several sub-CTBs, correcting a defect of $\calB$, each of these constructions being functorial. The first one is $\cocom\calB$, the greatest cocommutative sub-CTB (Proposition \ref{prop5.3.2}).
In the case of a non-degenerate CTDB, $\cocom\calB$ is a sub-CTDB.
We give a description of $\cocom\exch$ in terms of invariance under the action of symmetric groups (Definition \ref{defisigmainvariance}), which leads us to a description in terms of subsets of equivalences. 
A more convenient description of $\cocom\exch_+$ is lead in Chapter \ref{chap7}, through a family of hypergraphs called antichain covers. 

The second one is linked with the reciprocity principle of posets. As a well-known property of Ehrhart polynomials (which we prove here in Corollary \ref{cor5.3.5} in an algebraic way), for any poset $P$,
\[\phi_+^\calP(P)(-T)=(-1)^{|V(P)|}\phi^\calP_-(P)(T).\]
We shall say that the CTB $\calP$ of posets satisfy the reciprocity principle.
We show that any $\calB$ has a sub-CTB $\calB^{\mathrm{recip}}$, maximal for the inclusion, satisfying the reciprocity principle (Proposition \ref{prop5.3.6}).
By functoriality, all the polynomial invariants with the reciprocity principle factorize through $\exch_+^{\mathrm{recip}}$,
via the morphism $\Gamma_+^\calB$. This CTB $\exch_+^{\mathrm{recip}}$ is rather mysterious; we show that it is not restricted to the image of posets, with two examples of elements which do not correspond to an Ehrhart polynomial (Examples \ref{ex5.3.2} and \ref{ex5.3.3}).

A classical result on graphs is the property of alternating sign of their chromatic polynomial: if $G$ is a graph with $n$ vertices, 
we write its chromatic polynomial as
\[\phi_+^\calG+(G)=a_nX^n+\cdots+a_1X+a_0.\]
Then $(-1)^{n+k}a_k\geq 0$ for any $k\in \{0,\ldots,n\}$. This is classically shown by an induction on the number of edges of $G$, with the help of the deletion-contraction formula for the chromatic polynomial \cite{Rota1964}.
We shall say that $\calG$ satisfies the alternating sign condition (Definition \ref{defi5.3.8}).
We show that any CTB $\calB$ has a sub-CTB $\calB^{\mathrm{alt}}$, maximal for the inclusion, satisfying the alternating sign condition (Proposition \ref{prop5.3.9}). When $\calB$ is a cocommutative non-degenerate CTB, 
we show that $\calB$ satisfies the alternating condition if, and only if, $\phi_+^\calB$ factorizes through $\calG$, via a specific morphism $\theta_G^\calB$ of CTBs (Theorem \ref{thoe7.2.30}):
in other words, in the cocommutative case, the fundamental polynomial invariants with the alternating sign conditions are all chromatic polynomials of graphs.
We also show that the other fundamental polynomial invariants will have at least a zero coefficient between two non-zero coefficients.

\section*{More results for cocommutative CTBs}

The seventh chapter is devoted to the simplification of the preceding results in the cocommutative and non-degenerate case.
The first simplification is for the coproduct, which in this case is necessarily given by a restriction operation (Proposition \ref{prop7.1.1}):
If $X\subseteq Y$ are two finite sets, we have a map
\[\left\{\begin{array}{rcl}
\calB[Y]&\longrightarrow&\calB[X]\\
b&\longmapsto&b_{\mid X},
\end{array}\right.\]
compatible with the species structure, such that if $b\in \calB[Y]$,
\begin{align*}
\Delta_{X,Y\setminus X}(b)&=b_{\mid X}\otimes b_{\mid Y\setminus X}&\mbox{or}&&\Delta_{X,Y\setminus X}(b)&=0.
\end{align*}
This implies that if $X\subseteq Y\subseteq Z$ are finite sets, and $b\in \calB[Z]$,
\[\left(b_{\mid Y}\right)_{\mid X}=b_{\mid X}.\]
Moreover, the polynomial invariants $\phi_+^\calB$ and $\phi_-^\calB$ take their values in $\Z[T]$. By the way, in this case $\phi_-^\calB$ does not carry much information, as for any $b\in \calB[X]$, $\phi_-(b)=T^{|X|}$ (Proposition \ref{prop7.1.1}).
We have a partial dual implication for a non-degenerate CTDB $\calB$: in this case, if $\phi_+^\calB$ takes its values in $\Z[T]$, then $\calB$ is cocommutative (Theorem \ref{theo7.1.2}). This is not true if $\calB$ is not double, however.

We next show how to replace the CTB $\exch$, rather complicated, by the simpler CTB of hypergraphs $\calH\calG$. 
Firstly, we show how to associate to any convenient subspecies $\calC$ of $\calB$ (that is to say, a subspecies of the species of indecomposable elements of $\calB$) a morphism $\theta_\calC^\calB$ of CTB from $\calB$ to $\calH\calG$:
the hyperedges of the hypergraph $\theta_\calC^\calB(b)$ associated to an element $b\in \calB[X]$ are the subsets $Y$ of $X$ such that $b_{\mid Y}\in \calC[Y]$. 
We will be interested in two particular convenient subspecies: the first one is $\calC_{\CS}$, the subspecies of all indecomposable elements, the second one is $\calC_{\AC}$, the subspecies of minimal indecomposable elements.
This leads us to two sub-CTBs of $\calH\calG$, which we now describe.

A hypergraph $H$ is a connectedness system \cite[Entry A072446]{Sloane} if all singletons are hyperedges of $G$, and if the union of two intersecting hyperedges of $H$ is also a hyperedge of $H$.
A classical example of connectedness system is attached to a graph $G$: the corresponding hyperedges are the subsets $Y$ of $V(G)$ such that $G_{\mid Y}$ is connected -- hence, the name connectedness system. 
Connectedness systems defines a sub-CTB of $\calH\calG$ (Theorem \ref{theo7.2.7}). This is generalized as follows.
Let $\calB$ be a non-degenerate cocommutative CTB. The morphism $\theta_{\CS}^\calB$ associated to the convenient subspecies $\calC_{\CS}$ sends any $b\in \calB[X]$ to the hypergraph whose hyperedges are the subsets $Y$ of $X$ such that $b_{\mid Y}$ is indecomposable.
We prove that for any the morphism $\theta_{\CS}^\calB$ takes its values in $\calC\calS$, and that it is the unique morphism from $\calB$ to $\calC\calS$ that conserves the number of indecomposable components (Theorem \ref{theo7.2.8}).
This morphism $\theta_{\CS}^\calB$ is used to give a cancellation-free formula for the antipode (Proposition \ref{prop7.2.29}), with coefficients whose combinatorial interpretation (even of their sign, in fact) may be difficult.

A hypergraph $H$ is an antichain cover \cite[Entry A006126]{Sloane} if no singleton is a hyperedge of $H$ and if $Y\subseteq Z$ are two hyperedges of $H$, then $Y=Z$.
Antichain covers define a sub-CTB of $\calH\calG$, denoted by $\calA\calC$ (Theorem \ref{theo7.2.15}).
Let $\calB$ be a non-degenerate cocommutative CTB $\calB$. The morphism $\theta_{\AC}^\calB$ associated to the convenient subspecies $\calC_{\AC}$ sends any $b\in \calB[X]$ to the hypergraph whose hyperedges are the subsets $Y$ of $X$ such that $b_{\mid Y}$ is indecomposable, and minimal for this property. We prove in Theorem \ref{theo7.2.16} that 
$\theta_{\AC}^\calB$ takes its values in $\calA\calC$. Moreover,
\begin{align*}
\epsilon^{\calA\calC}_+\circ \theta^\calB_{\AC}&=\epsilon^\calB_+,&
\Gamma^{\calA\calC}_+\circ \theta^\calB_{\AC}&=\Gamma^\calB_+,&
\phi^{\calA\calC}_+\circ\theta^\calB_{\AC}&=\phi^\calB_+.
\end{align*}
A consequence is that any fundamental polynomial invariant of a cocommutative CTB is in fact the chromatic polynomial of an antichain cover. Another consequence is that the greatest cocommutative sub-CTB $\cocom\exch_+$ of $\exch_+$ is isomorphic to $\calA\calC$, through the morphism $\Gamma_+^{\calA\calC}$ (Corollary \ref{cor7.2.20}).

For the double-bialgebraic aspects, it turns out that $\calC\calS$ is a sub-CTDB of $\calH\calG$ (Theorem \ref{theo7.2.7}), whereas $\calA\calC$ is not.
Moreover, for any non-degenerate cocommutative CTDB $\calB$, $\theta_{\calC\calS}^\calB$ is the unique morphism of CTDBs from $\calB$ to $\calC\calS$ (Theorem \ref{theo7.2.8}).
Unfortunately, $\Gamma_+^{\calC\calS}$ does not induce a morphism of CTDBs from $\calC\calS$ to $\exch_+^{\double}$, and we introduce a sub-CTDB of $\calC\calS$, denoted by $\calC\calS^{\double}$, with this property and maximal for the inclusion (Proposition \ref{prop7.2.12}).
For $\calA\calC$, which is not double, we prove the existence of a sub-CTB $\calA\calC^{\double}$, which can be made double and is maximal for the inclusion (Proposition \ref{prop7.2.24}). 
Moreover, the double structure on it is unique, and $\Gamma_+^{\calA\calC^{\double}}$ induces a morphism of CTDBs from $\calA\calC^{\double}$ to $\exch_+^{\double}$ (Proposition \ref{prop7.2.24}), whose image is the greatest cocommutative sub-CTDB of $\exch_+^{\double}$ (Corollary \ref{cor7.2.25}). 
We give several equality or inclusion results on these CTBs in Theorem \ref{theo7.2.27}. 

\section*{Bialgebraic structure of $\exch$}

The last chapter is devoted to the study of $\exch$ and its sub-objects $\exch_\pm$ and $\exch_{\pp}$, as twisted bialgebras. We show that all of them are, in a non-trivial way, isomorphic to shuffle twisted bialgebras. 
To prove this, we use a twisted version of Loday and Ronco's rigidity theorem for cofree bialgebras \cite{Loday2006}. We firstly introduce twisted infinitesimal bialgebras (Definition \ref{defi8.1.3}), and prove a rigidity theorem for them (Theorem \ref{theo8.1.5}),with the help of the antipode. 
We also prove Theorem \ref{theo8.1.9}, which establishes that any commutative cofree twisted bialgebra is isomorphic to a shuffle algebra, thanks to the eulerian idempotent. 
These results are used on $\exch$ and its three sub-objects $\exch_+$, $\exch_-$ and $\exch_{\pp}$, thanks to a new concatenation-like product $\diamond$ on $\exch$, $\exch_+$ and $\exch_{\pp}$ (Theorem \ref{theo8.1.11}), and to another product $\blackdiamond$ on $\exch_-$ (Theorem \ref{theo8.1.14}).

\chapter*{Notations}.

\begin{enumerate}
\item $\K$ is a commutative field of characteristic zero. All the vector spaces in this text will be taken over $\K$.
\item For $n\in \N$, we put $[n]=\{1,\ldots,n\}$. In particular, $[0]=\emptyset$.
\item For $n\in \N$, the $n$-th symmetric group is denoted by $\sym_n$.
\item For $V,W$ two $\K$-vector spaces, then $\tau_{V,W}$ is the usual flip
\[\tau_{V,W}:\left\{\begin{array}{rcl}
V\otimes W&\longrightarrow&W\otimes V\\
v\otimes w&\longmapsto&w\otimes v.
\end{array}\right.\]\end{enumerate}

\chapter{Reminders and notations}

\section*{Introduction}

This first chapter is devoted to reminders and complements on the objects which will be used in the sequel: species, twisted algebras, coalgebras, bialgebras, Hopf algebras, etc. We start with species. 
These objects were introduced by Joyal \cite{Joyal1981,Joyal1986} in order to have a better understanding of the behavior of generating formal series associated to a family of combinatorial objects.
For example, graphs, hypergraphs, posets, partitions, equivalences, compositions, and many more, are organized as species.
Species are functors from the category of finite sets with bijections, to the category of sets with maps (set species) or to the category of vector spaces over $\K$ with linear maps (linear species).
The category of linear species is a symmetric tensor category with the Cauchy tensor product, and this allows to define algebras, coalgebras, bialgebras, etc, in this category: these are what we call twisted algebras, coalgebras, bialgebras, etc \cite{Patras2004}.
We give the definitions in the first section of this chapter, as well as examples which will play an important role in the sequel: the species of equivalences $\calE$ and the species of compositions $\comp$, see Example \ref{ex1.1.1}.

We generalize well-known results from classical objects to twisted ones in the second section. We firstly introduce the dual $\calQ^*$ of a species $\calQ$ (Definition \ref{defi1.2.1}). 
When $\calQ$ is a twisted coalgebra, then its dual is a (classical) algebra with the convolution product (Proposition \ref{prop1.2.2}).
The natural graduation of a species induces a metric on $\calQ^*$, making it a complete metric space (Propositions \ref{prop1.2.3}, and \ref{prop1.2.4}).
This allows to define the action of the algebra $\K[[T]]$ of formal series on $\calQ^*$, and more generally on the space of species morphisms between two species (Theorem \ref{theo1.2.6}).
This allows for example to determine the inverse for the convolution product of an element of $\calQ^*$ (Proposition \ref{prop4.2.12}), and to define the exponential and the logarithm of certain elements of $\calQ^*$ (Proposition \ref{prop1.2.8}).
We then define characters and infinitesimal characters of a twisted bialgebra $\calQ$ (Definition \ref{defi1.2.9}), generalizing the equivalent notions for bialgebras.
When $\calQ$ is connected, that is to say if $\calQ[\emptyset]$ is 1-dimensional, the exponential and the logarithm define bijections between infinitesimal characters and characters, inverse one from the other (Corollary \ref{cor1.2.12}), as for classical connected bialgebras.
We also give the definition of twisted Hopf algebras and give the twisted version of Takeuchi's formula \cite{Takeuchi1971} (Definition \ref{defi1.2.13} and Proposition \ref{prop1.2.14}).
We introduce the twisted version of the eulerian idempotent \cite{Bandiera2017,Burgunder2008,Foissy40} in Definition \ref{defi1.2.15}.
If $\calQ$ is a commutative or cocommutative twisted bialgebra, then the eulerian idempotent is a projection (which justifies its name), on the space of primitive elements when $\calQ$ is cocommutative, 
vanishing on the square of the augmentation ideal when $\calQ$ is commutative (Theorem \ref{theo1.2.19}).
The proof use some combinatorial interpretations in term of packed matrices. 

We then introduce polynomial invariants of a twisted bialgebra $\calQ$ (Definition \ref{defi1.2.20}). 
They correspond, at the classical level, to bialgebra morphisms taking their values in the polynomial bialgebra $(\K[T],m,\Delta)$, where $m$ is the usual product of polynomials and $\Delta$ the coproduct defined by
\[\Delta(T)=T\otimes 1+1\otimes T.\]
We show that these polynomial invariants are in a bijective correspondence with characters (Theorem \ref{theo1.2.23}). The character associated to a polynomial invariant is given by the evaluation at $T=1$,
and we give an explicit formula for the converse bijection, with the help of Hilbert polynomials. Using the exponential-logarithm bijection, we can also associate bijectively a polynomial invariant to any infinitesimal character (Proposition \ref{prop1.2.25}). Here, Hilbert polynomials are replaced by monomials $\dfrac{T^n}{n!}$. 

In the next section, we introduce the twisted version of double bialgebras. A twisted double bialgebra is a twisted bialgebra $\calQ$ with complementary coproducts $\delta_\sim$, associated to any equivalence of any finite set,
satisfying good compatibilities (Definition \ref{defi1.3.1}). We show that the characters of a twisted double bialgebra inherit a second convolution product, right-distributive on the first one (Proposition \ref{prop1.3.2}).
Among all the polynomial invariants of a twisted double bialgebra $\calQ$, only one is compatible with the coproduct $\delta$ and its associated counit, $\K[T]$ being given the second coproduct $\delta$ defined by 
\[\delta(T)=T\otimes T,\]
see Definition \ref{defi1.3.4} and Theorem \ref{theo1.3.6}. It will be called the fundamental polynomial invariant of $\calQ$. 
We show in Proposition \ref{prop1.3.7} that the monoid of characters of $\calQ$, with the product associated to $\delta$, acts (on the right) on the set of polynomial invariants of $\calQ$ and that this action is freely generated by the fundamental polynomial invariant:
in other words, the fundamental polynomial invariant allows to bijectively find all the other polynomial invariants by the action of characters.
We also show that the action of formal series of a connected double twisted bialgebra can be simplified thanks to the counit $\epsilon_\delta$ of the coproduct $\delta$ (Proposition \ref{prop1.3.9}). 
Consequently, cancellation-free formulas for the antipode or for the eulerian idempotent of a connected twisted double bialgebra can be given, with the help of the fundamental polynomial invariant (Theorem \ref{theo1.3.10}).

Reminders on the bosonic Fock functor are given in the last section of this chapter. This functor, defined in \cite{Aguiar2010}, associates to any linear species a graded vector space. As it is a functor of symmetric categories,
it sends any twisted algebra, coalgebra or bialgebra $\calQ$ to a classical algebra, coalgebra, or bialgebra $\calF[\calQ]$. 
We show that when $\calQ$ is a double twisted double bialgebra, then $\calF[\calQ]$ inherits a second coproduct $\delta$ making it a double bialgebra,
that is to say that $(\calF[\calQ],m,\Delta)$ is a bialgebra in the category of right comodules over the bialgebra $(\calF[\calQ],m,\delta)$, the coaction being $\delta$ itself.

The last section of this chapter is about the twisted quasi-shuffle algebra of compositions $\comp$.
We give the twisted version of Aguiar, Bergeron and Sottile's theorem on the terminal property of $\QSym$ \cite{Aguiar2006-2}:
we first introduce the quasi-shuffle product $\squplus$ on $\comp$ (Theorem \ref{theo1.5.5}), thanks to a cofreeness property (Proposition \ref{prop1.5.4}). 
Then the twisted bialgebra $(\ell\comp,\squplus,\Delta)$, together with a character $\epsilon$, satisfies a terminal property (Theorem \ref{theo1.5.6}).
it is noticeable that the application of the Fock functor sends $\ell\comp$ to $\QSym$.

\section{Twisted objects}

\subsection{Species}

\begin{defi}\label{defi1.1.1} \cite{Joyal1981,Joyal1986} \begin{enumerate}
\item A species is a functor $\calQ$ from the category of finite sets, with bijections, to the category of sets, with maps. In other words:
\begin{itemize}
\item For any finite set $X$, $\calQ[X]$ is a set.
\item For any finite sets $X$ and $Y$, and any bijection $\sigma:X\longrightarrow Y$, $\calQ[\sigma]$ is a map from $\calQ[X]$ to $\calQ[Y]$.
\item If $X$ is a finite set, $\calQ[\id_X]=\id_{\calQ[X]}$.
\item If $X$, $Y,Z$ are three finite sets, and $\sigma:X\longrightarrow Y$ and $\tau:Y\longrightarrow Z$ are bijections, then $\calQ[\tau\circ \sigma]=\calQ[\tau]\circ \calQ[\sigma]$.
\end{itemize}
\item If $\calQ$ and $\calR$ are two species, a species morphism from $\calQ$ to $\calR$ is a natural transformation $\phi$ between these two functors:
\begin{itemize}
\item For any finite set $X$, $\phi_X:\calQ[X]\longrightarrow \calR[Y]$ is a map.
\item If $\sigma:X\longrightarrow Y$ is a bijection between two finite sets, the following diagram commutes:
\[\xymatrix{\calQ[X]\ar[d]_{\phi_X} \ar[r]^{\calQ[\sigma]}&\calQ[Y]\ar[d]^{\phi_Y}\\
\calR[X]\ar[r]_{\calR[\sigma]}&\calR[Y]}\]
\end{itemize}
We shall simply write $\phi$ instead of $\phi_X$, in order to enlighten the notations.
\end{enumerate}\end{defi}

\begin{example}\label{ex1.1.1} \begin{enumerate}
\item For any finite set $X$, we denote by $\calE[X]$ the set of equivalences on $X$. Then $\calE$ is a species.
The sequence $(|\calE[n]|)_{n\geq 0}$ is the sequence of Bell numbers, see \cite[Entry A000110]{Sloane}.
\item For any finite set $X$, we denote by $\calS[X]$ the set of subsets of $X$. For any $k\geq 0$,
we denote by $\calS_k[X]$ the set of subsets of $X$ of cardinality $k$, and by $\calS_{\geq k}[X]$ the set of subsets of $X$ of cardinality $\geq k$. Then $\calS$ is a species and, for any $k\geq 0$, $\calS_k$ and $\calS_{\geq k}$ are subspecies of $\calS$. For any finite set $X$,
\begin{align*}
|\calS[X]|&=2^{|X|}, &|\calS_k[X]|&=\binom{|X|}{k}.
\end{align*}
\item A composition, or ordered partition of a finite set $X$, is a family $\pi=(X_1,\ldots,X_k)$ of nonempty subsets of $X$, such that $X=X_1\sqcup\cdots\sqcup X_k$. The length $\ell(\pi)$ is the integer $k$.
We denote by $\comp$ the species of compositions. For example,
\begin{align*}
\comp[\emptyset]&=\{\emptyset\},\\
\comp[1]&=\{(\{1\})\},\\
\comp[2]&=\{(\{1,2\}),\:(\{1\},\{2\}),\: (\{2\},\{1\})\},\\
\comp[3]&=\left\{\begin{array}{c}
(\{1,2,3\}),\:(\{1,2\},\{3\}),\:(\{1,3\},\{2\}),\:(\{2,3\},\{1\}),\\
(\{1\},\{2,3\}),\:(\{2\},\{1,3\}),\:(\{3\},\{1,3\}),\\
(\{1\},\{2\},\{3\}),\:(\{1\},\{3\},\{2\}),\:(\{2\},\{1\},\{3\}),\\
(\{2\},\{3\},\{1\}),\:(\{3\},\{1\},\{2\}),\:(\{3\},\{2\},\{1\})
\end{array}\right\}.
\end{align*}
In the sequel, we will shorten the notations by deleting some brackets and comas, as follows:
\begin{align*}
\comp[\emptyset]&=\{\emptyset\},\\
\comp[1]&=\{(1)\},\\
\comp[2]&=\{(12),\:(1,2),\: (2,1)\},\\
\comp[3]&=\left\{\begin{array}{c}
(123),\:(12,3),\:(13,2),\:(23,1),\\
(1,23),\:(2,13),\:(3,13),\\
(1,2,3),\:(1,3,2),\:(2,1,3),\\
(2,3,1),\:(3,1,2),\:(3,2,1)
\end{array}\right\}.
\end{align*}
The cardinality of $\comp[n]$ is given by the $n^{\mathrm{th}}$ Fubini number, see Entry\href{https://oeis.org/A000670}{A000670} of the OEIS \cite{Sloane}. 
For $k\in \N$, we denote by $\comp_k$ the subspecies of $\comp$ formed by compositions of length $k$. In particular, if $|X|<k$, then $\comp_k[X]=\emptyset$.
\item Other examples will be given in the next chapter (graphs, hypergraphs, posets, etc), see Examples \ref{ex2.1.1}-\ref{ex2.1.8}.
\end{enumerate}\end{example}

Similarly, a $\K$-linear species (shortly, linear species) is a functor from the category of finite sets, with bijections, to the category of $\K$-vector spaces, with linear maps.
Morphism of linear species are natural transformations $\phi$, such that for any finite set $X$, $\phi_X$ is linear.
There exists a functor $\ell$ from the category of species to the category of $\K$-linear species, sending a species $\calQ$ to the linear species $\ell\calQ$, defined as follows:
\begin{itemize}
\item For any finite set $X$, $\ell\calQ[X]$ is the vector space generated by the set $\calQ[X]$.
In other words, $\calQ[X]$ is a basis of $\ell\calQ[X]$.
\item For any bijection $\sigma:X\longrightarrow Y$ between two finite sets, $\ell\calQ[\sigma]:\ell\calQ[X]\longrightarrow\ell\calQ[Y]$ is the unique linear map such that $\ell\calQ[\sigma]_{\mid \calQ[X]}=\calQ[\sigma]$. 
\end{itemize}

\subsection{Twisted algebras, coalgebras and bialgebras}

The category of linear species is a symmetric monoidal category, for the tensor product defined as follows: if $\calQ$ and $\calR$ are two linear species, and $X$ is a finite set,
\[\calQ\otimes \calR[X]=\bigoplus_{X'\sqcup X''=X}\calQ[X']\otimes \calR[X''].\]
This will allows to define algebras, coalgebras, bialgebras, Hopf algebras, etc, in the category of linear species.

If $\calB$ and $\calC$ are species, we define similarly their cartesian product as follows: if $X$ is a finite set, 
\[\calB\times \calC[X]=\bigsqcup_{X=X'\sqcup X''} \calB[X']\times \calC[X''].\]
Note that $\ell\calB\otimes \ell \calC$ is naturally isomorphic to $\ell(\calB\times \calC)$.

We now give details on twisted algebras, coalgebras, and bialgebras.

\begin{defi}\label{defi1.1.2}
A twisted algebra is an algebra in the category of linear species, with the Cauchy tensor product. 
With more details, a twisted algebra is a pair $(\calQ,m)$, where:
\begin{itemize}
\item $\calQ$ is a linear species.
\item For any pair $(X$, $Y)$ of disjoint finite sets, $m_{X,Y}:\calQ[X]\otimes\calQ[Y]\longrightarrow\calQ[X\sqcup Y]$ is a linear map.
\end{itemize}
Such that:
\begin{itemize}
\item If $(X$, $Y)$ and $(X',Y')$ are pairs of disjoint finite sets, $\sigma:X\longrightarrow X'$ and $\tau:Y\longrightarrow Y'$ are bijections, then the following diagram commutes:
\[\xymatrix{\calQ[X]\otimes\calQ[Y]\ar[r]^{m_{X,Y}}\ar[d]_{\calQ[\sigma]\otimes\calQ[\tau]}&\calQ[X\sqcup Y]\ar[d]^{\calQ[\sigma\sqcup\tau]}\\
\calQ[X']\otimes\calQ[Y']\ar[r]_{\hspace{2mm}m_{X',Y'}}&\calQ[X'\sqcup Y']}\]
where $\sigma\sqcup\tau$ is the bijection defined by
\[\sigma\sqcup\tau:\left\{\begin{array}{rcl}
X\sqcup Y&\longrightarrow&X'\sqcup Y'\\
x&\longmapsto&\begin{cases}
\sigma(x)\mbox{ if }x\in X,\\
\tau(x)\mbox{ if }x\in Y.
\end{cases}\end{array}\right.\]
\item If $(X$, $Y,Z)$ is a triple of pairwise disjoint finite sets, then the following diagram commutes:
\[\xymatrix{\calQ[X]\otimes\calQ[Y]\otimes\calQ[Z]\ar[rr]^{m_{X,Y}\otimes\id_{\calQ[Z]}}
\ar[d]_{\id_{\calQ[Z]}\otimes m_{Y,Z}}&&\calQ[X\sqcup Y]\otimes\calQ[Z]\ar[d]^{m_{X\sqcup Y,Z}}\\
\calQ[X]\otimes\calQ[Y\sqcup Z]\ar[rr]_{m_{X,Y\sqcup Z}}&&\calQ[X\sqcup Y\sqcup Z]}\]
\item There exists $1_\calQ\in\calQ[\emptyset]$, such that for any finite set $X$, for any $x\in\calQ[X]$,
\[m_{\emptyset,X}(1_\calQ\otimes x)=m_{X,\emptyset}(x\otimes 1_\calQ)=x.\]
\end{itemize}
It is commutative if, for any finite disjoint finite sets $X$ and $Y$,
\[m_{X,Y}\circ \tau_{\calQ[Y],\calQ[X]}=m_{Y,X}. \] 
\end{defi}

Dually:

\begin{defi}\label{defi1.1.3}
A twisted coalgebra is a coalgebra in the category of linear species, with the Cauchy tensor product. 
With more details, a twisted coalgebra is a pair $(\calQ,\Delta)$, where:
\begin{itemize}
\item $\calQ$ is a linear species.
\item For any pair $(X$, $Y)$ of disjoint finite sets, $\Delta_{X,Y}:\calQ[X\sqcup Y]\longrightarrow\calQ[X]\otimes\calQ[Y]$ is a linear map.
\end{itemize}
Such that:
\begin{itemize}
\item If $(X$, $Y)$ and $(X',Y')$ are pairs of disjoint finite sets, $\sigma:X\longrightarrow X'$ and $\tau:Y\longrightarrow Y'$ are bijections, then the following diagram commutes:
\[\xymatrix{\calQ[X]\otimes\calQ[Y]\ar[d]_{\calQ[\sigma]\otimes\calQ[\tau]}&\calQ[X\sqcup Y]\ar[l]_{\Delta_{X,Y}}\ar[d]^{\calQ[\sigma\sqcup\tau]}\\
\calQ[X']\otimes\calQ[Y']&\calQ[X'\sqcup Y']\ar[l]^{\hspace{2mm}\Delta_{X',Y'}}}\]
\item If $(X$, $Y,Z)$ is a triple of pairwise disjoint finite sets, then the following diagram commutes:
\[\xymatrix{\calQ[X]\otimes\calQ[Y]\otimes\calQ[Z]
&&\calQ[X\sqcup Y]\otimes\calQ[Z]\ar[ll]_{\Delta_{X,Y}\otimes\id_{\calQ[Z]}}\\
\calQ[X]\otimes\calQ[Y\sqcup Z]\ar[u]^{\id_{\calQ[Z]}\otimes\Delta_{Y,Z}}&&\calQ[X\sqcup Y\sqcup Z]\ar[ll]^{\Delta_{X,Y\sqcup Z}}\ar[u]_{\Delta_{X\sqcup Y,Z}}}\]
\item There exists a linear map $\varepsilon_\Delta:\calQ[\emptyset]\longrightarrow\K$, such that for any finite set $X$, 
\[\left(\varepsilon_\Delta\otimes\id_{\calQ[X]}\right)\circ\Delta_{\emptyset,X}=\left(\id_{\calQ[X]}\otimes\varepsilon_\Delta\right)\circ\Delta_{X,\emptyset}=\id_{\calQ[X]}.\]
\end{itemize}
It is cocommutative if, for any finite disjoint finite sets $X$ and $Y$,
\[ \tau_{\calQ[Y],\calQ[X]}\circ \Delta_{X,Y}=\Delta_{Y,X}. \] 
\end{defi}

\begin{notation}
If $(\calQ,\Delta)$ is a coalgebra, we extend $\varepsilon_\Delta$ to the whole $\calQ$ by putting, for any nonempty finite set $X$ and any $x\in \calQ[X]$, $\varepsilon_\Delta(x)=0$.
\end{notation}

And finally:

\begin{defi}\label{defi1.1.4}
A twisted bialgebra is a bialgebra in the category of linear species, with the Cauchy tensor product. With more details, a twisted coalgebra is a triple $(\calQ,m,\Delta)$, such that:
\begin{itemize}
\item $(\calQ,m)$ is a twisted algebra.
\item $(\calQ,\Delta)$ is a twisted coalgebra.
\item $\varepsilon_\Delta(1_\calQ)=1$, and for any $x,y\in\calQ[\emptyset]$, $\varepsilon_\Delta\circ m_{\emptyset,\emptyset}(x\otimes y)=\varepsilon_\Delta(x)\varepsilon_\Delta(y)$. 
\item If $(X$, $Y)$ and $(X',Y')$ are two pairs of disjoint finite sets such that $X\sqcup Y=X'\sqcup Y'$, then
the following diagram commutes:
\[\xymatrix{\calQ[X]\otimes\calQ[Y]\ar[rr]^{m_{X,Y}}\ar[d]_{\Delta_{X\cap X',X\cap Y'}\otimes\Delta_{Y\cap X,Y\cap Y'}}&&\calQ[X\sqcup Y]\ar[d]^{\Delta_{X',Y'}}\\
\calQ[X\cap X']\otimes\calQ[X\cap Y']\otimes\calQ[Y\cap X']\otimes\calQ[Y\cap Y']\ar[d]_{\id_{\calQ[X\cap X']}\otimes\tau_{\calQ[X\cap Y'],\calQ[Y\cap X']}\otimes\id_{\calQ[Y\cap Y']}}&&\calQ[X']\otimes\calQ[Y']\\
\calQ[X\cap X']\otimes\calQ[Y\cap X']\otimes\calQ[X\cap Y']\otimes\calQ[Y\cap Y']\ar[rru]_{\hspace{1.6cm}m_{X\cap X',Y\cap X'}\otimes m_{X\cap Y',Y\cap Y'}}&&}\]
\end{itemize}\end{defi}

\begin{remark}
If $\calQ$ and $\calR$ are two twisted algebras, respectively coalgebras, then $\calQ\otimes\calR$ is also a twisted algebra, respectively coalgebra. 
In a twisted bialgebra $\calQ$, the product is a morphism of twisted coalgebras from $\calQ\otimes \calQ$ to $\calQ$, whereas the coproduct is a morphism of twisted algebras from $\calQ$ to $\calQ\otimes \calQ$.
\end{remark}

\begin{example}
For any finite set $X$, we put $\lcom[X]=\K$. For any bijection $\sigma:X\longrightarrow Y$ between two finite sets, 
$\lcom[\sigma]$ is the identity of $\K$. Then $\lcom$ is a species. It is a twisted bialgebra:
\begin{enumerate}
\item For any pair $(X,Y)$ of finite sets, any $f\in \lcom[X]$, $f\in \lcom[Y]$, $m_{X,Y}(f\otimes f)=ff$.
\item For any pair $(X,Y)$ of finite sets, any $f\in \lcom[X\sqcup Y]$, $\Delta_{X,Y}(f)=f (1\otimes 1)$.
\end{enumerate}
It is commutative and cocommutative. 
\end{example}

\begin{example}
We can make $\ell\calS$ a twisted bialgebra as follows:
\begin{enumerate}
\item For any pair $(X,Y)$ of disjoint finite sets, any $A\in \calS[X]$, $B\in \calS[Y]$, $m_{X,Y}(A\otimes B)=A\sqcup B$.
\item For any pair $(X,Y)$ of disjoint finite sets, any $A\in \calS[X\sqcup Y]$, $\Delta_{X,Y}(A)=(A\cap X)\otimes (A\cap Y)$. 
\end{enumerate}
It is commutative and cocommutative. \end{example}

\begin{example}
We can make $\ell\comp$ a twisted bialgebra as follows:
\begin{enumerate}
\item For any pair $(X,Y)$ of disjoint finite sets, any $(X_1,\ldots,X_k)\in \comp[X]$, $(Y_1,\ldots,Y_l)\in \comp[Y]$,
\[m_{X,Y}((X_1,\ldots,X_k)\otimes (Y_1,\ldots,Y_l))=(X_1,\ldots,X_k,Y_1,\ldots,Y_l).\]
\item For any pair $(X,Y)$ of disjoint finite sets, any $(Z_1,\ldots,X_n)\in \comp[X\sqcup Y]$,
\[\Delta_{X,Y}(Z_1,\ldots,Z_n)=\elim(Z_1\cap X,\ldots,Z_n\cap X)\otimes\elim(Z_1\cap Y,\ldots,Z_n\cap Y),\]
where $\elim$ means that one deletes the empty sets.
\end{enumerate}
It is cocommutative, but not commutative. 
\end{example}

\section{Antipode, characters and polynomial invariants}

\subsection{Convolution}

\begin{notation}
Let $\calQ$ and $\calR$ be linear species. We denote by $\Hom(\calQ,\calR)$ the vector space of morphisms of species from $\calQ$ to $\calR$. 
\end{notation}

In the particular case where $\calR=\lcom$, we rather talk of the dual of $\calQ$:

\begin{defi}\label{defi1.2.1}
Let $\calQ$ be a linear species. We denote by $\calQ^*$ the space $\Hom(\calQ,\lcom)$. 
In other words, an element $f\in \calQ^*$ is given by a map $f_X:\calQ[X]\longrightarrow \K$ for any finite set $X$, such that, if $\sigma:X\longrightarrow Y$ is a bijection between two finite sets, the following diagram commutes:
\[\xymatrix{\calQ[X]\ar[rr]^{\calQ[\sigma]}\ar[rd]_{f_X}&&\calQ[Y]\ar[ld]^{f_Y}\\ &\K&}\]
With an abuse of notations, we shall often simply denote $f(x)$ instead of $f_X(x)$ for $X$ a finite set and $x\in \calQ[X]$.
\end{defi}

\begin{prop}\label{prop1.2.2}
Let $(\calQ,\Delta)$ be a twisted coalgebra and $(\calR,m)$ be a twisted algebra. The convolution product $\star$ of $\Hom(\calQ,\calR)$ is defined as follows: for any $f,g\in \Hom(\calQ,\calR)$, any finite set $X$,
\[(f\star g)_X=\sum_{X=X'\sqcup X''}m_{X,X'}\circ (f_{X'}\otimes g_{X''})\circ \Delta_{X',X''}.\]
Then $(\Hom(\calQ,\calR),\star)$ is an associative and unitary algebra. 
Its unit is denoted by $\iota_{\calQ,\calR}$. For any finite set $X$, and any $x\in \calQ[X]$,
\[\iota_{\calQ,\calR}(x)=\begin{cases}
\varepsilon_\Delta(x)1_\calA\mbox{ if }X=\emptyset,\\
0\mbox{ otherwise}.
\end{cases}\]
When $\calQ=\calR$, we will write $\iota_\calQ$ instead of $\iota_{\calQ,\calQ}$.
\end{prop}

\begin{proof}
Let $f,g,h\in \Hom(\calQ,\calR)$. 
\begin{align*}
&(f\star g)\star h\\
&=\sum_{X=X'\sqcup X''\sqcup X'''}m_{X'\sqcup X'',X'''}\circ \left(m_{X',X''}\otimes \id_{\calQ[X''']}\right)(f_{X'}\otimes g_{X''}\otimes h_{X'''})\circ \left(\Delta_{X',X''}\otimes \id_{\calQ[X''']}\right)\circ \Delta_{X'\sqcup X'',X'''}\\
&=\sum_{X=X'\sqcup X''\sqcup X'''}m_{X',X''\sqcup X'''}\circ \left(\id_{\calQ[X']}\otimes m_{X'',X'''}\right)\circ
(f_{X'}\otimes g_{X''}\otimes h_{X'''})\circ \left( \id_{\calQ[X']}\otimes\Delta_{X'',X'''}\right)\circ \Delta_{X',X''\sqcup X'''}\\
&=f\star(g\star h).
\end{align*}
So $\star$ is associative. Let $f\in \Hom(\calQ,\calR)$.
\[(f\star\iota_{\calQ,\calR})_X=m_{X,\emptyset}\circ(f\otimes \iota_{\calQ,\calR})\circ \Delta_{X,\emptyset}+0=f.\]
Similarly, $\iota_{\calQ,\calR}\star f=f$. So $\iota_{\calQ,\calR}$ is the unit of $\Hom(\calQ,\calR)$.
\end{proof}

\begin{remark} 
In the particular case where $\calR=\lcom$, we have $\iota_{\calQ,\lcom}=\varepsilon_\Delta$.
\end{remark}

We now give $\Hom(\calQ,\calR)$ a metric.

\begin{prop}\label{prop1.2.3}
Let $\calQ$ and $\calR$ be linear species. For $f\in \Hom(\calQ,\calR)$, nonzero, we denote by $\val(f)$ the minimal integer $n$ such that for any finite set $X$ of cardinality $<n$, $f_X=0$.
By convention, $\val(0)=+\infty$. We put $||f||=2^{-\val(f)}$, with the convention that $2^{-\infty}=0$. Then:
\begin{enumerate}
\item For any $f,g\in \Hom(\calQ,\calR)$, $\val(f+g)\geq\min (\val(f),\val(g))$, and $||f+g||\leq \max(||f||,||g||)\leq ||f||+||g||$.
\item For any $f\in \Hom(\calQ,\calR)$, and any $\lambda\in \K$, nonzero, $\val(\lambda f)=\val(f)$ and $||\lambda f||=||f||$.
\end{enumerate}
Consequently, putting $d(f,g)=||f-g||=2^{-\val(f-g)}$ for any $f,g\in \Hom(\calQ,\calR)$, $(\Hom(\calQ,\calR),d)$ is a metric space.
\end{prop}

\begin{proof}
1. Let us assume that $X$ is a finite set, with $|X|<\min(\val(f),\val(g))$. Then $f_X=g_X=0$, so $(f-g)_X=0$. Therefore, $\val(f+g)\geq\min (\val(f),\val(g))$.

2. Immediate.
\end{proof}

\begin{prop}\label{prop1.2.4}
Let $\calQ$ and $\calR$ be linear species. The metric space $(\Hom(\calQ,\calR),d)$ is complete.
\end{prop}

\begin{proof}
Let $\left(f^{(n)}\right)_{n\geq 0}$ be a Cauchy sequence of $\Hom(\calQ,\calR)$. For any $N\in \N$, there exists $n_N\in \N$, such that if $p,q\geq n_N$, $d\left(f^{(p)},f^{(q)}\right)<2^{-N}$, that is to say $\val\left(f^{(p)}-f^{(q)}\right)>N$. 
Consequently, if $|X|<N$ and $p,q\geq n_N$, then $f^{(p)}_X=f^{(q)}_X$: the sequence $(f^{(n)}_X)_{n\geq 0}$ is stationary. Let us denote by $f_X$ the common value of $f^{(n)}_X$ for $n$ great enough. This defines an element $f\in \Hom(\calQ,\calR)$.
Moreover, if $p,q\geq n_N$, $f^{(q)}_X=f^{(p)}_X=f_X$ for any finite set $X$ of cardinality $\leq N$, so $d\left(f^{(p)},f\right)<2^{-N}$. Consequently, $\left(f^{(n)}\right)_{n\geq 0}$ converges to $f$.
\end{proof}

\begin{prop}
Let $\calQ$ and $\calR$ be linear species. For any $f,g\in \Hom(\calQ,\calR)$, $||f\star g||\leq ||f||\cdot ||g||$. Consequently, the product $\star$ is continuous.
\end{prop}

\begin{proof}
Let $X$ be a finite set of cardinality $<\val(f)+\val(g)$. Then
\[(f\star g)_X=\sum_{X=X'\sqcup X''} m_{X',X''}\circ(f_{X'}\otimes g_{X''})\circ \Delta_{X',X''}.\]
If $X=X'\sqcup X''$, as $|X|<\val(f)+\val(g)$, either $|X'|<\val(f)$ or $|X''|<\val(g)$, so either $f_{X'}=0$ or $g_{X''}=0$. Therefore, $(f\star g)_X=0$, and $\val(f\star g)\geq \val(f)+\val(g)$, which gives $||f\star g||\leq ||f||\cdot ||g||$.
\end{proof}

\begin{example}\label{ex1.2.1}
Let us consider $\calQ=\calR=\lcom$. An element $f$ of $\lcom^*$ is identified with a sequence $\theta(f)=(f_n)_{n\ge 0}\in \K^\N$,
with $f(\lambda)=n!f_n\lambda$ for any $\lambda \in \lcom[X]=\K$, when $|X|=n$. 
Moreover, if $f,g\in \lcom^*$, and if $\lambda \in \lcom[X]$ with $|X|=n$,
\begin{align*}
(f\star g)_n&=\frac{1}{n!}(f\star g)(1)
=\frac{1}{n!}\sum_{X=X'\sqcup X''} f_{X'}(1)g_{X''}(1)
=\frac{1}{n!}\sum_{k=0}^n \binom{n}{k} k!f_k (n-k)! g_{n-k}
=\sum_{k=0}^n f_kg_{n-k}.
\end{align*}
Therefore, as an algebra, $(\lcom,\star)$ is isomorphic to the algebra $\K[[T]]$ of formal series, through the map $\theta$. This map also conserves the valuation, so is an isometry from $(\lcom,d)$ to $(\K[[T]],d)$. 
\end{example}

\begin{theo}\label{theo1.2.6}
Let $\calQ$ be a twisted coalgebra, $\calR$ be a twisted algebra, and $f\in \Hom(\calQ,\calR)$, such that $\val(f)\geq 1$. For $P=\displaystyle \sum_{k=0}^\infty a_k T^k\in \K[[T]]$, we define $P(f)\in \Hom(\calQ,\calR)$ as follows: for any finite set $X$,
\begin{align}
\label{EQ1.1} P(f)_X&=\sum_{k=0}^{|X|} \sum_{(X_1,\ldots,X_k)\in \comp_k[X]}a_k m_{X_1,\ldots,X_k}\circ f^{\otimes k}\circ \Delta_{X_1,\ldots,X_k}.
\end{align}
Then:
\begin{enumerate}
\item The map $P\longrightarrow P(f)$ is a continuous algebra map from $\K[[T]]$ to $(\Hom(\calQ,\calR),\star)$.
\item Let $P$, $Q\in \K[[T]]$, such that $Q(0)=0$. Then $P(Q(f))=P\circ Q(f)$. 
\end{enumerate}\end{theo}

\begin{proof}
1. Let us consider the unique algebra map $\theta$ from $\K[T]$ to $\calQ$ which sends $T$ to $f$. As $\val(f)\geq 1$, $||f||\leq 2^{-1}$, and for any $n\geq 0$,
\[||\theta(T^n)||=||f^{*n}||\leq ||f||^n \leq2^{-n}.\]
As a consequence, if $P\in \K[T]$, $||P(f)|| \leq 2^{-\val(P)}=||P||$, so $\theta$ is continuous. As $\K[T]$ is dense in $\K[[T]]$ and $(\Hom(\calQ,\calR),d)$ is complete, $\theta$ can be uniquely extended as a continuous linear map from $\K[[T]]$ to $\Hom(\calQ,\calR)$.
As both products of $\K[[T]]$ and $\Hom(\calQ,\calR)$ are continuous, this extension, still denoted by $\theta$, is an algebra morphism.

Let $P=\displaystyle \sum_{k=0}^\infty a_k T^k\in \K[[T]]$, and $X$ be a finite set. Then, as $f_{\emptyset}=0$,
\begin{align*}
\theta(P)_X&=\sum_{k=0}^\infty \sum_{X=X_1\sqcup \ldots X_k} a_k m_{X_1,\ldots,X_k}\circ f^{\otimes k}\circ \Delta_{X_1,\ldots,X_k}\\
&=\sum_{k=0}^{|X|} \sum_{(X_1,\ldots,X_k)\in \comp_k[X]}a^k m_{X_1,\ldots,X_k}\circ f^{\otimes k}\circ \Delta_{X_1,\ldots,X_k}+0=P(f).
\end{align*}

2. As $Q(0)=0$, $Q(f)_\emptyset=0$, so $\val(Q(f))\geq 1$ and $P(Q(f))$ exists. By composition, $P\longrightarrow P\circ Q(f)$ is a continuous algebra map from $\K[[T]]$ to $\Hom(\calQ,\calR)$, which sends $T$ to $Q(f)$. 
By density of $\K[T]$ in $\K[[T]]$, it is $P\longrightarrow P(Q(f))$.
\end{proof}

\begin{notation}
Let $\calQ$ be a twisted coalgebra, and $\calR$ be a twisted algebra. We put
\begin{align*}
\Hom(\calQ,\calR)_0&=\{f\in \Hom(\calQ,\calR)\mid f_\emptyset=0\},&\Hom(\calQ,\calR)_0&=\{f\in \Hom(\calQ,\calR)\mid f_\emptyset=\iota_{\calQ,\calR}\}.
\end{align*}\end{notation}

Here is a first application:

\begin{prop}\label{prop1.2.7}
Let $\calQ$ be a twisted coalgebra, and $\calR$ be a twisted algebra. Any element $f$ of $\Hom(\calQ,\calR)_1$ is invertible for the convolution product. Its inverse is given as follows: for any finite set $X$,
\begin{align}
\label{EQ1.2} f^{\star-1}_X&=\sum_{k=0}^{|X|} (-1)^k\sum_{(X_1,\ldots,X_k)\in \comp_k[X]} m_{X_1\ldots,X_k}\circ f^{\otimes k} \circ \Delta_{X_1,\ldots,X_k}.
\end{align}\end{prop}

\begin{proof}
As $f\in \Hom(\calQ,\calR)_1$, $f-\iota_{\calQ,\calR}\in \Hom(\calQ,\calR)_0$. Let us consider $g=(1+(f-\iota_{\calQ,\calR}))^{-1}$. This is well-defined, and
\[g\star f=\left(\frac{1}{1+T}(1+T)\right)(f-\iota_{\calQ,\calR})=1(f-\iota_{\calQ,\calR})=\iota_{\calQ,\calR}.\]
Similarly, $f\star g=\iota_{\calQ,\calR}$, so $f$ is invertible, of inverse $g$. Then (\ref{EQ1.2}) is implied by (\ref{EQ1.1}), with $a_k=(-1)^k$ for any $k\geq 0$.
\end{proof}

Here is second application:

\begin{prop}\label{prop1.2.8}
Let $\calQ$ be a twisted coalgebra, and $\calR$ be a twisted algebra. The two following maps are bijections, inverse one from the other:
\begin{align*}
\exp_\calQ&:\left\{\begin{array}{rcl}
\Hom(\calQ,\calR)_0&\longrightarrow&\Hom(\calQ,\calR)_1\\
f&\longmapsto&\exp(f),
\end{array}\right.&
\ln_\calQ&:\left\{\begin{array}{rcl}
\Hom(\calQ,\calR)_1&\longrightarrow&\Hom(\calQ,\calR)_0\\
f&\longmapsto&\ln(1+(f-\iota_{\calQ,\calR})).
\end{array}\right.
\end{align*}\end{prop}

\begin{proof}
If $f\in \Hom(\calQ,\calR)_0$, then $\val(f)\geq 1$ and $\exp(f)$ is well-defined. Moreover, $\exp(f)_\emptyset= f^{\star 0}=\iota_{\calQ,\calR}$, so $\exp(f)\in \Hom(\calQ,\calR)_1$.
If $f\in \Hom(\calQ,\calR)_1$, then $\val(f-\iota_{\calQ,\calR})\geq 1$, and $\ln(1+(f-\iota_{\calQ,\calR}))$ is well-defined. Moreover, $\ln(1+(f-\iota_{\calQ,\calR}))_\emptyset=0$, so $\ln(1+(f-\iota_{\calQ,\calR}))\in \Hom(\calQ,\calR)_0$.

Let $f\in \Hom(\calQ,\calR)_0$. 
\begin{align*}
\ln_\calQ\circ \exp_\calQ(f)=\ln(1+(\exp(f)-\iota_{\calQ,\calR}))=\ln(1+T)\circ (\exp(T)-1)(f)=T(f)=f.
\end{align*}
Let $f\in \Hom(\calQ,\calR)_1$.
\begin{align*}
\exp_\calQ\circ \ln_\calQ(f)&=\exp(\ln(1+(f-\iota_{\calQ,\calR})))\\
&=\exp(T)\circ \ln(1+T)(f-\iota_{\calQ,\calR})\\
&=(1+T)(f-\iota_{\calQ,\calR})\\
&=\iota_{\calQ,\calR}+f-\iota_{\calQ,\calR}=f.
\end{align*}
So $\exp_\calQ$ and $\ln_\calQ$ are bijective, inverse one from the other. 
\end{proof}

\begin{defi}\label{defi1.2.9}
Let $\calQ$ be a twisted bialgebra, and let $f\in \calQ^*$.
\begin{enumerate}
\item We shall say that $f$ is a character of $\calQ$ if $f(1_\calQ)=1$ and for any pair $(X,Y)$ of disjoint finite sets, any $a\in \calQ[X]$, $b\in \calQ[Y]$, $f(ab)=f(a)f(b)$.
The set of characters of $\calQ$ is denoted by $\Char(\calQ)$. 
\item We shall say that $f$ is an infinitesimal character of $\calQ$ if for any pair $(X,Y)$ of disjoint finite sets, any $a\in \calQ[X]$, $b\in \calQ[Y]$, $f(ab)=\varepsilon_\Delta(a)f(b)+f(a)\varepsilon_\Delta(b)$.
The subspace of infinitesimal character of $\calQ$ is denoted by $\InfChar(\calQ)$.
\end{enumerate}\end{defi}

\begin{remark}
Let $f\in \InfChar(\calQ)$.
\begin{enumerate}
\item $f(1_\calQ)=f(1_\calQ 1_\calQ)=\varepsilon_\Delta(1_\calQ)f(1_\calQ)+f(1_\calQ)\varepsilon_\Delta(1_\calQ)=2f(1_\calQ)$, so $f(1_\calQ)=0$.
\item If $(X,Y)$ is a pair of disjoint nonempty finite sets, $a\in \calQ[X]$, and $b\in \calQ[Y]$, then $f(ab)=0$. 
\end{enumerate}\end{remark}

\begin{example}
With the identification of Example \ref{ex1.2.1}, a character of $\lcom$ corresponds to a formal series $\exp(\lambda T)$, with $\lambda\in \K$. An infinitesimal character of $\lcom$ corresponds to a formal series $\lambda T$, with $\lambda \in \K$. 
\end{example}

\begin{theo}
Let $\calQ$ be a twisted bialgebra. Then $\exp_\calQ$ induces a bijection from $\InfChar(\calQ)\cap\calQ^*_0$ to $\Char(\calQ)\cap\calQ^*_1$. Its inverse is induced by $\ln_\calQ$.
\end{theo}

\begin{proof}
We put
\[m^*:\left\{\begin{array}{rcl}
\calQ^*&\longrightarrow&(\calQ\otimes \calQ)^*\\
f&\longmapsto&f\circ m.
\end{array}\right.\]
Let $f\in \calQ^*$. If $X$ is a finite set, of cardinality $<\val(f)$, then
\[m^*(f)_X=(f\circ m)_X=\bigoplus_{X=X'\sqcup X''}f_X\circ m_{X',X''}=0,\]
so $\val(m^*(f))\geq \val(f)$, which gives $||m^*(f)||\leq ||f||$: the map $f^*$ is continuous. Moreover, as $\calQ$ is a twisted bialgebra, $m:\calQ\otimes \calQ \longrightarrow \calQ$ is a coalgebra map, so $m^*$ is an algebra map.

Let $f\in \Char(\calQ)\cap \calQ^*_1$. As $f$ is a character, $m^*(f)=f\otimes f$. Consequently, as $m*$ is a continuous algebra map,
\begin{align*}
m^*(\ln(1+(f-\varepsilon_\Delta)))&=\ln(1+(f\otimes f-\varepsilon_\Delta\otimes \varepsilon_\Delta))\\
&=\ln(1+((f\otimes \varepsilon_\Delta)*(\varepsilon_\Delta\otimes f)-\varepsilon_\Delta\otimes \varepsilon_\Delta)\\
&=\ln(1+(f\otimes \varepsilon_\Delta-\varepsilon_\Delta\otimes \varepsilon_\Delta))+\ln(1+(\varepsilon_\Delta\otimes f-\varepsilon_\Delta\otimes \varepsilon_\Delta))\\
&=\ln(1+(f-\varepsilon_\Delta))\otimes \varepsilon_\Delta+\varepsilon_\Delta\otimes \ln(1+(f-\varepsilon_\Delta)),
\end{align*}
so $\ln_\calQ(f)$ is an infinitesimal character.

Let $f\in \InfChar(\calQ)\cap \Q^*_0$. As $f$ is an infinitesimal character, $m^*(f)=\varepsilon_\Delta\otimes f+f\otimes \varepsilon_\Delta$. Therefore, as $m*$ is a continuous algebra map,
\begin{align*}
m^*(\exp(f))&=\exp(\varepsilon_\Delta\otimes f+f\otimes \varepsilon_\Delta)\\
&=\exp(\varepsilon_\Delta\otimes f)*\exp(f\otimes \varepsilon_\Delta)\\
&=(\varepsilon_\Delta\otimes \exp(f))*(\exp(f)\otimes \varepsilon_\Delta)\\
&=\exp(f)\otimes \exp(f),
\end{align*}
so $\exp(f)$ is a character. 
\end{proof}

\begin{defi}
Let $\calQ$ be a twisted bialgebra. We shall say that is is connected if $\calQ[\emptyset]=\vect(1_\calQ)$.
\end{defi}

If $\calQ$ is connected, then $\Char(\calQ)\subseteq \calQ^*_1$ and $\InfChar(\calQ)\subseteq \calQ^*_0$. This gives:

\begin{cor}\label{cor1.2.12}
Let $\calQ$ be a connected twisted bialgebra. Then $\exp_\calQ$ induces a bijection from $\InfChar(\calQ)$ to $\Char(\calQ)$. Its inverse is induces by $\ln_\calQ$.
\end{cor}

\subsection{Twisted Hopf algebras and antipodes}

\begin{defi}\label{defi1.2.13}
Let $\calQ$ be a twisted bialgebra. We shall say that it is a twisted Hopf algebra if $\id_\calQ$ has an inverse for the convolution product of $\Hom(\calQ,\calQ)$. If so, its inverse is called the antipode of $\calQ$ and denoted by $S_\calQ$. 
In other terms, for any finite set $X$, and any $x\in \calQ[X]$,
\begin{align*}
&\sum_{X=X'\sqcup X''} m_{X',X''}\circ (S_\calQ\otimes \id_{\calQ[X'']})\circ \Delta_{X',X'}(x)\\
&=\sum_{X=X'\sqcup X''} m_{X',X''}\circ (\id_{\calQ[X']}\otimes S_\calQ)\circ \Delta_{X',X'}(x)=\varepsilon_\Delta(x)1_\calQ.
\end{align*}
In particular, if $X=\emptyset$, this is equal to $0$. 
\end{defi}

Let us give the twisted version of the Takeuchi's formula for the antipode of a classical Hopf algebra \cite{Takeuchi1971}:

\begin{prop}[\textbf{Takeuchi's formula}]\label{prop1.2.14}
Let $\calQ$ be a connected twisted bialgebra. Then it is a twisted Hopf algebra, and its antipode is given as follows: for any nonempty finite set $X$, and any $x\in \calQ[X]$,
\[S_\calQ(x)=\sum_{k=1}^{|X|} \sum_{(X_1,\ldots,X_k)\in \comp[X]} (-1)^k m_{X_1,\ldots,X_k}\circ \Delta_{X_1,\ldots,X_k}(x).\]
\end{prop}

\begin{proof}
As $\calQ$ is connected, $\id_\calQ\in \Hom(\calQ,\calQ)_1$. By Proposition \ref{prop1.2.7}, it is invertible for the convolution, and its inverse is given by (\ref{EQ1.2}) with $f=\id_\calQ$. This gives Takeuchi's formula.
\end{proof}

\subsection{The eulerian idempotent}

\begin{defi}\label{defi1.2.15}
Let $\calQ$ be a connected twisted bialgebra. The eulerian idempotent of $\calQ$ is $\ln(\id_\calQ)$: as $\calQ$ is connected, $\id\in \Hom(\calQ,\calQ)_1$, so this is well-defined. For any nonempty finite set $X$, and any $x\in \calQ[X]$,
\[\ln(\id_\calQ)(x)=\sum_{k=1}^{|X|}\sum_{(X_1,\ldots,X_k)\in \comp_k[X]} \dfrac{(-1)^{k-1}}{k} m_{X_1,\ldots,X_k}\circ \Delta_{X_1,\ldots,X_k}(x).\]
\end{defi}

\begin{notation}
Let $k,l\geq 0$. A packed matrix of size $(k,l)$ is an matrix $M\in \calM_{k,l}(\Z)$, whose all entries are equal to $0$ or $1$, with no zero row nor zero column. 
The set of packed matrices of size $(k,l)$ is denoted by $\pack_{k,l}$. If $M\in \pack_{k,l}$ is 
\[||M||=\sum_{i=1}^k\sum_{j=1}^l M_{i,j}.\]
For $k,l,n\geq 1$, we denote by $p_{k,l}^n$ the number of elements of $\pack_{k,l}$ of weight $n$.
\end{notation}

\begin{example}
\begin{align*}
\pack_{1,1}&=\left\{\begin{pmatrix}
1
\end{pmatrix}\right\},\\
\pack_{1,2}&=\left\{\begin{pmatrix}
1&1
\end{pmatrix}\right\},\\
\pack_{2,1}&=\left\{\begin{pmatrix}
1\\1
\end{pmatrix}\right\},\\
\pack_{2,2}&=\left\{\begin{pmatrix}
1&1\\
1&1
\end{pmatrix},\begin{pmatrix}
1&1\\
1&0
\end{pmatrix},\begin{pmatrix}
1&1\\
0&1
\end{pmatrix},\begin{pmatrix}
0&1\\
1&1
\end{pmatrix},\begin{pmatrix}
1&0\\
1&1
\end{pmatrix},\begin{pmatrix}
1&0\\
0&1
\end{pmatrix},\begin{pmatrix}
0&1\\
1&0
\end{pmatrix}\right\}.
\end{align*}
The following array gives the cardinality of $\pack_{k,l}$ for small values of $k$ and $l$, see \cite[Entry A183109]{Sloane}.
\[\begin{array}{|c||c|c|c|c|c|}
\hline k/l&1&2&3&4&5\\
\hline\hline 1&1&1&1&1&1\\
\hline 2&1&7&25&79&241\\
\hline 3&1&25&265&2161&16081\\
\hline 4&1&79&2161&41503&693601\\
\hline 5&1&241&16081&693601&24997921\\
\hline\end{array}\]
For any $k,l\geq 0$, we put
\[P_{k,l}=\sum_{n=0}^\infty p_{k,l}^n T^n.\]
If $n>kl$, $p_{k,l}^n=0$: this belongs to $\K[T]$. We denote by $m_{k,l}^n$ the number of matrices (not necessarily packed) with $k$ rows and $l$ columns, of weight $n$, with coefficients in $\{0,1\}$. Choosing the position of the $n$ coefficients $1$ in such a matrix, we obtain
\[m_{k,l}^n=\binom{kl}{n}.\]
Choosing the nonzero rows and columns of such a matrix, we obtain
\[m_{k,l}^n=\sum_{i=0}^k \sum_{j=0}^l \binom{k}{i}\binom{l}{j} p_{i,j}^n.\]
Multiplying by $T^n$ and summing over all $n$, we obtain
\begin{align*}
\sum_{n=0}^{kl} m_{k,l}^n T^n&=\sum_{n=0}^{kl}\binom{kl}{n}T^n=(1+T)^{kl}\\
&=\sum_{n=0}^{kl}\sum_{i=0}^k \sum_{j=0}^\infty \binom{k}{i}\binom{l}{j} p_{i,j}^n X^n\\
&=\sum_{i=0}^k \sum_{j=0}^\infty \binom{k}{i}\binom{l}{j} P_{i,j}.
\end{align*} 
We consider the infinite matrices
\begin{align*}
M&=\left((1+T)^{kl}\right)_{k,l\geq 0},&P&=\left(P_{k,l}\right)_{k,l\geq 0},&C&=\left(\binom{k}{l}\right)_{k,l\geq 0}.
\end{align*}
We obtain that 
\[M=CPC^\top.\]
The inverse of $C$ is the matrix
\[C^{-1}=\left((-1)^{k-l}\binom{k}{l}\right)_{k,l\geq 0},\]
which gives $P=C^{-1}MC^{-1\top}$, which finally gives, for any $k,l\geq 0$,
\[P_{k,l}=\sum_{i=0}^k\sum_{j=0}^l (-1)^{k+l-i-j} \binom{k}{i}\binom{l}{j} (1+T)^{ij}.\]
\end{example}

\begin{prop}\label{prop1.2.16}
Let $\calQ$ be a commutative or cocommutative twisted bialgebra. We put $\rho_\calQ=\id_\calQ-\iota_\calQ$. For any $k,l\geq 1$,
\[\rho_\calQ^{\star k}\circ \rho_\calQ^{\star l}=\sum_{n=\max(k,l)}^{kl}p_{k,l}^n\rho^{\star n}.\]
\end{prop}

\begin{proof}
Let $X$ be a finite nonempty set. Then
\begin{align*}
&(\rho_\calQ^{\star k}\circ \rho_\calQ^{\star l})_X\\
&=\sum_{\substack{(X_1,\ldots,X_k)\in \comp_k[X],\\ (Y_1,\ldots,Y_l)\in \comp_l[X]}} m_{X_1,\ldots,X_k}\circ \Delta_{X_1,\ldots,X_k}\circ m_{Y_1,\ldots,Y_l}\circ \Delta_{Y_1,\ldots,Y_l}\\
&=\sum_{\substack{(X_1,\ldots,X_k)\in \comp_k[X],\\ (Y_1,\ldots,Y_l)\in \comp_l[X]}} m_{X_1\cap Y_1,\ldots,X_k\cap Y_1,\ldots,X_1\cap Y_l,\ldots,X_k\cap Y_l}\circ \Delta_{X_1\cap Y_1,\ldots,X_k\cap Y_1,\ldots,X_1\cap Y_l,\ldots,X_k\cap Y_l}\\
&=\sum_{\pi\in \comp[X]} \sum_{\substack{(X_1,\ldots,X_k)\in \comp_k[X],\\ (Y_1,\ldots,Y_l)\in \comp_l[X],\\ \elim(X_1\cap Y_1,\ldots,X_l\cap Y_l)=\pi}}m_\pi\circ \Delta_\pi. 
\end{align*}
We used the compatibility between $\Delta$ and $m$ and the commutativity or the cocommutativity of $\calQ$ for the second equality.
Let $\pi\in \comp[X]$, and $(X_1,\ldots,X_k)\in \comp_k[X]$, $(Y_1,\ldots,Y_l)\in \comp_l[X]$, such that $\elim(X_1\cap Y_1,\ldots,X_l\cap Y_l)=\pi$. We define $M\in \pack_{k,l}$ by
\[M_{i,j}=\begin{cases}
1\mbox{ if }X_i\cap Y_j\neq \emptyset,\\
0\mbox{ otherwise}.
\end{cases}\]
The weight of $M$ is the length of $\pi$. Note that if $M\in \pack_{k,l}$, as there exists at least a coefficient 1 on each row and column of $M$,
\[\max(k,l)\leq \mathrm{weight}(M)\leq kl.\]

Conversely, if $\pi\in \comp[X]$ and $M\in \pack_{k,l}$ has for weight the length of $\pi$, define $(X_1,\ldots,X_k)$ and $(Y_1,\ldots,Y_l)\in \comp[X]$ as follows:
\begin{enumerate}
\item Replace the coefficients 1 of $M$ by the components of $\pi=(X_1,\ldots,X_p)$, in this order, from the position $(1,1)$ to the position $(k,l)$, row after row. 
\item $X_i$ is the union of the sets on the $i$-th row for $i\in [k]$ and $Y_j$ is union of the sets on the $j$-th column for $j\in [l]$.
\end{enumerate}
As $M$ is packed, $(X_1,\ldots,X_k)\in \comp_k[X]$, $(Y_1,\ldots,Y_l)\in \comp_l[X]$ and $\elim(X_1\cap Y_1,\ldots,X_l\cap Y_l)=\pi$ by construction. Therefore, 
\begin{align*}
(\rho_\calQ^{\star k}\circ \rho_\calQ^{\star l})_X&=\sum_{n=\max(k,l)}^{kl} p^n_{kl}\sum_{\pi\in \comp_i[X]} m_\pi\circ \Delta_\pi=\sum_{n=\max(k,l)}^{kl} p^n_{kl}\rho_X^{\star n}. \qedhere
\end{align*}\end{proof}
The eulerian idempotent is usually not a projector. However, the following Proposition gives the result if $\calQ$ is commutative or cocommutative, which justify this name.

\begin{prop}\label{prop1.2.17}
Let $\calB$ be a commutative or cocommutative connected twisted bialgebra. Then $\ln(\id_\calQ)\circ \ln(\id_\calQ)=\ln(\id_\calQ)$.
\end{prop}

\begin{proof}
Let $X$ be a finite set. We use the notations of Proposition \ref{prop1.2.16}. We obtain that
\begin{align}
\nonumber \ln(\id_\calQ)\circ \ln(\id_\calQ)&=\left(\sum_{k=1}^\infty \frac{(-1)^{k-1}}k \rho_\calQ^{\star k}\right)\circ\left(\sum_{l=1}^\infty \frac{(-1)^{l-1}}l\rho_\calQ^{\star l}\right)\\
\nonumber &=\sum_{k=1}^\infty \sum_{l=1}^\infty \sum_{n=\max(k,l)}^{kl} p_{k,l}^n \frac{(-1)^{k+l}}{kl}\rho_\calQ^{\star n}\\
\label{EQ1.3}&=\sum_{n=1}^\infty \left(\underbrace{\sum_{k=1}^n\sum_{l=1}^n p_{k,l}^n\frac{(-1)^{k+l}}{kl}}_{=S_n}\right)\rho_\calQ^{\star n}.
\end{align}
Our aim is now to prove that $S_n=\dfrac{(-1)^{n-1}}{n}$ for any $n\geq 1$. For this, we shall study firstly the simple example of $\lcom$. Let $n\geq 1$. Then, if $X$ is a finite set, 
\[{\rho_{\lcom}^{\star k}}_X=\begin{cases}
0\mbox{ if }|X|<k,\\
k!\mbox{ if }|X|=k.
\end{cases}\]
Therefore, by triangularity, if $P\in \K[[T]]$, then $P(\rho_{\lcom})=0$ if, and only if, $P=0$.

Let now $X$ be a nonempty finite set.
\[\ln(\id_{\lcom})_X=\sum_{k=1}^{|X|} \frac{(-1)^{k-1}}{k} \sum_{\pi\in \comp_k[X]} m_\pi\circ \Delta_\pi
=\sum_{k=1}^{|X|} \frac{(-1)^{k-1}}{k} |\comp_k[X]| \id_{\lcom[X]},\]
as $m_\pi\circ \Delta\pi=\id_{\lcom[X]}$ for any $\pi\in \comp[X]$. 
Let us consider
\[\Phi_X(T)=\sum_{k=1}^{|X|} |comp_k(X)| \dfrac{T(T-1)\cdots (T-k+1)}{k!}. \]
Then 
\[\ln(\id_{\lcom})_X=\dfrac{\mathrm{d}\Phi_X}{\mathrm{d}T}(0)\id_{\lcom[X]}.\]
Let us fix $N\geq 1$. To any map $\sigma:X\longrightarrow [N]$, we associate:
\begin{enumerate}
\item An integer $k$, as $|\sigma(X)|$.
\item A subset $Y$ of $[N]$ of cardinality $k$, as $\sigma(X)$.
\item A composition $(X_1,\ldots,X_k)\in \comp_k[X]$, as $\elim(\sigma^{-1}(1),\ldots \sigma^{-1}(N))$.
\end{enumerate}
These three elements entirely determine $\sigma$, and to each such triplet $(k,Y,(X_1,\ldots,X_k))$ we can associate a map $\sigma:X\longrightarrow [N]$. Counting these maps into two ways, we obtain
\[N^{|X|}=\sum_{k=1}^{|X|}|\comp_k[X]| \binom{N}{k}=\Phi_X(N).\]
As this is true for any $N\geq 1$, we deduce that $\Phi_X=T^{|X|}$, which finally gives
\[\ln(\id_{\lcom})_X=\begin{cases}
\id_{\lcom[X]} \mbox{ if }|X|=1,\\
0\mbox{ otherwise}.
\end{cases}\]
As a consequence, $\ln(\id_{\lcom})\circ \ln(\id_{\lcom})=\ln(\id_{\lcom})$. Then (\ref{EQ1.3}) with $\calQ=\lcom$ gives
\[\sum_{n=1}^\infty S_n \rho_{\lcom}^{\star n}=\ln(\id_{\lcom})
=\sum_{n=1}^\infty \frac{(-1)^{n-1}}{n}\rho_{\lcom}^{\star n}.\]
The preceding observation gives that for any $n\geq 1$,
\[S_n=\frac{(-1)^{n-1}}{n}.\]
Then (\ref{EQ1.3}) for any $\calQ$ gives that $\ln(\id_\calQ)\circ \ln(\id_\calQ)=\ln(\id_\calQ)$.
\end{proof}

\begin{defi}
Let $\calQ$ be a twisted bialgebra.
\begin{enumerate}
\item Let $X$ be a nonempty finite set, and $x\in \calQ[X]$. We shall say that $x$ is primitive if for any $(X_1,X_2)\in \comp[X]$, $\Delta_{X_1,X_2}(x)=0$.
The linear subspecies of $\calQ$ of primitive element is denoted by $\prim(\calQ)$.
\item The subspecies $\calQ^2$ of $\calQ$ is defined as follows: $\calQ[\emptyset]=\vect(1_\calQ)$ and, if $X$ is a finite nonempty set,
\[\calQ^2[X]=\sum_{(X',X'')\in \comp[X]} m_{X',X''}(\calQ[X']\otimes \calQ[X'']).\]
\end{enumerate}\end{defi}

\begin{theo}\label{theo1.2.19}
Let $\calQ$ be a connected twisted bialgebra. 
\begin{enumerate}
\item $\calQ=\ln(\id_\calQ)(\calQ)+\calQ^2$.
\item The restriction of $\ln(\id_\calQ)$ to $\prim(\calQ)$ is the identity of $\prim(\calQ)$.
\item Let us assume that $\calQ$ is commutative. Then $\ln(\id_\calQ)$ is a projector of kernel equal to $\calQ^2$.
\item Let us assume that $\calQ$ is cocommutative. Then $\ln(\id_\calQ)$ is a projector on $\prim(\id_\calQ)$.
\end{enumerate}\end{theo}

\begin{proof}
1. As $\calQ$ is connected, $\calQ[\emptyset]=\ln(\id_\calQ)(\calQ)[\emptyset]+\calQ^2[\emptyset]$, as $1_\calQ\in \calQ^2[\emptyset]$.
Let $X$ be a finite nonempty set.
For any $x\in \calQ[X]$,
\[x=\ln(\id_\calQ)(x)-\sum_{k=2}^{|X|} \sum_{\pi\in \comp_k[X]} \frac{(-1)^{k-1}}{k} \underbrace{m_\pi\circ \Delta_\pi(x)}_{\in \calQ^2[X]},\]
so $\calQ[X]=\ln(\id_\calQ)(\calQ[X])+\calQ^2[X]$.

2. Let $X$ be a finite set, and $x\in \prim(\calQ)[X]$. By definition of a primitive element,
\[\ln(\id_\calQ)(x)=m_X\circ \Delta(x)+0=x.\]

3. Let us assume that $\calQ$ is commutative.
By definition of a twisted bialgebra, $m:\calQ\otimes \calQ\longrightarrow \calQ$ is a morphism of twisted coalgebras. We denote by $m'$ the product of $\calQ\otimes \calQ$.
If $X'$, $X''$, $Y'$, $Y''$ are pairwise disjoint finite sets, $x'\in \calQ[X']$, $x''\in \calQ[X'']$, $y'\in \calQ[Y']$, $y''\in \calQ[Y'']$,
\begin{align*}
&m_{X'\sqcup Y',X''\sqcup Y''}\circ m_{X'\sqcup X'',Y'\sqcup Y''}((x'\otimes x'') \otimes (y'\otimes y''))\\
&=m_{X'\sqcup Y',X''\sqcup Y''}(x'y'\otimes x''y'')\\
&=x'y'x''y''\\
&=x'x''y'y''\\
&=m'_{X'\sqcup X'',Y'\sqcup Y''}\circ (m_{X',X''}(x'\otimes x'')\otimes m_{Y',Y''}(y'\otimes y'')),
\end{align*}
so $m$ is a morphism of twisted bialgebras. Therefore, the following map is an algebra morphism:
\[m^*:\left\{\begin{array}{rcl}
\Hom(\calQ,\calQ)&\longrightarrow&\Hom(\calQ\otimes \calQ,\calQ)\\
f&\longmapsto&f\circ m.
\end{array}\right.\]
If $f\in \Hom(\calQ,\calQ)$, then $\val(m^*(f))\geq \val(f)$, so $m^*$ is continuous. Therefore, $m^*(\ln(\id_\calQ))=\ln(m^*(\id_\calQ))=\ln(m)$. Moreover, 
\[m=(\varepsilon_\Delta\otimes \id_\calQ)\star(\id_\calQ\otimes \varepsilon_\Delta)=(\id_\calQ\otimes \varepsilon_\Delta)\star (\varepsilon_\Delta\otimes \id_\calQ).\]
Therefore,
\[\ln(m)=\ln(\varepsilon_\Delta\otimes \id_\calQ)+\ln(\id_\calQ\otimes \varepsilon_\Delta)=\varepsilon_\Delta\otimes \ln(\id_\calQ)+\ln(\id_\calQ)\otimes \varepsilon_\Delta.\]
As $X$ and $Y$ are nonempty, $\varepsilon_\Delta(x)=\varepsilon_\Delta(y)=0$ and
\[\ln(\id_\calQ)(xy)=m^*(\ln(\id_\calQ))(x\otimes y)=\ln(m)(x\otimes y)=\varepsilon_\Delta(x)\ln(\id_\calQ)(y)+\ln(\id_\calQ)(x)\varepsilon_\Delta(y)=0.\]
As $\ln(\id_\calQ)(1_\calQ)=0$, we obtain that $\calQ^2\subseteq \ker(\ln(\id_\calQ))$.
As we already know that $\ln(\id_\calQ)$ is a projector (Proposition \ref{prop1.2.17}), the first item implies that $\calQ^2\subseteq \ker(\ln(\id_\calQ))$, necessarily $\calQ^2=\ker(\ln(\id_\calQ))$.

4. By definition of a twisted bialgebra, $\Delta:\calQ\longrightarrow \calQ\otimes \calQ$ is a morphism of twisted algebras. Let $X$, $X'$, $Y$, $Y'$ be pairwise disjoint finite sets. We denote by $\Delta'$ the coproduct of $\calQ\otimes \calQ$.
As $\calQ$ is cocommutative,
\begin{align*}
\Delta'_{X'\sqcup Y',X''\sqcup Y''}\circ\Delta_{X'\sqcup X'',Y'\sqcup Y''}&=(\id_{\calQ[X']} \otimes \tau_{\calQ[X''],\calQ[Y']}\otimes \id_{\calQ[Y'']})\circ \Delta_{X',X'',Y',Y''}\\
&=\Delta_{X',Y',X'',Y''}\\
&=(\Delta_{X',Y'}\otimes \Delta_{X'',Y''})\circ \Delta_{X'\sqcup Y',X''\sqcup Y''},
\end{align*}
so $\Delta$ is a morphism of twisted bialgebras. Consequently, Therefore, the following map is an algebra morphism:
\[\Delta^*:\left\{\begin{array}{rcl}
\Hom(\calQ,\calQ\otimes \calQ)&\longrightarrow&\Hom(\calQ,\calQ)\\
f&\longmapsto&f\circ \Delta.
\end{array}\right.\]
If $f\in \Hom(\calQ,\calQ\otimes \calQ)$, then $\val(\Delta^*(f))\geq \val(f)$, so $\Delta^*$ is continuous. Therefore, $\Delta^*(\ln(\id_\calQ))=\ln(\Delta^*(\id_\calQ))=\ln(\Delta)$. Moreover, 
\[\Delta=(\eta\otimes \id_\calQ)\star (\id_\calQ\otimes \eta)=(\id_\calQ\otimes \eta)\star (\eta\otimes \id_\calQ),\]
where $\eta:\K\longrightarrow \calQ$ is the species morphism sending $\lambda \in \K$ to $\lambda 1_\calQ\in \calQ[\emptyset]$. Here, $\K$ is seen as a linear species concentrated in degree $0$, that is to say $\K[\emptyset]=\K$ and $\K[X]=\{0\}$ if $X\neq \emptyset$.
Therefore,
\[\ln(\Delta)=\ln(\eta \otimes \id_\calQ)+\ln(\id_\calQ\otimes \eta)=\eta\otimes \ln(\id_\calQ)+\ln(\id_\calQ)\otimes \eta.\]
Therefore, if $X$ and $Y$ are non-empty, as $\eta$ is zero on $\calQ[X]$ and $\calQ[Y]$, $\ln(\Delta)_{X,Y}=\Delta_{X,Y}\circ \ln(\id_\calQ)_{X\sqcup Y}=0$. Therefore, the image of $\ln(\id_\calQ)$ is included in $\prim(\calQ)$. The two first item imply that it is a projector on $\prim(\calQ)$. 
\end{proof}

\subsection{Polynomial invariants}

\begin{defi}\label{defi1.2.20}
Let $\calQ$ be a twisted bialgebra. A polynomial invariant of $\calQ$ is given by maps $\phi_X:\calQ[X]\longrightarrow \K[T]$ for any finite set $X$, such that:
\begin{itemize}
\item If $\sigma:X\longrightarrow Y$ is a bijection between two finite sets, then the following diagram commutes:
\[\xymatrix{\calQ[X]\ar[rr]^{\calQ[\sigma]}\ar[rd]_{\phi_X}&&\calQ[Y]\ar[ld]^{\phi_Y}\\ &\K[T]&}\]
\item For any $x\in \calQ[\emptyset]$, $\phi_\emptyset(x)=\varepsilon_\Delta(x)$.
\item If $X$, $Y$ are disjoint finite sets, $x\in\calQ[X]$ and $y\in\calQ[Y]$, then $\phi(xy)=\phi(x)\phi(y)$.
\item If $X$ is a finite set, and $x\in\calQ[X]$, then
\begin{align*}
\phi_X(x)(T+T')&=\sum_{X=X'\sqcup X''}(\phi_{X'}\otimes \phi_{X''})\circ \Delta_{X',X''}(x)(T,T').
\end{align*}\end{itemize}
By abuse of notations, we shall often simply write $\phi(x)$ instead of $\phi_X(x)$, for $X$ a finite set and $x\in \calQ[X]$.
\end{defi}

\begin{remark}\begin{enumerate}
\item If $\phi$ is a polynomial invariant on a twisted bialgebra, then for any $\lambda\in \K$, we define a character of $\calQ$, by sending any $b\in \calQ[X]$, where $X$ is a finite set, to $\phi(b)(\lambda)$.
\item If $\calQ$ is connected, the second item is equivalent to $\phi(1_\calQ)=1$.
\end{enumerate}\end{remark}

\begin{notation}
For $k\geq 1$, we put
\[\binom{T}{k}=\frac{T(T-1)\cdots (T-k+1)}{k!}.\]
By convention, $\displaystyle \binom{T}{0}=1$.
\end{notation}

\begin{lemma}[\textbf{Vandermonde's formula}] \label{lem1.2.21}
For any $k\geq 0$,
\[\binom{T+T'}{k}=\sum_{i=0}^k \binom{T}{i}\binom{T'}{k-i}.\]
\end{lemma}

\begin{proof}
By the classical Vandermonde's formula, if $n,n'\in \N$,
\[\binom{n+n'}{k}=\sum_{i=0}^k \binom{n}{i}\binom{n'}{k-i}.\]
As this is true for any $(n,n')\in \N^2$, the polynomial version of Vandermonde's formula is also true.
\end{proof}

\begin{lemma}\label{lem1.2.22}
Let $\calQ$ be a twisted bialgebra, and $f\in \calQ^*_1$. There exists a unique family of maps $\phi:\calQ[X]\longrightarrow \K[T]$ for any finite set $X$, such that:
\begin{itemize}
\item If $\sigma:X\longrightarrow Y$ is a bijection between two finite sets, then the following diagram commutes:
\[\xymatrix{\calQ[X]\ar[rr]^{\calQ[\sigma]}\ar[rd]_{\phi_X}&&\calQ[Y]\ar[ld]^{\phi_Y}\\ &\K[T]&}\]
\item For any $x\in \calQ[\emptyset]$, $\phi_\emptyset(x)=\varepsilon_\Delta(x)$.
\item If $X$ is a finite set, and $x\in\calQ[X]$, then
\begin{align*}
\phi(x)(T+T')&=\sum_{X=X'\sqcup X''}(\phi\otimes \phi)\circ \Delta_{X',X''}(x)(T,T').
\end{align*}
\item For any finite set $X$, and any $x\in \calQ[X]$, $\phi_X(x)(1)=f(x)$.
\end{itemize}
For any nonempty finite set $X$,
\begin{align}
\label{eqphi}
&\forall x\in \calQ[X],&\phi(x)&=\sum_{k=1}^{|X|}\sum_{(X_1,\ldots,X_k)\in \comp[X]} f^{\otimes k}\circ \Delta_{X_1,\ldots,X_k}(x) \binom{T}{k}.
\end{align}\end{lemma}

\begin{proof}
\textit{Existence}. Let us define $\phi$ as follows:
\begin{itemize}
\item If $x\in \calQ[\emptyset]$, then $\phi(x)=\varepsilon_\Delta(x)$.
\item If $X$ is a nonempty finite set and $x\in \calQ[X]$, $\phi(x)$ is defined by (\ref{eqphi}).
\end{itemize}
The two first items are clearly satisfied. Let us prove the third one. Firstly, let us assume that $X=\emptyset$.
\begin{align*}
\phi(x)(T+T')&=\varepsilon_\Delta(x),\\
\sum_{X=X_1\sqcup X''}(\phi\otimes \phi)\circ \Delta_{X_1,X''}(x)(T,T')&=(\varepsilon_\Delta\otimes \varepsilon_\Delta)\circ\Delta_{\emptyset,\emptyset}(x)(T,T')=\varepsilon_\Delta(x).
\end{align*}
Let us now assume that $X\neq \emptyset$. By Lemma \ref{lem1.2.21},
\begin{align*}
&\phi(x)(T+T')\\
&=\sum_{k=1}^{|X|}\sum_{(X_1,\ldots,X_k)\in \comp[X]} \sum_{i=0}^k f^{\otimes k}\circ \Delta_{X_1,\ldots,X_k}(x) \binom{T}{i}\binom{T'}{k-i}\\
&=\sum_{X=X'\sqcup X''} \sum_{i=0}^{|X'|} \sum_{j=0}^{|X''|}\sum_{\substack{(X'_1,\ldots,X'_i)\in \comp[X'],\\ (X''_1,\ldots,X''_j)\in \comp[X'']}} f^{\otimes (i+j)}\circ \Delta_{X'_1,\ldots,X'_i,X''_1,\ldots,X''_j}(x)\binom{T}{i}\binom{T'}{k-i}\\
&=\sum_{X=X'\sqcup X''} \sum_{i=0}^{|X'|} \sum_{j=0}^{|X''|}\sum_{\substack{(X'_1,\ldots,X'_i)\in \comp[X'],\\ (X''_1,\ldots,X''_j)\in \comp[X'']}} \left(f^{\otimes i}\circ \Delta_{X'_1,\ldots,X'_i}\otimes f^{\otimes j}\circ \Delta_{X''_1,\ldots,X''_j}\right)
\circ \Delta_{X',X''}(x)\binom{T}{i}\binom{T'}{k-i}\\
&=\sum_{X=X'\sqcup X''}(\phi\otimes \phi)\circ \Delta_{X',X''}(x)(T,T'),
\end{align*}
where $j=k-i$.

Let us now prove the fourth item. If $X=\emptyset$, $\phi(x)(1)=\varepsilon_\Delta(x)=f(x)$, as $f\in \calQ^*_1$. If $X\neq \emptyset$,
\[\phi(x)(1)=f\circ \Delta_X(x)+0=f(x).\]

\textit{Unicity}. Let us assume that $\phi$ and $\psi$ both satisfy the four items. Let us show that for any finite set $X$, $\phi_X=\psi_X$. We proceed by induction on $X$. If $X=\emptyset$, by the second item, $\phi_\emptyset=\psi_\emptyset$.
Let us assume the result for any finite set $X'$ of cardinality $<|X|$. By the second and third items,
\begin{align*}
\phi(x)(T+T')&=(\phi\otimes \phi)\circ \Delta_{X,\emptyset}(x)(T,T')+(\phi\otimes \phi)\circ \Delta_{\emptyset,X}(x)(T,T')\\
&+\sum_{(X',X'')\in \comp[X]}(\phi\otimes \phi)\circ \Delta_{X',X''}(x)(T,T')\\
&=(\phi\otimes \varepsilon_\Delta)\circ \Delta_{X,\emptyset}(x)(T,T')+(\varepsilon_\Delta\otimes \phi)\circ \Delta_{\emptyset,X}(x)(T,T')\\
&+\sum_{(X',X'')\in \comp[X]}(\phi\otimes \phi)\circ \Delta_{X',X''}(x)(T,T')\\
&=\phi(T)+\phi(T')+\sum_{(X',X'')\in \comp[X]}(\phi\otimes \phi)\circ \Delta_{X',X''}(x)(T,T'),
\end{align*}
and similarly, 
\[\psi(T+T')=\psi(T)+\psi(T')+\sum_{(X',X'')\in \comp[X]}(\psi\otimes \psi)\circ \Delta_{X',X''}(x)(T,T').\]
By the induction hypothesis, if $(X',X'')\in \comp[X]$, then $|X'|$, $|X''|<|X|$, so $\phi_{X"}=\psi_{X'}$ and $\phi_{X''}=\psi_{X''}$, which gives
\[\sum_{(X',X'')\in \comp[X]}(\phi\otimes \phi)\circ \Delta_{X',X''}(x)(T,T')=\sum_{(X',X'')\in \comp[X]}(\psi\otimes \psi)\circ \Delta_{X',X''}(x)(T,T').\]
Finally, $(\phi-\psi)(x)(T+T')=(\phi-\psi)(x)(T)+(\phi-\psi)(x)(T')$, which implies that there exists a scalar $\lambda$, such that $(\phi-\psi)(x)=\lambda T$. Evaluating in $1$, we obtain $\lambda=f(x)-f(x)=0$, so $\phi(x)=\psi(x)$, and finally $\phi_X=\psi_X$.
\end{proof}

\begin{theo}\label{theo1.2.23}
Let $\calQ$ be a twisted bialgebra, and $f\in \Char(\calQ)\cap\calQ^*_1$. There exists a unique polynomial invariant of $\calQ$, such that for any finite set $X$, for any $x\in \calQ[X]$, $\phi(x)(1)=f(x)$. Moreover, if $X$ is nonempty,
\[\phi(x)=\sum_{k=1}^{|X|}\sum_{(X_1,\ldots,X_k)\in \comp[X]} f^{\otimes k}\circ \Delta_{X_1,\ldots,X_k}(x) \binom{T}{k}.\]
\end{theo}

\begin{proof}
Lemma \ref{lem1.2.22} gives the unicity of $\phi$, and the existence of a $\phi$ satisfying the first, third, and fourth items of Definition \ref{defi1.2.20}.
It remains to prove the second item. Let us consider the twisted bialgebra $\calQ\otimes \calQ$. As $m$ is a coalgebra morphism, we obtain two family of maps, given by
\begin{align*}
\phi^{(1)}_X=m\circ \phi\otimes (\phi\otimes \phi):&(\calQ\otimes\calQ)_X\longrightarrow \K[T],\\
\phi^{(2)}_X=\phi\circ m:&(\calQ\otimes\calQ)_X\longrightarrow \K[T],
\end{align*}
both satisfying the first, third, and fourth items of Definition \ref{defi1.2.20}. Let $X$ and $Y$ be disjoint finite sets, $x\in\calQ[X]$, and $y\in \calQ[Y]$. Then
\begin{align*}
\phi^{(1)}_X(x\otimes y)(1)&=\phi(x)(1)\phi_y(1)=f(x)f(y),\\
\phi^{(2)}_X(x\otimes y)(1)&=\phi(xy)(1)=f(xy).
\end{align*}
As $f$ is a character, $\phi^{(1)}_X(x\otimes y)(1)=\phi^{(2)}_X(x\otimes y)(1)$. By unicity in Lemma \ref{lem1.2.22}, $\phi^{(1)}=\phi^{(2)}$, which gives the second item of Definition \ref{defi1.2.20}.
\end{proof}

In the particular where $\calQ$ is connected:

\begin{cor}
Let $\calQ$ be a connected twisted bialgebra, and $f\in \Char(\calQ)$. There exists a unique polynomial invariant of $\calQ$, such that for any finite set $X$, for any $x\in \calQ[X]$, $\phi(x)(1)=f(x)$. Moreover, if $X$ is nonempty,
\[\phi(x)=\sum_{k=1}^{|X|}\sum_{(X_1,\ldots,X_k)\in \comp[X]} f^{\otimes k}\circ \Delta_{X_1,\ldots,X_k}(x) \binom{T}{k}.\]
\end{cor}

\begin{prop}\label{prop1.2.25}
Let $\calQ$ be a twisted bialgebra, and $g\in \InfChar(\calQ)\calQ^*_0$. The following defines a polynomial invariant of $\calQ$: for any finite set $X$, and any $x\in \calQ[X]$,
\[\phi(x)=\sum_{k=1}^{|X|}g^{\otimes k}\circ \Delta_{X_1,\dots,X_k}(x) \frac{T^k}{k!}.\]
It is the polynomial invariant associated to $f=\exp_\calQ(g)$.
\end{prop}

\begin{proof}
Let $X$ be a nonempty finite set, and $x\in \calQ[X]$.
\begin{align*}
&\phi(x)(T+T')\\
&=\sum_{k=1}^{|X|}g^{\otimes k}\circ \Delta_{X_1,\dots,X_k}(x) \left(\sum_{i=0}^k \frac{T^i}{i!}\frac{T'^{k-i}}{(k-i)!}\right)\\
&=\sum_{X=X'\sqcup X''} \sum_{i=0}^{|X'|} \sum_{j=0}^{|X''|}\sum_{\substack{(X'_1,\ldots,X'_i)\in \comp[X'],\\ (X''_1,\ldots,X''_j)\in \comp[X'']}} g^{\otimes (i+j)}\circ \Delta_{X'_1,\ldots,X'_i,X''_1,\ldots,X''_j}(x)\frac{T^i}{j!}\frac{T^j}{j!}\\
&=\sum_{X=X'\sqcup X''} \sum_{i=0}^{|X'|} \sum_{j=0}^{|X''|}\sum_{\substack{(X'_1,\ldots,X'_i)\in \comp[X'],\\ (X''_1,\ldots,X''_j)\in \comp[X'']}} \left(g^{\otimes i}\circ \Delta_{X'_1,\ldots,X'_i}\otimes g^{\otimes j}\circ \Delta_{X''_1,\ldots,X''_j}\right)
\circ \Delta_{X',X''}(x)\frac{T^i}{j!}\frac{T^j}{j!}\\
&=\sum_{X=X'\sqcup X''}(\phi\otimes \phi)\circ \Delta_{X',X''}(x)(T,T'),
\end{align*}
where $j=k-i$. Moreover, 
\[\phi(x)(1)=\sum_{k=1}^{|X|}\frac{1}{k!}g^{\otimes k}\circ \Delta_{X_1,\dots,X_k}(x)=\exp_\calQ(g)(x). \]
We denote by $\psi$ the polynomial invariant associated to $\exp_\calQ(g)$. By unicity in Lemma \ref{lem1.2.22}, $\phi=\psi$.
\end{proof}

\begin{cor}\label{cor1.2.26}
Let $\calQ$ be a twisted bialgebra, and $f\in \Char(\calQ)\cap\calQ^*_1$. We denote by $\phi$ the polynomial invariant associated to $f$. For any finite set $X$, and any $x\in \calQ[X]$,
\[\ln_\calQ(f)(x)=\dfrac{\mathrm{d}\phi(x)}{\mathrm{d}T}(0).\]
\end{cor}

\begin{proof}
We put $g=\ln_\calQ(f)$. The result is immediate if $X$ is empty. Let us now assume that $X\neq \emptyset$. By Proposition \ref{prop1.2.25}, 
\[\phi(x)=\sum_{k=1}^{|X|}g^{\otimes k}\circ \Delta_{X_1,\dots,X_k}(x) \frac{T^k}{k!},\]
so
\[\dfrac{\mathrm{d}\phi(x)}{\mathrm{d}T}(0)=g\circ\Delta_X(x)+0=g(x). \qedhere\]
\end{proof}

\begin{cor}\label{cor1.2.27}
Let $\calQ$ be a twisted bialgebra, and $f\in \Char(\calQ)\cap\calQ^*_1$. We denote by $\phi$ the polynomial invariant associated to $f$ by Theorem \ref{theo1.2.23}. For any $k\in \Z$, for any finite set $X$, and any $x\in \calQ[X]$,
\[\phi(x)(k)=f^{*k}(x).\]
\end{cor}

\begin{proof}
For any $k\in \Z$, we denote by $f_k$ the character of $\calQ$ defined as follows: for any finite set $X$, and any $x\in \calQ[X]$, $f_k(x)=\phi(x)(k)$. Then $f_0=\varepsilon_\Delta=f^{*0}$ and $f_1=f=f^{*1}$. 
Moreover, by definition of a polynomial invariant, for any $k,l\in \Z$,
\[\phi(x)(k+l)=\sum_{X=X'\sqcup X''}(\phi\otimes \phi)\circ \Delta_{X',X''}(x)(k,l)=
\sum_{X=X'\sqcup X''}(f_k\otimes f_l)\circ \Delta_{X',X''}(x)=f_k*f_l(x).\]
Therefore, for any $k,l\in \Z$, $f_k*f_l=f_{k+l}$. An easy induction on $k$ proves that for any $k\in \N$, $f_k=f^{*k}$. Moreover, if $k\in \N$, 
\[f_k*f_{-k}=f_{-k}*f_k=f_0=\varepsilon_\Delta,\]
so $f_{-k}=f_k^{*(-1)}=f^{*(-k)}$.
\end{proof}

\section{Twisted double bialgebras}

We refer to \cite{Abe1980,Cartier2021} for usual references and results on bialgebras and Hopf algebras.

\subsection{Definition}

\begin{notation}
Let $X$ be a finite set. Then $\calE[X]$ is partially ordered by the inclusion. 
\begin{enumerate}
\item Let $\sim\in \calE[X]$, we denote by $\varpi_\sim:X\longrightarrow X/{\sim}$ the canonical surjection.
\item Let $\sim\in \calE[X]$. We identify $\calE[X/{\sim}]$ with $\{\sim'\in \calE[X]\mid \sim\subseteq \sim'\}$:
if $\sim\subseteq \sim' \in \calE[X]$, the equivalence induced on $X/{\sim}$ by $\sim'$ is denoted by $\overline{\sim}'$. 
In other words, for any $\overline{x}$, $\overline{y}\in X/{\sim}$,
\[\overline{x}\overline{\sim}'\overline{y}\Longleftrightarrow x\sim' y.\]
Note that with this identification, $\varpi_{\sim'}\circ \varpi_\sim=\varpi_{\sim'}$.
\end{enumerate}\end{notation}

\begin{defi}\label{defi1.3.1}
A twisted double bialgebra is a family $(\calQ,m,\Delta,\delta)$ such that:
\begin{enumerate}
\item $(\calQ,m,\Delta)$ is a twisted bialgebra.
\item 
\begin{enumerate}
\item For any finite set $X$ and any $\sim\in\calE[X]$, we have a map
\[\delta_\sim:\calQ[X]\longrightarrow\calQ[X/{\sim}]\otimes\calQ[X].\]
For any bijection $\sigma:X\longrightarrow Y$ between two finite sets, the following diagram commutes:
\[\xymatrix{\calQ[X]\ar[d]_{\calQ[\sigma]}\ar[r]^{\hspace{-5mm}\delta_\sim}&\calQ[X/{\sim}]\otimes\calQ[X]\ar[d]^{\calQ[\sigma/{\sim}]\otimes\calQ[\sigma]}\\
\calQ[Y]\ar[r]_{\hspace{-5mm}\delta_{\sim_\sigma}}&\calQ[Y/{\sim}_\sigma]\otimes\calQ[Y]}\]
\item (Coassociativity of $\delta$). For any finite set $X$, and any $\sim,\sim'\in\calE[X]$, such that $\sim\subseteq\sim'$,
\[\left(\delta_{\overline{\sim'}}\otimes\id_{\calQ[X]}\right)\circ\delta_\sim=\left(\id_{\calQ[X/{\sim}']}\otimes\delta_\sim\right)\circ\delta_{\sim'}.\]
For any finite set $X$, and any $\sim,\sim'\in\calE[X]$, such that we do not have $\sim\subseteq\sim'$,
\[\left(\id_{\calQ[X/{\sim}']}\otimes\delta_\sim\right)\circ\delta_{\sim'}=0.\]
\item (Counity for $\delta$). There exists $\epsilon_\delta \in \Char(\calQ)$, such that for any finite set $X$, 
\begin{align*}
&\forall\sim\in\calE[X],&\left(\id_{\calQ[X/{\sim}]}\otimes\epsilon_\delta\right)\circ\delta_\sim
&=\begin{cases}
\id_{\calQ[X]}\mbox{ if $\sim$ is the equality of $X$},\\
0\mbox{ otherwise}.
\end{cases}\\
&&\sum_{\sim\in\calE[X]}\left(\epsilon_\delta\otimes\id_{\calQ[X]}\right)\circ\delta_\sim&=\id_{\calQ[X]}.
\end{align*}
\item (Compatibility with the product $m$). For any pair $(X$, $Y)$ of disjoint finite sets, for any $x\in\calQ[X]$, for any $y\in\calQ[Y]$, for any $\sim\in\calQ(X\sqcup Y)$, putting $\sim_X=\sim\cap X^2$ and $\sim_Y=\sim\cap Y^2$,
\[\delta_\sim (xy)
=\begin{cases}
\delta_{\sim_X}(x)\delta_{\sim_Y}(y)\mbox{ if }\sim=\sim_X\sqcup\sim_Y,\\
0\mbox{ otherwise}.
\end{cases}\]
Moreover, denoting by $\sim_\emptyset$ the unique equivalence on $\emptyset$, $\delta_{\sim_\emptyset}(1_\calQ)=1_\calQ\otimes 1_\calQ$. 
\item (Compatibility with the coproduct $\Delta$). For any pair $(X$, $Y)$ of disjoint finite sets, for any $\sim_X\in\calE[X]$ and $\sim_Y\in\calE[Y]$,
\[\left(\Delta_{X/{\sim_X},Y/{\sim_Y}}\otimes\id_{\calQ[X\sqcup Y]}\right)\circ\delta_{\sim_X\sqcup\sim_Y}=m_{1,3,24}\circ\left(\delta_{\sim_X}\otimes\delta_{\sim_Y}\right)\circ\Delta_{X,Y},\] 
where $m_{1,3,24}$ sends $x_1\otimes x_2\otimes x_3\otimes x_4$ to $x_1\otimes x_3\otimes m_{X,Y}(x_2\otimes x_4)$.
\end{enumerate}\end{enumerate}\end{defi}

\begin{example}
We give $\lcom$ a second coproduct as follows: for any finite set $X$, and $\sim\in \calE[X]$, any $\lambda \in \K=\lcom[X]$,
\[\delta_\sim(\lambda)=\begin{cases}
\lambda(1\otimes 1)\mbox{ if $\sim$ is the equality of $X$},\\
0\mbox{ otherwise}.
\end{cases}\]
\end{example}

\subsection{Characters of a twisted double bialgebra}

\begin{prop}\label{prop1.3.2}
Let $(\calQ,m,\Delta,\delta)$ be a twisted double bialgebra. Then $\Char(\calQ)$ inherits a second product $*$, defined as follows: if $f,g\in \Char(\calQ)$, for any finite set $X$, 
\[(f*g)_X=\sum_{\sim\in \calE[X]}(f_{X/{\sim}}\otimes g_X)\circ \delta_\sim.\]
The unit is $\epsilon_\delta$. For any $f,g,h\in \Char(\calQ)$, 
\[(f\star g)*h=(f*h)\star (g*h).\]
\end{prop}

\begin{proof}
Let $f,g\in \Char(\calQ)$. Let us first prove that $f*g\in \Char(\calQ)$. Let $(X,Y)$ be a pair of disjoint finite sets,
$x\in\calQ[X]$, $y\in \calQ[Y]$. By the compatibility between $\delta$ and $m$,
\begin{align*}
(f*g)(xy)&=\sum_{\sim\in \calE[X\sqcup Y]} (f\otimes g)\circ \delta_\sim(xy)\\
&=\sum_{\substack{\sim_X\in \calE[X],\\ \sim_Y\in \calE[Y]}} (f\otimes g)(\delta_{\sim_X}(x)\delta_{\sim_Y}(y))\\
&=\sum_{\substack{\sim_X\in \calE[X],\\ \sim_Y\in \calE[Y]}} (f\otimes g)\circ \delta_{\sim_X}(x)(f\otimes g)\circ \delta_{\sim_Y}(y)\\
&=(f*g)(x)(f*g)(y).
\end{align*}
We use the fact that $f$ and $g$ are characters for the third equality. So $f*g\in \Char(\calQ)$.

Let $f\in \calQ$. For any finite set $X$, by definition of the convolution,
\begin{align*}
(\epsilon_\delta*f)_X&=\sum_{\sim\in \calE[X]}(\epsilon_\delta\otimes f)\circ \delta_\sim\\
&=f\circ \left(\sum_{\sim\in \calE[X]}\left(\epsilon_\delta\otimes \id_{\calQ[X]}\right)\circ \delta_\sim\right)\\
&=f\circ \id_{\calQ[X]}\\
&=f_X,\\
(f*\epsilon_\delta)_X&=\sum_{\sim\in \calE[X]}(f\otimes\epsilon_\delta)\circ \delta_\sim\\
&=f\circ \left(\sum_{\sim\in \calE[X]}\left(\id_{\calQ[X]}\otimes \epsilon_\delta\right)\circ \delta_\sim\right)\\
&=f,
\end{align*}
so $\epsilon_\delta$ is a unit for $*$. Let $f,g,h\in \Char(\calQ)$. The coassociativity of $\delta$ gives
\begin{align*}
(f*g)*h_X&=\sum_{\substack{\sim \in \calE[X],\\ \overline{\sim}'\in \calE[X/{\sim}]}} (f\otimes g\otimes h)\circ\left(\delta_{\overline{\sim}'}\otimes \id_{\calQ[X]}\right)\circ \delta_\sim\\
&=\sum_{\sim\subseteq \sim'\in \calE[X]}(f\otimes g\otimes h)\circ(\id_{\calQ[X/{\sim'}]} \otimes \delta_\sim)\circ \delta_{\sim'}\\
&=\sum_{\sim,\sim'\in \calE[X]} (f\otimes g\otimes h)\circ(\id_{\calQ[X/{\sim'}]} \otimes \delta_\sim)\circ \delta_{\sim'}\\
&=f*(g*h)_X.
\end{align*}
So $*$ is associative. By the compatibility between $\delta$ and $\Delta$, 
\begin{align*}
(f\star g)*h&=\sum_{\sim\in \calE[X]} \sum_{X/{\sim}=X'/{\sim'}\sqcup X''/{\sim''}} (f\otimes g\otimes h)\circ \left(\Delta_{X'/{\sim'},X''/{\sim''}}\otimes\id_{\calQ[X]}\right)\circ \delta_\sim\\
&=\sum_{X=X'\sqcup X''}\sum_{\substack{\sim'\in \calE[X'],\\ \sim''\in \calE[X'']}}(f\otimes g\otimes h)\circ m_{1,3,24}\circ (\delta_{\sim'}\otimes \delta_{\sim''})\circ \Delta_{X',X''}\\
&=\sum_{X=X'\sqcup X''}\sum_{\substack{\sim'\in \calE[X'],\\ \sim''\in \calE[X'']}}(f\otimes h\otimes g\otimes h)\circ
(\delta_{\sim'}\otimes \delta_{\sim''})\circ \Delta_{X',X''}\\
&=\sum_{X=X'\sqcup X''} ((f*h)\otimes (g*h))\circ \Delta_{X',X''}\\
&=(f*h)\star (g*h).
\end{align*}
We used the fact that $h$ is a character for the third equality. So $(f\star g)*h=(f*h)\star (g*h)$. 
\end{proof}

\begin{lemma}\label{lem1.3.3}
Let $\calQ$ be a double twisted bialgebra, and $f\in \Char(\calQ)$, invertible in the monoid $(\Hom(\calQ,\lcom),*)$. Then its inverse is also a character of $\calG$.
\end{lemma}

\begin{proof}
The inverse of $f$ is denoted by $g$. Then $f*g(1_\calQ)=f(1_\calQ)g(1_\calQ)=1$. As $f(1_\calQ)=1$, $g(1_\calQ)=1$. 
Note that $\calQ\otimes \calQ$ is also a twisted double bialgebra, and it is not difficult to show that in $(\Hom(\calQ\otimes \calQ,\lcom),*)$:
\begin{enumerate}
\item The inverse of $m\circ (f\otimes f)$ is $m\circ (g\otimes g)$.
\item The inverse of $f\circ m$ is $g\circ m$.
\end{enumerate}
As $f\in \Char(\calQ)$, $m\circ (f\otimes f)=f\circ m$, so $m\circ (g\otimes g)=g\circ m$. Therefore, $g\in \Char(\calQ)$. 
\end{proof}

\subsection{The polynomial invariant}

\begin{defi}\label{defi1.3.4}
Let $\calQ$ be a twisted double bialgebra. Let $\phi$ be a polynomial invariant of $(\calQ,m,\Delta)$.
\begin{enumerate}
\item We shall say that is is compatible with $\delta$ if for any finite set $X$,
\begin{align*}
&\forall x\in \calQ[X],&\phi(x)(TT')&=\sum_{\sim\in \calE[X]} (\phi\otimes \phi)\circ \delta_\sim(x)(T,T').
\end{align*}
\item We shall say that is a double polynomial invariant if it is compatible with $\delta$ and if for any finite set $X$, \begin{align*}
&\forall x\in \calQ[X],&\phi(x)(1)&=\epsilon_\delta(x),
\end{align*}
where $\epsilon_\delta$ is the counit of $\delta$.
\end{enumerate}\end{defi}

\begin{lemma}\label{lem1.3.5}
Let $(\calQ,m,\Delta,\delta)$ be a double bialgebra, and let $\phi$ be a polynomial invariant of $\calQ$. 
We denote by $f$ the character of $\calQ$ defined as follows: for any finite set $X$, and any $x\in \calQ[X]$,
\[f(x)=\phi(x)(1).\]
The following assertions are equivalent:
\begin{enumerate}
\item $\phi$ is compatible with $\delta$.
\item For any finite set $X$, and any $x\in \calQ[X]$, $\displaystyle \phi(x)=\sum_{\sim\in \calE[X]} (f\otimes \phi)\circ \delta_\sim(x)$.
\end{enumerate}\end{lemma}

\begin{proof}
$\Longrightarrow$. If so, then 
\[\phi(TT')=\sum_{\sim\in \calE[X]} (\phi\otimes \phi)\circ \delta_\sim(x)(T,T').\]
Evaluating in $T=1$, we obtain
\[\phi(T')=\sum_{\sim\in \calE[X]} (f\otimes \phi)\circ \delta_\sim(x)(T').\]

$\Longleftarrow$. We have to show that $\phi$ is compatible with $\delta$. Let us show that for any finite set $X$, and any $x\in \calQ[X]$, 
\[\phi(x)(nT')=\sum_{\sim\in \calE[X]} (\phi\otimes \phi)\circ \delta_\sim(x)(n,T').\]
We proceed by induction on $n$. If $n=1$, as $\phi(y)(1)=\epsilon(y)$ for any $y\in \calQ[Y]$, where $Y$ is a finite set, by definition of the counit of $\delta$,
\begin{align*}
\sum_{\sim\in \calE[X]} (\phi\otimes \phi)\circ \delta_\sim(x)(1,T')&=\sum_{\sim\in \calE[X]} (f\otimes \phi)\circ \delta_\sim(x)(T')=\phi(x)(T').
\end{align*}
Let us assume the result at rank $n$. As $\phi$ is compatible with the product and the coproduct,
\begin{align*}
\phi(x)((n+1)T')&=\phi(x)(nT'+1T')\\
&=\sum_{X=X'\sqcup X''}(\phi\otimes \phi)\circ \Delta_{X',X'}(x)(nT',1T')\\
&=\sum_{X=X'\sqcup X''}\sum_{\substack{\sim'\in \calE[X'],\\ \sim''\in \calE[X'']}}(\phi\otimes\phi\otimes \phi\otimes \phi)\circ (\delta_{\sim'}\otimes \delta_{\sim''})\circ \Delta_{X',X''}(x)(n,T',1,T')\\
&=\sum_{X=X'\sqcup X''}\sum_{\substack{\sim'\in \calE[X'],\\ \sim''\in \calE[X'']}}m_{1,3,24}\circ (\phi\otimes\phi\otimes \phi\otimes \phi)\circ (\delta_{\sim'}\otimes \delta_{\sim''})\circ \Delta_{X',X''}(x)(n,1,T')\\
&=\sum_{X=X'\sqcup X''}\sum_{\substack{\sim'\in \calE[X'],\\ \sim''\in \calE[X'']}} (\phi\otimes\phi\otimes \phi)\circ m_{1,3,24}\circ (\delta_{\sim'}\otimes \delta_{\sim''})\circ \Delta_{X',X''}(x)(n,1,T')\\
&=\sum_{\sim \in \calE[X]}\sum_{X/{\sim}=X'\sqcup X''} (\phi\otimes\phi\otimes \phi)\circ \left(\Delta_{X',X''}\otimes \id_{\calQ[X]}\right)\circ \delta_\sim(x)(n,1,T')\\
&=\sum_{\sim \in \calE[X]}(\phi\otimes \phi)\circ \delta_\sim(x)(n+1,T').
\end{align*}
We used the induction hypothesis and the case $n=1$ for the second equality. As this equality is true for any $n\geq 1$, we obtain the result.
\end{proof}

\begin{remark}
If $\phi$ is a polynomial invariant compatible with $\delta$, considering the associated character $f$, evaluating $T=1$ in Lemma \ref{lem1.3.5}, we obtain that for any finite set $X$,
\[f=\sum_{\sim\in \calE[X]} (f\otimes f)\circ \delta_\sim.\]
\end{remark}

\begin{theo}\label{theo1.3.6}
Let $\calQ$ be a twisted double bialgebra. The counit of $\delta$ is denoted by $\epsilon_\delta$. The polynomial invariant of $\calQ$ is the polynomial invariant associated to $\epsilon_\delta$. We denote it by $\phi$. Then $\phi$ is the unique double polynomial invariant of $\calF[\calQ]$.
\end{theo}

\begin{proof}
\textit{Unicity}. Let $\psi$ be a double polynomial invariant of $\calQ$. 
Then, for any finite set $X$, and any $x\in\calQ[X]$, $\psi(x)(1)=\epsilon_\delta(x)$. 
By The unicity in Theorem \ref{theo1.2.23}, $\phi=\psi$.

\textit{Existence}.
We apply Lemma \ref{lem1.3.5} with $f=\epsilon_\delta$. Let $X$ be a finite set, and $x\in \calQ[X]$. Then
\[\sum_{\sim\in \calE[X]}\left(\epsilon_\delta\otimes \phi\right)\circ\delta_\sim(x)=\phi(x),\]
as $\epsilon_\delta$ is the counit of $\delta$. Therefore, $\phi$ is compatible with $\delta$. 
\end{proof}

\begin{remark}
There are other polynomial invariants, compatible with the coproduct $\delta$, which are not double polynomial invariants. For example, take $\phi$ defined as follows: for any finite set $X$, and any $x\in \calQ[X]$, $\phi(x)=\varepsilon_\Delta(x)$. This comes from Lemma \ref{lem1.3.5}: 
\[\sum_{\sim\in \calE[X]}\left(\varepsilon_\Delta\otimes \phi)\circ \delta_\sim(x)\right)=\varepsilon_\Delta(x).\]
\end{remark}

\begin{prop}\label{prop1.3.7}
Let $\calQ$ be a twisted double bialgebra. Its polynomial invariant is denoted by $\phi$.
We denote by $\PI(\calQ)$ the set of polynomial invariants of $\calQ$. The following defines a right action of the monoid $(\Char(\calQ),*)$ on $\PI(\calQ)$:
if $\psi\in \PI(\calQ)$, and $f\in \Char(\calQ)$, for any finite set $X$,
\begin{align*}
(\psi\curvearrowleft f)_X&=\sum_{\sim\in \calE[X]} (\psi\otimes f)\circ \delta_\sim.
\end{align*}
Moreover, the two following maps are bijections, inverse one from the other:
\begin{align*}
\theta:&\left\{\begin{array}{rcl}
\Char(\calQ)&\longrightarrow&\PI(\calQ)\\
f&\longmapsto&\phi\curvearrowleft f,
\end{array}\right.&
\theta':&
\left\{\begin{array}{rcl}
\PI(\calQ)&\longrightarrow&\Char(\calQ)\\
\psi&\longmapsto&\psi(\cdot)(1).
\end{array}\right.
\end{align*}
\end{prop}

\begin{proof}
Let $\psi\in \PI(\calQ)$, $f\in \Char(\calQ)$. Let us first prove that $\psi'=\psi\curvearrowleft f\in \PI(\calQ)$. 
As $\delta$, $\psi$ and $f$ are all compatible with the product, $\psi'$ is also compatible with the product.
Let $X$ be a finite set, and $x\in \calB[X]$. If $X\neq \emptyset$, as $X/{\sim}\neq \emptyset$ for any $\sim\in \calE[X]$,
\begin{align*}
\varepsilon_\Delta\circ \psi'(x)&=\sum_{\sim\in \calE[X]} (\varepsilon_\Delta\circ \psi\otimes f)\circ \delta_\sim(x)
=\sum_{\sim\in \calE[X]} (\varepsilon_\Delta\otimes f)\circ \delta_\sim(x)
=0=\varepsilon_\Delta(x).
\end{align*}
If $X=\emptyset$, denoting by $\sim_\emptyset$ the unique equivalence on $\emptyset$,
\begin{align*}
\varepsilon_\Delta\circ \psi'(x)&=(\varepsilon_\Delta\circ \psi \otimes f)\circ \delta_{\sim_\emptyset}(x)
=(\varepsilon_\Delta \otimes f)\circ \delta_{\sim_\emptyset}(x)
=\varepsilon_\Delta(x)f(1)
=\varepsilon_\Delta(x).
\end{align*}
So $\psi'(x)(0)=\varepsilon_\Delta\circ \psi'(x)=\varepsilon_\Delta(x)$ in any case.

Let $(X,Y)$ be a pair of disjoint finite sets, and $x\in \calQ[X\sqcup Y]$. By the compatibility between $\delta$ and $\Delta$, and as $\phi$ is a polynomial invariant,
\begin{align*}
\psi'(x)(T+T')&=\sum_{\sim\in \calE[X]} (\psi\otimes f)\circ \delta_\sim(x)(T+T')\\
&=\sum_{\sim\in \calE[X]}\sum_{X/{\sim}=X'/{\sim'}\sqcup X''/{\sim''}} (\psi\otimes \psi\otimes f)\circ \left(\Delta_{X'/{\sim'},X''/{\sim''}}\otimes \id_{\calQ[X]}\right)\circ \delta_\sim(x)(T,T')\\
&=\sum_{X=X'\sqcup X''}\sum_{\substack{\sim'\in \calE[X'],\\ \sim''\in \calE[X'']}}
(\psi\otimes \psi\otimes f)\circ m_{1,3,24}\circ (\delta_{sim'}\otimes \delta_{\sim''})\circ \Delta_{X',X''}(x)(T,T')\\
&=\sum_{X=X'\sqcup X''}\sum_{\substack{\sim'\in \calE[X'],\\ \sim''\in \calE[X'']}}(\psi\otimes f\otimes \psi\otimes f)\circ 
\circ (\delta_{\sim'}\otimes \delta_{\sim''})\circ \Delta_{X',X''}(x)(T,T')\\
&=\sum_{X=X'\sqcup X''}(\psi'\otimes \psi')\circ \Delta_{X',X''}(x)(T,T').
\end{align*}
We used that $f$ is a character for the fourth equality. Therefore, $\psi'\in \PI(\calQ)$.

Let $\psi\in \PI(\calQ)$, $g,h\in \Char(\calQ)$. 
\begin{align*}
(\psi\curvearrowleft g)\curvearrowleft h_X&=\sum_{\substack{\sim \in \calE[X],\\ \overline{\sim}'\in \calE[X/{\sim}]}} (\psi\otimes g\otimes h)\circ\left(\delta_{\overline{\sim}'}\otimes \id_{\calQ[X]}\right)\circ \delta_\sim\\
&=\sum_{\sim\subseteq \sim'\in \calE[X]}(\psi\otimes g\otimes h)\circ(\id_{\calQ[X/{\sim'}]} \otimes \delta_\sim)\circ \delta_{\sim'}\\
&=\sum_{\sim,\sim'\in \calE[X]} (\psi\otimes g\otimes h)\circ(\id_{\calQ[X/{\sim'}]} \otimes \delta_\sim)\circ \delta_{\sim'}\\
&=\psi\curvearrowleft (g*h)_X.
\end{align*}

Theorem \ref{theo1.2.23} gives that $\theta'$ is bijective. 
Let $f\in \Char(\calQ)$. Let $X$ be a finite set, and $x\in \calQ[X]$. 
\[\phi\curvearrowleft f(x)(1)=\sum_{\sim\in \calE[X]} (\phi\otimes f)\circ \delta_\sim(x)(1)
=f\left(\sum_{\sim\in \calE[X]}\left(\epsilon_\delta\otimes \id_{\calQ[X]}\right)\circ \delta_\sim(x)\right)=f(x),\]
which proves that $\theta'\circ\theta=\id_{\Char(\calQ)}$. Therefore, the inverse of the bijection $\theta'$ is $\theta$. 
\end{proof}

\begin{cor}
Let $\calQ$ be a twisted double bialgebra. Its polynomial invariant is denoted by $\phi$.
For any character $f\in \Char(\calQ)$, the polynomial invariant $\phi_f$ associated to $f$ by Theorem \ref{theo1.2.23} is given as follows: for any finite set $X$, 
\begin{align*}
{\phi_f}_X&=\sum_{\sim\in \calE[X]} (\phi\otimes f)\circ \delta_\sim.
\end{align*}\end{cor}

\begin{proof}
This comes form $\phi_f=\theta'^{-1}(f)=\theta(f)=\phi\curvearrowleft f$.
\end{proof}

\subsection{Action of formal series on twisted double bialgebras}

\begin{prop}\label{prop1.3.9}
Let $\calQ$ be a connected twisted double bialgebra, $\calR$ be a twisted algebra, and $f:\calQ\longrightarrow\calR$ be an algebra morphism. For any formal series $P\in \K[[T]]$, for any finite set $X$,
\[P(f-\iota_{\calQ,\calR})=\sum_{\sim\in \calE[X]}(P(\epsilon_\delta-\varepsilon_\Delta)\otimes f)\circ \delta_\sim.\]
\end{prop}

\begin{proof}
Note that $f(1_\calQ)=\iota_{\calQ,\calR}(1_\calQ)=1\calR$ and $\varepsilon_\Delta(1_\calQ)=\epsilon_\delta(1_\calQ)=1$. Therefore, as $\calQ$ is connected, $P(f-\iota_{\calQ,\calR})$ and $P(\epsilon_\delta-\varepsilon_\Delta)$ exist.
Let us consider the following linear map:
\begin{align*}
\theta&:\left\{\begin{array}{rcl}
\K[T]&\longrightarrow&\Hom(\calQ,\calR)\\
P&\longmapsto&\displaystyle\theta(P)_X=\sum_{\sim\in \calE[X]}(P(\epsilon_\delta-\varepsilon_\Delta)\otimes f)\circ \delta_\sim.
\end{array}\right.
\end{align*}
Let $P$, $Q\in \K[T]$. 
\begin{align*}
&\theta(P)\star \theta(Q)_X\\
&=\sum_{X=X'\sqcup X''}\sum_{\substack{\sim'\in \calE[X'],\\ \sim''\in \calE[X'']}} m\circ (P(\epsilon_\delta-\varepsilon_\Delta)\otimes f)\circ (Q(\epsilon_\delta-\varepsilon_\Delta)\otimes f)\circ (\delta_{\sim'}\otimes \delta_{\sim''})\circ \Delta_{X,'X''}\\
&=\sum_{X=X'\sqcup X''}\sum_{\substack{\sim'\in \calE[X'],\\ \sim''\in \calE[X'']}} (P(\epsilon_\delta-\varepsilon_\Delta)\otimes Q(\epsilon_\delta-\varepsilon_\Delta)\otimes f)\circ m_{1,3,24}\circ (\delta_{\sim'}\otimes \delta_{\sim''})\circ \Delta_{X,'X''}\\
&=\sum_{\sim\in \calE[X]}\sum_{X/{\sim}=X'\sqcup X''} (P(\epsilon_\delta-\varepsilon_\Delta)\otimes Q(\epsilon_\delta-\varepsilon_\Delta)\otimes f)\circ (\Delta_{X',X''}\otimes \id_{\calQ})\circ \delta_\sim\\
&=\sum_{\sim\in \calE[X]}(P(\epsilon_\delta-\varepsilon_\Delta)\star Q(\epsilon_\delta-\varepsilon_\Delta)\otimes f)\circ \delta_\sim\\
&=\sum_{\sim\in \calE[X]}(PQ(\epsilon_\delta-\varepsilon_\Delta)\otimes f)\circ \delta_\sim\\
&=\theta(PQ)_X.
\end{align*}
We used that $f$ is an algebra morphism for the second equality, and the compatibility between $\Delta$ and $\delta$ for the third one. We obtain that $\theta$ is an algebra morphism. 
If $|X|<\val(P)$, then $P(\epsilon_\delta-\varepsilon_\Delta)_{X/{\sim}}=0$ for any $\sim\in \calE[X]$, as $|X/{\sim}|\leq |X|$, so $\theta(P)_X=0$. This gives that $\val(\theta(P))\geq \val(P)$, and that $\theta$ is continuous. 
Finally, if $X$ is a non-empty finite set, as $\varepsilon_\Delta=0$ on $X/{\sim}$ for any $\sim\in \calE[X]$, 
\[\theta(T)_X=\sum_{\sim\in \calE[X]} (\epsilon_\delta\otimes f)\circ \delta_\sim=f_X=(f-\iota_{\calQ,\calR})_X,\]
by definition of the counit of $\delta$. Moreover, denoting by $\sim_0$ the unique element of $\calE[\emptyset]$,
\[\theta(T)_\emptyset(1_\calB)=(\epsilon_\delta-\varepsilon_\Delta)(1_\calQ)f(1_\calQ)=0=(f-\iota_{\calQ,\calR})(1_\calQ),\]
so $\theta(T)=f-\iota_{\calQ,\calR}$. 
Consequently, $\theta$ sends any $P\in \K[[T]]$ to $P(f-\iota_{\calQ,\calR})$.
\end{proof}

For $\calQ=\calR$ and $f=\id_\calQ$, we obtain:

\begin{theo}\label{theo1.3.10}
Let $\calQ$ be a connected twisted double bialgebra. Its polynomial invariant is denoted by $\phi$. We shall use the following Sweedler's notation: if $X$ is a finite set, $\sim\in \calE[X]$ and $x\in \calQ[X]$,
\[\delta_\sim(x)=x/{\sim}\otimes x\mid\sim.\]
\begin{enumerate}
\item The antipode of $(\calQ,m,\Delta)$ is given as follows: for any non-empty finite set $X$, and any $x\in \calQ[X]$,
\[S_\calQ(x)=\sum_{\sim\in \calE[X]} \phi(x/{\sim})(-1) x\mid\sim.\]
\item The eulerian idempotent of $(\calQ,m,\Delta)$ is given as follows: for any non-empty finite set $X$, and any $x\in \calQ[X]$,
\[\ln(\id_\calQ)(x)=\sum_{\sim\in \calE[X]} \dfrac{\mathrm{d}\phi(x/{\sim})}{\mathrm{d}T}(0) x\mid\sim.\]
\end{enumerate}\end{theo}

\begin{proof}
We apply Proposition \ref{prop1.3.9} with $\calQ=\calR$ and $f=\id_\calQ$.

1. We take $f=(1+T)^{-1}$. It is enough to prove that for any finite set $X$, and any $x\in \calQ[X]$, 
\[f(1+(\epsilon_\delta-\varepsilon_\Delta))=\phi(x)(-1).\]
For any $k\geq 0$,
\[\binom{-1}{k}=\dfrac{(-1)(-2)\cdots (-k)}{k!}=(-1)^k.\]
Therefore, 
\begin{align*}
\phi(x)(-1)&=\sum_{k=1}^{|X|}\sum_{(X_1,\ldots,X_k)\in \comp_k[X]} (-1^k)\epsilon_\delta^{\otimes k}\circ \Delta_{X_1,\ldots,X_k}(x)=f(1+(\epsilon_\delta-\varepsilon_\Delta))(x).
\end{align*}

2. We take $f=\ln(1+T)$. It is enough to prove that for any finite set $X$, and any $x\in \calQ[X]$, 
\[\ln(1+(\epsilon_\delta-\varepsilon_\Delta)) =\dfrac{\mathrm{d}\phi(x)}{\mathrm{d}T}(0).\]
For any $k\geq 0$,
\[\dfrac{\displaystyle \mathrm{d}\binom{T}{k}}{\mathrm{d}T}(0)=\dfrac{(-1)\cdots(-k+1)}{k!}=\dfrac{(-1)^{k-1}}{k}.\]
Therefore, 
\begin{align*}
\frac{\mathrm{d}\phi(x)}{\mathrm{d}T}(0)&=\sum_{k=1}^{|X|}\sum_{(X_1,\ldots,X_k)\in \comp_k[X]} \dfrac{(-1)^{k-1}}{k} \epsilon_\delta^{\otimes k}\circ \Delta_{X_1,\ldots,X_k}(x)=\ln(1+(\epsilon_\delta-\varepsilon_\Delta))(x). \qedhere
\end{align*}
\end{proof}

We now apply this with $\calR=\lcom$. 

\begin{theo}
Let $\calQ$ be a connected twisted double bialgebra, and $f\in \Char(\calQ)$. Then
\[\ln(f)=\sum_{\sim\in \calE[X]} \dfrac{\mathrm{d}\phi(x/{\sim})}{\mathrm{d}T}(0)f(x\mid\sim).\]
\end{theo}

\begin{proof}
Apply Proposition \ref{prop1.3.9} with $\calR=\lcom$, and $P=\ln(1+T)$, and use the proof of the second item of Theorem \ref{theo1.3.10}. 
\end{proof}

\section{The bosonic Fock functor}

\subsection{Definition and compatibility with the Cauchy product}

\begin{defi}\cite{Aguiar2010} 
Let $\calQ$ be a linear species. For any $n\geq 0$, the space $\coinv\calQ[n]$ is 
\[\coinv\calQ[n]=\frac{\calQ[n]}{\vect(x-\calQ[\sigma](x)\mid x\in \calQ[n],\: \sigma \in \sym_n)}.\]
The bosonic Fock functor $\calF$ is defined as follows:
\begin{itemize}
\item If $\calQ$ is a linear species, then 
\[\calF[\calQ]=\bigoplus_{n=0}^\infty \coinv\calQ[n].\]
\item If $f:\calQ\longrightarrow\calR$ is a linear species morphism, then $\calF[f]:\calF[\calQ]\longrightarrow \calF[\calR]$ is the linear map defined by
\begin{align*}
&\forall n\geq 0,\: \forall \overline{x}\in \coinv\calQ[n],&\calF[f](\overline{x})&=\overline{f(x)}\in \coinv\calR[n].
\end{align*} \end{itemize}
This defines a functor from the category of linear species to the category of vector spaces.
\end{defi}

\begin{defi} 
Let $\calQ$ be a linear species, and $X$ be a finite set, whose cardinality is denoted by $n$.
Let $\sigma:X\longrightarrow [n]$ be a bijection. For any $x\in \calQ[n]$, the element $\overline{\calQ[\sigma](x)}\in \coinv\calQ[n]$ does not depend on the choice of $\sigma$. This element will be called the isoclass of $x$.
\end{defi}

Therefore, if $\calQ$ is a linear species, for any finite set $X$, we obtain a map $x\longrightarrow \overline{x}$
from $\calQ[X]$ to $\calF[\calQ]$, sending $x$ to its isoclass. 
In particular, if $\calQ$ and $\calR$ are two linear species, if $X$ and $Y$ are two disjoint finite sets, we define a map $\calQ[X]\otimes \calR[Y]$ to $\calF[\calQ]\otimes \calF[\calR]$, sending $x\otimes y$ to $\overline{x}\otimes \overline{y}$. This induces an isomorphism 
\[\calF[\calQ\otimes\calR]\approx \calF[\calQ]\otimes \calF[\calR].\]
We now identify $\calF[\calQ\otimes\calR]$ and $\calF[\calQ]\otimes \calF[\calR]$. As a consequence:
\begin{enumerate}
\item If $(\calQ,m)$ is a twisted algebra, then $(\calF[\calQ],\calF[m])$ is an algebra.
If $X$, $Y$ are disjoint finite sets,
\begin{align*}
&\forall x\in \calQ[X],\: \forall y\in \calQ[Y],&\calF[m](\overline{x}\otimes \overline{y})&=\overline{m_{X,Y}(x\otimes y)}.
\end{align*}
The unit is $\overline{1_\calQ}$. With an abuse of notations, the product of $\calF[\calQ]$ will be denoted by $m$, as the product of $\calQ$. 
\item If $(\calQ,\Delta)$ is a twisted coalgebra, then $(\calF[\calQ],\calF[\Delta])$ is a coalgebra:
If $X$ is a finite set,
\begin{align*}
&\forall x\in \calQ[X],&\Delta(\overline{x})&=\sum_{X=X'\sqcup X''}\overline{\Delta_{X',X''}(x)}.
\end{align*}
The counit is given by
\begin{align*}
&\forall x\in \calQ[X],&\varepsilon_\Delta(\overline{x})=\calF[\varepsilon_\Delta](\overline{x})&=\begin{cases}
\varepsilon_\Delta(x)\mbox{ if }x\in \calQ[\emptyset],\\
0\mbox{ otherwise}.
\end{cases}\end{align*}
 With an abuse of notations, the coproduct of $\calF[\calQ]$ will be denoted by $\Delta$, as the coproduct of $\calQ$. 
\item If $(\calQ,m,\Delta)$ is a twisted bialgebra, then $(\calF[\calQ],m,\Delta)$ is a bialgebra.
\item If $(\calQ,m,\Delta)$ is a twisted Hopf algebra, then $(\calF[\calQ],m,\Delta)$ is a Hopf algebra. Its antipode is $\calF[S_\calQ]$.
\end{enumerate}

In the sequel, we shall consider linear species $\calQ$ which are linearizations of a set species $\calB$.
In other terms, for any finite set $X$, $\calQ[X]$ has a basis $\calB[X]$, such that for any bijection $\sigma:X\longrightarrow Y$ between two finite sets, $\calQ[\sigma](\calB[X])=\calB[Y]$.
In this case, a basis of $\calF[\calQ]$ is given by the set of isoclasses of elements $b\in \calB[X]$,
 where $X$ is any finite set.

\begin{example}\label{ex1.4.1} We consider $\calF[\lcom]$. For any finite set $X$, of cardinality denoted by $n$, we denote by $T_n$ the isoclass of $1\in \lcom[X]$. Then $(T_n)_{n\geq 0}$ is a basis of $\calF[\lcom]$. 
Moreover, for any $k,l\geq 0$,
\[T_k\cdots T_l=T_{k+l}.\]
For any $n\geq 0$,
\[\Delta(T_n)=\sum_{k=0}^n \binom{n}{k} T_k\otimes T_{n-k}.\] 
As a consequence, $\calF[\lcom]$ is isomorphic to the polynomial algebra $(\K[T],m,\Delta)$, with its usual product $m$ and its multiplicative coproduct $\Delta$ defined by
\[\Delta(T)=T\otimes 1+1\otimes T.\]
With the preceding notations, $T=T_1$ and $1=T_0$. 
In particular, if $\calQ$ is a twisted bialgebra, polynomial invariants of $\calQ$ are in fact morphism of bialgebras from $(\calF[\calQ],m,\Delta)$ to $(\K[T],m,\Delta)$.
\end{example}

\begin{example}
Let us consider $\calF[\ell\calS]$. For any finite set $X$, of cardinality denoted by $n$, and any $Y\in \calS[X]$, of cardinality denoted by $k$, we denote by $T_{n,k}$ the isoclass of $Y$ in $\cal[\ell\calS]$. 
Then $(T_{n,k})_{0\leq k\leq n}$ is a basis of $\calF[\ell\calS]$. Moreover, if $0\leq k\leq m$ and $0\leq l\leq n$,
\[T_{m,k}\cdots T_{n,l}=T_{m+n,k+l}.\]
In particular, if $0\leq k\leq n$,
\[T_{n,k}=T_{1,1}^k T_{1,0}^{n-k},\]
So $\calF[\ell\calS]=\K[T_{1,1},\: T_{1,0}]$, with the coproduct defined by
\begin{align*}
\Delta(T_{1,1})&=T_{1,1}\otimes 1+1\otimes T_{1,1},&
\Delta(T_{1,0})&=T_{1,0}\otimes 1+1\otimes T_{1,0}.
\end{align*}\end{example}

\begin{example}
Let us consider $\calF[\ell\comp]$. For any $n_1,\ldots,n_k\geq 1$, let us denote by $\Sigma_{n_1,\ldots,n_k}$
the isoclass of any composition $(X_1,\ldots,X_k)$, with $|X_i|=i$ for $i\in [k]$. Then $(\Sigma_{n_1,\ldots,n_k})_{k\geq 0,n_1,\ldots,n_k\geq 1}$ is a basis of $\calF[\ell\comp]$. Moreover,
if $m_1,\ldots,m_k,n_1,\ldots,n_l\geq 1$,
\[\Sigma_{m_1,\ldots,m_k}\cdot \Sigma_{n_1,\ldots,n_l}=\Sigma_{m_1,\ldots,m_k,n_1,\ldots,n_l}.\]
If $n\geq 1$,
\[\Delta(\Sigma_n)=1\otimes \Sigma_n+\sum_{k=1}^{n-1} \binom{n}{k}\Sigma_k\otimes \Sigma_{n-k}+\Sigma_n\otimes 1,\]
so $\calF[\ell\comp]$ is isomorphic to the Hopf algebra $\mathbf{NSym}$ of non-symmetric functions \cite{Duchamp2002,Gelfand1995,Hazewinkel2005}.
\end{example}

\subsection{Bosonic Fock functor of a twisted double bialgebra}

Let $(\calQ,m,\Delta,\delta)$ be a twisted double bialgebra. We define a second coproduct $\calF[\delta]$ on $\calF[\calQ]$ as follows: if $X$ is a finite set,
\begin{align*}
&\forall x\in \calQ[X],&\calF[\delta](\overline{x})&=\sum_{\sim\in \calE[X]} \overline{\delta(x)}.
\end{align*}
 With an abuse of notations, this coproduct of $\calF[\calQ]$ will be denoted by $\delta$, as in $\calQ$. 
Then $(\calF[\calQ],m,\Delta,\delta)$ is a double bialgebra, that is to say:
\begin{enumerate}
\item $(\calF[\calQ],m,\Delta)$ is a bialgebra.
\item $(\calF[\calQ],m,\delta)$ is a bialgebra.
\item $\Delta:\calF[\calQ]\longrightarrow \calF[\calQ]\otimes \calF[\calQ]$ is a right comodule morphism:
\[\left(\Delta\otimes \id_{\calF[\calQ]}\right)\circ \delta=m_{1,3,24}\circ (\delta\otimes \delta)\circ \Delta.\]
\item $\varepsilon_\Delta:\calF[\calQ]\longrightarrow \K$ is a right comodule morphism:
\begin{align*}
&\forall \overline{x}\in \calF[\calQ],&\left(\varepsilon_\Delta \otimes\id_{\calF[\calQ]}\right)\circ \delta(\overline{x})=\varepsilon_\Delta(\overline{x})1_{\calF[\calQ]}.
\end{align*}\end{enumerate} 

\begin{example}
Let us consider the second coproduct $\delta$ on $\K[T]$ induced by the twisted double structure of $\calB$, with the notations of Example \ref{ex1.4.1}. We obtain that $\delta$ is the multiplicative product given by
\[\delta(T)=T\otimes T.\]
It is a very classical example of double bialgebra.
\end{example}

\begin{remark}
Let $\calQ$ be a twisted bialgebra. We saw that a polynomial invariant is a morphism of bialgebras from $(\calF[\calQ],m,\Delta)$ to $(\K[T],m,\Delta)$.
If $\calQ$ is a double twisted bialgebra, then a double polynomial invariant is a morphism of double bialgebras from $(\calF[\calQ],m,\Delta,\delta)$ to $(\K[T],m,\Delta,\delta)$, and a polynomial invariant compatible with $\delta$ is 
a morphism of bialgebras from $(\calF[\calQ],m,\Delta)$ to $(\K[T],m,\Delta)$,
 which is also a morphism of non counitary coalgebras from $(\calF[\calQ],\delta)$ to $(\K[T],\delta)$. 
\end{remark}

\section{The twisted bialgebra of compositions and its terminal property}

\subsection{Coalgebraic structure on $\comp$}

The species $\comp$ is defined in Example \ref{ex1.1.1}. Its linearization is denoted by $\ell\comp$.

\begin{prop}
The linear species $\ell\comp$ is a twisted coalgebra, with the coproduct defined as follows:
if $X$ and $Y$ are disjoint finite sets, and $\pi=(X_1,\ldots,X_k)\in \comp[X\sqcup Y]$, then 
\[\Delta_{X,Y}(\pi)=\begin{cases}
(X_1,\ldots,X_i)\otimes (X_{i+1},\ldots,X_k)\mbox{ if $X_1\sqcup \cdots\sqcup X_i=X$},\\
0\mbox{ otherwise}.
\end{cases}\]
\end{prop}

\begin{proof}
Note that, as $X_1\sqcup\cdots\sqcup X_k=X\sqcup Y$, the coproduct can be rewritten as
\[\Delta_{X,Y}(\pi)=\begin{cases}
(X_1,\ldots,X_i)\otimes (X_{i+1},\ldots,X_k)\mbox{ if }X_1\sqcup \cdots\sqcup X_i=X\\
\hspace{5cm}\mbox{ and }X_{i+1}\sqcup\cdots\sqcup X_k=Y,\\
0\mbox{ otherwise}.
\end{cases}\]

Let $X$, $Y$, $Z$ be pairwise disjoint finite sets, and $\pi=(X_1,\ldots,X_k)\in \comp[X\sqcup Y\sqcup Z]$. Then 
\begin{align*}
\left(\Delta_{X,Y}\otimes \id_{\ell\comp[Z]}\right)\circ\Delta_{X\sqcup Y,Z}(\pi)&=\begin{cases}
(X_1,\ldots,X_i)\otimes (X_{i+1},\ldots X_j)\otimes (X_{j+1},\ldots,X_k)\\
\hspace{2cm}\mbox{ if }X_1\sqcup \cdots\sqcup X_i=X\\
\hspace{2cm}\mbox{ and }X_{i+1}\sqcup\cdots\sqcup X_j=Y\\
\hspace{2cm}\mbox{ and }X_{j+1}\sqcup\cdots\sqcup X_k=Z,\\
0\mbox{ otherwise}
\end{cases}\\
&=\left(\id_{\ell\comp[X]}\otimes \Delta_{Y,Z}\right)\circ\Delta_{X,Y\sqcup Z}(\pi).
\end{align*}
So $\Delta$ is coassociative. If $X$ is a finite set, and $\pi\in \comp[X]$,
\begin{align*}
\Delta_{\emptyset,X}(\pi)&=1\otimes \pi,&\Delta_{X,\emptyset}(\pi)&=\pi\otimes 1,
\end{align*}
where $1$ is the unique element of $\comp[\emptyset]$. So $\ell\comp$ has a counit, defined by $\varepsilon_\Delta(1)=1$.
\end{proof}

\begin{lemma}\label{lem1.5.2}
Let $\calQ$ a a twisted coalgebra, so that $\calQ[\emptyset]$ is one-dimensional. There exists a unique element $1_\calQ\in \calQ[\emptyset]$, such that $\varepsilon_\Delta(1_\calQ)=1$. Moreover, for any finite set $X$, and any $x\in \calQ[X]$.
\begin{align*}
\Delta_{X,\emptyset}(x)&=x\otimes 1_\calQ,&\Delta_{\emptyset,X}(x)&=1_\calQ\otimes x.
\end{align*}
\end{lemma}

\begin{proof}
As $\calQ$ is not empty, $\varepsilon_\Delta:\calQ[\emptyset]\longrightarrow\K$ is not zero. As $\calQ[\emptyset]$ is one-dimensional, $\varepsilon_\Delta$ is bijective. Therefore, there exists a unique element $1_\calQ\in \calQ[\emptyset]$, such that $\varepsilon_\Delta(1_\calQ)=1$. 

Let $X$ be a finite set, and $x\in \calQ[X]$, nonzero.
As $\calQ[\emptyset]$ is one-dimensional, there exists a scalar $\lambda\in \K$, such that
\[\Delta_{\emptyset,X}(x)=\lambda 1_\calQ\otimes x.\]
Therefore, 
\[x=\left(\varepsilon_\Delta\otimes \id_{\calQ[\emptyset]}\right)\circ\Delta_{\emptyset,X}(x)=\lambda x,\]
so $\lambda=1$. The proof is similar for $\Delta_{X,\emptyset}(x)$.
\end{proof}

\begin{remark}
In particular, $\Delta_{\emptyset,\emptyset}(1_\calQ)=1_\calQ\otimes 1_\calQ$.
\end{remark}

\begin{defi}
We define a species morphism $\epsilon:\ell\comp\longrightarrow\ell\Com$ as follows: for any finite set $X$,
\begin{align*}
&\forall \pi\in \comp[X],&\epsilon(\pi)&=\begin{cases}
1\mbox{ if }\ell(\pi)\leq 1,\\
0\mbox{ otherwise}.
\end{cases}\end{align*}\end{defi}

\begin{prop}\label{prop1.5.4}
Let $\calQ$ be a twisted coalgebra, and let $f:\calQ\longrightarrow\ell\Com$ be a species morphism such that $f(1_\calQ)=1$. There exists a unique twisted coalgebra morphism $\phi:\calQ\longrightarrow \ell\comp$ such that 
$\epsilon\circ\phi=f$. 
\end{prop}

\begin{proof}
\textit{Existence}. We put $\phi(1_\calQ)=1$ and, if $X$ is a nonempty finite set and $x\in \calQ[X]$. We put
\[\phi(x)=\sum_{k=0}^{|X|} \sum_{(X_1,\ldots,X_k)\in \comp_k[X]} f^{\otimes k}\circ \Delta_{X_1,\ldots,X_k}(x)(X_1,\ldots,X_k).\]
Let us prove that $\phi$ is a twisted coalgebra morphism. Let $(X,Y)$ be a pair of disjoint finite sets, and $x\in \calQ[X\sqcup Y]$. 
It $X=\emptyset$ then
\begin{align*}
\Delta_{\emptyset,Y}\circ \phi(x)&=1\otimes \phi(x)=\phi(1_\calQ)\otimes \phi(x)=(\phi\otimes \phi)\circ \Delta_{\emptyset,Y}(x).
\end{align*}
The proof is similar if $X=\emptyset$. Let us assume that both $X$ and $Y$ are nonempty.
\begin{align*}
&\Delta_{X,Y}\circ \phi(x)\\
&=\sum_{k=1}^{|X|+|Y|} \sum_{(X_1,\ldots,X_k)\in \comp_k[X]} f^{\otimes k}\circ \Delta_{X_1,\ldots,X_k}(x)\Delta_{X,Y}(X_1,\ldots,X_k)\\
&=\sum_{i=1}^{|X|}\sum_{j=1}^{|Y|} \sum_{\substack{(X_1,\ldots,X_i)\in \comp[X],\\
(Y_1,\ldots,Y_j)\in \comp[Y]}}f^{\otimes (i+j)}\circ \Delta_{X_1,\ldots,X_i,Y_1,\ldots,Y_j}(x)(X_1,\ldots,X_i)\otimes (Y_1,\ldots,Y_j)\\
&=\sum_{i=1}^{|X|}\sum_{j=1}^{|Y|} \sum_{\substack{(X_1,\ldots,X_i)\in \comp[X],\\
(Y_1,\ldots,Y_j)\in \comp[Y]}}f^{\otimes (i+j)}\circ \left(\Delta_{X_1,\ldots,X_i}\otimes\Delta_{Y_1,\ldots,Y_j}\right)\circ \Delta_{X,Y}(x)(X_1,\ldots,X_i)\otimes (Y_1,\ldots,Y_j)\\
&=(\phi\otimes \phi)\circ \Delta_{X,Y}(x).
\end{align*}

Let us now prove that $\epsilon\circ \phi=f$. Let $X$ be a finite set, and $x\in \calQ[X]$. If $X=\emptyset$, we can assume that $x=1_\calQ$. Then 
\[\epsilon\circ \phi(1_\calQ)=\epsilon(1)=1=f(1_\calQ).\]
Otherwise, 
\[\epsilon\circ\phi(x)=f\otimes \Delta_X(x)\epsilon(X)+0=f(x).\]

\textit{Unicity}. Let $\phi$, $\psi:\calQ\longrightarrow\ell\comp$ be two twisted coalgebra morphisms such that $\epsilon\circ \phi=\epsilon\circ\psi$. Let us prove that for any finite set $X$, and for any $x\in \calQ[X]$, $\phi(x)=\psi(x)$. 
We proceed by induction on $|X|$. If $X=\emptyset$, we can assume that $x=1_\calQ$. Then 
\[\varepsilon_\Delta\circ \phi(1_\calQ)=\varepsilon_\Delta(1_\calQ)=1,\]
so $\phi(1_\calQ)=1$. Similarly, $\psi(1_\calQ)=1$. Let us assume the result at all ranks $<|X|$. If $(X',X'')\in \comp_2[X]$, the induction hypothesis, gives that 
\begin{align*}
\Delta_{X',X''}\circ \phi(x)&=(\phi\otimes \phi)\circ \Delta_{X',X''}(x)=(\psi\otimes \psi)\circ \Delta_{X',X''}(x)=\Delta_{X',X''}\circ \psi(x),
\end{align*}
so $\Delta_{X,'X''}(\phi(x)-\psi(x))=0$. We deduce that there exists $\lambda \in \K$, such that $\phi(x)-\psi(x)=\lambda X$.
Applying $\epsilon$,
\[\lambda=\epsilon\circ\phi(x)-\epsilon\circ\psi(x)=0,\]
so $\phi(x)=\psi(x)$.
\end{proof}

\subsection{Quasi-shuffle product on $\ell\comp$}

\begin{theo}\label{theo1.5.5}
We define a species morphism $\epsilon_2:\ell\comp\otimes\ell\comp_2\longrightarrow\Com$ as follows:
if $(X,Y)$ is a pair of disjoint finite sets, $(X_1,\ldots,X_k)\in \comp[X]$ and $(Y_1,\ldots,Y_l)\in \comp[Y]$, then 
\[\epsilon_2((X_1,\ldots,X_k)\otimes (Y_1,\ldots,Y_l))=(\epsilon\otimes \epsilon)((X_1,\ldots,X_k)\otimes (Y_1,\ldots,Y_l))
=\begin{cases}
1\mbox{ if }k\leq 1\mbox{ and }l\leq 1,\\
0\mbox{ otherwise}.
\end{cases}\]
We denote by $\squplus:\ell\comp\otimes\ell\comp\longrightarrow\ell\comp$ the unique twisted coalgebra morphism such that $\epsilon\circ \squplus=\epsilon_2$ (Proposition \ref{prop1.5.4}).
Then $(\ell\comp,\squplus,\Delta)$ is a commutative twisted bialgebra.
Moreover, $\epsilon$ is a character of this twisted bialgebra. 
\end{theo}

\begin{proof}
Let us first prove the associativity of $\squplus$. We consider the species morphisms
\begin{align*}
\phi&=\squplus\circ (\id_{\ell\comp}\otimes \squplus),&\psi&=\squplus\circ(\squplus\otimes \id_{\ell\comp}).
\end{align*}
By composition, both are twisted coalgebra morphisms from $\ell\comp^{\otimes 3}$ to $\ell\comp$. 
Let $X$, $Y$, $Z$ be pairwise disjoint finite sets, $\pi\in \comp[X]$, $\pi'\in \comp[Y]$, and $\pi''\in \comp[Z]$,
of respective lengths $k$, $l$ and $m$. Then
\begin{align*}
\epsilon\circ \phi(\pi\otimes \pi'\otimes \pi'')&=\epsilon_2(\pi\squplus \pi'\otimes \pi'')
=\epsilon(\pi\squplus\pi')\epsilon(\pi'')=\epsilon_2(\pi\otimes \pi')\epsilon(\pi'')=\epsilon(\pi)\epsilon(\pi')\epsilon(\pi'').
\end{align*}
Similarly, $\epsilon\circ \psi(\pi\otimes \pi'\otimes \pi'')=\epsilon(\pi)\epsilon(\pi')\epsilon(\pi'')$.
By Proposition \ref{prop1.5.4} (unicity part), $\phi=\psi$, so $\squplus$ is associative.

Let us denote by $\tau:\ell\comp\otimes\ell\comp\longrightarrow\ell\comp\otimes\ell\comp$ be the usual flip.
We consider the species morphism $\squplus^{\mathrm{op}}=\squplus\circ \tau$. By composition, it is a twisted coalgebra morphism from $\ell\comp^{\otimes 2}$ to $\ell\comp$. Let $(X,Y)$ be a pair of disjoint finite sets, $\pi\in \comp[X]$, and $\pi'\in \comp[Y]$, of respective lengths $k$ and $k$. Then 
\begin{align*}
\epsilon\circ \squplus^{\mathrm{op}}(\pi\otimes \pi')&=\epsilon\circ \squplus(\pi'\otimes \pi)
=\epsilon_2(\pi'\otimes \pi)=\epsilon(\pi')\epsilon(\pi)=\epsilon_2(\pi\otimes \pi')=\epsilon\circ \squplus.
\end{align*}
By Proposition \ref{prop1.5.4} (unicity), $\squplus=\squplus^{\mathrm{op}}$, so $\squplus$ is commutative.

Let us consider $\theta:\ell\comp\longrightarrow \ell\comp$, which sends $x\in \ell\comp[X]$ to $1\squplus x$.
By composition, it is a twisted coalgebra endomorphism. Moreover, for any $x\in \ell\comp[X]$,
\[\epsilon\circ \theta(x)=\epsilon(1_\calQ)\epsilon(x)=\epsilon(x).\]
By Proposition \ref{prop1.5.4} (unicity), $\theta=\id_{\ell\comp}$, so $1$ is the unit of $\ell\comp$.
We proved that $(\ell\comp,\squplus,\Delta)$ is a commutative twisted bialgebra.

By construction, $\epsilon\circ \squplus=\epsilon\otimes\epsilon$, so $\epsilon$ is a character of $\ell\comp$.
\end{proof}

\begin{remark}
It is possible to combinatorially describe this product $\squplus$ with the help of quasi-shuffles, as we will do for $\exch$ later. 
For example, if $X$, $Y$, $Z$ are pairwise disjoint sets,
\begin{align*}
(X)\squplus (Y)&=(X,Y)+(Y,X)+(X\sqcup Y),\\
(X)\squplus (Y,Z)&=(X,Y,Z)+(X\sqcup Y,Z)+(Y,X,Z)+(Y,X\sqcup Z)+(Y,Z,X).
\end{align*}
\end{remark}

Let us now prove the terminal property of the twisted bialgebra $\ell\comp$.

\begin{theo}\label{theo1.5.6}
Let $\calQ$ be a twisted bialgebra, such that $\calQ[\emptyset]$ is one-dimensional, and let $f:\calQ\longrightarrow \Com$ be a species morphism, such that $f(1_\calQ)=1$. We denote by $\phi_f$ the unique twisted coalgebra morphism from $\calQ$ to $\ell\comp$ such that $\epsilon\circ\phi_f=f$. 
Then $\phi_f$ is a morphism of twisted bialgebras if, and only if, $f$ is a character of $\calQ$.
\end{theo}

\begin{proof}
Note that $1_\calQ$, defined as in Lemma \ref{lem1.5.2}, is necessarily the unit of $\calQ$.

$\Longrightarrow$. By composition, $\epsilon\circ \phi_f=f$ is a twisted algebra morphism from $\calQ$ to $\Com$, so is a character.

$\Longleftarrow$. Let us assume that $f$ is a character of $\calQ$. We consider the two species morphisms
\begin{align*}
\Phi&=\phi_f\circ m_{\calQ},&\Psi&=\squplus\circ (\phi_f\otimes \phi_f).
\end{align*}
By composition, both are twisted coalgebra morphisms from $\calQ\otimes \calQ$ to $\ell\comp$.
Let $(X,Y)$ be a pair of disjoint finite sets, $x\in \calQ[X]$ and $y\in \calQ[Y]$. 
\begin{align*}
\epsilon\circ \Phi(x\otimes y)&=\epsilon \circ \phi_f\circ m_{\calQ}(x\otimes y)\\
&=f \circ m_{\calQ}(x\otimes y)\\
&=f(xy),\\
\epsilon\circ \Psi(x\otimes y)&=\epsilon\circ\squplus\circ (\phi_f\otimes \phi_f)(x\otimes y)\\
&=(\epsilon\otimes \epsilon)\circ (\phi_f\otimes \phi_f)(x\otimes y)\\
&=f(x)f(y).
\end{align*}
As $f$ is a character, $\epsilon\circ \Phi=\epsilon\circ \Psi$. By Proposition \ref{prop1.5.4} (unicity), $\Phi=\Psi$,
so $\phi_f$ is a morphism of twisted bialgebras from $\calQ$ to $\ell\comp$.
\end{proof}

\begin{remark}
Applying the Fock functor, we obtain the bialgebra $\calF[\comp]=\QSym$ of quasi-symmetric functions, which the quasi-shuffle product and the deconcatenation coproduct, with a character $\epsilon$.
It satisfies also a terminal property, in the category of graded and connected bialgebras \cite{Aguiar2006-2}.
\end{remark}

\chapter{The category of combinatorial twisted bialgebras}

\section*{Introduction}

We now introduce combinatorial twisted bialgebras (shortly, CTBs). They are twisted bialgebras obtained by linearization of a set species $\calB$, with the following constraints:
\begin{itemize}
\item For any pair $(X,Y)$ of disjoint finite sets, for any $a \in \calB[X]$, and $b\in \calB[Y]$, the product of $ab$ belongs to $\calB[X\sqcup Y]$, or is zero.
\item For any pair $(X,Y)$ of disjoint finite sets, for any $b\in \calB[X\sqcup Y]$, the coproduct $\Delta_{X,Y}(b)$ is the tensor of an element of $\calB[X]$ and an element of $\calB[Y]$, or is zero.
\end{itemize}
See Definition \ref{defi2.1.1} for the full details. We introduce the main examples which will be used in this text:
the CTBs of graphs (Example \ref{ex2.1.2}), of oriented graphs (Example \ref{ex2.1.3}), of hypergraphs (Example \ref{ex2.1.4}, with four different coproducts), of posets and of quasi-posets (Examples \ref{ex2.1.5} and \ref{ex2.1.8}).
CTBs form a category, with the notion of morphism introduced in Definition \ref{defi2.1.2}: a morphism between two CTBs $\calB$ and $\calB'$ is a morphism of twisted bialgebras, sending any element of $\calB$ to an element of $\calB'$ or to $0$. We give several examples of morphisms, from oriented graphs to quasi-posets (Example \ref{ex2.1.13}), from quasi-posets to posets (Example \ref{ex2.1.14}), from hypergraphs to boolean functions (Example \ref{ex2.1.1.15}). We also introduce the notion of biideal of a CTB (Definition \ref{defi2.1.6}), which allows to define quotients and to associate a kernel to any morphism of CTBs, see Proposition \ref{prop2.1.8}.

We study the algebraic structure of a CTB in the second section. We introduce the notion of indecomposable elements (Definition \ref{defi2.2.1}), which will play a crucial role in the sequel, and show that any element of a CTB $\calB$ can be written as a product of indecomposable elements, uniquely up to the order of the factors (Theorem \ref{theo2.2.4}). 
For graphs, posets, quasi-posets, this corresponds to the well-known decomposition into connected components.
This also shows that, as a twisted algebra, any CTB is generated by its indecomposable elements, with commutation relations between some of them.
If $X$ is a nonempty finite set and $b\in \calB[X]$, then $1\leq \ic(b)\leq |X|$: elements for which $\ic$ is maximal are called discrete (Definition \ref{defi2.2.7}), and discrete elements form a sub-CTB $\bfD\calB$ of $\calB$, generated in degree 1. This sub-CTB will have an important role for double CTBs, where it will be related to the counit for the second product in most cases, as we shall see in Chapter \ref{chap4}.

\section{Combinatorial twisted bialgebras}

\subsection{Definition and examples}

\begin{defi}\label{defi2.1.1}
A combinatorial twisted bialgebra (briefly, CTB) is a triple $(\calB,m,\Delta)$ such that:
\begin{enumerate}
\item $\calB$ is a species.
\item $(\ell\calB,m,\Delta)$ is a twisted bialgebra.
\item $\calB[\emptyset]=\{1_\calB\}$, where $1_\calB$ is the unit of $\calB$.
\item If $X$ and $Y$ are two disjoint finite sets, and $a\in \calB[X]$, $b\in \calB[Y]$, then $m_{X,Y}(a\otimes b)=0$ or belongs to $\calB[X\sqcup Y]$.
 \item If $X$ and $Y$ are two disjoint finite sets, and $b\in \calB[X\sqcup Y]$, $\Delta_{X,Y}(b)=0$ or is of the form $b_{\mid_LX}\otimes b_{\mid_RY}$, with $b_{\mid_LX}\in \calB[X]$ and $b_{\mid_RY}\in \calB[Y]$.
\end{enumerate}\end{defi}

\begin{remark}
If $\calB$ is a CTB, as $\calB[\emptyset]=\{1_\calB\}$, its counit is uniquely determined by $\varepsilon_\Delta(1_\calB)=1$.
\end{remark}

\begin{example}\label{ex2.1.1}
For any finite set $X$, we put $\Com[X]=\{X\}$. It is a CTB with, for any disjoint finite sets $X$ and $Y$,
\begin{align*}
m_{X,Y}(X\otimes Y)&=X\sqcup Y,&\Delta_{X,Y}(X\sqcup Y)&=X\otimes Y.
\end{align*}
Its linearization is, up to an isomorphism sending $X\in \Com[X]$ to $1\in \K=\lcom[X]$, the twisted bialgebra $\lcom$ of Example \ref{ex1.2.1}.
\end{example}

\begin{example}\label{ex2.1.2} A graph $G$ is a pair $(V(G),E(G))$, where $V(G)$ is a finite set and $E(G)$ a subset of $\calS_2[V(G)]$. The elements of $E(G)$ are called the edges of $G$. For any finite set $X$, $\calG[X]$ is the set of graphs $G$ such that $V(G)=X$. 
This defines a species $\calG$. For example,
\begin{align*}
\calG[1]&=\{\grdun{$1$}\},\\
\calG[2]&=\left\{\begin{array}{c}
\grdun{$1$}\grdun{$2$}, \grddeux{$1$}{$2$}
\end{array}\right\},\\
\calG[3]&=\left\{\begin{array}{c}
\grdun{$1$}\grdun{$2$}\grdun{$3$},\grdun{$1$}\grddeux{$2$}{$3$},\grdun{$2$}\grddeux{$1$}{$3$},\grdun{$3$}\grddeux{$1$}{$2$},\grdtroisdeux{$1$}{$3$}{$2$},\grdtroisdeux{$2$}{$3$}{$1$},\grdtroisdeux{$3$}{$2$}{$1$},\grdtroisun{$1$}{$2$}{$3$}
\end{array}\right\}.
\end{align*}
 If $X$ is a finite set,
\[|\calG[X]|=2^{\binom{|X|}{2}}.\]
Let $X$ and $Y$ be disjoint finite sets. If $G\in \calG[X]$ and $H\in \calG[Y]$, we define $GH\in \calG[X\sqcup Y]$ by
\[E(GH)=E(G)\sqcup E(H).\]
This defines a product $m$ on $\calG$. If $Y\subseteq X$ are finite sets and $G\in \calG[X]$, we define $G_{\mid Y}\in \calG[Y]$ by
\[E(G_{\mid Y})=\{e\in E(G)\mid e\subseteq Y\}.\]
We then define a coproduct on $\calG$ as follows: if $X$, $Y$ are disjoint finite sets,
\begin{align*}
&\forall G\in \calG[X\sqcup Y],&\Delta_{X,Y}(G)&=G_{\mid X}\otimes G_{\mid Y}.
\end{align*}
Then $(\calG,m,\Delta)$ is a CTB. It is commutative and cocommutative.
The bialgebra $\calF[\calG]$ is studied in \cite{Foissy36}. 
\end{example}

\begin{example}\label{ex2.1.3}
 An oriented graph $G$ is a pair $(V(G),A(G))$, where $V(G)$ is a finite set and $A(G)$ a subset of 
\[\{(x,y)\in V(G)^2\mid x\neq y\}.\]
The elements of $A(G)$ are called the arcs of $G$. For any finite set $X$, $\calO\calG[X]$ is the set of oriented graphs $G$ such that $V(G)=X$. This defines a species $\calO\calG$. If $X$ is a finite set,
\[|\calO\calG[X]|=2^{|X|(|X|-1)}.\]
Let $X$ and $Y$ be disjoint finite sets. If $G\in \calO\calG[X]$ and $H\in \calO\calG[Y]$, we define $GH\in \calO\calG[X\sqcup Y]$ by
\[A(GH)=A(G)\sqcup A(H).\]
This defines a product on $\calO\calG$. If $Y\subseteq X$ are finite sets, and $G\in \calO\calG[X]$, we define $G_{\mid Y}\in \calO\calG[Y]$ by
\[A(G_{\mid Y})=\{e\in A(G)\mid e\in Y^2\}.\]
We shall say that $Y$ is a lower part of $G$ if 
\begin{align*}
&\forall x,x'\in X,&(x,x')\in A(G)\mbox{ and }x'\in X&\Longrightarrow x\in X.
\end{align*}
We then define a coproduct on $\calO\calG$ as follows: if $X$, $Y$ are disjoint finite sets,
\begin{align*}
&\forall G\in \calO\calG[X\sqcup Y],&\Delta_{X,Y}(G)&=\begin{cases}
G_{\mid X}\otimes G_{\mid Y}\mbox{ if $X$ is a lower part of $G$},\\
0\mbox{ otherwise}.
\end{cases}\end{align*}
Then $(\calO\calG,m,\Delta)$ is a CTB. It is commutative and not cocommutative.
The bialgebra $\calF[\calO\calG]$ is studied in \cite{Manchon2012}.
\end{example}

\begin{example}\label{ex2.1.4} A hypergraph $G$ is a pair $(V(G),\HE(G))$, where $V(G)$ is a finite set and $\HE(G)$ a subset of $\calS_{\geq 1}[V(G)]$, see \cite{Berge89,Bretto2013} for more details and applications. 
The elements of $\HE(G)$ are called the hyperedges of $G$, and the elements of $\HE(G)$ of cardinality $\geq 2$ are called the nontrivial hyperedges of $G$. For any finite set $X$, $\calH\calG[X]$ is the set of hypergraphs $G$ such that $V(G)=X$. 
This defines a species $\calH\calG$. If $X$ is a finite set,
\[|\calH\calG[X]|=2^{2^{|X|}-1}.\]
Let $X$ and $Y$ be disjoint finite sets. If $G\in \calH\calG[X]$ and $H\in \calH\calG[Y]$, we define $GH\in \calH\calG[X\sqcup Y]$ by
\[\HE(GH)=\HE(G)\sqcup \HE(H).\]
This defines a product $m$ on $\calH\calG$. If $Y\subseteq X$ are finite sets, and $G\in \calH\calG[X]$, we define $G_{\mid Y}$ and $G_{\mid_\cap Y}\in \calH\calG[Y]$ by
\begin{align*}
\HE(G_{\mid Y})&=\{e\in \HE(G)\mid e\subseteq Y\},&\HE(G_{\mid_\cap Y})&=\{e\cap Y\mid e\in \HE(G)\}\setminus\{\emptyset\}.
\end{align*}
We then define four coproducts on $\calH\calG$ as follows: if $X$, $Y$ are disjoint finite sets, 
\begin{align*}
&\forall G\in \calH\calG[X\sqcup Y],&\Delta_{X,Y}^{\subset,\subset}(G)&=G_{\mid X}\otimes G_{\mid Y},&\Delta_{X,Y}^{\cap,\cap}(G)&=G_{\mid_\cap X}\otimes G_{\mid_\cap Y},\\
&&\Delta_{X,Y}^{\subset,\cap}(G)&=G_{\mid X}\otimes G_{\mid_\cap Y},&\Delta_{X,Y}^{\cap,\subset}(G)&=G_{\mid_\cap X}\otimes G_{\mid Y}.
\end{align*}
Then $(\calH\calG,m,\Delta^{\subset,\subset})$, $(\calH\calG,m,\Delta^{\cap,\cap})$, $(\calH\calG,m,\Delta^{\subset,\cap})$, and $(\calH\calG,m,\Delta^{\cap,\subset})$,
are CTBs. All are commutative. The two first ones are cocommutative. The two last ones are not, and are co-opposite one from the other.
The bialgebras $\calF[\calH\calG]$, with it various coproducts, is studied in \cite{Ebrahimi-Fard2022,Foissy44}.
\end{example}

\begin{example}\label{ex2.1.5} A poset is a pair $P=(V(P),\leq_P)$ where $V(P)$ is a finite set and $\leq_P$ is an order on $V(P)$, that is to say a transitive, reflexive, and antisymmetric relation on $V(P)$.
For any finite set $X$, $\calP[X]$ is the set of posets $P$ such that $V(P)=X$. This defines a species $\calP$. 
For example, representing posets by their Hasse diagrams,
\begin{align*}
\calP[1]&=\{\tdun{$1$}\},\\
\calP[2]&=\left\{\begin{array}{cl}
\tdun{$1$}\tdun{$2$},\tddeux{$1$}{$2$},\tddeux{$2$}{$1$}
\end{array}\right\},\\
\calP[3]&=\left\{\begin{array}{cl}
\tdun{$1$}\tdun{$2$}\tdun{$3$},\tdun{$1$}\tddeux{$2$}{$3$},\tdun{$1$}\tddeux{$3$}{$2$},\tdun{$2$}\tddeux{$1$}{$3$},\tdun{$2$}\tddeux{$3$}{$1$},\tdun{$3$}\tddeux{$1$}{$2$},\tdun{$3$}\tddeux{$2$}{$1$},\\
\tdtroisun{$1$}{$3$}{$2$},\tdtroisun{$2$}{$3$}{$1$},\tdtroisun{$3$}{$2$}{$1$},
\pdtroisun{$1$}{$2$}{$3$},\pdtroisun{$2$}{$1$}{$3$},\pdtroisun{$3$}{$1$}{$2$},
\tdtroisdeux{$1$}{$2$}{$3$},\tdtroisdeux{$1$}{$3$}{$2$},\tdtroisdeux{$2$}{$1$}{$3$},
\tdtroisdeux{$2$}{$3$}{$1$},\tdtroisdeux{$3$}{$1$}{$2$},\tdtroisdeux{$3$}{$2$}{$1$}
\end{array}\right\}.
\end{align*}

The following array gives the first values of $|\calP[X]|$, see Entry A001035 of the OEIS \cite{Sloane}.
\[\begin{array}{|c||c|c|c|c|c|c|c|c|c|c|}
\hline |X|&1&2&3&4&5&6&7&8&9&10\\
\hline\hline |\calP[X]|&1&3&19&219&4231&130023&6129859&431723379&44511042511&6611065248783\\
\hline\end{array}\]
Let $X$ and $Y$ be disjoint finite sets. If $P\in \calP[X]$ and $Q\in \calP[Y]$, we define $PQ\in \calP[X\sqcup Y]$ by
\[\leq_{PQ}=\leq_P\sqcup \leq_Q.\]
This defines a product $m$ on $\calP$. If $Y\subseteq X$ are finite sets, and $P\in \calP[X]$, we define $P_{\mid Y}\in \calP[Y]$ by
\[\leq_{P_{\mid Y}}=\leq_P \cap Y^2.\]
We shall say that $Y$ is a lower part of $P$ if for any $x,x'\in X$, 
\[x\leq_P x'\mbox{ and }x'\in X\Longrightarrow x\in X.\]
We then define a coproduct on $\calP$: if $X$, $Y$ are disjoint finite sets,
\begin{align*}
&\forall P\in \calP[X\sqcup Y],&\Delta_{X,Y}(P)&=\begin{cases}
P_{\mid X}\otimes P_{\mid Y}\mbox{ if $X$ is a lower part of $P$},\\
0\mbox{ otherwise}.
\end{cases}\end{align*}
Then $(\calP,m,\Delta)$ is a CTB. It is commutative and not cocommutative. The bialgebra $\calF[\calP]$ is studied in \cite{Foissy37}.
\end{example}

\begin{example}
A poset $P$ is a rooted forest if it does not contain a sub-poset isomorphic to the poset $\ptroisun$. The subspecies of rooted forests is a sub-CTDB of $\calP$, denoted by $\calF$.
It is the twisted version of the Butcher-Connes-Kreimer Hopf algebra \cite{Brouder2004,Connes1998,Calaque2011}. It is commutative and not cocommutative. For example,
\begin{align*}
\calF[1]&=\{\tdun{$1$}\},\\
\calF[2]&=\left\{\begin{array}{cl}
\tdun{$1$}\tdun{$2$},\tddeux{$1$}{$2$},\tddeux{$2$}{$1$}
\end{array}\right\},\\
\calF[3]&=\left\{\begin{array}{cl}
\tdun{$1$}\tdun{$2$}\tdun{$3$},\tdun{$1$}\tddeux{$2$}{$3$},\tdun{$1$}\tddeux{$3$}{$2$},\tdun{$2$}\tddeux{$1$}{$3$},\tdun{$2$}\tddeux{$3$}{$1$},\tdun{$3$}\tddeux{$1$}{$2$},\tdun{$3$}\tddeux{$2$}{$1$},\\
\tdtroisun{$1$}{$3$}{$2$},\tdtroisun{$2$}{$3$}{$1$},\tdtroisun{$3$}{$2$}{$1$},
\tdtroisdeux{$1$}{$2$}{$3$},\tdtroisdeux{$1$}{$3$}{$2$},\tdtroisdeux{$2$}{$1$}{$3$},
\tdtroisdeux{$2$}{$3$}{$1$},\tdtroisdeux{$3$}{$1$}{$2$},\tdtroisdeux{$3$}{$2$}{$1$}
\end{array}\right\}.
\end{align*}
For any finite set $X$, $|\calF[X]|=(|X|+1)^{|X|-1}$, see Entry A000272 in \cite{Sloane}.
\end{example}

\begin{example}
A planar rooted forest is a pair formed by a rooted forest and an embedding of its Hasse graph in the plane (that is to say, the roots are totally ordered from left to right, and the descendants of any vertex are also totally ordered from left to right). This defines a species denoted by $\calP\calF$. For example,
 \begin{align*}
\calP\calF[1]&=\{\tdun{$1$}\},\\
\calP\calF[2]&=\left\{\begin{array}{cl}
\tdun{$1$}\tdun{$2$},\tdun{$2$}\tdun{$1$},\tddeux{$1$}{$2$},\tddeux{$2$}{$1$}
\end{array}\right\},\\
\calP\calF[3]&=\left\{\begin{array}{cl}
\tdun{$1$}\tdun{$2$}\tdun{$3$},\tdun{$1$}\tdun{$3$}\tdun{$2$},\tdun{$2$}\tdun{$1$}\tdun{$3$},\tdun{$2$}\tdun{$3$}\tdun{$1$},\tdun{$3$}\tdun{$1$}\tdun{$2$},\tdun{$3$}\tdun{$2$}\tdun{$1$},\\
\tdun{$1$}\tddeux{$2$}{$3$},\tdun{$1$}\tddeux{$3$}{$2$},\tdun{$2$}\tddeux{$1$}{$3$},\tdun{$2$}\tddeux{$3$}{$1$},\tdun{$3$}\tddeux{$1$}{$2$},\tdun{$3$}\tddeux{$2$}{$1$},
\tddeux{$2$}{$3$}\tdun{$1$},\tddeux{$3$}{$2$}\tdun{$1$},\tddeux{$1$}{$3$}\tdun{$2$},\tddeux{$3$}{$1$}\tdun{$2$},\tddeux{$1$}{$2$}\tdun{$3$},\tddeux{$2$}{$1$}\tdun{$3$},\\
\tdtroisun{$1$}{$3$}{$2$},\tdtroisun{$1$}{$2$}{$3$},\tdtroisun{$2$}{$3$}{$1$},\tdtroisun{$2$}{$1$}{$3$},\tdtroisun{$3$}{$2$}{$1$},\tdtroisun{$3$}{$1$}{$2$},
\tdtroisdeux{$1$}{$2$}{$3$},\tdtroisdeux{$1$}{$3$}{$2$},\tdtroisdeux{$2$}{$1$}{$3$},
\tdtroisdeux{$2$}{$3$}{$1$},\tdtroisdeux{$3$}{$1$}{$2$},\tdtroisdeux{$3$}{$2$}{$1$}
\end{array}\right\}.
\end{align*}
For any finite set $X$, $|\calP\calF[X]|=\dfrac{(2|X|)!}{(|X|+1)!}$, see Entry A001761 of \cite{Sloane}.
The product of two planar forests is their concatenation. If $F$ is a planar forest and $I\subseteq V(F)$ is a lower part or an upper part of $F$, then $F_{\mid I}$ is naturally a planar forest, and we define the coproduct of $\calP\calF$ as follows: if $X$, $Y$ are disjoint finite sets, 
\begin{align*}
&\forall F\in \calP\calF[X\sqcup Y],&\Delta_{X,Y}(F)&=\begin{cases}
F_{\mid X}\otimes F_{\mid Y}\mbox{ if $X$ is a lower part of $F$},\\
0\mbox{ otherwise}.
\end{cases}\end{align*}
This is the twisted version of the noncommutative Butcher-Connes-Kreimer bialgebra, see \cite{Foissy1,Foissy2,Holtkamp2003}. 
\end{example}

\begin{example}\label{ex2.1.8} A quasi-poset is a pair $P=(V(P),\leq_P)$ where $V(P)$ is a finite set and $\leq_P$ is a quasi-order on $V(P)$, that is to say a transitive, and reflexive relation on $V(P)$.
For any finite set $X$, $\calQ\calP[X]$ is the set of quasi-posets $P$ such that $V(P)=X$. This defines a species $\calQ\calP$. By Alexandroff's theorem \cite{Alexandroff1937}, this species is isomorphic to the species of finite topologies. For example,
\begin{align*}
\calQ\calP[1]&=\calP[1],\\
\calQ\calP[2]&=\calP[2]\sqcup\{\tun_{\tiny 12}\:\},\\
\calQ\calP[3]&=\calP[3]\sqcup\left\{\begin{array}{c}
\tun_{123},\tun_1\tun_{23},\tun_2\tun_{13},\tun_3\tun_{13},\tdeux_{12}^3,\tdeux_{13}^2,\tdeux_{23}^1,\tdeux_{13}^2,
\tdeux_1^{23},\tdeux_2^{13},\tdeux_3^{12}
\end{array}\right\}.
\end{align*}
In this representation, a vertex with several indices $x_1,\ldots,x_k$ means that for any $p,q\in [k]$, $x_p\leq_P x_q$.
The following array gives the first values of $|\calQ\calP[X]|$, see Entry \href{https://oeis.org/A000798}{A000798} of the OEIS \cite{Sloane}.
\[\begin{array}{|c||c|c|c|c|c|c|c|c|c|c|}
\hline |X|&1&2&3&4&5&6&7&8&9&10\\
\hline\hline |\calQ\calP[X]|&1&4&29&355&6942&209527&9535241&642779354&63260289423&8977053873043\\
\hline\end{array}\]
Let $X$ and $Y$ be disjoint finite sets. If $P\in \calQ\calP[X]$ and $Q\in \calQ\calP[Y]$, we define $PQ\in \calQ\calP[X\sqcup Y]$ by
\[\leq_{PQ}=\leq_P\sqcup \leq_Q.\]
This defines a product $m$ on $\calQ\calP$. If $Y\subseteq X$ are finite sets, and $P\in \calQ\calP[X]$, we define $P_{\mid Y}\in \calQ\calP[Y]$ by
\[\leq_{P_{\mid Y}}=\leq_P \cap Y^2.\]
We shall say that $Y$ is a lower part of $P$ if for any $x,x'\in X$, 
\[x\leq_P x'\mbox{ and }x'\in X\Longrightarrow x\in X.\]
We then define a coproduct on $\calQ\calP$: if $X$, $Y$ are disjoint finite sets,
\begin{align*}
&\forall P\in \calQ\calP[X\sqcup Y],&\Delta_{X,Y}(P)&=\begin{cases}
P_{\mid X}\otimes P_{\mid Y}\mbox{ if $X$ is a lower part of $P$},\\
0\mbox{ otherwise}.
\end{cases}\end{align*}
Then $(\calQ\calP,m,\Delta)$ is a CTB. It is commutative and not cocommutative. The bialgebra $\calF[\calQ\calP]$ is studied in \cite{Foissy37}.
\end{example}

\begin{example}
Let $X$ be a finite set. A boolean function on $X$ is a map $f:\calS[X]\longrightarrow \Z$, such that $f(\emptyset)=0$. This defines a species $\calB\calF$. To the difference with the preceding examples,
if $X$ is not empty, then $\calB\calF[X]$ is infinite.
If $X$, $Y$ are disjoint finite sets, $f\in \calB\calF[X]$ and $g\in \calB\calF[Y]$, we define $fg\in \calB\calF[X\sqcup Y]$ by
\begin{align*}
&\forall A\subseteq X\sqcup Y,& fg(A)&=f(A\cap X)+g(A\cap Y).
\end{align*}
If $h\in \calB\calF[X\sqcup Y]$, we put
\[\Delta_{X,Y}(f)=f_{\mid \calS[X]}\otimes g_{\mid \calS[Y]}.\]
With these product and coproduct, $\calB\calF$ is a CTB, which has been studied in \cite{Foissy57}. Identifying a matroid with its rank function, the species of matroids becomes a sub-CTB of $\calB\calF$. 
\end{example}

Other examples (or at least, their Fock functors) can be found on mixed graphs \cite{Foissy45}, and noncrossing partitions \cite{Ebrahimi-Fard2020,Foissy43}.

\subsection{Morphisms}

\begin{defi}\label{defi2.1.2}
Let $\calB$ and $\calB'$ be CTBs.
A morphism of CTBs from $\calB$ to $\calB'$ is a morphism $\phi:\ell\calB\longrightarrow \ell\calB'$ of twisted bialgebras such that for any finite set $X$, any $b\in \calB[X]$, $\phi(b)\in \calB'[X]\sqcup \{0\}$.
\end{defi}

In this way, CTBs form a category. 

\begin{example}\label{ex2.1.10}
The following maps are morphisms of CTBs:
\begin{align*}
\phi_{\calG,\Com}&:\left\{\begin{array}{rcl}
\calG&\longrightarrow&\Com\\
G\in \calG[X]&\longmapsto&X,
\end{array}\right.&\phi_{\calH\calG,\Com}&:\left\{\begin{array}{rcl}
\calH\calG&\longrightarrow&\Com\\
G\in \calH\calG[X]&\longmapsto&X.
\end{array}\right.
\end{align*}
\end{example}

\begin{example}
Graphs are particular hypergraphs, and this induces a canonical injection $\iota_\calG:\calG\longrightarrow \calH\calG$. It is not difficult to see that is is a morphism of CTBs from $\calG$ to $(\calH\calG,m,\Delta^{\subset,\subset})$.
Another injection $\phi_\calG:\calG\longrightarrow \calH\calG$ is defined by sending any graph $G\in \calG[X]$ to the hypergraph on $X$ defined by
\[\HE(\phi_\calG(G))=E(G)\sqcup\{\{x\}\mid x\in X\}.\]
This is a morphism of CTBs from $\calG$ to $(\calH\calG,m,\Delta^{\subset,\subset})$, $(\calH\calG,m,\Delta^{\subset,\cap})$, $(\calH\calG,m,\Delta^{\cap,\subset})$, and $(\calH\calG,m,\Delta^{\cap,\cap})$.
\end{example}

\begin{example}
Posets are quasi-posets: we have a canonical injection $\phi_{\calP}:\calP\longrightarrow\calQ\calP$. It is not difficult to see that $\phi_\calP$ is a morphism of CTBs from $(\calP,m,\Delta)$ to $(\calQ\calP,m,\Delta)$.
\end{example}

\begin{example}\label{ex2.1.13}
Let$G$ be an oriented graph. We define the quasi-order $\leq_G$ on $V(G)$ as the transitive and reflexive closure of $A(G)$. In other terms, $x\leq_G y$ if $x=y$ or if there exist arcs $(x_0,,x_1)$, $(x_1,x_2),\ldots$, $(x_{k-1},x_k)\in A(G)$, with $x_0=x$ and $x_k=y$. 
This defines a non-injective, surjective species map $\phi_{\calO\calG}$ from $\calO\calG$ to $\calQ\calP$, sending $G$ to $(V(G),\leq_G)$. It is not difficult to see that it is a morphism of CTBs from $(\calO\calG,m,\Delta)$ to $(\calQ\calP,m,\Delta)$.
\end{example}

\begin{example}\label{ex2.1.14}
We define a species morphism from $\calQ\calP$ to $\calP$ as follows. If $X$ is a finite set, and $P\in \calQ\calP[X]$, then
\begin{align*}
\phi_{\calQ\calP}(P)&=\begin{cases}
P\mbox{ if }P\in \calP[X],\\
0\mbox{ otherwise}.
\end{cases}
\end{align*}
Then $\phi_{\calQ\calP}$ is a morphism of CTBs from $(\calQ\calP,m,\Delta)$ to $(\calP,m,\Delta)$. Note that
\[\phi_{\calQ\calP}\circ \phi_{\calP}=\id_{\calP}.\]
 If $G$ is an oriented graph, then
$\phi_{\calQ\calP}\circ \phi_{\calO\calG}(G)=\begin{cases}
(V(G),\leq_G)\mbox{ if $G$ is acyclic},\\
0\mbox{ otherwise}.
\end{cases}$
\end{example}

\begin{example} \label{ex2.1.1.15}
Let $H\in\calH\calG[X]$ be a hypergraph. We define a boolean function 
\[\phi_{\calH\calG,\calB\calF}:\left\{\begin{array}{rcl}
\calS[X]&\longrightarrow&\Z\\
Y&\longmapsto&\calS[Y]\cap \HE(G).
\end{array}\right.\]
This defines an injective morphism $\phi_{\calH\calG,\calB\calF}$ of CTBs from $(\calH\calG,m,\Delta^{\subset,\subset})$ to $\calB\calF$.
\end{example}

\begin{lemma}\label{lem2.1.3}
Let $(\calQ,m,\Delta)$ be a twisted bialgebra, such that $\calQ[\emptyset]$ is 1-dimensional. 
\begin{enumerate}
\item For any finite set $X$, for any $x\in \calQ[X]$,
\[\Delta_{\emptyset,\ldots,\emptyset,X,\emptyset,\ldots,\emptyset}(x)=1_\calQ\otimes \cdots \otimes 1_\calQ\otimes x\otimes 1_\calQ\otimes \cdots \otimes 1_\calQ.\]
\item For any $X_1,\ldots,X_k$ pairwise disjoint finite sets, for any $b_i\in \calQ[X_i]$ with $i\in \calQ[X]$, for any $\sigma\in \sym_k$,
\[\Delta_{X_1,\ldots,X_k}(b_{\sigma(1)}\cdots b_{\sigma(k)})=b_1\otimes\cdots\otimes b_k.\]
\end{enumerate}\end{lemma}

\begin{proof}
1. As $1_\calQ\in \calQ[\emptyset]$, which is one-dimensional, $(1_\calQ)$ is a basis of $\calQ[\emptyset]$. Therefore, $\Delta_{X,\emptyset}(x)=y\otimes 1_\calQ$ for a certain $y\in \calQ[X]$. Moreover,
\[x=(\id_{\calQ[X]}\otimes \varepsilon_\Delta)\circ \Delta_{X,\emptyset}(x)=y\varepsilon_\Delta(1_\calQ)=y.\]
Similarly, $\Delta_{\emptyset,X}(x)=1_\calQ\otimes x$. As $\Delta_{\emptyset,\emptyset}(1_\calQ)=1_\calQ\otimes 1_\calQ$, we obtain the result for $\Delta_{\emptyset,\ldots,\emptyset,X,\emptyset,\ldots,\emptyset}(x)$ by induction on the number of $\emptyset$.

2. From the first item, we obtain
\begin{align*}
\Delta_{X_1,\ldots,X_k}(b_{\sigma(1)}\cdots b_{\sigma(k)})&=
\Delta_{X_1\cap X_{\sigma(1)},\ldots,X_k\cap X_{\sigma(1)}}(b_{\sigma(1)})\cdots\Delta_{X_1\cap X_{\sigma(k)},\ldots,X_k\cap X_{\sigma(1)}}(b_{\sigma(k)})\\
&=\Delta_{\emptyset,\ldots,X_{\sigma(1)},\ldots,\emptyset}(b_{\sigma(1)})\cdots
\Delta_{\emptyset,\ldots,X_{\sigma(k)},\ldots,\emptyset}(b_{\sigma(k)})\\
&=(1_\calQ\otimes \cdots\otimes b_{\sigma(1)}\otimes\cdots\otimes 1_\calQ)\cdots(1_\calQ\otimes \cdots\otimes b_{\sigma(k)}\otimes\cdots\otimes 1_\calQ)\\
&=b_1\otimes\cdots\otimes b_k. \qedhere
\end{align*}
\end{proof}

This lemma will be applied to CTBs $\calB$, as $\ell\calB[\emptyset]=\vect(1_\calB)$. Here is a first example of application:

\begin{prop}\label{prop2.1.4}
Let $(\calB,m,\Delta)$ be a CTB. Let $X$, $Y$ be two disjoint finite sets and $a\in \calB[X]$, $b\in \calB[Y]$. Then $ab\neq 0$. 
\end{prop}

\begin{proof}
We obtain $\Delta_{X,Y}(ab)=a\otimes b\neq 0$, so $ab\neq 0$. 
\end{proof}

Let us formulate the coassociativity and multiplicativity of $\Delta$ with the preceding notations.

\begin{prop}
Let $\calB$ be a CTB.
\begin{enumerate}
\item Let $X$, $Y$ be disjoint finite sets, $a\in \calB[X]$, and $b\in \calB[Y]$. Let $X'$, $Y'$ be disjoint finite sets such that $X\sqcup Y=X'\sqcup Y'$. Then $\Delta_{X',Y'}(ab)\neq 0$ if, and only if, $\Delta_{X\cap X',X\cap Y'}(a)\neq 0$ and $\Delta_{Y\cap X',Y\cap Y'}(b)\neq 0$. If this holds, 
\begin{align*}
(ab)_{\mid_LX'}&=a_{\mid_LX\cap X'}b_{\mid_LY\cap X'},&(ab)_{\mid_RY'}&=a_{\mid_RX\cap Y'}b_{\mid_RY\cap Y'}.
\end{align*}
\item Let $X$, $Y$, $Z$ be pairwise disjoint finite sets, and $b\in \calB[X\sqcup Y\sqcup Z]$. 
The following assertions are equivalent:
\begin{itemize}
\item $\Delta_{X\sqcup Y,Z}(b)\neq 0$ and $\Delta_{X,Y}(b_{\mid_LX\sqcup Y})\neq 0$.
\item $\Delta_{X,Y\sqcup Z}(b)\neq 0$ and $\Delta_{Y,Z}(b_{\mid_RY\sqcup Z})\neq 0$.
\end{itemize}
If these conditions hold,
\begin{align*}
\left(b_{\mid_LX\sqcup Y}\right)_{\mid_LX}&=b_{\mid_L X},&
\left(b_{\mid_LX\sqcup Y}\right)_{\mid_RY}&=\left(b_{\mid_RY\sqcup Z}\right)_{\mid_LY},&
b_{\mid_RZ}&=\left(b_{\mid_RY\sqcup Z}\right)_{\mid_RZ}.
\end{align*}\end{enumerate}\end{prop}

\begin{proof}
1. We have
\begin{align*}
\Delta_{X',Y'}(ab)&=\begin{cases}
(ab)_{\mid_LX'}\otimes (ab)_{\mid_RY'}\mbox{ if }\Delta_{X',Y'}(ab)\neq 0,\\
0\mbox{ otherwise}
\end{cases}\\
=\Delta_{X\cap X',X\cap Y'}(a)\Delta_{Y\cap X',Y\cap Y'}(b)&=\begin{cases}
a_{\mid_LX\cap X'}b_{\mid_LY\cap X'}\otimes a_{\mid_RX\cap Y'}b_{\mid_RY\cap Y'}\\
\hspace{1cm} \mbox{ if }\Delta_{X\cap X',X\cap Y'}(a)\neq 0\mbox{ and }\Delta_{Y\cap X',Y\cap Y'}(b)\neq 0,\\
 0\mbox{ otherwise}.
\end{cases}
\end{align*}
As the product of two elements of $\calB$ is never zero, we obtain the result.

2. We have
\begin{align*}
\left(\Delta_{X,Y}\otimes \id_{\calB[Z]}\right)\circ \Delta_{X\sqcup Y,Z}(b)&=\begin{cases}
\left(b_{\mid_LX\sqcup Y}\right)_{\mid_LX}\otimes\left(b_{\mid_LX\sqcup Y}\right)_{\mid_RY}\otimes b_{\mid_RZ}\\
\hspace{1cm} \mbox{ if } \Delta_{X\sqcup Y,Z}(b)\neq 0\mbox{ and }\Delta_{X,Y}(b_{\mid_LX\sqcup Y})\neq 0,\\
0\mbox{ otherwise}
\end{cases}\\
=\left(\id_{\calB[X]}\otimes \Delta_{Y,Z}\right)\circ\Delta_{X,Y\sqcup Z}(b)&=\begin{cases}
b_{\mid_L X}\otimes \left(b_{\mid_RY\sqcup Z}\right)_{\mid_LY}\otimes \left(b_{\mid_RY\sqcup Z}\right)_{\mid_RZ}\\
\hspace{1cm}\mbox{ if }\Delta_{X,Y\sqcup Z}(b)\neq 0\mbox{ and }\Delta_{Y,Z}(b_{\mid_RY\sqcup Z})\neq 0,\\
0\mbox{ otherwise}.
\end{cases}\end{align*}
The result immediately follows.
\end{proof}

\subsection{Quotients}

\begin{defi}\label{defi2.1.6}
Let $(\calB,m,\Delta)$ be a CTB, and $\calI$ be a subspecies of $\calB$. We shall say that $\calI$ is a biideal of $\calB$ if:
\begin{enumerate}
\item $\calI[\emptyset]=\emptyset$.
\item For any pair $(X$, $Y)$ of disjoint sets, any $a\in \calB[X]$, $b\in \calB[Y]$, if $a\in \calI[X]$ or $b\in \calI[Y]$, then $ab\in \calI[X\sqcup Y]$.
\item For any pair $(X$, $Y)$ of disjoint sets, and $b\in \calI[X\sqcup Y]$ such that $\Delta_{X,Y}(b)=b_{\mid_LX}\otimes b_{\mid_RY}$, with $b_{\mid_LX}\in \calB[X]$ and $b_{\mid_RY}\in \calB[Y]$, then $b_{\mid_LX}\in \calI[X]$ or $b_{\mid_RY}\in \calI[Y]$.
\end{enumerate}
\end{defi}

\begin{prop}
Let $(\calB,m,\Delta)$ be a CTB, and $\calI$ be a biideal of $\calB$. Then $\calB\setminus \calI$ naturally inherits a CTB structure.
\end{prop}

\begin{proof}
By definition, $\ell\calI$ is a biideal of the twisted bialgebra $\ell\calB$. The quotient $\ell\calB/\ell\calI$ is naturally a twisted bialgebra. 
For any finite set $X$, the classes of the elements of $\calB[X]\setminus \calI[X]$ gives a basis of $\ell\calB/\ell\calI[X]$. Identifying these elements with their classes, this gives $\calB\setminus \calI$ a structure of CTB.
\end{proof}

\begin{example}
The subspecies $\calI=\calQ\calP\setminus\calP$ is a biideal of $\calQ\calP$: 
\begin{enumerate}
\item $\calI[\emptyset]=\emptyset$.
\item If $P$ or $Q$ are quasi-posets, and $P$ or $Q$ is not a poset, then $PQ$ is not a poset.
\item Let $P\in \calQ[X\sqcup Y]\setminus \calP[X\sqcup Y]$, such that $\Delta_{X,Y}(P)\neq 0$. Therefore, $X$ is a lower part of $P$. As $P$ is not a poset, there exists $x\neq y\in X\sqcup Y$, with $x\leq_P y$ and $y\leq_P x$.
As $X$ is a lower part of $P$, both $x,y\in X$, or both $x,y\in Y$. Consequently, either $P_{\mid X}$ is not a poset, or $P_{\mid Y}$ is not a poset.
\end{enumerate}
The quotient $\calQ\calP\setminus\calI$ is $\calP$, as a CTB.
\end{example}

\begin{prop}\label{prop2.1.8}
Let $\calB$ and $\calB'$ be CTBs, and $\phi:\calB\longrightarrow \calB'$ be a morphism of CTBs. For any finite set $X$, we put
\[\calI[X]=\{b\in \calB[X]\mid \phi(x)=0\}.\]
Then $\calI$ is a biideal of $\calB$. Moreover, $\phi$ induces a morphism CTBs $\overline{\phi}:\calB\setminus \calI\longrightarrow \calB'$, and for any finite set $X$, any $b\in \calB\setminus \calI[X]$, $\overline{\phi}(b)\neq 0$.
\end{prop}

\begin{proof}
As $\phi(1_\calB)=1_{\calB'}\neq 0$, $1_\calB\notin \calI[\emptyset]$, so $\calI[\emptyset]=\emptyset$. Let $(X$, $Y)$ be pair $(X$, $Y)$ of disjoint sets. Let $a\in \calI[X]$ and $b\in \calB[X]$. Then $\phi(ab)=\phi(a)\phi(b)=0$, so $ab\in \calI[X\sqcup Y]$.
Similarly, if $a\in \calB[X]$ and $b\in \calI[X]$, $ab\in \calI[X\sqcup Y]$.
Let $b\in \calB[X\sqcup Y]$, such that $\Delta_{X,Y}(b)=b_{\mid_LX}\otimes b_{\mid_RY}$, with $b_{\mid_LX}\in \calB[X]$ and $b_{\mid_RY}\in \calB[Y]$. We assume that both $\phi(b_{\mid_LX})$ and $\phi(b_{\mid_RY})$ are nonzero. Then 
\[\Delta_{X,Y}(\phi(b))=(\phi\otimes \phi)\circ \Delta_{X,Y}(b)=\phi(b_{\mid_LX})\otimes \phi(b_{\mid_RY})\neq 0,\]
so $\phi(b)\neq 0$. Dually, if $\phi(b)=0$, then either $\phi(b_{\mid_LX})=0$ or $\phi(b_{\mid_RY})=0$, so $\calI$ is a biideal.

By definition, $\ell\calI$ is included in the kernel of the linearization $\ell\phi$ of $\phi$, so $\ell\phi$ induces a morphism of twisted bialgebras $\ell\overline{\phi}:\ell\calB/\ell\calI\longrightarrow \ell\calB'$. 
This morphism is precisely the linearization of a morphism $\overline{\phi}:\calB\setminus \calI\longrightarrow \calB'$ of CTBs. By construction, if $b\in \calB\setminus \calI[X]$, $\overline{\phi}(b)=\phi(b)\neq 0$. 
\end{proof}

\begin{remark}
The morphism $\overline{\phi}$ is not necessarily injective, as shown in the next example. 
\end{remark}

\begin{example}
We consider the morphism $\phi=\phi_{\calQ\calP}\circ \phi_{\calO\calG}:\calO\calG\longrightarrow \calP$ of Example \ref{ex2.1.14}. For any finite set $X$,
\[\calI[X]=\{G\in \calO\calG[X]\mid\mbox{$G$ contains a cycle}\}.\]
The quotient is the species $\calO\calG_{\mathrm{ac}}$ of acyclic oriented graphs, which is in this way a CTB. The map $\overline{\phi}$ sends any acyclic oriented graph $G$ to the poset $(V(G),\leq_G)$. It is not injective: for example, if $G,H\in \calO\calG_{\mathrm{ac}}[3]$ are defined by
\begin{align*}
A(G)&=\{(1,2),(2,3)\},&A(H)&=\{(1,2),(1,3),(2,3)\},
\end{align*}
then $\phi(G)=\phi(H)$ is the usual order of $[3]$.
\end{example}

\section{The algebraic structure of CTBs}

\subsection{Indecomposable elements}

\begin{defi}\label{defi2.2.1}
Let $(\calB,m,\Delta)$ be a CTB, $X$ be a nonempty finite set, and $b\in \calB[X]$. We shall say that $b$ is indecomposable if, for any disjoint finite sets $X',X''$ such that $X=X'\sqcup X''$, for any $b'\in \calB[X']$ and $b''\in \calB[X'']$, 
\[b=b'b''\Longrightarrow X'=\emptyset \mbox{ or }X''=\emptyset.\] 
\end{defi}

\begin{example}\begin{enumerate}
\item For $\Com$ (Example \ref{ex2.1.1}), the indecomposable elements are the sets of cardinality one.
\item We consider the CTB $\calG$ of graphs (Example \ref{ex2.1.2}). A graph $G\in \calG[X]$ is connected if for any composition $(X',X'')\in \comp[X]$, there exists $x'\in X'$, $x''\in X''$, such that $\{x',x''\}\in E(G)$.
Then $G$ is indecomposable in $(\calG,m,\Delta)$ if, and only if, $G$ is connected. 
Similar results can be proved for oriented graphs, hypergraphs, posets and quasi-posets (Examples \ref{ex2.1.3}--\ref{ex2.1.8}). In the sequel, we shall indifferently use the words indecomposable or connected for these objects.
\end{enumerate}\end{example}

\begin{prop}\label{prop2.2.2}
Let $(\calB,m,\Delta)$ be a CTB. Let $X$ be a finite set and $b\in \calB[X]$. There exists a composition $(X_1,\ldots,X_k)$ of $X$, and, for any $i\in [k]$, $b_i\in \calB[X_i]$, indecomposable, such that $b=b_1\cdots b_k$.
\end{prop}

\begin{proof}
By induction on $|X|$. If $X=\emptyset$, then $b=1_\calB$, and we take, by convention, the empty composition of $\emptyset$. Let us assume the result at all ranks $<|X|$. If $b$ is indecomposable, we take the composition $(X)$ and $b_1=b$.
Otherwise, there exist a pair $(X',X'')$ of disjoint nonempty finite sets, $b'\in \calB[X']$, $b''\in \calB[X'']$, such that $b=b'b''$. 
By the induction hypothesis applied to $b'$ and $b''$, there exist compositions $(X'_1,\ldots,X'_p)\in \comp[X']$ and $(X''_1,\ldots,X''_q)\in \comp[X'']$, and indecomposable elements $b'_i\in \calB[X'_i]$, $b''_j\in \calB[X''_j]$, such that 
\begin{align*}
b'&=b'_1\cdots b'_p,&b''&=b''_1\cdots b''_q.
\end{align*}
Then $b=b'_1\cdots b'_pb''_1\cdots b''_q$: we take the composition $(X'_1,\ldots,X'_p,X''_1,\ldots,X''_q)$.
\end{proof}

\begin{lemma}\label{lem2.2.3}
Let $(\calB,m,\Delta)$ be a CTB. Let $X$ be a nonempty finite set and $b\in \calB[X]$, with two decompositions into indecomposable elements:
\[b=b_1\cdots b_k=b'_1\cdots b'_l.\]
Then $k=l$ and there exists $\sigma\in \sym_k$ such that $b_i=b'_{\sigma(i)}$ for any $i\in [k]$.
\end{lemma}

\begin{proof}
Let us assume that $b_i\in \calB[X_i]$ and $b'_j\in \calB[X'_j]$. Then, by Lemma \ref{lem2.1.3},
\begin{align}
\label{EQ2.1}\Delta_{X_1,\ldots,X_k}(b)&=b_1\otimes \cdots \otimes b_k\neq 0.
\end{align}
On the other hand,
\begin{align}
\label{EQ2.2}\Delta_{X_1,\ldots,X_k}(b)&=\Delta_{X_1\cap X'_1,\ldots,X_k\cap X'_1}(b'_1)\cdots \Delta_{X_1\cap X'_k,\ldots,X_k\cap X'_1}(b'_k).
\end{align}
Therefore, for any $j$, $\Delta_{X_1\cap X'_j,\ldots,X_k\cap X'_j}(b'_j)\neq 0$, and we can write 
\[\Delta_{X_1\cap X'_j,\ldots,X_k\cap X'_j}(b'_j)=b'_{1,j}\otimes \cdots \otimes b'_{k,j},\]
with $b'_{j,i}\in \calB[X_i\cap X'_j]$. We deduce from (\ref{EQ2.1}) and (\ref{EQ2.2}) that for any $i\in [k]$,
\[b_i=b'_{i,1}\cdots b'_{i,l}.\]
As $b_i$ is indecomposable, all the $X_i\cap X'_j$, for $j\in [l]$, are empty, except one, which is equal to $X_i$: there exists a unique $\sigma(i)\in [l]$, such that $X_i\subseteq X'_{\sigma(i)}$. 
By symmetry, $X_i=X'_{\sigma(i)}$, and $\sigma$ is a bijection (which implies that $k=l$). Moreover, we obtain
\[b_i=1_\calB\cdots b'_{\sigma(i)}\cdots 1_\calB=b'_{\sigma(i)}. \qedhere\]
\end{proof}

\begin{theo}\label{theo2.2.4}
Let $(\calB,m,\Delta)$ be a CTB. Let $X$ be a finite set and $b\in \calB[X]$, which we decompose in two ways into indecomposable elements:
\[b=b_1\cdots b_k=b'_1\cdots b'_l.\]
Then $k=l$, and it is possible to go from the first decomposition to the second one by permutations of consecutive commuting elements of $\{b_1,\ldots,b_k\}$.
\end{theo}

\begin{proof}
Lemma \ref{lem2.2.3} gives $k=l$ and the existence of $\sigma \in \sym_k$ such that $b'_i=b_{\sigma(i)}$ for any $i$. Therefore, 
\[b=b_1\cdots b_k=b_{\sigma(1)}\cdots b_{\sigma(k)}.\]
Let us consider a decomposition of $\sigma$ into a product of elementary transpositions $(i,i+1)$ of minimal length:
\[\sigma=t_1\circ \cdots \circ t_p.\]
If $(i,i+1)$ appears in this decomposition, as it is of minimal length, $i+1$ appears before $i$ in $\sigma$, seen as the word $\sigma(1)\cdots\sigma(k)$. Therefore, by Lemma \ref{lem2.1.3},
\begin{align*}
\Delta_{X\setminus (X_i\sqcup X_{i+1}),X_i\sqcup X_{i+1}}(b)&=\Delta_{X\setminus (X_i\sqcup X_{i+1}),X_i\sqcup X_{i+1}}(b_1\ldots b_k)=x\otimes b_ib_{i+1},\\
&=\Delta_{X\setminus (X_i\sqcup X_{i+1}),X_i\sqcup X_{i+1}}(b_{\sigma(i)}\ldots b_{\sigma(k)})=y\otimes b_{i+1}b_i,
\end{align*} 
for certain $x$, $y\in \calB[X\setminus (X_i\sqcup X_{i+1})]$. Therefore, $b_ib_{i+1}=b_{i+1}b_i$. Using the decomposition of $\sigma$, one can go from the first decomposition to the second one by permutations of consecutive commuting elements of $b_1,\ldots,b_k$.
\end{proof}

\begin{prop}
Let $(\calB,m,\Delta)$ be a CTB. For any finite set $X$, and any $b\in \calB[X]$, we denote by $\ic(b)$ the number of factors in any decomposition of $b$ into indecomposable elements. Then:
\begin{enumerate}
\item $\ic(1_\calB)=0$. 
\item If $X\neq \emptyset$, $1\leq \ic(b)\leq |X|$.
\item If $X$, $Y$ are disjoint finite sets, $b\in \calB[X]$, and $b'\in \calB[Y]$, then $\ic(bb')=\ic(b)+\ic(b')$.
\end{enumerate}
\end{prop}

\begin{proof}
Immediate.
\end{proof}

\begin{example}
In $\Com$ (Example \ref{ex2.1.1}), for any finite set $X$, $\ic(X)=|X|$.
For graphs (Example \ref{ex2.1.2}), if $G\in \calG[X]$, then $\ic(G)$ is the number of connected components of $G$. Similar results hold for oriented graphs, hypergraphs, posets, quasi-posets
(Examples \ref{ex2.1.3}--\ref{ex2.1.8}).
\end{example}

\begin{prop}\label{prop2.2.6}
Let $\calB$ and $\calB'$ be two CTBs and $\phi:\calB\longrightarrow\calB'$ be a morphism
of CTBs.
\begin{enumerate}
\item For any finite set $X$, and any $b\in \calB[X]$ such that $\phi(b)\neq 0$, $\ic(\phi(b))\geq \ic(b)$.
\item Let us assume that $\phi$ is injective. For any finite set $X$, and any $b\in \calB[X]$, $\ic(\phi(b))=\ic(b)$.
\end{enumerate}\end{prop}

\begin{proof}
1. Let us decompose $b$ into indecomposable elements:
\[b=b_1\cdots b_k,\]
with $k=\ic(b)$. Then $\phi(b)=\phi(b_1)\cdots \phi(b_k)$, so $\ic(\phi(b))\geq k=\ic(b)$.

2. Let us decompose $\phi(b)$ into indecomposable elements:
\[\phi(b)=b'_1\cdots b'_l,\]
with $l=\ic(\phi(b))$, and $b'_i\in \calB'[X_i]$ for $i\in [k]$. By Lemma \ref{lem2.1.3},
\[\phi^{\otimes l}\circ \Delta_{X_1,\ldots,X_l}(b)=\Delta_{X_1,\ldots,X_l}(\phi(b))=b'_1\otimes \cdots\otimes b'_l\neq 0,\]
so $\Delta_{X_1,\ldots,X_l}(b)\neq 0$. We then put $\Delta_{X_1,\ldots,X_l}(b)=b_1\otimes\cdots \otimes b_l$, and we obtain that $\phi(b_i)=b'_i$ for any $i\in [k]$. This gives
\[\phi(b)=b'_1\cdots b'_l=\phi(b_1)\cdots \phi(b_l)=\phi(b_1\cdots b_l).\]
As $\phi$ is injective, $b=b_1\cdots b_l$, so $\ic(b)\geq l=\ic(\phi(b))$. Combined with the first item, $\ic(b)=\ic(\phi(b))$.
\end{proof}

\begin{remark}\begin{enumerate}
\item Non-injective morphisms of CTBs do not necessarily conserve $\ic$. Take the morphism $\phi_\Com$ (Example \ref{ex2.1.10}). If $G$ is a graph, then $\ic(G)$ is the number of connected components of $G$, whereas $\ic(\phi_\Com(G))$ is the number of vertices of $G$. 
\item A morphism of CTBs may conserve $\ic$ and not be injective: this is the case of the morphism $\phi_{\calO\calG}$ (Example \ref{ex2.1.13}).
\end{enumerate}

\end{remark}

\subsection{Discrete elements}

\begin{defi}\label{defi2.2.7}
Let $\calB$ be a CTB. For any finite set $X$, we put
\[\bfD\calB[X]=\{b\in \calB[X]\mid \ic(b)=|X|\}.\]
This defines a subspecies of $\calB$.
\end{defi}

\begin{remark}
If $X=\emptyset$ or if $|X|=1$, then $\bfD\calB[X]=\calB[X]$.
\end{remark}

\begin{example}\begin{enumerate}
\item $\bfD\Com=\Com$.
\item For graphs (Example \ref{ex2.1.2}), for any finite set $X$, there exists a unique element in $\bfD\calG[X]$, which is the graph $(X,\emptyset)$.
\item For oriented graphs (Example \ref{ex2.1.3}), for any finite set $X$, there exists a unique element in $\bfD\calO\calG[X]$, which is the oriented graph $(X,\emptyset)$.
\item For hypergraphs (Example \ref{ex2.1.4}), for any finite set $X$, there exist $2^{|X|}$ elements in $\bfD\calH\calG[X]$. They are the hypergraphs $G\in \calH\calG[X]$ such that $\HE(G)\subseteq \calS_1[X]$.
\item For posets (Example \ref{ex2.1.5}), for any finite set $X$, there exists a unique element in $\bfD\calP[X]$, which is poset $P\in \calG[X]$ such that 
\begin{align*}
&\forall x,y\in X,&x\leq_P y&\Longleftrightarrow x=y.
\end{align*}
This poset is also the unique element of $\bfD\calQ\calP[X]$ (Example \ref{ex2.1.8}).
\end{enumerate}\end{example}

\begin{prop}\label{prop2.2.8}
Let $(\calB,m,\Delta)$ be a CTB. Let $X$, $Y$ be disjoint finite sets.
\begin{enumerate}
\item If $a\in \calB[X]$ and $b\in \calB[Y]$, then $ab\in \bfD\calB[X\sqcup Y]$, if, and only if, $a\in \bfD\calB[X]$ and $b\in \bfD\calB[Y]$.
\item If $b\in \bfD\calB[X\sqcup Y]$, then $\Delta_{X,Y}(b)\neq 0$ and belongs to $\bfD\calB[X]\otimes \bfD\calB[Y]$.
\end{enumerate}
\end{prop}

\begin{proof}
1. Note that $\ic(a)\leq |X|$, $\ic(b)\leq |Y|$, and $\ic(ab)=\ic(a)+\ic(b)$. Therefore, $\ic(ab)=|X|+|Y|$ if, and only if, $\ic(a)=|X|$ and $\ic(b)=|Y|$. This implies the result.

2. As $\ic(b)=|X|$, there exists an indexation $x_1,\ldots,x_k$ of the elements of $X\sqcup Y$, and $b_i\in \calB[\{x_i\}]$ for $i\in [k]$, such that 
\[b=b_1\cdots b_k.\]
Lemma \ref{lem2.1.3} gives
\[\Delta_{X,Y}(b)=\underbrace{\prod_{b_i\in X}b_i}_{=b_{\mid_LX}}\otimes\underbrace{\prod_{b_i\in Y}b_i}_{=b_{\mid_RY}}.\] 
This belongs to $\bfD\calB[X]\otimes \bfD\calB[Y]$, as $\ic \left(b_{\mid_LX}\right)=|X|$ and $\ic\left(b_{\mid_RY}\right)=|Y|$.
\end{proof}

\begin{remark}
Proposition \ref{prop2.2.8} implies that $\bfD\calB$ is a sub-CTDB of $\calB$.
\end{remark}

\chapter{The combinatorial twisted bialgebra $\exch$}

\section*{Introduction}

We introduce in this chapter the CTB $\exch$ and its special character $\epsilon$, and explore its universal property. We start with the species of set compositions $\comp$, defined in Example \ref{ex1.1.1}.
It is a twisted algebra, with the concatenation product: if $X$ and $Y$ are disjoint finite sets, $\pi=(X_1,\ldots,X_k)\in \comp[X]$ and $\pi'=(Y_1,\ldots,Y_l)\in \comp[Y]$, then 
\[(\pi,\pi')=(X_1,\ldots,X_k,Y_1,\ldots,Y_l).\]
A nonempty subset $I$ of $\comp[X]$ satisfies the exchange condition if for any $(X_1,X_2)\in \comp[X]$, and any $\pi_i,\pi'_i\in \comp[X_i]$ for $i\in [2]$,
\[(\pi_1,\pi_2)\mbox{ and }(\pi'_1,\pi'_2)\in I\Longrightarrow (\pi_1',\pi_2)\in I,\]
see Definition \ref{defi3.1.1}. These subsets of $\comp[X]$ form a subset $\exch[X]$ of $\calS[\comp[X]]$, and this defines a species $\exch$. 
We give it a structure of CTB as follows:
\begin{itemize}
\item If $X'\subseteq X''$ are two finite sets and $\pi=(X_1,\ldots,X_k)\in \comp[X'']$, then $\pi_{\mid X'}\in \comp[X']$ is obtained by deleting the $\emptyset$'s in the sequence $(X_1\cap X',\ldots,X_k\cap X'')$.
If $X$, $Y$ are two disjoint finite sets, $I\subseteq \comp[X]$, and $J\subseteq \comp[Y]$, then
\[I\squplus J=\{\pi \in \comp[X\sqcup Y] \mid \pi_{\mid X}\in I, \pi_{\mid Y}\in J\},\]
see Definition \ref{defi3.1.4}. For example, if $(X_1,X_2)\in \comp[X]$,
\begin{align*}
\{(X_1,X_2),\: (X)\}\squplus \{(Y)\}&=\left\{\begin{array}{c}
(X_1,X_2,Y),\:(X_1,Y,X_2),\:(Y,X_1,X_2),\\
(X_1\sqcup Y,X_2),\:(X_1,X_2\sqcup Y),\:(X$, $Y),\:(Y,X),\:(X\sqcup Y)
\end{array}\right\}.
\end{align*}
This product can be also defined with the help of quasi-shuffles (Proposition \ref{prop3.1.8}), justifying the notation $\squplus$. 
\item If $(X,Y)$ is a pair of disjoint finite sets, and $I\in \exch[X\sqcup Y]$, we put
\begin{align*}
I_{\mid_L X}&=\{\pi\in \comp[X]\mid \exists \pi'\in \comp[Y],\: (\pi,\pi')\in I\},\\
I_{\mid_R Y}&=\{\pi'\in \comp[Y]\mid \exists \pi\in \comp[X],\: (\pi,\pi')\in I\},
\end{align*}
and
\[\Delta_{X,Y}(I)=\begin{cases}
I_{\mid_LX}\otimes I_{\mid R_Y}\mbox{ if $I_{\mid_LX}\neq \emptyset$ and $I_{\mid_RY}\neq \emptyset$},\\
0\mbox{ otherwise}.
\end{cases}\]
\end{itemize}
The exchange condition insures that $(\exch,\squplus,\Delta)$ is a CTB (Theorem \ref{theo3.1.16}).

The second section is devoted to the terminal property of $\exch$. If $\calB$ is a CTB, a combinatorial character of $\calB$ is a character of the twisted algebra $\ell\calB$, taking only the values $0$ or $1$ on $\calB$ (Definition \ref{defi3.2.1}).
If $f$ is a combinatorial character of $\calB$, then the elements $b$ of $\calB$ such that $f(b)=1$ define a twisted subalgebra of $\ell\calB$, which we call special (Definition \ref{defi3.2.2}).
This defines a one-to-one correspondence between combinatorial characters and special twisted subalgebras of $\calB$ (Proposition \ref{prop3.2.5}). 
In particular, combinatorial twisted subbialgebras of $\calB$ are special twisted subalgebras, but the converse is generally false.
Two particular special twisted subalgebras immediately appear: $\calB$ itself and $\bfD\calB$, corresponding to the combinatorial characters $\epsilon^\calB_+$ and $\epsilon^\calB_-$, given as follows: for any finite set $X$,
\begin{align*}
&\forall b\in \calB[X],&\epsilon^\calB_+(b)&=1,&\epsilon^\calB_-(b)&=\begin{cases}
1\mbox{ if }\ic(b)=|X|,\\
0\mbox{ otherwise}.
\end{cases}\end{align*}
In the particular case of $\exch$, we simply denote $\epsilon_\pm^\exch=\epsilon_\pm$. We introduce a third character $\epsilon$ of $\exch$, defined by
\begin{align*}
&\forall I\in \exch[X],&\epsilon(I)&=\begin{cases}
1\mbox{ if }(X)\in I,\\
0\mbox{ otherwise},
\end{cases}
\end{align*}
see Proposition \ref{prop3.2.7}. The universal property of $(\exch,\epsilon)$ is stated in Theorem \ref{theo3.2.9}: if $\calB$ is a CTB and $f$ is a combinatorial character of $\calB$,
there exists a unique morphism $\Gamma_f:\calB\longrightarrow \exch$ such that $\epsilon\circ \Gamma_f=f$. This is a combinatorial version of Aguiar, Bergeron and Sottile's result on the Hopf algebra of quasi-symmetric functions \cite{Aguiar2006-2}.
Moreover, for any finite set $X$, and any $b\in \calB[X]$, if $\Gamma_f(b)\neq 0$, then it is equal to the set of compositions $(X_1,\ldots,X_k)\in \comp[X]$ such that
\[f^{\otimes k}\circ \Delta_{X_1,\ldots,X_k}(b)=1.\]

To push forward this result, we shall need two sub-CTBs of $\exch$. For any finite set $X$, $\comp[X]$ is naturally ordered by the refinement. For example, if $(X_1,X_2,X_3)\in \comp[X]$,
\[(X_1\sqcup X_2\sqcup X_3)\leq (X_1\sqcup X_2,X_3),\: (X_1,X_2\sqcup X_3)\leq (X_1,X_2,X_3).\]
A subset $I$ of $\comp[X]$ is an upper part of $\comp[X]$ if
\begin{align*}
&\forall \pi,\pi'\in \comp[X],&(\pi\leq \pi'\mbox{ and }\pi\in I)&\Longrightarrow \pi'\in I.
\end{align*}
It is a lower part of $\comp[X]$ if
\begin{align*}
&\forall \pi,\pi'\in \comp[X],&(\pi\leq \pi'\mbox{ and }\pi'\in I)&\Longrightarrow \pi\in I.
\end{align*}
The subspecies $\exch_+$ and $\exch_-$ are defined as, respectively, the elements of $\exch$ which are upper parts, respectively lower parts (Definition \ref{defi3.1.1}). 
It turns out that $\exch_+$ and $\exch_-$ are sub-CTBs of $\exch$ (Theorem \ref{theo3.1.16}). We prove in Proposition \ref{prop3.2.7} that
\begin{align*}
\epsilon_+^{\exch_+}&={\epsilon_+}_{\mid \exch_+}=\epsilon_{\mid \exch_+},&\epsilon_-^{\exch_-}&={\epsilon_-}_{\mid \exch_-}=\epsilon_{\mid \exch_-}.
\end{align*}
Moreover, for any CTB $\calB$, $\Gamma_{\epsilon_\pm^\calB}$, which we simply denote by $\Gamma_\pm^\calB$, takes in values in $\exch_\pm$ (Corollary \ref{cor3.2.11} for $-$ and Theorem \ref{theo3.2.13} for $+$).
We can apply this for $\calB=\exch$: we simply write $\Gamma_\pm^{\exch}=\Gamma_\pm$, and this gives projections on $\exch_\pm$. If $I\in \exch[X]$,
\begin{itemize}
\item $\Gamma_-(I)=m(I)=\{\pi\in \comp[X]\mid \exists \pi'\in I,\: \pi\leq \pi'\}$ (Theorem \ref{theo3.2.15}).
\item $\Gamma_+(I)=\{\pi\in m(I)\mid \forall \pi'\in \comp[X],\: \pi\leq \pi'\Longrightarrow \pi'\in m(I)\}$ if this set is not empty, and $0$ otherwise (Theorems \ref{theo3.2.13} and \ref{theo3.2.16}).
\end{itemize}
This leads to the important notion of non-degenerate CTB (Definition \ref{defi3.2.14}): $\calB$ is non-degenerate if $\Gamma_+^\calB$ does not take the value $0$ on $\calB$.
For example, graphs, posets, hypergraphs form non-degenerate CTBs, but not quasi-posets.
We show in Theorem \ref{theo3.2.17} how to reduce ourselves to the non-degenerate case, without changing $\Gamma_+^\calB$, with the help of a combinatorial quotient $\nd\calB$ of $\calB$.

Another interesting subspecies $\exch_{\pp}$ of $\exch$ is given by subsets of $\comp[X]$, satisfying the exchange condition, and whose all elements are of maximal length: this is stable under $\Delta$, but not under $\squplus$.
However, it can be given another product $\shuffle$, making it a CTB (Theorem \ref{theo3.4.1}). We show (Proposition \ref{prop3.4.3}) that both morphisms $\Gamma^{\exch_{\pp}}_\pm:\exch_{\pp}\longrightarrow \exch_\pm$ are injective.
This CTB $\exch_{\pp}$ is used to construct a polynomial invariant $\phi_{\pp}^\calB$ for any CTB $\calB$ (Proposition \ref{prop3.5.1}), which turns out to be homogeneous, and the associated character $f_{\pp}^\calB$, taking its values in $[0;1]\cap \Q$.
More generally, we use the results of the first chapter, especially Theorem \ref{theo1.2.23}, on a CTB $\calB$, in order to obtain particular polynomial invariants.
Firstly, the polynomial invariant of $\exch$ associated to $\epsilon$ is denoted by $\Phi$. From Proposition \ref{prop3.4.3}, 
\begin{align*}
&\forall I\in \exch[X],&\Phi(I)&=\sum_{k=0}^{|X|} |I\cap \comp_k[X]|\binom{T}{k}.
\end{align*}
 The combinatorial characters $\epsilon_\pm^\calB$ induce two polynomial invariants $\phi_\pm^\calB$ (Proposition \ref{prop3.5.13}), satisfying 
\begin{align}
\label{EQ3.1}\phi^\calB_+&=\Phi\circ \Gamma_+^\calB,&\Phi_-&=\Phi\circ \Gamma_-^\calB.
\end{align}
These polynomial invariants admit a combinatorial interpretation (Proposition \ref{prop3.5.6}). In particular:
\begin{itemize}
\item For graphs, $\phi_+^\calG$ is the chromatic polynomial \cite{Harary1969}. For hypergraphs, with their different CTB structures, it gives one of the various chromatic polynomials of \cite{Ebrahimi-Fard2022,Aval2019,Borowiecki2000}. The invariants $\phi_-^\calG$ and $\phi_-^{\calH\calG}$ are far less interesting, as they only depend on the number of vertices.
\item For posets, $\phi_+^\calP$ is the strict Ehrhart polynomial, whereas $\phi_-^\calP$ is the (non strict) Ehrhart polynomial.
\end{itemize}
We then detail several properties of these polynomial invariants, including results on their degrees and their leading coefficient (given by the character $f_{\pp}^\calB$, Proposition \ref{prop3.5.8}).
A consequence of (\ref{EQ3.1}) is that the set of polynomials $\phi_\pm^\calB(b)$'s, associated to objects $b$ of a fixed degree $n$ of any CTB $\calB$, is finite. We give these sets for $n\leq 3$ in Remark \ref{rk3.5.3}. 
They include chromatic polynomials of graphs and hypergraphs, Ehrhart polynomials, and many more. We also detail a result on the possible leading coefficients $f_{\pp}^\calB$ of these polynomials in Proposition \ref{prop3.5.17}.

An important point is that the twisted bialgebras $\ell\exch$ is not isomorphic to $\ell\comp$. Both satisfies a terminal property, but not in the same category: the category of CTB for the first one, the category of twisted bialgebras $\calB$ such that $\calB[\emptyset]$ is one-dimensional for the second one.
However, as CTBs are also twisted bialgebras, there exists a natural morphism of twisted bialgebras $\Theta$ from $\exch$ to $\ell\comp$, which sends $I\in \exch[X]$ to the sum of its elements in $\comp[X]$ (Proposition \ref{prop3.3.1}).
This morphism relates the two terminal properties of $\exch$ and $\comp$, as shown in Corollary \ref{cor3.3.2}. More precisely, if $\calB$ is a CTB and $f$ is a combinatorial character of $\calB$, then:
\begin{itemize}
\item There exists a unique morphism of CTBs$\Phi_f:\calB\longrightarrow\exch$, such that $\epsilon\circ\Phi_f=f$.
\item There exists a unique morphism of twisted bialgebras $\phi_f:\ell\calB\longrightarrow\ell\comp$ such that $\epsilon\circ \phi_f=f$.
\end{itemize}
Then $\phi_f=\Theta\circ\Phi_f$.
The morphism allows to see $\exch$ as a lifting of the quasi-shuffle product of compositions to the level of subsets of compositions. In fact, $\exch$ is maximal for such a construction, as it is shown in Proposition \ref{prop3.3.3}.
Applying the full Fock functor $\calF_F$, we obtain a morphism from $\calF_F[\exch]$ to the Hopf algebra of packed words $\WQSym$ \cite{Novelli2006-2}, sending the class of an element of $\exch[n]$ to the sum of its elements, seen as packed words rather than compositions.
A similar result is also proved for $\exch_{\pp}$, with $\WQSym$ replaced by $\FQSym$, the Malvenuto-Reutenauer Hopf algebra of permutations \cite{Malvenuto1995,Duchamp2002}, see Proposition \ref{propMalvenutoReutenauer}.

\section{The exchange condition}

\subsection{Definition}

\begin{notation}\begin{enumerate}
\item For any finite set $X$, $\comp[X]$ is partially ordered by the refinement $\leq$: for $\pi,\pi'\in \comp[X]$, $\pi\leq \pi'$ if each block of $\pi$ is the union of consecutive blocks of $\pi'$.
The minimal elements of $\comp[X]$ is $(X)$, its maximal elements are the $|X|!$ elements of $\comp_{|X|}[X]$.
\item If $X$ and $Y$ are two disjoint finite sets, $\pi=(X_1,\ldots,X_k)\in \comp[X]$, and $\pi'=(Y_1,\ldots,Y_l)\in \comp[Y]$, we put
\[(\pi,\pi')=(X_1,\ldots,X_k,Y_1,\ldots,Y_l)\in \comp[X\sqcup Y].\]
\item Let $(B,\Delta)$ be a twisted coalgebra. For any finite set $X$, and any $\pi=(X_1,\ldots,X_k)\in \comp[X]$, we put
\[\Delta_\pi=\Delta_{X_1,\ldots,X_k}.\]
Note that if $X$, $Y$ are two disjoint finite sets, $\pi\in \comp[X]$ and $\pi'\in \comp[Y]$, then
\[\Delta_{(\pi,\pi')}=(\Delta_\pi\otimes\Delta_{\pi'})\circ\Delta_{X,Y}.\]
\end{enumerate}\end{notation}

\begin{example}
Here is the Hasse graph of $\comp[2]$. 
\begin{align*}
\xymatrix{(1,\:2)\ar@{-}[rd]&&(2,\:1)\ar@{-}[ld]\\
&(12)&}\\
\end{align*}
Here is the Hasse graph of $\comp[3]$.
\begin{align*}
\xymatrix{(1,3,2)&(1,2,3)&(2,1,3)&(2,3,1)&(3,2,1)&(3,1,2)\\
(1,23)\ar@{-}[u]\ar@{-}[ru]&(12,3)\ar@{-}[u]\ar@{-}[ru]&(2,13)\ar@{-}[u]\ar@{-}[ru]&(23,1)\ar@{-}[u]\ar@{-}[ru]&
(3,12)\ar@{-}[u]\ar@{-}[ru]&(13,2)\ar@{-}[u]\ar@{-}[lllllu]&\\
&&&(123)\ar@{-}[lllu]\ar@{-}[llu]\ar@{-}[lu]\ar@{-}[u]\ar@{-}[ru]\ar@{-}[rru]&&}
\end{align*}\end{example}

\begin{defi}\label{defi3.1.1}\begin{enumerate}
\item Let $X$ be a finite set and $I\subseteq \comp[X]$, nonempty. We shall say that $I$ satisfies the exchange condition if for any $(X_1,X_2)\in \comp_2[X]$, for any $\pi_i,\pi'_i\in \comp[X_i]$ for $i\in [2]$,
\[(\pi_1,\pi_2)\mbox{ and }(\pi'_1,\pi'_2)\in I\Longrightarrow (\pi_1',\pi_2)\in I.\]
The set of nonempty subsets of $\comp[X]$ satisfying the exchange condition is denoted by $\exch[X]$. This defines a species $\exch$.
\item Let $X$ be a finite set and $I\in \exch[X]$. We shall say that $I\in \exch_+[X]$ if
\begin{align*}
&\forall \pi,\:\pi'\in \comp[X],&\pi\in I\mbox{ and }\pi\leq \pi'&\Longrightarrow \pi'\in I.
\end{align*}
This defines a subspecies $\exch_+$ of $\exch$.
\item Let $X$ be a finite set and $I\in \exch[X]$. We shall say that $I\in \exch_-[X]$ if \begin{align*}
&\forall \pi,\:\pi'\in \comp[X],&\pi'\in I\mbox{ and }\pi\leq \pi'&\Longrightarrow \pi'\in I.
\end{align*}
This defines a subspecies $\exch_-$ of $\exch$.
\item Let $X$ be a finite set and $I\in \exch[X]$. We shall say that $I\in \exch_{\pp}[X]$ if $I\subseteq \comp_{|X|}[X]$. 
This defines a subspecies $\exch_{\pp}$ of $\exch$.
\end{enumerate}\end{defi}

\begin{remark}
If $|X|\leq 3$, then any nonempty $I\subseteq \comp[X]$ belongs to $\exch[X]$. This allows to obtain these cardinalities:
\[\begin{array}{|c||c|c|c|}
\hline n&1&2&3\\
\hline\hline |\exch[n]|&1&7&8191\\
\hline |\exch_+[n]|&1&4&322\\
\hline |\exch_-[n]|&1&4&322\\
\hline |\exch_{\pp}[n]|&1&3&63\\
\hline\end{array}\]
This does not work anymore if $|X|\geq 4$, as shown by the following example:
\[I=\{(12,34),(1,2,3,4)\}\subset \comp[4].\]
As $(1,2,34)\notin I$, $I\notin \exch[4]$.
It is noticeable that if $n\leq 3$, $|\exch_+[n]|=|\exch_-[n]|$. This is due to the fact that the complement of an element of $\exch_\pm[n]$ belongs to $\exch_\mp[n]$ or is empty. This is no more true if $n\geq 4$, as the complement of an element of $\exch[n]$ does not necessarily satisfies the exchange condition. For example, take
\[I=\{(1,2,3,4)\}\in \exch_+[4].\]
Its complement contains $(2,1,3,4)$ and $(1,2,4,3)$, but not $(1,2,3,4)$,
so does not belong to $\exch[4]$. 
\end{remark}

\begin{prop}
Let $X$ be a finite set, and $I\in \exch[X]$. We put
\[m(I)=\{\pi \in \comp[X]\mid \exists \pi'\in I,\:\pi\leq \pi'\}.\]
Then $m(I)\in \exch_-[X]$. 
\end{prop}

\begin{proof} 
Note that $I\subseteq m(I)$, so $m(I)\neq \emptyset$.
Let $\pi'\in m(I)$ and $\pi\in \comp[X]$ such that $\pi\leq \pi'$. By definition of $m(I)$, there exists $\pi''\in I$, such that $\pi'\leq \pi''$. Therefore, $\pi\leq \pi'\leq \pi''$, so $\pi \in m(I)$. It remains to show that $m(I)$ satisfies the exchange condition. 
Let $(X_1,X_2)\in \comp[X]$, $\pi_1,\pi'_1\in \comp[X_1]$, $\pi_2,\pi'_2\in \comp[X_2]$, such that $(\pi_1,\pi_2)$ and $(\pi_1',\pi_2')\in m(I)$. By definition of $m(I)$, there exists $\varpi,\varpi'\in I$, such that $(\pi_1,\pi_2)\leq \varpi$ and $(\pi_1',\pi_2')\leq \varpi'$.
By definition of the order on compositions, $\varpi=(\varpi_1,\varpi_2)$ and $\varpi'=(\varpi'_1,\varpi'_2)$, with $\varpi_1,\varpi'_1\in \comp[X_1]$ and $\varpi_2,\varpi'_2\in \comp[X_2]$. Moreover, $\pi_i\leq \varpi_i$ and $\pi'_i\leq \varpi'_i$ for $i\in [2]$.
By the exchange condition for $I$, $(\varpi'_1,\varpi_2)\in I$. As $(\pi'_1,\pi_2)\leq (\varpi'_1,\varpi_2)$, $(\pi'_1,\pi_2)\in m(I)$.
\end{proof}

\begin{prop}\label{prop3.1.3}
Let $X$ be a finite set, and $I\in \exch[X]$. We put
\[M(I)=I\cap \comp_{|X|}[X].\]
Then $M(I)=\emptyset$ or $M(I)\in \exch_{\pp}[X]$.
In particular, if $I\in \exch_+[X]$, then $M(I)\in \exch_{\pp}[X]$ and $m(I)=m(M(I))$.
\end{prop}

\begin{proof}
Let us assume that $M(I)\neq \emptyset$. Let $(X_1,X_2)\in \comp[X]$, $\pi_i,\pi'_i\in \comp[X_i]$ for $i\in [2]$, such that $(\pi_1,\pi_2)\in M(I)$ and $(\pi'_1,\pi'_2)\in M(I)$. Then 
\[\ell((\pi_1,\pi_2))=\ell(\pi_1)+\ell(\pi_2)=|X|=|X_1|+|X_2|.\]
As $\ell(\pi_i)\leq |X_i|$ for $i\in [2]$, necessarily $\ell(\pi_i)=|X_i|$ for $i\in [2]$, and similarly $\ell(\pi'_i)=|X_i|$ for $i\in [2]$. By the exchange condition for $I$, $(\pi'_1,\pi_2)\in I$, and 
\[\ell((\pi'_1,\pi_2))=\ell(\pi'_1)+\ell(\pi_2)=|X_1|+|X_2|=|X|,\]
so $(\pi'_1,\pi_2)\in M(I)$. Therefore, $M(I)\in \exch_{\pp}[X]$.

Let us assume that $I\in \exch_+[X]$. As $I$ is nonempty, let us take $\pi\in I$. Let $\pi'\in \comp[X]$,
maximal for $\leq$, such that $\pi\leq \pi'$. Then $\pi'\in \comp_{|X|}[X]$ and $\pi'\in I$,
so $\pi'\in M(I)$, which is consequently nonempty. As $M(I)\subseteq I$, $m(M(I))\subseteq m(I)$.
Let $\pi\in m(I)$. There exists $\pi'\in I$, such that $\pi\leq \pi'$. Let us consider a maximal element $\pi''$ of $\comp[X]$, such that $\pi'\leq \pi''$. Then $\pi''\in \comp_{|X|}[X]$ As $I\in \exch_+[X]$, $\pi''\in I$, so $\pi''\in M(I)$. Therefore, $\pi\leq\pi''$, so $\pi\in m(M(I))$. 
\end{proof}

\subsection{Product on $\exch$}

\begin{defi}\label{defi3.1.4}
\begin{enumerate}
\item Let $X$ be a finite set, $Y\subseteq X$ and $\pi=(X_1,\ldots,X_k)\in \comp[X]$. We put
\[\pi_{\mid Y}=\elim(X_1\cap Y,\ldots,X_k\cap Y)\in \comp[Y],\]
where $\elim$ means that one deletes the empty $X_i\cap Y$'s in the sequence $(X_1\cap Y,\ldots,X_k\cap Y)$.
\item Let $X$, $Y$ be two disjoint finite sets, and let $I\subseteq \comp[X]$, and $J\subseteq \comp[Y]$. We put
\[I\squplus J=\{\pi \in \comp[X\sqcup Y] \mid \pi_{\mid X}\in I, \pi_{\mid Y}\in J\}.\]
\end{enumerate}
\end{defi}

\begin{example}
Let $X$, $Y$ be disjoint finite sets. 
\begin{align*}
\{(X)\}\squplus\{(Y)\}&=\{(X,Y),(Y,X),(X\sqcup Y)\}.
\comp[X]\squplus\comp[Y]&=\comp[X\sqcup Y].
\end{align*}
If $(X_1,X_2)\in \comp[X]$,
\begin{align*}
\{(X_1,X_2),\: (X)\}\squplus \{(Y)\}&=\left\{\begin{array}{c}
(X_1,X_2,Y),\:(X_1,Y,X_2),\:(Y,X_1,X_2),\\
(X_1\sqcup Y,X_2),\:(X_1,X_2\sqcup Y),\:(X$, $Y),\:(Y,X),\:(X\sqcup Y)
\end{array}\right\}.
\end{align*}
\end{example}

\begin{remark}
Obviously, if $Z\subseteq Y\subseteq X$ are finite sets, and $\pi\in \comp[X]$, then
\[(\pi_{\mid Y})_{\mid Z}=\pi_{\mid Z}.\]
If $Y\subseteq X$ are finite sets, and $\pi\leq \pi'\in \comp[X]$, then
\[\pi_{\mid Y}\leq \pi'_{\mid Y}.\]
If $X$, $Y$ are disjoint finite sets, for any $I\subseteq \comp[X]$, and any $J\subseteq \comp[Y]$,
\[\emptyset\squplus J=I\squplus\emptyset=\emptyset.\]
Moreover, by definition of $\squplus$,
\[I\squplus J=J\squplus I.\]
\end{remark}

\begin{lemma}
Let $X$, $Y,Z$ be pairwise disjoint finite sets, $I\subseteq \comp[X]$, $J\subseteq \comp[Y]$ and $K\subseteq \comp[Z]$. Then $(I\squplus J)\squplus K=I\squplus (J\squplus K)$. 
Moreover, $\{\emptyset\}\subseteq \comp[\emptyset]$ is a unit for the product $\squplus$.
\end{lemma}

\begin{proof}
Let $\pi\in \comp[X\sqcup Y\sqcup Z]$.
\begin{align*}
\pi \in (I\squplus J)\squplus K&\Longleftrightarrow \pi_{\mid X\sqcup Y}\in I\squplus J,\: \pi_{\mid Z}\in K\\
&\Longleftrightarrow (\pi_{\mid X\sqcup Y})_{\mid X}\in I,\: (\pi_{\mid X\sqcup Y})_{\mid Y}\in J,\: \pi_{\mid Z}\in K\\
&\Longleftrightarrow \pi_{\mid X}\in I,\: \pi_{\mid Y}\in J,\: \pi_{\mid Z}\in K\\
&\Longleftrightarrow \pi_{\mid X}\in I,\: (\pi_{\mid Y\sqcup Z})_{\mid Y}\in J,\: (\pi_{\mid Y\sqcup Z})_{\mid Z}\in K\\
&\Longleftrightarrow \pi_{\mid X}\in I,\: (\pi_{\mid Y\sqcup Z})\in J\squplus K\\
&\Longleftrightarrow \pi \in I\squplus(J\squplus K).
\end{align*}
So $(I\squplus J)\squplus K=I\squplus (J\squplus K)$. 

Let $X$ be a finite set, and $I\subseteq \comp[X]$. Let $\pi\in \comp[X]$.
\begin{align*}
\pi \in \{\emptyset\}\squplus I&\Longleftrightarrow \pi_{\mid \emptyset}=\emptyset\mbox{ 	and }\pi_{\mid X}\in I\Longleftrightarrow \emptyset=\emptyset\mbox{ and }\pi\in I\Longleftrightarrow \pi\in I.
\end{align*}
So $ \{\emptyset\}\squplus I=I$. Similarly, $I\squplus \{\emptyset\}=I$.
\end{proof}

\begin{lemma}
Let $X$, $Y$ be disjoint finite sets, $I\subseteq \comp[X]$, and $J\subseteq \comp[Y]$.
\begin{enumerate}
\item If $I\in \exch[X]$ and $J\in \exch[Y]$, then $I\squplus J\in \exch[X\sqcup Y]$.
\item If $I\in \exch_-[X]$ and $J\in \exch_-[Y]$, then $I\squplus J\in \exch_-[X\sqcup Y]$.
\item If $I\in \exch_+[X]$ and $J\in \exch_+[Y]$, then $I\squplus J\in \exch_+[X\sqcup Y]$.
\end{enumerate}\end{lemma}

\begin{proof}
1. Firstly, observe that $I$ and $J$ are nonempty. Let us take $\pi\in I$ and $\pi'\in J$, then $(\pi,\pi')\in I\squplus J$, so $I\squplus J$ is nonempty.
Let $(X_1,X_2)\in \comp[X\sqcup Y]$, $\pi_i,\pi'_i\in \comp[X_i]$ for $i\in [2]$, such that $(\pi_1,\pi_2)\in I\squplus J$, $(\pi'_1,\pi'_2)\in I\squplus J$. Then
\begin{align*}
(\pi_1,\pi_2)_{\mid X}=({\pi_1}_{\mid X},{\pi_2}_{\mid X})&\in I,&(\pi_1,\pi_2)_{\mid Y}=({\pi_1}_{\mid Y},{\pi_2}_{\mid Y})&\in J,\\
(\pi'_1,\pi'_2)_{\mid X}=({\pi'_1}_{\mid X},{\pi'_2}_{\mid X})&\in I,&(\pi'_1,\pi'_2)_{\mid Y}=({\pi'_1}_{\mid Y},{\pi'_2}_{\mid Y})&\in J.
\end{align*}
By the exchange condition for $I$ and $J$,
\begin{align*}
(\pi'_1,\pi_2)_{\mid X}=({\pi'_1}_{\mid X},{\pi_2}_{\mid X})&\in I,&(\pi'_1,\pi_2)_{\mid Y}=({\pi'_1}_{\mid Y},{\pi_2}_{\mid Y})&\in J.
\end{align*}
So $(\pi'_1,\pi_2)\in I\squplus J$. Therefore, $I\squplus J\in \exch[X\sqcup Y]$.

2. Let $\pi'\in I\squplus J$ and $\pi\in \comp[X\sqcup Y]$, such that $\pi\leq \pi'$. Then $\pi_{\mid X}\leq \pi'_{\mid X}$: as $I\in \exch_-[X]$, $\pi_{\mid X}\in I$. Similarly, $\pi_{\mid Y}\in J$, so $\pi\in I\squplus J$.

3. Similar proof. \end{proof}

As a direct consequence:

\begin{prop}\label{prop3.1.7}
The bilinear extension of the product $\squplus$ makes $\ell\exch$ a commutative twisted algebra, which unit is $1=\{\emptyset\}\in \exch[\emptyset]$. Moreover, $\ell\exch_+$ and $\ell\exch_-$ are twisted subalgebras of $\ell\exch$.
\end{prop}

Let us now give a more concrete description of this product $\squplus$.

\begin{notation}\begin{enumerate}
\item For $k,l\geq 0$, we denote by $\qsh(k,l)$ the set of $(k,l)$-quasi-shuffles, that is to say surjective maps $\sigma:[k+l]\longrightarrow [\max(\sigma)]$ such that
\begin{align*}
&\sigma(1)<\cdots<\sigma(k),&&\sigma(k+1)<\cdots<\sigma(k+l).
\end{align*}
\item Let $X$ be a finite set, $\pi=(X_1,\ldots,X_k)\in \comp_k[X]$, and $\sigma:[k]\longrightarrow [l]$ be a surjective map.
We put
\[\sigma\cdot \pi=\left(\bigsqcup_{\sigma(i)=1}X_i,\ldots, \bigsqcup_{\sigma(i)=l}X_i\right)\in \comp_l[X].\]
\end{enumerate}\end{notation}

\begin{prop}\label{prop3.1.8}

Let $X$, $Y$ be disjoint finite sets. If $I\in \exch[X]$ and $Y\in \exch[Y]$, then
\[I\squplus J=\{\sigma\cdot(\pi,\pi')\mid \pi\in I,\:\pi'\in J,\: \sigma\in \qsh(\ell(\pi),\ell(\pi'))\}.\]
\end{prop}

\begin{proof}
$\subseteq$. Let $(X_1,\ldots,X_n)\in I\squplus J$. Then $\pi=\elim(X_1\cap X,\ldots,X_k\cap n)=(X'_1,\ldots,X'_k)\in I$ and $\pi=\elim(X_1\cap Y,\ldots,X_Y\cap n)=(X'_{k+1},\ldots,X'_{k+l})\in J$.
Let $\sigma:[k+l]\longrightarrow [n]$ be defined by $X'_i\subseteq X_{\sigma(i)}$ for any $i\in [k+l]$. Then $\sigma\in \qsh(k,l)$ and $(X_1,\ldots,X_n)=\sigma\cdot (\pi,\pi')$.

$\supseteq$. Let $\pi=(X_1,\ldots,X_k)\in I$, $\pi'=(X_{k+1},\ldots,X_{k+l})\in J$ and $\sigma\in \qsh(k,l)$. Then, as $\sigma(1)<\cdots<\sigma(k)$,
\begin{align*}
\sigma\cdot(\pi,\pi')_{\mid X}&=\elim \left(\bigsqcup_{\sigma(i)=1}X_i\cap X,\ldots, \bigsqcup_{\sigma(i)=\max(\sigma)}X_i\cap X\right)\\
&=\elim\left(\emptyset,\ldots,\emptyset,X_1,\ldots,X_k,\emptyset,\ldots,\emptyset\right)\\
&=\pi.
\end{align*}
Similarly, $\sigma\cdot(\pi,\pi')_{\mid Y}=\pi'$. So $\sigma\cdot(\pi,\pi')_{\mid X}\in I\squplus J$. 
\end{proof}

\subsection{Coproduct on $\exch$}

\begin{defi}
Let $X$, $Y$ be disjoint finite sets, and $I\subseteq \comp[X\sqcup Y]$. We put
\begin{align*}
I_{\mid_L X}&=\{\pi\in \comp[X]\mid \exists \pi'\in \comp[Y],\: (\pi,\pi')\in I\},\\
I_{\mid_R Y}&=\{\pi'\in \comp[Y]\mid \exists \pi\in \comp[X],\: (\pi,\pi')\in I\}.
\end{align*}
\end{defi}

\begin{lemma}\label{lem3.1.10}
Let $X$, $Y$ be disjoint finite sets, and $I\in \exch[X\sqcup Y]$. 
\begin{enumerate}
\item The following conditions are equivalent:
\begin{enumerate}
\item $(X,Y)\in m(I)$.
\item $I_{\mid_L X}\neq \emptyset$ and $I_{\mid_RY}\neq \emptyset$.
\item $I_{\mid_LX}\neq \emptyset$ or $I_{\mid_RY}\neq \emptyset$.
\end{enumerate}
If this holds, $\{(\pi,\pi')\mid \pi\in I_{\mid_LX},\: \pi'\in I_{\mid_RY}\}\subseteq I$.
\item We assume that $(X$, $Y)\in m(I)$. 
\begin{enumerate}
\item If $I\in \exch[X\sqcup Y]$, then $I_{\mid_L X}\in \exch[X]$ and $I_{\mid_R Y}\in \exch[Y]$.
\item If $I\in \exch_-[X\sqcup Y]$, then $I_{\mid_L X}\in \exch_-[X]$ and $I_{\mid_R Y}\in \exch_-[Y]$.
\item If $I\in \exch_+[X\sqcup Y]$, then $I_{\mid_L X}\in \exch_+[X]$ and $I_{\mid_R Y}\in \exch_+[Y]$.
\end{enumerate}\end{enumerate}\end{lemma}

\begin{proof}
1. $(a)\Longrightarrow (b)$. There exists $\pi\in I$, such that $(X,Y)\leq \pi$. By definition of the order, $\pi=(\pi',\pi'')$, with $\pi'\in \comp[X]$ and $\pi''\in \comp[Y]$. Therefore, $\pi'\in I_{\mid_LX}$ and $\pi''\in I_{\mid_RY}$.

$(b)\Longrightarrow (c)$. Obvious.

$(c)\Longrightarrow (a)$. Let us for example take $\pi\in I_{\mid_LX}$. By definition, there exists $\pi''\in \comp[Y]$, such that $(\pi,\pi')\in I$. Then $(X,Y)\leq (\pi,\pi')$, so $(X,Y)\in m(I)$.

Let us take $\pi\in I_{\mid_LX}$ and $\pi'\in I_{\mid_RY}$. There exist $\varpi\in \comp[X]$, $\varpi'\in \comp[Y]$, such that $(\pi,\varpi)$ and $(\pi',\varpi')\in I$. By the exchange condition, $(\pi,\pi')\in I$.

2. By the first item, $I_{\mid_LX}$ and $I_{\mid_RY}$ are nonempty.

$(a)$ Let $(X_1,X_2)\in \comp[X]$, $\pi_i,\pi'_i\in \comp[X_i]$ for $i\in [2]$, such that both $(\pi_1,\pi_2)$ and $(\pi'_1,\pi'_2)$ belong to $I_{\mid_LX}$. By definition of $I_{\mid_LX}$, there exists $\pi_3,\pi'_3\in \comp[Y]$ such that 
\begin{align*}
(\pi_1,\pi_2,\pi_3)&\in I,&(\pi'_1,\pi'_2,\pi'_3)&\in I.
\end{align*}
By the exchange condition for $I$, $(\pi'_1,\pi_2,\pi_3)\in I$, 
so $(\pi_1',\pi_2)\in I_{\mid_LX}$: $I_{\mid_LX}\in \exch[X]$. The proof is similar for $I_{\mid_RY}$.

2. $(b)$ Let $\pi'\in I_{\mid_LX}$ and $\pi\in \comp[X]$, such that $\pi\leq \pi'$. There exists $\varpi\in \comp[Y]$, such that $(\pi',\varpi)\in I$. Moreover, $(\pi,\varpi)\leq (\pi',\varpi)$. As $I\in \exch_-[X]$, $(\pi,\varpi)\in I$, so $\pi\in I_{\mid_LX}$: $I_{\mid_LX}\in \exch_-[X]$.
The proof is similar for $I_{\mid_RY}$.

2. $(c)$ Similar proof.
\end{proof}

\begin{lemma}\label{lem3.1.11}
Let $X$, $Y,Z$ be pairwise disjoint finite sets, and $I\in \exch[X\sqcup Y\sqcup Z]$. The following assertions are equivalent:
\begin{enumerate}
\item $(X$, $Y\sqcup Z)\in m(I)$ and $(Y,Z)\in m(I_{\mid_RY\sqcup Z})$.
\item $(X\sqcup Y,Z)\in m(I)$ and $(X$, $Y)\in m(I_{\mid_LX\sqcup Y})$.
\item $(X$, $Y,Z)\in m(I)$.
\end{enumerate}
If this holds,
\begin{align*}
\left(I_{\mid_L X\sqcup Y}\right)_{\mid_L X}&=I_{\mid_LX},&\left(I_{\mid L X\sqcup Y}\right)_{\mid_RY}&=\left(I_{\mid_R Y\sqcup Z}\right)_{\mid_L Y},&\left(I_{\mid_R Y\sqcup Z}\right)_{\mid_R Z}&=I_{\mid_R Z}.
\end{align*}
\end{lemma}

\begin{proof}
$1.\Longrightarrow 3$. As $(Y,Z)\in m(I_{\mid_RY\sqcup Z})$, there exists $\varpi\in I_{\mid_RY\sqcup Z}$, $(Y,Z)\leq \varpi$. By definition of $\leq$, $\varpi=(\pi',\pi'')$, with $\pi'\in \comp[Y]$ and $\pi''\in \comp[Z]$. As $(\pi',\pi'')\in I_{\mid_RY\sqcup Z}$,
there exists $\pi\in \comp[X]$, $(\pi,\pi',\pi'')\in I$. As $(X$, $Y,Z)\leq (\pi,\pi',\pi'')$, $(X$, $Y,Z)\in m(I)$.

$3.\Longrightarrow 1$. There exists $\varpi\in I$, such that $(X$, $Y,Z)\leq \varpi$. By definition of $\leq$, $\varpi=(\pi,\pi',\pi'')$, with $\pi\in \comp[X]$, $\pi'\in \comp[Y]$ and $\pi''\in \comp[Z]$.
Then $(X$, $Y\sqcup Z)\leq (\pi,\pi',\pi'')$, so $(X$, $Y\sqcup Z)\in m(I)$. Moreover, $(\pi',\pi'')\in I_{\mid_RY\sqcup Z}$ and $(Y,Z)\leq (\pi',\pi'')$, so $(Y,Z)\in m(I_{\mid_RY\sqcup Z})$.

$2.\Longleftrightarrow 3$. Similar proof.

Let us assume that these conditions hold. Then none of $\left(I_{\mid_L X\sqcup Y}\right)_{\mid_L X}$, $I_{\mid_LX}$, $\left(I_{\mid L X\sqcup Y}\right)_{\mid_RY}$, $\left(I_{\mid_R Y\sqcup Z}\right)_{\mid_L Y}$, $\left(I_{\mid_R Y\sqcup Z}\right)_{\mid_R Z}$ and $I_{\mid_R Z}$ is empty.

$\left(I_{\mid_L X\sqcup Y}\right)_{\mid_L X}\subseteq I_{\mid_LX}$. Let $\pi\in \left(I_{\mid_L X\sqcup Y}\right)_{\mid_L X}$. There exists $\pi'\in \comp[Y]$, such that $(\pi,\pi')\in I_{\mid_L X\sqcup Y}$. Therefore, there exists $\pi''\in \comp[Z]$, such that $(\pi,\pi',\pi'')\in I$. 
Hence, $\pi\in I_{\mid_LX}$.

$I_{\mid_LX}\subseteq \left(I_{\mid_L X\sqcup Y}\right)_{\mid_L X}$. Let $\pi \in I_{\mid_LX}$. As $(Y,Z)\in m(I_{\mid_R Y\sqcup Z})$, there exist $\pi'\in \comp[Y]$, $\pi''\in \comp[Z]$, such that $(\pi,\pi')\in I_{\mid_R Y\sqcup Z}$. By Lemma \ref{lem3.1.10}, $(\pi,\pi',\pi'')\in I$, so $(\pi,\pi')\in I_{\mid_LX\sqcup Y}$ and
 finally $\pi \in \left(I_{\mid_L X\sqcup Y}\right)_{\mid_L X}$.

$\left(I_{\mid_R Y\sqcup Z}\right)_{\mid_R Z}=I_{\mid_R Z}$: similar proof.

$\left(I_{\mid L X\sqcup Y}\right)_{\mid_RY}=\left(I_{\mid_R Y\sqcup Z}\right)_{\mid_L Y}$. Let $\pi'\in \comp[Y]$.
\begin{align*}
\pi' \in \left(I_{\mid L X\sqcup Y}\right)_{\mid_RY}&\Longleftrightarrow \exists \pi\in \comp[X],\: (\pi,\pi')\in I_{\mid_L X\sqcup Y}\\
&\Longleftrightarrow \exists \pi\in \comp[X],\:\exists \pi''\in \comp[Z],\: (\pi,\pi',\pi'')\in I\\
&\Longleftrightarrow \exists \pi''\in \comp[Z],\: (\pi',\pi'')\in I_{\mid_RY\sqcup Z}\\
&\Longleftrightarrow \pi'\in \left(I_{\mid_R Y\sqcup Z}\right)_{\mid_L Y}. \qedhere
\end{align*}\end{proof}

\begin{prop}\label{prop3.1.12}
We define a coproduct $\Delta$ on $\exch$ as follows: if $X$, $Y$ are disjoint finite sets, and $I\in \exch[X\sqcup Y]$,
\begin{align*}
\Delta_{X,Y}(I)&=\begin{cases}
I_{\mid_LX}\otimes I_{\mid_RY}\mbox{ if }(X$, $Y)\in m(I),\\
0\mbox{ otherwise}.
\end{cases}\end{align*}
Then $(\ell\exch,\Delta)$ is a twisted coalgebra, and $\ell\exch_+$ and $\ell\exch_-$ are twisted subcoalgebras of $\ell\exch$. 
\end{prop}

\begin{proof}
By Lemma \ref{lem3.1.10}, $\Delta$ is well-defined, and $\exch_+$ and $\exch_-$ are subcoalgebras.
Let us prove the coassociativity of $\Delta$. Let $X$, $Y,Z$ be pairwise disjoint finite sets, and let $I\in \exch[X\sqcup Y\sqcup Z]$. By Lemmas \ref{lem3.1.11},
\begin{align*}
&\left(\Delta_{X,Y}\otimes \id_{\exch[Z]}\right)\circ \Delta_{X\sqcup Y,Z}(I)\\
&=\begin{cases}
\left(I_{\mid_LX\sqcup Y}\right)_{\mid_LX}\otimes \left(I_{\mid_LX\sqcup Y}\right)_{\mid_RY}\otimes I_{\mid_RZ}\mbox{ if }(X\sqcup Y,Z)\in m(I)\mbox{ and }(X$, $Y)\in m\left(I_{\mid_LX\sqcup Y}\right),\\
0\mbox{ otherwise}
\end{cases}\\
&=\begin{cases}
I_{\mid_LX}\otimes \left(I_{\mid_RY\sqcup Z}\right)_{\mid_LY}\otimes \left(I_{\mid_RY\sqcup Z}\right)_{\mid_RZ}\mbox{ if }(X$, $Y\sqcup Z)\in m(I)\mbox{ and }(Y,Z)\in m \left(I_{\mid_RY\sqcup Z}\right),\\
0\mbox{ otherwise}
\end{cases}\\
&=\left(\id_{\exch[X]}\otimes \Delta_{Y,Z}\right)\circ \Delta_{X,Y\sqcup Z}(I).
\end{align*}

Let $X$ be a finite set, and $I\in \exch[X]$. Then 
\begin{align*}
\Delta_{\emptyset,X}(I)&=\{\emptyset\}\otimes X,&\Delta_{X,\emptyset}(I)&=X\otimes \{\emptyset\}.
\end{align*}
So $\Delta$ is counitary, and its counit is defined by $\varepsilon_\Delta(\{\emptyset\})=1$.
\end{proof}

\begin{remark}
The twisted coalgebras $\ell\exch$, $\ell\exch_+$ and $\ell\exch_-$ are not cocommutative. Take $I=\{(1,2,3)\}\in \exch_+[3]$, then
\begin{align*}
\Delta_{\{1\},\{2,3\}}(I)&=\{(1)\}\otimes \{(2,3)\},&\Delta_{\{2,3\},\{1\}}(I)&=0.
\end{align*}
Take $J=\{(123),\:(1,2,3)\}\in \exch_-[3]$, then
\begin{align*}
\Delta_{\{1\},\{2,3\}}(J)&=\{(1)\}\otimes \{(2,3)\},&\Delta_{\{2,3\},\{1\}}(J)&=0.
\end{align*}\end{remark}

\subsection{Bialgebraic structure on $\exch$}

\begin{lemma}\label{lem3.1.13}
Let $X$ be a finite set, and $I\in \exch[X]$. Then 
\[m(I)=\{\pi\in \comp[X]\mid \Delta_\pi(I)\neq 0\}.\]
\end{lemma}

\begin{proof}
$\supseteq$. Let $\pi=(X_1,\ldots,X_k)\in \comp[X]$, such that $\Delta_{X_1,\ldots,X_k}(I)\neq \emptyset$. We proceed by induction on $k$. If $k=1$, then $\pi=(X)$. As $I\neq \emptyset$, choosing any $\pi\in I$, $(X)\leq \pi$, so $(X)\in m(I)$.
Let us assume the result at rank $k-1$, with $k\geq 2$. 
\[\Delta_{X_1,\ldots,X_k}(I)=\left(\id_{\exch[X_1]}\otimes \Delta_{X_2,\ldots,X_k}\right)\circ \Delta_{X_1,X_2\sqcup \cdots \sqcup X_k}(I)\neq 0,\]
so $\Delta_{X_1,X_2\sqcup \cdots \sqcup X_k}(I)\neq 0$ and $\Delta_{X_2,\ldots,X_k}\left(I_{\mid_RX_2\sqcup \cdots \sqcup X_k}\right)\neq 0$. By definition of $\Delta$, $(X_1,X_2\sqcup \cdots \sqcup X_k)\in m(I)$, and, by the induction hypothesis, 
$(X_2,\ldots,X_k)\in m(I_{\mid_RX_2\sqcup \cdots \sqcup X_k})$. Therefore, there exists $\pi'\in I_{\mid_RX_2\sqcup \cdots \sqcup X_k}$, such that $(X_2,\ldots,X_k)\leq \pi'$. Finally, there exists $\pi\in \comp[X]$, such that $(\pi,\pi')\in I$. 
As $(X_1,\ldots,X_k)\leq (\pi,\pi')$, $(X_1,\ldots,X_k)\in m(I)$.

$\subseteq$. Let $\pi=(X_1,\ldots,X_k)\in m(I)$. We proceed by induction on $k$. If $k=1$, $\pi=(X)$, and $\Delta_X(I)=I\neq \emptyset$. Let us assume the result at rank $k-1$, with $k\geq 2$.
There exists $\pi'\in I$, such that $\pi\leq \pi'$. By definition of $\leq$, $\pi'=(\pi'_1,\ldots,\pi'_k)$, with $\pi'_i\in \comp[X_i]$ for any $i\in [k]$. Moreover, 
\[(X_1,X_2\sqcup \cdots \sqcup X_k)\leq \pi\leq \pi',\]
so $(X_1,X_2\sqcup \cdots X_k)\in m(I)$. Therefore, by definition of $\Delta$,
\[\Delta_{X_1,X_2\sqcup \cdots \sqcup X_k}(I)=I_{\mid_LX_1}\otimes I_{\mid_RX_2\sqcup \cdots \sqcup X_k}\neq 0.\]
 Moreover, as $\pi'\in I$, $(\pi'_2,\ldots,\pi'_k)\in I_{\mid_RX_2\sqcup \cdots \sqcup X_k}$, so $(X_2,\ldots,X_k)\in m(I_{\mid_RX_2\sqcup \cdots \sqcup X_k})$. By the induction hypothesis,
\[\Delta_{X_1,\ldots,X_k}(I)=\left(\id_{\exch[X_1]}\otimes \Delta_{X_2,\ldots,X_k}\right)\circ \Delta_{X_1,X_2\sqcup \cdots \sqcup X_k}(I)\neq 0. \qedhere\]
\end{proof}

\begin{lemma}\label{lem3.1.14}
Let $X$, $Y$ be disjoint finite sets, $I\in \exch[X]$ and $J\in \exch[Y]$. Then 
\[m(I\squplus J)=m(I)\squplus m(J).\] 
\end{lemma}

\begin{proof}
$\subseteq$. Let $\pi\in m(I\squplus J)$. There exists $\pi'\in I\squplus J$, such that $\pi\leq \pi'$. Therefore, $\pi_{\mid X}\leq \pi'_{\mid X}$. As $\pi'\in I\squplus J$, $\pi'_{\mid X}\in I$, so $\pi_{\mid X}\in m(I)$. Similarly, $\pi_{\mid Y}\in m(I)$, and finally, $\pi\in m(I)\squplus m(J)$.

$\supseteq$. Let $\pi\in m(I)\squplus m(J)$. There exist $\pi'\in I$, $\pi''\in J$, $\pi_{\mid X}\leq \pi'$ and $\pi_{\mid Y}\leq \pi''$. Moreover, by Proposition \ref{prop3.1.8}, $\pi$ can be obtained by the action of a quasi-shuffle $\sigma$ on $(\pi_{\mid X},\pi_{\mid Y})$, which we write
\[\pi=\sigma\cdot (\pi_{\mid X},\pi_{\mid Y}),\]
see Proposition \ref{prop3.1.8}. Applying an extension $\tau$ of this quasi-shuffle to $(\pi',\pi'')$, we obtain $\varpi\in \comp[X\sqcup Y]$ such that $\pi\leq \varpi$, $\varpi_{\mid X}=\pi'$, and $\varpi_{\mid Y}=\pi''$. Therefore, $\varpi\in I\squplus J$ and $\pi\in m(I\squplus J)$. 
\end{proof}

\begin{lemma}\label{lem3.1.15}
Let $(X$, $Y)$ and $(X',Y')$ be two pairs of disjoint finite sets, such that $X\sqcup Y=X'\sqcup Y'$. Let $I\in \exch[X]$ and $J\in \exch[Y]$. The following assertions are equivalent:
\begin{enumerate}
\item $(X',Y')\in m(I\squplus J)$.
\item $(X',Y')_{\mid X}\in m(I)$ and $(X',Y')_{\mid Y}\in m(J)$.
\end{enumerate}
If this holds, 
\begin{align*}
(I\squplus J)_{\mid_LX'}&=I_{\mid_LX\cap X'}\squplus J_{\mid_LY\squplus X'},&
(I\squplus J)_{\mid_RY'}&=I_{\mid_RX\cap Y'}\squplus J_{\mid_RY\squplus Y'}.
\end{align*}\end{lemma}

\begin{proof}
By Lemma \ref{lem3.1.14},
\begin{align*}
(X',Y')\in m(I\squplus J)&\Longleftrightarrow (X',Y')\in m(I)\squplus m(J)\\
&\Longleftrightarrow (X',Y')_{\mid X}\in m(I) \mbox{ and }(X',Y')_{\mid Y}\in m(J).
\end{align*}
Let us assume that this holds. Let $\pi\in (I\squplus J)_{\mid_LX'}$. There exists $\pi'\in \comp[Y']$, $(\pi,\pi')\in I\squplus J$. Then $(\pi_{\mid X\cap X'},\pi'_{\mid X\cap Y'})\in I$, so $\pi_{\mid X\cap X'}\in I_{\mid_LX\cap X'}$. 
Similarly, $\pi_{\mid Y\cap X'}\in I_{\mid_LY\cap X'}$, so $\pi\in I_{\mid_LX\cap X'}\squplus J_{\mid_LY\squplus X'}$.

Let $\pi\in I_{\mid_LX\cap X'}\squplus J_{\mid_LY\squplus X'}$. Then $\pi_{\mid X\cap X'}\in I_{\mid_LX\cap X'}$, so there exists $\varpi\in \comp[X\cap Y']$, $(\pi_{\mid X\cap X'},\varpi)\in I$. 
Similarly, there exists $\varpi'\in \comp[Y\cap Y']$, $(\pi_{\mid Y\cap X'},\varpi')\in J$. Then
\begin{align*}
(\pi,\varpi,\varpi')_{\mid X}&=(\pi_{\mid X'\cap X},\varpi)\in I,&(\pi,\varpi,\varpi')_{\mid Y}&=(\pi_{\mid X'\cap Y},\varpi')\in J,
\end{align*}
so $(\pi,\varpi,\varpi')\in I\squplus J$, and finally $\pi\in (I\squplus J)_{\mid X'}$.

The proof of $(I\squplus J)_{\mid_RY'}=I_{\mid_RX\cap Y'}\squplus J_{\mid_RY\squplus Y'}$ is similar.
\end{proof}

\begin{theo}\label{theo3.1.16}
$(\exch,\squplus,\Delta)$ is a CTB. Moreover, $\exch_+$ and $\exch_-$ are sub-CTBs of $\exch$.
\end{theo}

\begin{proof}
We already proved that $(\ell\exch,\squplus)$ is a twisted algebra and $(\ell\exch,\Delta)$ a twisted coalgebra. 
It remains to prove the compatibility between $\squplus$ and $\Delta$. 
Let $(X$, $Y)$ and $(X',Y')$ be two pairs of disjoint finite sets, such that $X\sqcup Y=X'\sqcup Y'$, and let $I\in \exch[X]$ and $J\in \exch[Y]$. By Lemma \ref{lem3.1.15},
\begin{align*}
&\Delta_{X',Y'}(I\squplus J)\\
&=\begin{cases}
(I\squplus J)_{\mid_LX'}\otimes (I\squplus J)_{\mid_RY'}\mbox{ if }(X',Y')\in m(I\squplus J),\\
0\mbox{ otherwise}
\end{cases}\\
&=\begin{cases}
I_{\mid_LX\cap X'}\squplus J_{\mid_LY\squplus X'}\otimes I_{\mid_RX\cap Y'}\squplus J_{\mid_RY\squplus Y'}\mbox{ if $(X',Y')_{\mid X}\in m(I)$ and $(X',Y')_{\mid Y}\in m(J)$},\\
0\mbox{ otherwise}
\end{cases}\\
&=\Delta_{X'\cap X,Y'\cap X}(I)\squplus \Delta_{X'\cap Y,Y'\cap Y}(J). 
\end{align*}
By Propositions \ref{prop3.1.7} and \ref{prop3.1.12}, $\exch_+$ and $\exch_-$ are sub-CTBs of $\exch$.
\end{proof}

Let us now describe the discrete elements of $\exch$.

\begin{prop}\label{prop3.1.17}
Let $X$ be a finite set.
\begin{enumerate}
\item The unique element of $\bfD\exch[X]$ is $\comp[X]$. The same result holds for $\exch_+$ and $\exch_-$. 
\item $\exch_+[X]\cap \exch_-[X]=\{\comp[X]\}=\bfD\exch[X]$.
\end{enumerate}\end{prop}

\begin{proof}
1. Let us put $X=\{x_1,\ldots,x_n\}$. Putting $\pi_i=\{(\{x_i\})\}$ for any $i\in [n]$, it is immediate that
\[\pi_1\squplus\cdots \squplus \pi_n=\comp[X].\]
The result follows.

2. Let $I\in \exch_+[X]\cap \exch_-[X]$. As $I\in \exch_-[X]$, it contains the unique minimal element of $\comp[X]$, that is to say $(X)$. As $I\in \exch_+[X]$, $I=\comp[X]$. 
\end{proof}

\section{Terminal property of $\exch$}

In this section, we fix a CTB $\calB$.

\subsection{Characters and special twisted subalgebras}

The definition of the character of a twisted bialgebra is given in Definition \ref{defi1.2.9}.
In the case of a CTB $\calB$, we shall write "$f$ is a character of $\calB$" for "$f$ is a character of $\ell\calB$". 

\begin{defi}\label{defi3.2.1}
Let $\calB$ be a CTB, and $f$ be a character of $\calB$. 
We shall say that $\calB$ is a combinatorial character of $\calB$ if for any finite set $X$, and any $b\in \calB[X]$, $f(b)\in \{0,1\}$. 
\end{defi}
\begin{remark}
In other words, a combinatorial character of $\calB$ is a morphism of twisted algebras from $\calB$ to $\Com$.
\end{remark}

This notion of combinatorial character is related to particular twisted subalgebras:

\begin{defi}\label{defi3.2.2}
Let $\calB$ be a CTB and $ \calA$ be a subspecies of $\calB$. We shall say that $ \calA$ is a special twisted subalgebra of $\calB$ if:
\begin{itemize}
\item $1_\calB\in \calA[\emptyset]$.
\item For any pair $(X$, $Y)$ of disjoint finite sets, for any $a\in \calB[X]$, $b\in \calB[Y]$,
\[ab\in \calA[X\sqcup Y]\Longleftrightarrow a\in \calA[X]\mbox{ and }b\in \calA[Y].\]
\end{itemize}\end{defi}

\begin{prop}
Let $\calB$ be a CTB and $ \calA$ be a sub-CTB of $\calB$. Then $ \calA$ is a special twisted subalgebra of $\calB$.
\end{prop}

\begin{proof}
Let $X$, $Y$ be disjoint finite sets, $a\in \calB[X]$ and $b\in \calB[Y]$. As $ \calA$ is a twisted subalgebra of $\calB$, if $a\in \calA[X]$ and $b\in \calA[Y]$, then $ab\in \calA[X\sqcup Y]$.
Conversely, as $ \calA$ is a twisted subcoalgebra, by Lemma \ref{lem2.1.3}, $\Delta_{X,Y}(ab)=a\otimes b\in \calA[X]\otimes \calA[Y]$, so $a\in \calA[X]$ and $b\in \calA[Y]$.
\end{proof}

\begin{defi}\label{defi3.2.4}
Let $\calB$ be a CTB. A convenient subspecies of $\calB$ is a subspecies $\calC$ of the subspecies of indecomposable elements of $\calB$.
\end{defi}

There is a bijective correspondence between combinatorial characters, special twisted subalgebras, and convenient subspecies:

\begin{prop}\label{prop3.2.5}
Let $\calB$ be a CTB. 
\begin{enumerate}
\item Let $ \calA$ be a special twisted subalgebra $\calB$. For any finite set $X$, for any $b\in \calB[X]$, we put
\[\epsilon_ \calA(b)=\begin{cases}
1\mbox{ if }b\in \calA[X],\\
0\mbox{ otherwise}.
\end{cases}\]
Then $\epsilon_ \calA$ is a combinatorial character of $\calB$.
\item Let $f$ be a character of $\calB$. For any finite set $X$, we put
\[ \calB_f[X]=\{b\in \calB[X]\mid f(b)\neq 0\}.\]
Then $ \calB_f$ is a special twisted subalgebra of $\calB$. 
\item If $\calA$ is a special twisted subalgebra of $\calB$, there exists a unique combinatorial character $f$ of $\calB$ such that $\calA=\calB_f$.
\item Let $\calC$ be a convenient subspecies of $\calB$. For any finite set $X$, we denote by $\calB_\calC[X]$ the set of $b\in \calB[X]$ whose indecomposable components all belong to $\calC$.
Then $\calB_\calC$ is a special twisted subalgebra of $\calB$.
\item If $\calA$ is a special twisted subalgebra of $\calB$, there exists a unique convenient subspecies $\calC$ of $\calB$, such that $\calA=\calB_\calC$. 
\end{enumerate}\end{prop}

\begin{proof}
1. As $1_\calB\in \calA[\emptyset]$, $\epsilon_ \calA(1_\calB)=1$. Let $X$, $Y$ be disjoint finite sets, $a\in \calB[X]$, $b\in \calB[Y]$.
As $ \calA$ is a special twisted subalgebra of $\calB$,
\begin{align*}
\epsilon_ \calA(ab)=1&\Longleftrightarrow ab\in \calA[X\sqcup Y]\\
&\Longleftrightarrow a\in \calA[X]\mbox{ and }b\in \calA[Y]\\
&\Longleftrightarrow \epsilon_ \calA(a)=\epsilon_ \calA(b)=1.
\end{align*}
As $\epsilon_ \calA$ only takes the values $0$ and $1$, $\epsilon_ \calA(ab)=\epsilon_ \calA(a)\epsilon_ \calA(b)$. Therefore, $\epsilon_ \calA$ is a character of $\calB$.

2. As $f(1_\calB)=1$, $1_\calB\in \calA[\emptyset]$. Let $X$, $Y$ be disjoint finite sets, $a\in \calB[X]$, $b\in \calB[Y]$.
\begin{align*}
ab\in \calB_f[X\sqcup Y]&\Longleftrightarrow f(ab)\neq 0\\
&\Longleftrightarrow f(a)f(b)\neq 0\\
&\Longleftrightarrow f(a)\neq 0\mbox{ and }f(b)\neq 0\\
&\Longleftrightarrow a\in \calB_f[X]\mbox{ and } \calB_f[Y].
\end{align*}
Consequently, $ \calA$ is a special twisted subalgebra of $\calB$.

3. \textit{Existence}. Let us prove that $\calA=\calB_{\epsilon_\calA}$. Let $X$ be a finite set, and $b\in \calB[X]$. As $\epsilon_\calA$ only takes the values $0$ and $1$ on $\calB$,
\begin{align*}
b\in \calA[X]&\Longleftrightarrow \epsilon_\calA(b)=1\Longleftrightarrow \epsilon_\calA(b)\neq 0\Longleftrightarrow b\in \calB_{\epsilon_\calA}[X].
\end{align*}

\textit{Unicity}. Let $f$ be a combinatorial character of $\calB$, such that $\calA=\calB_f$. Let $X$ be a finite set, and $b\in \calB[X]$.
\begin{align*}
f(b)=1&\Longleftrightarrow f(b)\neq 0\Longleftrightarrow b\in \calB_f[X]\Longleftrightarrow b\in \calA[X]\Longleftrightarrow \epsilon_\calA(b)=1.
\end{align*}
As both $f$ and $\epsilon_\calA$ are combinatorial characters, $f=\epsilon_\calA$.

4. Let $X$, $Y$ be disjoint finite sets, $a\in \calB[X]$ and $b\in \calB[Y]$. As the indecomposable components of $ab$ are the ones of $a$ and of $b$, $ab\in \calB_\calC[X\sqcup Y]$ if, and only if, $a\in \calB_\calC[X]$ and $b\in \calB_\calC[Y]$.

5. \textit{Existence}. For any finite set $X$, let us denote by $\calC[X]$ the set of indecomposable elements of $\calB[X]$ which belong to $\calA[X]$. 

$\calA[X]\subseteq \calB_\calC[X]$. Let $b\in \calA[X]$, which we decompose into indecomposable elements:
\[b=b_1\cdots b_k,\]
with $b_i\in \calB[X_i]$, indecomposable. As $\calA$ is a special twisted subalgebra of $\calB$, for any $i\in [k]$, $b_i\in \calA[X_i]$, so $b_i\in \calC[X_i]$. As a consequence, $b\in \calB_\calC[X]$.

$\calB_\calC[X]\subseteq \calA[X]$. Let $b\in \calB_\calC[X]$, which we decompose into indecomposable elements:
\[b=b_1\cdots b_k,\]
with $b_i\in \calB[X_i]$, indecomposable. By definition of $\calB_\calC$, $b_i\in \calC[X_i]\subseteq \calA[X_i]$ for any $i\in [k]$. As $\calA$ is a twisted subalgebra of $\calB$, $b\in \calA[X]$.

\textit{Unicity}. Let us assume that $\calB_\calC=\calB_\calD$, where $\calC$ and $\calD$ are convenient subspecies of $\calB$. Let $X$ be a finite set. If $b\in \calD[X]$, then $b$ is indecomposable and belongs to $\calB_\calD[X]=\calB_\calC[X]$, so $b\in \calC[X]$.
Symmetrically, if $b\in \calC[X]$, then $b\in \calD[X]$, so $\calC=\calD$. 
\end{proof}

Not all special twisted subalgebras are twisted subbialgebras. Let us give a criterion on convenient subspecies to obtain twisted subbialgebras.

\begin{prop}
Let $\calB$ be a CTB, and $\calC$ a convenient subspecies of $\calB$. Then $\calB_\calC$ is a twisted subbialgebra of $\calB$ if, and only if, the two following conditions hold:
\begin{enumerate}
\item For any disjoint finite set $X$, $Y$, for any $b\in \calC[X\sqcup Y]$, such that $\Delta_{X,Y}(b)\neq 0$ and $b_{\mid_LX}$ is indecomposable, $b_{\mid_LX}\in \calC[X]$.
\item For any disjoint finite set $X$, $Y$, for any $b\in \calC[X\sqcup Y]$, such that $\Delta_{X,Y}(b)\neq 0$ and $b_{\mid_RY}$ is indecomposable, $b_{\mid_RY}\in \calC[X]$.
\end{enumerate}\end{prop}

\begin{proof}
$\Longleftrightarrow$. Let us prove 1. As $\calB_\calC$ is a twisted subbialgebra, $b_{\mid_LX}\in \calB_\calC[X]$. As it is indecomposable, it belongs to $\calC[X]$. The proof of 2. is similar.

$\Longleftarrow$. Let $X$, $Y$ be disjoint finite sets, $b\in \calB_\calC[X\sqcup Y]$, such that $\Delta_{X,Y}(b)\neq 0$. Let us prove that $b_{\mid_LX}\in \calB_\calC[X]$. We decompose $b$ into indecomposable elements:
\[b=b_1\cdots b_k,\]
with, for any $i\in [k]$, $b_i\in \calC[X_i]$. Then 
\[b_{\mid_LX}={b_1}_{\mid_L X_1\cap X}\cdots {b_k}_{\mid_L X_i\cap X}.\]
Let $b'\in \calB[X']$ be an indecomposable component of $b_{\mid_LX}$. There exists $i\in [k]$, such that $X'\subseteq X_i\cap X$. Then, using the coassociativity of $\Delta$,
\begin{align*}
b'&=({b_i}_{\mid_L X_1\cap X})_{\mid_KX'}={b_i}_{\mid_LX'}.
\end{align*} 
By hypothesis 1., $b'\in \calC[X']$. Therefore, $b_{\mid_LX}\in \calB_\calC[X]$. Using hypothesis 2., we prove that $b_{\mid_RY}\in \calB_\calC[Y]$. 
\end{proof}

\begin{example}\label{ex3.2.1}\begin{enumerate}
\item The combinatorial character associated to the twisted subbialgebra $\calB$ is denoted by $\epsilon^\calB_-$:
\begin{align*}
&\forall b\in \calB[X],&\epsilon^\calB_-(b)&=1.
\end{align*}
The associated convenient subspecies is the subspecies $\calB_-$ of all indecomposable elements of $\calB$.
\item 
The combinatorial character associated to the twisted subbialgebra $\bfD\calB$ is denoted by $\epsilon^\calB_+$:
\begin{align*}
&\forall b\in \calB[X],&\epsilon^\calB_+(b)=\begin{cases}
1\mbox{ if }\ic(b)=|X|,\\
0\mbox{ otherwise}.
\end{cases}\end{align*}
The associated convenient subspecies $\calC_{\calB_+}$ is concentrated in degree 1: for any finite set $X$,
\[\calC_{\calB_+}[X]=\begin{cases}
\calB[X]\mbox{ if }|X|=1,\\
\emptyset\mbox{ otherwise}.
\end{cases}\]
\item For $\calB=\exch$, $\exch_+$ or $\exch_-$, we shall simply write $\epsilon_+$ and $\epsilon_-$.
\end{enumerate}\end{example}

Here are other examples.

\begin{example}\begin{enumerate}
\item Let us take the CTB $\calQ\calP$ of quasi-posets (Example \ref{ex2.1.8}). For any finite set $X$, the set of equivalences $\calE[X]$ on $X$ is a subset of $\calQ\calP[X]$, and it is not difficult to prove that $\calE$ is a twisted subbialgebra of $\calQ\calP$. 
The associated character is given as follows: for any finite set $X$,
\begin{align*}
&\forall P\in \calQ\calP[X],&\epsilon_\calE[P]&=\begin{cases}
1\mbox{ if $\leq_P$ is an equivalence},\\
0\mbox{ otherwise}.
\end{cases}
\end{align*}
The associated convenient subspecies is denoted by $\calQ\calP_\calE$. For any finite set $X$, $\calQ\calP_\calE[X]$ is reduced to a single element, the equivalence $X^2$. 
\item The subspecies of oriented acyclic graphs $\calO \calA_{\mathrm{ac}}$ is a twisted subbialgebra of $\calO \calA$ of Example \ref{ex3.2.4}. The associated character is given as follows: for any finite set $X$, 
\begin{align*}
&\forall G\in \calO \calA[X],&\epsilon_{\mathrm{ac}}(G)&=\begin{cases}
1\mbox{ if $G$ is acyclic},\\
0\mbox{ otherwise}.
\end{cases}\end{align*}\end{enumerate}\end{example}

\begin{prop}\label{prop3.2.7}
For any finite set $X$, we put
\begin{align*}
&\forall I\in \exch[X],&\epsilon(I)&=\begin{cases}
1\mbox{ if $(X)\in I$},\\
0\mbox{ otherwise}.
\end{cases}\end{align*}
Then $\epsilon$ is a combinatorial character of $\exch$. Moreover, for any finite set $X$,
\begin{align*}
&\forall I\in \exch_+[X],&\epsilon(I)&=\epsilon_+(I),\\
&\forall I\in \exch_-[X],&\epsilon(I)&=\epsilon_-(I).
\end{align*}\end{prop}

\begin{proof}
Let us put, for any finite set $X$,
\[ \calA[X]=\{\pi\in \comp[X]\mid (X)\in I\}.\]
Let $X$, $Y$ be disjoint finite sets, $I\in \exch[X]$, $J\in \exch[Y]$. If $(X)\in I$ and $(Y)\in J$, then $(X\sqcup Y)\in I\squplus J$, as $(X\sqcup Y)_{\mid X}=(X)$ and $(X\sqcup Y)_{\mid Y}=(Y)$. Conversely, if $(X\sqcup Y)\in I\squplus J$, then, by definition of $\squplus$, $(X\sqcup Y)_{\mid X}=(X)\in I$ and $(X\sqcup Y)_{\mid Y}=(Y)\in J$.
So $ \calA$ is a special twisted subalgebra of $\exch$, and $\epsilon$ is the combinatorial character associated to $ \calA$, by Proposition \ref{prop3.2.5}. Consequently, $\epsilon$ is a combinatorial character of $\exch$.

Let us assume that $I\in \exch_-[X]$. As $I\neq \emptyset$, it contains $\pi\in \comp[X]$. As $(X)\leq \pi$, $(X)\in I$, so $\epsilon(I)=1=\epsilon_-(I)$. Let us assume that $I\in \exch_+[X]$. Then
\begin{align*}
\epsilon(I)=1&\Longleftrightarrow (X)\in I\Longleftrightarrow I=\comp[X]\Longleftrightarrow \epsilon_+(I)=1.
\end{align*} 
As both $\epsilon$ and $\epsilon_+$ are combinatorial characters of $\exch$, they coincide on $\exch_+$. 
\end{proof}

\begin{remark}
The twisted subalgebra $ \calA$ introduced in this proof is not a twisted subbialgebra: take for example $I=\{(123), (1,2,3)\}\in \calA[3]$. We obtain $I_{\mid_L[2]}=\{(1,2)\}\notin \calA[2]$.
However, the intersection of $ \calA$ with $\exch_+$ and $\exch_-$ are twisted subbialgebras: as $\epsilon_{\mid \exch_\pm}=\epsilon_\pm$,
\begin{align*}
 \calA\cap \exch_+&=\bfD\exch_+,& \calA\cap \exch_-&=\exch_-.
\end{align*}
\end{remark}

We shall need the following definition:

\begin{defi}
Let $X$ be a finite set, and $b\in \calB[X]$. We put
\[\calC_-(b)=\{\pi\in \comp[X]\mid \Delta_\pi(b)\neq 0\}.\]
\end{defi}

\begin{theo}\label{theo3.2.9}
Let $\calB$ be a CTB.
\begin{enumerate}
\item Let $f$ be a combinatorial character of $\calB$. We denote by $ \calB_f$ the associated special twisted subalgebra of $\calB$ associated to $f$ by Proposition \ref{prop3.2.5}: for any finite set $X$,
\[\calB_f[X]=\{b\in \calB[X]\mid f(b)=1\}.\]
Let $X$ be a finite set, and $b\in \calB_f[X]$. We denote by $\calC_f(b)$ the set of compositions $(X_1,\ldots,X_k)\in \calC_-(b)$, such that $\Delta_{X_1,\ldots,X_k}(b)\in \calB_f[X_1]\otimes\cdots\otimes \calB_f[X_k]$. We put
\[\Gamma_f:\left\{\begin{array}{rcl}
\calB&\longrightarrow&\exch\\
b\in \calB[X]&\longmapsto&\begin{cases}
\calC_f(b)\mbox{ if }\calC_f(b)\neq\emptyset,\\
0\mbox{ otherwise}.
\end{cases}\end{array}\right.\] 
Then $\Gamma_f$ is the unique morphism of CTBs such that $\epsilon\circ \Gamma_f=f$.
\item Let $\Gamma:\calB\longrightarrow \exch$ be a morphism of CTBs. There exists a unique combinatorial character $f$ of $\calB$, such that $\Gamma=\Gamma_f$.
\end{enumerate}\end{theo}

\begin{proof}
1. \textit{First step}. Let us prove that for finite set $X$, and any $b\in \calB[X]$, $\calC_f(b)=\emptyset$ or belongs to $\exch[X]$. Let us assume that $\calC_f(b)\neq \emptyset$. 
Let $(X_1,X_2)\in \comp[X]$, $\pi_i,\pi'_i\in \comp[X_i]$ for $i\in [2]$, such that $(\pi_1,\pi_2)$ and $(\pi'_1,\pi'_2)$ belong to $\calC_f(b)$. 
Then both $\Delta_{\pi_1,\pi_2}(b)$ and $\Delta_{\pi'_1,\pi'_2}(b)$ are tensors of elements of $ \calB_f$. Moreover,
\[\Delta_{\pi_1,\pi_2}=(\Delta_{\pi_1}\otimes \Delta_{\pi_2})\circ \Delta_{X_1,X_2}(b).\]
We deduce that $\Delta_{X_1,X_2}(b)\neq 0$. Then 
\[\Delta_{\pi_1,\pi_2}(b)=\Delta_{\pi_1}(b_{\mid_L X_1})\otimes \Delta_{\pi_2}(b_{\mid_R X_2}).\]
As a conclusion, $\Delta_{\pi_1}(b_{\mid_L X_1})$ and $\Delta_{\pi_2}(b_{\mid_R X_2})$ are both tensors of elements of $ \calB_f$. By symmetry, the same holds for $\Delta_{\pi_1'}(b_{\mid_L X_1})$ and $\Delta_{\pi_2'}(b_{\mid_R X_2})$. Finally,
\[\Delta_{\pi'_1,\pi_2}(b)=(\Delta_{\pi'_1}\otimes \Delta_{\pi_2})\circ \Delta_{X_1,X_2}(b)=\Delta_{\pi'_1}(b_{\mid_L X_1})\otimes \Delta_{\pi_2}(b_{\mid_R X_2}),\]
so this is a tensor of elements of $ \calB_f$: we proved that $(\pi'_1,\pi_2)\in \calC_f(b)$. Hence, $\calC_f(b)\in \exch[X]$.

\textit{Second step}. Let Let $X$, $Y$ be disjoint finite sets, $a\in \calB[X]$, and $b\in \calB[Y]$. Let us prove that
\[\calC_f(ab)=\calC_f(a)\squplus\calC_f(b).\]
$\subseteq$. Let $\pi=(X_1,\ldots,X_k)\in \calC_f(ab)$. Then 
\[\Delta_\pi(ab)=\Delta_{X_1\cap X,\ldots,X_k\cap X}(a)\Delta_{X_1\cap Y,\ldots,X_k\cap Y}(b),\]
and this a tensor of elements of $ \calB_f$. By Proposition \ref{prop3.2.5}, $\Delta_{X_1\cap X,\ldots,X_k\cap X}(a)$ and $\Delta_{X_1\cap Y,\ldots,X_k\cap Y}(b)$ are tensors of elements of $ \calB_f$.
Applying $\varepsilon_\Delta$ on the components in $\calB[\emptyset]$, we obtain that $\Delta_{\pi_{\mid X}}(a)$ and $\Delta_{\pi_{\mid Y}}(b)$ are tensors of elements of $ \calB_f$, so $\pi_{\mid X}\in \calC_f(a)$ and $\pi_{\mid Y}\in \calC_f(b)$.
As a conclusion, $\pi\in \calC_f(a)\squplus \calC_f(b)$.

$\supseteq$. Let $\pi=(X_1,\ldots,X_k)\in \calC_f(a)\squplus\calC_f(b)$. Then $\Delta_{\pi_{\mid X}}(a)$ and $\Delta_{\pi_{\mid Y}}(b)$ are tensors of elements of $ \calB_f$.
Inserting the unit $1_\calB$ in the required places, by Lemma \ref{lem2.1.3}, $\Delta_{X_1\cap X,\ldots,X_k\cap X}(a)$ and $\Delta_{X_1\cap Y,\ldots,X_k\cap Y}(b)$ are tensors of elements of $ \calB_f$.
Multiplying these elements together, by Proposition \ref{prop2.2.8}, $\Delta_\pi(ab)$ is a tensor of elements of $ \calB_f$. So $\pi\in \calC_f(b)$.

\textit{Third step}. Let $X$, $Y$ be disjoint finite sets, and $b\in \calB[X\sqcup Y]$, such that $(X$, $Y)\in \calC_-(b)$. Let us prove the following items:
\begin{enumerate}
\item If $\calC_f(b)=\emptyset$, then $\calC_f(b_{\mid_LX})=\emptyset$ or $\calC_f(b_{\mid_RY})=\emptyset$.
\item If $\calC_f(b)\neq \emptyset$, then 
$\calC_f(b_{\mid_LX})=\calC_f(b)_{\mid_LX}$, and $\calC_f(b_{\mid_RY})=\calC_f(b)_{\mid_RY}$.
\end{enumerate}
1. Let us assume that $\calC_f(b_{\mid_LX})\neq \emptyset$ and $\calC_f(b_{\mid_RY})\neq \emptyset$, and let us choose $\pi\in \calC_f(b_{\mid_LX})$ and $\pi'\in \calC_f(b_{\mid_RY})$. Then
\[\Delta_{(\pi,\pi')}(b)=\Delta_{\pi}(b_{\mid_LX})\otimes \Delta_{\pi'}(b_{\mid_RY}),\]
so this is a tensor of elements of $ \calB_f$: $(\pi,\pi')\in \calC_f(b)$, which is consequently nonempty.

2. $\calC_f(b_{\mid_LX})\subseteq \calC_f(b)_{\mid_LX}$. Let $\pi\in \calC_f(b_{\mid_LX})$. By the first item, $\calC_f(b_{\mid_RY})\neq \emptyset$: let us pick $\pi'\in \calC_f(b_{\mid_RY})$. Then 
\[\Delta_{(\pi,\pi')}(b)=\Delta_{\pi}(b_{\mid_LX})\otimes \Delta_{\pi'}(b_{\mid_RY}),\]
and this is a tensor of elements of $ \calB_f$: $(\pi,\pi')\in \calC_f(b)$, and $\pi\in \calC_f(b)_{\mid_LX}$.

$\calC_f(b)_{\mid_LX}\subseteq \calC_f(b_{\mid_LX})$. Let $\pi\in \calC_f(b)_{\mid_LX}$. There exists $\pi'\in \comp[Y]$, such that $(\pi,\pi')\in \calC_f(b)$. Then
\[\Delta_{(\pi,\pi')}(b)=\Delta_{\pi}(b_{\mid_LX})\otimes \Delta_{\pi'}(b_{\mid_RY})\]
is a tensor of elements of $ \calB_f$. Consequently, $\Delta_\pi(b)$ is a tensor of elements of $ \calB_f$, so $\pi\in \calC_f(b_{\mid_LX})$.

$\calC_f(b_{\mid_RY})=\calC_f(b)_{\mid_RY}$: similar proof.

\textit{Fourth step}. Let $X$, $Y$ be disjoint finite sets, and $b\in \calB[X\sqcup Y]$, such that $(X$, $Y)\notin \calC_-(b)$.

Let us assume that $(X$, $Y)\in m(\calC_f(b))$. There exists $\pi\in \calC_f(b)$, such that $(X$, $Y)\leq \pi$. By definition, $\calC_f(b)\subseteq \calC_-(b)$. As $\calC_-(b)\in \exch_-[X\sqcup Y]$, $(X$, $Y)\in \calC_-(b)$. 
We proved that $(X$, $Y)\notin m(\calC_f(b))$.

\textit{Fifth step}. Let us prove that $\Gamma_f$ is a twisted algebra morphism. 
 Let $X$, $Y$ be two disjoint sets, $a\in \calB[X]$, $b\in \calB[Y]$. 
\begin{itemize}
\item If $\calC_f(a)=\emptyset$, by the second step, $\calC_f(ab)=\emptyset\squplus \calC_f(a)=\emptyset$, so $\Gamma_f(ab)=0=\Gamma_f(a)\squplus \Gamma_f(b)$.
\item If $\calC_f(b)=\emptyset$, similarly, $\Gamma_f(ab)=0=\Gamma_f(a)\squplus \Gamma_f(b)$.
\item If both $\calC_f(a)$ and $\calC_f(b)$ are nonempty, by the second step, $\calC_f(ab)=\calC_f(a)\squplus \calC_f(b)$, so $\Gamma_f(ab)=\Gamma_f(a)\squplus \Gamma_f(b)$.
\end{itemize}

Let us now prove that $\Gamma_f$ is a twisted coalgebra morphism. Let $X$, $Y$ be two disjoint sets, and $b\in \calB[X\sqcup Y]$.
If $(X$, $Y)\notin \calC_-(b)$, then $\Delta_{X,Y}(b)=0$. By the fourth step, $(X$, $Y)\notin m(\calC_f(b))$, so $\Delta_{X,Y}(\calC_f(b))=0$. We obtain that in this case,
\[\Delta_{X,Y}\circ \Gamma_f(b)=(\Gamma_f\otimes \Gamma_f)\circ \Delta_{X,Y}(b)=0.\]
Let us now assume that $(X$, $Y)\in \calC_-(b)$. We put $\Delta_{X,Y}(b)=b_{\mid_LX}\otimes b_{\mid_RY}$. 
If $\calC_f(b)=\emptyset$, by the third step (first item), $\calC_-(b_{\mid_LX})=\emptyset$ or $\calC_-(b_{\mid_RY})=\emptyset$, so 
\[\Delta\circ \Gamma_f(b)=(\Gamma_f\otimes \Gamma_f)\circ \Delta_{X,Y}(b)=0.\]
Let us now assume that $\calC_f(b)\neq \emptyset$. By the third step, 
\begin{align*}
\Delta_{X,Y}\circ \Gamma_f(b)&=\calC_f(b)_{\mid_LX}\otimes \calC_f(b)_{\mid_RY}=\calC_f(b_{\mid_LX})\otimes \calC_f(b_{\mid_RY})=(\Gamma_f\otimes \Gamma_f)\circ \Delta_{X,Y}(b).
\end{align*}
So $\Gamma_f$ is a coalgebra morphism.

Let $X$ be a finite set, and $b\in \calB[X]$.
\begin{align*}
\epsilon\circ \Gamma_f(b)=1&\Longleftrightarrow (X)\in \calC_f(b)\Longleftrightarrow b\in \calB_f[X]\Longleftrightarrow f(b)=1.
\end{align*} 
As both $\epsilon\circ \Gamma_f$ and $f$ are combinatorial characters of $\calB$, they are equal.

 \textit{Final step}. Let us prove the unicity. Let $\phi,\psi:\calB\longrightarrow \exch$ be two morphisms of CTBs, such that $\epsilon\circ \phi=\epsilon\circ \psi$. 
Let us prove that for any finite set $X$, and any $b\in \calB[X]$, $\phi(b)=\psi(b)$. We proceed by induction on $|X|$. If $X=\emptyset$, then $b=1_\calB$, and $\phi(b)=\psi(b)=1_{\exch}$. 
Let us assume the result at all ranks $<|X|$. Let us firstly assume that $\phi(b)\neq 0$. Then $m(\phi(b))\neq \emptyset$. Let $\pi\in m(\phi(b))$. 
\begin{itemize}
\item If $\pi=(X)$, then $\epsilon\circ \phi(I)=\epsilon\circ \psi(b)=1$, so $\psi(b)\neq \emptyset$.
\item Otherwise, we put $\pi=(X_1,\ldots,X_k)$, with $k\geq 2$. By the induction hypothesis applied to $X_i$ for $i\in [k]$,
\[\Delta_{X_1,\ldots,X_k}\circ \psi(b)=\psi^{\otimes k}\circ \Delta_{X_1,\ldots,X_k}(b)=\phi^{\otimes k}\circ \Delta_{X_1,\ldots,X_k}(b)=\Delta_{X_1,\ldots,X_k}\circ \phi(b)\neq 0,\]
so $\psi(b)\neq 0$.
\end{itemize}
We proved that $\psi(b)\neq 0$. By symmetry, $\phi(b)=0$ if, and only if, $\psi(b)=0$. We now assume that $\phi(b)=I$ and $\psi(b)=J$ are both non-zero. Let $(X_1,\ldots,X_k)\in I$.
\begin{itemize}
\item If $k\geq 2$, then $\Delta_{X_1,\ldots,X_k}(I)\neq 0$. We put
\[\Delta_{X_1,\ldots,X_k}(I)=I_1\otimes \cdots \otimes I_k.\]
Then $(X_i)\in I_i$ for $i\in [k]$. By the induction hypothesis applied to $X_i$,
\begin{align*}
\Delta_{X_1,\ldots,X_k}(J)&=\Delta_{X_1,\ldots,X_k}\circ \psi(b)\\
&=\psi^{\otimes k}\circ \Delta_{X_1,\ldots,X_k}(b)\\
&=\phi^{\otimes k}\circ \Delta_{X_1,\ldots,X_k}(b)\\
&=\Delta_{X_1,\ldots,X_k}\circ \phi(b)\\
&=\Delta_{X_1,\ldots,X_k}(I)=I_1\otimes \cdots \otimes I_k.
\end{align*}
By definition of $\Delta$, for any $i\in [k]$, there exists $(\pi_{i,1},\ldots,\pi_{i,k})\in J$, with $\pi_{i,j}\in \comp[X_j]$ and $\pi_{i,i}=(X_i)$. By the exchange condition, $(X_1,\ldots,X_k)\in J$.
\item If $k=1$, then $(X_1,\ldots,X_k)=(X)$, and $\epsilon(I)=1$. As $\epsilon\circ\phi=\epsilon\circ \psi$, $\epsilon(J)=1$, so $(X)\in J$.
\end{itemize}
We proved that $I\subseteq J$. By symmetry, $I=J$, so $\phi(b)=\psi(b)$.

2. By composition, $f=\epsilon\circ \Gamma$ is a combinatorial character of $\calB$. By the unicity part of the first item, $\Gamma=\Gamma_f$. 
\end{proof}

\begin{lemma}\label{lem3.2.10}
Let $f$ be a combinatorial character of $\calB$. 
\begin{enumerate}
\item $\Gamma_f$ takes its values in $\exch_+$ if, and only if, the two following conditions are satisfied:
\begin{itemize}
\item $ \calB_f$ is a sub-CTB of $\calB$.
\item For any finite set $X$, and any $b\in \calB_f[X]$, $\calC_-(b)=\comp[X]$.
\end{itemize}
\item $\Gamma_f$ takes its values in $\exch_-$ if, and only if, for any pair $(X$, $Y)$ of disjoint finite set $X$, and for any $b\in \calB[X\sqcup Y]$,
\[\calC_f(b)\neq \emptyset\Longleftrightarrow f(b)=1.\]
\end{enumerate}\end{lemma}

\begin{proof}
1. $\Longrightarrow$. Let $X$ be a finite set, and $b\in \calB_f[X]$. Then $f(b)=1$, so $(X)\in \calC_f(b)$. As $\calC_f(b)\in \exch_+[X]$, $\calC_f(b)=\comp[X]$, so $\calC_-(b)=\comp[X]$.

1. $\Longleftarrow$. Let $X$ be a finite set, $b\in \calB[X]$, $\pi\in \calC_f(b)$, and $\pi'\in \comp[X]$, with $\pi\leq \pi'$. We put $\pi=(X_1,\ldots,X_k)$. Then $\pi'=(\pi'_1,\ldots,\pi'_k)$, with $\pi'_i\in \comp[X_i]$ for $i\in [k]$, by definition of $\leq$. We put
\[\Delta_{X_1,\ldots,X_k}(b)=b_1\otimes\cdots\otimes b_k,\]
with $b_i\in \calB_f[X_i]$ for $i\in [k]$. As $ \calC_-(b_i)=\comp[X_i]$, $\Delta_{\pi'_i}(b_i)$ is non-zero. As $\calB_f$ is a twisted subbialgebra of $\calB$, it is a tensor of elements of $ \calB_f$, and 
\[\Delta_{\pi'}(b)=\Delta_{\pi'_1}(b_1)\otimes\cdots\otimes\Delta_{\pi'_k}(b_k)\]
is consequently a tensor of elements of $ \calB_f$. So $\pi'\in \calC_f(b)$. We proved that, if is is not empty, $\calC_f(b)\in \exch_+[X]$.

2. Let $X$ be a finite set and $b\in \calB[X]$. If $f(b)=1$, then $(X)\in \calC_f(b)$, so $\calC_f(b)\neq \emptyset$.

$\Longrightarrow$. Let $X$ be a finite set, and $b\in \calB[X]$ such that $\calC_f(b)\neq \emptyset$. As $\calC_f(b)\in \exch_-[X]$, it contains the unique minimal element of $\comp[X]$, that is to say $(X)$. Therefore, $f(b)=1$.

$\Longrightarrow$. Let $X$ be a finite set, and $b\in\calB[X]$, such that $\calC_f(b)\neq \emptyset$. Let us prove that $\calC_f(b)\in \exch_-[X]$. Let $\pi'\in \calC_-(b)$, and $\pi=(X_1,\ldots,X_k)\in \comp[X]$, such that $\pi\leq \pi'$.
We put $\pi'=(\pi'_1,\ldots,\pi'_k)$, with $\pi'_i\in \comp[X_i]$ for $i\in [k]$, and $\Delta_{X_1,\ldots,X_k}(b)=b_1\otimes\cdots\otimes b_k$. Then
\[\Delta_{\pi'}(b)=\Delta_{\pi_1}(b_1)\otimes\cdots\otimes \Delta_{\pi_k}(b_k),\]
and this is a tensor of elements of $\calB_f$. Consequently, for any $i\in [k]$, $\Delta_{\pi_i}(b_i)$ is a tensor of elements of $\calB_f$, so $\pi_i\in \calC_f(b_i)$. By hypothesis, $f(b_i)=1$, so $\pi\in \calC_f(b)$. We obtained that $\calC_f(b)\in \exch_-[X]$.
\end{proof}

With the two characters $\epsilon^\calB_-$ and $\epsilon^\calB_-$, we obtain two particular morphisms $\Gamma^\calB_+$ and $\Gamma^\calB_-$.

\begin{cor}\label{cor3.2.11}
Let $\calB$ be a CTB. The following is a morphism of CTBs:
\[\Gamma^\calB_-:\left\{\begin{array}{rcl}
\calB&\longrightarrow&\exch_-\\
b\in \calB[X]&\longmapsto&\calC_-(b).
\end{array}\right.\]
It is the unique morphism of CTBs from $\calB$ to $\exch_-$ such that $\epsilon_-\circ \Gamma^\calB_-=\epsilon^\calB_-$.
\end{cor}

\begin{proof}
This is exactly the morphism $\Gamma^\calB_{\epsilon_-}$. Note that for any $b\in \calB[X]$, $(X)\in \calC_-(b)$, so $\calC_-(b)\neq \emptyset$.
The character $\epsilon^\calB_-$ obviously satisfies the condition of the second item of Lemma \ref{lem3.2.10}, so $\Gamma^\calB_-$ takes its values in $\exch_-$. 
\end{proof}

For $\exch$, $\exch_-$, and $\exch_+$:

\begin{lemma}
Let $X$ be a finite set, and $I\in \exch[X]$. Then $\calC_+(I)\neq\emptyset$ if, and only if, $M(I)\neq \emptyset$.
\end{lemma}

\begin{proof}
$\Longrightarrow$. Let us assume that $\calC_+(I)\neq \emptyset$. As $\calC_+(I)\in \exch_+[X]$, it contains maximal elements of $\comp[X]$, that is to say elements of $\comp_{|X|}[X]$. Let $(\{x_1\},\ldots,\{x_n\})\in \calC_+(I)$.
Then $\Delta_{\{x_1\},\ldots,\{x_n\}}(I)\neq 0$. As $\comp[\{x_i\}]$ is reduced to a single element,
\[\Delta_{\{x_1\},\ldots,\{x_n\}}(I)=\{(\{x_1\})\}\otimes\cdots\otimes \{(\{x_n\})\}.\]
By definition of $\Delta$, $(\{x_1\},\ldots,\{x_n\})\in I$, so $M(I)\neq \emptyset$.

$\Longleftarrow$. Let us assume that $M(I)\neq\emptyset$, and let $(\{x_1\},\ldots,\{x_n\})\in M(I)$. Then 
\[\Delta_{\{x_1\},\ldots,\{x_n\}}(I)=\{(\{x_1\})\}\otimes\cdots\otimes \{(\{x_n\})\}\in \bfD\exch[\{x_1\}]\otimes\cdots\otimes \bfD\exch[\{x_n\}],\]
so $(\{x_1\},\ldots,\{x_n\})\in \calC_+(I)$. 
\end{proof}

\begin{theo}\label{theo3.2.13}
Let $\calB$ be a CTB. 
Let $X$ be a finite set. We put
\begin{align*}
&\forall b\in \calB[X],&\calC_+(b)&=\{(X_1,\ldots,X_k) \in \comp[X]\mid 0\neq\Delta_{X_1,\ldots,X_k}(b)\in \bfD\calB[X_1]\otimes \cdots \otimes \bfD\calB[X_k]\}.
\end{align*}
The following is a morphism of CTBs:
\[\Gamma^\calB_+:\left\{\begin{array}{rcl}
\calB&\longrightarrow&\exch_+\\
b\in \calB[X]&\longmapsto&\begin{cases}
\calC_+(b)\mbox{ if }\calC_+(b)\neq\emptyset,\\
0\mbox{ otherwise}.
\end{cases}\end{array}\right.\]
It is the unique morphism of CTBs from $\calB$ to $\exch_+$ such that $\epsilon_+\circ \Gamma^\calB_+=\epsilon^\calB_+$.
 \end{theo}

\begin{proof}
This is exactly the morphism $\Gamma^\calB_{\epsilon_+}$. Proposition \ref{prop2.2.8}-2 easily implies that if $\epsilon^\calB_+(b)=1$, then $\calC_+(b)=\calC_{\epsilon_+}(b)=\comp[X]$.
By Lemma \ref{lem3.2.10}, $\calC_+$ takes its values in $\exch_+$.
\end{proof}

The following notion will be important in the sequel:

\begin{defi}\label{defi3.2.14}
Let $\calB$ be a CTB. 
We shall say that $\calB$ is non-degenerate if for any finite set $X$, and any $b\in \calB[X]$, $\calC_+(b)\neq \emptyset$. 
\end{defi}

\begin{example}\begin{enumerate}
\item The CTBs of graphs, of posets, and of hypergraphs (whatever the choice of the coproduct is for the last one) are non-degenerate.
\item If $\calB$ is a CTB, then $\bfD\calB$ is a non-degenerate sub-CTB of $\calB$.
\end{enumerate}\end{example}

\begin{example}\label{ex3.2.4}\begin{enumerate}
\item The CTB of quasi-posets $\calQ\calP$ is degenerate: if $X$ is a finite set, and $P\in \calQ\calP[X]$, then $\calC_+(P)\neq \emptyset$ if, and only if, $P$ is a poset. Therefore, $\nd\calQ\calP=\calP$.
\item The CTB of oriented graphs $\calO\calG$ is degenerate: if $X$ is a finite set, and $G\in \calO\calG[X]$, then $\calC_+(G)\neq \emptyset$ if, and only if, $G$ is acyclic. 
Therefore, $\nd\calO\calG$ is the CTB of acyclic oriented graphs $\calO\calG_{\mathrm{ac}}$.
\end{enumerate}\end{example}

\begin{example}\label{ex3.2.5}
The CTBs $\exch$ and $\exch_-$ are degenerate: consider for example $I=\{(X)\}\in \exch_-[X]\subseteq \exch[X]$, with $|X|\geq 2$. Then $\calC_-(I)=\emptyset$, as $M(I)=\emptyset$. Moreover, for any finite set $X$,
\begin{align*}
\nd\exch[X]&=\{I\in \exch[X]\mid M(I)\neq \emptyset\},&\nd\exch_-[X]&=\{I\in \exch_-[X]\mid M(I)\neq \emptyset\}.
\end{align*}
At the opposite, $\exch_+$ is non-degenerate: if $I\in \exch_+[X]$, then it contains maximal elements of $\comp[X]$, so $M(I)\neq\emptyset$, and $\calC_+(I)\neq \emptyset$.
\end{example}

\subsection{Application to $\exch$}

Let us apply these results with $\calB=\exch$. In this case, we shall simply denote $\Gamma_\pm$ instead of $\Gamma^\exch_\pm$.

\begin{theo}\label{theo3.2.15}
For any finite set $X$, and any $I\in \exch[X]$, $\Gamma_-(I)=m(I)$. In particular, if $I\in \exch_-[X]$, then $\Gamma_-(I)=I$. Therefore, $\Gamma_-$ is a projection from $\exch$ to $\exch_-$.
\end{theo}

\begin{proof}
Lemma \ref{lem3.1.13} gives immediately that $\calC_-(I)=m(I)$. The result follows.
\end{proof}

\begin{theo}\label{theo3.2.16}
Let $X$ be a finite set, and $I\in \exch[X]$. Then $\calC_+(I)$ is the set of compositions $\pi\in m(I)$ such that
\begin{align*}
&\forall \pi'\in \comp[X],&\pi\leq \pi'\Longrightarrow \pi'\in m(I).
\end{align*}
If $I\in \exch_+[X]$, then $\calC_+(I)=I$. In other words, $\Gamma_+$ is a projection from $\exch$ to $\exch_+$. 
\end{theo}

\begin{proof}
Obviously, $\calC_+(I)\subseteq \calC_-(I)=m(I)$.
Let $\pi \in \calC_+(I)$. If $\pi\leq \pi'$, as $\calC_+(I)\in \exch_+[X]$, $\pi'\in \calC_+(I)$.

Let us assume that $\pi\in m(I)$ and that if $\pi'\geq \pi$, then $\pi'\in m(I)$. Let us put 
\[\Delta_\pi(I)=I_1\otimes\cdots\otimes I_k,\]
with $I_i\in \exch[X_i]$. If $\pi_i\in \comp[X_i]$ for $i\in [k]$, then $(\pi_1,\ldots,\pi_k)\geq \pi$. By hypothesis, $(\pi_1,\ldots,\pi_k)\in m(I)$, so
\[\Delta_{(\pi_1,\ldots,\pi_k)}(I)=\Delta_{\pi_1}(I_1)\otimes\cdots\otimes \Delta_{\pi_k}(I_k)\neq 0,\]
which implies that for any $i\in [k]$, for any $\pi_i\in \comp[X_i]$, $\Delta_{\pi_i}(I_i)\neq 0$. We conclude that $m(I_i)=\comp[X_i]$, and that $I_i\in \bfD\calB[X_i]$, by Proposition \ref{prop3.1.17}. In other words, $\pi\in \calC_+(I)$.

Let us consider the morphism of CTBs
\[\Gamma_+:\left\{\begin{array}{rcl}
\exch_+&\longmapsto&\exch_+\\
I\in \exch_+[X]&\longmapsto&\calC_+(I).
\end{array}\right.\]
Then $\epsilon_+\circ \Gamma_+=\epsilon_+$. By the unicity in Theorem \ref{theo3.2.13}, $\Gamma_+=\id_{\exch_+}$. 
\end{proof}

\subsection{Reduction to the non-degenerate case}

\begin{theo}\label{theo3.2.17}
Let $\calB$ be a CTB. For any finite set $X$, we denote by $\calI[X]$ the set of elements $b\in \calB[X]$ such that $\calC_+(b)=\emptyset$. Then $\calI$ is a biideal of $\calB$. The quotient is denoted by $\nd\calB=\calB\setminus \calI$. For any finite set $X$,
\[\nd\calB[X]=\{b\in \calB[X]\mid \calC_+(b)\neq \emptyset\}.\]
Moreover, $\nd\calB$ is non-degenerate, and for any finite set $X$, and any $\overline{b}\in \nd\calB[X]$,
\begin{align*}
\ic(\overline{b})&=\ic(b),&\calC_+(\overline{b})&=\calC_+(b).
\end{align*} \end{theo}

\begin{proof}
By Proposition \ref{prop2.1.8} with the morphism $\Gamma_+:\calB\longrightarrow \exch_+$ of Theorem \ref{theo3.2.13}, $\calI$ is a biideal of $\calB$.

Let $X$ be a finite set, and $b\in \calB[X]\setminus \calI[X]$, indecomposable. Let us assume that $\overline{b}=\overline{b'}\overline{b''}$, with $b'\in \calB[X']\setminus \calI[X']$, $b''\in \calB[X'']\setminus \calI[X']$, and $X=X'\sqcup X''$. We can write in $\calB[X]$ that
\[b=b'b''+\sum_{a\in \calI[X]} x_aa,\]
where $x_a\in \K$ for any $a\in\calI[X]$. As $b\notin \calI[X]$, identifying in the basis $\calB[X]$, necessarily $b=b'b''$. As $b$ is indecomposable, $X'=\emptyset$ or $X''=\emptyset$, so $\overline{b}$ is indecomposable.

Let now $X$ be a finite set, and $b\in \calB[X]\setminus \calI[X]$, not necessarily indecomposable. Let us decompose $b$ into indecomposable elements:
\[b=b_1\cdots b_k,\]
with $k=\ic(b)$. Then
\[\overline{b}=\overline{b_1}\cdots \overline{b_k},\]
and the preceding observation proves that $\overline{b_i}$ is indecomposable for any $i\in [k]$. Therefore, this is a decomposition of $\overline{b}$ into indecomposable elements, so $\ic(\overline{b})=k=\ic(b)$.

Let $X$ be a finite set, and $b\in \calB[X]\setminus \calI[X]$. We consider $(X_1,\ldots,X_k)\in \calC_+(b)$. We put
\[\Delta_{X_1,\ldots,X_k}(b)=b_1\otimes\cdots\otimes b_k.\]
For any $i\in [k]$, $b_i\in \bfD\calB[X_i]$, so $\ic(b_i)=|X_i|$. Then
\[\Delta_{X_1,\ldots,X_k}(\overline{b})=\overline{b_1}\otimes\cdots\otimes \overline{b_k},\]
and for any $i\in [k]$, $\ic(\overline{b_i})=\ic(b_i)=|X_i|$, so $\overline{b_i}\in \bfD\nd\calB[X_i]$. We obtain that $(X_1,\ldots,X_k)\in \calC_+(\overline{b})$.

We now consider $(X_1,\ldots,X_k)\in \calC_+(\overline{b})$. We put
\[\Delta_{X_1,\ldots,X_k}(\overline{b})=\overline{b_1}\otimes\cdots\otimes \overline{b_k},\]
with $b_i \in\calB[X_i]\setminus\calI[X_i]$ for any $i\in [k]$. Then $\overline{b_i}\in \bfD\nd\calB[X_i]$,
so $\ic(\overline{b_i})=\ic(b_i)=|X_i|$. Finally, $b_i\in \bfD\calB[X_i]$, so $(X_1,\ldots,X_k)\in \calC_+(b)$. We proved that $\calC_+(\overline{b})=\calC_+(b)$.

By definition, for any finite set $X$, for any $b\in \calB[X]\setminus \calI[X]$, $\calC_+(\overline{b})=\calC_+(b)\ne \emptyset$, so $\nd\calB$ is non-degenerate. 
\end{proof}

\begin{remark}
The conservation of the number of indecomposable components implies that \[\bfD\nd\calB=\bfD\calB.\]
\end{remark}

\begin{prop}\label{prop3.2.18}
Let $\calB$ and $\calB'$ be two twisted CTBs, and $\phi:\calB\longrightarrow\calB'$ be a morphism of CTBs. We assume that $\calB'$ is non-degenerate. 
Then $\phi$ induces a morphism of CTBs from $\nd\calB$ to $\calB'$.
\end{prop}

\begin{proof}
It is enough to prove that if $X$ is a finite set, and $b\in \calB[X]$ such that $\calC_+(b)=\emptyset$, then $\phi(b)=0$. By duality, let us prove that if $\phi(b)\neq 0$, then $\calC_+(b)\neq 0$. 
If $\phi(b)\neq 0$, then $\phi(b)\in \calB'[X]$. As $\calB'$ is non-degenerate, there exists an indexation $\{x_1,\ldots,x_k\}$ of $X$, such that $(\{x_1\},\ldots,\{x_k\})\in \calC_+(\phi(b))$. Then
\[\phi^{\otimes k}\circ \Delta_{\{x_1\},\ldots,\{x_k\}}(b)=\Delta_{\{x_1\},\ldots,\{x_k\}}\circ \phi(b)\neq 0,\]
so $(\{x_1\},\ldots,\{x_k\})\in \calC_-(b)$. As its parts are of cardinality 1, $(\{x_1\},\ldots,\{x_k\})\in \calC_+(b)$,
so $\calC_+(b)\neq 0$.
\end{proof}

\subsection{From $\ell\exch$ to $\ell\comp$}

Let us go back to the (not combinatorial) twisted bialgebra $\ell\comp$:

\begin{prop}\label{prop3.3.1}
The following defines a morphism of twisted bialgebras from $\ell\exch$ to $\ell\comp$: for any finite set $X$,
\begin{align*}
&\forall I\in \exch[X],&\Theta(I)&=\sum_{\pi\in I}\pi.
\end{align*}\end{prop} 

\begin{proof}
By the terminal property of $\ell\comp$ (Theorem \ref{theo1.5.6}), there exists a unique morphism of twisted bialgebras $\Theta:\ell\exch\longrightarrow \ell\comp$, such that $\epsilon\circ\Theta=\epsilon$. Moreover, if $X$ is a nonempty finite set, and $I\in \exch[X]$, the proof of Proposition \ref{prop1.5.4} gives that
\begin{align*}
\Theta(I)&=\sum_{k=1}^{|X|}\sum_{(X_1,\ldots,X_k)\in \comp_k[X]} \epsilon^{\otimes k}\circ \Delta_{X_1,\ldots,X_k}(I)(X_1,\ldots,X_k).
\end{align*}
Observe firstly that $(X_1,\ldots,X_k)\in \comp_k[X]$, then $\epsilon^{\otimes k}\circ \Delta_{X_1,\ldots,X_k}(I)=0$ or $1$.
Let us assume that $\epsilon^{\otimes k}\circ \Delta_{X_1,\ldots,X_k}(I)=1$. Then $(X_1,\ldots,X_k)\in m(I)$, and for any $i\in [k]$, 
$\epsilon\left(I_{\mid_i X_i}\right)=1$, so $(X_i)\in I_{\mid_i X_i}$. 
The exchange condition gives that $(X_1,\ldots,X_k)\in I$. Conversely, if $(X_1,\ldots,X_k)\in I$, then $(X_1,\ldots,X_k)\in m(I)$, and for any $i\in [k]$, $(X_i)\in I_{\mid_iX_i}$, so $\epsilon^{\otimes k}\circ \Delta_{X_1,\ldots,X_k}(I)=1$. We obtain that 
\begin{align*}
\Theta(I)&=\sum_{k=1}^{|X|}\sum_{(X_1,\ldots,X_k)\in \comp_k[X]\cap I}(X_1,\ldots,X_k)=\sum_{\pi\in I}I. \qedhere
\end{align*}\end{proof}

\begin{cor}\label{cor3.3.2}
Let $\calB$ be a CTB and $f$ be a combinatorial character of $\calB$. 
The unique morphism of twisted bialgebras $\theta_f:\ell\calB\longrightarrow\ell\comp$ such that $\epsilon\circ \theta_f=f$ of Theorem \ref{theo1.5.6} is $\Theta\circ \Gamma_f^\calB$. For any finite set $X$, and any $b\in \calB[X]$,
\[\Theta\circ \Gamma_f^\calB(b)=\sum_{\pi\in \calC_f(b)}\pi.\]
\end{cor}

\begin{proof}
By composition, $\Theta\circ \Gamma_f^\calB$ is a morphism of twisted bialgebras. Moreover,
\[\epsilon\circ \Theta\circ \Gamma_f^\calB=\epsilon\circ \Gamma_f^\calB=f.\]
The formula for $\Theta\circ \Gamma_f^\calB(b)$ is immediate.
\end{proof}

As a consequence, applying the Fock functors of \cite{Aguiar2010}, we obtain interesting morphisms of Hopf algebras. Applying the bosonic Fock functor $\calF$, we obtain a morphism $\calF[\Theta]:\calF[\ell\exch]\longrightarrow \QSym$.
If $X$ is a finite set,
\begin{align*}
&\forall I\in \exch[X],&\calF[\Theta](\overline{I})&=\sum_{\pi=(X_1,\ldots,X_k)\in I} (|X_1|,\ldots,|X_k|).
\end{align*}
Applying the full Fock functor $\calF_F$, we obtain a morphism $\calF_F[\Theta]=\calF_F[\ell\exch]\longrightarrow \WQSym$, where $\WQSym$ is the Hopf algebra of packed words of \cite{Novelli2006-2}. Let us describe this morphism. 
If $\pi=(X_1,\ldots,X_p)\in \comp[n]$, we associate to $\pi$ the packed word $w_\pi=k_1\ldots k_n$ of length $n$, such that for any $i\in [n]$, $i\in X_{k_i}$. Then
\begin{align*}
&\forall I\in \exch[n],&\calF_F[\Theta](I)&=\sum_{\pi\in I} w_\pi.
\end{align*}

To end the study of the links between $\exch$ and $\comp$, let us give a maximality result:

\begin{prop}\label{prop3.3.3}
Let $\calB$ be a subspecies of $\calS(\comp)$, such that $\emptyset \notin \calB[X]$ for any finite set $X$.
We assume that $\calB$ is given a CTB structure, such that the following map is a morphism of CTB:
\[\Theta:\left\{\begin{array}{rcl}
\calB&\longrightarrow&\comp\\
I\in \calB[X]&\longmapsto&\displaystyle \sum_{\pi\in I}\pi.
\end{array}
\right.\]
Then $\calB$ is a sub-CTB of $\exch$.
\end{prop}

\begin{proof}
Let us prove that the coproduct of $\calB$, which we denote by $\blacktriangle$, is the restriction of $\Delta$. 
Let $(X,Y)$ be a pair of disjoint finite sets, and $I\in \calB[X\sqcup Y]$.

Let us assume that $(X,Y)\in \calC_-(I)$. We put $\blacktriangle(I)=I_1\otimes I_2$, with $I_1\in \calB[X]$ and $I_2\in \calB[Y]$.
Let $\pi_1\in I_1$, and $\pi_2\in I_2$. Then $\pi_1\otimes \pi_2$ appears in $(\Theta\otimes \Theta)\circ\blacktriangle(I)=\Delta\circ \Theta(I)$, so $(\pi_1,\pi_2)$ appears in $\Theta(I)$, and $(\pi_1,\pi_2)\in I$.
We obtain that $\pi_1\in I_{\mid_LX}$ and $\pi_2\in I_{\mid_RY}$. 
Conversely, if $\pi_1\in I_{\mid_LX}$, there exists $\pi_1\in \comp[Y]$, such that $(\pi_1,\pi_2)\in I$. So $(\pi_1,\pi_2)$ appears in $\Theta(I)$,
and $\pi_1\otimes \pi_2$ appears in $\Delta\circ \Theta(I)=(\Theta\otimes \Theta)\circ \blacktriangle(I)$.
As a consequence, $\pi_1\in I_1$. We proved that $I_1=I_{\mid_LX}$. Similarly, $I_2=I_{\mid_RY}$.
Finally, in this case, $\blacktriangle(I)=\Delta(I)$.

Let us assume that $(X,Y)\notin \calC_-(I)$. Therefore, $\blacktriangle_{X,Y}(I)=0$, so $(\Theta\otimes \Theta)\circ \blacktriangle(I)=\Delta\circ \Theta(I)=0$. As a consequence, for any $\pi\in I$, $\Delta_{X,Y}(\pi)=0$.
This gives that $I_{\mid_LX}=\emptyset$ and $I_{\mid_RY}=\emptyset$, so $\Delta_{X,Y}(I)=0$. 
In this case also, $\blacktriangle(I)=\Delta(I)$.

Let $X$ be a finite set, and $I\in \calB[X]$. Let us show that $I$ satisfies the exchange condition. Let $(X_1,X_2)\in \comp[X]$,
$\pi_i,\pi'_i\in \comp[X_i]$ for $i\in [2]$, such that $(\pi_1,\pi_2)$ and $(\pi'_1,\pi'_2)$ belong to $I$.
Therefore, $(\pi_1,\pi_2)$ and $(\pi'_1,\pi'_2)$ both appear in $\Theta(I)$, so $\pi_1\otimes \pi_2$ and $\pi'_1\otimes \pi'_2$ both appear in $\Delta_{X_1,X_2}\circ \Theta(I)=(\Theta\otimes \Theta)\circ \Delta_{X_1,X_2}(I)$.
Consequently, $(X_1,X_2)\in \calC_-(I)$, 
$\pi_1$, $\pi'_1\in I_{\mid_L X_1}$ and $\pi_2$, $\pi'_2\in I_{\mid_R X_2}$. As a consequence, $\pi_1\otimes \pi'_2$ appears in $(\Theta\otimes \Theta)\circ \Delta_{X_1,X_2}(I)=\Delta_{X_1,X_2}\circ \Theta(I)$, so $(\pi_1,\pi'_2)$ appears in $\Theta(I)$ and finally $(\pi_1,\pi'_2)\in I$.

Let us finally prove that the product of $\calB$ is the restriction of $\squplus$. We denote by $\cdot$ the product of $\calB$.
Let $(X,Y)$ be a pair of disjoint finite sets, $I\in \calB[X]$, and $J\in \calB[Y]$. As $\Theta$ is a twisted algebra morphism, 
\begin{align*}
\Theta(I\cdot J)&=\Theta(I)\squplus \Theta(J)=\sum_{\pi\in I}\sum_{\pi'\in J} \pi\squplus \pi'
=\sum_{k=0}^{|X|}\sum_{l=0}^{|Y|} \sum_{\substack{\pi\in \comp_k[X]\cap I,\\ \pi'\in \comp_l[Y]\cap J,\\ \sigma\in \qsh(k,l)}}
\sigma\cdot(\pi\otimes \pi').
\end{align*}
As $I\cdot J\in \calB[X\sqcup Y]$, $\calB$ being a CTB,
\[I\cdot J=\{\sigma\cdot(\pi\otimes \pi')\mid k\leq |X|,\: l\leq |Y|,\: \pi\in \comp_k[X]\cap I,\: \pi'\in \comp_l[Y]\cap J,\: \sigma\in \qsh(k,l)\}=I\squplus J.\]
Therefore, $\cdot$ is the restriction of $\squplus$ to $\calB$.\end{proof}

\section{Structure on $\exch_{\pp}$}

The subspecies $\exch_{\pp}$ is stable under $\Delta$, but not under $\squplus$. However, it has a structure of CTB, which we now detail.

\subsection{The shuffle product}

\begin{theo}\label{theo3.4.1}\begin{enumerate}
\item Let $X$, $Y$ be disjoint finite sets. For $I\in \exch_{\pp}[X]$ and $J\in \exch_{\pp}[J]$, we put 
\[I\shuffle J=(I\squplus J)\cap \comp_{|X|+|Y|}[X\sqcup Y] \in \exch_{\pp}[X\sqcup Y].\]
Then $(\exch_{\pp},\shuffle,\Delta)$ is a CTB. 
\item Let $\calB$ be a CTB. For any finite set $X$, and any $b\in \calB[X]$, we put
\[\calC_{\pp}(b)=\calC_-(b)\cap \comp_{|X|}[X].\]
The following is a morphism of CTBs:
\[\Gamma^\calB_{\pp}:\left\{\begin{array}{rcl}
\calB&\longrightarrow&(\exch_{\pp},\shuffle, \Delta)\\
b\in \calB[X]&\longmapsto&\begin{cases}
\calC_{\pp}(b)\mbox{ if }\calC_{\pp}(b)\neq \emptyset,\\
0\mbox{ otherwise}.
\end{cases}
\end{array}\right.\]
\end{enumerate}\end{theo}

\begin{proof}
2. We firstly prove that $\Gamma^\calB_{\pp}$ is compatible with the products and the coproducts.
Let $X$, $Y$ be disjoint finite sets, $a\in \calB[X]$ and $b\in \calB[Y]$. As $\Gamma^\calB_-$ is a morphism of twisted algebras,
\[\calC_-(ab)=\calC_-(a)\squplus \calC_-(b).\]
By definition of $\squplus$, intersecting with $\comp_{|X|+|Y|}[X\sqcup Y]$,
\[\calC_{\pp}(ab)=\calC_{\pp}(a)\shuffle \calC_{\pp}(b).\]
So $\Gamma^\calB_{\pp}$ is compatible with the products.

Let $b\in \calB[X\sqcup Y]$. As $\Gamma^\calB_-$ is a morphism of twisted algebras,
\[\Delta_{X,Y}\circ \Gamma^\calB_-(b)=\left(\Gamma^\calB_-\otimes \Gamma^\calB_-\right)\circ \Delta_{X,Y}(b).\]
Intersecting the tensors on both sides with $\comp_{|X|}[X]$ on the left, and $\comp_{|Y|}[Y]$ on the right, we obtain
\[\Delta_{X,Y}\circ \Gamma^\calB_{\pp}(b)=\left(\Gamma^\calB_{\pp}\otimes \Gamma^\calB_{\pp}\right)\circ \Delta_{X,Y}(b).\]
So $\Gamma^\calB_{\pp}$ is a compatible with the coproducts.

1. Let us firstly apply the second point with $\calB=\exch$. Then for any finite set $X$, any $I\in \exch[I]$,
\[\calC_{\pp}(I)=\calC_-(I)\cap \comp_{|X|}[X]=I\cap\comp_{|X|}[X].\]
Consequently, if $I\in \exch_{\pp}[X]$, $\Gamma^{\exch}_{\pp}(I)=I$, so $\Gamma^{\exch}_{\pp}$ is surjective.
As $\exch$ is a twisted bialgebra and $\Gamma^{\exch}_{\pp}$ is compatible with the products and the coproducts,
by taking the quotient, $(\exch_{\pp},\shuffle,\Delta)$ is a twisted bialgebra. It is then immediate that it is a CTB.
\end{proof}

\begin{remark}\label{rk3.4.1}
If $X$ and $Y$ are disjoint finite sets, $I\in \exch_{\pp}[X]$ and $J\in \exch_{\pp}[Y]$,
\[I\shuffle J=\{\sigma\cdot (\pi,\pi')\mid \pi\in I,\: \pi'\in J,\: \sigma \in \sh(|X|,|Y|)\},\]
where for any $k,l\geq 0$, $\sh(k,l)$ is the set of $(k,l)$-shuffles, that is to say of bijective $(k,l)$-quasi-shuffles. 
\end{remark}

\begin{example}
For the CTB $\calG$ of graphs (Example \ref{ex2.1.2}), if $X$ is a finite set and $G\in \calG[X]$, then $\calC_-(G))=\comp[X]$, so $\Gamma^\calB_{\pp}(G)=\comp_{|X|}[X]$.
The same results holds for hypergraphs (Example \ref{ex2.1.4}).
\end{example}

\begin{example}
For the CTB $\calP$ of posets (Example \ref{ex2.1.5}), if $X$ is a finite set, of cardinality denoted by $n$, and $P\in \calP[X]$, then 
\[(\{x_1\},\ldots,\{x_n\})\in \calC_{\pp}(P))\Longleftrightarrow (\forall i,j\in [n], x_i\leq_P x_j\Longrightarrow i\leq j).\]
Therefore, $\calC^\calP_{\pp}(P)$ is in bijection with the set of total orders on $X$ which refine $\leq_P$ (they are called linear extensions of $P$).
This is easily generalized to quasi-posets: if $P\in \calQ\calP[X]$, then
\[\Gamma^{\calQ\calP}_{\pp}(P)=\begin{cases}
\Gamma^\calP_{\pp}(P)\mbox{ if }P\in \calP[X],\\
0\mbox{ otherwise}.
\end{cases}\]\end{example}

\begin{example}
Let us consider the CTB $\calO\calG$ of oriented graphs (Example \ref{ex2.1.3}). Similarly to the case of posets, if $X$ is a finite set, of cardinality denoted by $n$, and $G\in \calO\calG[X]$, then $\calC^{\calO\calG}_{\pp}(G)$ is in bijection with the set of bijections $c:X\longrightarrow [n]$ such that
\begin{align*}
&\forall x,x'\in X,&(x,x')\in A(G)&\Longrightarrow c(x)<c(x').
\end{align*}
In particular, if $G$ is not acyclic, $\Gamma^{\calO\calG}_{\pp}(G)=0$. Moreover, 
$\Gamma^{\calO\calG}_{\pp}=\Gamma^{\calQ\calP}_{\pp}\circ \Gamma_{\calO\calG}$.
\end{example}

\subsection{Morphisms $\Gamma^{\exch_{\pp}}_\pm$}

\begin{lemma}\label{lem3.4.2}
For any finite set $X$, the unique element of $\bfD\exch_{\pp}[X]$ is $\comp_{|X|}[X]$.
\end{lemma}

\begin{proof}
If $I\in \bfD\exch_{\pp}[X]$, putting $X=\{x_1,\ldots,x_k\}$, then 
\[I=\{(\{x_1\})\}\shuffle\cdots\shuffle\{(\{x_k\})\}=\comp_k[X],\]
which gives the result.
\end{proof}

\begin{prop}\label{prop3.4.3}\begin{enumerate}
\item Let $X$ be a finite set, and $I\in \exch_{\pp}[X]$. Then 
\[\calC_+(I)=\{(X_1,\ldots,X_k)\in \comp[X]\mid \comp_{|X_1|}[X_1]\cdot \:\cdots\: \cdot \comp_{|X_k|}[X_k]\subseteq I\}.\]
Consequently, $\calC_+(I)\cap \comp_{|X|}[X]=I$, and $\Gamma^{\exch_{\pp}}_+$ is injective.
\item Let $X$ be a finite set, and $I\in \exch_{\pp}[X]$. Then 
\[\calC_-(I)=\{(X_1,\ldots,X_k)\in \comp[X]\mid \comp_{|X_1|}[X_1]\cdot \:\cdots\: \cdot \comp_{|X_k|}[X_k]\cap I\neq \emptyset\}.\]
Consequently, $\calC_-(I)\cap \comp_{|X|}[X]=I$, and $\Gamma^{\exch_{\pp}}_-$ is injective.
\end{enumerate}\end{prop}

\begin{proof}
1. $\subseteq$. If $(X_1,\ldots,X_k)\in \calC_+(I)$, putting
\[\Delta_{X_1,\ldots,X_k}=I_{\mid_1 X_1}\otimes\cdots \otimes I_{\mid_k X_k},\]
then for any $i\in [k]$, by Lemma \ref{lem3.4.2}, $I_{\mid_i X_i}=\exch_{\pp}[X_i]$. Let us choose $\varpi_i\in \exch_{\pp}[X_i]$ for any $i\in [k]$. By definition of the coproduct,
there exist $\pi_{i,j}\in \comp[X_j]$ for $j\in [k]$, such that $\pi_{i,1}\cdot\:\cdots\:\cdot \pi_{i,k}\in I$ and $\pi_{i,i}=\varpi_i$.
By the exchange condition, $\varpi_1\cdot\:\cdots\:\cdot \varpi_k\in I$, so
 $\comp_{|X_1|}[X_1]\cdot \:\cdots\: \cdot \comp_{|X_k|}[X_k]\subseteq I$.

$\supseteq$. If $\comp_{|X_1|}[X_1]\cdot \:\cdots\: \cdot \comp_{|X_k|}[X_k]\subseteq I$, by definition of the coproduct,
\[\Delta_{X_1,\ldots,X_k}(I)=\comp_{|X_1|}[X_1]\otimes\cdots\otimes\comp_{|X_k|}[X_k].\]
By Lemma \ref{lem3.4.2}, $(X_1,\ldots,X_k)\in \calC_+(I)$.

If $(X_1,\ldots,X_k)\in \calC_+(I)\cap \comp_{|X|}(I)$, then, as $|X_i|=1$ for any $i\in [k]$,
\[\Delta_{X_1,\ldots,X_k}(I)=\{(X_1)\}\otimes\cdots\otimes \{(X_k)\}.\]
By definition of the coproduct, $(X_1,\ldots,X_k)\in I$. Conversely, if $(X_1,\ldots,X_k)\in I$,
then $(X_1,\ldots,X_k)\in m(I)=\calC_+(I)$, so $(X_1,\ldots,X_k)\in \calC_+(I)\cap\comp_{|X|}[X]$.

2. $\subseteq$. Let $(X_1,\ldots,X_k)\in \calC_-(I)=m(I)$. There exists $\pi\in I$, such that $(X_1,\ldots,X_k)\leq\pi$.
By definition of the order on $\comp[X]$, $\pi=(\pi_1,\ldots,\pi_k)$, with $\pi_i\in \comp[X_i]$. Moreover, as $\pi\in \comp_{|X|}[X]$, for any $i\in [k]$, $\pi_i\in \comp_{|X_i|}[X_i]$. So 
\[\pi\in \comp_{|X_1|}[X_1]\cdot \:\cdots\: \cdot \comp_{|X_k|}[X_k]\cap I. \]

$\supseteq$. Let $\pi \in \comp_{|X_1|}[X_1]\cdot \:\cdots\: \cdot \comp_{|X_k|}[X_k]\cap I$. Then $\pi=(\pi_1,\ldots,\pi_k)$,
with $\pi_i\in \comp[X_i]$. Therefore, $\pi\geq (X_1,\ldots,X_k)$, so $(X_1\ldots,X_k)\in m(I)=\calC_-(I)$.

The proof of $\calC_-(I)\cap \comp_{|X|}[X]=I$ is similar to the end of the first item, and we leave it to the reader.
\end{proof}

\begin{prop}\label{prop3.4.4}
Let $\calB$ be a CTB. Then
\begin{align*}
\Gamma^{\exch_{\pp}}_+\circ \Gamma^\calB_{\pp}&=\Gamma^\calB_{\epsilon'^\calB_+},&
\Gamma^{\exch_{\pp}}_-\circ \Gamma^\calB_{\pp}&=\Gamma^\calB_{\epsilon'^\calB_-},
\end{align*}
where $\epsilon'^\calB_+$ and $\epsilon'^\calB_-$ are the combinatorial characters of $\calB$ defined by
\begin{align*}
&\forall b\in \calB[X],&\epsilon'^\calB_+(b)&=\begin{cases}
1\mbox{ if }\calC_-(b)=\comp[X],\\
0\mbox{ otherwise},
\end{cases}&\epsilon'^\calB_-(b)&=\begin{cases}
1\mbox{ if }\calC_+(b)\neq \emptyset,\\
0\mbox{ otherwise}.
\end{cases}\end{align*}\end{prop}

\begin{proof}
Let $X$ be a finite set, and $b\in \calB[X]$. 
\begin{align*}
\epsilon\circ \Gamma^{\exch_{\pp}}_+\circ \Gamma^\calB_{\pp}(b)&=\epsilon_+\circ \Gamma^\calB_{\pp}(b)
=\begin{cases}
1\mbox{ if }\calC_{\pp}(b)=\comp_{|X|}[X],\\
0\mbox{ otherwise}.
\end{cases}\end{align*}
If $\calC_{\pp}(b)=\comp_{|X|}[X]$, then for any $\pi\in \comp[X]$, there exists $\pi'\in \comp_{|X|}[X]\subseteq \calC_-(b)$, such that $\pi\leq\pi'$. As $\calC_-(b)\in \exch_-[X]$, $\pi\in \calC_-(b)$, so $\calC_-(b)=\comp[X]$.
Conversely, if $\calC_-(b)=\comp[X]$, then $\calC_{\pp}(b)=\comp[X]\cap \comp_{|X|}[X]=\comp_{|X|}[X]$.
Therefore, $\epsilon\circ \Gamma^{\exch_{\pp}}_+\circ \Gamma^\calB_{\pp}=\epsilon'^\calB_+$.
By unicity in Theorem \ref{theo3.2.13}, $ \Gamma^{\exch_{\pp}}_+\circ \Gamma^\calB_{\pp}=\Gamma^\calB_{\epsilon'^\calB_+}$.

On the other hand,
\begin{align*}
\epsilon\circ \Gamma^{\exch_{\pp}}_-\circ \Gamma^\calB_{\pp}(b)&=\epsilon_-\circ \Gamma^\calB_{\pp}(b)=\begin{cases}
1\mbox{ if }\calC_{\pp}(b)\neq \emptyset,\\
0\mbox{ otherwise}.
\end{cases}\end{align*}
If $\calC_{\pp}(b)\neq \emptyset$, as $\calC_{\pp}(b)\subseteq \calC_+(b)$, $\calC_+(b)\neq \emptyset$.
Conversely, if $\calC_+(b)\neq \emptyset$, as it belongs to $\exch_+[X]$, it contains maximal elements of $\comp[X]$, that is to say elements of $\comp_{|X|}[X]$. Therefore, $\calC_{\pp}(b)\neq \emptyset$.
Consequently, $\epsilon\circ \Gamma^{\exch_{\pp}}_-\circ \Gamma^\calB_{\pp}=\epsilon'^\calB_-$.
By unicity in Theorem \ref{theo3.2.13}, $ \Gamma^{\exch_{\pp}}_-\circ \Gamma^\calB_{\pp}=\Gamma^\calB_{\epsilon'^\calB_-}$.
\end{proof}

\begin{remark}
If $\calB$ is non-degenerate, then $\epsilon'^\calB_-=\epsilon^\calB_-$, so $\Gamma^{\exch_{\pp}}_-\circ \Gamma^\calB_{\pp}=\Gamma_-^\calB$.
\end{remark}

\begin{example}
In $\calG$, as it is non-degenerate, $\Gamma^{\exch_{\pp}}_+\circ \Gamma^\calG_{\pp}=\Gamma^\calG_-$.
Moreover, for any finite set $X$, for any $G\in \calG[X]$, $\calC_-(G)=\comp[X]$, so $\epsilon'^\calG_+(G)=1$,
and finally $\epsilon'^\calG_+=\epsilon^\calG_-$ and $\Gamma^{\exch_{\pp}}_+\circ \Gamma^\calG_{\pp}=\Gamma^\calG_-$. Similar results hold for $\calH\calG$.
\end{example}

\begin{example}
In $\calP$, as it is non-degenerate, $\Gamma^{\exch_{\pp}}_-\circ \Gamma^\calP_{\pp}=\Gamma_-^\calP$.
for any finite set $X$, and any $b\in \calP[X]$, $\calC_-(b)=\comp[X]$ if, and only if, $b\in \bfD\calP[X]$. Therefore, $\epsilon'^\calP_+=\epsilon^\calP_+$, and $\Gamma^{\exch_{\pp}}_+\circ \Gamma^\calP_{\pp}=\Gamma^\calP_+$. 
\end{example}

\begin{example}
In $\calQ\calP$, for any finite set $X$, and any $b\in \calQ\calP[X]$, $\calC_-(b)=\comp[X]$ if, and only if, $b\in \bfD\calQ\calP[X]$. Therefore, $\epsilon'^{\calQ\calP}_+=\epsilon^{\calQ\calP}_+$, and $\Gamma^{\exch_{\pp}}_+\circ \Gamma^{\calQ\calP}_{\pp}=\Gamma^{\calQ\calP}_+$. Moreover, $\calC_+(b)\neq \emptyset$, if, and only if, $P\in \calP[X]$. Therefore,
\[\Gamma^{\exch_{\pp}}_-\circ \Gamma^{\calQ\calP}_{\pp}(P)=\begin{cases}
\Gamma^\calP_-(P)\mbox{ if }P\in \calP[X],\\
0\mbox{ otherwise}.
\end{cases}\]
Similar results for $\calO\calG$, the condition "$P\in \calP[X]$" being replaced by "$G$ is acyclic".
\end{example}

\subsection{A link with the Malvenuto-Reutenauer Hopf algebra}

For any finite set $X$, we denote by $\calT\calO[X]$ the set of total orders on $X$. This defines a species $\calT\calO$. Let $(X,Y)$ be a pair of disjoint finite sets. For $\leq\in \calT\calO[X\sqcup Y]$, we put
\[\Delta_{X,Y}(\leq)=\begin{cases}
\leq_{\mid X}\otimes \leq_{\mid Y}\mbox{ if }\forall x\in X,\: \forall y\in Y,\: x\leq y,\\
0\mbox{ otherwise}.
\end{cases}\]
For $\leq_X\in \calT\calO[X]$ and $\leq_Y\in \calT\calO[Y]$, we put
\[\leq_X\shuffle \leq_Y=\sum_{\substack{\leq \in \calT\calO[X\sqcup Y],\\ \leq_{\mid X}=\leq_X,\\ \leq_{\mid_Y}=\leq_Y}}\leq.\]
Then $\ell\calT\calO$ is a twisted bialgebra. Note that $\calT\calO$ is not a CTB, as the product of two elements of $\calT\calO[X]$ and $\calT\calO[Y]$ is generally not an element of $\calT\calO|X\sqcup Y]$, but a sum of different elements of $\calT\calO[X\sqcup Y]$.

\begin{prop}\label{propMalvenutoReutenauer}
Let $X$ be a finite set, and $\pi \in \comp_{|X|}[X]$. Then $\pi$ induces an order $\leq_\pi$ on $X$, defined as follows: for $x,y\in X$, $x\leq_\pi y$ if $\{x\}$ appears on the left of $\{y\}$ in $\pi$.
The following map defines a morphism of twisted bialgebras from $(\ell\exch_{\pp},\shuffle,\Delta)$ to $(\ell\calT\calO,\shuffle,\Delta)$: for any finite set $X$,
\begin{align*}
&\forall I\in \exch_{\pp}[X],&\Theta(I)&=\sum_{\pi\in I}\leq_\pi.
\end{align*}\end{prop}

\begin{proof}
Let $(X,Y)$ be a pair of disjoint finite sets, $I\in \exch_{\pp}[X]$ and $J\in \exch_{\pp}[Y]$.
\begin{align*}
\Theta(I\shuffle J)&=\theta\left(\{\pi\in \comp_{|X|+|Y|}[X\sqcup Y]\mid \pi_{\mid X}\in I,\: \pi_{\mid Y}\in J\}\right)\\
&=\sum_{\substack{\pi \in \comp_{|X|+|Y|}[X\sqcup Y],\\ \pi_{\mid X}\in I,\\ \pi_{\mid Y}\in J}}\leq_\pi\\
&=\sum_{\substack{\pi\in I,\\ \pi'\in I}}\sum_{\substack{\leq \in \calT\calO[X\sqcup Y],\\ \leq_{\mid X}=\leq_\pi,\\ \leq_{\mid Y}=\leq_{\pi'}}}\leq\\
&=\sum_{\substack{\pi\in I,\\ \pi'\in I}} \leq_\pi \shuffle \leq_{\pi'}\\
&=\Theta(I)\shuffle \Theta(J).
\end{align*}
So $\Theta$ is a morphism of twisted algebras. Let now $I\in \exch_{\pp}[X\sqcup Y]$.
\begin{align*}
(\Theta\otimes \Theta)\circ \Delta_{X,Y}(I)&=\sum_{\substack{\pi\in \comp_{|X|}[X],\\ \pi'\in \comp_{|Y|}[Y],\\ (\pi,\pi')\in I}}\leq_{\pi}\otimes \leq_{\pi'}\\
&=\sum_{\substack{\pi \in I,\\ \forall (x,y)\in X\times Y,\: x\leq_\pi y}}{\leq_\pi}_{\mid X}\otimes {\leq_\pi}_{\mid Y}\\
&=\Delta_{X,Y}\circ \Theta(I).
\end{align*}
We used the exchange condition for the first equality. Therefore, $\Theta$ is a morphism of twisted coalgebras.
\end{proof}

We now apply the Fock functors of \cite{Aguiar2010}. The bosonic Fock functor sends $\calF[\calT\calO]$ to $\K[T]$, trough the isomorphism sending any $\overline{\leq}$ to $\dfrac{T^n}{n!}$, for $\leq \in \calT\calO[X]$ and $|X|=n$.
Moreover, $\calF[\Theta]$ is $\phi_{\pp}^{\exch_{\pp}}$. Let us now describe $\calF_F[\calT\calO]$. For any $n\geq 0$, we identify $\calT\calO[n]$ with $\sym_n$, through the bijection sending $\sigma\in \sym_n$ to the order $\leq_\sigma\in \calT\calO[n]$ defined by
\begin{align*}
&\forall i,j\in [n],&i\leq_\sigma j\Longleftrightarrow \sigma^{-1}(i)\leq \sigma^{-1}(j).
\end{align*}
In other words, $\sigma(1)\leq_\sigma\cdots\leq_\sigma \sigma(n)$. Therefore, we obtain the identification
\[\calF_F[\calT\calO]=\bigoplus_{n=0}^\infty \ell\calT\calO[n]=\bigoplus_{n=0}^\infty\K\sym_n.\]
The product and coproduct of $\calF_F[\calT\calO]$ are given as follows.
\begin{itemize}
\item For $\sigma\in \sym_k$ and $\tau\in \sym_l$, 
\[\sigma\shuffle \tau=\sum_{\alpha\in \sh(k,l)} \alpha^{-1}\circ (\sigma\otimes\tau).\]
For example, representing the mutation $\sigma\in \sym_n$ by the word $\sigma(1)\cdots\sigma(n))$,
\begin{align*}
(12)\shuffle(12)&=(1234)+(1324)+(3124)+(1342)+(3142)+(3412),\\
(21)\shuffle(12)&=(2134)+(2314)+(321+(2341)+(3241)+(3421),\\
(12)\shuffle(21)&=(1243)+(1423)+(4123)+(1432)+(4132)+(4312),\\
(21)\shuffle(21)&=(2143)+(2413)+(4213)+(2431)+(4231)+(4321).
\end{align*}
\item For $\sigma\in \sym_n$,
\[\Delta(\sigma)=\sum_{k=0}^n \mathrm{std}(\sigma(1)\cdots\sigma(k))\otimes \mathrm{std}(\sigma(k+1)\cdots \sigma(n)),\]
where $\mathrm{std}$ is the standardization. For example,
\begin{align*}
\Delta((521634))&=(521634)\otimes 1+\mathrm{std}(52163)\otimes \mathrm{std}(4)+\mathrm{std}(5216)\otimes \mathrm{std}(34)+\mathrm{std}(521)\otimes \mathrm{std}(634)\\
&+\mathrm{std}(52)\otimes \mathrm{std}(1634)+\mathrm{std}(5)\otimes \mathrm{std}(21634)+1\otimes (521634)\\
&=(521634)\otimes 1+(42353)\otimes (1)+(3214)\otimes (12)+(321)\otimes(312)\\
&+(21)\otimes(1423)+(1)\otimes(21534)+1\otimes (521634).
\end{align*}\end{itemize}
We recognize the Malvenuto-Reutenauer Hopf algebra $\FQSym$ \cite{Malvenuto1995,Duchamp2002}. Therefore, we obtain a morphism Hopf algebra $\calF_F[\Theta]:\calF_F[\ell\exch_{\pp}]\longrightarrow \FQSym$, such that for any $n\geq 0$,
\begin{align*}
&\forall I\in \exch_{\pp}[n],&\calF_F[\Theta](I)&=\sum_{(\{\sigma(1)\},\ldots,\{\sigma(n)\})\in I} \sigma.
\end{align*}

\section{Characters and polynomial invariants}

\subsection{An example for $\exch_{\pp}$}

\begin{prop}\label{prop3.5.1}
Let $\calB$ be a CTB. The following defines a polynomial invariant of $\calB$: for any finite set $X$, 
\begin{align*}
&\forall b\in \calB[X],&\phi_{\pp}^\calB(b)&=\frac{|\calC_{\pp}(b)|}{|X|!} T^{|X|}.
\end{align*}
The associated character (Theorem \ref{theo1.2.23}) is denoted by $f_{\pp}^\calB$: for any finite set $X$,
\begin{align*}
&\forall b\in \calB[X],&f^\calB_{\pp}(b)&=\frac{|\calC_{\pp}(b)|}{|X|!}.
\end{align*}\end{prop}

\begin{proof}
We shall use Proposition \ref{prop1.2.25}, with the following infinitesimal character of $\calB$: for any finite set $X$,
\begin{align*}
&\forall b\in \calB[X],&g^\calB_{\pp}(b)&=\begin{cases}
1\mbox{ if }|X|=1,\\
0\mbox{ otherwise}.
\end{cases}\end{align*}
The associated polynomial invariant of $\calB$ (Proposition \ref{prop1.2.25}) is denoted by $\phi_{\pp}^\calB$. For any finite set $X$, and any $b\in \calB[X]$,
\begin{align*}
\phi_{\pp}^\calB(x)&=\sum_{k=1}^{|X|} g_{\pp}^\calB(x)\circ \Delta_{X_1,\ldots,X_k}(b) \dfrac{T^k}{k!}\\
&=\sum_{k=1}^{|X|} \sum_{(X_1,\ldots,X_k)\in \calC_-(I)\cap \comp_k[X]} g_{\pp}^\calB(x)\circ \Delta_{X_1,\ldots,X_k}(b) \dfrac{T^k}{k!}\\
&=\sum_{(X_1,\ldots,X_k)\in \calC_-(I)\cap \comp_{|X|}[X]}\dfrac{T^{|X|}}{|X|!}+0\\
&=\frac{\calC_{\pp}(b)}{|X|!} T^{|X|}. \qedhere
\end{align*}\end{proof}

\begin{remark}
By Proposition \ref{prop1.2.25}, $f^\calB_{\pp}=\exp\left(g^\calB_{\pp}\right)$.
\end{remark}

\begin{notation}
For $\calB=\exch_{\pp}$, we shall simply write $\Phi_{\pp}$ for $\phi^{\exch_{\pp}}$, and $\Lambda_{\pp}$ for $f_{\pp}^{\exch_{\pp}}$. For any finite set $X$, and any $I\in \exch_{\pp}[X]$, as $\calC_{\pp}(I)=I$,
\begin{align*}
\Phi_{\pp}(I)&=\dfrac{|I|}{|X|!}T^{|X|},&\Lambda_{\pp}(I)&=\dfrac{|I|}{|X|!}.
\end{align*}
\end{notation}

We immediately obtain:

\begin{prop}
Let $\calB$ be a CTB. Then $\phi_{\pp}^\calB=\Phi_{\pp}\circ \Gamma^\calB_{\pp}$ and $f_{\pp}^\calB=\Lambda_{\pp}\circ \Gamma^\calB_{\pp}$.
\end{prop}

\subsection{Polynomial invariants of a CTB}

Let $\calB$ be a CTB, and let $f$ be a character of $\calB$.
We shall denote by $\phi^\calB_f$ the polynomial invariant associated to $f$ in Theorem \ref{theo1.2.23}. if $X$ is a nonempty finite set,
\begin{align*}
&\forall b\in \calB[X],&\phi^\calB_f(b)=\sum_{k=1}^{|X|} \sum_{(X_1,\ldots,X_k)\in \calC_-(b)} f^{\otimes k}\circ \Delta_{X_1,\ldots,X_k}(b)\binom{T}{k}.
\end{align*}

\begin{prop}\label{prop3.5.3}
Let us put $\Phi=\phi^\exch_\epsilon$. For any finite set $X$,
\begin{align*}
&\forall I\in \exch[X],&\Phi(I)&=\sum_{k=0}^{|X|} |I\cap \comp_k[X]|\binom{T}{k}.
\end{align*}\end{prop}\begin{proof}
By definition, 
\[\Phi(I)=\sum_{k=0}^{|X|} a_k(I)\binom{T}{k},\]
where $a_k(I)$ is the number of compositions $(X_1,\ldots,X_k)\in \comp_k[X]$ such that $\epsilon^{\otimes k}\circ \Delta_{X_1,\ldots,X_k}(I)\neq 0$.
If $(X_1,\ldots,X_k)$ is such a composition, putting 
\[\Delta_{X_1,\ldots,X_k}(I)=I_{\mid_1 X_1}\otimes \cdots \otimes I_{\mid_k X_k},\]
then for any $i\in [k]$, $(X_i)\in I_{\mid_i X_i}$. By definition of $\Delta$, there exist $\pi_{i,j}\in \comp[X_j]$ for $j\in [k]$,
such that $(\pi_{i,1},\ldots,\pi_{i,k})\in I$ and $\pi_{i,i}=(X_i)$. By the exchange condition, $(X_1,\ldots,X_k)\in I$.
Conversely, if $(X_1,\ldots,X_k)\in I$, Then $(X_1,\ldots,X_k)\in \calC_-(I)$ and $(X_i)\in I_{\mid_i X_i}$ for any $i\in [k]$,
so $\epsilon^{\otimes k}\circ \Delta_{X_1,\ldots,X_k}(I)\neq 0$. Therefore, $a_k(I)=|I\cap \comp_k[X]|$. 
\end{proof}

\begin{prop}\label{prop3.5.4}
Let $\calB$ be a CTB and $f$ be a combinatorial character of $\calB$. 
Then $\phi^\calB_f=\Phi\circ \Gamma^\calB_f$.
\end{prop}

\begin{proof}
As $\Phi$ is a polynomial invariant of $\exch$ and $\Gamma^\calB_f$ is a morphism of CTBs,
by composition $\Phi\circ\Gamma^\calB_f$ is a polynomial invariant of $\calB$. Moreover, for any finite set $X$, and any $b\in \calB[X]$,
\[\Phi\circ\Gamma^\calB_f(b)(1)=\epsilon\circ \Gamma^\calB_f(b)=f(b).\]
So $\Phi\circ \Gamma^\calB_f=\phi^\calB_f$.
\end{proof}

\begin{notation}
We shall consider in particular
\begin{align*}
\phi^\calB_+&=\phi_{\epsilon^\calB_+},&\phi^\calB_-&=\phi_{\epsilon^\calB_-}.
\end{align*}
When $\calB=\exch_\pm$, we shall write $\Phi_{\pm}$ instead of $\phi_\pm^{\exch_\pm}$. Then
\begin{align*}
\Phi_+&=\Phi\circ \Gamma_+^{\exch},&\Phi_-&=\Phi\circ \Gamma_-^{\exch}.
\end{align*}
Therefore, for any non-empty finite set $X$, and any $I\in \exch[X]$,
\begin{align*}
\Phi_\pm(I)&=\sum_{k=1}^{|X|}|\calC_\pm(I)\cap \comp_k[X]|\binom{T}{k}.
\end{align*}
In particular, if $I\in \exch_{\pm}[X]$,
\[\Phi_\pm(I)=\Phi(I)=\sum_{k=1}^{|X|}|I\cap \comp_k[X]|\binom{T}{k}.\]
\end{notation}

\begin{cor}
Let $X$ be a non-empty finite set and $I\in \exch[X]$. Then
\begin{align*}
\epsilon^{*(-1)}(I)&=\sum_{k=1}^{|X|}(-1)^k |I\cap \comp_k[X]|,&
\epsilon_\pm^{*(-1)}(I)&=\sum_{k=1}^{|X|}(-1)^k |\calC_\pm(I)\cap \comp_k[X]|.
\end{align*}
\end{cor}

\begin{proof}
We apply Proposition \ref{prop4.1.23}, with $f=\epsilon$ or $\epsilon_+$ and $k=-1$.
As $\displaystyle \binom{-1}{k}=(-1)^k$ for any $k\in \N$, we obtain the result.
\end{proof}

\begin{prop}\label{prop3.5.6}
Let $\calB$ be a CTB. For any finite set $X$, for any $b\in \calB[X]$, for any $N\in \N_{>0}$:
\begin{enumerate}
\item $\phi^\calB_+(b)(N)$ is the number of maps $c:X\longrightarrow [N]$ such that
\[0\neq \Delta_{c^{-1}(1),\ldots,c^{-1}(N)}(b)\in \bfD\calB[c^{-1}(1)]\otimes \cdots \otimes \bfD\calB[c^{-1}(N)].\]
\item $\phi^\calB_-(b)(N)$ is the number of maps $c:X\longrightarrow [N]$ such that
$\Delta_{c^{-1}(1),\ldots,c^{-1}(N)}(b)\neq 0$.
\end{enumerate}\end{prop}

\begin{proof}
Observe that
\[\phi^\calB_+(b)(N)=\sum_{k=0}^\infty \sum_{(X_1,\ldots,X_k)\in \calC_-(b)}\underbrace{{\epsilon^\calB_+}^{\otimes k}\circ \Delta_{X_1,\ldots,X_k}(b)}_{\in \{0,1\}}\binom{N}{k}.\]
If $k\geq 0$, we put
\[a_k=\sum_{(X_1,\ldots,X_k)\in \calC_-(b)}{\epsilon^\calB_+}^{\otimes k}\circ \Delta_{X_1,\ldots,X_k}(b).\]
By definition of $\epsilon^\calB_+$, ${\epsilon^\calB_+}^{\otimes k}\circ \Delta_{X_1,\ldots,X_k}(b)=1$, if, and only if, $\Delta_{X_1,\ldots,X_k}(b)\in \bfD\calB[X_1]\otimes\cdots\otimes \bfD[X_k]$. 
So $a_k$ is the number of surjective maps $c:X\longrightarrow [k]$ such that
\[\Delta_{c^{-1}(1),\ldots,c^{-1}(k)}(b)\in \bfD\calB[c^{-1}(1)]\otimes \cdots \otimes \bfD\calB[c^{-1}(k)].\]
This implies the result. The proof is similar for $\phi^\calB_-$. 
\end{proof}

\begin{example}
We consider the case $\calG$ of graphs (Example \ref{ex2.1.2}). 
For any finite set $X$, and any $G\in \calG[X]$, if $X_1\sqcup\cdots\sqcup X_n=X$, then 
\[\Delta_{X_1,\ldots,X_n}(G)=G_{\mid X_1}\otimes \cdots \otimes G_{\mid X_k}\neq 0.\]
Therefore, $\phi^\calG_-(G)=T^{|X|}$. As the unique element of $\bfD\calG[X]$ is the graph with no edge, for any $N\in \N$, $\phi^\calG_-(G)(N)$ is the number of maps $c:X\longrightarrow [N]$ such that if $\{x,x'\}\in E(G)$, then $c(x)\neq c(y)$. 
In other words, $\phi^\calG_-(G)$ is the chromatic polynomial of $G$ \cite{Harary1969}, see also \cite{Foissy36}.
\end{example}

\begin{example}
We consider the case $\calH\calG$ of hypergraphs (Example \ref{ex2.1.4}), with one of the four coproducts we defined. This coproduct is denoted by $\Delta$.
For any finite set $X$, and any $G\in \calH\calG[X]$, if $X_1\sqcup\cdots\sqcup X_n=X$, then 
\[\Delta_{X_1,\ldots,X_n}(G)\neq 0.\]
Therefore, $\phi^{\calH\calG}_-(G)=T^{|X|}$. The unique element of $\bfD\calH\calG[X]$ is the hypergraph on $X$ with no non-trivial hyperedge. Therefore, for any $N\in \N$:
\begin{itemize}
\item If $\Delta=\Delta^{\subset,\subset}$, $\phi^{\calH\calG^{\subset,\subset}}_+(G)(N)$ is the number of maps $c:X\longrightarrow [N]$
such that for any $e\in \HE(G)$, of cardinality $\geq 2$, $c_{\mid e}$ is not constant.
\item If $\Delta=\Delta^{\cap,\cap}$, $\phi^{\calH\calG^{\cap,\cap}}_+(G)(N)$ is the number of maps $c:X\longrightarrow [N]$
such that for any $e\in \HE(G)$, $c_{\mid e}$ is injective.
\item If $\Delta=\Delta^{\cap,\subset}$ or $\Delta^{\subset,\cap}$, $\phi^{\calH\calG^{\cap,\subset}}_+(G)(N)$ is the number of maps $c:X\longrightarrow [N]$
such that for any $e\in \HE(G)$, $c_{\mid e}$ reaches its minimum only once.
\end{itemize}
As observed in \cite{Ebrahimi-Fard2022,Foissy44}, we obtain the various chromatic polynomials of hypergraphs, including the one of \cite{Aval2019,Borowiecki2000}.
\end{example}

\begin{example}
We now consider the case $\calP$ of posets (Examples \ref{ex2.1.5}). For any finite set $X$, and any $P\in \calP[X]$, if $X_1\sqcup\cdots\sqcup X_n=X$, the following assertions are equivalent:
\begin{enumerate}
\item $\Delta_{X_1,\ldots,X_n}(G)\neq 0$.
\item For any $x,x'\in X$, assuming that $x\in X_i$ and $x'\in X_j$, then
\[x\leq_P x'\Longrightarrow i\leq j.\]
\end{enumerate}
As a consequence, for any $N\in\N$, $\phi^\calP_+(P)(N)$ is the number of maps $c:X\longrightarrow [N]$ such that
\begin{align*}
&\forall x,x'\in X,&x\leq_Px'&\Longrightarrow c(x)\leq c(x').
\end{align*}
So $\phi^\calP_-(P)$ is the Ehrhart polynomial associated to $P$. More precisely, if $P\in \calP[n]$ is a poset,
we define a polytope $\mathrm{pt}(P)$ by
\[\mathrm{pt}(P)=\{(x_1,\ldots,x_n)\in [0;1]^2\mid \forall i,j\in [n],\: i\leq_p j\Longrightarrow x_i\leq x_j\}.\]
Then $\phi^\calP_-(P)(N)$ is the number of integer points in $h_{O,N}(\mathrm{pt}(P))$, where $h_{O,N}$ is the homothety of center $O$ and ratio $N$. 

The unique element of $\bfD\calP[X]$ is the equality of $X$. Therefore, the following assertions are equivalent:
\begin{enumerate}
\item $\Delta_{X_1,\ldots,X_n}(G)\in \bfD\calP[X_1]\otimes\cdots\otimes\bfD\calP[X_k]$.
\item For any $x,x'\in X$, assuming that $x\in X_i$ and $x'\in X_j$, then
\[x\leq_P x' \mbox{ and }x\neq x'\Longrightarrow i<j.\]
\end{enumerate}
As a consequence, for any $N\in\N$, $\phi^\calP_-(P)(N)$ is the number of maps $c:X\longrightarrow [N]$ such that
\begin{align*}
&\forall x,x'\in X,&x\leq_P x' \mbox{ and }x\neq x'&\Longrightarrow c(x)< c(x').
\end{align*}
So $\phi^\calP_+(P)$ is the strict Ehrhart polynomial associated to $P$.
More precisely, $\phi^\calP_+(P)(N)$ is the number of integer points in the interior of $h_{O,N}(\mathrm{pt}(P))$.

These results on posets can be extended to quasi-posets (Example \ref{ex2.1.8}) without any problem.
In particular, if $P\in \calQ\calP[X]$ is not a poset, there exists $x,x'\in X$, such that $x\neq x'$, $x\leq_P x'$ and $x'\leq_P x$. If $c:X\longrightarrow [N]$ is a map counted by $\phi^{\calQ\calP}_-(P)(N)$, then 
\begin{align*}
c(x)&<c(x'),&c(x')&<c(x).
\end{align*} 
Therefore, such a map does not exist: $\phi^{\calQ\calP}_-(P)=0$. 
\end{example}

\begin{example}
We now consider the case $\calO\calG$ of oriented graphs (Example \ref{ex2.1.3}).
Let $X$ be a finite set, $G\in \calO\calG[X]$, and $N\in \N$. We obtain similarly that
\begin{enumerate}
\item $\phi^{\calO\calG}_-(G)(N)$ is the number of maps $c:X\longrightarrow [N]$ such that
\begin{align*}
&\forall x,x'\in X,&(x,x')\in A(G)&\Longrightarrow c(x)\leq c(x').
\end{align*}
\item $\phi^{\calO\calG}_+(G)(N)$ is the number of maps $c:X\longrightarrow [N]$ such that
\begin{align*}
&\forall x,x'\in X,&(x,x')\in A(G)&\Longrightarrow c(x)<c(x').
\end{align*}
\end{enumerate}
Note that $\phi^{\calO\calG}_{\pm}=\phi^{\calQ\calP}_{\pm}\circ \phi_{\calO\calG}$. 
This can be generalized, as we shall see in Proposition \ref{prop3.5.11}.
\end{example}

\begin{cor}\label{cor3.5.7}
Let $\calB$ be a CTB, $X$ be a non-empty finite set, and $b\in \calB[X]$. Then
\begin{align*}
{\epsilon_\pm^\calB}^{*(-1)}(b)&=\sum_{k=1}^{|X|}(-1)^k|\calC_\pm(b)\cap \comp_k[X]|.
\end{align*}\end{cor}

\begin{proof}
We use Proposition \ref{prop4.1.23}, with $f=\epsilon_\pm$. Then
\begin{align*}
{\epsilon_\pm^\calB}^{*(-1)}(b)&=\phi_\pm^\calB(-1)\\
&=\Phi_\pm\circ \Gamma_\pm^\calB(b)(-1)\\
&=\Phi\circ \Gamma_\pm^{\exch}\circ \Gamma_\pm^\calB(b)(-1)\\
&=\Phi\circ \Gamma_\pm^\calB(b)(-1)\\
&=\epsilon^{*(-1)}\left(\calC_\pm(b)\right).
\end{align*}
For the fourth equality, we used that $\Gamma_\pm^\calB$ takes its values in $\exch_\pm$, on which $\Gamma_\pm^{\exch}$ is the identity. We conclude with Proposition \ref{prop3.4.3}.
\end{proof}

\begin{prop}\label{prop3.5.8}
Let $\calB$ be a CTB, and $f$ be a character of $\calB$. Let $X$ be a finite set, of cardinality denoted by $n$, and $b\in \calB[X]$.
\begin{enumerate}
\item $\phi^\calB_f$ is of degree $\leq n$.
\item Let us assume that for any $b\in \calB[1]$, $f(b)=1$. The coefficient of $T^n$ in both $\phi^\calB_f(b)$ is $f^\calB_{\pp}(b)$. 
\end{enumerate}\end{prop}

\begin{proof}
Observe that if $k>n$, $\comp_k[X]=\emptyset$: this implies that $\phi^\calB_f(b)$ is of degree $\leq n$. 
The coefficient of $T^n$ in $\phi^\calB_+(b)$ is
\[\frac{1}{n!}\sum_{\pi\in \calC_{\pp}(b)} f^{\otimes n}\circ \Delta_\pi(b).\]
If the hypothesis of the second item is satisfied, then for any $\pi\in \calC_{\pp}(b)$, $f^{\otimes n}\circ \Delta_\pi(b)=1$, and the coefficient of $T^n$ is
\[\frac{|\calC_{\pp}(b))|}{n!}=f^\calB_{\pp}(b). \qedhere\]
\end{proof}

Proposition \ref{prop3.5.8} can be applied to $\phi_\pm$. For $\phi_+$, we obtain more:

\begin{prop}
Let $\calB$ be a CTB, $X$ be a finite set, of cardinality denoted by $n$, and $b\in \calB[X]$. Then $f^\calB_{\pp}(b)=0$ if, and only if, $\phi^\calB_+(b)=0$.
\end{prop}

\begin{proof}
By Proposition \ref{prop3.5.8}, if $\phi^\calB_+(b)=0$, then $f^\calB_{\pp}(b)=0$.
Let us assume that $f^\calB_{\pp}(b)=0$. Then $\calC_+(b)\cap \comp_{|X|}[X]=\emptyset$. If $\calC_+(b)\neq \emptyset$, as it belongs to $\exch_+[X]$, it should contain maximal elements of $\comp[X]$, that is to say elements of $\comp_n[X]$: this is a contradiction. Therefore, $\calC_+(b)=\emptyset$, so $\phi^\calB_+(b)=0$. 
\end{proof}

\begin{example}
We may have $f^\calB_{\pp}(b)=0$ and $\phi^\calB_-(b)\neq \emptyset$. Let us take a finite set $X$ with at least two elements, and $I=\{(X)\}\in \exch_-[X]$. Then $f^\calB_{\pp}(b)=0$, as $|X|\geq 2$, and $\Phi_-(I)=T$. 
\end{example}

\begin{prop}\label{prop3.5.10}
Let $\calB$ be a CTB. Let $X$ be a finite set and $b\in \calB[X]$. The following assertions are equivalents:
\begin{enumerate}
\item $\phi^\calB_+(b)=\phi^\calB_-(b)$.
\item $\phi^\calB_+(b)=T^{|X|}$.
\item $b\in \bfD\calB[X]$.
\end{enumerate}\end{prop}

\begin{proof}
For any $k\in [|X|]$, we put $a^+_k=\calC_+(b)\cap \comp_k[X]$ and $a^-_k=\calC_-(b)\cap \comp_k[X]$. Then
\begin{align*}
\phi^\calB_+(b)&=\sum_{k=0}^{|X|}a^+_k \binom{T}{k},&\phi^\calB_-(b)&=\sum_{k=0}^{|X|}a^-_k \binom{T}{k}.
\end{align*}

$1.\Longrightarrow 3$. Then, if $0\leq k\leq |X|$, $a_k^+=a_k^-$. As $\calC_+(b)\subseteq \calC_-(b)$, we deduce that $\calC_+(b)=\calC_-(b)$. As $\calC_-(b)\in \exch_-[X]$, $(X)\in \calC_-(b)=\calC_+(b)$. Consequently, $\Delta_X(b)=b\in \bfD\calB[X]$.

$3.\Longrightarrow 2$. Then, we can write $X=\{x_1,\ldots,x_n\}$, and $b=b_1\cdots b_n$, with $b_i\in \calB[\{x_i\}]$. Therefore,
\[\phi^\calB_+(b)=\prod_{i=1}^n \phi^\calB_+(b_i)=T^n=T^{|X|}.\]

$2\Longrightarrow 1$. As 
\[T^{|X|}=\sum_{k=0}^{|X|} |\comp_k[X]| \binom{T}{k},\]
we deduce that if $0\leq k\leq |X|$, $a_k^+=|\comp_k[X]|$. Therefore, $\calC_+(b)=\comp[X]$. Consequently, as $\calC_+(b)\subseteq \calC_-(b)$, $\calC_-(b)=\calC_+(b)=\comp[X]$. This implies that $\phi^\calB_-(b)=T^{|X|}$. 
\end{proof}

\begin{prop}\label{prop3.5.11}
Let $\calB,\calB'$ be two CTBs, and $\psi:\calB\longrightarrow \calB'$ be a morphism of CTBs.
\begin{enumerate}
\item The following conditions are equivalent:
\begin{enumerate}
\item For any finite set $X$, and any $b\in \calB[X]$, $\psi(b)\neq 0$.
\item $\epsilon^{\calB'}_-\circ\psi=\epsilon^\calB_-$.
\item $\phi^{\calB'}_-\circ \psi=\phi^\calB_-$.
\end{enumerate}
\item The following conditions are equivalent:
\begin{enumerate}
\item For any finite set $X$, for any $b\in \calB[X]$, $b\in \bfD\calB[X]$ if, and only if, $\psi(b)\in \bfD\calB'[X]$.
\item $\epsilon^{\calB'}_+\circ\psi=\epsilon^\calB_+$.
\item $\phi^{\calB'}_+\circ \psi=\phi^\calB_+$.
\end{enumerate}
\item Let us assume that for any finite set $X$, and any $b\in \calB[X]$, $\psi(b)\neq 0$.Then $\phi^{\calB'}_{\pp}\circ \psi=\phi^\calB_{\pp}$.
\end{enumerate}\end{prop}

\begin{proof}
1. $(a)\Longrightarrow (b)$. For any finite set $X$, and any $b\in \calB[X]$, as $\psi(b)\in \calB'[X]$ by hypothesis, $\epsilon^{\calB'}_-\circ \psi(b)=1=\epsilon^\calB_-(b)$.

2. $(a)\Longrightarrow (b)$. Let $X$ be a finite set, and $b\in \calB[X]$. By hypothesis, 
\begin{align*}
\epsilon^\calB_+(b)=1&\Longleftrightarrow b\in \bfD\calB[X]\Longleftrightarrow \psi(b)\in \bfD\calB'[X]\Longleftrightarrow \epsilon^{\calB'}_+\circ \psi(b)=1.
\end{align*}
As both $\epsilon^\calB_+$ and $\epsilon^{\calB'}_+\circ\psi$ are combinatorial characters on $\calB$, they are equal.

1. and 2. $(b)\Longrightarrow (c)$. For any finite set $X$, and any $b'\in \calB'[X]$, $\phi^\calB_\pm(b')(1)=\epsilon_\pm(b')$. Therefore, or any finite set $X$, and any $b\in \calB[X]$,
\[\phi^{\calB'}_\pm\circ \psi(b)(1)=\epsilon^{\calB'}_\pm(\psi(b))=\epsilon^\calB_\pm(b).\]
Therefore, by Theorem \ref{theo1.2.23} (unicity), $\phi^{\calB'}_\pm\circ\psi=\phi^\calB_\pm$.

1. $(c)\Longrightarrow (a)$. Let $X$ be a finite set, and $b\in \calB[X]$. Then $\phi^\calB_-(b)\neq 0$, so $\phi^{\calB'}_-\circ \psi(b)\neq 0$, which implies that $\psi(b)\neq 0$.

2. $(c)\Longrightarrow (a)$. Let $X$ be a finite set, and $b\in \calB[X]$. The cardinality of $X$ is denoted by $n$.
By Proposition \ref{prop3.5.10},
\begin{align*}
b\in \bfD\calB[X]&\Longleftrightarrow \phi^\calB_+(b)=T^n\Longleftrightarrow \phi^{\calB'}_+\circ \psi(b)=T^n\Longleftrightarrow \psi(b)\in \bfD\calB'[X]. 
\end{align*}

3. By the first item, $\phi^{\calB'}_+\circ\psi=\phi^\calB_+$. By Proposition \ref{prop3.5.8}-1, if $b\in \calB[X]$,
the leading term of $\phi^{\calB'}_+\circ \psi(b)=\phi^\calB_+(b)$ is $f^\calB_{\pp}(b)=f^{\calB'}_{\pp}\circ \psi(b)$.
By Proposition \ref{prop3.5.8}-2, 
\[\phi^{\calB'}_{\pp}\circ \psi(b)=f_{\calB'}\circ\psi(B)T^{|X|}=f^\calB_{\pp}(b)T^{|X|}=\phi^\calB_{\pp}(b). \qedhere\]
\end{proof}

For the third item, the converse is false, as we shall see in the next Example.

\begin{example}
We consider the morphism $\phi_{\calQ\calP}:\calQ\calP\longrightarrow \calP$ of Example \ref{ex2.1.14}.
If $P$ is a poset, then $f_{\pp}^{\calP\calQ}(P)=f_{\pp}^\calP(P)$. If $P$ is not a poset, then $f_{\pp}^{\calP\calQ}(P)=0$, so
\[\phi_{\pp}^{\calQ\calP}(P)=0=\phi^\calP_{\pp}\circ \phi_{\calQ\calP}(P).\]
Therefore, $\phi_{\pp}^{\calQ\calP}=\phi^\calP_{\pp}\circ \phi_{\calQ\calP}$.
\end{example}

\begin{prop}\label{prop3.5.12}
Let $\calB$ be a cocommutative CTB, and $f$ a character of $\calB$. 
Let us assume that for any finite set $X$, for any $b\in \calB[X]$, $f(b)\in \Z$. Then for any finite set $X$, for any $b\in \calB[X]$, $\phi_f(b)\in \Z[T]$.
\end{prop}

\begin{proof}
We have
\[\phi_f(b)=\sum_{k=0}^\infty\sum_{(X_1,\ldots,X_k)\in \calC_-(b)}f^{\otimes k}\circ \Delta_{X_1,\dots,X_k}(b)\binom{T}{k}.\]
As $\Delta$ is cocommutative, if $(X_1,\ldots,X_k)\in \calC_-(b)$, then for any $\sigma\in \sym_k$, $(X_{\sigma(1)},\ldots,X_{\sigma(k)})\in \calC_-(b)$: $\sym_k$ acts freely on the elements of $\calC_-(b)$ of length $k$.
Choosing a transversal $A_k$ of the orbits, we obtain
\begin{align*}
\phi_f(b)&=\sum_{k=0}^\infty\sum_{(X_1,\ldots,X_k)\in A_k}k!f^{\otimes k}\circ \Delta_{X_1,\dots,X_k}(b)\binom{T}{k}\\
&=\sum_{k=0}^\infty\sum_{(X_1,\ldots,X_k)\in A_k}\underbrace{f^{\otimes k}\circ \Delta_{X_1,\dots,X_k}(b)}_{\in \Z}T(T-1)\cdots (T-k+1)\in \Z[T]. \qedhere
\end{align*}\end{proof}

\begin{prop}\label{prop3.5.13}
Let $\calB$ be a CTB. Then $\phi^\calB_\pm=\Phi_\pm\circ \Gamma^\calB_\pm$.
\end{prop}

\begin{proof}
For any finite set $X$, for any $b\in \calB[X]$, $\calC_-(b)\neq \emptyset$, as it contains $(X)$. 
So $\Gamma^\calB_-(b)\neq 0$.
By Proposition \ref{prop3.5.11}-1, $\Phi_-\circ \Gamma^\calB_-=\phi^\calB_-$.
Let $X$ be a finite set, and $b\in \calB[X]$.
\begin{align*}
b\in \bfD\calB[X]&\Longleftrightarrow (X)\in \calC_+(b)\Longleftrightarrow \calC_+(b)\in \bfD\exch_+[X]
\Longleftrightarrow \Gamma^\calB_+(b)\in \bfD\exch_+[X].
\end{align*}
By Proposition \ref{prop3.5.11}, $\Phi_+\circ \Gamma^\calB_+=\phi^\calB_+$.
\end{proof}

\begin{remark}
By construction, we also have $\phi^\calB_{\pp}=\Phi_{\pp}\circ \Gamma^\calB_+$.
\end{remark}

\begin{remark}\label{rk3.5.3}
As a consequence, for a given integer $n$, for any finite set $X$ of cardinality $n$, and any $b\in \calB[X]$, $\phi^\calB_\pm(b)=\Phi_\pm(\calC_\pm(b))$. 
This implies that the set of non-zero polynomials $\calP_\pm(n)$ which are obtained in this way from an object of a given degree $n$ of any CTB is equal to $\Phi_\pm(\exch_\pm[n])$. For example, 
\begin{align*}
\calP_+(1)&=\{T\},\\
\calP_+(2)&=\left\{T^2,\:(T-1)T,\:\dfrac{(T-1)T}{2},\right\},\\
\calP_+(3)&=\left\{\begin{array}{l}
T^3,\:(T-1)T^2,\:(T-1)^2T,\:(T+1)(T-1)T,\:(T-2)(T-1)T,\:\\[3mm]
\dfrac{(T-1)^2T}{2},\:\dfrac{(T-1)T^2}{2},\:\dfrac{(T-2)(T-1)T}{2},\\[3mm]
\:\dfrac{(2T-1)(T-1)T}{2},\:\dfrac{(2T-3)(T-1)T}{2},\:\dfrac{(2T+1)(T-1)T}{2},\:\\[3mm]
\dfrac{(T-2)(T-1)T}{3},\:\dfrac{2(T-2)(T-1)T}{3},\:\dfrac{(2T-1)(T-1)T}{3},\:\\[3mm]
\dfrac{(T-2)(T-1)T}{6},\:\dfrac{(2T-1)(T-1)T}{6},\:\dfrac{(4T+1)(T-1)T}{6},\:\\[3mm]
\dfrac{(4T-5)(T-1)T}{6},\:\dfrac{5(T-2)(T-1)T}{6},\:\dfrac{(5T-1)(T-1)T}{6},\:\\[3mm]
\dfrac{(5T+2)(T-1)T}{6},\:\dfrac{(5T-7)(T-1)T}{6},\:\dfrac{(5T-4)(T-1)T}{6}
\end{array}\right\}.
\end{align*}
Among them, here are the chromatic polynomials of graphs of degree 3:
\[\calP_{\mathrm{chr}}(3)=\left\{T^3,\:(T-1)T^2,\:(T-1)^2T,\:(T-2)(T-1)T\right\}.\]
It is also set of chromatic polynomials of hypergraphs of degree 3 attached to the coproduct $\Delta^{\cap,\cap}$. The set of chromatic polynomials of hypergraphs of degree 3 attached to the coproduct $\Delta^{\subset,\subset}$ is
\[\calP_{\mathrm{chr}}(3)\sqcup\{(T+1)(T-1)T\}.\]
The set of chromatic polynomials of hypergraphs of degree 3 attached to the coproduct $\Delta^{\subset,\cap}$ or $\Delta^{\cap,\subset}$ is
\[\calP_{\mathrm{chr}}(3)\sqcup\left\{\dfrac{(2T-1)(T-1)T}{2},\:\dfrac{(2T-3)(T-1)T}{2}\right\}.\]
The set of strict Ehrhart polynomials of degree 3 attached to a poset or a quasi-poset is
\[\calP_{\mathrm{Ehr}}(3)=\left\{T^3,\:\dfrac{(T-1)T^2}{2},\:\dfrac{(T-2)(T-1)T}{6},\:\dfrac{(2T-1)(T-1)T}{6}\right\}.\]
The set of polynomials attached to a mixed graph on three vertices is 
\[\calP_{\mathrm{chr}}(3)\cup \calP_{\mathrm{Ehr}}(3)\cup\left\{\frac{(T-1)^2T}{2},\:\frac{(T-2)(T-1)T}{2},\:\frac{(T-2)(T-1)T}{3}\right\}.\]
The set of polynomials obtained from a noncrossing partition on three blocks is
\[\left\{\frac{(2T-3)(T-1)T}{2},\:\frac{(T-1)^2T}{2},\: \frac{(T-2)(T-1)T}{3},\: \frac{(2T-1)(T-1)T}{6},\: \frac{(T-2)(T-1)T}{6}\right\}.\]
For $\calP_-$, we obtain
\begin{align*}
\calP_-(1)&=\{T\},\\
\calP_-(2)&=\left\{T,\:T^2,\:\dfrac{(T+1)T}{2},\right\},\\
\calP_-(3)&=\left\{\begin{array}{l}
T,\:T^2,\:(2T-1)T,\:(3T-2)T,\: \dfrac{(T+1)T}{2},\:\dfrac{(3T-1)T}{2},\:\dfrac{(5T-3)T}{2},\\[3mm]
T^3,\:(T^2+2T-1)T,\:\dfrac{(T+1)T^2}{2},\:\dfrac{(T^2+3T-2)T}{2},\\[3mm]
\dfrac{(T^2+6T-4)T}{3},\:\dfrac{(2T-1)(T+2)T}{3},\:\dfrac{(T^2+3T-1)T}{3},\\[3mm]
\dfrac{(T+2)(T+1)T}{6},\:\dfrac{(T^2+12T-7)T}{6},\:\dfrac{(T^2+15T-10)T}{6},\\[3mm]
\dfrac{(T^2+6T-1)T}{6},\:\dfrac{(T^2+9T-4)T}{6},\:\dfrac{(2T+1)(T+1)T}{6},\\[3mm]
\dfrac{(2T-1)(T+5)T}{6},\:\dfrac{(4T-1)(T+1)T}{6},\:\dfrac{(5T-2)(T+1)}{6}
\end{array}\right\}.\end{align*}\end{remark}

\begin{remark}
There is no way to obtain $\Phi_-(I)$ from $\Phi_+(I)$ uniquely. Take for example, in $\exch_+[3]$,
\begin{align*}
I&=\{(1,2,3),\: (1,3,2)\},&
J&=\{(1,2,3),\: (2,3,1)\},
\end{align*}
then $\Phi_+(I)=\Phi_+(J)=\dfrac{(T-2)(T-1)T}{3}$. However,
\begin{align*}
m(I)&=I\sqcup\{(1,23),\:(12,3),\:(13,2),\:(123)\},\\
m(J)&=I\sqcup\{(1,23),\:(2,13),\:(12,3),\:(23,1),\:(123)\},
\end{align*}
so 
\begin{align*}
\Phi_-(I)&=2\dfrac{T(T-1)(T-2)}{6}+3\dfrac{T(T-1)}{2}+T=\dfrac{(2T+1)(T+1)T}{6},\\
\neq\Phi_-(J)&=2\dfrac{T(T-1)(T-2)}{6}+4\dfrac{T(T-1)}{2}+T=\dfrac{(T^2+3T-1)T}{3}.
\end{align*}\end{remark}

\begin{prop}
For any $n\geq 1$, $\calP_-(n)\subsetneq \calP_-(n+1)$.
\end{prop}

\begin{proof}
Let $X$ be a finite set of cardinality $n$, and $I\in \exch_-[X]$. Let $X'=X\sqcup\{*\}$ be a set of cardinality $n+1$ that contains $X$. We put
\[I'=\{(X_1,\ldots,X_k\sqcup\{*\})\mid (X_1,\ldots,X_k)\in I\}.\]
Then $I'\in \exch_-[X']$ and $\Phi_-(I')=\Phi_-(I)$. So $\calP_-(n)\subseteq \calP_-(n+1)$. The inclusion is strict: the elements of $\calP_-(n)$ are all of degree $\leq n$, and $T^{n+1}=\Phi_-(\comp[n+1])\in \calP_-(n+1)$.
\end{proof}

\begin{prop}\label{prop3.5.15}
Let $\calB$ be a CTB. The following assertions are equivalent:
\begin{enumerate}
\item $\calB$ is non-degenerate.
\item For any finite set $X$, and any $b\in \calB[X]$, $f^\calB_{\pp}(b)\neq 0$.
\item For any finite set $X$, and any $b\in \calB[X]$, $\phi^\calB_+(b)$ is of degree $|X|$.
\item For any finite set $X$, and any $b\in \calB[X]$, $\phi^\calB_+(b)\neq 0$.
\item For any finite set $X$, and any $b\in \calB[X]$, $\phi^\calB_-(b)$ is of degree $|X|$.
\end{enumerate}\end{prop}

\begin{proof}
$1.\Longrightarrow 2$. Let $X$ be a finite set, and $b\in \calB[X]$. As $\calB$ is non-degenerate, $\calC_+(b)\neq \emptyset$. As $\calC_+(b)\in \exch_+[X]$, $\calC_{\pp}(b))\neq \emptyset$.
Therefore, $f^\calB_{\pp}(b)\neq 0$.

$2.\Longrightarrow 3$. Immediate consequence of Proposition \ref{prop3.5.8}.

$3.\Longrightarrow 4$. Immediate.

$4.\Longrightarrow 5$. By hypothesis, for any finite set $X$, and any $b\in \calB[X]$, $\calC_+(b)\neq \emptyset$. As $\calC_+(b)\in \exch_+[X]$, $\calC_{\pp}(b))\neq \emptyset$,
so $\calC_+(b)$ contains compositions of length $|X|$. Moreover, $\calC_+(b)\subseteq \calC_-(b)$, so $\calC_-(b)$ contains compositions of length $|X|$, which implies 5.

$5.\Longrightarrow 1$. Then, for any finite set $X$, and any $b\in \calB[X]$, $M(\calC_-(b))\neq \emptyset$. Any element of $M(\calC_-(b))$ belongs to $\calC_+(b)$, so $\calC_+(b)\neq \emptyset$. Therefore, $\calB$ is non-degenerate.
\end{proof}

\begin{lemma}\label{lem3.5.16}
Let $\calB$ be a non-degenerate CTB, $X$ a finite set, and $b\in \calB[X]$. The following assertions are equivalent:
\begin{enumerate}
\item $\phi^\calB_+(b)$ is unitary.
\item $\phi^\calB_-(b)$ is unitary.
\item $\phi^\calB_-(b)\in \Z[T]$.
\item $\phi^\calB_-(b)=T^{|X|}$.
\item $f^\calB_{\pp}(b)=1$.
\end{enumerate}\end{lemma}

\begin{proof}
As $\calB$ is non-degenerate, from Proposition \ref{prop3.5.15}, the degree of $\phi^\calB_+(b)$ and $\phi^\calB_-(b)$ is $|X|$. 
By Proposition \ref{prop3.5.8}, they share the same leading coefficient $f^\calB_{\pp}(b)\in ]0;1]$. This immediately implies that $1.\Longleftrightarrow 2.\Longleftrightarrow 5$.

$4.\Longrightarrow 3$. Obvious.

$3.\Longrightarrow 5$. Then the leading coefficient $f^\calB_{\pp}(b)$ of $\phi^\calB_-(b)$ belongs to $\Z$. Moreover, by definition, $0<f^\calB_{\pp}(b)\leq 1$, so $f^\calB_{\pp}(b)=1$.

$5.\Longrightarrow 4$. Then $\calC_-(b)\cap \comp_{|X|}[X]=\comp_{|X|}[X]$,
so $\calC_-(b)$ contains all the maximal elements of $\comp[X]$. Therefore, $\phi^\calB_-(b)=T^{|X|}$.
\end{proof}

\begin{remark}\label{rk3.5.5} 
\begin{enumerate}
\item If $\calB$ is non-degenerate and cocommutative, by Proposition \ref{prop3.5.12} for $f=\epsilon_-$, for any finite set $X$ and any $b\in \calB[X]$, the 5 assertions of Lemma \ref{lem3.5.16} are satisfied. Still by Proposition \ref{prop3.5.12}, but for $f=\epsilon^\calB_+$, $\phi^\calB_+(b)\in \Z[T]$. 
\item In general, given a CTB $\calB$, a finite set $X$ and $b\in \calB[X]$, $\phi^\calB_\pm(b)\in \Z[T]$ if, and only if, for any $k\in [|X|]$,
\[k!\mid |\calC_\pm(b)\cap \comp_k[X]|.\]
\end{enumerate}\end{remark}

\begin{example}
It may happen that $\phi^\calB_+(b)$ is unitary, but not in $\Z[T]$. Take for example
\begin{align*}
J&=\comp[4]\sqcup\{(12,3)\}\in \exch_+[3].
\end{align*}
 Then 
\[\Phi_+(J)=6\dfrac{T(T-1)(T-2)}{6}+\dfrac{T(T-1)}{2}=\dfrac{(2T-3)(T-1)T}{2},\]
which is unitary, but not in $\Z[T]$. 
\end{example}

\subsection{Forbidden values for $f^\calB_{\pp}$}

Let $\calB$ be a CTB. For any finite set $X$, of cardinality $n$, for any $b\in \calB[X]$, $n!f^\calB_{\pp}(b)\in \{0,\dots,n!\}$.
Let us show that it cannot take all possible values in $\{0,\ldots,n!\}$.

\begin{prop}\label{prop3.5.17}
Let $\calB$ be a CTB, $X$ be a finite set of cardinality $n$, $b\in \calB[X]$. We put $N=n!f^\calB_{\pp}(b)$. Then, if $k+l=n$ and $k,l\geq 2$, then the coefficient of $T^N$ in the following polynomial is non-zero:
\[P_{k,l}=\left(1+\sum_{p=1}^{k!}\sum_{q=1}^{l!} T^{pq}\right)^{\binom{k+l}{k}}.\]
\end{prop} 

\begin{proof}
Up to a bijection, we assume that $X=\{1,\ldots,n\}$, and put
\[\Lambda(b)=\{\sigma \in \sym_n\mid (\{\sigma(1)\},\ldots,\{\sigma(n)\})\in \calC_{\pp}(b)\}.\]
Then $n!f^\calB_{\pp}(b)=|\Lambda(b)|$.
Any $\sigma\in \sym_n$ can be uniquely written as
\[\sigma=\alpha\circ (\sigma_1\otimes \sigma_2),\]
with $\sigma_1\in \sym_k$, $\sigma_2\in \sym_l$, and $\alpha$ in the set $\sh(k,l)$ of $(k,l)$-shuffles. The exchange condition implies that for any $\sigma_1,\sigma'_1\in \sym_k$, $\sigma_2,\sigma'_2\in \sym_l$, and $\alpha\in \sh(k,l)$,
\begin{align*}
&(\alpha\circ (\sigma_1\otimes \sigma_2)\in \Lambda(b) \mbox{ and }(\alpha\circ (\sigma'_1\otimes \sigma'_2)\in \Lambda(b))\\
\Longleftrightarrow&(\alpha\circ (\sigma'_1\otimes \sigma_2)\in \Lambda(b) \mbox{ and }(\alpha\circ (\sigma_1\otimes \sigma'_2)\in \Lambda(b)).
\end{align*} 
Let us decompose $\Lambda(b)$ according to $\alpha$:
\[\Lambda(I)=\bigsqcup_{\alpha\in \sh(k,l)}\Lambda_\alpha(I),\]
where $\Lambda_\alpha(I)$ is the set of elements of $\Lambda(I)$ of the form $\alpha\circ (\sigma_1\otimes \sigma_2)$, with $\sigma_1\in \sym_1$ and $\sigma_2\in \sym_2$. 
The exchange condition implies that there exists $S_1\subseteq \sym_k$, $S_2\subseteq \sym_l$, such that
\[\Lambda_\alpha(b)=\{\alpha\circ (\sigma_1\otimes \sigma_2)\mid \sigma_1\in S_1,\l \sigma_2\in S_2\}.\]
Therefore, $|\Lambda_\alpha(b)|=0$ or of the form $pq$, with $p\in [k!]$ and $q\in [l!]$. As $N=\Lambda_\alpha(I)$ is a sum of $\binom{k+l}{k}$ such numbers, $T^N$ appears in $P_{k,l}$. 
\end{proof}

\begin{example}
For $n=4$, taking $(k,l)=(2,2)$, we obtain that $4!f^\calB_{\pp}(b)\neq 23$. For $n=5$, taking $(k,l)=(2,3)$, we obtain that 
\[5!f^\calB_{\pp}(b)\notin\{115,117,119\}.\]
For $n=6$, taking $(k,l)=(2,4)$ and $(3,3)$, we obtain that
\[6!f^\calB_{\pp}(b)\notin\{695, 697, 699, 701, 703, 705, 706, 707, 709, 710, 711, 712, 713, 715, 716, 717, 718, 719\}.\]
(There are 18 forbidden values). 
For $n=7$, taking $(k,l)=(2,5)$ and $(3,4)$, we obtain that
\[7!f^\calB_{\pp}(b)\notin\left\{\begin{array}{c}
4921, 4923, 4925, 4927, 4929, 4931, 4933, 4935, 4937, 4939, 4941, 4943, \\
4945, 4947, 4949, 4951, 4953, 4955, 4957, 4959, 4961, 4963, 4965, 4967, \\
4969, 4971, 4973, 4975, 4977, 4979, 4981, 4983, 4985, 4987, 4989, 4991, \\
4993, 4995, 4997, 4999, 5001, 5002, 5003, 5005, 5007, 5008, 5009, 5011, \\
5012, 5013, 5014, 5015, 5017, 5018, 5019, 5020, 5021, 5023, 5024, 5025, \\
5026, 5027, 5029, 5030, 5031, 5032, 5033, 5035, 5036, 5037, 5038, 5039
\end{array}\right\}.\]
(There are 72 forbidden values).
\end{example}

\chapter{Double-bialgebraic aspects}

\label{chap4}

We now introduce combinatorial twisted double bialgebras (briefly, CTDBs), that are both CTBs and twisted double bialgebras, following the definition of the first chapter, with a combinatorial compatibility detailed in Definition \ref{defi4.1.1}.
Examples are given by graphs (Example \ref{ex4.1.1}), where this second coproduct is given by contraction of connected subgraphs, hypergraphs (Example \ref{ex4.1.2}) with a similar construction, and posets (Example \ref{ex4.1.3}). 
A CTDB is always commutative (Proposition \ref{prop4.1.9}), but not necessarily cocommutative, as it is shown by the CTDB of posets $\calP$.
To any element $b\in \calB[X]$ of a CTDB $\calB$ is associated a set of equivalences $\calE[b] \subseteq \calE[X]$ (Notation \ref{notation4.1.1}). We show in Proposition \ref{prop4.1.7} that $\calE[b]$ has always a smaller element $\sim_b^{(m)}$ (which is always the equality of $X$)
and a greater element $\sim_b^{(M)}$, which is a refinement of the equivalence $\sim_b^{(\ic)}$, whose classes are the indecomposable components of $b$. The most interesting cases of CTDBs are when $\sim_b^{(M)}=\sim_b^{(\ic)}$ for any $b\in \calB$: they are called regular.
We give eight characterizations of regular CTDBs in Proposition \ref{prop4.1.15}; one of the most useful is that the counit of the second coproduct $\delta$ is the character $\epsilon_+^\calB$ of the third chapter.
Graphs, hypergraphs, posets, and more generally any non-degenerate CTDB is regular (Proposition \ref{prop4.1.16}). 
A non regular CTDB has, in some sense, a defect of gradation, and we explain how to correct this, and to reduce ourselves to the case of regular CTDBs.
Regular CTDBs are closely related to operads: we show in Theorem \ref{theo4.1.22} that their subspecies of indecomposable elements are given an operadic structure, inducing the second coproduct $\delta$.
If $\calB$ is a regular CTDB, then its fundamental polynomial invariant, as defined in Theorem \ref{theo1.3.6}, is $\phi_+^\calB$. 
The double structure allows to have an expression for the coefficients of this polynomial invariant with the help of a character $g^\calB_{\pp}$ (Theorem \ref{theo4.1.24}).

We introduce in the second section the weaker notion of quasi-cocommutative CTBs. The definition of these objects is rather technical (see Definition \ref{defi4.2.1}) and we won't detail it now. 
They generalize non-degenerate twisted double bialgebras, and several results on the latter can be extended to them, as we shall see. Generally, a quasi-cocommutative CTB $\calB$ does not have a second coproduct $\delta$,
but a coaction-like map $\xi$ can be defined from $\calB$ to $\exch_+\otimes \calB$. 
First of all, a set of equivalences $\calE[b]$ is associated to any element $b\in \calB$ (Definition \ref{defi4.2.5}), and an element $\calC_+(b,\sim)\in \exch_+$ is associated to any $\sim\in \calE[b]$ (Theorem \ref{theo4.2.10}). In the case of a CTDB, $\calC_+(b,\sim)=\calC_+(b/{\sim})$; 
however, when $\calB$ is only quasi-cocommutative, $b/{\sim}$ is not necessarily well-defined. 
This map $\xi$ is sufficient to find a cancellation-free formula for the antipode and the eulerian idempotent (Corollary \ref{cor4.2.13}). 
The fundamental polynomial invariant is also generalized to the level of quasi-cocommutative twisted combinatorial bialgebras (Proposition \ref{prop4.2.14}), and it is proved to be compatible with $\xi$. A formula for the coefficients of this polynomial invariant can be found in Theorem \ref{theo4.2.17}. 
These results will be used in the forthcoming chapter \ref{chap6}, on a twisted bialgebra on trees introduced in \cite{Xiao2026}. 
We also prove that any commutative twisted combinatorial bialgebra $\calB$ contains a maximal quasi-cocommutative combinatorial twisted subbialgebra $\calB^{\qc}$, and that this construction is functorial (Theorem \ref{theo4.2.18} and Proposition \ref{prop4.2.19}).

In the third section of this chapter, we turn back to $\exch_+$ and look for a second coproduct $\delta$ on it, hoping to make it a CTDB. As it is not quasi-cocommutative, the results implies that this does not exist, and we have to restrict ourselves, at least to $\exch_+^{\qc}$. 
In fact, this is not enough, and we prove in Theorem \ref{theo4.3.1} the existence of a twisted combinatorial twisted subbialgebra $\exch_+^{\double}$, such that:
\begin{itemize}
\item $\exch_+^{\double}$ can be given a unique coproduct $\delta$, making it double.
\item If $\calB$ is a twisted combinatorial twisted subbialgebra of $\exch_+$ which can be given a coproduct $\delta'$, making it a CTDB, then $\calB\subseteq \exch_+^{\double}$ and $\delta'=\delta_{\calB}$.
\end{itemize}
A consequence of this result is that if $\calB$ is a CTB such that $\Gamma_+^\calB$ is injective, then $\calB$ has at most one coproduct $\delta$ making it double (Corollary \ref{cor4.3.2}). This includes the CTDB $\calG$ of graphs and $\calP$ of posets.
Unfortunately, in general, if $\calB$ is a CTDB, then $\Gamma_+^\calB$ does not induce a morphism of CTDBs from $\calB$ to $\exch_+^{\double}$, and we give a necessary and sufficient condition for $\Gamma_+^\calB$ to be so in Theorem \ref{theo4.3.5}, when $\calB$ is non-degenerate 
-- note that in this case, the only potential morphism of CTDBs from $\calB$ to $\exch_+^{\double}$ is $\Gamma_+^\calB$. This condition is about the conservation of indecomposable components, and is always satisfied when $\Gamma_+^\calB$ is injective (Corollary \ref{cor4.3.6}). 
For the case of a quasi-cocommutative CTB $\calB$, we can ask ourselves is the map $\xi$ is a left coaction from $\exch_+^{\double}$ onto $\calB$. This is not the case in general, and we give a necessary and sufficient condition for this in Proposition \ref{prop4.3.9}.

\section{Combinatorial twisted double bialgebras}

\subsection{Definition and examples}

\begin{defi}\label{defi4.1.1}
A combinatorial twisted double bialgebra (briefly, CTDB) is a family $(\calB,m,\Delta,\delta)$ such that:
\begin{enumerate}
\item $(\calB,m,\Delta)$ is a CTB.
\item $(\ell\calB,m,\Delta,\delta)$ is a twisted double bialgebra.
\item For any finite set $X$, and $\sim\in [X]$, and any $b\in \calB[X]$, either $\delta_\sim(b)=0$ or is of the form
\[\delta_\sim(b)=b/{\sim}\otimes b\mid \sim,\]
with $b/{\sim}\in \calB[X/{\sim}]$ and $b\mid\sim\in \calB[X]$.
\end{enumerate}\end{defi}

\begin{remark}
As $\ell\calB[\emptyset]$ is 1-dimensional, generated by $1_\calB$, we have for free that for any $x\in \ell\calB[\emptyset]$,
\[\left(\varepsilon_\Delta\otimes\id_{\calB(\emptyset)}\right)\circ\delta_{\sim_\emptyset}(x)=\varepsilon_\Delta(x)1_\calB.\]
\end{remark}

If $(\calB,m,\Delta,\delta)$ is a CTDB, then $\calF[\ell\calB]$ is a double bialgebra in the classical sense, with the coproduct $\delta$ defined by
\begin{align*}
&\forall b\in \calB[n],&\delta(\overline{b})&=\sum_{\sim\in \calE[b]}\overline{b/{\sim}}\otimes \overline{b\mid\sim}.
\end{align*}

If $\calB$ is a CTDB, then $\calF[\calB]$ inherits a second coproduct $\delta$, making it a double bialgebra, in the classical meaning: for any finite set $X$,
\begin{align*}
&\forall b\in \calB[X],&\delta(\overline{b})&=\sum_{\substack{\sim\in \calE[X],\\ \delta_\sim(b)=b/{\sim}\otimes b\mid\sim\neq 0}}\overline{b/{\sim}}\otimes \overline{b\mid\sim},&\epsilon_\delta(\overline{b})&=\epsilon_\delta(b)
\end{align*}

\begin{notation}\label{notation4.1.1}
Let $\calB$ be a CTDB.
Let $X$ be a finite set, and $b\in \calB[X]$. 
We denote by $\calE[b]$ the set of equivalences $\sim\in \calE[X]$ such that $\delta_\sim(b)\neq 0$.
\end{notation}

Let us reformulate the coassociativity of $\delta$:

\begin{prop}\label{lem4.1.2}
Let $X$ be a finite set, $b\in \calB[X]$, and $\sim'\subseteq \sim\in \calE[X]$. We denote by $\overline{\sim}$ the equivalence on $X/{\sim}'$ induced by $\sim$. The following assertions are equivalent:
\begin{enumerate}
\item $\sim'\in \calE[b]$ and $\overline{\sim}\in \calE[b/{\sim}']$.
\item $\sim\in \calE[b]$ and $\sim'\in \calE[b\mid \sim]$.
\end{enumerate}
If so, then
\begin{align*}
(b/{\sim}')/\overline{\sim}&=b/{\sim},&b/{\sim}'\mid \overline{\sim}&=(b\mid \sim)/{\sim}',&(b\mid \sim)\mid \sim'&=b\mid \sim'.
\end{align*}\end{prop}

\begin{proof}
As $\delta$ is coassociative,
\begin{align*}
(\delta_{\overline{\sim}} \otimes \id_{\calB[X]})\circ \delta_{\sim'}(b)&=\begin{cases}
(b/{\sim}')/\overline{\sim}\otimes (b/{\sim}')\mid \overline{\sim}\otimes b\mid \sim'\mbox{ if }\sim'\in \calE[b]\mbox{ and }\overline{\sim}\in \calE[b/{\sim}'],\\
0\mbox{ otherwise}
\end{cases}\\
=(\id_{\calB[X/{\sim}]}\otimes \delta_{\sim'})\circ \delta_\sim(b)&=\begin{cases}
b/{\sim}\otimes (b\mid \sim)/{\sim}'\otimes (b\mid \sim)\mid \sim'\mbox{ if }\sim \calE[b]\mbox{ and }\sim'\in \calE[b/{\sim}],\\
0\mbox{ otherwise}.
\end{cases}
\end{align*}
The result follows.
\end{proof}

Let us now reformulate the compatibility between $\delta$ and the product:

\begin{prop}\label{prop4.1.3}
Let $X$, $Y$ be a disjoint finite sets, $a\in \calB[X]$, and $b\in \calB[Y]$. Then
\[\calE[ab]=\{\sim_X\sqcup \sim_Y\mid\sim_X\in \calE[a],\:\sim_Y\in \calE[b]\}.\] 
Moreover, if $\sim_X\in \calE[a]$, $\sim_Y\in \calE[b]$, and $\sim=\sim_X\sqcup\sim_Y$, then
\begin{align*}
(ab)/{\sim}&=(a/{\sim}_X)(b/{\sim}_Y),&(ab)\mid\sim=(a\mid\sim_X)(b\mid\sim_Y).
\end{align*}
\end{prop}

\begin{proof}
$\subseteq$. Let $\sim\in \calE[ab]$. We put $\sim_X=\sim\cap X^2$ and $\sim_Y=\sim\cap Y^2$. If $\sim\neq \sim_X\sqcup\sim_Y$, as $\delta$ is compatible with $m$, $\delta_\sim(ab)=0$.
So $\sim=\sim_X\sqcup \sim_Y$, and $\delta_\sim(ab)=\delta_{\sim_X}(a)\delta_{\sim_Y}(b)\neq 0$.
Consequently, $\delta_{\sim_X}(a)\neq 0$, and $\delta_{\sim_Y}(b)\neq 0$: $\sim_X\in \calE[a]$, and $\sim_Y\in \calE[b]$.

$\supseteq$. Let $\sim_X\in \calE[a]$, and $\sim_Y\in \calE[b]$. We put $\sim=\sim_X\sqcup\sim_Y$. As $\delta$ is compatible with $m$,
\[\delta_\sim(ab)=\delta_{\sim_X}(a)\delta_{\sim_Y}(b)=(a/{\sim}_X)(b/{\sim}_Y)\otimes (a\mid\sim_X)(b\mid\sim_Y).\]
As the product of two elements of $\calB$ is never zero (Proposition \ref{prop2.1.4}), this is nonzero. So $\sim\in \calE[ab]$, and we obtain the formulas for $(ab)/{\sim}$ and $(ab)\mid\sim$.
\end{proof}

Let us reformulate the compatibility between $\delta$ and $\Delta$:

\begin{prop}\label{prop4.1.4}
Let $X$ be a finite set, $b\in \calB[X]$, $\sim\in \calE[X]$, and $(\overline{X'},\overline{X''})\in \comp[X/{\sim}]$. We put $X'=\varpi_\sim^{-1}(\overline{X'})$, $X''=\varpi_\sim^{-1}(\overline{X'})$, $\sim'=\sim\cap X'^2$, and $\sim''=\sim\cap X''^2$.
The following conditions are equivalent:
\begin{enumerate}
\item $\sim\in \calE[b]$ and $(\overline{X'},\overline{X''})\in \calC_-(b/{\sim})$.
\item $(X',X'')\in \calC_-(b)$, $\sim'\in \calE[b_{\mid_LX'}]$, and $\sim''\in \calE[b_{\mid_RX''}]$.
\end{enumerate}
Moreover, if this holds,
\begin{align*}
(b/{\sim})_{\mid_L\overline{X'}}&=b_{\mid_LX'}/{\sim}',&(b/{\sim})_{\mid_R\overline{X''}}&=b_{\mid_RX''}/{\sim}'',&b\mid\sim&=(b_{\mid_LX'}\mid\sim')(b_{\mid_RX''}\mid\sim'').
\end{align*}\end{prop}

\begin{proof}
As $\Delta$ and $\delta$ are compatible, 
\begin{align*}
m_{1,3,24}\circ(\delta_{\sim'}\otimes \delta_{\sim''})\circ \Delta_{X',X''}(b)&=(\Delta_{\overline{X'},\overline{X''}}\otimes \id_{\calB[X]})\circ \delta_\sim(b).
\end{align*}
The left side is nonzero if, and only if, (2) holds. The right side is nonzero if, and only if, (1) holds. Consequently, (1) and (2) are equivalent. If they both hold,
\begin{align*}
m_{1,3,24}\circ(\delta_{\sim'}\otimes \delta_{\sim''})\circ \Delta_{X',X''}(b)&=(b_{\mid LX'}/{\sim}'\otimes b_{\mid_RX''}/{\sim}'' \otimes (b_{\mid_LX'}\mid\sim')(b_{\mid_RX''}\mid\sim''),\\
=(\Delta_{\overline{X'},\overline{X''}}\otimes \id_{\calB[X]})\circ \delta_\sim(b)&=(b/{\sim})_{\mid_L\overline{X'}}\otimes (b/{\sim})_{\mid_R\overline{X''}}\otimes b\mid\sim,
\end{align*}
which gives the three equalities by identification. \end{proof}

\begin{example}\label{ex4.1.1}
Let us consider again the CTB $\calG$ of graphs (Example \ref{ex2.1.2}). Let $X$ be a finite set, $G\in \calG[X]$, and $\sim\in \calE[X]$. We define $G/{\sim}\in \calG[X/{\sim}]$ and $G\mid\sim \in \calG[X]$ by
\begin{align*}
E(G/{\sim})&=\{\{p_\sim(x),p_\sim(y)\}\mid \{x,y\}\in E(G),\: x\nsim y\},\\
E(G\mid\sim)&=\{\{x,y\}\in E(G)\mid x\sim y\}.
\end{align*}
We shall say that $\sim\in \calE[G]$ if for any $Y\in X/{\sim}$, $G_{\mid Y}$ is connected. We put
\[\delta_\sim(G)=\begin{cases}
G/{\sim}\otimes G\mid \sim \mbox{ if }\sim\in\calE[G],\\
0\mbox{ otherwise}.
\end{cases}\]
Then $(\calG,m,\Delta,\delta)$ is a CTDB.
\end{example}

\begin{example}\label{ex4.1.2}
Let us consider the CTB $\calH\calG$ of hypergraphs (Example \ref{ex2.1.4}). Let $X$ be a finite set, $G\in \calH\calG[X]$, and $\sim\in \calE[X]$. We define $G/{\sim}\in \calH\calG[X/{\sim}]$, 
and $G\mid_\cap\sim, G\mid\sim \in \calH\calG[X]$ by
\begin{align*}
\HE(G/{\sim})&=\{p_\sim(e)\mid e\in \HE(G)\},\\
\HE(G\mid\sim)&=\{e\in \HE(G)\mid \forall x,y\in e,\: x\sim y\},\\
\HE(G\mid_\cap\sim)&=\{e\cap Y\mid e\in \HE(G),\: Y\in X/{\sim}\}\setminus \{\emptyset\}.
\end{align*}
We shall say that $\sim\in \calE^\cap(G)$ if for any $Y\in X/{\sim}$, $G_{\mid_\cap Y}$ is connected. We shall say that $\sim\in \calE^\subset(G)$ if for any $Y\in X/{\sim}$, $G_{\mid Y}$ is connected. We then define $\delta^\cap$ and $\delta^\subset$ by
\begin{align*}
\delta^\cap_\sim(G)&=\begin{cases}
G/{\sim}\otimes G\mid_\cap \sim \mbox{ if }\sim\in\calE^\cap(G),\\
0\mbox{ otherwise},
\end{cases}&
\delta^\subset_\sim(G)&=\begin{cases}
G/{\sim}\otimes G\mid \sim \mbox{ if }\sim\in\calE^\subset(G),\\
0\mbox{ otherwise}.
\end{cases}
\end{align*}
Then $(\calH\calG,m,\Delta^{\cap,\cap},\delta^\cap)$ and $(\calH\calG,m,\Delta^{\subset,\subset},\delta^\subset)$ are CTDBs.
 We shall come back later to the case of $\Delta^{\cap,\subset}$ and $\Delta^{\subset,\cap}$ in Example \ref{ex4.2.3}.
\end{example}

\begin{example}\label{ex4.1.3}
Let us consider the species $\calP$ of posets (Example \ref{ex2.1.5}). Let $X$ be a finite set, $P\in \calP[X]$, and $\sim\in \calE[X]$. We define a relation $\leq_{P/{\sim}}$ on $X/{\sim}$, as the transitive closure of $\sim\cup \leq_P$.
This is a quasi-order, but usually not an order. We define an order $\leq_{P\mid\sim}$ on $X$, by
\[\leq_{P\mid\sim}=\leq_P\cap \sim.\]
We shall say that $\sim\in \calE[P]$ if the two following conditions are satisfied:
\begin{itemize}
\item $P/{\sim}$ is a poset.
\item For any $Y\in X/{\sim}$, $P_{\mid Y}$ is an indecomposable poset.
\end{itemize}
We then put
\[\delta_\sim(P)=\begin{cases}
P/{\sim}\otimes P\mid\sim\mbox{ if }\sim\in \calE[P],\\
0\mbox{ otherwise}.
\end{cases}\]
Then $(\calP,m,\Delta,\delta)$ is a CTDB. 
A similar structure of CTDB can be defined on the species $\calO\calG_{\mathrm{ac}}$ of acyclic oriented graphs of Example \ref{ex3.2.4}.
\end{example}

\begin{remark}
In all these examples, $\epsilon_\delta=\epsilon^\calB_+$. This is quite a common phenomenon, see Propositions \ref{prop4.1.15} and \ref{prop4.1.16}.
\end{remark}

\subsection{Morphisms}

\begin{defi}
Let $\calB$ and $\calB'$ be two CTDBs.
A morphism of CTDBs from $\calB$ to $\calB'$ is morphism of CTBs $\phi:\calB\longrightarrow \calB'$ such that for any finite set $X$, any $\sim\in \calE[X]$, $\delta'_\sim\circ \phi=(\phi\otimes \phi)\circ \delta_\sim$.
\end{defi}

\begin{example}
The morphism $\theta:\calO\calG_{\mathrm{ac}}\longrightarrow \calP$, sending an acyclic oriented graph $G$ to the order $\leq_G$, is a morphism of CTDBs. 
\end{example}

\begin{example}
For any finite set $X$, any graph $G\in \calG[X]$, we define a hypergraph $\theta_{\calG\rightarrow\calH\calG}(G)\in \calH\calG[X]$ by
\[\HE(\theta_{\calG\rightarrow\calH\calG}(G))=E(G)\sqcup\{\{x\}\mid x\in X\}.\]
Then $\theta_{\calG\rightarrow\calH\calG}$ is a morphism of CTDBs from $(\calG,m,\Delta,\delta)$ to $(\calH\calG,m,\Delta,\delta)$. 
\end{example}

\begin{example}
For any finite set $X$, any hypergraph $G\in \calH\calG[X]$, we define a graph $\theta_{\calH\calG\rightarrow\calG}(G)\in \calG[X]$ by
\[E(\theta_{\calH\calG\rightarrow\calG}(G))=\{Y\in \calS_2[X]\mid \exists Z\in \HE(G),Y\subseteq Z\}.\]
Then $\theta_{\calH\calG\rightarrow\calG}$ is a morphism of CTDBs from $(\calH\calG,m,\Delta^{\cap,\cap},\delta^\cap)$ to $(\calG,m,\Delta,\delta)$.
\end{example}

We have introduced in Theorem \ref{theo1.3.6}. Let us detail some examples for CTBs.

\begin{example}\begin{enumerate}
\item For $(\calG,m,\Delta,\delta)$, $\epsilon_\delta=\epsilon^\calG_+$, so the fundamental polynomial invariant of $\calG$ is the chromatic polynomial.
\\item For both $(\calH\calG,m,\Delta^{\cap,\cap},\delta^\cap)$ and $(\calH\calG,m,\Delta^{\subset,\subset},\delta^\subset)$, $\epsilon_\delta=\epsilon^{\calH\calG}_+$, so the fundamental polynomial invariants are the chromatic polynomials attached to $\Delta^{\cap,\cap}$ and $\Delta^{\subset,\subset}$.
\item For $(\calP,m,\Delta,\delta)$, $\epsilon_\delta=\epsilon^\calP_+$, so the fundamental polynomial invariant is the strict Ehrhart polynomial.
\end{enumerate}\end{example}

\subsection{Structure of the set of equivalences associated to an element}

We now consider a CTDB $(B,m,\Delta,\delta)$.

\begin{defi}
Let $X$ be a finite set, and $b\in \calB[X]$. 
\begin{enumerate}
\item We denote by $\sim_b^{(m)}$ the equality of $X$. 
\item We denote by $\sim_b^{(\ic)}$ the equivalence on $X$ whose classes are the elements $X_1,\ldots,X_k$ which underlie any decomposition of $b$ into indecomposable elements (Proposition \ref{prop2.2.8} and Lemma \ref{lem2.2.3}).
\end{enumerate}\end{defi}

\begin{prop}\label{prop4.1.7}
Let $X$ be a finite set, and $b\in \calB[X]$.
\begin{enumerate}
\item $\sim_b^{(m)}\in \calE[b]$, and there exists $\sim_b^{(M)}\in \calE[b]$, such that for any $\sim\in \calE[b]$, $\sim_b^{(m)}\subseteq \sim\subseteq \sim_b^{(M)}$.
Moreover, $\sim_b^{(M)}\subseteq \sim_b^{(\ic)}$.
\item For any $\sim\in \calE[b]$, $\ic(b/{\sim})\geq \ic(b)$. 
\item The unique $\sim\in \calE[b]$ such that $\epsilon_\delta(b\mid \sim)=1$ is $\sim_b^{(m)}$. Moreover, $b/{\sim}_b^{(m)}=b$.
\item The unique $\sim\in \calE[b]$ such that $\epsilon_\delta(b/{\sim})=1$ is $\sim_b^{(M)}$. 
Moreover, $b\mid\sim_b^{(M)}=b$.
\item $\epsilon_\delta(b)=\begin{cases}
1\mbox{ if }\sim_b^{(m)}=\sim_b^{(M)},\\
0\mbox{ otherwise}.
\end{cases}$
\item If $b\in \bfD\calB[X]$, then $\sim_b^{(m)}=\sim_b^{(M)}=\sim_b^{(\ic)}$, and $\epsilon_\delta(b)=1$.
\end{enumerate}
\end{prop}

\begin{proof}
2. Let us first prove that for any $\sim\in \calE[b]$, $\sim_b^{(m)}\subseteq \sim\subseteq \sim_b^{(\ic)}$.
The first inclusion is obvious. Let us decompose $b$ into indecomposable elements:
\[b=b_1\cdots b_k,\]
with $b_i\in \calB[X_i]$ for any $i\in [k]$. Let $\sim\in \calE[b]$. We put $\sim_i=\sim\cap X_i^2$, for $i\in [k]$.
By Proposition \ref{prop4.1.3},
\[\sim=\sim_1\sqcup\cdots\sqcup\sim_k\subseteq X_1^2\sqcup\cdots\sqcup X_k^2=\sim_b^{(\ic)}.\]
Proposition \ref{prop4.1.3} gives
\[b/{\sim}=(b_1/{\sim}_1)\cdots (b_k/{\sim}_k),\]
so $\ic(b/{\sim})\geq k=\ic(b)$, which proves the second item.

3. Let $\sim\in \calE[X]$. Then, as $\epsilon_\delta$ is a right counit,
\[(\id_{\calB[X/{\sim}]}\otimes \epsilon_\delta)\circ \delta_\sim(b)=\epsilon_\delta(b\mid\sim)b/{\sim}=\begin{cases}
b\mbox{ if }\sim=\sim_b^{(m)},\\
0\mbox{ otherwise}.
\end{cases}\]
Therefore, if $\sim\neq \sim_b^{(m)}$, then $\epsilon_\delta(b\mid\sim)=0$. Moreover, $\delta_{\sim_b^{(m)}}(b)\neq 0$, so $\sim_b^{(m)}\in \calE[b]$, and
 \[(\id_{\calB[X]} \otimes \epsilon_\delta)\circ \delta_{\sim_b^{(m)}}(b)=\epsilon_\delta\left(b\mid\sim_b^{(m)}\right)b/ \sim_b^{(m)}=b,\]
so $\epsilon_\delta\left(b\mid\sim_b^{(m)}\right)=1$ and $b/ \sim_b^{(m)}=b$.

4. As $\epsilon_\delta$ is a left counit,
\begin{align}
\label{EQ4.1}\sum_{\sim\in\calE[b]} \epsilon_\delta(b/{\sim})b\mid \sim&=b.
\end{align}
Therefore, there exists $\sim\in \calE[b]$, such that $b\mid\sim=b$. Applying this result for $b/{\sim}$ with $\sim\in \calE[b]$, there exists $\sim'\in \calE[X/{\sim}]$ such that $(b/{\sim})\mid\sim'=b/{\sim}$. 
By the third item applied to $b/{\sim}$, $\epsilon_\delta(b/{\sim})=\epsilon_\delta((b/{\sim})\mid\sim')\in \{0,1\}$.
Consequently, (\ref{EQ4.1}) implies that there exists exactly one $\sim\in \calE[b]$, such that $\epsilon_\delta(b/{\sim})=1$. We denote it by $\sim_b^{(M)}$. 
Then (\ref{EQ4.1}) implies that $b\mid\sim_b^{(M)} \sim=b$.

5. Let us assume that $\sim_b^{(m)}=\sim_b^{(M)}$. By the third item, 
$b/{\sim}_b^{(M)}=b$. By the fourth item, $\epsilon_\delta(b)=1$. Otherwise, by the third item, $b/{\sim}_b^{(m)}=b$. By the fourth item, as $\sim_b^{(m)}\neq \sim_b^{(M)}$, $\epsilon_\delta(b)=0$.

6. Let us assume that $b\in \bfD\calB[X]$. Then, as $\ic(b)=|X|$, $\sim_b^{(m)}=\sim_b^{(\ic)}$. By our first observation, $\sim_b^{(m)}=\sim_b^{(M)}=\sim_b^{(\ic)}$. By the third item, $b/{\sim}_b^{(m)}=b$.
By the fourth item, $\epsilon_\delta\left(b/{\sim}_b^{(m)}\right)=\epsilon_\delta(b)=1$.

1. It remains to prove that, if $\sim\in \calE[b]$, then $\sim\subseteq \sim_b^{(M)}$. Let $\sim\in \calE[b]$, such that $\sim$ is not included in $\sim_b^{(M)}$. By the coassociativity of $\delta$, and as $b\mid \sim_b^{(M)}=b$,
\[\left(\id_{\calB\left[X/{\sim}_b^{(M)}\right]}\otimes \delta_\sim\right)\circ \delta_{\sim_b^{(M)}}(b)
=b/{\sim}_b^{(M)}\otimes \delta_\sim(b)=0,\]
so $\delta_\sim(b)=0$ and $\sim\notin \calE[b]$. 
 \end{proof}

\begin{remark}
We obtain that $\epsilon_\delta$ is a combinatorial character of $\calB$, which was not part of the definition of a CTDB (we only asked for $\epsilon_\delta$ to be a character).
\end{remark}

\begin{cor}\label{cor4.1.8}
Let $X$, $Y$ be disjoint finite sets.
\begin{enumerate}
\item If $a\in \calB[X]$ and $b\in \calB[Y]$, then $\sim_{ab}^{(M)}=\sim_a^{(M)}\sqcup \sim_b^{(M)}$,
and $\sim_{ab}^{(m)}=\sim_a^{(m)}\sqcup \sim_b^{(m)}$.
\item For any finite set $X$, we put 
\[\bfC\calB[X]=\{b\in \calB[X]\mid \epsilon_\delta(b)=1\}.\]
Then $\bfC\calB$ is a twisted subbialgebra of $\calB$ which contains $\bfD\calB$.
\end{enumerate}
\end{cor}

\begin{proof}
1. $\sim_{ab}^{(m)}=\sim_a^{(m)}\sqcup\sim_b^{(m)}$ is obvious.
By Proposition \ref{prop4.1.3}, the maximal element of $\calE[ab]$ is $\sim_a^{(M)}\sqcup \sim_b^{(M)}$: this is $\sim_{ab}^{(M)}$.

2. $\bfC\calB$ is the subspecies $\calB_{\epsilon_\delta}$ associated to the combinatorial character $\epsilon_\delta$: by Proposition \ref{prop3.2.5}, it is a twisted subalgebra of $\calB$. 
By Proposition \ref{prop4.1.7}-6, it contains $\bfD\calB$. Let us prove that it is a twisted subcoalgebra. Let $X$, $Y$ be disjoint finite sets, and $b\in \calB[X\sqcup Y]$, such that $(X$, $Y)\in \calC_-(b)$. We put
\[\Delta_{X,Y}(b)=b_{\mid_LX}\otimes b_{\mid_RY}.\]
Note that $\calE[b]$ is reduced to the single element $\sim_b^{(m)}$, by Proposition \ref{prop4.1.7}-1 and 5. We put $\sim_X=\sim_{b_{\mid_LX}}^{(M)}$ and $\sim_Y=\sim_{b_{\mid_RY}}^{(M)}$. 
By Proposition \ref{prop4.1.4}, as $\sim_b^{(m)}=\sim_X\sqcup\sim_Y$, 
\[b\mid\sim=b=(b_{\mid_LX}\mid\sim_X) (b_{\mid_RY}\mid \sim_Y)=b_{\mid_LX}b_{\mid_RY}.\]
Therefore, $1=\epsilon_\delta(b)=\epsilon_\delta(b_{\mid_LX})\epsilon_\delta(b_{\mid_RY})$, so $\epsilon_\delta(b_{\mid_LX})=\epsilon_\delta(b_{\mid_RY})=1$, as $\epsilon_\delta$ is a combinatorial character of $\calB$.
\end{proof}

\begin{prop}\label{prop4.1.9}
Let $\calB$ be a CTDB. Then it is commutative. 
\end{prop}

\begin{proof}
Let $X$, $Y$ be disjoint finite sets, $a\in \calB[X]$, $b\in \calB[Y]$. We take $\sim_X=\sim_a^{(M)}$, and $\sim_Y=\sim_b^{(M)}$, and $\sim=\sim_X\sqcup\sim_Y$. By Proposition \ref{prop4.1.3} with $(X',X'')=(X,Y)\in \calC_-(ab)$,
\[(ab)\mid \sim=(a\mid\sim_X)(b\mid \sim_Y)=ab.\]
Still by Proposition \ref{prop4.1.3}, with $(X',X'')=(Y,X)\in \calC_-(ab)$,
\[(ab)\mid \sim=(b\mid\sim_Y)(a\mid \sim_X)=ba.\] 
So $ab=ba$.
\end{proof}

\begin{prop}\label{prop4.1.10}
Let $X$ be a finite set, $b\in \calB[X]$ and $\sim\in \calE[b]$. Then
\begin{align*}
\sim_{b/{\sim}}^{(m)}&=\overline{\sim},& \sim_{b/{\sim}}^{(M)}&=\overline{\sim_b^{(M)}},&
\sim_{b\mid\sim}^{(m)}&=\sim_b^{(m)},&\sim_{b\mid\sim}^{(M)}&=\sim.
\end{align*}
\end{prop}

\begin{proof}
As $\sim_{b/{\sim}}^{(m)}$ is the equality of $X/{\sim}$, it is $\overline{\sim}$. By Proposition \ref{lem4.1.2}, $\overline{\sim_b^{(M)}}\in \calE[b/{\sim}]$, and
\[(b/{\sim})/\overline{\sim_b^{(M)}}=b/{\sim}_b^{(M)}.\]
Therefore, $\epsilon_\delta\left((b/{\sim})/\overline{\sim_b^{(M)}}\right)=\epsilon_\delta\left(b/{\sim}_b^{(M)}\right)=1$. Hence, $\overline{\sim_b^{(M)}}=\sim_{b/{\sim}}^{(M)}$.

Both $\sim_b^{(m)}$ and $\sim_{b\mid\sim}^{(m)}$ are the equality of $X$. As $\overline{\sim}=\sim_{b/{\sim}}^{(m)}$, by Proposition \ref{lem4.1.2}, $\sim\in \calE[b\mid\sim]$, and
\[(b\mid\sim)/{\sim}=(b/{\sim})\mid \overline{\sim}=(b/{\sim})\mid \sim_{b/{\sim}}^{(m)},\]
which implies that $\epsilon_\delta\left((b\mid\sim)/{\sim}\right)=1$. Therefore, $\sim=\sim_{b\mid\sim}^{(M)}$. \end{proof}

\begin{defi}
Let $X$ be a nonempty finite set, and $b\in \calB[X]$. We shall say that $b$ is primitive if for any $(X_1,X_2)\in \comp[X]$, $\Delta_{X_1,X_2}(b)=0$.
We shall say that $b$ is group-like if $\delta_{\sim_b^{(m)}}(b)=b\otimes b$.
\end{defi}

\begin{remark}
In other words, $b\in \calB[X]$ is primitive if, and only if, $\calC_-(b)=\{(X)\}$.
\end{remark}

\begin{lemma}\label{lem4.1.13}
Let $X$ be a finite set, and $b\in \calB[X]$. The following assertions are equivalent:
\begin{enumerate}
\item $b\in \bfC\calB[X]$.
\item $\calE[b]=\left\{\sim_b^{(m)}\right\}$.
\item $b$ is group-like.
\end{enumerate}
\end{lemma}

\begin{proof}
$1.\Longrightarrow 2$. By Proposition \ref{prop4.1.7}-5, $\sim_b^{(m)}=\sim_b^{(M)}$.
By Proposition \ref{prop4.1.7}-1, $\calE[b]=\left\{\sim_b^{(m)}\right\}$.

$2.\Longrightarrow 3$. Then $\sim_b^{(m)}=\sim_b^{(M)}$. By Proposition \ref{prop4.1.7}-3 and 4, $b/{\sim}_b^{(m)}=b\mid\sim_b^{(m)}=b$.

$3.\Longrightarrow 1$. By Proposition \ref{prop4.1.7}-3, $\epsilon_\delta(b)=\epsilon_\delta\left(b\mid \sim_b^{(m)}\right) =1$.
\end{proof}

\begin{lemma}\label{lem5.14}
Let $X$ be a finite set, and $b\in \bfC\calB[X]$. Then $b$ is indecomposable if, and only if, $b$ is primitive.
\end{lemma}

\begin{proof}
Let us assume that $b$ is not primitive. Let $(X',X'')\in \calC_-(b)$. We put $\sim'=\sim_{b_{\mid_LX'}}^{(m)}$ and $\sim''=\sim_{b_{\mid_RX''}}^{(m)}$. Then $\sim=\sim'\sqcup \sim''=\sim_b^{(m)}$, and, by Proposition \ref{prop4.1.4},
\[b=b\mid\sim=(b_{\mid_LX'}\mid\sim')(b_{\mid_RX''}\mid\sim''),\]
so $b$ is decomposable.

Let us assume that $b$ is decomposable. There exists $(X',X'')\in \comp[X]$, $b'\in \calB[X']$, $b''\in \calB[X'']$, such that $b=b'b''$. Then, by Lemma \ref{lem2.1.3},
\[\Delta_{X',X''}(b)=b'\otimes b''\neq 0,\]
so $b$ is not primitive.. 
\end{proof}

\begin{prop}\label{prop4.1.14}
 As a twisted algebra, $\bfC\calB$ is generated by the subspecies of primitive group-like elements.
\end{prop}

\begin{proof}
Let $X$ be a finite set, $b\in \bfC\calB[X]$, which we decompose into indecomposable elements:
\[b=b_1\cdots b_k.\]
By Proposition \ref{prop3.2.5}, $b_i,\ldots,b_k\in \bfC\calB[X_i]$ for any $i\in [k]$. By Lemma \ref{lem5.14}, $b_1,\ldots,b_k$ are primitive (and group-like, as they belong to $\bfC\calB$), which gives the result. 
\end{proof}

\begin{prop}\label{prop4.1.15}
Let $\calB$ be a CTDB. The following assertions are equivalent:
\begin{enumerate}
\item For any finite set $X$, for any $b\in \calB[X]$, $\sim_b^{(\ic)}=\sim_b^{(M)}$.
\item For any finite set $X$, for any $b\in \calB[X]$, $\sim_b^{(\ic)}\in \calE[b]$.
\item For any finite set $X$, for any $b\in \calB[X]$, for any $\sim\in \calE[b]$, $\ic(b\mid\sim)\geq \cl(\sim)$.
\item For any finite set $X$, for any $b\in \calB[X]$, for any $\sim\in \calE[b]$, 
$\left\{\begin{array}{rcl}
\ic(b\mid\sim)&=&\cl(\sim),\\ 
\ic(b/{\sim})&=&\ic(b).
\end{array}\right.$
\item Let $X$ be a finite set. If $\calB[X]$ contains a primitive group-like element, then $|X|=1$.
\item $\epsilon_\delta=\epsilon^\calB_+$.
\item $\bfC\calB\subseteq \bfD\calB$.
\item $\bfC\calB=\bfD\calB$.
\end{enumerate}
Such a CTDB will be called regular.
\end{prop}

\begin{proof}
$1.\Longrightarrow 2$. Immediate.

$2.\Longrightarrow 3$. Let us consider $\sim'=\sim_{b\mid\sim}^{(\ic)}$. By hypothesis, $\sim'\in \calE[b\mid\sim]$, so
\[(\id_{\calB[X/{\sim}]} \otimes \delta_{\sim'})\circ \delta_\sim(b)=b/{\sim}\otimes \delta_{\sim'}(b\mid\sim)\neq 0.\]
Consequently, by the coassociativity of $\delta$, $\sim'\subseteq \sim$. Therefore, $\cl(\sim')=\ic(b\mid\sim)\geq \cl(\sim)$.

$3.\Longrightarrow 1$. Let $X$ be a finite set, $b\in \calB[X]$, and $\sim\in \calE[b]$. Then $\sim\subseteq \sim_b^{(\ic)}$, by Proposition \ref{prop4.1.7}-1. 
If $\sim \neq \sim_b^{(\ic)}$, then $\cl(\sim)>\cl\left(\sim_b^{(\ic)}\right)=\ic(b)$. By hypothesis, $\ic(b\mid\sim)\geq \cl(\sim)>\ic(b)$, so $b\mid\sim\neq b$. By Proposition \ref{prop4.1.7}-4, $\sim\neq \sim_b^{(M)}$.
Consequently, $\sim_b^{(M)}=\sim_b^{(\ic)}$.

We proved that $1$., $2$. and $3$. are equivalent.

$4.\Longrightarrow 3$. Obvious.

$1. \Longrightarrow 4$. Let $\sim\in \calE[b]$. By Proposition \ref{prop4.1.10}, $\sim=\sim_{b\mid\sim}^{(M)}$, so $\sim=\sim_{b\mid\sim}^{(\ic)}$. Consequently, 
\[\cl(\sim)=\cl\left(\sim_{b\mid\sim}^{(\ic)}\right)=\ic(b\mid\sim).\]
By Proposition \ref{prop4.1.10}, $\sim_{b/{\sim}}^{(\ic)}=\sim_{b/{\sim}}^{(M)}=\overline{\sim_b^{(M)}}=\overline{\sim_b^{(\ic)}}$, so
\[\ic(b/{\sim})=\cl\left(\sim_{b/{\sim}}^{(\ic)}\right)=\cl\left(\overline{\sim_b^{(\ic)}}\right)=\cl\left(\sim_b^{(\ic)}\right)=\ic(b).\]

We proved that $1.$, $2$, $3$. and $4$ are equivalent.

$1.\Longrightarrow 6$. Let $X$ be a finite set, and $b\in \calB[X]$. If $b\in \bfD\calB[X]$, then $\sim_b^{(m)}=\sim_b^{(\ic)}$. By Proposition \ref{prop4.1.7}-5, $\epsilon_\delta(b)=\epsilon^\calB_+(b)=1$. Otherwise, $\sim_b^{(m)}\neq \sim_b^{(\ic)}=\sim_b^{(M)}$. By Proposition \ref{prop4.1.7}-4 and 5,
\[\epsilon_\delta(b)=\epsilon_\delta\left(b/{\sim}_b^{(m)}\right)=0=\epsilon^\calB_+(b).\]
Consequently, $\epsilon_\delta=\epsilon^\calB_+$.

$6.\Longrightarrow 7$. Let $X$ be a finite set, and $b\in \bfC\calB[X]$. Then $\epsilon_\delta(b)=\epsilon^\calB_+(b)=1$, so $b\in \bfD\calB[X]$.

$7.\Longrightarrow 8$. As $\bfD\calB\subseteq \bfC\calB$, by Proposition \ref{prop4.1.7}-6.

$8.\Longrightarrow 5$. Let $X$ be a finite set, and $b\in \calB[X]$ be a primitive group-like element. Then $b\in \bfD\calB[X]$, so $\ic(b)=|X|$. By Lemma \ref{lem5.14}, $b$ is indecomposable, so $\ic(b)=1$. Therefore, $|X|=\ic(b)=1$.

$5.\Longrightarrow 1$. Let us assume that there exists a finite set $X$, and $b\in \calB[X]$, such that $\sim_b^{(M)}\neq \sim_b^{(\ic)}$. We choose $X$ and $b$, such that $|X|$ is minimal. Let us decompose $b$ into indecomposable elements:
\[b=b_1\cdots b_k.\]
Then
\begin{align*}
\sim_b^{(M)}&=\sim_{b_1}^{(M)}\sqcup\cdots\sqcup\sim_{b_k}^{(M)}\neq\sim_b^{(\ic)}=\sim_{b_1}^{(\ic)}\sqcup\cdots\sqcup\sim_{b_k}^{(\ic)}.
\end{align*}
Therefore, there exists $i\in [k]$, such that $\sim_{b_i}^{(M)}\neq\sim_{b_i}^{(\ic)}$: by minimality of $|X|$, $X=X_i$, so $b$ is indecomposable. As $\sim_b^{(M)}\neq \sim_b^{(\ic)}=X^2$, $\left|X/{\sim}_b^{(M)}\right|\geq 2$.
Note that $b/{\sim}_b^{(M)}$ is a group-like element of $X/{\sim}_b^{(M)}$. Let us show that it is primitive. Let $X'/{\sim}'$, $X''/{\sim}''$ be subsets of $X/{\sim}_b^{(M)}$, with $X/{\sim}_b^{(M)}=(X'/{\sim}')\sqcup (X''/{\sim}'')$, such that $\Delta_{X'/{\sim}',X''/{\sim}''}\left(b/{\sim}_b^{(M)}\right)\neq 0$. 
By Proposition \ref{prop4.1.4}, $b=b\mid \sim_b^{(M)}=\left(b_{\mid_LX'}\mid\sim'\right)\left(b_{\mid_RX''}\mid\sim''\right)$.
As $b$ is indecomposable, $X'=\emptyset$ or $X''=\emptyset$. Therefore, $X'/{\sim}'=\emptyset$ or $X''/{\sim}''=\emptyset$: $b/{\sim}_b^{(M)}$ is a primitive group-like element of $\calB\left[X/{\sim}_b^{(M)}\right]$, and $\left|X/{\sim}_b^{(M)}\right|\geq 2$. 
\end{proof}

\begin{example}
The CTDBs $\calG$, $\calO\calG$, $\calH\calG$, and $\calP$ are regular. Indeed, in all cases, by definition of $\delta$, for any $b$, $\ic(b\mid\sim)=\cl(\sim)$. So the third item is satisfied. 
\end{example}

\begin{example}\label{ex4.1.11}
Let us give an example of non regular CTDB. For any finite set $X=\{x,y\}$ of cardinality two, we consider $G_X\in \calO\calG[X]$, such that $A(G_X)=\{(x,y),\:(y,x)\}$.
We denote by $\calB_i$ the twisted subalgebra of $\calO\calG$ generated by these oriented graphs $G_X$. If $X=\{x,y\}$,
\[\Delta_{\{x\},\{y\}}(G_X)=\Delta_{\{y\},\{x\}}(G_X)=0,\]
so $\calB_i$ is a sub-CTDB of $\calO\calG$. We give it a coproduct $\delta$ as follows: if $X=\{x,y\}$ is a set of cardinality 2, 
\begin{align*}
&\forall \sim\in \calE[X],&\delta_\sim(G_X)&=\begin{cases}
G_X\otimes G_X\mbox{ if }\sim=\sim^{(m)}_X,\\
0\mbox{ otherwise}.
\end{cases}
\end{align*}
Then $\calB_i$ is a CTDB. It is not regular: indeed, 
\begin{align*}
\epsilon_\delta(G_X)&=1,&\epsilon^{\calB_i}_+(G_X)&=0.
\end{align*}\end{example}

\begin{prop}\label{prop4.1.16}
Let $\calB$ a CTDB. 
\begin{enumerate}
\item If $(\calB,m,\Delta)$ is non-degenerate, then $(\calB,m,\Delta,\delta)$ is regular.
\item If $(\calB,m,\Delta,\delta)$ is regular, then $\nd\calB$ inherits a contraction-restriction coproduct $\delta$, making $(\nd\calB,m,\Delta,\delta)$ a non-degenerate CTDB.
\end{enumerate}\end{prop}

\begin{proof}
1. Let $X$ be a finite set, $b\in \calB[X]$, and $\sim\in \calE[b]$. As $\calB$ is non-degenerate, $\calC_{\pp}(b/{\sim})\neq \emptyset$. Let us put $X/{\sim}=\{X_1,\ldots,X_k\}$, such that 
$(\{X_1\},\ldots,\{X_k\})\in \calC_{\pp}(b/{\sim})$, with $k=\cl(\sim)$.Then, putting $\sim_i=X_i^2$ for $i\in [k]$,
\begin{align*}
&m_{1,3,\ldots,2k-1,2\ldots 2k}\circ (\delta_{\sim_1}\otimes\cdots\otimes\delta_{\sim_k})\circ \Delta_{X_1,\ldots,X_k}(b)=\left(\Delta_{\{X_1\},\ldots,\{X_k\}}\otimes \id_{\calB[X]}\right)\circ \delta_\sim(b)\neq 0.
\end{align*}
Therefore, $(X_1,\ldots,X_k)\in \calC_-(b)$. Identifying the last tensor, we obtain
\[b\mid\sim=\left(b_{\mid_1 X_1}\mid\sim_1\right)\cdots \left(b_{\mid_k X_k}\mid\sim_k\right).\]
Consequently, $\ic(b\mid\sim)\geq k=\cl(\sim)$. We obtain that $\calB$ is regular.

2. Let us show that the biideal $\calI$ defining $\nd\calB$ is also a coideal for $\delta$. Let $X$ be a finite set. Recall that $\calI[X]$ is formed by the elements $b\in \calB[X]$, such that $\calC_+(b)=\emptyset$. For any $b\in \calB[X]$, 
\[\phi^\calB_+(b)=\sum_{k=1}^{|X|} |\comp_k[X]\cap \calC_+(b)|\frac{T(T-1)\cdots(T-k+1)}{k!}.\]
Therefore, $b\in \calI[X]$ if, and only if, $\phi^\calB_+(b)=0$.

Let $b\in \calI[X]$. Then $\phi^\calB_+(b)=0$. This implies that 
\[\epsilon_\delta(b)=\epsilon^\calB_+(b)=\phi^\calB_+(b)(1)=0.\]
As $\epsilon_\delta=\epsilon^\calB_+$, $\phi^\calB_+$ is the fundamental polynomial invariant of $\calB$. We obtain
\[\phi^\calB_+(b)(TT')=\sum_{\sim\in \calE[b]} \phi^\calB_+(b/{\sim})(b)(T)\phi^\calB_+(b\mid\sim)(T').\]
If $\phi^\calB_+(b/{\sim})(T)\neq 0$, and $\phi^\calB_+(b/{\sim})(T')\neq 0$, the leading monomial of $\phi^\calB_+(b/{\sim})(b)(T)\phi^\calB_+(b\mid\sim)(T')$ (for both $T$ and $T'$) has a strictly positive coefficients A sum of such polynomials cannot be zero. 
So, for any $\sim\in \calE[b]$, either $\phi^\calB_+(b/{\sim})(T)=0$, and $\phi^\calB_+(b/{\sim})(T')=0$, so either $b/{\sim}\in \calI[X/{\sim}]$ or $b\mid\sim \in \calI[X]$. We proved that $\calI$ is also a coideal for $\delta$. \end{proof}

\begin{defi}
Let $\calB$ be a CTDB. For any finite set $X$ and any $b\in \calB[X]$, we denote by $\calC_\delta(b)$ the sets of compositions $(X_1,\ldots,X_k)\in \comp[X]$ such that
\[\epsilon_\delta^{\otimes k}\circ \Delta_{X_1,\ldots,X_k}(b)=1.\]
If $\calB$ is regular, then $\calC_\delta(b)=\calC_+(b)$.
\end{defi}

\begin{notation}
Let $X$ be a finite set, and $\sim\in \calE[X]$. For any $(X_1,\ldots,X_k)\in \comp[X/{\sim}]$, we put $\upsilon_\sim(X_1,\ldots,X_k)=(\varpi_\sim^{-1}(X_1),\ldots,\varpi_\sim^{-1}(X_k))\in \comp[X]$. This defines an injective map $\upsilon_\sim:\comp[X/{\sim}]\longrightarrow \comp[X]$.
\end{notation}

\begin{prop}
Let $\calB$ a CTDB, $X$ a finite set, and $b\in \calB[X]$. Then
\[\calC_-(b)=\bigsqcup_{\sim\in \calE[b]} \upsilon_\sim(\calC_\delta(b/{\sim})).\]
\end{prop}

\begin{proof}
$\subseteq$. Let $(X_1,\ldots,X_k)\in \calC_-(b)$. We put
\[\Delta_{X_1,\ldots,X_k}(b)=b_1\otimes \cdots\otimes b_k,\]
with $b_i\in \calB[X_i]$ for $i\in [k]$. We put $\sim_i=\sim_{b_i}^{(M)}$ for $i\in [k]$ and $\sim=\sim_1\sqcup\cdots\sqcup\sim_k\in \calE[X]$. Then $\epsilon_\delta({b_i}/{\sim_i})=1$. Moreover, as $\delta$ and $\Delta$ are compatible, 
\begin{align*}
\left(\Delta_{X_1/{\sim_1},\ldots, X_k/{\sim_k}}\otimes \id_{\calB[X]}\right)\circ \delta_\sim(b)&=m_{1,3,\ldots,2k-1,2\cdots 2k}\circ(\delta_{\sim_1}\otimes\cdots \otimes \delta_{\sim_k}) \circ \Delta_{X_1,\ldots,X_k}(b)\\
&={b_1}/{\sim_1}\otimes\cdots\otimes {b_k}/{\sim_k}\otimes b_1\mid\sim_1\cdots b_k\mid\sim_k\neq 0,
\end{align*}
so $\delta_\sim(b)\neq 0$, which gives $\sim\in \calE[b]$, and $(X_1/{\sim_1},\ldots, X_k/{\sim_k})\in \calC_-(b/{\sim})$. Moreover, identifying,
\[\Delta_{X_1/{\sim_1},\ldots, X_k/{\sim_k}}(b/{\sim})={b_1}/{\sim_1}\otimes\cdots\otimes {b_k}/{\sim_k},\]
so $\epsilon_\delta^{\otimes k}\circ \Delta_{X_1/{\sim_1},\ldots, X_k/{\sim_k}}(b/{\sim})=1$, which gives $(X_1/{\sim_1},\ldots, X_k/{\sim_k})\in \calC_\delta(b/{\sim})$ and finally $(X_1,\ldots,X_k)\in \upsilon_\sim(\calC_\delta(b/{\sim}))$.

$\supseteq$. Let $\sim\in \calE[b]$ and $(X_1/{\sim_1},\ldots, X_k/{\sim_k})\in \calC_\delta(b/{\sim})$. Then 
\begin{align*}
\sim&=\sim_1\sqcup\cdots\sqcup\sim_k,&\upsilon_\sim(X_1/{\sim_1},\ldots, X_k/{\sim_k})&=(X_1,\ldots,X_k). 
\end{align*}
As $\delta$ and $\Delta$ are compatible,
\begin{align*}
m_{1,3,\ldots,2k-1,2\cdots 2k}\circ(\delta_{\sim_1}\otimes\cdots \otimes \delta_{\sim_k}) \circ \Delta_{X_1,\ldots,X_k}(b)&=\left(\Delta_{X_1/{\sim_1},\ldots, X_k/{\sim_k}}\otimes \id_{\calB[X]}\right)\circ \delta_\sim(b)\neq 0,
\end{align*}
so $\Delta_{X_1,\ldots,X_k}(b)\neq 0$, and $(X_1,\ldots,X_k)\in \calC_-(b)$.

Let us finally prove that the union is disjoint. Let $\sim\in \calE[b]$ and $(X_1/{\sim_1},\ldots, X_k/{\sim_k})\in \calC_\delta(b/{\sim})$. Then $\sim=\sim_1\sqcup\cdots\sqcup\sim_k$.
We put $\Delta_{X_1/{\sim_1},\ldots, X_k/{\sim_k}}(b/{\sim}))=b'_1\otimes\cdots\otimes b'_k$, with $b'_i\in \calB[X_i/{\sim_i}]$ for $i\in [k]$. With the notations of the $\subseteq$, for any $i\in [k]$, $b'_i=b_i/{\sim_i}$. As $\epsilon_\delta(b'_i)=1$, $\sim_i=\sim_{b_i}^{(M)}$, and we obtain that
$\sim=\sim_{b_1}^{(M)}\sqcup\cdots\sqcup \sim_{b_k}^{(M)}$. Therefore, $\sim$ is entirely determined by $(X_1,\ldots,X_k)$, so the union is disjoint.
\end{proof}

\begin{cor}
Let $\calB$ a regular CTDB, $X$ a finite set, and $b\in \calB[X]$. Then
\[\calC_-(b)=\bigsqcup_{\sim\in \calE[b]} \upsilon_\sim(\calC_+(b/{\sim})).\]
\end{cor}

\subsection{Reduction the regular case}

Let $\calB$ be a CTDB.

\begin{defi}
Let $X$ be a finite set, $b\in \calB[X]$. We put $\displaystyle \sim_b^{(\epsilon)}=\sim_{b\mid\sim_b^{(m)}}^{(\ic)}$.
\end{defi}

\begin{remark}
If $\calB$ is regular, then $\sim_b^{(\epsilon)}=\sim_b^{(m)}$, as $b\mid\sim_b^{(m)}\in \bfC\calB[X]=\bfD\calB[X]$.
\end{remark}

\begin{prop}\label{prop4.1.21}\begin{enumerate}
\item Let $X$, $Y$ be disjoint finite sets, $a\in \calB[X]$ and $b\in \calB[Y]$. Then $\sim_{ab}^{(\epsilon)}=\sim_a^{(\epsilon)}\sqcup \sim_b^{(\epsilon)}$.
\item Let $X$, $Y$ be disjoint finite sets, and $b\in \calB[X\sqcup Y]$. We assume that $(X$, $Y)\in \calC_-(b)$. Then $\sim_b^{(\epsilon)}=\sim_{b_{\mid_LX}}^{(\epsilon)}\sqcup \sim_{b_{\mid_RY}}^{(\epsilon)}$.
\item Let $X$ be a finite set, $b\in \calB[X]$, and $\sim\in \calE[b]$. Then $\sim_{b\mid\sim}^{(\epsilon)}=\sim_b^{(\epsilon)}$ and $\sim_{b/{\sim}}^{(\epsilon)}=\overline{\sim_b^{(\epsilon)}}$.
\end{enumerate}\end{prop}

\begin{proof}
1. We have
\[(ab)\mid\sim_{ab}^{(m)}=\left(a\mid\sim_a^{(m)}\right) \left(b\mid\sim_b^{(m)}\right),\]
so
\[\sim_{ab}^{(\epsilon)}=\sim_{\left(a\mid\sim_a^{(m)}\right) \left(b\mid\sim_b^{(m)}\right)}^{(\ic)}=\sim_{a\mid\sim_a^{(m)}}^{(\ic)}\sqcup \sim_{b\mid\sim_b^{(m)}}^{(\ic)}=\sim_a^{(\epsilon)}\sqcup \sim_b^{(\epsilon)}.\]

2. We have
\begin{align*}
m_{1,3,24}\circ \left(\delta_{\sim_{b_{\mid_LX}}^{(m)}}\otimes \delta_{\sim_{b_{\mid_RY}}^{(m)}}\right)\circ \Delta_{X,Y}(b)
&=b_{\mid_LX}/{\sim}_{b_{\mid_LX}}^{(m)}\otimes b_{\mid_RY}/{\sim}_{b_{\mid_RY}}^{(m)}\otimes \left(b_{\mid_LX}\mid \sim_{b_{\mid_LX}}^{(m)}\right) \left(b_{\mid_RY}\mid \sim_{b_{\mid_RY}}^{(m)}\right)\\
=\left(\Delta_{X/{\sim}_{b_{\mid_LX}}^{(m)},Y/{\sim}_{b_{\mid_RY}}^{(m)}}\otimes \id_{\calB[X]}\right)\circ \delta_{\sim_b^{(m)}}(b)
&=\left(b/{\sim}_b^{(m)}\right)'\otimes \left(b/{\sim}_b^{(m)}\right)''\otimes b\mid \sim_b^{(m)}. 
\end{align*}
This gives
\[b\mid \sim_b^{(m)}= \left(b_{\mid_LX}\mid \sim_{b_{\mid_LX}}^{(m)}\right) \left(b_{\mid_RY}\mid \sim_{b_{\mid_RY}}^{(m)}\right).\]
Therefore,
\[\sim_{b\mid \sim_b^{(m)}}^{(\ic)}=\sim_{b_{\mid_LX}\mid \sim_{b_{\mid_LX}}^{(m)}}^{(\ic)}\sqcup \sim_{b_{\mid_RY}\mid \sim_{b_{\mid_RY}}^{(m)}}^{(\ic)},\]
which gives the result.

3. As $\sim_{b\mid\sim}^{(m)}=\sim_b^{(m)}$, 
\[(b\mid\sim)\mid \sim_{b\mid\sim}^{(m)}=(b\mid\sim)\mid\sim_b^{(m)}=b\mid\sim_b^{(m)}.\]
This gives
\[\sim_{b\mid\sim}^{(\epsilon)}=b_{(b\mid\sim)\mid \sim_{b\mid\sim}^{(m)}}^{(\ic)}=
\sim_{b\mid\sim_b^{(m)}}^{(\ic)}=\sim_b^{(\epsilon)}.\]

We now prove that $\sim_{b/{\sim}}^{(\epsilon)}\subseteq \overline{\sim_b^{(\epsilon)}}$. We denote by $\sim'$ the unique element of $\calE[X]$ containing $\sim$, such that $\overline{\sim'}=\sim_{b/{\sim}}^{(\ic)}$.
By Proposition \ref{lem4.1.2},
\[(b/{\sim})\mid \sim_{b/{\sim}}^{(m)}=(b/{\sim})\mid\overline{\sim}=(b\mid\sim)/{\sim}.\]
Let us denote by $b\mid\sim=b_1\cdots b_k$ the decomposition of $b\mid\sim$ into indecomposable elements, with $b_i\in\calB[X_i]$, for $i\in [k]$. We put $\sim_i=\sim\cap X_i^2$.
As $\overline{\sim}=\sim_{b/{\sim}}^{(m)}\in \calE[b/{\sim}]$, by Proposition \ref{lem4.1.2}, $\sim\in \calE[b\mid\sim]$.By Proposition \ref{prop4.1.3}, $\sim=\sim_1\sqcup\cdots\sqcup\sim_k$, and
\[(b\mid\sim)/{\sim}=(b_1/{\sim}_1)\cdots (b_k/{\sim}_k).\]
Therefore, 
\[\sim_b^{(\epsilon)}=\sim_{(b\mid\sim)/{\sim}}^{(\ic)}=\sim_{b_1/{\sim}_1}^{(\ic)}\sqcup\cdots\sqcup \sim_{b_k/{\sim}_k}^{(\ic)}.\]
So $\sim_b^{(\epsilon)}\subseteq (X_1/{\sim}_1)^2\sqcup\cdots\sqcup (X_k/{\sim}_k)^2$, which gives
\[\sim'\subseteq X_1^2\sqcup\cdots\sqcup X_k^2=\sim_b^{(\epsilon)}.\]
Consequently, $\sim_{b/{\sim}}^{(\epsilon)}\subseteq \overline{\sim_b^{(\epsilon)}}$.

Let us assume that there exist a finite set $X$, $b\in \calB[X]$, and $\sim\in \calE[b]$, such that $\sim_{b/{\sim}}^{(\epsilon)}\subsetneq \overline{\sim_b^{(\epsilon)}}$.
We choose such a triple $(X,b,\sim)$, with $|X|$ minimal. Let $\sim'\in \calE[b]$, such that $\sim\subseteq \sim'$. By Proposition \ref{lem4.1.2}, $(b/{\sim})/\overline{\sim'}=b/{\sim}'$, so
\[\sim_{b/{\sim}'}^{(\epsilon)}=\sim_{(b/{\sim})/\overline{\sim'}}^{(\epsilon)}\subseteq \overline{\sim_{b/{\sim}}^{(\epsilon)}}\subsetneq \overline{\overline{\sim_b^{(\epsilon)}}}=\overline{\sim_b^{(\epsilon)}}.\]
Therefore, we can replace $\sim$ by $\sim'$, without loss of generality. We choose to replace $\sim$ by $\sim_b^{(M)}$: this allows to assume that $b/{\sim}$ is group-like and that $b\mid\sim=b$. 
Let us decompose $b/{\sim}$ into indecomposable elements:
\[b/{\sim}=b'_1\cdots b'_k,\]
with $b'_i\in \calB[X_i/{\sim}_i]$ for $i\in [k]$. Then $(X_1/{\sim}_1,\ldots,X_k/{\sim}_k)\in \calC_-(b/{\sim})$, and
\begin{align*}
&\left(\Delta_{X_1/{\sim}_1,\ldots,X_k/{\sim}_k}\otimes \id_{\calB[X]}\right)\circ \delta_\sim(b)\\
&=(b/{\sim})_{\mid_1 X_1/{\sim}_1}\otimes\cdots\otimes (b/{\sim})_{\mid_k X_k/{\sim}_k}\otimes b\\
&=m_{1,3,\ldots,2k-1,2\ldots 2k}\circ (\delta_{\sim_1}\otimes\cdots\otimes \delta_{\sim_k})\circ\Delta_{X_1,\ldots,X_k}(b)\\
&=\left(b_{\mid_1 X_1}\right)/{\sim}_1\otimes\cdots\otimes=\left(b_{\mid_k X_k}\right)/{\sim}_k\otimes \left(b_{\mid_1 X_1}\right)\mid \sim_1 \cdots \left(b_{\mid_k X_k}\right)\mid \sim_k.
\end{align*}
We put $b_i=\left(b_{\mid_i X_i}\right)\mid \sim_i\in \calB[X_i]$ for $i\in [k]$. Then, $b=b_1\cdots b_k$. By the first item,
\[\sim_b^{(\epsilon)}=\sim_{b_1}^{(\epsilon)}\sqcup\cdots\sqcup \sim_{b_k}^{(\epsilon)}.\]
Moreover, $b/{\sim}=b_1/{\sim}_1\cdots b_k/{\sim}_k$, by Proposition \ref{prop4.1.3}, so, by the first item,
\[\sim_{b/{\sim}}^{(\epsilon)}=\sim_{b_1/{\sim}_1}^{(\epsilon)}\sqcup\cdots\sqcup \sim_{b_k/{\sim}_k}^{(\epsilon)}.\]
As $\sim_{b/{\sim}}^{(\epsilon)}\subsetneq \overline{\sim_b^{(\epsilon)}}$, there exists $i\in [k]$, such that $\sim_{b_i/{\sim}_i}^{(\epsilon)}\subsetneq \overline{\sim_{b_i}^{(\epsilon)}}$.
By minimality of $|X|$, necessarily $k=1$: we obtain that $b/{\sim}$ is indecomposable. By Proposition \ref{lem4.1.2},
\[(b/{\sim})\mid \sim_{b/{\sim}}^{(m)}=(b/{\sim})\mid\overline{\sim}=(b\mid\sim)/{\sim}=)b/{\sim},\]
so $\sim_{b/{\sim}}^{(\epsilon)}=\sim_{b/{\sim}}^{(\ic)}=(X/{\sim})^2$, as $b/{\sim}$ is indecomposable. Proposition \ref{prop4.1.3} implies that $\ic(b/{\sim})\geq \ic(b)$. 
As $\ic(b/{\sim})=1$, $\ic(b)=1$, so $b$ is indecomposable. Therefore,
\[\sim_b^{(\epsilon)}=\sim_{b\mid\sim}^{(\ic)}=\sim_b^{(\ic)}=X^2,\]
and finally $\sim_{b/{\sim}}^{(\epsilon)}=\overline{\sim_b^{(\epsilon)}}=(X/{\sim})^2$: this is a contradiction with the choice of $b$. Therefore, for any finite set $X$, any $b\in \calB[X]$, and any $\sim\in \calE[b]$, $\sim_{b/{\sim}}^{(\epsilon)}=\overline{\sim_b^{(\epsilon)}}$.
\end{proof}

Let us consider the double bialgebra $\calF[\calB]$. For any finite set $X$, and any $b\in \calB[X]$, we denote by $\tilde{b}$ is isoclass, that it to say its class in $\calF[\calB]$. The bialgebra $\calF[\calB]$ is naturally graded by the cardinality of the underlying finite set $X$:
\begin{enumerate}
\item If $\tilde{a}\in \calF[\calB]_k$, and $\tilde{b}\in \calF[\calB]_l$, then $\tilde{a}\tilde{b}\in \calF[\calB]_{k+l}$.
\item If $\tilde{b}\in \calF[\calB]_n$, then $\displaystyle \Delta(\tilde{b}) \in \bigoplus_{k=0}^n \calF[\calB]_k\otimes \calF[\calB]_{n-k}$, and $\displaystyle \delta(\tilde{b}) \in \bigoplus_{k=0}^n \calF[\calB]_k\otimes \calF[\calB]_n$.
\end{enumerate}
The group-like elements (for $\delta)$ of $\calF[\calB]$ are the isoclasses of elements of $\bfC\calB$.

Let us give another graduation to $\calF[\calB]$. For any finite set $X$, and any $b\in \calB[X]$, we put
\[\deg_\reg(\tilde{b})=\cl\left(\sim_b^{(\epsilon)}\right).\]
We denote by $\calF[\calB]^{(n)}$ the subspace of $\calF[\calB]$ of homogeneous elements of degree $n$, for this degree $\deg_\reg$. Proposition \ref{prop4.1.21} implies that:
\begin{enumerate}
\item If $\tilde{a}\in \calF[\calB]^{(k)}$, and $\tilde{b}\in \calF[\calB]^{(l)}$, then $\tilde{a}\tilde{b}\in \calF[\calB]^{(k+l)}$.
\item If $\tilde{b}\in \calF[\calB]^{(n)}$, then $\displaystyle \Delta(\tilde{b}) \in \bigoplus_{k=0}^n \calF[\calB]^{(k)}\otimes \calF[\calB]^{(n-k)}$, and $\displaystyle \delta(\tilde{b}) \in \bigoplus_{k=0}^n \calF[\calB]^{(k)}\otimes \calF[\calB]^{(n)}$.
\end{enumerate}
It is then possible to define a CTDB $\reg\calB$, such that, as a graded bialgebra, $\calF[\reg\calB]$ is $\calF[\calB]$, graded with the graduation induced by $\deg_\reg$. In other words,
\begin{itemize}
\item For any $n\in\N$, the space of coinvariants of $\calB[n]$ is generated by isoclasses of elements $b\in \calB[n]$.
\item For any $n\in \N$, the space of coinvariants of $\reg\calB[n]$ is generated by isoclasses of elements $b$ such that $\cl\left(\sim_b^{(\epsilon)}\right)=n$.
\end{itemize}
Let $b\in\bfC\reg\calB[X]$. Note that $|X|=\cl\left(\sim_b^{(\epsilon)}\right)$, by definition of $\reg\calB$. Then $\tilde{b}$ is a group-like element of $\calF[\reg\calB]=\calF[\calB]$. Therefore, $b$ is a group-like element of $\calB[Y]$, for a certain set $Y$. 
Therefore, 
\[b\mid\sim_b^{(m)}=b,\]
so $\sim_b^{(\epsilon)}=\sim_b^{(\ic)}$. Therefore, 
\[|X|=|X|=\cl\left(\sim_b^{(\epsilon)}\right)=\cl\left(\sim_b^{(\ic)}\right)=\ic(b).\]
We conclude that $b\in \bfD\reg\calB[X]$. By the characterization of regularity (item 7), $\reg\calB$ is regular.

\begin{remark}
The unique double bialgebra morphism from $\calF[\calB]$ to $\K[T]$ is $\phi_{\epsilon_\delta}$. The unique double bialgebra morphism from $\calF[\reg\calB]$ to $\K[T]$ is $\phi^\calB_+$, as $\reg\calB$ is regular. 
As $\calF[\calB]=\calF[\reg\calB]$, if $b\in \calB[X]$, denoting by $b'$ a lifting of $\tilde{b}$ to $\reg\calB$, we obtain
\[\phi_{\epsilon_\delta}(b)=\phi^\calB_+(b').\]
As a consequence, the polynomials obtained by application of the unique invariant compatible with the whole structure of a CTDB, regular or not, always belong to one of the $\calP_+(n)$'s.
\end{remark}

\begin{example}
Let us consider $\calB_i$, defined in Example \ref{ex4.1.11}. For any finite set $X$, of cardinality 1, let us consider an element $H_X$, which can be seen as the oriented graph with a loop on a single vertex.
Then $\reg\calB_i$ is the free commutative twisted algebra generated by these elements, with the coproducts defined as follows: if $|X|=1$, denoting by $\sim_X$ the unique element of $\calE[X]$,
\begin{align*}
\Delta_{X,\emptyset}(H)&=H\otimes 1,&\Delta_{\emptyset,X}(H)&=1\otimes H,&\delta_{\sim_X}(H_X)&=H_X\otimes H_X.
\end{align*}
Consequently, $\reg\calB_i$ is isomorphic to $\Com$. 
\end{example}

\subsection{Links with operads}

Let us first give a few reminders on compositions of species. Let $\calQ$ and $\calR$ be $\K$-linear species, with $\calR[\emptyset]=(0)$. The composition of species $\calQ\circ\calR$ is defined by
\[\calQ\circ\calR[X]=\bigoplus_{\sim\in \calE[X]} \calQ[X/{\sim}] \otimes \bigotimes_{Y\in X/{\sim}} \calR[Y].\]
This makes the category of connected $\K$-linear species (that are, species $\calQ$ with $\calQ[\emptyset]=(0)$) a monoidal (non symmetric) category. The unit is the species $\calU$ defined by
\[\calU[X]=\begin{cases}
\K\mbox{ if }|X|=1,\\
0\mbox{ otherwise}.
\end{cases}\]
An operad is an associative algebra in this monoidal category, see \cite{Dotsenko2016,Loday2012,Markl2002,Mendez2015} for more details and results.

\begin{theo}\label{theo4.1.22}
Let $\calB$ be a CTDB. We assume that for any finite set $X$, $\calB[X]$ is finite. 
\begin{enumerate}
\item Let us assume that $\calB$ is regular. We define a non unitary operad on $\calB$ as follows: for any finite set $X$, $\sim\in \calE[X]$, for any $b_Y\in \calB[Y]$ for $Y\in X/{\sim}$, and any $b\in \calB[X/{\sim}]$,
\[b\circ\left(\bigotimes_{Y\in X/{\sim}} b_Y\right)=\sum_{\substack{b'\in \calB[X],\\ \sim\in \calE[b'],\\ b'/{\sim}=b,\\ b'\mid\sim=\prod b_Y}}b'.\]
\item Let us assume that $\calB$ is regular. This operadic structure is compatible with the product, in the following sense: if $X$, $Y$ are disjoint finite sets, $\sim_X\in \calE[X\sim_X]$, $\sim_Y\in \calE|Y/{\sim}_Y]$, for any $a\in \calB[X]$, $b\in \calB[Y]$, $b'_Z\in \calB[Z]$ for $Z\in (X/{\sim}_X)\sqcup(Y/{\sim}_Y)$,
\[ab\circ\left(\bigotimes_{Z\in (X/{\sim}_X)\sqcup (Y/{\sim}_Y)} b'_Z\right)=\left(a\circ\left(\bigotimes_{Z\in X/{\sim}_X} b'_Z\right)\right)\left(b\circ\left(\bigotimes_{Z\in Y/{\sim}_Y} b'_Z\right)\right).\]
\item Let us assume that $\calB$ is regular. We denote by $\ind\calB$ the subspecies of indecomposable elements of $\calB$. If $\calB$ is regular, then $\ind\calB$ is a suboperad of $\calB$. It has a unit, defined by
\[U_X=\sum_{b\in \calB[X]} b,\]
for any set $X$ of cardinality $1$. 
\end{enumerate}\end{theo}

\begin{proof}
Let $X$ be a finite set. 
\[\calB\circ\calB\circ\calB[X]=\bigoplus_{\sim'\subset \sim\in \calE[X]}\calB[X/{\sim}]\otimes \bigotimes_{Y\in X/{\sim}}\calB[Y/{\sim}'\cap Y^2]\otimes\bigotimes_{Y\in X/{\sim}'}\calB[Y].\]
Let $\sim'\subseteq\sim \in \calE[X]$, $b\in \calB[X/{\sim}]$, $b_Y\in \calB[X/{\sim}'\cap Y^2]$ for $Y\in X/{\sim}$, and $b'_Z\in \calB[Z]$ for $Z\in X/{\sim}'$. 
\begin{align*}
\left(b\circ \bigotimes_{Y\in X/{\sim}} b_Y\right)\circ \bigotimes_{Z\in X/{\sim}'} b'_Z&=\sum_{\substack{b'\in \calB[X/{\sim}'],\\ \overline{\sim}\in \calE[b'],\\ b'/\overline{\sim}=b,\\ b'\mid\overline{\sim}=\prod b_Y}}
\sum_{\substack{b''\in \calB[X],\\\sim'\in \calE[b''],\\ b''/{\sim}'=b',\\ b''\mid\sim'=\prod b'_Z}} b''\\
&=\sum_{\substack{b''\in \calB[X],\\ \sim'\in \calE[b''],\: \overline{\sim}\in \calE[b''/{\sim}'],\\ (b''/{\sim}')/\overline{\sim}=b,\\ (b''/{\sim}')\mid\overline{\sim}=\prod b_Y,\\b''\mid \sim'=\prod b'_Z}}b''\\
&=\sum_{\substack{b''\in \calB[X],\\ \sim\in \calE[b''],\: \sim'\in \calE[b''\mid\sim],\\ b''/{\sim}=b,\\ (b''\mid\sim)/{\sim}'=\prod b_Y,\\ ((b''\mid\sim)\mid\sim')=\prod b'_Z}}b''\\
&=\sum_{\substack{b'_Y\in \calB[Y],\\ \sim'\cap Y^2\in \calE[b'_Y],\\ b'_Y/{\sim}_Y\cap Y^2=b_Y,\\ b'_Y\mid \sim_Y\cap Y^2=\prod_{Z\subseteq Y}b'_Z}}
\sum_{\substack{b''\in \calB[X],\\ \sim\in \calE[b''],\\ b''/{\sim}=b,\\ b''\mid\sim=\prod b'_Y}}b''\\
&=b\circ\left( \bigotimes_{Y\in X/{\sim}} b_Y\circ\left(\bigotimes_{Z\in X/{\sim}',\: Z\subseteq Y} b'_Z\right)\right).
\end{align*}
We used Proposition \ref{lem4.1.2} for the third equality. So $\circ$ is associative.

We now assume $\calB$ is regular. Let us prove the compatibility with the product.
\begin{align*}
ab\circ\left(\bigotimes_{Z\in (X/{\sim}_X)\sqcup (Y/{\sim}_Y)} b'_Z\right)&=\sum_{\substack{b''\in \calB[X\sqcup Y],\\ \sim_X\sqcup\sim_Y\in \calE[b''],\\ b''/{\sim}_X\sqcup \sim_Y=ab,\\ b''\mid \sim_X\sqcup\sim_Y=\prod b'_Z}} b''.
\end{align*}
If $b''/{\sim}_X\sqcup\sim_Y=ab$, as $\ic(b/{\sim}_X\sqcup\sim_Y)=\ic(ab)$ (item 4), necessarily $\overline{\sim_{b''}}^{(\ic)}=\sim_{b''/{\sim}}^{(\ic)}$, so $b''$ can be written under the form $b''=a'b'$, with $a'\in \calB[X]$ and $b'\in \calB[Y]$. As $\delta$ is compatible with the product, we obtain
\begin{align*}
ab\circ\left(\bigotimes_{Z\in (X/{\sim}_X)\sqcup (Y/{\sim}_Y)} b'_Z\right)&=\sum_{\substack{a'\in \calB[X],\\ \sim_X\in \calE[a'],\\ a'/{\sim}_X=a,\\ a'\mid\sim_X=\prod_{Z\in X/{\sim}_X} b'_Z}}
\sum_{\substack{b'\in \calB[Y],\\ \sim_Y\in \calE[b'],\\ b'/{\sim}_Y=b,\\ b'\mid\sim_Y=\prod_{Z\in Y/{\sim}_Y} b'_Z}} a'b'\\
&=\left(a\circ\left(\bigotimes_{Z\in X/{\sim}_X} b'_Z\right)\right)\left(b\circ\left(\bigotimes_{Z\in Y/{\sim}_Y} b'_Z\right)\right).
\end{align*}

Let $X$ be a finite set, $b\in \ind\calB[X/{\sim}]$, $b_Y\in \ind\calB[Y]$ for $Y\in X/{\sim}$, and finally $b'\in \calB[X]$, such that $b'$ appears in $\displaystyle b\circ \left(\bigotimes_{Y\in X/{\sim}}b_Y\right)$.
Then $\sim\in \calE[b']$, and $b'/{\sim}=b$. By definition of the regularity (item 4), $\ic(b')=\ic(b)=1$, so $b'\in \ind\calB[X]$: $\ind\calB$ is a suboperad.

Note that $\ind\calB[X]=\calB[X]$ if $|X|=1$, so $U_X\in \ell\ind\calB[X]$.
Let $X$ be a finite set, $b\in \ind\calB[X]$, and $\sim=X^2$. 
\[U_{X/{\sim}}\circ b=\sum_{b'\in \calB[X/{\sim}]} b'\circ b.\]
As $b$ is indecomposable, $\sim=\sim_b^{(\ic)}$. As $\calB$ is regular, $\sim\in \calE[b]$ (second item), and $b\mid\sim=b$. Therefore, if $b'\in \calB[X/{\sim}]$,
\[b'\circ b=\begin{cases}
b\mbox{ if }b'=b/{\sim},\\
0\mbox{ otherwise}.
\end{cases}\]
Therefore, $U_{X/{\sim}}\circ b=b$, and $U$ is a left unit on $\ind\calB$. Moreover, if $b\in \calB[X]$, denoting by $\sim$ the equality of $X$,
\begin{align*}
b\circ \left(\bigotimes_{x\in X} U_{\{x\}}\right)&=\prod_{b_x \in \calB[\{x\}]}\sum_{\substack{b'\in \calB[X],\\ \sim \in \calE[b'],\\ b'/{\sim}=b,\\ b'\mid\sim=\prod b_x }}b'.
\end{align*}
For any $b'\in \calB[X]$, $\sim\in \calE[b]$ and $b'/{\sim}=b'$. As $\calB$ is regular (item 4), $\ic(b\mid\sim)=\cl(\sim)=|X|$, so there exists a unique family $(b_x)_{x\in X}$, with $b_x\in \calB[\{x\}]$ for $x\in X$, such that $b'\mid\sim=\prod b_x$. Consequently, 
\[b\circ \left(\bigotimes_{x\in X} U_{\{x\}}\right)=b.\]
So $U$ is a right counit for the whole $\calB$, so for $\ind\calB$ also by restriction.
\end{proof}

\begin{remark}
$U$ is not a left counit on $\calB$: if $b\in \calB[X]$ is not indecomposable, denoting $\sim=X^2$, then $\sim\notin \calE[b]$, so $U_{X/{\sim}}\circ b=0$.
\end{remark}

\begin{example}
For the CTDB of of posets $\calP$, the underlying operadic structure on connected posets can be found in \cite{Foissy53}.
\end{example}

\subsection{Characters and polynomial invariants}

Let us first characterize characters invertible for the convolution product associated to the coproduct $\delta$. Recall that if $\calB$ is a CTDB, then for any finite set $X$,
\[\bfC\calB[X]=\{b\in \calB[X]\mid \epsilon_\delta(b)=1\}.\]

\begin{prop}\label{prop4.1.23}
Let $\calB$ be a regular CTDB, and let $f\in \Char(\calB)$. Then $f$ is invertible for the convolution product $*$ associated to $\delta$ if, and only if, for any $b\in \calB[1]$, $f(b)\neq 0$.
\end{prop}

\begin{proof}
$\Longrightarrow$. We denote by $g$ the inverse of $f$. Let $b\in \bfC\calB[1]$. Then 
\[f*g(b)=(f\otimes g)\circ \delta_{\sim_1}(b)=f(b)g(b)=\epsilon_\delta(b)=1,\]
where $\sim_1$ is the unique equivalence on $[1]$. So $f(b)\neq 0$.

$\Longleftarrow$. Let $X$ be a finite set, and $b\in \bfD\calB[X]$. Let $b_1,\ldots,b_n$ be the indecomposable components of $b$, with $b_i\in \calB[X_i]$ for $i\in [n]$. As $b\in \bfD\calB[X]$, $|X_i|=1$ for any $i\in [n]$. 
Therefore, $f(b)=f(b_1)\ldots f(b_n)$ and, by hypothesis on $f$, $f(b_i)\neq 0$ for any $i\in [n]$, so $f(b)\neq 0$. 

Let $X$ be a finite set, and $b\in \calB[X]$. We define $g(b)$ by a decreasing induction on $\ic(b)$ as follows:
\begin{enumerate}
\item if $\ic(b)=|X|$, then $b\in \bfC\calB[X]$, and we put $g(b)=\dfrac{1}{f(b)}$, which exists as we already noticed that $f(b)\neq 0$. 
Moreover, 
\begin{align*}
f*g(b)&=\sum_{\sim\in \calE[b]} (f\otimes g)\circ \delta_\sim(b)=(f\otimes g)\circ\delta_{\sim^{(\ic)}_b}(b)=f(b)g(b)=1=\epsilon_\delta(b).
\end{align*}
\item If $\ic(b)<|X|$, then $\epsilon_\delta(b)=0$. Moreover,
\begin{align*}
\sum_{\sim\in \calE[X]} \delta_\sim(x)&=\sum_{\sim\in \calE[b]\setminus \left\{\sim_b^{(\ic)}\right\}} b/{\sim}\otimes b\mid \sim+b/{\sim_b^{(\ic)}}\otimes b.
\end{align*}
If $\sim\in \calE[b]\setminus \left\{\sim_b^{(\ic)}\right\}$, then $\sim\subsetneq \sim_b^{(\ic)}$, so $\cl(\sim)>\cl\left(\sim_b^{(\ic)}\right)=\ic(b)$. As $\ic(b\mid\sim)=\ic(\sim)$ by regularity of $\calB$, $\ic(b\mid\sim)>\ic(b)$. 
Moreover, $b/{\sim_b^{(\ic)}} \in \bfD\calB\left[X/{\sim_b^{(\ic)}}\right]$, so $f\left(b/{\sim_b^{(\ic)}}\right)\neq 0$. Therefore, we can put
\[g(b)=-\frac{1}{f\left(b/{\sim_b^{(\ic)}}\right)}\sum_{\sim\in \calE[b]\setminus \left\{\sim_b^{(\ic)}\right\}} f(b/{\sim})g(b\mid \sim).\]
By construction, $f*g(b)=0=\epsilon_\delta(b)$.
\end{enumerate}
Therefore, $g$ is a right inverse of $f$ in $(\Hom(\calB,\lcom),*)$. We now define $h(b)$ by induction on $|X|$ as follows:
\begin{enumerate}
\item If $X=\emptyset$, $b=1_\calB$ and we put $h(b)=1$.
\item If $X\neq \emptyset$, then 
\[\sum_{\sim\in \calE[X]} \delta_\sim(b)=\sum_{\sim\in \calE[b]\setminus\left\{\sim_b^{(m)}\right\}}\delta_\sim(b)+b\otimes b\mid \sim_b^{(m)}.\]
If $\sim\in \calE[b]\setminus\left\{\sim_b^{(m)}\right\}$, then $|X/{\sim}|<|X|$. Moreover, $b\mid \sim_b^{(m)}\in \bfD\calB[X]$, so $f\left(b\mid \sim_b^{(m)}\right)\neq 0$. We then put
\[h(b)=\frac{1}{f\left(b\mid \sim_b^{(m)}\right)}\left(\epsilon_\delta(b)-\sum_{\sim\in \calE[b]\setminus\left\{\sim_b^{(m)}\right\}}(h\otimes f)\circ \delta_\sim(b)\right).\]
\end{enumerate}
By construction, $h\circ f(b)=\epsilon_\delta(b)$, so $h$ is a left inverse of $f$ in $(\Hom(\calB,\lcom),*)$. As $*$ is associative,
\[g=\epsilon_\delta*g=(h*f)*g=h*(f*g)=h*\epsilon_\delta=h,\]
so $f$ is invertible in in $(\Hom(\calB,\lcom),*)$. As the inverse of a character is a character (Lemma \ref{lem1.3.3}), $f$ is invertible in $(\Char(\calB),*)$.\end{proof}

\begin{theo}\label{theo4.1.24}
Let $\calB$ be a regular CTDB. We define a character of $\calB$ as follows: for any finite $X$, and any $b\in \calB[X]$, indecomposable,
\[g^\calB_{\pp}(b)=\ln_\calB(\epsilon_\delta)(b)=\sum_{k=1}^{|X|}\frac{(-1)^{k-1}}{k} |\calC_+(b)\cap \comp_k[X]|.\]
Then 
\[\phi^\calB_+(b)=\sum_{k=1}^{|X|}\left(\sum_{\substack{\sim\in \calE[b],\\ \cl(\sim)=k}} \frac{|\calC_{\pp}(b/{\sim})|}{k!}g^\calB_{\pp}(b\mid\sim)\right)X^k.\]
\end{theo}

\begin{proof}
Let us consider the polynomial invariant $\phi_{\pp}^\calB$. By Proposition \ref{prop1.3.7}, $\phi_{\pp}^\calB=\phi^\calB_+\curvearrowleft f_{\pp}^\calB$, as $\phi^\calB_+$ is the fundamental polynomial invariant of $\calB$, which is regular.
By Proposition \ref{prop4.1.23}, as $f_{\pp}^\calB(b)=1$ for any $b\in \calB[1]$, $f^\calB_{\pp}$ is invertible for the convolution product $*$. Its inverse is denoted by $g^\calB_{\pp}$. As $\curvearrowleft$ is an action,
\[\phi^\calB_+=\phi_{\pp}^\calB\curvearrowleft g^\calB_{\pp}.\]
Therefore, for any finite set $X$, and any $b\in \calB[X]$,
\[\phi^\calB_+(b)=\sum_{\sim\in \calE[b]} \phi_{\pp}^\calB(b/{\sim})g_{\pp}^\calB(b\mid\sim)=\sum_{\sim\in \calE[b]} \dfrac{|\calC_{\pp}(b/{\sim})|}{|\cl(\sim)|!} g_{\pp}^\calB(b\mid\sim) T^{\cl(\sim)}.\]

Let us now assume that $b$ is indecomposable. Then $X^2=\sim_b^{(\ic)}$ is the unique element of $\calE[b]$ with only one class. Therefore, the coefficient of $X$ in $\phi^\calB_+(b)$ is 
\begin{align*}
\frac{\mathrm{d}\phi_+^\calB(b)}{\mathrm{d}T}(0)&=\dfrac{\left|\calC_{\pp}\left(b/{\sim_b^{(\ic)}}\right)\right|}{\left|\cl\left(\sim_b^{(\ic)}\right)\right|!} g_{\pp}^\calB\left(b\mid\sim_b^{(\ic)}\right)=g_{\pp}^\calB(b). 
\end{align*}
Moreover, by Corollary \ref{cor1.2.26}, $\dfrac{\mathrm{d}\phi_+^\calB(b)}{\mathrm{d}T}(0)=\ln_\calB(\epsilon_\delta)(b)$.
\end{proof}

\begin{remark}
Note that $\ln_\calB(\epsilon_\delta)$ is not a character, but an infinitesimal character. It just coincide with $g_{\pp}^\calB$ on indecomposable elements. The unicity of the decomposition in indecomposable elements (up to the order)
implies that this is a coherent way to define $g_{\pp}^\calB$.
\end{remark}

\begin{notation}
For $\calB=\exch_+$, we shall simply note $\Upsilon_+$ instead of $g_{\pp}^{\exch_+}$. For any finite set $X$, and any $I\in \exch_+[X]$, indecomposable, as $\calC_+(I)=I$,
\[\Upsilon_+(I)=\sum_{k=1}^{|X|}\frac{(-1)^{k-1}}{k} |I\cap \comp_k[X]|.\]
\end{notation}

\begin{remark}
We do not necessarily have $\Upsilon_+\circ \Gamma_+^\calB=g_{\pp}^\calB$. However, this is satisfied if, for any finite set $X$, for any $b\in \calB[X]$, $\ic(\calC_+(b))=\ic(b)$. 
\end{remark}

\section{Quasi-cocommutative CTBs}

We now introduce quasi-cocommutative CTBs. The definition is rather technical, but its interest will appear in the next section, when we will try to introduce a structure of CTDB on $\exch_+$.

\subsection{Definition and examples}

\begin{defi}\label{defi4.2.1}
Let $\calB$ be a CTB. We shall say that it is quasi-cocommutative if it is commutative and if for any finite set $X$, any $b\in \calB[X]$, and any $(X_1,\ldots,X_k)$ and any $(X_{\sigma(1)}),\ldots,X_{\sigma(k)}) \in \calC_-(b)$ with $\sigma \in \sym_k$ such that $b_{\mid_i X_i}$ and $b_{\mid_i X_{\sigma(i)}}$ are indecomposable for any $i\in [k]$, then
\begin{align*}
&\forall i\in [k],& b_{\mid_i X_i}&=b_{\mid_{\sigma^{-1}(i)} X_i}.
\end{align*}\end{defi}

\begin{remark}
If $\calB$ is a cocommutative CTB, then it is quasi-cocommutative.
Therefore, $\calG$, $(\calH\calG,m,\Delta^{\subset,\subset},\Delta^{\cap,\cap})$ are quasi-cocommutative.
\end{remark}

\begin{example}
The CTBs $\calP$ and $\calQ\calP$ are quasi-cocommutative, but not cocommutative, as the coproduct is given by restriction.
\end{example}

\begin{example}\begin{enumerate}
\item The CTB $(\calH\calH,m,\Delta^{\cap,\subset})$ is not quasi-cocommutative. Let $G$ be the hypergraph on $[4]$ whose hyperedges are given by
\[\HE(G)=\{\{1\},\{2\},\{3\},\{4\},\{1,2\},\{1,3\},\{1,2,3,4\}\}.\]
We take $X=[3]$ and $Y=\{4\}$. Then 
\begin{align*}
\HE(G_{\mid X})&=\{\{1\},\{2\},\{3\},\{1,2\},\{1,3\}\},&\HE(G_{\mid_\cap X})&=\{\{1\},\{2\},\{3\},\{1,2\},\{1,3\},\{1,2,3\}\},\\
\HE(G_{\mid Y})&=\{\{4\}\},&\HE(G_{\mid_\cap Y})&=\{\{4\}\}.
\end{align*}
Hence, $G_{\mid X}$, $G_{\mid_\cap X}$, $G_{\mid Y}$, and $G_{\mid_\cap Y}$ are connected, so indecomposable. Moreover, 
\begin{align*}
\Delta_{X,Y}^{\cap,\subset}(G)&=G_{\mid_\cap X}\otimes G_{\mid Y},&\Delta_{Y,X}^{\cap,\subset}(G)&=G_{\mid_\cap Y}\otimes G_{\mid X},
\end{align*} 
and $G_{\mid_\cap X}\neq G_{\mid X}$. 
\item The CTB $\exch_+$ is not quasi-cocommutative. Let 
\[I=\{(1,2,3),(3,1,2),(3,2,1)\}\in \exch_+[3].\]
We take $X_1=\{1,2\}$ and $X_2=\{3\}$. Then 
\begin{align*}
\Delta_{X_1,X_2}(I)&=\{(1,2)\}\otimes \{(3)\},\\
\Delta_{X_2,X_1}(I)&=\{(3)\}\otimes \{(1,2),(2,1)\},
\end{align*}
and $\{(1,2)\}$, $\{(1,2),(2,1)\}$, and $\{(3)\}$ are indecomposable. 
\end{enumerate}\end{example}

\begin{prop}\label{cor4.2.2}
Let $\calB$ be a non-degenerate CTDB.
Then it is quasi-cocommutative. 
\end{prop}

\begin{proof}
Let $X$ be a finite set, and $b\in \calB[X]$. We assume that $(X_1,\ldots,X_k)$ and $(X_{\sigma(1)},\ldots X_{\sigma(k)})$ both belong to $\calC_-(b)$, with $\sigma\in \sym_k$.
Finally, we assume that for any $i\in [k]$, $b_{\mid_i X_i}$ and $b_{\mid {\sigma^{-1}(i)} X_i}$ are indecomposable.
Let us prove that for any $i\in [k]$, $b_{\mid_i X_i}=b_{\mid {\sigma^{-1}(i)} X_i}$.

Let us put $\sim=X_1^2\sqcup\cdots\sqcup X_k^2$. Then $X/{\sim}=\{X_1,\ldots,X_k\}$, and, by Lemma \ref{lem4.2.3}, $\sim\in \calC_-(b)$, and
\[b\mid \sim=b_{\mid_1X_1}\cdots b_{\mid_kX_k}.\]
Similarly,
\[b\mid \sim=b_{\mid_1 X_{\sigma(1)}}\cdots b_{\mid_k X_{\sigma(k)}}.\]
By the unicity of the decomposition into indecomposable elements, for any $i\in [k]$, 
$b_{\mid_iX_i}=b_{\mid {\sigma^{-1}(i)} X_i}$.
\end{proof}

\begin{example} \label{ex4.2.3}
As a consequence, nor $\exch_+$, nor $(\calH\calG,m,\Delta^{\subset,\cap})$ can be made CTDBs, as they are both non-degenerate but not quasi-cocommutative.
\end{example}

\subsection{Two lemmas on non-degenerate CTDBs}

\begin{lemma}\label{lem4.2.3}
Let $\calB$ be a non-degenerate CTDB.
Let $X$ be a finite set, $b\in \calB[X]$ and $\sim\in \calE[X]$. The following assertions are equivalent.
\begin{enumerate}
\item $\sim\in \calE[b]$.
\item There exists $(X_1,\ldots,X_k)\in \calC_-(b)$, such that:
\begin{itemize}
\item $X/{\sim}=\{X_1,\ldots,X_k\}$.
\item Putting $\Delta_{X_1,\ldots,X_k}(b)=b_{\mid_1X_1}\otimes \ldots \otimes b_{\mid_kX_k}$, then for any $i\in [k]$, $b_{\mid_iX_i}$ is indecomposable. 
\end{itemize}
Moreover, if this holds, $b\mid \sim=b_{\mid_1X_1}\cdots b_{\mid_kX_k}$. 
\end{enumerate}\end{lemma}

\begin{proof}
$\Longrightarrow$. As $\calB$ is non-degenerate, $\calC_{\pp}(b/{\sim})\neq \emptyset$: there exists an indexation $\{X_1,\ldots,X_k\}$ of $X/{\sim}$,
such that $(\{X_1\},\ldots,\{X_k\})\in \calC_{\pp}(b/{\sim})$. We obtain
\begin{align*}
&\left(\Delta_{\{X_1\},\ldots,\{X_k\}}\otimes \id_{\calB[X]}\right)\circ \delta_\sim(b)\\
&=(b/{\sim})_{\mid_1\{X_1\}}\otimes \cdots \otimes (b/{\sim})_{\mid_k\{X_k\}}\otimes b\mid\sim\\
&=m_{1,3,\ldots,2k-1,2\ldots 2k}\circ (\delta_{\sim_1}\otimes \cdots \otimes \delta_{\sim_k})\circ \Delta_{X_1,\ldots,X_k}(b)\\
&=(b_{\mid_1X_1})/{\sim}_1\otimes\cdots\otimes(b_{\mid_kX_k})/{\sim}_k\otimes (b_{\mid_1X_1})\mid\sim_1\cdots(b_{\mid_kX_k})\mid \sim_k,
\end{align*}
where $\sim_i=\sim\cap X_i^2=X_i^2$ for $i\in [k]$. As this is nonzero, $(X_1,\ldots,X_k)\in \calC_-(b)$ and, for any $i\in [k]$, $\sim_i\in \calE[b_{\mid_iX_i}]$. So $\sim_i=X_i^2\subseteq \sim_{b_{\mid_iX_i}}^{(\ic)}$,
which gives $\sim_{b_{\mid_iX_i}}^{(\ic)}=X_i^2$, and finally that $b_{\mid_iX_i}$ is indecomposable. Moreover, for $i\in [k]$,
\[(b_{\mid_i X_i})\mid \sim_i=(b_{\mid_i X_i})\mid \sim_{b_{\mid_i X_i}}^{(\ic)}=b_{\mid_i X_i},\]
so
\[b\mid \sim=(b_{\mid_1 X_1})\mid\sim_1\cdots(b_{\mid_kX_k})\mid \sim_k=b_{\mid_1 X_1}\cdots b_{\mid_k X_k}.\]

$\Longleftarrow$. As $b_{\mid_i X_i}$ is indecomposable, $\sim_{b_{\mid_i X_i}}^{(\ic)}=X_i^2\in \calE[b_{\mid_i X_i}]$ for $i\in [k]$. We put $\sim_i=X_i^2$ and 
\[\sim=\sim_1\sqcup \cdots \sqcup \sim_k.\]
Then
\begin{align*}
\left(\Delta_{\{X_1\},\ldots,\{X_k\}}\otimes \id_{\calB[X]}\right)\circ \delta_\sim(b)&=
m_{1,3,\ldots,2k-1,2\ldots 2k}\circ (\delta_{\sim_1}\otimes \cdots \otimes \delta_{\sim_k})\circ \Delta_{X_1,\ldots,X_k}(b)\neq 0,
\end{align*}
so $\delta_\sim(b)\neq 0$, and $\sim\in \calE[b]$. 
\end{proof}

\begin{lemma}\label{lem4.2.4}
Let $\calB$ a non-degenerate CTDB. Let $X$ be a finite set, $\sim\in \calE[b]$, $(X_1/{\sim}_1,\dots,X_k/{\sim}_k)\in \comp[X/{\sim}]$. The following assertions are equivalent:
\begin{enumerate}
\item $(X_1/{\sim}_1,\ldots,X_k/{\sim}_k)\in \calC_+(b/{\sim})$.
\item $(X_1,\ldots,X_k)\in \calC_-(b)$ and $\sim_i=\sim_{b_{\mid_i X_i}}^{(\ic)}$ for any $i\in [k]$.
\end{enumerate}
\end{lemma}

\begin{proof}
Observe that $\sim=\sim_1\sqcup \cdots\sqcup \sim_k$. Moreover,
\begin{align*}
&\left(\epsilon_+^{\otimes k}\otimes \id_{\calB[X]}\right)\circ\left(\Delta_{X_1/{\sim}_1,\ldots,X_k/{\sim}_k}\otimes\id_{\calB[X]}\right)\circ \delta_\sim(b)\\
&=\left(\epsilon_+^{\otimes k}\otimes \id_{\calB[X]}\right)\circ m_{1,3,\ldots,2k-1,2\ldots 2k}\circ (\delta_{\sim_1}\otimes\cdots\otimes\delta_{\sim_k})\circ \Delta_{X_1,\ldots,X_k}(b).
\end{align*}
Then
\begin{align*}
1.\Longleftrightarrow &\left(\epsilon_+^{\otimes k}\otimes \id_{\calB[X]}\right)\circ\left(\Delta_{X_1/{\sim}_1,\ldots,X_k/{\sim}_k}\otimes\id_{\calB[X]}\right)\circ \delta_\sim(b)\neq 0,\\
2.\Longleftrightarrow&\left(\epsilon_+^{\otimes k}\otimes \id_{\calB[X]}\right)\circ m_{1,3,\ldots,2k-1,2\ldots 2k}\circ (\delta_{\sim_1}\otimes\cdots\otimes\delta_{\sim_k})\circ \Delta_{X_1,\ldots,X_k}(b)\neq 0.
\end{align*}
The result follows, as $\epsilon_+=\epsilon_\delta$.
\end{proof}

\subsection{Equivalences in a quasi-cocommutative CTB}

Mimicking Lemma \ref{lem4.2.3}, we can associate a set of equivalences to any element of a quasi-cocommutative CTB.

\begin{defi}\label{defi4.2.5}
Let $\calB$ be a quasi-cocommutative CTB, $X$ be a finite set, and $b\in \calB[X]$.
Let $\sim\in \calE[X]$. We shall say that $\sim\in \calE[b]$ if there exists $(X_1,\ldots,X_k)\in \calC_-(b)$, such that:
\begin{itemize}
\item $X/{\sim}=\{X_1,\ldots,X_k\}$.
\item Putting $\Delta_{X_1,\ldots,X_k}(b)=b_{\mid_1X_1}\otimes \ldots \otimes b_{\mid_kX_k}$, then for any $i\in [k]$, $b_{\mid_iX_i}$ is indecomposable. 
\end{itemize}
If this holds, we put $b\mid \sim=b_{\mid_1X_1}\cdots b_{\mid_kX_k}$. The definition of the quasi-cocommutativity implies that this does not depend on the choice of $(X_1,\ldots,X_k)$. 
\end{defi}

\begin{remark}
If $\calB$ is a non-degenerate CTDB, Lemma \ref{lem4.2.3} insures that this definitions of $\calE[b]$ and $b\mid\sim$ coincide with the preceding ones. 
\end{remark}

It is not easy to define $b/{\sim}$ at this level of generality as an element of $\calB$. However, it is possible to define it as an element of $\exch_+[X/{\sim}]$.
Mimicking Lemma \ref{lem4.2.4}:

\begin{defi}\label{defi4.2.6}
Let $\calB$ be a quasi-cocommutative CTB, $X$ be a finite set, and $b\in \calB[X]$. Let $\sim\in \calE[X]$ and $\pi=(X_1/{\sim_1},\ldots,X_k/{\sim_k})\in \comp[X/{\sim}]$. We shall say $\pi\in\calC_+(b,\sim)$ if $(X_1,\ldots,X_k)\in \calC_-(b)$ and $\sim_i=\sim_{b_{\mid_i X_i}}^{(\ic)}$ for any $i\in [k]$. Then:
\begin{enumerate}
\item $\calC_+(b,\sim)\neq \emptyset$ if, and only if, $\sim\in \calE[b]$.
\item If $\sim\in \calE[b]$, then $\calC_+(b,\sim)\in \exch_+[X/{\sim}]$. 
\end{enumerate}\end{defi}
 
\begin{proof}
We put $I=\calC_+(b,\sim)$, in order to enlighten the writing of this proof. 

1. $\Longrightarrow$. Let us assume that $\sim\in \calE[b]$, and let us first prove that $I\neq \emptyset$. Let $(X_1,\ldots,X_k)$ associated to $\sim$, as in Definition \ref{defi4.2.5}. For $i\in [k]$, as $b_{\mid_iX_i}$ is indecomposable, $\sim^{(\ic)}_{b_{\mid_i X_i}}=X_i^2=\sim\cap X_i^2=\sim$. Therefore, by definition, $(X_1/{\sim_1},\ldots,X_k/{\sim_k})\in I$.

1. $\Longleftarrow$. Let us assume that $I\neq \emptyset$, and let us take $(X_1/{\sim_1},\ldots,X_k/{\sim_k})\in I$. Then $(X_1,\ldots,X_k)\in \calC_-(b)$ and $\sim_i=\sim_{b_{\mid_i X_i}}^{(\ic)}$ for any $i\in [k]$. Let us denote by $X_{i_1},\ldots,X_{i,l_i}$ the connected components of $b_{\mid_i X_i}$, and put $\sim_{i,j}=X_{i,j}$ for any $i\in [k]$ and $j\in [l_i]$. As $\sim_i=\sim_{b_{\mid_i X_i}}$,
the $X_{i,j}$'s are the classes of $\sim$. By construction of the $X_{i,j}$'s,
$\pi=(X_{1,1},\ldots,X_{1,l_1},\ldots,X_{k,1},\ldots,X_{k,l_k})\in \calC_-(b)$ and $\Delta_{\pi}(b)$ is a tensor of indecomposable elements. Therefore, $\sim\in \calE[b]$.

2. Let us now prove the exchange condition for $I$. 
Let $(X'/{\sim}',X''/{\sim}'')\in \comp[X]$, $\pi',\varpi'\in \comp[X']$, $\pi'',\varpi''\in \comp[X'']$, such that $(\pi',\pi'')$ and $(\varpi',\varpi')$ are both elements of $I$. We denote by $\overline{\pi}'$, $\overline{\pi}''$, $\overline{\varpi}'$ and $\overline{\varpi}''$ the lifting of these compositions to compositions of $X'$ or $X''$. 
By definition of $I$, $(X_1,\ldots,X_k)=(\overline{\pi}',\overline{\pi}'')\in \calC_-(b)$, $(Y_1,\ldots,Y_l)=(\overline{\varpi}',\overline{\varpi}'')\in \calC_-(b)$, and:
\begin{itemize}
\item For any $i\in [k]$, $\sim\cap X_i^2=\sim^{(\ic)}_{b_{\mid_iX_i}}$.
\item For any $j\in [l]$, $\sim\cap Y_j^2=\sim^{(\ic)}_{b_{\mid_jY_j}}$.
\end{itemize}
By the exchange condition for $\calC_-(b)$, $(\overline{\varpi}',\overline{\pi}'')\in \calC_-(b)$. Moreover, the factors of $\Delta_{(\overline{\varpi}',\overline{\pi}'')}(b)$ are, either some $b_{\mid_iX_i}$, or some $b_{\mid_jY_j}$.
So the second condition of Definition \ref{defi4.2.6} is satisfied for $(\varpi',\pi'')$, which indeed belong to $I$.

Let $\pi\leq \pi'\in \comp[X/{\sim}]$, with $\pi\in I$. We put $\pi=(X_1/{\sim_1},\ldots,X_k/{\sim_k})$ and
$\pi'=(X'_1/{\sim'_1},\ldots,X'_k/{\sim'_l})$. By hypothesis, $(X_1,\ldots,X_k)\in \calC_-(b)$ and $(X_1,\ldots,X_k)\leq (X'_1,\ldots,X'_l)$. By definition of the order, $\pi=(\pi'_1,\ldots,\pi'_k)$, with $\pi_i\in \comp[X_i]$ for $i\in [k]$.
Let us take $i\in [k]$. As $\sim_i=\sim^{(\ic)}_{b_{\mid_i X_i}}$, any $Y\subseteq X_i$, union of classes of $\sim_i$, is a disjoint union of connected components of $b_{\mid_i X_i}$. 
Therefore, any composition $\varpi$ of $\comp[X_i]$ made of unions ok classes of $\sim_i$ belong to $\calC_-(b_{\mid_i X_i})$, and if $b'\in \calB[Y]$ is a factor of $\Delta_\varpi(b_{\mid_i X_i})$, then $\sim^{(\ic)}_{b'}=\sim_i\cap Y$. 
Therefore, $\varpi'\in \calC_+(b_{\mid_i X_i},\sim_i)$. We obtain in particular that for any $i\in [k]$, $\pi'_i\in \calC_+(b_{\mid_i X_i},\sim_i)$. This implies that $\pi'\in\calC_+(b,\sim)$,
so $\calC_+(b,\sim)\in \exch_+[X/{\sim}]$. \end{proof}

\begin{remark}
By construction, if $\calB$ is a non-degenerate CTDB, for any finite set $X$, and any $b\in \calB[X]$, $\calC_+(b,\sim)=\calC_+(b/{\sim})$.
\end{remark}

\begin{lemma}\label{lem4.2.7}
Let $\calB$ be a quasi-commutative CTB, $X$ and $Y$ be disjoint finite sets, $a\in \calB[X]$ and $b\in \calB[Y]$. Then
\begin{align*}
\calE[ab]&=\{\sim_X\sqcup \sim_Y\mid \sim_X\in \calE[a],\: \sim_Y\in \calE[b]\}.
\end{align*}
Moreover, if $\sim_X\in \calE[a]$ and $\sim_Y\in \calE[b]$, putting $\sim=\sim_X\sqcup \sim_Y$, then
\begin{align*}
\calC_+(ab,\sim)&=\calC_+(a,\sim_X)\squplus \calC_+(b,\sim_Y),&
(ab)\mid\sim&=(a\mid\sim_X) (b\mid\sim_Y).
\end{align*}\end{lemma}

\begin{proof}
$\subseteq$. Let $\sim\in \calE[ab]$. There exists $(X_1,\ldots,X_k)\in \calC_-(ab)$, such that $X/{\sim}=\{X_1,\ldots,X_k\}$, and 
\[\Delta_{X_1,\ldots,X_k}(ab)=(ab)_{\mid_1 X_1}\otimes \cdots \otimes (ab)_{\mid_k X_k},\]
with, for any $i\in [K]$, $(ab)_{\mid_i X_i}$ indecomposable. As $\Delta$ is compatible with the product, for any $i\in [k]$,
\[(ab)_{\mid_i X_i}=a_{\mid_i X_i\cap X}b_{\mid_i X_i\cap Y}.\]
As $(ab)_{\mid X_i}$ is indecomposable, either $X_i\subseteq X$ or $X_i\subseteq Y$. Consequently, $\sim=\sim_X\sqcup \sim_Y$, with $\sim_X=\sim\cap X^2$ and $\sim_Y=\sim\cap Y^2$.
Moreover, if $X_i\subseteq X$, $(ab)_{\mid_i X_i}=a_{\mid_i X_i}$ is indecomposable. This implies that $\sim_X\in \calE[a]$. Similarly, $\sim_Y\in \calE[b]$.

$\supseteq$. Let $\sim_X\in \calE[a]$, $\sim_Y\in \calE[b]$, and $\sim=\sim_X\sqcup \sim_Y$. There exists $(X_1,\ldots,X_k)\in \calC_-(a)$, $(Y_1,\ldots,Y_l)\in \calC_-(b)$, such that
\begin{align*}
X/{\sim}_X&=\{X_1,\ldots,X_k\},&Y/{\sim}_Y&=\{Y_1,\ldots,Y_l\},
\end{align*}
and, for any $i\in [k]$, $a_{\mid_i X_i}$ is indecomposable, and for any $j\in [l]$, $b_{\mid_j Y_j}$ is indecomposable.
Then $(X\sqcup Y)/{\sim}=\{X_1,\ldots,X_k,Y_1,\ldots,Y_l\}$, and
\begin{align*}
\Delta_{X_1,\ldots,X_k,Y_1,\ldots,Y_l}(ab)&=\left(\Delta_{X_1,\ldots,X_k}\otimes \Delta_{Y_1,\ldots,Y_l}\right)\circ \Delta_{X,Y}(ab)\\
&=\left(\Delta_{X_1,\ldots,X_k}\otimes \Delta_{Y_1,\ldots,Y_l}\right)(a\otimes b)\\
&=a_{\mid_1 X_1}\otimes \cdots \otimes a_{\mid_k X_k}\otimes b_{\mid_1 Y_1}\otimes \cdots \otimes b_{\mid_j Y_l}
\end{align*}
so $\sim\in \calE[ab]$.

Moreover, by definition,
\[(ab)\mid \sim=\left(\prod_{i=1}^k a_{\mid_i X_i}\right)\left(\prod_{j=1}^j b_{\mid_j Y_j}\right)=(a\mid \sim_X)(b\mid\sim_Y).\]

We put $I=\calC_+(a,\sim_X)$, $J=\calC_+(b,\sim_Y)$, and $K=\calC_+(ab,\sim)$, in order to enlighten the writing of the end of this proof.

$K\subseteq I\squplus J$. Let $(Z_1/{\sim}_1,\ldots,Z_n/{\sim}_n)\in K$, with $(Z_1,\ldots,Z_n)\in \calC_-(ab)$, and $\sim_i=\sim_{(ab)_{\mid_i Z_i}}^{(\ic)}$ for any $i\in [n]$.
For any $i\in [n]$, we put $Z_i/{\sim}_i=(X_i/{\sim}'_i)\sqcup (Y_i/{\sim}''_i)$, with $X_i\subseteq X$ and $Y_i\subseteq Y$. Then 
\[\Delta_{Z_1,\ldots,Z_l}(ab)=\Delta_{X_1,\ldots,X_n}(a)\Delta_{Y_1,\ldots,Y_n}(b)\neq 0,\]
so $\elim (X_1,\ldots,X_n)\in \calC_-(a)$ and $\elim (Y_1,\ldots,Y_n)\in \calC_-(b)$. Moreover, 
\[\sim_i=\sim_{(ab)_{\mid_i Z_i}}^{(\ic)}=\sim_{(a_{\mid_i X_i})(b_{\mid_i Y_i})}^{(\ic)}=\sim_{a_{\mid_i X_i}}^{(\ic)}\sqcup \sim_{b_{\mid_i Y_i}}^{(\ic)},\]
so $\sim'_i=\sim_{a_{\mid_i X_i}}^{(\ic)}$ and $\sim''_i=\sim_{b_{\mid_i Y_i}}^{(\ic)}$. This gives that $\elim (X_1,\ldots,X_n)\in a$ and $\elim (Y_1,\ldots,Y_n)\in b$,
so finally $(Z_1/{\sim}_1,\ldots,Z_n/{\sim}_n)\in I\squplus J$.

$I\squplus J\subseteq K$. We keep the same notations. We obtain that
\[\Delta_{Z_1,\ldots,Z_l}(ab)=\Delta_{X_1,\ldots,X_n}(a)\Delta_{Y_1,\ldots,Y_n}(b)\neq 0,\]
as $\Delta_{X_1,\ldots,X_n}(a)$ and $\Delta_{Y_1,\ldots,Y_n}(b)$ are both tensor products of elements of $\calB$, and the product of two elements of $\calB$ is never zero.
Moreover, for any $i\in [l]$, 
\[\sim_i=\sim'_i\sqcup \sim''_i=\sim_{a_{\mid_i X_i}}^{(\ic)}\sqcup \sim_{b_{\mid_i Y_i}}^{(\ic)}=\sim_{(ab)_{\mid_i Z_i}}^{(\ic)},\]
so $(Z_1,\ldots,Z_n)\in K$.
\end{proof}

\begin{lemma}\label{lem4.2.8}
Let $\calB$ be a quasi-commutative CTB, $X$ be a finite set, and $b\in \calB[X]$. 
\begin{enumerate}
\item $\sim^{(\ic)}_b\in \calE[b]$. Moreover, for any $\sim\in \calE[b]$, the following assertions are equivalent:
\begin{enumerate}
\item $\calC_+(b,\sim)=\comp[X/{\sim}]$.
\item $b\mid\sim=b$.
\item $\sim=\sim^{(\ic)}$.
\end{enumerate}
\item $\sim_b^{(m)}\in \calE[b]$ if, and only if, $\calC_+(b)\neq \emptyset$. Moreover, for any $\sim\in \calE[b]$, the following assertions are equivalent:
\begin{enumerate}
\item $\calC_+(b,\sim)=\calC_+(b)$.
\item $b\mid\sim\in \bfD\calB[X]$.
\item $\sim=\sim_b^{(m)}$.
\end{enumerate}\end{enumerate}\end{lemma}

\begin{proof}
1. Let $b=b_1\cdots b_k$ be the decomposition of $b$ into indecomposable elements, with $b_i\in \calB[X_i]$ for $i\in [k]$. Then
\[\Delta_{X_1,\ldots,X_k}(b)=b_1\otimes\cdots\otimes \cdots b_k,\]
and, as $b_i$ is indecomposable for any $i\in [k]$, $\sim_b^{(\ic)}\in \calE[b]$. Moreover, $b\mid\sim=b_1\cdots b_k=b$.

We have just proved that $(c)\Longrightarrow (b)$. Let us now prove that $(b)\Longrightarrow (c)$. If $\sim\in \calE[b]$, by construction of $b\mid\sim$, the connected components of $b\mid \sim$ are the classes of $\sim$. Therefore, if $b\mid\sim=b$, then $\sim=\sim_b^{(\ic)}$.

$(a)\Longrightarrow (c)$. If $\calC_+(b,\sim)=\comp[X/{\sim}]$, then $(x/{\sim})\in \calC_+(b,\sim)$, so $(X)\in \calC_-(b)$ (which is obvious) and $\sim=\sim_b^{(\ic)}$, by definition of $\calC_+(b,\sim)$.

$(c)\Longrightarrow (a)$. Let $(X_1/{\sim_1},\ldots,X_k/{\sim_k})\in \comp[X/{\sim}]$.
Then any $X_i$ is a union of connected components of $b$, so $(X_1,\ldots,X_k)\in \calC_-(b)$ and $b_{\mid_iX_i}$ is the product of the connected components of $b$ included in $X_i$. Therefore, $\sim_i=\sim^{(\ic)}_b\cap X_i^2=\sim^{(\ic)}_{b_{\mid_iX_i}}$. As a conclusion, $(X_1/{\sim_1},\ldots,X_k/{\sim_k})\in \calC_+(b,\sim)$.

2. Let us assume that $\sim_b^{(m)}\in \calE[b]$. There exists $(X_1,\ldots,X_k)\in \calC_-(b)$, such that $X/{\sim_b^{(m)}}=\{X_1,\ldots,X_k\}$. So the $X_i$'s are singletons: $(X_1,\ldots,X_k)\in \calC_+(b)$.

Let us assume that $\calC_+(b)\neq \emptyset$. As $\calC_+(b)\in \exch_+[X]$, it contains a maximal element $(X_1,\ldots,X_k)$ of $\comp[X]$, that is to say $|X_i|=1$ for any $i\in [k]$, which implies that $b_{\mid_iX_i}$ is indecomposable. Therefore, the equivalence whose classes are $X_1,\ldots,X_k$, which is $\sim_b^{(m)}$, belongs to $\calE[b]$.

$(c)\Longleftrightarrow (b)$. As $\sim_{b\mid \sim}^{(\ic)}=\sim$, 
\[b\mid\sim\in \bfD\calB[X]\Longleftrightarrow \sim_{b\mid \sim}^{(\ic)}=\sim_b^{(m)}\Longleftrightarrow \sim=\sim_b^{(m)}.\]

$(a)\Longrightarrow (c)$. If $\calC_+(b,\sim)=\calC_+(b)$, then $X/{\sim}=X$, so $\sim=\sim_b^{(m)}$.

$(c)\Longrightarrow (a)$. By definition, the elements of $\calC_+\left(b,\sim_b^{(m)}\right)$ are the $(X_1,\ldots,X_k)\in \calC_-(b)$, such that 
\begin{align*}
&\forall i\in [k],&\sim_b^{(m)}\cap X_i^2&=\sim_{b_{\mid_i X_i}}^{(\ic)},
\end{align*}
that is to say, such that $b_{\mid_i X_i}\in \bfD\calB[X_i]$ for any $i\in i\in [k]$. In other words, $\calC_+\left(b,\sim_b^{(m)}\right)=\calC_+(b)$. 
\end{proof}

\begin{lemma}\label{lem4.2.9}
Let $X$, $Y$ be disjoint finite sets, $b\in\calB[X\sqcup Y]$, $\sim' \in \calE[X]$ and $\sim''\in \calE[Y]$. We put $\sim=\sim'\sqcup \sim''$. The following assertions are equivalent:
\begin{enumerate}
\item $\sim\in \calE[b]$, and $(X/{\sim}',Y/{\sim}'')\in m(\calC_+(b,\sim))$.
\item $(X$, $Y)\in \calC_-(b)$, and $\sim'\in \calE[b_{\mid_LX}]$ and $\sim''\in \calE[b_{\mid_RY}]$.
\end{enumerate}
If so, $\Delta_{X/{\sim}',Y/{\sim}''}(\calC_+(b,\sim))=\calC_+\left(b_{\mid_LX},\sim'\right)\otimes \calC_+\left(b_{\mid_RY},\sim''\right)$, and $b\mid\sim=(b_{\mid_LX}\mid\sim')(b_{\mid_RY}\mid \sim'')$.
\end{lemma}

\begin{proof}
$2.\Longrightarrow 1$. By definition of $\calE[b_{\mid_LX}]$ and $\calE[b_{\mid_RY}]$, there exist $(X_1,\ldots, X_k)\in \calC_-(b_{\mid_LX})$, $(Y_1,\ldots,Y_l)\in \calC_-(b_{\mid_RY})$, such that $\Delta_{X_1,\ldots,X_k}(b_{\mid_LX})$ and $\Delta_{Y_1,\ldots,Y_l}(b_{\mid_RY})$ are tensors of indecomposable elements, and
\begin{align*}
X/{\sim}'&=\{X_1,\ldots,X_k\},&Y/{\sim}''&=\{Y_1,\ldots,Y_k\}.
\end{align*} 
As $\Delta_{X_1,\ldots,X_k,Y_1,\ldots,Y_l}(b)=\Delta_{X_1,\ldots,X_k}(b_{\mid_LX})\otimes \Delta_{Y_1,\ldots,Y_l}(b_{\mid_RY})$, this is a tensor of indecomposable elements, and
\[(X\sqcup Y)/{\sim}=(X/{\sim}')\sqcup (Y/{\sim}'')=\{X_1,\ldots,X_k,Y_1,\ldots,Y_k\}.\]
So $\sim\in \calE[b]$. Moreover, as $(X_1,\ldots,X_k,Y_1,\ldots,Y_l)\in \calC_-(b)$, we obtain an element of $\calC_+(b,\sim)$ of the form $(\pi',\pi'')$, with $\pi'\in \comp[X/{\sim}]$ and $\pi''\in \comp[Y/{\sim}'']$. As $(X/{\sim}',Y/{\sim}'')\leq (\pi,\pi')$, $(X/{\sim}',Y/{\sim}'')\in m(\calC_+(b,\sim))$.

$1.\Longrightarrow 2$. As $(X/{\sim}',Y/{\sim}'')\in m(\calC_+(b,\sim))$, there exists $\pi\in \calC_+(b,\sim)$, such that $(X/{\sim}',Y/{\sim}'')\leq \pi$. By definition of the order, $\pi=(\pi',\pi'')$, with $\pi'\in \comp[X]$ and $\pi''\in \comp[Y]$. 
We denote by $\varpi'\in \comp[X]$ and $\varpi''\in \comp[Y]$ the liftings of $\pi'$ and $\pi''$. Then $(\varpi',\varpi'')\in \calC_-(b)$. As $\calC_-(b)\in \exch_-[X]$ and $(X$, $Y)\leq (\varpi',\varpi'')$, $(X$, $Y)\in \calC_-(b)$.

As $(X/{\sim}',Y/{\sim}'')\in m(\calC_+(b,\sim))$, there exists $\pi'\in \comp[X/{\sim}']$, $\pi''\in \comp[Y/{\sim}'']$, such that $(\pi',\pi'')\in \calC_+(b,\sim)$. We denote by $\varpi'\in \comp[X]$ and $\varpi''\in\comp[Y]$ the liftings of, respectively, $\pi$ and $\pi'$.
By definition of $\calC_+(b,\sim)$, $(\varpi',\varpi'') \in \calC_-(b)$. Moreover, as $(X$, $Y)\leq (\varpi',\varpi'')$, $(X$, $Y)\in \calC_-(b)$.

We put 
\[\Delta_{X/{\sim}',Y/{\sim}'}(\calC_+(b,\sim))=J'\otimes J''.\]
$J'\subseteq \calC_+\left(b_{\mid_LX},\sim'\right)$. Let $\pi'\in J'$. There exists $\pi''\in \comp[Y/{\sim}]$, such that $(\pi',\pi'')\in \calC_+(b,\sim)$. Let us denote by $\varpi'$ and $\varpi''$ the liftings of $\pi'$ and $\pi''$. Then $(\varpi',\varpi'')\in \calC_-(b)$. Hence, $\varpi'\in \calC_-(b_{\mid_LX})$, and, as the restriction ${b_{\mid_LX}}_{\mid_i X_i'}$ of $b$
to an element $X'_i\in X/{\sim}$ are also restrictions of $b$, we obtain that the restriction of $\sim'$ to $X'_i$ is $\sim_{{b_{\mid_LX}}_{\mid_i X_i'}}^{(\ic)}$: therefore, $\pi'\in \calC_+\left(b_{\mid_LX},\sim'\right)$.

$\calC_+\left(b_{\mid_LX},\sim'\right)\subseteq J'$. Let $\pi\in b_{\mid_LX}/{\sim}$. We denote by $\varpi'\in \comp[X]$ the lifting of $\pi'$. Then $\varpi'\in \calC_-(b_{\mid_LX})$. Taking any $\varpi''\in \calC_-(b_{\mid_RY})$, then $(\varpi',\varpi'')\in \calC_-(b)$. 
As $(X$, $Y)\in \calC_-(b)$, if $(Y_1,\ldots,Y_l)$ is a composition such that the indecomposable components of $b_{\mid_LX}'$ are $Y_1,\ldots,Y_l$, then $(X$, $Y_1,\ldots,Y_l)\in \calC_-(b)$. 
By the exchange condition for $\calC_-(b)$, $(\varpi',Y_1,\ldots,Y_l)\in \calC_-(b)$. As $\sim''\in \calE[b_{\mid_LX}']$, the $Y_j$'s are $\sim''$-saturated, by Lemma \ref{lem4.2.7}. Therefore, there exists $\pi''\in \comp[Y/{\sim}]$, such that the lifting of $\pi''$ is $(Y_1,\ldots,Y_l)$.
By definition of $\calC_+(b,\sim)$, $(\pi',\pi'')\in \calC_+(b,\sim)$, so $\pi'\in J'$. We prove similarly that $\calC_+\left(b_{\mid_RY},\sim''\right)=J''$.

$b\mid\sim=(b_{\mid_LX}\mid\sim')(b_{\mid_RY}\mid \sim'')$. Let $(X_1,\ldots,X_k)\in \calC_-(b_{\mid_LX})$ and $(Y_1,\ldots,Y_l)\in \calC_-(b_{\mid_RY})$, such that
\begin{align*}
X/{\sim}'&=\{X_1,\ldots,X_k\},&Y/{\sim}''&=\{Y_1,\ldots,Y_k\}.
\end{align*}
We put
\begin{align*}
\Delta_{X_1,\ldots,X_k}(b_{\mid_LX})&=b'_1\otimes\cdots \otimes b'_k,&\Delta_{Y_1,\ldots,Y_l}(b_{\mid_RY})&=b''_1\otimes\cdots \otimes b''_l.
\end{align*}
Then, by definition, $b_{\mid_LX}\mid\sim'=b'_1\cdots b'_k$ and $b_{\mid_RY}\mid\sim''=b''_1\cdots b''_l$. Moreover,
\[\Delta_{X_1,\ldots,X_k,Y_1,\ldots,Y_l}(b)=\Delta_{X_1,\ldots,X_k}(b_{\mid_LX})\otimes \Delta_{Y_1,\ldots,Y_l}(b_{\mid_RY})=b'_1\otimes\cdots \otimes b'_k\otimes b''_1\otimes\cdots \otimes b''_l.\]
By definition, $b\mid\sim=b'_1\cdots b'_k b''_1\cdots b''_l=(b_{\mid_LX}\mid\sim') (b_{\mid_RY}\mid\sim'')$.
\end{proof}

\begin{theo}\label{theo4.2.10}
Let $\calB$ be a quasi-cocommutative CTB. For any finite set $X$, any $\sim\in \calE[X]$, and any $b\in \calB[X]$, we define $\xi_\sim(b)\in \exch_+[X/{\sim}]\otimes \calB[X]$ by
\[\xi_\sim(b)=\begin{cases}
\calC_+(b,\sim)\otimes b\mid \sim\mbox{ if }\sim\in \calE[b],\\
0\mbox{ otherwise}.
\end{cases}\]
Then:
\begin{enumerate}
\item For any pair $(X,Y)$ of disjoint finite sets, for any $a\in \calB[X]$ and $b\in \calB[Y]$, for any $\sim\in \calE[X\sqcup Y]$, putting $\sim_X=\sim\cap X^2$ and $\sim_Y\cap Y^2$,
\[\xi_\sim(ab)=\begin{cases}
\xi_\sim(a)\xi_\sim(b)\mbox{ if }\sim=\sim_X\sqcup \sim_Y,\\
0\mbox{ otherwise}.
\end{cases}\]
\item For any pair $(X,Y)$ of disjoint finite sets, for any $\sim_X\in \calE[X]$ and $\sim_Y\in \calE[Y]$,
\[\left(\Delta_{X/{\sim_X},Y/{\sim_Y}}\otimes \id_{\calB[X\sqcup Y]}\right)\circ \xi_{\sim_X\sqcup \sim_Y}
=m_{1,3,24}\circ (\xi_{\sim_X}\otimes \xi_{\sim_Y})\circ \Delta_{X,Y}.\]
\item For any finite set $X$, and any $\sim,\sim'\in \calE[X]$, such that we do not have $\sim\subseteq \sim'$,
\[\left(\id_{\exch_+[X/{\sim'}]}\otimes \xi_\sim\right)\circ \xi_{\sim'}=0.\]
\item For any finite set $X$, and any $b\in \calB[X]$,
\begin{align*}
\sum_{\sim\in \calE[X]} \left(\id_{\exch_+[X/{\sim}]}\otimes \epsilon_+^\calB\right)\circ \xi_\sim(b)&=\begin{cases}
\calC_+(b)\mbox{ if }\calC_+(b)\neq \emptyset,\\
0\mbox{ otherwise}.
\end{cases}\\
\sum_{\sim\in \calE[X]} \left(\epsilon_+\otimes \id_{\calB[X]}\right)\circ \xi_\sim(b)&=b.
\end{align*}
\end{enumerate}\end{theo}

\begin{proof}
1. This is a direct consequence of Lemma \ref{lem4.2.7}.

2. This is a direct consequence of Lemma \ref{lem4.2.9}.

3. If $\sim'\notin \calE[b]$, then $\xi_{\sim'}(b)=0$. Otherwise, let us decompose $b\mid\sim'$ into indecomposable elements:
\[b\mid\sim'=b_1\cdots b_k,\]
with $b_i\in \calB'[X_i]$ for $i\in [k]$. By definition of $\calE[b]$ and $b\mid\sim'$, $X_1,\ldots,X_k$ are the classes of $\sim$.
As $\sim$ is not included in $\sim'$, $\sim$ is not included in $X_1^2\sqcup\cdots\sqcup X_k^2$.
By Lemma \ref{lem4.2.7}, $\sim\notin \calE[b\mid\sim']$, so $\xi_\sim(b\mid\sim')=0$, which gives the result.

4. Both formulas come from Lemma \ref{lem4.2.8}, and the fact that $\epsilon_+(I)=1$ if, and only if, $I=\comp[X]$, for any $I\in \exch_+[I]$, and $\epsilon_+^\calB(a)=1$ if, and only if, $a\in \bfD\calB[X]$ for any $b\in \calB[X]$.
\end{proof}

\begin{remark}
If $\calB$ is a CTDB, then it is quasi-cocommutative and $\xi=\left(\Gamma_\calB^+\otimes \id_{\calB}\right)\circ \delta$, by Lemmas \ref{lem4.2.3} and \ref{lem4.2.4}.
\end{remark}

Not all the morphisms between two quasi-cocommutative CTBs are compatible with $\xi$, as shown in the next result.

\begin{prop}
Let $\calB$ and $\calB'$ be two quasi-cocommutative CTBs, and let $\phi:\calB\longrightarrow \calB'$ be a morphism of CTBs. We assume that for any finite set $X$, for any $b\in \calB[X]$, $\phi(x)\neq 0$. The following assertions are equivalent:
\begin{enumerate}
\item For any finite set $X$, for any $\sim\in \calE[X]$, for any $b\in \calB[X]$,
\[\xi_\sim\circ \phi(b)=\left(\id_{\exch_+[X/{\sim}]}\otimes \phi\right)\circ \xi_\sim(b).\]
\item For any finite set $X$, for any $\sim\in \calE[X]$, $\ic(\phi(b))=\ic(b)$.
\end{enumerate}\end{prop}

\begin{proof}
$1.\Longrightarrow 2$. By Proposition \ref{prop2.2.6}, $\ic(\phi(b))\geq \ic(b)$. Let us take $\sim=\sim_b^{(\ic)}$.
Then $\delta_\sim(b)=\calC_-(b,\sim)\otimes b$, so 
\[\delta_\sim(\phi(b))=\calC_-(b,\sim)\otimes \phi(b)\neq 0.\]
As a consequence, $\sim\in \calE[\phi(b)]$. Therefore, there exists $(X_1,\ldots,X_k)\in \calC_-(\phi(b))$,
such that $X/{\sim}=\{X_1,\ldots,X_k\}$ and for any $i\in [k]$, $\phi(b)_{\mid_i X_i}$ is indecomposable. 
Therefore, by multiplicativity of the coproduct of $\calB'$, for any $i\in [k]$, $X_i$ is included in an indecomposable component of $\phi(b)$. As a consequence, $\sim\subseteq \sim_{\phi(b)}^{(\ic)}$, so 
\[\ic(b)=\cl(\sim)\geq \cl\left(\sim_{\phi(b)}^{(\ic)}\right)=\ic(\phi(b)).\]

$2.\Longrightarrow 1$. Let $\sim\in \calE[X]$. If $\sim\in \calE[b]$, there exists $(X_1,\ldots,X_k)\in \calC_-(b)$,
such that $X/{\sim}=\{X_1,\ldots,X_k\}$ and, $b_{\mid_i X_i}$ is indecomposable for any $i\in [k]$.
As $\phi$ is a morphism of CTBs, 
\[\Delta_{X_1,\ldots,X_k}(\phi(b))=\phi\left(b_{\mid_1 X_1}\right)\otimes\cdots\otimes\phi\left(b_{\mid_k X_k}\right).\]
As $\phi$ does not take the value $0$ on $\calB$, this is non-zero, so $(X_1,\ldots,X_k)\in \calC_-(\phi(b))$. 
Moreover, for any $i\in [k]$, $\phi(b)_{\mid_iX_i}=\phi\left(b_{\mid_i X_i}\right)$ is indecomposable, as $\phi$ conserves $\ic$. Therefore, $\sim\in \calE[\phi(b)]$ and
\[\phi(b)\mid\sim=\phi(b)_{\mid_1 X_1}\cdots \phi(b)_{\mid_k X_k}
=\phi\left(b_{\mid_1 X_1}\right)\cdots \phi\left(b_{\mid_k X_k}\right)
=\phi\left(b_{\mid_1 X_1}\cdots b_{\mid_k X_k}\right)=\phi(b\mid \sim).\]
Conversely, if $\sim\in \calE[\phi(b)]$, there exists $(X_1,\ldots,X_k)\in \calC_-(\phi(b))$, such that $X/{\sim}=\{X_1,\ldots,X_k\}$ and, for any $i\in [k]$, $\phi(b)_{\mid_i X_i}$ is indecomposable. 
\[\Delta_{X_1,\ldots,X_k}(\phi(b))=\phi\left(b_{\mid_1 X_1}\right)\otimes\cdots\otimes\phi\left(b_{\mid_k X_k}\right)\neq 0,\]
so necessarily $(X_1,\ldots,X_k)\in \calC_-(b)$. Moreover, for any $i\in [k]$, $\phi(b)_{\mid_iX_i}=\phi\left(b_{\mid_i X_i}\right)$. As $\ic\left(\phi\left(b_{\mid_i X_i}\right)\right)\leq \ic\left(\phi(b)_{\mid_iX_i}\right)=1$,
$\phi\left(b_{\mid_i X_i}\right)$ is indecomposable, so $\sim\in \calE[b]$. We obtain that $\calE[b]=\calE[\phi(b)]$.

Let us now assume that $\sim\in \calE[b]=\calE[\phi(b)]$. Let $(X_1/{\sim_1},\ldots,X_k/{\sim_k})\in \comp[X/{\sim}]$.
As $\phi$ conserves $\ic$, and as for any $b\in \calB[X]$, $\sim_{\phi(b)}^{(\ic)}\subseteq \sim_b^{(\ic)}$, necessarily
$\sim_{\phi(b)}^{(\ic)}=\sim_b^{(\ic)}$. Therefore, 
\begin{align*}
(X_1/{\sim_1},\ldots,X_k/{\sim_k})\in \calC_+(b,\sim)&\Longleftrightarrow (X_1,\ldots,X_k)\in \calC_-(b)\mbox{ and }
\forall i\in [k],\: \sim_i=\sim_b^{(\ic)}\\
&\Longleftrightarrow (X_1,\ldots,X_k)\in \calC_-(\phi(b))\mbox{ and }
\forall i\in [k],\: \sim_i=\sim_{\phi(b)}^{(\ic)}\\
&\Longleftrightarrow (X_1/{\sim_1},\ldots,X_k/{\sim_k})\in \calC_+(\phi(b),\sim).
\end{align*}
We use the fact that, as $\phi$ does not take the value $0$ on $\calB$, $\calC_-(\phi(b))=\calC_-(b)$ for any $b\in \calB[X]$. As a conclusion, $\calC_-(b,\sim)=\calC_-(\phi(b),\sim)$, and
\begin{align*}
\xi_\sim\circ \phi(b)&=\begin{cases}
\calC_-(\phi(b),\sim)\otimes \phi(b)\mid\sim \mbox{ if }\sim\in \calE[\phi(b)],\\
0\mbox{ otherwise}
\end{cases}\\
&=\begin{cases}
\calC_-(b,\sim)\otimes \phi(b\mid\sim) \mbox{ if }\sim\in \calE[b],\\
0\mbox{ otherwise}
\end{cases}\\
&=\left(\id_{\exch_+[X/{\sim}]}\otimes \phi\right)\circ \xi_\sim(b). \qedhere
\end{align*}\end{proof}

\subsection{Convolution and applications}

\begin{prop}\label{prop4.2.12}
Let $\calB$ be a quasi-cocommutative CTB. We put $\rho_\calB=\id_\calB-\iota_\calB$.
For any $k\geq 1$, for any finite nonempty set $X$ and any $b\in \calB[X]$,
\[\rho_\calB^{\star k}(b)=\sum_{\sim\in \calE[b]} |\calC_+(b,\sim)\cap \comp_k[X/{\sim}]|b\mid\sim.\]
\end{prop}

\begin{proof}
We proceed by induction on $k$. If $k=1$, observe that for any $I\in \exch_+[X]$,
\[\epsilon_+(I)=|I\cap \comp_1[X]|.\]
Therefore, with Theorem \ref{theo4.2.10}-4,
\[\sum_{\sim\in \calE[b]} |\calC_+(b,\sim)\cap \comp_1[X/{\sim}]|b\mid\sim=\sum_{\sim\in \calE[X]}\left(\epsilon_+\otimes \id_{\calB[X]}\right)\circ \xi_\sim(b)=b=\rho_\calB(b),\]
the last equality coming from the fact that $X$ is nonempty. Let us assume the result at rank $k$.
\begin{align*}
&\rho_\calB^{\star(k+1)}(b)\\
&=\sum_{X=X'\sqcup X''}m_{X',X''}\circ (\rho_\calB^{\star k}\otimes \rho_\calB)\circ \Delta_{X',X''}(b)\\
&=\sum_{(X',X'')\in \calC_-(b)}m_{X',X''}\circ (\rho_\calB^{\star k}\otimes \rho_\calB)\circ \Delta_{X',X''}(b)+0\\
&=\sum_{\substack{(X',X'')\in \calC_-(b),\\ \sim'\in \calE[b_{\mid_LX'}],\\ \sim''\in \calE[b_{\mid_RX''}]}}
|\calC_+(b_{\mid_LX'},\sim')\cap \comp_k[X'/{\sim'}]|\times |\calC_+(b_{\mid_RX''},\sim'')\cap \comp_1[X''/{\sim''}]| 
\left(b_{\mid_LX'}\right)\mid\sim' \left(b_{\mid_RX''}\right)\mid\sim''\\
&=\sum_{\substack{\sim\in \calE[b],\\ (X'/{\sim'},X''/{\sim''})\in m(\calC_+(b,\sim))}}
|(\calC_+(b,\sim))_{\mid_LX'/{\sim'}}\cap \comp_k[X'/{\sim'}]|\times |(\calC_+(b,\sim))_{\mid_RX''/{\sim}}\cap \comp_1[X''/{\sim}'']|b\mid\sim.
\end{align*}
We used the induction hypothesis and the case $k=1$ for the third equality, Theorem \ref{theo4.2.10}-2 for the fourth one. Moreover, by definition of the coproduct on $\exch$, for any $I\in \exch[X]$, for any $k,l\geq 1$,
\[\sum_{(X',X'')\in m(I)}|I_{\mid_LX'}\cap \comp_k[X']|\times |I_{\mid_RX''}\cap\comp_l[X'']|=|I\cap \comp_{k+l}[X]|,\]
as any element of $\comp_{k+l}[X]$ can be uniquely written as the concatenation of a composition of length $k$ with a composition of length $l$. Therefore,
\[\rho_\calB^{\star(k+1)}(b)=\sum_{\sim\in \calE[b]}|\calC_+(b,\sim)\cap \comp_{k+1}[X/{\sim}]| b\mid\sim. \qedhere\]
\end{proof}

\begin{cor}\label{cor4.2.13}
Let $\calB$ be a quasi-cocommutative CTB, $X$ be a finite nonempty set, and $b\in \calB[X]$.
\begin{enumerate}
\item The antipode $S_\calB$ of $\calB$ is given as follows:
\[S_\calB(b)=\sum_{\sim\in \calE[b]} \Phi_+(\calC_+(b,\sim))(-1) b\mid \sim=\sum_{\sim\in \calE[b]}\left(\sum_{k=1}^{|X/{\sim}|}(-1)^k \left|\calC_+(b,\sim)\cap \comp_k[X/{\sim}]\right|\right) b\mid \sim.\]
\item The eulerian idempotent $\ln(\id_\calB)$ is given as follows:
\[\ln(\id_\calB)(x)=\sum_{\sim\in \calE[b]} \frac{\mathrm{d}\Phi_+(\calC_+(b,\sim))}{\mathrm{d}T}(0) b\mid\sim.\]
\end{enumerate}\end{cor}

\begin{proof}
Let $P=\displaystyle \sum_{k=0}a_kT^k\in \K[[T]]$. Then, by Proposition \ref{prop4.2.12}, 
\[P(\rho_\calB)=\sum_{k=1}^{|X|} a_k \rho_\calB^{\star k}(b)=\sum_{\sim\in \calE[X]} \left(\sum_{k=1}^{|X|} a_k|\calC_+(b,\sim)\cap \comp_k[X]|\right)b\mid\sim.\]
Recall that for any $I\in \exch_+[X]$,
\[\Phi_+(I)=\sum_{k=1}^{|X|} |I\cap \comp_k[X]| \binom{T}{k}.\]
Consequently, 
\begin{align*}
\Phi_+(I)(-1)&=\sum_{k=1}^{|X|}(-1)^k |I\cap \comp_k[X]|,&
\frac{\mathrm{d}\Phi_+}{\mathrm{d}T}(0)&=\sum_{k=1}^{|X|}\frac{(-1)^{k-1}}{k} |I\cap \comp_k[X]|.
\end{align*}
Applying these observations for $P=\dfrac{1}{1+T}$, that is to say $a_k=(-1)^k$ for any $k$, we obtain the formula for $P(\rho_\calB)=S_\calB$, by Takeuchi's formula. Taking $P=\ln(1+T)$, that is to say $a_k=\dfrac{(-1)^{k-1}}{k}$ for any $k\geq 1$, we obtain the formula for $\ln(\id_\calB)$. 
\end{proof}

\subsection{Polynomial invariants}

\begin{prop}\label{prop4.2.14}
Let $\calB$ be a quasi-cocommutative CTB. Let $X$ be a finite set, and $b\in \calB[X]$. 
\[\phi_+^\calB(b)(TT')=\sum_{\sim\in \calE[b]} \left(\Phi_+\otimes \phi_+^\calB\right)\circ \xi_\sim(b)(T,T').\]
\end{prop}

\begin{proof}
Let us prove that for any $n\geq 1$, 
\[\phi_+^\calB(b)(nT)=\sum_{\sim\in \calE[b]} \left(\Phi_+\otimes \phi_+^\calB\right)\circ \xi_\sim(b)(n,T).\]
For $n=1$,
\begin{align*}
\sum_{\sim\in \calE[b]} \left(\Phi_+\otimes \phi_+^\calB\right)\circ \xi_\sim(b)(n,T)&=\sum_{\sim\in \calE[b]} \left(\epsilon_+\otimes \phi_+^\calB\right)\circ \xi_\sim(b)=\phi_+^\calB(T).
\end{align*}
Let us assume the result at rank $n$. 
\begin{align*}
&\phi_+^\calB(b)((n+1)T')\\
&=\phi_+^\calB(b)(nT'+1T')\\
&=\sum_{X=X'\sqcup X''}(\phi_+^\calB\otimes \phi_+^\calB)\circ \Delta_{X',X'}(b)(nT',1T')\\
&=\sum_{X=X'\sqcup X''}\sum_{\substack{\sim'\in \calE[X'],\\ \sim''\in \calE[X'']}}(\Phi_+\otimes\phi_+^\calB\otimes \Phi_+\otimes \phi_+^\calB)\circ (\delta_{\sim'}\otimes \delta_{\sim''})\circ \Delta_{X',X''}(b)(n,T',1,T')\\
&=\sum_{X=X'\sqcup X''}\sum_{\substack{\sim'\in \calE[X'],\\ \sim''\in \calE[X'']}}m_{1,3,24}\circ (\Phi_+\otimes\phi_+^\calB\otimes \Phi_+\otimes \phi_+^\calB)\circ (\delta_{\sim'}\otimes \delta_{\sim''})\circ \Delta_{X',X''}(b)(n,1,T')\\
&=\sum_{X=X'\sqcup X''}\sum_{\substack{\sim'\in \calE[X'],\\ \sim''\in \calE[X'']}} (\Phi_+\otimes\Phi_+\otimes \phi_+^\calB)\circ m_{1,3,24}\circ (\delta_{\sim'}\otimes \delta_{\sim''})\circ \Delta_{X',X''}(b)(n,1,T')\\
&=\sum_{\sim \in \calE[X]}\sum_{X/{\sim}=X'\sqcup X''}(\Phi_+\otimes\Phi_+\otimes \phi_+^\calB)\circ \left(\Delta_{X',X''}\otimes \id_{\calQ[X]}\right)\circ \delta_\sim(b)(n,1,T')\\
&=\sum_{\sim \in \calE[X]}(\Phi_+\otimes \phi_+^\calB)\circ \delta_\sim(b)(n+1,T').
\end{align*}
We used the induction hypothesis and the case $n=1$ for the second equality. As this equality is true for any $n\geq 1$, we obtain the result.
\end{proof}

\begin{lemma}
Let $\calB$ be a quasi-cocommutative CTB. Let $\psi$ be a polynomial invariant of $\exch_+$ and $f\in \Char(\calB)$. The following is a polynomial invariant of $\calB$: for any finite set $X$,
\begin{align*}
&\forall b\in \calB[X],&\psi\leftsquigarrow f(b)&=\sum_{\sim\in \calE[b]} (\psi\otimes f)\circ \xi_\sim(b).
\end{align*}\end{lemma}

\begin{proof}
Let $(X,Y)$ be a pair of disjoint finite sets, $a\in \calB[X]$, $b\in \calB[Y]$. By the compatibility between $\xi$ and $m$,
\begin{align*}
\psi\leftsquigarrow f(ab)&=\sum_{\substack{\sim_X\in \calE[a],\\ \sim_X\in \calE[b]}} (\psi\otimes f)\circ \xi_{\sim_X\sqcup \sim_Y}(ab)\\
&=\sum_{\substack{\sim_X\in \calE[a],\\ \sim_X\in \calE[b]}} (\psi\otimes f)\left(\xi_{\sim_X}(a)\xi_{\sim_Y}(b)\right)\\
&=\psi\leftsquigarrow f(a)\psi\leftsquigarrow f(b).
\end{align*}
We used that $\psi$ respects $\squplus$ and $f$ respects $m$ for the last equality.

Let $X$ be a finite set, and $b\in \calB[X]$. By the compatibility between $\xi$ and $\Delta$,
\begin{align*}
&\psi\leftsquigarrow f(b)(T+T')\\
&=\sum_{\sim\in \calE[b]} (\psi\otimes f)\circ \delta_\sim(b)(T+T')\\
&=\sum_{\sim\in \calE[b]}\sum_{X/{\sim}=X'/{\sim'}\sqcup X''/{\sim''}} (\psi\otimes \psi\otimes f)\circ \left(\Delta_{X'/{\sim'},X''/{\sim''}}\otimes \id_{\calB[X]}\right)\circ \delta_\sim(b)(T,T')\\
&=\sum_{X=X'\sqcup X''}\sum_{\substack{\sim_X\in \calE[b_{\mid_LX'}],\\ \sim_Y\in \calE[b_{\mid_RX''}]}}(\psi\otimes \psi\otimes f)\circ m_{1,3,24}\circ \left(\delta_{\sim_X}\otimes \delta_{\sim_Y}\right)\circ\Delta_{X',X''}(b)(T,T')\\
&=\sum_{X=X'\sqcup X''}\sum_{\substack{\sim_X\in \calE[b_{\mid_LX'}],\\ \sim_Y\in \calE[b_{\mid_RX''}]}}(\psi\otimes f\otimes \psi\otimes f)\circ\left(\delta_{\sim_X}\otimes \delta_{\sim_Y}\right)\circ\Delta_{X',X''}(b)(T,T')\\
&=\sum_{X=X'\sqcup X''}(\psi\leftsquigarrow f\otimes\psi\leftsquigarrow f)\circ \Delta_{X',X''}(b)(T,T').
\end{align*}
We used that $\psi$ is a polynomial invariant of $\exch_+$ for the second equality, and that $f\in \Char(\calB)$ for the fourth one. So $\psi$ is a polynomial invariant of $\calB$.
\end{proof}

\begin{lemma}\label{lem4.2.16}
Let $\calB$ be a quasi-cocommutative CTB, and $\psi$ be a polynomial invariant of $\exch_+$. We denote by $I_1=\{([1])\}$ the unique element of $\exch_+[1]$.
If $\psi(I_1)\neq 0$, there exists $f \in \Char(\calB)$, such that $\phi_+^\calB=\psi\leftsquigarrow f$.
\end{lemma}

\begin{proof}
 Let $b\in \calB[X]$. By Lemma \ref{lem4.2.8}, $\calC_+\left(b,\sim_b^{(\ic)}\right)=\comp[X/{\sim}]$, so the indecomposable components of $\calC_+\left(b,\sim_b^{(\ic)}\right)$ are singletons. 
Therefore, by multiplicativity of $\psi$, $\psi\left(\calC_+\left(b,\sim_b^{(\ic)}\right)\right)(1)\neq 0$.
We define a character $f$ of $\calB$ as follows. Let $X$ be a finite set and $b\in \calB[X]$. We proceed by induction on $|X|$. If $X=\emptyset$, then $b=1_\calB$ and we put $f(1_\calB)=1$. Otherwise, let $b=b_1\cdots b_k$ the decomposition of $b$ into indecomposable elements. 
\begin{enumerate}
\item If $k\geq 2$, we put $f(b)=f(b_1)\cdots f(b_k)$.
\item If $k=1$, then we put
\[f(b)=\frac{1}{\psi\left(\calC_+\left(b,\sim_b^{(\ic)}\right)\right)(1)}\left(\epsilon_+(b)-\sum_{\sim\in \calE[b]\setminus \left\{\sim_b^{(\ic)}\right\}} \psi\left(\calC_+(b,\sim)\right)(1) f(b\mid\sim)\right).\]
Note that if $\sim\in \calE[b]\setminus \left\{\sim_b^{(\ic)}\right\}$, $b\mid \sim$ is not indecomposable, so $f(b\mid\sim)$ is already defined.
\end{enumerate} 
By the unicity of the decomposition into indecomposable components (up to the order), $f$ is well-defined and is a character of $\calB$. For any $b\in \calB[X]$, by construction,
\[\psi\leftsquigarrow f(b)(1)=\epsilon_+(b).\]
By multiplicativity, this is in fact true for any $b\in \calB[X]$. Therefore, by unicity in Theorem \ref{theo1.2.23}, $\psi\leftsquigarrow f=\phi_+^\calB$.
\end{proof}

We can apply this with $\psi=\Phi_{\pp}$, and obtain a quasi-cocommutative version of Theorem \ref{theo4.1.24}.

\begin{theo}\label{theo4.2.17}
Let $\calB$ be a non-quasi-cocommutative CTDB. We define a character of $\calB$ as follows: for any finite $X$, and any $b\in \calB[X]$, indecomposable,
\[g^\calB_{\pp}(b)=\ln_\calB(\epsilon_\delta)(b)=\sum_{k=1}^{|X|}\frac{(-1)^{k-1}}{k} |\calC_+(b)\cap \comp_k[X]|.\]
Then 
\[\phi^\calB_+(b)=\sum_{k=1}^{|X|}\left(\sum_{\substack{\sim\in \calE[b],\\ \cl(\sim)=k}} \frac{|\calC_+(b,\sim)\cap \comp_k[X]|}{k!}g^\calB_{\pp}(b\mid\sim)\right)X^k.\]
\end{theo}

\begin{proof}
We consider $\psi=\Phi_{\pp}$. Then $\Phi_{\pp}(I_1)=X$: we can apply Lemma \ref{lem4.2.16}. There exists $f\in \Char(\calB)$, such that $\phi_+^\calB=\Phi\leftsquigarrow f$. Therefore, for any finite set $X$, and any $b\in \calB[X]$,
\[\phi^\calB_+(b)=\sum_{\sim\in \calE[b]} \Phi_{\pp}(\calC_+(b,\sim))f(b\mid\sim)=\sum_{\sim\in \calE[b]} \dfrac{|\calC_+(b,\sim)\cap \comp_{\cl(\sim)}[X]|}{|\cl(\sim)|!} f(b\mid\sim) T^{\cl(\sim)}.\]

Let us now assume that $b$ is indecomposable. Then $X^2=\sim_b^{(\ic)}$ is the unique element of $\calE[b]$ with only one class. Therefore, the coefficient of $X$ in $\phi^\calB_+(b)$ is 
\begin{align*}
\frac{\mathrm{d}\phi_+^\calB(b)}{\mathrm{d}T}(0)&=\dfrac{|I_1|}{1!}f\left(b\mid\sim_b^{(\ic)}\right)=f(b). 
\end{align*}
Moreover, by Corollary \ref{cor1.2.26}, $\dfrac{\mathrm{d}\phi_+^\calB(b)}{\mathrm{d}T}(0)=\ln_\calB(\epsilon_\delta)(b)$.
\end{proof}

\subsection{Maximal quasi-cocommutative subbialgebra}

\begin{theo}\label{theo4.2.18}
Let $\calB$ be a commutative CTB. For any finite set $X$, we denote by $\calB^{\qc}[X]$ the subset of $\calB[X]$ of elements $b$ such that for any $(X_1,\ldots,X_k)$ and $(X_{\sigma(1)}),\ldots,X_{\sigma(k)}) \in \calC_-(b)$ with $\sigma \in \sym_k$, 
such that $b_{\mid_i X_i}$ and $b_{\mid_i X_{\sigma(i)}}$ are indecomposable for any $i\in [k]$, then
\begin{align*}
&\forall i\in [k],& b_{\mid_i X_i}&=b_{\mid_{\sigma^{-1}(i)} X_i}.
\end{align*}
Then $\calB^{\qc}$ is the greatest quasi-cocommutative twisted subbialgebra of $\calB$.
\end{theo}

\begin{proof}
Let us first prove that $\calB^{\qc}$ is a twisted subbialgebra of $\calB$. 
Let $X$, $Y$ be disjoint finite sets, $a\in \calB^{\qc}[X]$, and $b\in \calB^{\qc}[Y]$. Let us assume that $(X_1,\ldots,X_k)$, $(X_{\sigma(1)},\ldots,X_{\sigma(k)})\in \calC_-(ab)$, with, for $i\in [k]$, $(ab)_{\mid_i X_i}$ and $(ab)_{\mid {\sigma^{-1}(i)} X_i}$ indecomposable. If $i\in [k]$, 
\[(ab)_{\mid_i X_i}=a_{\mid_i X_i\cap X}b_{\mid_i X_i\cap Y}.\]
As $(ab)_{\mid_i X_i}$ is indecomposable, then $X_i\subseteq X$ or $X_i\subseteq Y$. We denote by $(X_{i_1},\ldots,X_{i_p})$ the sequence of $X_i$ included in $X$, and by $(X_{j_1},\ldots,X_{j_q})$ the sequence of $X_i$ included in $Y$. 
Then $\Delta_{X_1,\ldots,X_k}(ab)$ is obtained by permutation of the tensors of
\[\Delta_{X_{i_1},\ldots,X_{i_p}}(a)\otimes \Delta_{X_{j_1},\ldots,X_{j_q}}(b).\]
Similarly, $\Delta_{X_{\sigma(1)},\ldots,X_{\sigma(k)}}(ab)$ is obtained by a permutation of tensors
of $\Delta_{X_{\sigma(i_1)},\ldots,X_{\sigma(i_p)}}(a)\otimes\Delta_{X_{\sigma(j_1)},\ldots,X_{\sigma(j_q)}}(b)$. As $a\in \calB^{\qc}[X]$ and $b\in \calB^{\qc}[Y]$, $\Delta_{X_{\sigma(i_1)},\ldots,X_{\sigma(i_p)}}(a)\otimes\Delta_{X_{\sigma(j_1)},\ldots,X_{\sigma(j_q)}}(b)$
is obtained from $\Delta_{X_{i_1},\ldots,X_{i_p}}(a)\otimes \Delta_{X_{j_1},\ldots,X_{j_q}}(b)$ by permutations of the tensors. Therefore, $ab\in \calB^{\qc}[X\sqcup Y]$.

Let $X$, $Y$ be disjoint finite sets, and $b\in \calB^{\qc}[X\sqcup Y]$, such that $(X$, $Y)\in \calC_-(b)$. Let us assume that $(X_1,\ldots,X_k)$, $(X_{\sigma(1)},\ldots,X_{\sigma(k)})\in \calC_-(b_{\mid_LX})$, with, 
for $i\in [k]$, ${b_{\mid_LX}}_{\mid_i X_i}$ and ${b_{\mid_LX}}_{\mid {\sigma^{-1}(i)} X_i}$ indecomposable. Let $(Y_1,\ldots,Y_l)\in \calC_-(b_{\mid_RY})$, such that $\Delta_{Y_1,\ldots,Y_l}(b_{\mid_RY})$ is a tensor of indecomposable elements: this exists, take for example $(Y_1,\ldots,Y_l)$ to be the indecomposable components of $b_{\mid_RY}$. Then
\begin{align*}
\Delta_{X_1,\ldots,X_k,Y_1,\ldots,Y_l}(b)&=\Delta_{X_1,\ldots,X_k}(b_{\mid_LX})\otimes \Delta_{Y_1,\ldots,Y_l}(b_{\mid_RY}),\\
\Delta_{X_{\sigma(1)},\ldots,X_{\sigma(k)},Y_1,\ldots,Y_l}(b)&=\Delta_{X_{\sigma(1)},\ldots,X_{\sigma(k)}}(b_{\mid_LX})\otimes \Delta_{Y_1,\ldots,Y_l}(b_{\mid_RY}),
\end{align*}
and both are tensors of indecomposable elements. As $b\in \calB^{\qc}[X\sqcup Y]$, the second one is obtained from the first one by permutations of tensors. Consequently, the same holds for $\Delta_{X_1,\ldots,X_k}(b_{\mid_LX})$ and $\Delta_{X_{\sigma(1)},\ldots,X_{\sigma(k)}}(b_{\mid_RY})$:
$b_{\mid_LX}\in \calB^{\qc}[X]$. Similarly, $b_{\mid_RY}\in \calB^{\qc}[Y]$.
As $\calB^{\qc}$ is commutative (as a twisted subalgebra of $\calB$), it is quasi-cocommutative.

Let $\calB'$ be a quasi-cocommutative twisted subbialgebra of $\calB$. By definition of quasi-cocommutativity, $\calB'\subseteq \calB^{\qc}$. 
\end{proof}

\begin{prop}\label{prop4.2.19}
Let $\calB$ and $\calB'$ be two commutative CTB, and $\phi:\calB\longrightarrow\calB'$ be a morphism of CTBs. Then $\phi(\calB^{\qc})\subseteq \calB'^{\qc}\sqcup \{0\}$.
\end{prop}

\begin{proof}
It is enough to prove that the image of $\phi(\calB^{\qc})$ (except 0) is quasi-commutative.
Let $X$ be a finite set, $b\in \calB[X]$ such that $\phi(b)\neq 0$, and $(X_1,\ldots,X_k)$, $(X_{\sigma(1)},\ldots,X_{\sigma(k)})$ two elements of $\calC_-(\phi(b))$, such that for any $i\in [k]$, $\phi(b)_{\mid_i X_i}$ and $\phi(b)_{\mid_{\sigma^{-1}(i) X_i}}$ are indecomposable. 
As $\phi$ is a morphism of twisted coalgebras, for any $i\in [k]$,
\begin{align*}
\phi(b)_{\mid_i X_i}&=\phi\left(b_{\mid_i X_i}\right),&
\phi(b)_{\mid_{\sigma^{-1}(i) X_i}}&=\phi\left(b_{\mid_{\sigma^{-1}(i) X_i}}\right).
\end{align*}
By Proposition \ref{prop2.2.6}, $\ic\left(\phi(b)_{\mid_i X_i}\right)=1\geq \ic\left(b_{\mid_i X_i}\right)$,
so $b_{\mid_i X_i}$ is indecomposable. Similarly, $b_{\mid_{\sigma^{-1}(i) X_i}}$ is indecomposable. 
As $b\in \calB^{\qc}[X]$, for any $i\in [k]$, $b_{\mid_i X_i}=b_{\mid_{\sigma^{-1}(i) X_i}}$, which implies that
$\phi(b)_{\mid_i X_i}=\phi(b)_{\mid_{\sigma^{-1}(i) X_i}}$. As a conclusion, $\phi(b)\in \calB'^{\qc}[X]$.
\end{proof}

\begin{prop}\label{prop4.2.20}
Let $\calB$ be a quasi-cocommutative twisted bialgebra, $X$ a finite set and $b\in \calB[X]$. 
For any $\sim\in \calE[b]$, such that $\calC_+(b,\sim)\neq \emptyset$, $\calC_+(b,\sim)\in \exch_+^{\qc}[X/{\sim}]$.
\end{prop}

\begin{proof}
Let $(X_1/{\sim_1},\ldots,X_k/{\sim_k})$ and $(X_{\sigma(1)}/{\sim_{\sigma(1)}},\ldots,X_{\sigma(k)}/{\sim_{\sigma(k)}})$ be elements of $m(\calC_+(b,\sim))$, such that for any $i\in [k]$, $\calC_+(b,\sim)_{\mid_i X_i}$ and $\calC_+(b,\sim)_{\mid_{\sigma^{-1}(i)} X_i/{\sim_i}}$ are indecomposable.
Then $(X_1,\ldots,X_k)$ and $(X_{\sigma(1)},\ldots,X_{\sigma(k)})$ both belong to $\calC_-(b)$ and for any $i\in [k]$,
$\sim_i=\sim^{(\ic)}_{b_{\mid_i X_i}}=\sim^{(\ic)}_{b_{\mid_{\sigma^{-1}(i)}X_i}}$. We arbitrarily index the elements of $X/{\sim}_i$ as $(X_{i,1},\ldots,X_{i,l_i})$. Then 
\[\pi=(X_{1,1},\ldots,X_{1,l_1},\ldots,X_{k,1},X_{k,l_k})\in \calC_-(b),\]
and the factors of $\Delta_\pi(b)$ are the indecomposable components of the $b_{\mid_iX_i}$'s.
Similarly, $\pi'=(X_{\sigma(1),1},\ldots,X_{\sigma(1),l_{\sigma(1)}},\ldots,X_{\sigma(k),1},\ldots,X_{\sigma(k),l_{\sigma(k)}})\in \calC_-(b)$ and the factors of $\Delta_{\pi'}(b)$ are the indecomposable components of the $b_{\mid_{\sigma^{-1}(i)} X_i}$'s. 
As $\calB$ is quasi-cocommutative, the factors of $\Delta_\pi(b)$ and $\Delta_{\pi'}(b)$ are the same, so for any $i\in [k]$, $b_{\mid_iX_i}$ and $b_{\mid_{\sigma^{-1}(i)} X_i}$ have the same indecomposable components: 
as a consequence of the commutativity of $\calB$, $b_{\mid_iX_i}=b_{\mid_{\sigma^{-1}(i)} X_i}$. The compatibility of $\Delta$ and $\xi$ gives that for any $i\in [k]$,
\begin{align*}
(b_{\mid_iX_i})/{\sim_i}&(\calC_+(b,\sim))_{\mid_i X_i/{\sim_i}},&
(b_{\mid_{\sigma^{-1}(i)}X_i})/{\sim_i}&=(\calC_+(b,\sim))_{\mid_{\sigma^{-1}(i)} X_i/{\sim_i}}
\end{align*}
which finally gives $(\calC_+(b,\sim))_{\mid_i X_i/{\sim_i}}=(\calC_+(b,\sim))_{\mid_{\sigma^{-1}(i)} X_i/{\sim_i}}$. So $\calC_+(b,\sim)\in \exch_+^{\qc}[X/{\sim}]$.
\end{proof}

\section{Second coproduct on $\exch_+$}

We know (Example \ref{ex4.2.3}) that the whole $\exch_+$ cannot be made a CTDB, as it is not quasi-cocommutative. We now study $\exch_+^{\qc}$, and look for a possible second coproduct on it.

\subsection{Second coproduct on $\exch_+^{\qc}$}

Let us apply the results and definitions of the preceding section for $\calB=\exch_+$.We obtain a twisted subbialgebra $\exch_+^{\qc}$. We denote $\delta$ instead of $\xi$ for $\exch_+^{\qc}$, and for any finite set $X$, and any $\sim\in \exch_+^{\qc}[X]$,
$\delta_\sim:\exch_+^{\qc}[X]\longrightarrow \exch_+^{\qc}[X/{\sim}]\otimes\exch_+^{\qc}[X]$. Let us give more details on these objects.

For any finite set $X$, $\exch_+^{\qc}[X]$ is the set of elements $I\in \exch_+[X]$ such that, if $(X_1,\ldots,X_k)$ and $(X_{\sigma(1)},\ldots X_{\sigma(k)})$ both belong to $m(I)$, for $\sigma\in \sym_k$, 
with, for any $i\in [k]$, $I_{\mid_i X_i}$ and $I_{\mid_{\sigma^{-1}(i)} X_i}$ indecomposable, then, for any $i\in [k]$, $I_{\mid_i X_i}=I_{\mid_{\sigma^{-1}(i)} X_i}$. 

Let $I\in \exch+^{\qc}[X]$, and $\sim\in \calE[X]$. Then $\sim\in \calE[I]$ if there exists $(X_1,\ldots,X_k)\in m(I)$, such that:
\begin{itemize}
\item $X/{\sim}=\{X_1,\ldots,X_k\}$.
\item Putting $\Delta_{X_1,\ldots,X_k}(I)=I_{\mid_1X_1}\otimes \ldots \otimes I_{\mid_kX_k}$, then for any $i\in [k]$, $I_{\mid_iX_i}$ is indecomposable. 
\end{itemize}
If $\sim\in \calE[I]$, we put $I\mid \sim=I_{\mid_1X_1}\squplus\cdots\squplus I_{\mid_kX_k}$. Moreover, 
for any $\pi=(X_1/{\sim}_1,\ldots,X_k/{\sim}_k)\in \comp[X/{\sim}]$, $\pi\in I/{\sim}$ if $(X_1,\ldots,X_k)\in m(I)$ and $\sim_i=\sim_{I_{\mid_i X_i}}^{(\ic)}$ for any $i\in [k]$. 
By Definition \ref{defi4.2.6} and Proposition \ref{prop4.2.20} for $\calB=\exch_+^{\qc}$, $I/{\sim}\in \exch_+^{\qc}[X/{\sim}]$.
The coproduct $\delta_\sim$ of Theorem \ref{theo4.2.10} is defined from $\exch_+^{\qc}[X]$ to $\exch_+^{\qc}[X/{\sim}]\otimes \exch^{\qc}[X]$. Moreover, by Theorem \ref{theo4.2.10}:
\begin{enumerate}
\item $\delta$ is compatible with the product $\squplus$.
\item $\delta$ is compatible with the coproduct $\Delta$.
\item $\delta$ satisfies this partial coassociativity: for any finite set $X$, and any $\sim,\sim'\in \calE[X]$, such that we do not have $\sim\subseteq \sim'$,
\[\left(\id_{\exch_+[X/{\sim'}]}\otimes \delta_\sim\right)\circ \delta_{\sim'}=0.\]
\item $\epsilon_+$ is a counit for $\delta_\sim$.
\end{enumerate}
Unfortunately, $\delta$ is not fully coassociative. We shall have to restrict again to a smaller subspecies. 

\begin{theo} \label{theo4.3.1}
There exists a sub-CTB $\exch_+^{\double}$ of $\exch_+^{\qc}$, such that:
\begin{enumerate}
\item $(\exch_+^{\double},\squplus,\delta,\Delta)$ is a CTDB.
\item If $(\calA,\squplus,\Delta)$ is a sub-CTB of $\exch_+$ which can be given a coproduct $\delta'$ which makes it double, then $\calA\subseteq \exch_+^{\double}$ and $\delta'$ is the restriction of $\delta$ to $\calA$. 
\end{enumerate} \end{theo}

\begin{proof}
1. We define $\exch_+^{\double}[X]$ by induction on $|X|$ as follows. 
\begin{enumerate}
\item $\exch_+^{\double}[\emptyset]=\{1_{\exch_+}\}$.
\item If $|X|\geq 1$, then $\exch_+^{\double}[X]$ is the set of $I\in \exch+^{\qc}[X]$ such that:
\begin{enumerate}
\item If $(X',X'')\in m(I)$, then $I_{\mid_LX'}\in \exch_+^{\double}[X']$ and $I_{\mid_RX''}\in \exch_+^{\double}[X'']$.
\item If $\sim \in \calE[I]$, different from $\sim_I^{(m)}$, then $I/{\sim}\in \exch_+^{\double}[X/{\sim}]$.
\item For any $\sim\subseteq \sim'\in \calE[X]$, 
\[\left(\delta_{\overline{\sim'}}\otimes \id_{\exch+^{\qc}[X]}\right)\circ \delta_\sim(I)=\left(\id_{\exch+^{\qc}[X/{\sim}]}\otimes \delta_\sim\right)\circ \delta_{\sim'}(I).\]
\end{enumerate}\end{enumerate}
This is well-defined, as if $(X',X'')\in m(I)$, then $|X'|<|X|$ and $|X''|<|X|$, and, if $\sim \neq \sim_I^{(m)}$, then $|X/{\sim}|<|X|$. By 2.(a)-(b), $\exch_+^{\double}$ is stable under $\Delta$. By 2.(c), the restriction of $\delta$ to $\exch_+^{\double}$ is coassociative.
Moreover, as $\delta$ is compatible with $\squplus$ and $\Delta$, it is not difficult to obtain that if $X$, $Y$ are disjoint finite sets, $I\in \exch_+^{\double}[X]$ and $J\in \exch_+^{\double}[Y]$, then $I\squplus J\in \exch_+^{\double}[X\sqcup Y]$ (proceed by induction on $|X|+|Y|$). Consequently, as $\exch_+^{\double}$ is stable under $\Delta$ and $\squplus$, if $\sim\in \calE[I]$ and $I\in \exch_+^{\double}[X]$, then $I\mid\sim\in \exch_+^{\double}[X]$.
Therefore, $(\exch_+^{\double},\squplus,\Delta,\delta)$ is a CTDB.

Let $\calA$ be a subspecies of $\exch+$, such that $(\calA,\squplus,\Delta,\delta')$ is a CTDB for a certain coproduct $\delta$'. Then $\calA$ is quasi-cocommutative. 
By Proposition \ref{prop4.2.14}, the canonical injection takes its values in $\exch_+^{\qc}$, so $\calA\subseteq \exch_+^{\qc}$. 
Let us show that $\calA[X]\subseteq \exch_+^{\double}[X]$ by induction on $|X|$. As $\exch_+^{\double}[\emptyset]=\exch_+[\emptyset]$, this is obvious if $X=\emptyset$. Let us assume the result at all ranks $<|X|$. Let $I\in \calA[X]$. As $(\calA,\squplus,\Delta,\delta)$ is a CTDB, 2.(c) is satisfied. Let $(X',X'')\in m(I)$. The induction hypothesis insures that $I_{\mid_LX'}\in \exch_+^{\double}[X']$ and $I_{\mid_RX''}\in \exch_+^{\double}[X'']$, as $\calA$ is a twisted subcoalgebra for $\Delta$. Hence, 2. (a) is satisfied. Let $\sim\in \calE[I]$, different from $\sim_I^{(m)}$. Then $I/{\sim}\in \calA[X/{\sim}]\subseteq \exch_+^{\double}[X/{\sim}]$, by the induction hypothesis. Hence, 2. (b) is satisfied. Therefore, $I\in \exch_+^{\double}[X]$.

It remains to show that $\delta'$ is the restriction of $\delta$ to $\calA$. 
Let $X$ be a finite set, and $I\in \calA[X]$. As $\exch_+$ is non-degenerate, $\calA$ is non-degenerate. By Lemma \ref{lem4.2.3} for $b=I$, for any $\sim\in \calE[X]$, $\sim\in \calE'[I]$ if, and only if, $\sim\in \calE[I]$ by definition of $\calE[I]$, as $\calC_-(I)=m(I)$. Moreover, if $\sim\in \calE[I]$, Lemma \ref{lem4.2.3} gives that $I\mid'\sim=I\mid\sim$. Lemma \ref{lem4.2.4} gives that $I/'\sim=I/{\sim}$. Finally, $\delta'=\delta$.
\end{proof}

\begin{cor}\label{cor4.3.2}
Let $\calB$ be a CTB, such that $\Gamma^\calB_+:\calB\longrightarrow \exch_+$ is injective. Then there exists at most one coproduct $\delta$ on $\calB$, making $\calB$ a CTDB.
\end{cor}

\begin{proof}
Let $\delta$ be a coproduct on $\calB$, making $\calB$ a CTDB.
We define a coproduct $\delta'$ on $\Gamma^\calB_+(\calB)$ by 
\[\delta'=\left(\Gamma^\calB_+\otimes \Gamma^\calB_+\right)\circ \delta \circ {\Gamma^\calB_+}^{-1},\]
 thanks to the injectivity of $\Gamma^\calB_+$. Then $\Gamma^\calB_+(\calB)$ is a twisted subbialgebra of $\exch_+$, and a CTDB. 
By Theorem \ref{theo4.3.1}, $\Gamma^\calB_+(\calB)\subseteq \exch_+^{\double}$, and $\delta'=\delta_{\mid\Gamma^\calB_+(\calB)}$. Therefore, 
\[\delta=\left({\Gamma^\calB_+}^{-1}\otimes {\Gamma^\calB_+}^{-1}\right)\circ \delta_{\mid \exch_+^{\double}}\circ \Gamma^\calB_+.\]
So, if it exists, $\delta$ is unique. 
\end{proof}

Let us apply this to graphs and posets. 

\begin{lemma}\label{lem4.3.3}
Let $X$ be a finite set. 
\begin{enumerate}
\item Let $G\in \calG[X]$. Then $E(G)$ is the set of pairs $\{x,y\}$ such that for any $\pi\in \calC_+(G)$, $x$ and $y$ are in two different blocks of $\pi$.
Consequently, $\Gamma^\calB_+:\calG\longrightarrow \exch_+$ is injective.
\item Let $P\in \calP[X]$. For any $x\neq y\in X$, $x\leq_P y$ if, and only if, in any $\pi\in \calC_+(P)$, the block containing $x$ appears before the block containing $y$.
Consequently, $\Gamma^\calB_+:\calP\longrightarrow \exch_+$ is injective.
\end{enumerate}
\end{lemma}

\begin{proof}
$1$. Let $\{x,y\}\in E(G)$. For any $\pi=(X_1,\ldots,X_k)\in \calC_+(G)$, for any $i\in [k]$, $G_{\mid X_i}$ does not have any edge, so $x,y$ cannot be in the same block of $\pi$. 
Dually, let us assume that $\{x,y\}\notin E(G)$. Then $G_{\mid \{x,y\}}$ has no edge, so belong to $\bfD\calG[\{x,y\}]$. Taking a composition $\pi$ of $X$ of the form $(\{x,y\},\{x_1\},\ldots,\{x_n\})$, $\pi\in \calC_+(G)$, and $x,y$ are in the same block of $\pi$.

$2.\Longrightarrow$. Let $\pi=(X_1,\ldots,X_k)\in \calC_+(P)$. As it belongs also to $\calC_-(P)$, by definition of $\Delta$ on $\calP$, assuming that $x\in X_i$ and $y\in X_j$, necessarily $i\leq j$. 
Moreover, $P_{\mid X_i}\in \bfD\calP[X_i]$, so $i<j$. 

$2.\Longleftarrow$. We assume that we do not have $x\leq_P y$. We consider
\begin{align*}
X''&=\{x'\in X\mid x\leq_P x''\},&X'&=X\setminus X''.
\end{align*}
If $x''\in X''$ and $x'''\in X$, such that $x''\leq_P x'''$, then $x\leq_P x''\leq_P x'''$, so $x'''\in X''$: we obtain that $(X',X'')\in \calC_-(P)$. Moreover, by hypothesis, $y\in X'$, and obviously, $x\in X''$. 
As $\calP$ is non-degenerate, let us take $\pi'\in \calC_+(P_{\mid X')})$ and $\pi''\in \calC_+(P_{\mid X''})$. Then $\pi=(\pi',\pi'')\in \calC_+(P)$. As $x\in X''$ and $y\in X'$, the block containing $y$ appears before the block containing $x$ in $\pi$. 
\end{proof}

\begin{example}\label{ex4.3.1}\begin{enumerate}
\item As $\Gamma^\calG_+$ and $\Gamma^\calP_+$ are injective, $\delta$ described in Examples \ref{ex4.1.1} and \ref{ex4.1.3} are the unique coproducts making $\calG$ and $\calP$ CTDBs.
\item $\Gamma^{\calO\calG}_+$ is not injective. Consider the two oriented graphs $G$, $G'\in \calO\calG[3]$ such that
\begin{align*}
A(G)&=\{(1,2),\:(2,3)\},& A(G')&=\{(1,2),\:(1,3),\:(2,3)\}.
\end{align*}
Then $\calC_+(G)=\calC_+(G')=\{(1,2,3)\}$.
\item $\Gamma^{\calH\calG}_+$ is not injective. Consider the two hypergraphs $G$ and $G'\in \calH\calG[3]$, given by
\begin{align*}
\HE(G)&=\{\{1\},\:\{2\},\:\{3\},\:\{1,2\}\},&E(G')&=\{\{1\},\:\{2\},\:\{3\},\:\{1,2\},\:\{1,2,3\}\}.
\end{align*}
Then $\calC_+(G)=\calC_+(G')$ is the set of compositions $\pi\in \comp[3]$ such that $1$ and $2$ do not belong to the same block of $\pi$. 
\end{enumerate}\end{example}

\subsection{$\Gamma^\calB_+$ as a double bialgebra morphism}

\begin{prop}\label{prop4.3.4}
Let $\calB$ and $\calB'$ be two regular CTDBs. Let $\phi:\calB\longrightarrow \calB'$ be a morphism of CTDBs. Then, for any finite set $X$, for any $b\in \calB[X]$ such that $\phi(b)\neq 0$, 
\[\ic(\phi(b))=\ic(b).\]
\end{prop}

\begin{proof}
By Proposition \ref{prop2.2.6}-2, $\ic(\phi(b))\geq \ic(b)$.

Let us prove that $\ic(\phi(b))=\ic(b)$ by induction on $|X|$. If $X=\emptyset$, then $b=1_\calB$, $\phi(b)=1_{\calB'}$, and $\ic(\phi(b))=\ic(b)=0$. If $|X|=1$, then both $b$ and $\phi(b)$ are indecomposable, so $\ic(\phi(b))=\ic(b)=1$.
Let us assume the result at all ranks $<|X|$, with $|X|\geq 2$.
\begin{itemize}
\item If $\sim_b^{(\ic)}\neq \sim_b^{(m)}$, let us write $\sim=\sim_b^{(\ic)}$ for short. Then $|X/{\sim}|<|X|$. Moreover,
\[\epsilon_\delta(b/{\sim})=1=\epsilon_\delta(\phi(b/{\sim})),\]
so $\phi(b/{\sim})\neq 0$. As $\calB$ and $\calB'$ are regular, and by the induction hypothesis,
\[\ic(b)=\ic(b/{\sim})=\ic(\phi(b/{\sim}))=\ic(\phi(b)/{\sim})=\ic(\phi(b)).\]
\item If $\sim_b^{(\ic)}=\sim_b^{(m)}$, then $b\in \bfD\calB[X]$. Therefore, $\ic(b)=|X|$, and
\[|X|\geq \ic(\phi(b))\geq \ic(b)=|X|,\]
so $\ic(\phi(b))=\ic(b)=|X|$. \qedhere
\end{itemize}\end{proof}

\begin{theo}\label{theo4.3.5}
Let $\calB$ be a non-degenerate CTDB. 
\begin{enumerate}
\item If $\Psi:\calB\longrightarrow \exch_+^{\double}$ is a morphism of CTDBs, then $\Psi=\Gamma^\calB_+$.
\item The following assertions are equivalent:
\begin{enumerate}
\item $\Gamma^\calB_+:\calB\longrightarrow \exch_+$ takes its values in $\exch_+^{\double}$ and is a morphism of CTDBs.
\item For any finite set $X$, for any $b\in \calB[X]$, $\ic(\calC_+(b))=\ic(b)$.
\end{enumerate}\end{enumerate}\end{theo}

\begin{proof}
1. Then $\epsilon_+\circ \Psi=\epsilon^\calB_+$, as $\Psi$ is compatible with the counits of the coproduct $\delta$. By Theorem \ref{theo3.2.13}, $\Psi=\Gamma^\calB_+$.

2. $(a)\Longrightarrow (b)$. As $\calB$ and $\exch_+^{\double}$ are both non-degenerate, they are regular (Proposition \ref{prop4.1.16}). By Proposition \ref{prop4.3.4}, for any finite set $X$, for any $b\in \calB[X]$, $\ic(\calC_+(b))=\ic(b)$.

2. $(b)\Longrightarrow (a)$. Let $X$ be a finite set, $b\in \calB[X]$, and $\sim \in \calE[X]$. Then $\calC_-(b)=m(\calC_+(b))$.
The hypothesis (b) implies that for any $b'\in \calB$, $b'$ is indecomposable if, and only if, $\calC_+(b')$ is indecomposable. Therefore, Lemma \ref{lem4.2.3} implies that $\calE[b]=\calE[\calC_+(b)]$,
and that $\calC_+(b\mid\sim)=\calC_+(b)\mid\sim$ for any $\sim\in \calE[b]$.
Lemma \ref{lem4.2.4} and the hypothesis $(b)$ implies that $\calC_+(b/{\sim})=\calC_+(b)/{\sim}$. Therefore, 
\[\delta\circ \Gamma^\calB_+=(\Gamma^\calB_+\otimes \Gamma^\calB_+)\circ \delta.\]
Consequently, $\Gamma^\calB_+(\calB)$ is a sub-CTDB of $\exch_+$ with the coproduct $\delta$.
By Theorem \ref{theo4.3.1}, $\Gamma^\calB_+(\calB)\subseteq \exch_+^{\double}$. 
\end{proof}

\begin{cor}\label{cor4.3.6}
Let us assume that $\calB$ is a CTDB, such that $\Gamma^\calB_+$ is injective.
Then $\Gamma^\calB_+$ takes its values in $\exch_+^{\double}$ and is a morphism of CTDBs.
\end{cor}

\begin{proof}
By Proposition \ref{prop2.2.6}, as $\Gamma^\calB_+$ is injective, is conserves the number of indecomposable components, so it a morphism of CTDBs.
\end{proof}

\begin{example}\label{ex4.3.2}\begin{enumerate}
\item For the CTDBs $\calG$ of graphs and $\calP$ of posets, $\Gamma^\calB_+$ is injective (Lemma \ref{lem4.3.3}). 
Therefore, in both cases, $\Gamma^\calB_+$ is a morphism of CTDBs from $\calG$ or $\calP$ to $\exch_+^{\double}$. 
\item Let us consider again the hypergraphs $G$ and $G'$ on $[3]$ of Example \ref{ex4.3.1}:
\begin{align*}
\HE(G)&=\{\{1\},\:\{2\},\:\{3\},\:\{1,2\}\},&\HE(G')&=\{\{1\},\:\{2\},\:\{3\},\:\{1,2\},\:\{1,2,3\}\}.
\end{align*}
We already observed that $\calC_+(G)=\calC_+(G')$. As $G'$ is connected, $\ic(G')=1$. As $G$ is not connected, $\ic(\calC_+(G'))=\ic(\calC_+(G))\geq 2$. 
Therefore, $\Gamma^\calB_+:\calH\calG\longrightarrow\exch$ is not a morphism of CTDBs. 
\item Let us consider the CTDB of acyclic oriented graphs $\calO\calG_{\mathrm{ac}}$ (Example \ref{ex4.1.3}). 
As the map $\phi_{\calO\calG}$ is a morphism of CTDBs, by Proposition \ref{prop4.3.4}, for any finite set $X$ and any $G\in \calO\calG_{\mathrm{ac}}[X]$,
\[\ic(\phi_{\calO\calG}(G))=\ic(G).\]
Moreover, as $\Gamma^\calB_+$ is injective for posets, 
\[\ic(\Gamma^\calB_+(G))=\ic(\Gamma^\calB_+(\theta(G)))=\ic(\theta(G))=\ic(G).\]
So $\Gamma^\calB_+:\calO\calG_{\mathrm{ac}}\longrightarrow \exch_+^{\double}$ is a morphism of CTDBs.
\end{enumerate}\end{example}

\begin{cor}
Let $\calB$ be a regular CTDB, such that for any finite set $X$, for any $b\in \calB[X]$, $\calC_+(b)=\emptyset$ or $\ic(\calC_+(b))=\ic(b)$. 
Then $\Gamma^\calB_+:\calB\longrightarrow \exch_+$ takes its values in $\exch_+^{\double}$ and is a morphism of CTDBs. 
\end{cor}

\begin{proof}
This comes the fact that $\Gamma^\calB_+=\Gamma^{\calB'}_+\circ \varpi$, where $\varpi:\calB\longrightarrow \nd\calB$ is the canonical surjection (which is a morphism of CTDBs), 
and $\Gamma^{\nd\calB}_+:\nd\calB\longrightarrow \exch_+$ is a morphism of CTDBs, from Theorem \ref{theo4.3.5}.
\end{proof}

\subsection{The map $\zeta$ as a coaction}

\begin{defi}
Let $\calB$ be a quasi-cocommutative CTB. We shall say that $\zeta$ is a coaction if:
\begin{enumerate}
\item For any finite set $X$, and $b\in \calB[X]$, and any $\sim\in \calB[b]$, $\calC_+(b,\sim)\in \exch_+^{\double}[X/{\sim}]$.
\item For any finite set $X$, and $b\in \calB[X]$, and any $\sim\subseteq \sim'\in \calB[b]$,
\[(\delta_{\overline{\sim}} \otimes \id_{\calB[X]})\circ \zeta_{\sim'}(b)
=(\id_{\exch_+[X/{\sim}]}\otimes \zeta_{\sim'})\circ \zeta_\sim(b).\]
\end{enumerate}\end{defi}

If $\zeta$ is a coaction, we can define $\zeta:\calF[\calB]\longrightarrow \calF[\exch_+^{\double}]\otimes \calF[\calB]$ by
\[\zeta(\overline{b})=\sum_{\sim\in \calE[b]} \overline{\calC_+(b,\sim)} \otimes \overline{b\mid\sim}.\]
Then $\zeta$ is a left coaction of the bialgebra $(\calF[\exch_+^{\double}],\squplus,\delta)$ onto $\calF[\calB]$,
and $(\calF[\calB],m,\Delta)$ is a bialgebra in the category of left comodules over the bialgebra
$(\calF[\exch_+^{\double}],\squplus,\delta)$.

\begin{prop}\label{prop4.3.9}
Let $\calB$ be a quasi-cocommutative CTB. Then $\zeta$ is a coaction if, and only if, the two following properties are satisfied:
\begin{enumerate}
\item For any finite set $X$, and $b\in \calB[X]$, and any $\sim\subseteq \sim'\in \calE[X]$,
\[ (\sim'\in \calE[b]\mbox{ and }\overline{\sim}\in \calE[\calC_+(b,\sim')])\Longrightarrow (\sim \calE[b]\mbox{ and }\sim'\in \calE[b/{\sim}]).\]
\item For any finite set $X$, and $b\in \calB[X]$, and any $\sim\in \calE[b]$, $\ic(\calC_+(b,\sim))=\ic(b)$.
\end{enumerate}\end{prop}

\begin{proof}
We start with a series of observations.

\textit{First observation}. Let $X$ be a finite set, and $b\in \calB[X]$. For any $\sim\subseteq \sim'\in \calE[X]$,
\begin{align*}
(\delta_{\overline{\sim}} \otimes \id_{\calB[X]})\circ \zeta_{\sim'}(b)&=\begin{cases}
\calC_+(b,\sim')/\overline{\sim}\otimes \calC_+(b,\sim')\mid \overline{\sim}\otimes b\mid \sim'\mbox{ if }\sim'\in \calE[b]\mbox{ and }\overline{\sim}\in \calE[\calC_+(b,\sim')],\\
0\mbox{ otherwise}
\end{cases}\\
(\id_{\exch_+[X/{\sim}]}\otimes \zeta_{\sim'})\circ \zeta_\sim(b)&=\begin{cases}
\calC_+(b,\sim)\otimes \calC_+(b\mid \sim,\sim')\otimes (b\mid \sim)\mid \sim'\mbox{ if }\sim \calE[b]\mbox{ and }\sim'\in \calE[b/{\sim}],\\
0\mbox{ otherwise}.
\end{cases}\end{align*}
Therefore, $\zeta$ is a coaction if, and only if, the two following conditions hold:
\begin{enumerate}
\item For any finite set $X$, any $b\in \calB[X]$, and $\sim\subseteq \sim'\in \calE[X]$,
\[(\sim \calE[b]\mbox{ and }\sim'\in \calE[b/{\sim}])\Longleftrightarrow (\calE[b]\mbox{ and }\overline{\sim}\in \calE[\calC_+(b,\sim')]).\]
\item For any finite set $X$, any $b\in \calB[X]$, and $\sim\subseteq \sim'\in \calE[X]$, such that the conditions of the first item hold, 
\begin{align*}
\calC_+(b,\sim')/\overline{\sim}&=\calC_+(b,\sim),& \calC_+(b,\sim')\mid\overline{\sim}=\calC_+(b\mid\sim,\sim'),&b\mid\sim'&=(b\mid\sim)\mid\sim'.
\end{align*}\end{enumerate}

\textit{Second observation}. Let $X$ be a finite set, $b\in \calB[X]$, and $\sim\in \calE[b]$. 
Let us decompose $b$ into indecomposable elements:
\[b=b_1\cdots b_k,\]
with $b_i\in \calB[X_i]$ for $i\in [k]$. By the compatibility between $\zeta$ and the product, putting $\sim=\sim\cap X_i^2$ for $i\in [k]$, $\sim=\sim_1\sqcup\cdots\sqcup \sim_k$ and
\[\calC_+(b,\sim)=\calC_+(b_1,\sim_1)\squplus\cdots \squplus \calC_+(b_k,\sim_k).\]
Therefore, $\sim_{\calC_+(b,\sim)}^{(\ic)}\subseteq \overline{\sim_b^{(\ic)}}$, so $\ic\left(\calC_+(b,\sim)\right)\geq \ic(b)$, with an equality if, and only if, $\sim_{\calC_+(b,\sim)}^{(\ic)}=\overline{\sim_b^{(\ic)}}$.

\textit{Third observation}. Let $X$ be a finite set, $b\in \calB[X]$, and $\sim\in \calE[b]$. Let us prove that
\[m(\calC_+(b,\sim))=\{(X_1/{\sim_1},\ldots,X_k/{\sim_k})\in \comp[X/{\sim}]\mid (X_1,\ldots,X_k)\in \calC_-(b)\}.\]
$\subseteq$. Let $(X_1/{\sim_1},\ldots,X_k/{\sim_k})\in m(\calC_+(b,\sim))$. There exists $(X_1'/{\sim_1'},\ldots,X'_l/{\im'_l})\in \calC_+(b,\sim)$, such that
\[(X_1/{\sim_1},\ldots,X_k/{\sim_k})\leq (X_1'/{\sim_1'},\ldots,X'_l/{\im'_l}).\]
By definition of $\calC_+(b,\sim)$, $(X'_1,\ldots,X'_l)\in \calC_b(b)$. Moreover, $(X_1,\ldots,X_k)\leq (X'_1,\ldots,X'_l)$. 
As $\calC_-(b)\in \exch_-[X]$, $(X_1,\ldots,X_k)\in \calC_-(b)$.

$\supseteq$. Let $(X_1/{\sim_1},\ldots,X_k/{\sim_k})\in \comp[X/{\sim}]$, such that $(X_1,\ldots,X_k)\in \calC_-(b)$.
By the compatibility between $\zeta$ and $\Delta$, for any $i\in [k]$,
\[\calC_+\left(b_{\mid_i X_i},\sim_i)\right)=\calC_+(b,\sim)_{\mid_i X_i}.\]
Let $X_{i,1}/{\sim_{i,1}},\ldots,X_{i,l_i}/{\sim_{i,l_i}}$ be the indecomposable components of $\calC_+\left(b_{\mid_i X_i},\sim_i\right)$. Then 
\[(X_{1,1}/{\sim_{1,1}},\ldots,X_{k,l_k}/{\sim_{k,l_k}})\in m(\calC_+(b,\sim)).\]
By the compatibility of $\zeta$ with $\Delta$, $(X_{1,1}/{\sim_{1,1}},\ldots,X_{k,l_k}/{\sim_{k,l_k}})\in \calC_-(b)$.
We denote by $\alpha(i,j)$ the position of $X_{i,j}$ in this sequence. For any $i\in [k]$, any $j\in [l_i]$,
by the hypothesis of conservation of $\ic$ and the second observation,
\begin{align*}
\overline{\sim_{\left(b_{\mid_i X_i}\right)_{\mid_j X_{i,j}}}^{(\ic)}}&=\overline{\sim_{b_{\mid_{\alpha(i,j)} X_{i,j}}}^{(\ic)}}\\
&=\sim^{(\ic)}_{\calC_+\left(b_{\mid_{\alpha(i,j)}X_{i,j}},\sim\right)}\\
&=\sim^{(\ic)}_{\calC_+(b,\sim)_{\mid_{\alpha(i,j)}X_{i,j}}}\\
&=\sim^{(\ic)}_{\left(\calC_+(b,\sim)_{\mid_i X_i}\right)_{\mid_j X_{i,j}}}.
\end{align*}
As $\left(\calC_+(b,\sim)_{\mid_i X_i}\right)_{\mid_j X_{i,j}}$ is indecomposable by construction, $\sim_{\left(b_{\mid_i X_i}\right)_{\mid_j X_{i,j}}}^{(\ic)}$ has a unique class, so $\left(b_{\mid_i X_i}\right)_{\mid_j X_{i,j}}$ is indecomposable. Therefore,
\[(X_{1,1}/{\sim_{1,1}},\ldots,X_{k,l_k}/{\sim_{k,l_k}})\in\calC_+(b,\sim),\]
and finally, $(X_1/{\sim_1},\ldots,X_k/{\sim_k})\in m(\calC_+(b,\sim))$, as
\[(X_1/{\sim_1},\ldots,X_k/{\sim_k})\leq (X_{1,1}/{\sim_{1,1}},\ldots,X_{k,l_k}/{\sim_{k,l_k}}).\]

\textit{Fourth observation}. Let $X$ be a finite set, $b\in \calB[X]$, $\sim\subseteq \sim'\in \calE[X]$. 
\begin{align*}
&\sim\in \calE[b]\mbox{ and }\overline{\sim'}\in \calE[\calC_+(b,\sim)]\\
&\Longleftrightarrow \exists (X_1,\ldots,X_l)\in \calC_-(b),\:X/{\sim}=\{X_1,\ldots,X_l\},\: \forall i\in [l],\: b_{\mid_i X_i}\mbox{ is indecomposable}\\
&\mbox{and }\exists (X'_1/{\sim'_1},\ldots,X'_k/{\sim'_k})\in m(\calC_+(b,\sim)),\: (X/{\sim})/{\sim'}=\left\{X'_1/{\sim'1},\ldots,X'_k/{\sim'_k}\right\},\\
&\forall i\in [k],\: \calC_+(b,\sim)_{\mid_i X'_i/{\sim'_i}}\mbox{ is indecomposable}.
\end{align*}
By the compatibility between $\Delta$ and $\delta$,
\[ \calC_+(b,\sim)_{\mid_i X'_i/{\sim'_i}}=\calC_+\left(b_{\mid_i X'_i},\sim'_i\right).\]
By the conservation of the number of indecomposable components, $\calC_+(b,\sim)_{\mid_i X'_i/{\sim'_i}}$ is indecomposable if, and only if, $b_{\mid_i X'_i}$ is indecomposable.
Combining with the third observation, we obtain
\begin{align*}
&\sim\in \calE[b]\mbox{ and }\overline{\sim'}\in \calE[\calC_+(b,\sim)]\\
&\Longleftrightarrow \exists (X_1,\ldots,X_l)\in \calC_-(b),\:X/{\sim}=\{X_1,\ldots,X_l\},\: \forall i\in [l],\: b_{\mid_i X_i}\mbox{ is indecomposable}\\
&\mbox{and }\exists (X'_1,\ldots,X'_k)\in \calC_-(b),\: X/{\sim'}=\{X'_1,\ldots,X'_k\},\: \forall i\in [k],\: b_{\mid_i X'_i} \mbox{ is indecomposable}.
\end{align*}
On the other hand,
\begin{align*}
&\sim'\in \calE[b]\mbox{ and }\sim\in \calE[b\mid\sim']\\
&\Longleftrightarrow
\exists (X'_1,\ldots,X'_k)\in \calC_-(b),\: X/{\sim'}=\{X'_1,\ldots,X'_k\},\: \forall i\in [k],\: b_{\mid_i X'_i} \mbox{ is indecomposable},\\
&\mbox{and }\exists (X_1,\ldots,X_k)\in \calC_-(b\mid\sim'),\: X/{\sim}=\{X_1,\ldots,X_l\},\: \forall i\in [l],\: (b\mid\sim')_{\mid_i X_i}\mbox{ is indecomposable}.
\end{align*}

Let us assume that $\sim'\in \calE[b]$ and $\sim\in \calE[b\mid\sim']$. Let $(X_1,\ldots,X_l)$ and $(X'_1,\ldots,X'_l)$ as before. For $i\in [k]$, let us put 
\[\{X_j\mid X_j\subseteq X'_i\}=\{X_{i,1},\ldots,X_{i,l_i}\},\]
where $X_{i,1},\ldots,X_{i,l_i}$ are indexed in the same order as they appear in $(X_1,\ldots,X_l)$. 
As $b\mid\sim'=b_{\mid_1 X'_1}\cdots b_{\mid_k X'_k}$, by multiplicativity of $\Delta$, we obtain that
\[\left(\Delta_{X_{1,1},\ldots,X_{i,l_1}}\otimes\cdots \Delta_{X_{k,1},\ldots,X_{k,l_k}}\right)\circ \Delta_{X'_1,\ldots,X'_k}(b)\]
is obtained from $\Delta_{X_1,\ldots,X_k}(b\mid\sim')$ by permutation of tensors. Therefore, it is a tensor of indecomposable elements. As $\{X_{1,1},\ldots,X_{k,l_k}\}=\{X_1,\ldots,X_l\}=X/{\sim}$, we obtain that
$\sim\in \calE[b]$ and $\overline{\sim'}\in \calE[\calC_+(b,\sim)]$.

Moreover, for certain indices $\alpha(i,j)$, giving the position of $X_{i,j}$ in $(X_{1,1},\ldots,X_{k,l_k})$,
by coassociativity of $\Delta$,
\[(b\mid\sim')\mid\sim=\prod_{i=1}^k \prod_{j=1}^{l_i} \left(b_{\mid_i X'_i}\right)_{\mid_j X_{i,j}}=
\prod_{i=1}^k \prod_{j=1}^{l_i} b_{\mid_{\alpha(i,j)} X_{i,j}}=b\mid\sim.\]

$\Longrightarrow$. Let us assume that $\zeta$ is a coaction. The first observation implies the first item. 
Let $X$ be a finite set, $b\in \calB[X]$ and $\sim\in \calE[b]$. We suppose that $\ic\left(\calC_+(b,\sim)\right)\neq \ic(b)$. The second observation implies that $\sim^{(\ic)}_{\calC_+(b,\sim)}\subsetneq \overline{\sim_b^{(\ic)}}$.
Let us put $\sim'=\sim_b^{(\ic)}$, then $\sim\subseteq \sim'$, and
\begin{align*}
\left(\id_{\exch_+[X/{\sim'}]}\otimes \zeta_\sim\right)\circ \zeta_{\sim'}(b)&=\calC_+(b,\sim')\otimes \delta_\sim(b)\neq 0,
\end{align*}
whereas
\begin{align*}
\left(\delta_{\overline{\sim'}}\otimes \id_{\calB[X]}\right)\circ \zeta_\sim(b)&=\delta_{\overline{\sim'}}(\calC_+(b,\sim))\otimes b\mid\sim=0.
\end{align*}
Indeed, we do not have $\overline{\sim'}\subseteq \sim_{\calC_+(b,\sim)}^{(\ic)}$, as $\sim^{(\ic)}_{\calC_+(b,\sim)}\subsetneq \overline{\sim_b^{(\ic)}}$, so, by the compatibility between $\delta$ and $\squplus$, $\delta_{\overline{\sim'}}(\calC_+(b,\sim))=0$. This contradicts the fact that $\zeta$ is a coaction.

$\Longleftarrow$. Let $X$ be a finite set, $b\in \calB[X]$, $\sim\subseteq \sim'\in \calE[X]$. By hypothesis and the fourth observation we obtain that 
\[(\sim \calE[b]\mbox{ and }\sim'\in \calE[b/{\sim}])\Longleftrightarrow (\calE[b]\mbox{ and }\overline{\sim}\in \calE[\calC_+(b,\sim')]).\]
Let us assume that these conditions hold. By the fourth observation, $(b\mid\sim')\mid\sim=b\mid\sim$.
Let $\pi=(X_1/{\sim'_1},\ldots,X_n/{\sim'_n})\in \comp[X/{\sim'}]$, that is to say $(X_1,\ldots,X_n)\in \comp[X]$ and $\sim'_i\in \calE[X_i]$ for $i\in [n]$, with $\sim'=\sim'_1\sqcup\cdots\sqcup \sim'_n$. 
For $i\in [n]$, we put $\sim_i=\sim\cap X_i^2$. As $\sim\subseteq \sim'$, $\sim=\sim_1\sqcup\cdots \sqcup \sim_n$. 
Then
\begin{align*}
\pi \in \calC_+(b,\sim)/{\overline{\sim'}}&\Longleftrightarrow(X_1/{\sim_1},\ldots,X_n/{\sim_n})\in m(\calC_+(b,\sim))
\mbox{ and }\forall i\in [n],\: \sim^{(\ic)}_{\calC_+(b,\sim)_{\mid_{i X_i}}}=\overline{\sim'_i}.
\end{align*}
By the compatibility between $\Delta$ and $\zeta$, $\calC_+(b,\sim)_{\mid_i X_i}=\calC_+\left(b_{\mid_i X_i},\sim_i'\right)$. By the second observation, 
\[\sim^{(\ic)}_{\calC_+(b,\sim)_{\mid_{i X_i}}}=\overline{\sim_{b_{\mid_i X_i}}^{(\ic)}}.\]
Combining with the third observation,
\begin{align*}
\pi \in \calC_+(b,\sim)/{\overline{\sim'}}&\Longleftrightarrow (X_1,\ldots,X_n)\in \calC_-(b)\mbox{ and }
\forall i\in [n],\: \sim_{b_{\mid_i X_i}}^{(\ic)}=\sim'_i\\
&\Longleftrightarrow \pi\in \calC_+(b,\sim').
\end{align*}
Therefore, $\calC_+(b,\sim)/{\overline{\sim'}}= \calC_+(b,\sim')$.

Let $\pi=(X_1/{\sim_1},\ldots,X_n/{\sim_n})\in \comp[X/{\sim}]$, that is to say $(X_1,\ldots,X_n)\in \comp[X]$ and $\sim_i\in \calE[X_i]$ for $i\in [n]$, with $\sim=\sim_1\sqcup\cdots\sqcup \sim_n$. 
As $\sim'\in \calE[b]$, there exists $(X'_1,\ldots,X'_k)\in \calC_-(b)$, such that $b_{\mid_i X_i}$ is indecomposable for $i\in [k]$ and $X/{\sim'}=\{X'_1,\ldots,X'_k\}$.
\begin{align*}
\pi\in \calC_+(b,\sim)\mid\overline{\sim'}&\Longleftrightarrow \pi \in \calC_+(b,\sim)_{\mid_1 X'_1/{\sim'_1}}\sqcup\cdots\sqcup\calC_+(b,\sim)_{\mid_k X'_k/{\sim'_k}}\\
&\Longleftrightarrow \forall j\in [k],\: \elim\left(X_1/{\sim_1}\cap X'_j/{\sim'_j},\ldots,X_n/{\sim_n}\cap X'_j/{\sim'_j}\right)\in \calC_+(b,\sim)_{\mid_j X'_j/{\sim'_j}}\\
&\Longleftrightarrow \forall j\in [k],\: \elim\left(X_1/{\sim_1}\cap X'_j/{\sim'_j},\ldots,X_n/{\sim_n}\cap X'_j/{\sim'_j}\right)\in \calC_+\left(b_{\mid_j X'_j},\sim'_j\right)\\
&\Longleftrightarrow \forall j\in [k],\:\elim(X_1\cap X'_j,\ldots,X_n\cap X'_j)\in \calC_-\left(b_{\mid_j X'_j}\right) \\
&\mbox{and }\forall i\in [n], \forall j\in [k],\:\sim\cap (X_i\cap X'_j)^2=\sim^{(\ic)}_{\left(b_{\mid_j X'_j}\right)_{\mid_i X_i\cap X'_j}}\\
&\Longleftrightarrow (X_1,\ldots,X_n)\in \calC_-\left(\prod_{i=1}^k b_{\mid_j X'_j}\right)\\
&\mbox{and }\forall i\in [n], \forall j\in [k],\:\sim\cap (X_i\cap X'_j)^2=\sim^{(\ic)}_{\left(b_{\mid_j X'_j}\right)_{\mid_i X_i\cap X'_j}}\\
&\Longleftrightarrow (X_1,\ldots,X_n)\in \calC_-(b\mid\sim')\mbox{ and }\forall i\in [n], \sim_i=\sim^{(\ic)}_{(b\mid\sim')_{\mid_i X_i}}\\
&\Longleftrightarrow \pi \in \calC_+(b\mid\sim',\sim).
\end{align*}
We used the compatibility between $\zeta$ and $\Delta$ for the third equivalence. Therefore, $\calC_+(b,\sim)\mid\overline{\sim'}=\calC_+(b\mid\sim',\sim)$.

It remains to show that for any finite set $X$, any $b\in \calB[X]$, any $\sim\in \calE[b]$, $\calC_+(b,\sim)\in \exch_+^{\double}[X/{\sim}]$.
Let $X$ be a finite set. We consider the set $\calD[X]$ of elements of the form $\exch_+[\sigma](\calC_+(b,\sim))$, where $Y$ is a finite set, $b\in \calB[Y]$, $\sim\in \calE[b]$, and $\sigma:Y/{\sim}\longrightarrow X$ a bijection. 
This defines a subspecies of $\exch_+^{\qc}$.
Let $X$ and $X'$ be disjoint finite sets, $I=\exch_+[\sigma](\calC_+(b,\sim))\in \calD[X]$, $I'=\exch_+[\sigma'](\calC_+(b',\sim'))\in \calD[X']$, with $b\in \calB[Y]$ and $b'\in \calB[Y']$. Up to an isomorphism, we can assume that $Y$ and $Y'$ are disjoint. Then
\begin{align*}
I\sqcup I'&=\exch_+[\sigma](\calC_+(b,\sim))\sqcup \exch_+[\sigma'](\calC_+(b',\sim'))\\
&=\exch_+[\sigma \sqcup \sigma'](\calC_+(b,\sim)\sqcup \calC_+(b',\sim'))\\
&=\exch_+[\sigma \sqcup \sigma'](\calC_+(bb',\sim\sqcup\sim'))\\
&\in \calD[X\sqcup X'].
\end{align*}
So $\calD$ is a twisted subalgebra of $\exch_+^{\qc}$.

Let $X$ and $X'$ be disjoint finite sets, $I\in \calD[X\sqcup X']$, such that $(X,X')\in m(I)$. We put $I=\exch_+[\sigma\sqcup \sigma'](\calC_+(b,\sim))$, where $\sigma:Y/{\sim}\longrightarrow X$ and $\sigma':Y/{\sim'}\longrightarrow X'$ are two bijections. Then
\begin{align*}
\Delta_{X,X'}(I)&=(\exch_+[\sigma]\otimes \exch_+[\sigma'])\circ \Delta_{Y/{\sim},Y'/{\sim'}}(\calC_+(b,\sim))\\
&=(\exch_+[\sigma]\otimes \exch_+[\sigma'])\left(\calC_+\left(b_{\mid Y},\sim\right)\otimes \calC_+\left(b_{\mid Y'},\sim'\right)\right)\\
&\in \calD[X]\otimes \calD'[X'].
\end{align*} 
We used the compatibility between $\delta$ and $\Delta$ for the second equality. Consequently, $\calD$ is a twisted subbialgebra of $\exch_+^{\qc}$.

Let $I=\exch_+[\sigma](\calC_+(b,\sim))\in \calD[X]$, $\sim''\in \calE[I]$. We put $\overline{\sim'}=\calE[\sigma^{-1}](\sim'')\in \calE[Y/{\sim}]$ and by $\sim'\in \calE[Y]$ the corresponding equivalence such that $\sim\subseteq \sim'$. 
For a well-chosen bijection $\overline{\sigma}:Y/{\sim'}\longrightarrow X/{\sim''}$,
\begin{align*}
\delta_{\sim''}(I)\otimes b\mid \sim&=\left(\exch_+[\overline{\sigma}]\otimes \exch_+[\sigma]\otimes \id_{\calB[Y]}\right)\circ \left(\delta_{\overline{\sim'}}\otimes \id_{\calB[Y]}\right)\circ \zeta_\sim(b)\\
&=\left(\exch_+[\overline{\sigma}]\otimes \exch_+[\sigma]\otimes \id_{\calB[Y]}\right)\circ \left(\id_{\exch_+[Y/{\sim'}]}\otimes \zeta_\sim \right)\circ \zeta_{\sim'}(b)\\
&=\left(\exch_+[\overline{\sigma}]\otimes \exch_+[\sigma]\otimes \id_{\calB[Y]}\right)\left(\calC_+(b,\sim')\otimes \calC_+(b\mid\sim',\sim)\otimes (b\mid\sim')\mid\sim\right).
\end{align*}
We deduce that $I/{\sim''}=\exch_+[\overline{\sigma}](\calC_+(b,\sim'))\in \calD[X/{\sim''}]$. So $\calD$ is stable under $\delta$.

Let us now prove the coassociativity of $\delta$ on $\calD$. Let $I=\exch_+[\sigma](\calC_+(b,\sim))\in \calD[X]$, $\sim'\subseteq \sim''\in \calE[X]$. For a good choice of bijections $\sigma_i$,
\begin{align*}
\left(\delta_{\overline{\sim''}}\otimes \id_{\exch_+[X]}\right)\circ \delta_{\sim'}(I)\otimes b\mid\sim&=\left(\exch_+[\sigma_1]\otimes\exch_+[\sigma_2]\otimes \exch_+[\sigma_3] \otimes \calB[Y]\right)\\
&\circ\left(\delta_{\overline{\sim''}}\otimes \id_{\calB[Y/{\sim}]}\otimes \id_{\calB[Y]}\right)\circ \left(\delta_{\overline{\sim'}}\otimes \id_{\calB[Y]}\right)\circ \zeta_\sim(b)\\
\end{align*}
Moreover, 
\begin{align*}
&\left(\delta_{\overline{\sim''}}\otimes \id_{\calB[Y/{\sim}]}\otimes \id_{\calB[Y]}\right)\circ \left(\delta_{\overline{\sim'}}\otimes \id_{\calB[Y]}\right)\circ \zeta_\sim(b)\\
&=\left(\delta_{\overline{\sim''}}\otimes \id_{\calB[Y/{\sim}]}\otimes \id_{\calB[Y]}\right)\circ \left(\id_{\exch_+[Y/{\sim'}]}\otimes \zeta_\sim\right)\circ \zeta_{\sim'}(b)\\
&=\left(\id_{\exch_+[Y/\sim'']}\otimes \id_{\exch_+[Y/{\sim'}]}\otimes \zeta_\sim\right)\circ \left(\delta_{\overline{\sim''}}\otimes \id_{\calB[Y]}\right)\circ \zeta_{\sim'}(b)\\
&=\left(\id_{\exch_+[Y/\sim'']}\otimes \id_{\exch_+[Y/{\sim'}]}\otimes \zeta_\sim\right)\circ \left(\id_{\calB[Y/{\sim''}]}\otimes \zeta_{\sim'}\right)\circ \zeta_{\sim''}(b)\\
&=\left(\id_{\exch_+[Y/\sim'']}\otimes\delta_{\overline{\sim'}}\otimes \id_{\calB[Y]}\right)\circ \left(\id_{\calB[Y/{\sim''}]}\otimes \zeta_\sim\right)\circ \zeta_{\sim''}(b)\\
&=\left(\id_{\exch_+[Y/\sim'']}\otimes\delta_{\overline{\sim'}}\otimes \id_{\calB[Y]}\right)\circ \left(\delta_{\sim''}\otimes \id_{\calB[Y]}\right)\circ \zeta_\sim(b).
\end{align*}
Therefore,
\begin{align*}
\left(\delta_{\overline{\sim''}}\otimes \id_{\exch_+[X]}\right)\circ \delta_{\sim'}(I)\otimes b\mid\sim&=\left(\exch_+[\sigma_1]\otimes\exch_+[\sigma_2]\otimes \exch_+[\sigma_3] \otimes \calB[Y]\right)\\
&=\circ\left(\id_{\exch_+[Y/\sim'']}\otimes\delta_{\overline{\sim'}}\otimes \id_{\calB[Y]}\right)\circ \left(\delta_{\sim''}\otimes \id_{\calB[Y]}\right)\circ \zeta_\sim(b)\\
&=\left(\id_{\exch_+[X/\sim'']}\otimes\delta_{\overline{\sim'}}\otimes \id_{\calB[Y]}\right)\circ \left(\delta_{\sim''}\otimes \id_{\calB[Y]}\right)\circ \zeta_\sim(b)\\
&=\left(\id_{\exch_+[X/\sim'']}\otimes\delta_{\overline{\sim'}}\right)\circ delta_{\sim''}(I)\otimes b\mid \sim.
\end{align*}
Therefore, $\left(\delta_{\overline{\sim''}}\otimes \id_{\exch_+[X]}\right)\circ \delta_{\sim'}(I)=\left(\id_{\exch_+[X/\sim'']}\otimes\delta_{\overline{\sim'}}\right)\circ \delta_{\sim''}(I)$. As a conclusion, $(\calD,\squplus,\Delta,\delta)$ is a CTDB.
By Theorem \ref{theo4.3.1}, $\calD\subseteq \exch_+^{\double}$. Therefore, for any finite set $X$, and $b\in \calB[X]$, and any $\sim\in \calE[b]$, $\calC_+(b,\sim)\in \exch_+^{\double}[X/{\sim}]$.
\end{proof}

\begin{remark}
If $\calB$ is cocommutative, then the first item of Proposition \ref{prop4.3.9} is automatically satisfied, if the second item is satisfied. We use the notation of the proof of the fourth observation of the proof of Proposition \ref{prop4.3.9}.
Let us assume that $\sim\in \calE[b]$ and $\overline{\sim'}\in \calE[\calC_+(b,\sim)]$.
By cocommutativity, up to a permutation of the $X_i$'s, we can assume that 
\begin{align*}
X'_1&=X_1\sqcup \cdots \sqcup X_{l_1},&&\cdots&X'_k=X_{l_1+\cdots+l_{k-1}+1}\sqcup \cdots \sqcup X_{l_1+\cdots+l_k}.
\end{align*}
Then, by multiplicativity of $\Delta$,
\begin{align*}
\Delta_{X_1,\ldots,X_k}(b\mid \sim')&=\Delta_{X_1,\ldots,X_k}(b_{\mid_1X'_1}\cdots b_{\mid_k X'_k})\\
&=\left(\Delta_{X_1,\ldots,X_{l_1}}\otimes\cdots\otimes \Delta_{X_{l_1+\cdots+l_{k-1}+1},\ldots,X_{l_1+\cdots+l_k}}\right)\circ \Delta_{X'_1,\ldots,X'_l}(b)\\
&=\Delta_{X_1,\ldots,X_k}(b).
\end{align*}
So $\sim'\in \calE[b]$ and $\sim\in \calE[b\mid\sim']$. 
\end{remark}

\chapter{Internal products and sub-objects}

\section*{Introduction}

We introduce in this chapter more exotic structures on CTBs, using the Hadamard product of species. 
If $\calP$ and $\calQ$ are two species, the Hadamard product of $\calP$ and $\calQ$ is given on any finite set $X$ by
\[\calP\boxtimes \calQ[X]=\calP[X]\times \calQ[X].\]
The Hadamard product of two CTBs is also a CTB (Proposition \ref{prop5.1.2}), but this is not true for CTDBs.
An internal product on a CTB $\calB$ is defined as an associative morphism of CTBs from $\calB\boxtimes \calB$ to $\calB$ (Definition \ref{defi5.1.4}). 
For example, $\exch$ and its sub-CTBs $\exch_+$ and $\exch_-$ have an internal product, given by the intersection (Proposition \ref{prop5.1.5});
the CTBs of graphs, of hypergraphs, and of oriented graphs have also internal products, given by the union of set of hyperedges or arcs (Propositions \ref{prop5.1.6} and \ref{prop5.1.7}). Another example is given by quasi-posets, with the transitive closure (Proposition \ref{prop5.1.8}).
A major difference of this notion with the preceding ones is that it is not compatible with the bosonic Fock functor.

Another application of the Hadamard product is the notion of coloration of a CTB (Definition \ref{defi5.2.1}):
if $\Omega$ is a set and $\calB$ is a CTB, we can define another CTB $\calF_\Omega[\calB]=\calB\boxtimes \calD_\Omega$,
whose objects are the objects of $\calB$, given an extra coloration by elements of $\Omega$ (Proposition \ref{prop5.2.2}). 
If $\calB$ is a CTDB and $\Omega$ is a commutative semigroup, then $\calF_\Omega[\calB]$ is a CTDB (Proposition \ref{prop5.2.5}). An example is given by quasi-posets, which can be obtained from posets by a coloration by elements of $\N_{>0}$.
This formalism is used to study the formal flow associated to a CTB $\calB$. A similar task has been done in \cite{Bandiera2025} for free Novikov algebras (objects which can also be lifted to the level of species, but not as a CTB) and for posets.
We decorate $\calB$ with elements $\omega_1,\ldots,\omega_k$ of a set $\Omega$, and define in this way $k$ infinitesimal characters $\iota_i$ of $\calB$. We then define a character of $\calF_\Omega[\calB]$ by
\[f_\Omega=\exp(\iota_1)*\cdots*\exp(\iota_k),\]
and the formal flow is the polynomial invariant of $\calF_\Omega[\calB]$ associated to $f_\Omega$, that is to say $F_\Omega^\calB=\phi_{f_\Omega}^{\calF_\Omega[\calB]}$ (Definition \ref{defi5.2.8}).
As a consequence, the infinitesimal generator of $F_\Omega^\calB$ is 
\begin{align*}
\dfrac{\mathrm{d}F_\Omega^\calB}{\mathrm{d}T}(0)&=\ln\left(\exp(\iota_1)*\cdots*\exp(\iota_k)\right),
\end{align*}
which shows the strong link between this construction and the Campbell-Baker-Hausdorff formula in the group of characters of $\calB$.
We give a combinatorial interpretation of the values of the formal flows when evaluated to an element of $\N$,
leading to the notion of $\calB$-partitions associated to a decorated element of $\calB$ (Theorem \ref{theo5.2.9}).
In the case of posets, we recover the classical notion of $\calP$-partitions \cite{Malvenuto1993}, see Example \ref{ex5.2.3}.
We finally prove that the formal flow is compatible with $\Gamma_-^\calB$, and with $\Gamma_+^\calB$ when $\calB$ is non-degenerate (Corollary \ref{cor5.2.11}).

The second section is devoted to the definition of sub-objects of a CTB. The first one is the greatest cocommutative sub-CTB $\cocom\calB$ of a given CTB $\calB$, whose existence is proved in Proposition \ref{prop5.3.2}.
This can be rather small: for example, for the CTB $\calP$ of posets, we only obtain the discrete elements.
It is more interesting for $\exch$, and we give a description of $\cocom\exch$ with the help of an invariance condition under the action of symmetric groups (Definition \ref{defisigmainvariance}), which leads to a description of $\cocom\exch$ in terms of sets of equivalences. 
The two other sub-objects are inspired by remarkable properties of famous families of polynomials. 
The first example is given by the reciprocity principle for Ehrhart polynomials: it is well-known that the Ehrhart polynomial $\mathrm{Ehr}(P)$ and the strict Ehrhart polynomial $\mathrm{Ehr}^s(P)$ of a poset $P$ on a set of cardinality $n$
(and more generally, of certain polytopes of dimension $n$) are related by the formula
\[\mathrm{Ehr}^s(P)(-T)=(-1)^n \mathrm{Ehr}(P)(T).\]
In other words, for the CTB of posets $\calP$, for any finite set $X$,
\begin{align*}
&\forall P\in \calP[X],&\phi_+^\calP(P)(-T)&=(-1)^{|X|}\phi_-^\calP(P).
\end{align*}
We shall say that $\calP$ satisfies the reciprocity principle (Definition \ref{defi5.3.3}).
We prove that any CTB $\calB$ has a sub-CTB $\calB^{\mathrm{recip}}$ satisfying the reciprocity principle, maximal for the inclusion, and that this construction is functorial (Proposition \ref{prop5.3.6}). 
As a consequence, the sub-CTB $\exch_+^{\mathrm{recip}}$ gives all the polynomial invariants with the reciprocity principle,
and we show that there are more such examples than the Ehrhart polynomials of posets, see Examples \ref{ex5.3.2} and \ref{ex5.3.3}.
The second example is given by the alternating sign condition, satisfied by chromatic polynomials of graphs (Definition \ref{defi5.3.8}). 
We prove that any CTB $\calB$ has a sub-CTB $\calB^{\mathrm{alt}}$ satisfying the alternating sign condition, maximal for the inclusion, and that this construction is functorial (Proposition \ref{prop5.3.9}). This condition will be explored with more details in the cocommutative case in the forthcoming chapter \ref{chap7}.

\section{Hadamard product of CTBs}

\subsection{Construction}

\begin{defi}
Let $\calP$ and $\calQ$ be two species. The Hadamard product of $\calQ$ and $\calQ$ is the species $\calP\boxtimes\calQ$ defined as follows: for any finite set $X$,
\[\calP\boxtimes\calQ[X]=\calP[X]\times \calQ[X].\]
\end{defi}

\begin{prop} \label{prop5.1.2}
Let $\calB$ and $\calB'$ be two CTBs. Then $\calB\boxtimes\calB'$ is also a CTB, with the following product and coproduct:
\begin{itemize}
\item If $X$, $Y$ are disjoint finite sets, for any $a\in \calB[X]$, $a'\in \calB'[X]$, $b\in \calB[Y]$, $b'\in \calB[Y]$,
\[m_{X,Y}((a,b)\otimes (a',b'))=m_{X,Y}(a\otimes b),m_{X,Y}(a',b').\]
\item If $X$, $Y$ are disjoint finite sets, for any $b\in \calB[X\sqcup Y]$, $b'\in \calB'[X\sqcup Y]$,
\[\Delta_{X,Y}(b,b')=\begin{cases}
\left(b_{\mid_LX},b'_{\mid_LX}\right)\otimes\left(b_{\mid_RY},b'_{\mid_RY}\right)\mbox{ if }(X$, $Y)\in \calC_-(b)\cap \calC_-(b'),\\
0\mbox{ otherwise}.
\end{cases}\]\end{itemize}
\end{prop}

\begin{proof}
Direct verifications.
\end{proof}

\begin{example}
The CTB of mixed graphs which can be deduced from the construction of \cite{Foissy45} is a sub-object of $\calG\boxtimes \calO\calG$: a mixed graph $G=(V(G),E(G),A(G))$ is identified with the pair $((V(G),E(G)),(V(G),A(G))$. 
We do not obtain the whole $\calG\boxtimes\calO\calG$, as we impose constraints on arcs and edges for mixed graphs (namely, we cannot have an edge $\{x,y\}$ and an arc $(x,y)$ in the same mixed graph). This injection is a morphism of CTBs.
\end{example}

\begin{example}\label{ex5.1.2}
Let $\calT\calO$ be the species of total orders, with the product and coproduct defined as follows:
\begin{itemize}
\item If $X$, $Y$ are disjoint finite sets, $\leq_X\in \calT\calO[X]$, and $\leq_Y\in \calT\calO[Y]$, then $\leq_X\cdot \leq_Y$ is the ordinal sum of $\leq_X$ and $\leq_Y$:
\[\leq_X\cdot \leq_Y=\leq_X \sqcup \leq_Y\sqcup (X\times Y).\]
\item If $X$, $Y$ are disjoint finite sets, and $\leq\in \calT\calO[X\sqcup Y]$, then 
\[\Delta_{X,Y}(\leq)=\leq_{\mid X}\otimes \leq_{\mid Y}.\]
\end{itemize}
If$ \calB$ is a CTB, then for any finite set $X$, $\calB\boxtimes \calT\calO[X]$ is made of pairs $(b,\leq)$, where $b\in \calB[X]$ and $\leq$ is a total order on $X$. Therefore, $\calF[\calB\boxtimes\calT\calO]$ has for basis 
\[\bigsqcup_{n\geq 0} \calB[n],\]
identifying $b\in \calB[n]$ with the isoclass of $(b,\leq)$, where $\leq$ is the usual order on $[n]$. Therefore, $\calF[\calB\boxtimes\calT\calO]=\calF_F[\calB]$ is the full Fock functor of \cite{Aguiar2010} applied to $\calB$.
For example, for any finite set $X$, $\calP\boxtimes \calT\calO[X]$ is the set of pairs $(\leq_P,\leq_T)$, where $\leq_P$ is a partial order on $X$ and $\leq_T$ is a total order on $X$. We obtain the CTB of special posets of \cite{Foissy20}.
\end{example}

\begin{remark}
If $\calB$ and $\calB'$ are CTDBs, one can define a coproduct $\delta$ on $\calB\boxtimes\calB'$ as follows: if $X$ is a finite set, $\sim\in \calE[X]$, $b\in \calB[X]$ and $b'\in \calB'[X]$,
\[\delta_\sim(b,b')=\begin{cases}
(b/{\sim},b'/{\sim})\otimes (b\mid\sim,b'\mid\sim)\mbox{ if }\sim\in \calE[b]\cap \calE[b'],\\
0\mbox{ otherwise}.
\end{cases}\]
Unfortunately, this generally does not give $\calB\boxtimes \calB'$ a structure of CTDB. For example, take $\calB=\calB'=\calG$, $X=[2]$, and $G,H\in \calG[2]$ defined by
\begin{align*}
E(G)&=\emptyset,&E(H)=\{[2]\}.
\end{align*}
For any $\sim\in \calE[2]$,
\[\delta_\sim(G,H)=\begin{cases}
(G,H)\otimes (G,G)\mbox{ if }\sim=\sim_G^{(m)},\\
0\mbox{ otherwise}.
\end{cases}\]
So $\delta$ has no right counit. 
\end{remark}

\begin{prop}\label{prop5.1.3}
Let $\calB$ and $\calB'$ be two CTBs. Let $X$ be a finite set, $b\in \calB[X]$, and $b'\in \calB'[X]$.
\begin{enumerate}
\item $\calC_-(b,b')=\calC_-(b)\cap\calC_-(b')$.
\item Let us assume that $\calB$ or $\calB'$ is commutative. Then $\sim^{(\ic)}_{(b,b')}$ is the transitive closure of $\sim_b^{(\ic)}\cup \sim_{b'}^{(\ic)}$.
\item Let us assume that $\calB$ or $\calB'$ is commutative. Then $\calC_+(b,b')=\calC_+(b)\cap\calC_+(b')$.
\end{enumerate}\end{prop}

\begin{proof}
1. Comes directly of the definition of $\Delta$ on $\calB\boxtimes\calB'$.

2. We assume for example that $\calB$ is commutative. We denote by $\sim$ the transitive closure of $\sim_b^{(\ic)}\cup \sim_{b'}^{(\ic)}$. As $\sim_{b'}^{(\ic)}\subseteq \sim$, we can write $b'=b'_1\cdots b'_n$, with $X/{\sim}=\{X_1,\ldots,X_n\}$,
and $b'_i\in \calB'[X_i]$ for $i\in [k]$. By commutativity of $\calB$, as $\sim_b^{(\ic)}\subseteq \sim$, we can also write
$b=b_1\cdots b_n$, with $b_i\in \calB[X_i]$ for $i\in [n]$. Then
\[(b,b')=(b_1\cdots b_n,b'_1\cdots b'_n)=(b_1,b'_1)\cdots (b_n,b'_n),\]
so
\[\sim_{(b,b')}^{(\ic)}=\bigcup_{k=1}^n \sim_{(b_k,b'_k)}^{(\ic)}\subseteq \bigcup_{k=1}^n X_k^2=\sim.\]
So $\sim_{(b,b')}^{(\ic)}\subseteq \sim$. Moreover, we can write
\[(b,b')=(a_1,a'_1)\cdots (a_m,a'_m),\]
with $X/{\sim}_{(b,b')}^{(\ic)}=\{X'_1,\ldots,X'_m\}$ and $(a_i,a'_i)\in \calB\boxtimes\calB'[X'_i]$ for $i\in [m]$. Then
\begin{align*}
b&=a_1\cdots a_m,&b'&=a'_1\cdots a'_m.
\end{align*}
Therefore, $\sim_b^{(\ic)}\subseteq \sim_{(b,b')}^{(\ic)}$ and $\sim_{b'}^{(\ic)}\subseteq \sim_{(b,b')}^{(\ic)}$.
As $\sim_{(b,b')}^{(\ic)}$ is transitive, it contains the transitive closure of $\sim_b^{(\ic)}\cup \sim_{b'}^{(\ic)}$.
Therefore, $\sim\subseteq \sim_{(b,b')}^{(\ic)}$.

3. \textit{First step}. Let us prove that $\bfD(\calB\boxtimes \calB')=\bfD\calB\boxtimes \bfD\calB'$.
Let $X$ be a finite set. We denote by $\sim_X^{(m)}$ the equality of $X$. Let $b\in \calB[X]$ and $b'\in \calB'[X]$.
By the second item,
\begin{align*}
(b,b')\in \bfD(\calB\boxtimes \calB')[X]&\Longleftrightarrow \sim_{(b,b')}^{(\ic)}=\sim_X^{(m)}\\
&\Longleftrightarrow\mbox{the transitive closure of $\sim_b^{(\ic)}\cup \sim_{b'}^{(\ic)}$ is $\sim_X^{(m)}$}\\
&\Longleftrightarrow \sim_b^{(\ic)}=\sim_{b'}^{(\ic)}=\sim_X^{(m)}\\
&\Longleftrightarrow b\in \bfD\calB[X]\mbox{ and }b'\in \bfD\calB'[X]\\
&\Longleftrightarrow (b,b')\in \bfD\calB\boxtimes \bfD\calB'[X].
\end{align*}

\textit{Second step}. Let $(X_1,\ldots,X_k)\in \comp[X]$. By the first step,
\begin{align*}
(X_1,\ldots,X_k)\in \calC_+(b,b')&\Longleftrightarrow (X_1,\ldots,X_k)\in \calC_-(b,b')\\
&\mbox{ and }\Delta_{X_1,\ldots,X_k}(b,b')\in \bfD(\calB\boxtimes\calB')[X_1]\otimes\cdots\otimes \bfD(\calB\boxtimes\calB')[X_k]\\
&\Longleftrightarrow (X_1,\ldots,X_k)\in \calC_-(b)\cap \calC_-(b')\\
&\mbox{ and }\Delta_{X_1,\ldots,X_k}(b,b')\in \bfD\calB\boxtimes\bfD\calB'[X_1]\otimes\cdots\otimes \bfD\calB\boxtimes\bfD\calB'[X_k]\\
&\Longleftrightarrow (X_1,\ldots,X_k)\in \calC_-(b)\mbox{ and }(X_1,\ldots,X_k)\in \calC_-(b')\\
&\mbox{ and }\Delta_{X_1,\ldots,X_k}(b)\in \bfD\calB[X_1]\otimes\cdots\otimes\bfD\calB[X_n]\\
&\mbox{ and }\Delta_{X_1,\ldots,X_k}(b')\in \bfD\calB'[X_1]\otimes\cdots\otimes\bfD\calB'[X_n]\\
&\Longleftrightarrow (X_1,\ldots,X_k)\in \calC_+(b)\cap \calC_+(b'). \qedhere
\end{align*}\end{proof}

\begin{example}
If $\calB$ and $\calB'$ are not commutative, it may happen that $\calC_+(b,b')\subsetneq\calC_+(b)\cap\calC_+(b')$.
Here is an example. Let us consider the CTB $\calT\calO$ of Example \ref{ex5.1.2}. 
Consider $\leq,\leq'\in \calT\calO[2]$, defined by
\begin{align*}
&1\leq 2,&& 2\leq' 1.
\end{align*}
Then 
\begin{align*}
\leq&=\leq_{\mid \{1\}}\cdot \leq_{\mid \{2\}}\in \bfD\calT\calO[2],&\leq'&=\leq'_{\mid \{2\}}\cdot \leq'_{\mid \{1\}}\in \bfD\calT\calO[2],
\end{align*}
but $(\leq,\leq')\notin \bfD(\calT\calO\boxtimes \calT\calO)[2]$: otherwise, we would have
\begin{align*}
(\leq,\leq')&=(\leq,\leq')_{\mid \{1\}}\cdot (\leq,\leq')_{\mid \{2\}}=(\leq,\leq),\\
\mbox{ or }(\leq,\leq')&=(\leq,\leq')_{\mid \{2\}}\cdot (\leq,\leq')_{\mid \{1\}}=(\leq',\leq'),
\end{align*}
which is not true, as $\leq\neq\leq'$. 
\end{example}

\subsection{Internal products}

\begin{defi}\label{defi5.1.4}
Let $\calB$ be a CTB. An internal product is an associative morphism $\divideontimes :\calB\boxtimes \calB\longrightarrow \calB$ of twisted CTBs, that is to say:
\begin{itemize}
\item For any finite set $X$, $\divideontimes :\ell\calB[X]\otimes \ell\calB[X]\longrightarrow \ell \calB[X]$ is an associative and unitary product. It unit will be denoted by $I_X$. If $\sigma:X\longrightarrow Y$ is a bijection between two finite sets, then the following diagram commutes:
\[\xymatrix{\ell\calB[X]\otimes \ell\calB[X]\ar[r]^{\hspace{6mm}\divideontimes }\ar[d]_{\ell\calB[\sigma]\otimes \ell\calB[\sigma]}&\ell\calB[X]\ar[d]^{\ell\calB[\sigma]}\\
\ell\calB[Y]\otimes \ell\calB[Y]\ar[r]_{\hspace{6mm}\divideontimes }&\ell\calB[Y]}\]
Moreover, $1_\calB\divideontimes 1_\calB=1_\calB$.
\item For any finite set $X$, and any $a,b\in \calB[X]$, $a\divideontimes b=0$ or belongs to $\calB[X]$. 
\item $\divideontimes $ is compatible with the product $m$: for any disjoint finite sets $X$ and $Y$, $a,b\in \calB[X]$, $a',b'\in \calB[Y]$,
\[m_{X,Y}(a\divideontimes b\otimes a'\divideontimes b')=m_{X,Y}(a\otimes a')\divideontimes m_{X,Y}(b\otimes b').\]
\item $\divideontimes $ is compatible with the coproduct: for any disjoint finite sets $X$ and $Y$, $a,b\in \calB[X\sqcup Y]$,
\[\Delta_{X,Y}(a\divideontimes b)=\Delta_{X,Y}(a)\divideontimes \Delta_{X,Y}(b).\]
\end{itemize}\end{defi}

\begin{prop}\label{prop5.1.5}
Let $(X,Y)$ be a pair of disjoint finite sets. We put
\begin{align*}
&\forall I\in \exch[X],\: \forall J\in \exch[Y],&I\divideontimes J&=\begin{cases}
I\cap J\mbox{ if }I\cap J\neq \emptyset,\\
0\mbox{ otherwise.}
\end{cases}
\end{align*}
Then $\divideontimes$ is an internal product on $\exch$. Moreover, $\exch_+$ and $\exch_-$ are stable under $\divideontimes$.
\end{prop}

\begin{proof}
Let us consider the character $\epsilon\boxtimes\epsilon$ of $\exch\boxtimes\exch$: if $X$ is a finite set,
\begin{align*}
&\forall I,J\in \exch[X],&\epsilon\boxtimes\epsilon(I,J)&=\begin{cases}
1\mbox{ if }(X)\in I\cap J,\\
0\mbox{ otherwise}.
\end{cases}
\end{align*}
We denote by $\divideontimes:\exch\boxtimes \exch\longrightarrow \exch$ the unique morphism of CTBs such that $\epsilon\circ \divideontimes=\epsilon\boxtimes \epsilon$.
Let $X$ be a finite set, and $(I,J)\in \exch\boxtimes\exch[X]$. Let us show that $\calC_{\epsilon\boxtimes\epsilon}(I,J)=I\cap J$. 

$\subseteq$. Let $(X_1,\ldots,X_k)\in \calC_{\epsilon\boxtimes\epsilon}(I,J)$. Then
\[(\epsilon\boxtimes \epsilon)^{\otimes k}\circ \Delta_{X_1,\ldots,X_k}(I,J)=1.\]
Therefore, $(X_1,\ldots,X_k)\in \calC_-(I,J)$ and for any $i\in [k]$, $(X_i)\in I_{\mid_i X_i}\cap J_{\mid_i X_i}$. By the exchange condition applied to $I$ and $J$, $(X_1,\ldots,X_k)\in I\cap J$. 

$\supseteq$. Let $(X_1,\ldots,X_k)\in I\cap J$. Then $(X_1,\ldots,X_k)\in \calC_-(I)$, and $(X_i)\in I_{\mid_i X_i}$ for any $i\in [k]$. Similarly, $(X_1,\ldots,X_k)\in \calC_-(J)$, and $(X_i)\in J_{\mid_i X_i}$ for any $i\in [k]$. 
As a consequence, $(\epsilon\boxtimes\epsilon)\left(I_{\mid_i X_i}\otimes J_{\mid_i X_i}\right)=1$, so $(X_1,\ldots,X_k)\in \calC_{\epsilon\boxtimes \epsilon}(I,J)$. 

Therefore, $\calC_[\epsilon\boxtimes \epsilon](I,J)=I\cap J$, and
\[I\divideontimes J=\begin{cases}
I\cap J\mbox{ if }I\cap J\neq \emptyset,\\
0\mbox{ otherwise.}
\end{cases}\]
For any $I,J,K\in \exch[X]$,
\[(I\divideontimes J)\divideontimes K=I\divideontimes(J\divideontimes K)=\begin{cases}
I\cap J\cap K\mbox{ if }I\cap J\cap K\neq \emptyset,\\
0\mbox{ otherwise}.
\end{cases}\]
So $\divideontimes$ is associative.

If $I,J\in \exch_\pm[X]$, it is not difficult to show that $I\cap J$ is, either empty, or in $\exch_\pm[X]$. So $\exch_+$ and $\exch_-$ are stable under $\divideontimes$.
\end{proof}

\begin{remark}
The product $*$ is usually not defined on the classical bialgebra $\calF[\calB]$. Here is an example. We consider
\begin{align*}
I&=\{(1,2)\}\in \exch[2],&J&=\{(2,1)\}\in \exch[2].
\end{align*}
Then $I$ and $J$ are isomorphic, through the transposition $(1\:2)$. Therefore, in $\calF[\exch]$, $\overline{I}=\overline{J}$. But $I*J=0$ and $I*I=I$. 
\end{remark}

\begin{prop}\label{prop5.1.6}
Let $X$ be a finite set, $G,H\in \calH\calG[X]$. We define a hypergraph $G\divideontimes H\in \calH\calG[X]$ by
\[\HE(G\divideontimes H)=\HE(G)\cup \HE(H).\]
Then $\divideontimes $ is an internal product of $(\calH\calG,m,\Delta^{\subset,\subset})$. Moreover, $\Gamma^{\calH\calG}_+:\calH\calG\longrightarrow\exch_+$ and $\Gamma^{\calH\calG}_-:\calH\calG\longrightarrow\exch_-$ compatible with $\divideontimes $.
\end{prop}

\begin{proof}
The product $\divideontimes $ is obviously associative. For any finite set $X$, the unit is the hypergraph on $X$ with no hyperedge. Let $X$, $Y$ be disjoint finite sets, $G,H\in \calH\calG[X]$, $G',H'\in \calH\calG[Y]$.
\begin{align*}
\HE((G\divideontimes H)(G'\divideontimes H'))&=\HE(G\divideontimes H)\sqcup \HE(G'\divideontimes H')\\
&=(\HE(G)\cup \HE(H))\sqcup (\HE(G')\cup \HE(H'))\\
&=(\HE(G)\sqcup \HE(G'))\cup (\HE(H)\sqcup \HE(H'))\\
&=\HE(GG')\cup \HE(HH')\\
&=\HE((GG')\divideontimes (HH')),
\end{align*}
so $(G\divideontimes H)(G'\divideontimes H')=(GG')\divideontimes (HH')$: the product $\divideontimes $ is compatible with the product of hypergraphs.

Let $X$, $Y$ be disjoint finite sets, and $G,H\in \calH\calG[X\sqcup Y]$. 
\begin{align*}
\HE((G\divideontimes H)_{\mid Y})&=\{Z\in \HE(G\divideontimes H)\mid Z\subseteq Y\}\\
&=\{Z\in \HE(G)\cup \HE(H)\mid Z\subseteq Y\}\\
&=\{Z\in \HE(G)\mid Z\subseteq Y\}\cup \{Z\in \HE(H)\mid Z\subseteq Y\}\\
&=\HE(G_{\mid Y})\cup \HE(H_{\mid Y})\\
&=\HE(G_{\mid Y}\divideontimes H_{\mid Y}).
\end{align*}
So $(G\divideontimes H)_{\mid Y}=G_{\mid Y}\divideontimes H_{\mid Y}$. This implies the compatibility of $\divideontimes $ with $\Delta^{\subset,\subset}$.

Let $X$ be a finite set, $G,H\in \calH\calG[X]$. Let $(X_1,\dots,X_k)\in \comp[X]$.
\begin{align*}
&(X_1,\ldots,X_k)\in \calC_+(G\divideontimes H)\\
&\Longleftrightarrow \forall i\in [k],\: X_i\mbox{ does not contain any non-trivial edge of $G\divideontimes H$}\\
&\Longleftrightarrow \forall i\in [k],\: X_i\mbox{ does not contain any non-trivial edge of $G$, nor of $H$}\\
&\Longleftrightarrow (X_1,\ldots,X_k)\in \calC_+(G)\cap \calC_+(H)\\
&\Longleftrightarrow (X_1,\ldots,X_k)\in \calC_+(G)\divideontimes \calC_+(H).
\end{align*}
Therefore, $\Gamma^{\calH\calG}(G\divideontimes H)=\calC_+(G\divideontimes H)=\calC_+(G)\divideontimes \calC_+(H)=\Gamma^{\calH\calG}_+(G)\divideontimes \Gamma^{\calH\calG}_+(H)$.

Let $X$ be a finite set, and $G,H\in \calH\calG[X]$. Then
\[\Gamma^{\calH\divideontimes\calG}_-(GH)=\comp[X]=\comp[X]\cap \comp[X]=\comp[X]\divideontimes\comp[X]=\Gamma^{\calH\calG}(G)\divideontimes \Gamma^{\calH\calG}(H),\]
so $\Gamma^{\calH\divideontimes\calG}_-$ is also compatible with $\divideontimes$.
\end{proof}

\begin{remark}
The CTB $\calG$ of graphs, seen as a sub-CTDB of $\calH\calG$, is stable under $\divideontimes $, so also inherits an internal product, compatible with $\Gamma^\calG_+$ and $\Gamma^\calG_-$.
\end{remark}

\begin{prop}\label{prop5.1.7}
Let $X$ be a finite set, $G,H\in \calO\calG[X]$. We define an oriented graph $G\divideontimes H\in \calO\calG[X]$ by
\[A(G\divideontimes H)=A(G)\cup A(H).\]
Then $\divideontimes $ is an internal product of $\calO\calG$. Moreover, $\Gamma^{\calO\calG}_+:\calO\calG\longrightarrow\exch_+$ and $\Gamma^{\calO\calG}_-:\calO\calG\longrightarrow\exch_-$ are compatible with $\divideontimes $.
\end{prop}

\begin{proof}
The proof that $\divideontimes $ is an internal product on $\calO\calG$ is very similar to the proof for $\calH\calG$, and we let it to the reader. For any finite set $X$, $I_X$ is the oriented graph on $X$ with no arc. Let us prove the compatibility with $\Gamma^{\calO\calG}_+$. Let $X$ be a finite set, and $G,H\in \calO\calG[X]$.
Then $\calC_+(G\divideontimes H)$ is the set of compositions $(X_1,\ldots,X_k)\in \comp[X]$ such that, for any $x,y\in X$,
\[(x,y)\in A(G\divideontimes H)\Longrightarrow x\in X_i,\: y\in X_j,\: \mbox{ with }i<j.\]
As $A(G\divideontimes H)=A(G)\cup A(H)$, $\calC_+(G\divideontimes H)$ is the set of compositions $(X_1,\ldots,X_k)\in \comp[X]$ such that, for any $x,y\in X$,
\begin{align*}
(x,y)\in A(G)&\Longrightarrow x\in X_i,\: y\in X_j,\: \mbox{ with }i<j\\
\mbox{ and }(x,y)\in A(H)&\Longrightarrow x\in X_i,\: y\in X_j,\: \mbox{ with }i<j,
\end{align*}
that is to say $\calC_+(G)\cap \calC_+(H)$. Therefore, 
\begin{align*}
\Gamma^{\calO\calG}_+(G\divideontimes H)&=\begin{cases}
\calC_+(G)\cap \calC_+(H)\mbox{ if }\calC_+(G)\cap\calC_+(H)\neq \emptyset,\\
0\mbox{ otherwise},
\end{cases}\\
&=\Gamma^{\calO\calG}_+(G)\divideontimes \Gamma^{\calO\calG}_+(H). 
\end{align*}
Similarly, $\calC_-(G\divideontimes H)$ is the set of compositions $(X_1,\ldots,X_k)\in \comp[X]$ such that, for any $x,y\in X$,
\begin{align*}
(x,y)\in A(G)&\Longrightarrow x\in X_i,\: y\in X_j,\: \mbox{ with }i\leq j\\
\mbox{ and }(x,y)\in A(H)&\Longrightarrow x\in X_i,\: y\in X_j,\: \mbox{ with }i\leq j,
\end{align*}
that is to say $\calC_-(G)\cap \calC_-(H)$. Therefore, $\Gamma_-^{\calO\calG}$ is compatible with $\divideontimes$.
\end{proof}

\begin{prop}\label{prop5.1.8}
Let $X$ be a finite set, $P,Q\in \calQ\calP[X]$. We define a quasi-poset $P\divideontimes Q\in \calQ\calP[X]$ by
\[\leq_{P\divideontimes Q}=\leq_P\vee \leq_Q,\]
where $\vee$ is the transitive closure of the union. Then $\divideontimes $ is an internal product of $\calQ\calP$. 
Moreover, $\phi_{\calO\calG\longrightarrow \calQ\calP}$, $\Gamma^{\calQ\calP}_+:\calQ\calP\longrightarrow\exch_+$, and $\Gamma^{\calQ\calP}_-:\calQ\calP\longrightarrow\exch_-$ are all compatible with $\divideontimes $. 
\end{prop}

\begin{proof}
Let us first prove that for any finite set $X$, any $G,H\in \calO\calG[X]$,
\[\phi_{\calO\calG}(G\divideontimes H)=\phi_{\calO\calG}(G)\divideontimes \phi_{\calO\calG}(H).\]
Firstly, $\phi_{\calO\calG}(G\divideontimes H)$ is the quasi-poset given by the transitive closure of $A(G\divideontimes H)=A(G)\cup A(H)$. Secondly, $\phi_{\calO\calG}(G)\divideontimes \phi_{\calO\calG}(H)$ is the quasi-poset given by the transitive closure of the union of transitive closures of $A(G)$ and $A(H)$. 
It is not difficult to show that they are equal.

As $\phi_{\calO\calG}:\calO\calG\longrightarrow \calQ\calP$ is a surjective morphism of CTBs and is compatible with $\divideontimes $, Proposition \ref{prop5.1.7} implies that $\divideontimes $ is an internal product for $\calQ\calP$.
For any finite set $X$, $I_X$ is the equality of $X$. As $\Gamma^{\calQ\calP}_\pm\circ \phi_{\calO\calG}=\Gamma^{\calO\calG}_\pm$ and $\Gamma^{\calO\calG}_\pm$ is compatible with $\divideontimes $, $\Gamma^{\calQ\calP}_\pm$ is also compatible with $\divideontimes $. 
\end{proof}

\begin{remark}
$\calP$ is not stable under $\divideontimes $. For example, consider $P,Q\in \calP[2]$, defined by
\begin{align*}
&1\leq_P 2,&&2\leq_Q 1.
\end{align*}
Then $\leq_P\vee \leq_Q=[2]^2\notin \calP[2]$. 
\end{remark}

\section{Colorations and formal flow}

\subsection{Coloration}

Let us give an example of endofunctor constructed from the Hadamard product of species.

\begin{defi}\label{defi5.2.1}
Let $\Omega$ be a non-empty set. For any finite set $X$, let us denote by $\calD_\Omega[X]$ the set of maps $\ell:X\longrightarrow \Omega$. This defines a species $\calD_\Omega$, which is made a CTB as follows:
\begin{enumerate}
\item Let $X$ and $Y$ be disjoint finite sets, $\ell_X\in \calD_\Omega[X]$ and $\ell_Y\in \calD_\Omega[Y]$. We define $\ell_X\cdot\ell_Y:X\sqcup Y\longrightarrow \Omega$ by
\begin{align*}
(\ell_X\cdot \ell_Y)_{\mid X}&=\ell_X,&(\ell_X\cdot \ell_Y)_{\mid Y}&=\ell_Y.
\end{align*}
The unit is the unique element of $\calD_\Omega[\emptyset]$.
\item Let $X$ and $Y$ be disjoint finite sets, and $\ell\in \calD_\Omega[X\sqcup Y]$. We put $\Delta_{X,Y}(\ell)=\ell_{\mid X}\otimes \ell_{\mid Y}$.
\end{enumerate} \end{defi}

\begin{proof}
Direct verifications.
\end{proof}

\begin{remark}
\begin{enumerate}
\item If $|\Omega|=1$, then $\calD_\Omega$ is isomorphic to $\Com$.
\item $\nd\calD_\Omega=\calD_\Omega$.
\end{enumerate}
\end{remark}
This allows to define the endofunctor associated to $\Omega$:

\begin{prop}\label{prop5.2.2}
Let $\calB$ be a CTB, and let $\Omega$ be a non-empty set. We put $\calF_\Omega[\calB]=\calB\boxtimes\calD_\Omega$: this defines an endofunctor of the category of CTBs. More precisely:
\begin{enumerate}
\item For any finite set $X$, $\calF_\Omega[\calB][X]$ is the set of pairs $(b,\ell)$, where $b\in \calB[X]$ and $\ell:X\longrightarrow \Omega$.
\item For any pair $(X,Y)$ of disjoint finite sets, and any $a\in \calB[X]$, $b\in \calB[Y]$, $\ell_X\in \calD_\Omega[X]$, $\ell_Y\in \calD_\Omega[Y]$,
\[(a,\ell_X)(b,\ell_Y)=(ab,\ell_X\cdot\ell_Y).\]
\item For any pair $(X,Y)$ of disjoint finite sets, and any $b\in \calB[X\sqcup Y]$, $\ell\in \calD_\Omega[X\sqcup Y]$,
\[\Delta_{X,Y}(b,\ell)=\begin{cases}
(b_{\mid_LX},\ell_{\mid X})\otimes (b_{\mid_RY},\ell_{\mid Y})\mbox{ if }(X,Y)\in \calC_-(b),\\
0\mbox{ otherwise}.
\end{cases}\]
\item We denote by $S_\calB$ the antipode of $\calB$. Then the antipode of $\calF_\Omega[\calB]$ is $S_\calB\boxtimes \id_{\calD_\Omega}$. 
\end{enumerate} \end{prop}

\begin{proof}
It remains to prove the formula for the antipode. We put $T=S_\calB\boxtimes \id_{\calD_\Omega}$. Note that $T\left(1_{\calF_\Omega[\calB]}\right)=1_{\calF_\Omega[\calB]}$. 
Let $X$ be a nonempty finite set, and $(b,\ell)\in \calF_\Omega[\calB]$. Then
\begin{align*}
T*\id_{\calF_\Omega[S_\calB]}(b,\ell)&=T(b,f)+(b,f)+\sum_{(X',X'')\in \calC_-(b)} T\left(b_{\mid_LX'},\ell_{\mid X''}\right)\left(b_{\mid_RX''},\ell_{\mid X''}\right)\\
&=(S_\calB(b),f)+(b,f)+\sum_{(X',X'')\in \calC_-(b)}\left(S_\calB\left(b_{\mid_LX'}\right)b_{\mid_RX''},f_{\mid X'}f_{\mid X''}\right)\\
&=(S_\calB(b),f)+(b,f)+\sum_{(X',X'')\in \calC_-(b)}\left(S_\calB\left(b_{\mid_LX'}\right)b_{\mid_RX''},f\right)\\
&=\left(S_\calB(b)+b+\sum_{(X',X'')\in \calC_-(b)}S_\calB\left(b_{\mid_LX'}\right)b_{\mid_RX''},f\right)\\
&=0.
\end{align*}
Similarly, $\id_{\calF_\Omega[\calB]}*T(b,\ell)=0$. So $T$ is the antipode of $\calF_\Omega[\calB]$.
\end{proof}

\begin{example}\label{ex5.2.1}
Let us consider $\calP\boxtimes \calD_{\mathbb{N}_*}$. Then $\calF[\calP\boxtimes \calD_{\mathbb{N}_*}]$ has for basis the set of isoclasses of $(P,\ell)$, where $P=(V(P),\leq_P)$ is a poset and $\ell:V(P)\longrightarrow \N_{>0}$ is a map.
Such an element is identified with the isoclass of a quasi-poset $Q$, with the same Hasse graph, the decorations on the vertices indicating the cardinalities of the classes of $\sim_Q$. As a conclusion, $\calF_{\N_{>0}}[\calP]$ is isomorphic, as a bialgebra, to $\calF[\calQ\calP]$.
\end{example}

If $\calB$ is a CTDB, then it is possible to make $\calF_\Omega[\calB]$ also a CTDB, under the condition to give $\Omega$ a commutative semigroup. Let us detail this construction. We shall need a definition and a technical Lemma.

\begin{defi}
We assume that $(\Omega,+)$ is a commutative semigroup. Let $X$ be a finite set, $\sim\in \calE[X]$, and $\ell\in \calD_\Omega[X]$. We define $\ell/{\sim}\in \calD_\Omega[X/{\sim}]$ by
\begin{align*}
&\forall Y\in X/{\sim},&\ell/{\sim}(Y)&=\sum_{y\in Y}\ell(y).
\end{align*} 
This sum has to be understood in $(\Omega,+)$, and is well-defined as $+$ is associative and commutative.
\end{defi}

\begin{lemma}\label{lem5.2.4}
Let $X$ be a finite set, $\sim,\sim'\in \calE[X]$, such that $\sim\subseteq \sim'$, and $\ell\in \calD_\Omega[X]$. Then
\[(\ell/{\sim})/{\overline{\sim'}}=\ell/{\sim'}.\]
\end{lemma}

\begin{proof}
Let $Z\in (X/{\sim})/{\overline{\sim}'}$. We identify it with $Z'\in X/{\sim'}$.
\begin{align*}
(\ell/{\sim})/{\overline{\sim'}}(Z)&=\sum_{Y\in Z}\ell/{\sim}(Y)=\sum_{Y\in Z}\sum_{y\in Y}\ell(y)=\sum_{y\in Z'}\ell(y)=\ell/{\sim}(Z').
\end{align*}
So $(\ell/{\sim})/{\overline{\sim'}}=\ell/{\sim'}$.
\end{proof}

\begin{prop}\label{prop5.2.5}
Let $(\Omega,+)$ be commutative semigroup, and $\calB$ be a CTDB. Then $\calF_\Omega[\calB]$ is also a CTDB, with the coproduct $\delta$ defined as follows:
for any finite set $X$, and $\sim\in \calE[X]$, $b\in \calB[X]$, $\ell\in \calD_\Omega[X]$,
\[\delta_\sim(b,\ell)=\begin{cases}
(b/{\sim},\ell/{\sim})\otimes (b\mid\sim,\ell)\mbox{ if }\sim\in \calE[b],\\
0\mbox{ otherwise}.
\end{cases}\]
We denote respectively by $\phi_+^{calB,\Omega}$ and $\phi_+^{\calB}$ the fundamental polynomial invariants by $\calF_\Omega[\calB]$ and of $\calB$. For any finite set $X$, and $\ell:X\longrightarrow \Omega$ and any $b\in \calB[X]$,
\[\phi_+^{\calB,\Omega}(b,\ell)=\phi_+^\calB(b).\]
\end{prop}

\begin{proof}
Let us prove the coassociativity of $\delta$. Let $X$ be a finite set, $\sim,\sim'\in \calE[X]$, and $b\in \calB[X]$. If $\sim$ is not included in $\sim'$, then $(\id_{\calB[X/{\sim'}]}\otimes \delta_\sim)\circ\delta_{\sim'}(b)=0$. This implies that
\[\left(\id_{\calF_\Omega[\calB][X/{\sim'}]}\otimes \delta_\sim\right)\circ\delta_{\sim'}(b,\ell)=0.\]
Let us now assume that $\sim\subseteq\sim'$. By the coassociativity of $\delta$ on $\calB$ and Lemma \ref{lem5.2.4},
\begin{align*}
\left(\delta_{\overline{\sim}'}\otimes \id_{\calF_\Omega[\calB][X]}\right)\circ \delta_\sim(b,\ell)&=\begin{cases}
((b/{\sim})/{\overline{\sim}'},(\ell/{\sim})/{\overline{\sim}'})\otimes (b/{\sim}\mid \overline{\sim'},\ell/{\sim})\otimes (b\mid\sim,\ell)\\
\hspace{1cm}\mbox{ if }\sim\in \calE[b]\mbox{ and }\overline{\sim'}\in \calE[b/{\sim}],\\
0\mbox{ otherwise}
\end{cases}\\
&=\begin{cases}
(b/{\sim'},\ell/{\sim'})\otimes ((b\mid\sim')/{\sim},\ell/{\sim})\otimes ((b\mid\sim')\mid\sim,\ell)\\
\hspace{1cm}\mbox{ if }\sim'\in \calE[b]\mbox{ and }\sim\in \calE[b\mid\sim'],\\
0\mbox{ otherwise}
\end{cases}\\
&=\left(\id_{\calF_\Omega[\calB][X/{\sim'}]}\otimes \delta_\sim\right)\circ\delta_{\sim'}(b,\ell).
\end{align*}

Let us put, for any finite set $X$, $b\in \calB[X]$, $\ell\in \calD_\Omega[X]$, $\epsilon(b,\ell)=\epsilon(b)$. Then, for any $\sim\in \calE[X]$, as $\calB$ is a twisted double bialgebra,
\begin{align*}
\left(\id_{\calF_\Omega[\calB][X]}\otimes \epsilon\right)\circ \delta_\sim(b,\ell)&=\begin{cases}
(b,\ell)\mbox{ if $\sim$ is the equality of $X$},\\
0\mbox{ otherwise}.
\end{cases}\\
\sum_{\sim\in \calE[X]}\left(\epsilon\otimes \id_{\calF_\Omega[\calB][X]}\right)\circ \delta_\sim(b,\ell)&=\sum_{\sim\in \calE[b]}\epsilon(b/{\sim})(b,\ell)=(b,\ell).
\end{align*}
So $\epsilon$ is a counit for $\calF_\Omega[\calB]$.

Let us finally prove the compatibility between $\delta$ and $\Delta$. Let $X$, $Y$ be disjoint finite sets, $\sim_X\in \calE[X]$, $\sim_Y\in \calE[Y]$, $b\in \calB[X\sqcup Y]$ and $\ell\in \calD_\Omega[X\sqcup Y]$. We put $\sim=\sim_X\sqcup \sim_Y$. 
As $\calB$ is a twisted double bialgebra,
\begin{align*}
&\left(\Delta_{X/{\sim_X},Y/{\sim_Y}}\otimes \id_{\calF_\Omega[\calB][X]}\right)\circ \delta_{\sim_X\sqcup \sim_Y}(b,\ell)\\
&=\begin{cases}
\left((b/{\sim})_{\mid_L X/{\sim_X}},(\ell/{\sim})_{\mid X/{\sim_X}}\right)\otimes \left((b/{\sim})_{\mid_R Y/{\sim_Y}},(\ell/{\sim})_{\mid Y/{\sim_Y}}\right)\otimes (b\mid\sim,\ell)\\
\hspace{1cm}\mbox{ if $\sim\in \calE[b]$ and }(X/{\sim_X},Y/{\sim_Y})\in \calC_-(b/{\sim}),\\
0\mbox{ otherwise}
\end{cases}\\
&=\begin{cases}
\left((b_{\mid_LX})/{\sim_X},(\ell_{\mid X})/{\sim_X}\right)\otimes \left((b_{\mid_RY})/{\sim_Y},(\ell_{\mid Y})/{\sim_Y}\right)\otimes \left((b_{\mid_LX})\mid\sim_X(b_{\mid_RY})\mid\sim_Y,\ell\right)\\
\hspace{1cm}\mbox{ if $(X,Y)\in \calC_-(b)$, $\sim_X\in \calE[b_{\mid_LX}]$ and $\sim_Y\in \calE[b_{\mid_RY}]$},\\
0\mbox{ otherwise}
\end{cases}\\
&=\begin{cases}
\left((b_{\mid_LX})/{\sim_X},(\ell_{\mid X})/{\sim_X}\right)\otimes \left((b_{\mid_RY})/{\sim_Y},(\ell_{\mid Y})/{\sim_Y}\right)\otimes \left((b_{\mid_LX})\mid\sim_X,\ell_{\mid X}\right)\left((b_{\mid_RY})\mid\sim_Y,\ell_{\mid Y}\right)\\
\hspace{1cm}\mbox{ if $(X,Y)\in \calC_-(b)$, $\sim_X\in \calE[b_{\mid_LX}]$ and $\sim_Y\in \calE[b_{\mid_RY}]$},\\
0\mbox{ otherwise}
\end{cases}\\
&=m_{1,3,24}\circ(\delta_{\sim_X}\otimes \delta_{\sim_Y})\circ \Delta_{X,Y}(b).
\end{align*}
So $\calF_\Omega[\calB]$ is a CTDB. 

Let us observe that for any finite set $X$, any $\ell:X\longrightarrow \Omega$, and any $b\in \calB[X]$,
\[\epsilon_\delta(b,\ell)=\epsilon_\delta(b).\]
Therefore, the polynomial invariant defined by
\[\phi(b,\ell)=\phi_+^\calB(b)\]
is a polynomial invariant of $\calF_\Omega[\calB]$, such that for any finite set $X$, any $\ell:X\longrightarrow \Omega$, and any $b\in \calB[X]$,
\[\phi(b,\ell)(1)=\phi_+^\calB(b)(1)=\epsilon_\delta(b)=\epsilon_\delta(b,\ell).\]
Therefore, $\phi=\phi_+^{\calB,\Omega}$.
\end{proof}

\begin{remark}
Let $\phi:\calB\longrightarrow \calB'$ be a morphism of CTBs. Then $\calF_\Omega[\phi]$ send $(b,\ell)$ to $(\phi(b),\ell)$, for any finite set $X$, $b\in \calB[X]$ and $\ell\in \calD_\Omega[X]$.
\end{remark}

\begin{example}
Going back to Example \ref{ex5.2.1}, the isomorphism of bialgebras from $\calF[\calP\boxtimes \calD_{\N_{>0}}]$ to $\calF[\calQ\calP]$ is in fact an isomorphism of double bialgebras, when $\N_{>0}$ is given its usual addition. 
\end{example}

\subsection{Formal flow and $\calB$-partitions}

\begin{defi}\label{defi5.2.6}
Let $\calB$ be a CTB. Let $X$ be a finite set, of cardinality denoted by $n$, $b\in \calB[X]$, $f:X\longrightarrow [n]$ be a bijection, and $g:X\longrightarrow [N]$ be a map. 
We define a relation on $X$ as the lexicographic order associated to the pair $(g,f)$:
\begin{align*}
&\forall x,y\in X,& x\leq_{g,f} y\mbox{ if } &g(x)<g(y)\mbox{ or }(g(x)=g(y)\mbox{ and }f(x)\leq f(y)).
\end{align*}
As $f$ is a bijection, $\leq_{g,f}$ is a total order on $X$. 
We then write $X=\{x_1\leq_{g,f}\cdots \leq_{g,f} x_n\}$. We shall say that $g$ is a $\calB$-partition of $(b,f)$ if $(\{x_1\},\ldots,\{x_n\})\in \calC_{\pp}(b)$.
\end{defi}

\begin{example}\label{ex5.2.3}
Let us consider the CTB of posets $\calP$. Let $X$ be a finite set, and $P\in \calP[X]$. 
If $X=\{x_1,\ldots,x_n\}$, then $(\{x_1\},\ldots,\{x_n\})\in \calC_{\pp}(P)$ if, and only if,
\begin{align*}
&\forall i,j\in [n],&x_i\leq_P x_j&\Longrightarrow i\leq j.
\end{align*}
Let $f:X\longrightarrow [n]$ a bijection and $g:X\longrightarrow [N]$ be a map. Reformulating Definition \ref{defi5.2.6}, $g$ is a $\calP$-partition of $(P,f)$, if, and only if
\begin{align*}
&\forall x,y\in X,&x\leq_P y\Longrightarrow&g(x)<g(y)\\
&&&\mbox{ or }g(x)=g(y)\mbox{ and }f(x)\leq f(y). 
\end{align*}
This is equivalent to
\begin{align*}
&\forall x,y\in X,&x\leq_P y&\Longrightarrow g(x)\leq g(y),\\
&&(f(x)>f(y) \mbox{ and }x\leq_Py)&\Longrightarrow g(x)<g(y).
\end{align*}
We recover the notion of $\calP$-partitions of an indexed poset \cite{Malvenuto1993,Foissy27}.
\end{example}

\begin{example}
Let $X$ be a finite set, $I\in \exch[X]$ and $f:X\longrightarrow [n]$ be a bijection. 
An $\exch$-partition of $(I,f)$ is a map $g:X\longrightarrow [n]$ such that there exists $(\{x_1\},\ldots,\{x_n\})\in I$, with
\begin{align*}
&\forall i,j\in [n],& i\leq j&\Longrightarrow g(x_i)\leq g(x_j),\\
&&(i\leq j\mbox{ and }g(x_i)=g(x_j))&\Longrightarrow f(x_i)\leq f(x_j).
\end{align*}
This is equivalent to 
\begin{align*}
&\forall i,j\in [n],& i\leq j&\Longrightarrow g(x_i)\leq g(x_j),\\
&&(i\leq j\mbox{ and }f(x_i)>f(x_j))&\Longrightarrow g(x_i)<f(x_j).
\end{align*}
\end{example}

\begin{example}\label{ex5.2.5}
Let $b\in \calB[X]$, such that $\calC_{\pp}(b)=\comp_{|X|}[X]$, and let $f:X\longrightarrow [n]$ be a bijection. Then any map $g:X\longrightarrow [N]$ is a $\calB$-partition of $(b,f)$. 
\end{example}

Let us now introduce the formal flow. For this, we consider a finite, totally ordered set $\Omega=\{\omega_1<\cdots<\omega_k\}$. 

\begin{defi}\label{defi5.2.7}
Let $\calB$ be a CTB. For any $\omega\in \Omega$, we consider the infinitesimal character $\iota_\omega$ of $\calF_\Omega[\calB]$ defined as follows: for any finite set $X$, any $b\in \calB[X]$, and any $\ell\in \calD_\Omega[X]$,
\[\iota_\omega(b,\ell)=\begin{cases}
1\mbox{ if $|X|=1$ and $\ell(X)=\{\omega\}$},\\
0\mbox{ otherwise}.
\end{cases}\]
We then consider the character $f_\Omega=\exp(\iota_{\omega_1})*\cdots*\exp(\iota_{\omega_k})$. The formal flow $F_\Omega^\calB$ is the polynomial invariant $\phi_{f_\Omega}$ of $\calF_\Omega[\calB]$ associated to $f_\Omega$. 
\end{defi}

\begin{remark}
As a consequence, if $\lambda$ is the unique infinitesimal character of $\calB$ such that
\[\exp(\iota_{\omega_1})*\cdots*\exp(\iota_{\omega_k})=\exp(\lambda),\]
then 
\[\lambda(b,\ell)=\dfrac{\mathrm{d}F_\Omega^\calB(b,\ell)}{\mathrm{d}T}(0).\]
\end{remark}

\begin{defi}\label{defi5.2.8}
Let $X$ be a finite set, of cardinality denoted by $n$, $\ell\in \calD_\Omega[X]$, and $f:X\longrightarrow [n]$ a bijection.
We shall say that $f$ is $\ell$-compatible if 
\begin{align*}
&\forall x,y\in X,&\ell(x)<\ell(y)\Longrightarrow f(x)<f(y).
\end{align*}
The set of $\ell$-compatible bijections is denoted by $\CB(\ell)$.
\end{defi}

\begin{theo}\label{theo5.2.9}
Let $\calB$ be a CTB. Let $X$ be a finite set, of cardinality denoted by $n$, $b\in \calB[X]$, and $\ell\in \calD_\Omega[X]$. For any $N\in \N$,
\[F_\Omega^\calB(b,\ell)(N)=\frac{1}{\ell!}\sum_{f\in \CB(\ell)}|\{g:X\longrightarrow[N]\mid\mbox{$g$ is a $\calB$-partition of $(b,f)$}\}|,\]
where
\[\ell!=\prod_{\omega\in \Omega}|\ell^{-1}(\omega)|!\] 
\end{theo}

\begin{proof}
Notice that
\begin{align*}
&F_\Omega^\calB(b,\ell)\\
&=\sum_{(X_1,\ldots,X_l)\in\comp[X]} f_\Omega^{\otimes l}\circ \Delta_{X_1,\ldots,X_l}(b,\ell)\binom{T}{l}\\
&=\sum_{\substack{(X_1,\ldots,X_l)\in\comp[X],\\ f_i:X_i\longrightarrow[n_i],\mbox{\scriptsize bijective},\\m_{i,1}+\cdots+m_{i,k}=n_i}}
\left(\prod_{i=1}^l\prod_{k=1}^k m_{i,j}!\right)^{-1}\underbrace{(\iota_{\omega_1}^{\otimes m_{1,1}}\otimes\cdots\otimes \iota_{\omega_k}^{\otimes m_{1,k}}\otimes\cdots\otimes \iota_{\omega_1}^{\otimes m_{l,1}}\otimes\cdots\otimes \iota_{\omega_k}^{\otimes m_{l,k}})
\circ\Delta_f(b,\ell)}_{a_f}\binom{T}{l}
\end{align*}
where
\[\Delta_f=\Delta_{\{f_1^{-1}(1)\},\ldots,\{f_1^{-1}(n_1)\},\ldots,\{f_l^{-1}(1)\},\ldots,\{f_l^{-1}(n_l)\}.}\]
The term $a_f$ is 0, except if $(\{f_1^{-1}(1)\},\ldots,\{f_1^{-1}(n_1)\},\ldots,\{f_l^{-1}(1)\},\ldots,\{f_l^{-1}(n_l)\})\in \calC_-(b)$ and 
\[\left(\ell(f^{-1}_1(1)),\ldots,\ell(f^{-1}_l(m_l))\right)=(\underbrace{\omega_1,\ldots,\omega_1}_{m_{1,1}},\ldots,\underbrace{\omega_k,\ldots,\omega_k}_{m_{1,k}},\ldots,
\underbrace{\omega_1,\ldots,\omega_1}_{m_{l,1}},\ldots,\underbrace{\omega_k,\ldots,\omega_k}_{m_{l,k}}).\] 
In this case, $a_F=1$, and
\[\prod_{i=1}^l\prod_{k=1}^k m_{i,j}!=\prod_{\omega\in \Omega}\prod_{i=1}^l |\ell^{-1}(\omega)\cap X_i|!.\]
Consequently, for $a_F$ to be nonzero, it is necessary that 
\begin{align*}
&\forall i\in [l],\: \forall x,y\in X_i,&\ell(x)<\ell(y)\Longrightarrow f_i(x)<f_i(y).
\end{align*}
For such a family $(f_1,\ldots,f_l)$, the number of $f\in \CB(\ell)$, such that the standardization of $f_{\mid X_i}$ is equal to $f_i$ for any $i\in [l]$ is 
\[\ell!\left(\prod_{\omega\in \Omega}\prod_{i=1}^l |\ell^{-1}(\omega)\cap X_i|!\right)^{-1}.\]
Therefore, we can rewrite
\begin{align*}
F_\Omega^\calB(b,\ell)&=\frac{1}{\ell!}\sum_{f\in \CB(\ell)} \sum_{(X_1,\ldots,X_l)\in \calC_-(b,f)}\binom{T}{l},
\end{align*}
where $\calC_-(b,f)$ is the set of compositions $(X_1,\ldots,X_l)\in \calC_-(b)$, such that, putting $f(X_i)=\{x_{i,1}<\cdots<x_{i,n_i}\}$ for $i\in [l]$, then
\[(\{f^{-1}(x_{1,1})\},\ldots,\{f^{-1}(x_{1,n_1})\},\ldots,\{f^{-1}(x_{N,1})\},\ldots,\{f^{-1}(x_{N,n_N})\})\in \calC_-(b).\]
Sending a composition $(X_1,\ldots,X_l)$ to the map $g:X\longrightarrow [l]$ such that $g^{-1}(i)=X_i$ for any $i\in [l]$, the set $\calC_-(b,f)$ is now identified with the set of surjective $\calB$-partitions of $(b,f)$. Therefore, 
\begin{align*}
F_\Omega^\calB(b,\ell)&=\frac{1}{\ell!}\sum_{f\in \CB(\ell)}\sum_{\mbox{\scriptsize $g$ surjective $\calB$-partition of $(b,f)$}} \binom{T}{\max(g)}.
\end{align*}
The result immediately follows. 
\end{proof}

\begin{remark}
We obtain that for any finite set $X$, any map $\ell:X\longrightarrow \Omega$, and any $b\in \calB[X]$,
\[F_\Omega^\calB(b,\ell)=\frac{n!}{\ell!}\sum_{f\in \CB(\ell)} F_\Omega^\calB(b,f).\]
\end{remark}

\begin{example}
Let us apply this to the CTB of posets $\calP$. Let $P$ be a finite poset, of cardinality $n$, and $\ell:V(P)\longrightarrow[k]$ be a map. According to Example \ref{ex5.2.3}, for any $N\in \N$, $\ell!F_{[k]}^\calP(P,\ell)(N)$ is the number of pairs $(f,g)$ such that:
\begin{enumerate}
\item $f:V(P)\longrightarrow [n]$ is a bijection such that
\begin{align*}
&\forall x,y\in V(P),&\ell(x)<\ell(y)&\Longrightarrow f(x)<f(y).
\end{align*}
\item $g:V(P)\longrightarrow [N]$ is a $P$-partition of $(P,f)$:
\begin{align*}
&\forall x,y\in V(P),& x\leq_P y&\Longrightarrow g(x)\leq g(y),\\
&&(x\leq_P y\mbox{ and }f(x)>f(y))&\Longrightarrow g(x)<g(y).
\end{align*}
\end{enumerate}
In particular, if $\ell$ is bijective, $n!F_{[n]}^\calP(P,\ell)(N)$ is the number of $P$-partitions $g:V(P)\longrightarrow [N]$ of $(P,\ell)$, as the unique $\ell$-compatible map is $\ell$ itself.
This description of $F_\Omega^\calP$ has been firstly done by Ruggero Bandiera. It does not seem that his work on this subject is available yet, but this result has been the object of a talk at the workshop \emph{Algebraic Combinatorics and Finite Groups IV}, see
\href{https://sites.google.com/uniroma1.it/claudia-malvenuto/activities/cetraro-2026}{the website of the meeting}.
\end{example}

\begin{example}
Let $X$ be a finite set, $I\in \exch[X]$ and $\ell:X\longrightarrow \Omega$ be a map. Then, for any $N\geq 1$, $\ell!F_\Omega(I,\ell)(N)$ is the number of pairs $(f,g)$ such that:
\begin{enumerate}
\item $f:X\longrightarrow [n]$ is a bijection and $g:X\longrightarrow [n]$ is a map.
\item For any $x,y\in X$, if $\ell(x)<\ell(y)$, then $f(x)<f(y)$.
\item There exists $(\{x_1\},\ldots,\{x_n\})\in I$ such that
\begin{align*}
&\forall i,j\in [n],& i\leq j&\Longrightarrow g(x_i)\leq g(x_j),\\
&&(i\leq j\mbox{ and }f(x_i)>g(x_j))&\Longrightarrow g(x_i)<f(x_j).
\end{align*}\end{enumerate}\end{example}

\begin{example}\label{ex5.2.6}
Let us assume that $\calC_{\pp}(b)=\comp_n[X]$, where $n=|X|$. For any labeling $\ell$ of $X$,
\begin{align*}
F_\Omega^\calB(b,l)(N)&=\frac{1}{\ell!}\sum_{f\in \CB(\ell)} N^n=\frac{1}{\ell!}|\CB(\ell)|N^n=N^n,
\end{align*}
as $|\CB(\ell)|=\ell!$. Therefore,
\[F_\Omega^\calB(b,l)=T^n.\]
Consequently, if $\lambda$ is the unique infinitesimal character of $\calB$ such that
\[\exp(\iota_{\omega_1})*\cdots*\exp(\iota_{\omega_k})=\exp(\lambda),\]
then 
\[\lambda(b,\ell)=\dfrac{\mathrm{d}F_\Omega^\calB(b,\ell)}{\mathrm{d}T}(0)=\delta_{n,1}.\]
In other words, 
\[\lambda(b,\ell)=\left(\sum_{\omega\in \Omega} \iota_\omega\right)(b,\ell).\]
This is the case for any graph in $\calG$, or any hypergraph in $\calH\calG$, whatever the choice of the coproduct is. 
\end{example}

\begin{cor}
Let $\calB$ be a CTB, $X$ be a finite set, and $b\in \calB[X]$. If $\ell:X\longrightarrow \Omega$ is a constant map, then
\[F_\Omega(b,\ell)=\phi_{\pp}^\calB(b).\]
\end{cor}

\begin{proof}
Any bijection $f:X\longrightarrow [n]$ is $\ell$-compatible. As a consequence, for any $N\geq 1$,
\[F^\calB_\Omega(b,\ell)=\frac{1}{n!}\sum_{f:X\longrightarrow [n],\:\mbox{\scriptsize bijective}}|\{g:X\longrightarrow[N]\mid\mbox{$g$ is a $\calB$-partition of $(b,f)$}\}|.\]
Therefore, $n! F^\calB_\Omega(b,\ell)(N)$ is the number of pairs $(f,g)$, such that:
\begin{enumerate}
\item $f:X\longrightarrow [n]$ is a bijection and $g:X\longrightarrow [N]$ is a map.
\item Putting $X=\{x_1\leq_{g,f}\cdots \leq_{g,f} x_n\}$, then $(\{x_1\},\ldots,\{x_n\})\in \calC_{\pp}(b)$.
\end{enumerate}
For any $(\{x_1\},\ldots,\{x_n\})\in \comp_n[X]$, we denote by $I(\{x_1\},\ldots,\{x_n\})$ the set of pairs $(f,g)$ such that $x_1\leq_{g,f}\cdots \leq_{g,f}$. The cardinality of this set only depends on $N$ and $n$ and is denoted by $P_n(N)$. Reformulating, we obtain 
\[F^\calB_\Omega(b,\ell)(N)=\frac{1}{n!}\sum_{(\{x_1\},\ldots,\{x_n\})\in \calC_{\pp}(b)}
|I(\{x_1\},\ldots,\{x_n\})|
=\frac{1}{n!}|\calC_{\pp}(b)| P_n(N)=f_{\pp}^{\calB}(b)P_n(N).\]
Let us now compute $P_n(N)$. Let us apply what precedes to the CTB $\calG$ and any graph $G$ in $\calG[X]$. By Example \ref{ex5.2.6}, as $\calC_{\pp}(G)=\comp_n[X]$,
\[F^\calB_\Omega(G,\ell)=P_n(N)=N^n.\]
As a consequence, $F^\calB_\Omega(b,\ell)=f_{\pp}^\calB(b)T^n=\phi_{\pp}^\calB(b)$.
\end{proof}

\begin{lemma}
Let $\calB$ be a CTB, and $\omega\in \Omega$. Then $\iota_\omega\circ \calF_\Omega[\Gamma^\calB_-]=\iota_\omega$. If $\calB$ is non-degenerate, $\iota_\omega\circ \calF_\Omega[\Gamma^\calB_+]=\iota_\omega$. 
\end{lemma}

\begin{proof}
Let $X$ be a finite set, $b\in \calB[X]$ and $\ell\in \calD_\Omega[X]$. As $\calC_-(b)\in \exch_-[X]$,
\begin{align*}
\iota_\omega\circ \calF_\Omega[\Gamma^\calB_-](b,\ell)&=\begin{cases}
1\mbox{ if }|X|=1\mbox{ and }\ell(X)=\{\omega\},\\
0\mbox{ otherwise}
\end{cases}\\
&=\iota_\omega(b,\ell).
\end{align*}
If $\calB$ is non-degenerate, the proof is similar, as then $\Gamma^\calB_+(b)=\calC_+(b)\in \exch_+[X]$. 
\end{proof}

As a consequence, by unicity in Theorem \ref{theo1.2.23}:

\begin{cor}\label{cor5.2.11}
Let $\calB$ be a CTB. Then $F_\Omega^{\exch_-}\circ \calF_\Omega[\Gamma^\calB_-]=F_\Omega^\calB$. If $\calB$ is non-degenerate, $F_\Omega^{\exch_+}\circ \calF_\Omega[\Gamma^\calB_+]=F_\Omega^\calB$.
\end{cor}

\begin{prop}
Let us assume that $\calB$ is a CTDB, and let us denote by $\phi_+^\calB$ its fundamental polynomial invariant.
For any finite set $X$, any $\ell:X\longrightarrow \Omega$ and any $b\in \calB[X]$, we denote by $\lambda^\calB_{F,\Omega}(b,\ell)$ then number of elements $(\{x_1\},\ldots,\{x_n\})\in \calC_{\pp}(b)$ such that
\begin{align*}
&\forall i,j\in [n],&i\leq j&\Longrightarrow \ell(i)\leq\ell(j).
\end{align*}
Then
\[F_\Omega^\calB(b,\ell)=\sum_{\sim\in \calE[b]} \lambda^\calB_{F,\Omega}(b\mid \sim,\ell)\phi_+^\calB(b/{\sim})(T).\]
\end{prop}

\begin{proof}
Let us choose an associative and commutative product on $\Omega$. Then $\calF_\Omega[\calB]$ is a CTDB, 
by Proposition \ref{prop5.2.5}. The fundamental polynomial invariant of $\calF_\Omega[\calB]$ is given by
\[\phi_+^{\calB,\Omega}(b,\ell)=\phi_+^\calB(b).\]
Let us use then Proposition \ref{prop1.3.7}. We put $\lambda^\calB_{F,\Omega}(b,\ell)=F_\Omega^\calB(b,\ell)(1)$. By definition of a $\calB$-partition with $N=1$,
for any bijection $f:X\longrightarrow [n]$, the constant map $X\longrightarrow [1]$ is a $\calB$-partition of $(b,\ell)$ if, and only if, 
\[(\{f^{-1}(1)\},\ldots,\{f^{-1}(n)\})\in \calC_{\pp}(b).\]
Therefore, $\lambda^\calB_{F,\Omega}(b,\ell)$ is the number of $\ell$-compatible bijections $f$ such that 
\[(\{f^{-1}(1)\},\ldots,\{f^{-1}(n)\})\in \calC_{\pp}(b).\]
It is equal to the number of elements $(\{x_1\},\ldots,\{x_n\})\in \calC_{\pp}(b)$, such that
\begin{align*}
&\forall i,j\in [n],&i\leq j&\Longrightarrow \ell(i)\leq\ell(j).
\end{align*}
Then Proposition \ref{prop1.3.7} gives that
\[F_\Omega^\calB(b,\ell)=(\phi_+^{\calB,\Omega}\curvearrowleft \lambda^\calB_{F,\Omega})(b,\ell)=
\sum_{\sim\in \calE[b]} \lambda^\calB_{F,\Omega}(b\mid \sim,\ell)\phi_+^\calB(b/{\sim})(T). \qedhere\]
\end{proof}

\begin{remark}
We put $\Omega=\{\omega_1<\cdots<\omega_k\}$. Then $\lambda^\calB_{F,\Omega}(b,\ell)$ is the number of elements $\pi\in \calC_{\pp}(b)$, such that 
\[\pi\geq \elim(\ell^{-1}(\omega_1),\ldots,\ell^{-1}(\omega_k)).\]
\end{remark}

\begin{example}
Let $P\in \calP[X]$ and $\ell:X\longrightarrow [\Omega]$. Then $\lambda^\calB_{F,\Omega}(P,\ell)$ is the number of total orders $\leq'_P$ containing $\leq_P$ such that
\begin{align*}
&\forall x,y\in X,&x\leq'_P y&\Longrightarrow \ell(x)\leq \ell(y).
\end{align*}
If $\sim\in \calE[P]$, then $\lambda^\calB_{F,\Omega}(P\mid \sim, \ell)$ is the number of total orders $\leq'_P$ containing $\sim\cap\leq_P$ such that
\begin{align*}
&\forall x,y\in X,&x\leq'_P y&\Longrightarrow \ell(x)\leq \ell(y).
\end{align*}\end{example}

\section{Particular sub-objects}

\subsection{A lemma on twisted subcoalgebras}

\begin{lemma}\label{lem5.3.1}
Let $\calB$ be a CTB, and let $\calB'$ be a subspecies of $\calB$. 
There exists a twisted subcoalgebra $\co\calB'$ of $\calB$, included in $\calB'$, and maximal for the inclusion. If $\calB'$ is a twisted subalgebra of $\calB$, then $\co\calB'$ is a sub-CTB of $\calB$.
\end{lemma}

\begin{proof}
Let us put 
\[\bfC=\{\calC\subseteq \calB'\mid \calC\mbox{ is a twisted subcoalgebra of $\calB$}\}.\]
Then $\bfC$ is not empty: it contains $\calB'$ defined by $\calB'[X]=\emptyset$ for any finite set $X$.
Let $(\calC_i)_{i\in I}$ be a family of elements of $\bfC$. Then 
\[\bigcup_{i\in I} \calC_i\in \bfC. \]
Therefore, we can consider 
\[\co\calB'=\sum_{\calC\in C}\calC.\]
This is an element of $\bfC$, obviously maximal for the inclusion.

Let us now assume that $\calB'$ is a twisted subbialgebra of $\calB$. 
Then $1_\calB\in \calB'$. Therefore, the subspecies defined by
\[\calC[X]=\begin{cases}
\{1_\calB\}\mbox{ if }X=\emptyset,\\
\emptyset \mbox{ otherwise}
\end{cases}\]
belong to $\bfC$, so $1_\calB$ belong to $\co\calB'$.
Let $\calC_1$ and $\calC_2$ in $C$. We consider the subspecies $\calC_1\calC_2$ defined as follows: for any finite set $X$,
\[\calC_1\calC_2[X]=\bigcup_{X=X_1\sqcup X_2}m_{X_1,X_2}(\calC_1[X_1]\otimes \calC_2[X_2]).\]
As $m$ and $\Delta$ are compatible, $\calC_1\calC_2\in \bfC$. 
Applying this to $\calC_1=\calC_2=\co\calB'$, we obtain, by maximality of $\co\calB'$, that 
$\co\calB'\co\calB'\subseteq \co\calB'$. Therefore, $\co\calB'$ is a twisted subalgebra of $\calB$.
\end{proof}

\begin{remark}
Here is a way to compute (more or less explicitly) $\co\calB'[X]$, by induction on $X$.
\begin{itemize}
\item $\co\calB'[\emptyset]=\{1_\calB\}$ if $1_\calB\in \calB'[\emptyset]$, and $\emptyset$ otherwise.
\item If $X\neq \emptyset$, then $\co\calB'[X]$ is the set of elements $b\in \calB'[X]$ such that for any $(X',X'')\in \calC_-(b)$, $\Delta_{X',X''}(b)\in \co\calB'[X']\otimes \co\calB'[X'']$.
\end{itemize}
By construction, the species constructed in this way is a twisted subcoalgebra of $\calB$ included in $\calB'$, and is maximal for the inclusion.
\end{remark}

\subsection{Cocommutative sub-objects}

\begin{prop}\label{prop5.3.2}\begin{enumerate}
\item Let $\calB$ be a CTB. It contains a cocommutative twisted subcoalgebra $\cocom\calB$, maximal for the inclusion. Moreover, $\cocom\calB$ is a sub-CTDB of $\calB$.
\item Let $\calB$ be a non-degenerate CTDB. Then $\cocom\calB$ is a sub-CTDB of $\calB$.
\end{enumerate}\end{prop}

\begin{proof}
1. We put, for any finite set $X$,
\[\calB'[X]=\{b\in \calB[X]\mid \forall (X',X'')\in \comp[X],\: \Delta_{X'',X'}(b)=\tau_{\calB[X'],\calB[X'']}\circ \Delta_{X',X''}(b)\}.\]
Applying Lemma \ref{lem5.3.1}, $\cocom\calB=\co\calB'$ is a cocommutative twisted subcoalgebra of $\calB$, and is maximal for the inclusion. Let $X$, $Y$ be disjoint finite sets, and $a\in \calB[X]$, $b\in \calB[Y]$. For $(X',X'')\in \comp[X\sqcup Y]$,
\begin{align*}
\tau_{\calB[X'],\calB[X'']}\circ\Delta_{X',X''}(ab)&=\tau_{\calB[X'],\calB[X'']}\left(\Delta_{X'\cap X,X''\cap X}(a)\Delta_{X'\cap Y,X'\cap Y}(b)\right)\\
&=\Delta_{X''\cap X,X'\cap X}(a)\Delta_{X''\cap Y,X'\cap Y}(b)\\
&=\Delta_{X'',X'}(ab),
\end{align*}
so $\calB'$ is a twisted subalgebra of $\calB$. By Lemma \ref{lem5.3.1}, $\cocom\calB$ is a twisted subbialgebra of $\calB$.

2. For any finite set $X$, we denote by $\calB'[X]$ the set of elements $b\in \calB[X]$, such that there exists a finite set $X'$, $b'\in \cocom\calB[X]$, $\sim\in \calE[b]$, and $\sigma:X'/{\sim}\longrightarrow X$ a bijection, such that $b=\calB[\sigma](b'/{\sim})$.
Taking $b'=b$ and $\sim=\sim_b^{(m)}$, we observe that $\cocom\calB\subseteq \calB'$. Let $b\in \calB'[X\sqcup Y]$, such that $(X$, $Y)\in \calC_-(b)$. Let $Z$ be a finite set, $b'\in \cocom\calB[Z]$, $\sim\in \calE[b]$, and $\sigma:Z/{\sim}\longrightarrow X\sqcup Y$ a bijection, such that $b=b'/{\sim}$.
Up to this bijection, we can assume that $X\sqcup Y=X'/{\sim}$ and $\sigma=\id_{X\sqcup Y}$. We put $X'=p_\sim^{-1}(X)$, $Y'=p_\sim^{-1}(Y)$, $\sim_X=\sim\cap X'^2$, $\sim_Y=\sim\cap Y'^2$, such that $Z=X'\sqcup Y'$ and $X=X'/{\sim_X}$, $Y=Y'/{\sim_Y}$. Then, as $b'\in \cocom\calB[X\sqcup Y]$,
\begin{align*}
&(\tau_{\calB[Y],\calB[X']}\otimes \id_{\calB[X]})\circ (\Delta_{Y,X}\otimes \id_{\calB[Z]})\circ \delta_\sim(b')\\
&=(\tau_{\calB[Y],\calB[X]}\otimes \id_{\calB[Z]})\circ m_{1,3,24}\circ (\delta_{\sim'}\otimes \delta_{\sim''})\circ\Delta_{Y',X'}(b')\\
&=(\tau_{\calB[Y],\calB[X]}\otimes \id_{\calB[Z]})\circ m_{1,3,24}\circ (\delta_{\sim'}\otimes \delta_{\sim''})\circ\tau_{\calB[X'],\calB[Y'']}\circ \Delta_{X',Y'}(b')\\
&=(\tau_{\calB[Y],\calB[X]}\otimes \id_{\calB[Z]})\circ(\tau_{\calB[X],\calB[Y]}\otimes \id_{\calB[Z]})\circ m_{1,3,24}\circ (\delta_{\sim'}\otimes \delta_{\sim''})\circ \Delta_{X',Y'}(b')\\
&=(\Delta_{X,Y}\otimes \id_{\calB[Z]})\circ \delta_\sim(b').
\end{align*}
Identifying the two first components,
\[\tau_{\calB[Y],\calB[X']}\circ \Delta_{Y,X}(b)=\Delta_{X,Y}(b),\]
so $\calB'$ is a cocommutative twisted subcoalgebra of $\calB$, which contains $\cocom\calB$. By maximality of $\cocom\calB$, $\calB'=\cocom\calB$. Therefore, if $X$ is a finite set, $b\in \cocom\calB[X]$, and $\sim\in \calE[b]$, then $b/{\sim}\in \calB'[X/{\sim}]=\cocom\calB[X/{\sim}]$.
By Lemma \ref{lem4.2.3}, $b\mid \sim\in \cocom\calB[X]$, so $\cocom\calB$ is stable under $\delta$.
\end{proof}

\begin{example}\begin{enumerate}
\item $\cocom\calP=\bfD\calP$: if $X$ is a finite set and $P\in \cocom\calP[X]$, then for any $(X',X'')\in \comp[X]$, 
\[\tau_{\calP[X''],\calP[X']}\circ \Delta_{X'',X'}(P)=\Delta_{X',X''}(P).\]
As a consequence, any lower part of $P$ is also the complement of a lower part of $P$. This implies that $P\in \bfD\calP[X]$.
Conversely, $\bfD\calP$ is a cocommutative twisted subcoalgebra of $\calP$, so is included in $\cocom\calP$.
\item Similarly, $\cocom\calP\calQ$ is the subspecies of equivalences, and $\cocom\calO\calG$ is the subspecies of graphs such that any vertex is part of a cycle. 
\item Consider $\calB=\cocom(\calH\calG,m,\Delta^{\subset,\cap})$. If $G\in \calB[X]$, then for any $Y\subseteq X$, $G_{\mid Y}=G_{\mid_\subset Y}$. Therefore, $\calB$ is the subspecies of hypergraphs $G$ such that the connected components $G'$ of $G$ are of one of the following form:
\begin{itemize}
\item $V(G')=\{x\}$ and $\HE(G')=\emptyset$ or $\HE(G')=\{\{x\}\}$.
\item For any $x\in V(G')$, $\{x\}\in \HE(G')$ and the non-trivial hyperedges of $G$ are of cardinality 2.
\end{itemize}
So $\cocom(\calH\calG,m,\Delta^{\subset,\cap})$ is isomorphic to the CTB of graphs, with the addition of another graph of degree 1, with no edge.
\end{enumerate}\end{example}

\begin{lemma}\label{lemsubcocom}
Let $\calC$ be a CTB and $\calB$ be a sub-CTB of $\calC$. Then $\cocom\calB=\calB\cap\cocom\calC$.
\end{lemma}

\begin{proof}
$\subseteq$. The sub-CTB $\cocom\calB$ of $\calC$ is cocommutative. By maximality of $\cocom\calC$, $\cocom\calB\subseteq \calB\cap \cocom\calC$.

$\supseteq$. By intersection, $\calB\cap \cocom\calC$ is a sub-CTB of $\calB$, and it is cocommutative, as $\cocom\calC$.
By maximality of $\cocom\calB$, $\calB\cap\cocom\calC\subseteq \cocom\calB$. 
\end{proof}

\subsection{Description of $\cocom\exch$}

Let us consider the case of $\exch$ and its sub-CTBs $\exch_\pm$ and $\exch_{\pp}$.
We shall need the following definition:

\begin{defi}\label{defisigmainvariance}
Let $X$ be a finite set, and $I\subseteq \comp[X]$. We shall say that $I$ is $\sym$-invariant if
\begin{align*}
&\forall (X_1,\ldots,X_k)\in I,\: \forall \sigma\in \sym_k,\: (X_{\sigma(1)},\ldots, X_{\sigma(k)})\in I.
\end{align*}
The subspecies of elements $I\in \exch$ which are $\sym$-invariant is denoted by $\sym\exch$.
\end{defi}

\begin{prop}\label{propsymivariance}\begin{align*}
\cocom\exch&=\sym\exch,&\cocom\exch_\pm&=\sym\exch\cap \exch_\pm,\\
&&\cocom\exch_{\pp}&=\sym\exch\cap \exch_{\pp}.
\end{align*}\end{prop}

\begin{proof}
$\cocom\exch\subseteq\sym\exch$. Let $X$ be a finite set, and $I\in \cocom\exch[X]$.
Let $(X_1,\ldots,X_k)\in I$. Then $(X_1,\ldots,X_k)\in m(I)$ and $(X_i)\in I_{\mid_i X_i}$ for any $i\in [k]$. 
As $\cocom\exch$ is cocommutative, for any $\sigma\in \sym_k$, $I_{\mid_i X_{\sigma(i)}}=I_{\mid_{\sigma(i)}X_{\sigma(i)}}$ for any $i\in [k]$ , so $(X_{\sigma(i)}) \in I_{\mid_i X_{\sigma(i)}}$. 
By the exchange condition, $(X_{\sigma(1)},\ldots,X_{\sigma(k)})\in I$.

$\sym\exch\subseteq \cocom\exch$. Let $(X,Y)$ be a pair of finite disjoint sets, and $I\in \sym\exch[X]$, such that $(X,Y)\in m(I)$. 
Let $(X_1,\ldots,X_k)\in I_{\mid_L X}$, and $\sigma\in \sym_k$. By definition, there exists $(Y_1,\ldots,Y_l)\in \comp[Y]$, such that $(X_1,\ldots,X_k,Y_1,\ldots,Y_l)\in I$. As $I$ is $\sym$-invariant, the composition $(X_{\sigma(1)},\ldots,X_{\sigma(k)},Y_1,\ldots,Y_l)$ belongs to $I$,
so $(X_{\sigma(1)},\ldots,X_{\sigma(k)})\in I_{\mid_LX}$. We proved that $I_{\mid_LX}\in \sym\exch[X]$. Similarly, $I_{\mid_RY}\in \sym\exch[Y]$. We obtain that $\sym\exch$ is a twisted subcoalgebra of $\exch$. 
Moreover, still by $\sym$-invariance of $I$, $(Y_1,\ldots,Y_l,X_1,\ldots,X_k)\in I$, so $(X_1,\ldots,X_k)\in I_{\mid_RX}$.
We proved that $I_{\mid_LX}\subseteq I_{\mid_RX}$. By symmetry, $I_{\mid_LX}=I_{\mid_RX}$, and, similarly, 
$I_{\mid_RY}=I_{\mid_LY}$, which gives that
\[\Delta_{Y,X}(I)=\tau_{\exch[X],\exch[Y]}\circ\Delta_{Y,X}(I).\]
We obtain that $\sym\exch$ is cocommutative. By maximality, $\sym\exch\subseteq \cocom\exch$.

The other equalities comes from Lemma \ref{lemsubcocom}. 
\end{proof}

This leads to the following description of $\cocom\exch=\sym\exch$. 

For any finite set $X$, we denote by $\calS\calE[X]$ the set of nonempty subsets of $\calE[X]$. This defines a species $\calS\calE$. For any finite set $X$ and any $I\in \\calS\calE[X]$, we put
\[\iota(I)=\{(X_1,\ldots,X_k)\in \comp[X]\mid \exists \sim\in I,\: X/{\sim}=\{X_1,\ldots,X_k\}\}.\]
Then $\iota$ defines a bijection from $\calS\calE[X]$ to the set of nonempty subsets of $\comp[X]$ which are $\sym$-invariant. Let us precise when $\iota$ gives elements of $\sym\exch$.

\begin{defi}
Let $X$ be a finite set, and $I\in \calS\calE[X]$. We shall say that $I$ satisfies the exchange condition if for any $(X_1,X_2)\in \comp[X]$,
\begin{align*}
&\forall \sim_1,\sim'_1\in \calE[X_1],\: \forall \sim_2,\sim'_2\in \calE[X_2],& \sim_1\sqcup \sim_2\in I\mbox{ and }\sim'_1\sqcup \sim'_2 \in I\Longrightarrow \sim'_1\sqcup \sim_2\in I.
\end{align*}
This defines a subspecies $\calE\exch$ of $\calS\calE$.
\end{defi}

By construction, $\iota$ induces an isomorphism of species between $\calE\exch$ and $\sym\exch=\cocom\exch$.
Therefore, $\calE\exch$ can be made a CTB, isomorphic to $\cocom\exch$. Let us describe this structure.
\begin{itemize}
\item Let $(X,Y)$ be a pair of disjoint finite sets, and $I\in \calE\exch[X\sqcup Y]$. We put
\begin{align*}
I_{\mid X}&=\{\sim\in \calE[X]\mid \exists \sim'\in \calE[Y],\: \sim\sqcup \sim'\in I\},\\
I_{\mid Y}&=\{\sim\in \calE[Y]\mid \exists \sim'\in \calE[X],\: \sim\sqcup \sim'\in I\}.
\end{align*}
Then
\[\Delta_{X,Y}(I)=\begin{cases}
I_{\mid X}\otimes I_{\mid Y}\mbox{ if }I_{\mid X}\neq \emptyset \mbox{ and }I_{\mid_Y}\neq \emptyset,\\
0\mbox{ otherwise}.
\end{cases}\]
\item Let $(X,Y)$ be a pair of disjoint finite sets, $I\in \calE\exch[X]$, and $J\in \calE\exch[Y]$. Then
\[I\squplus J=\{\sim \in \calE[X\sqcup Y]\mid \sim\cap X^2\in I,\: \sim\cap Y^2\in J\}.\]
\end{itemize}

Let us now describe $\cocom\exch_\pm$. Recall that $\calE$ is partially ordered by the inclusion. 
For $X$ a finite set and $I\subseteq \calE[X]$, we shall say that $I$ is an upper part of $\calE[X]$ if
\begin{align*}
&\forall \sim,\sim'\in \calE[X],\: \sim\subseteq \sim'\mbox{ and }\sim'\in I&\Longrightarrow \sim\in I.
\end{align*}
The set of elements $I\in \calE\exch[X]$ which are upper parts of $\calE[X]$ is denoted by $\calE\exch_+[X]$. This defines a species $\calE\exch_+$, isomorphic to $\sym\exch_+$ by construction.
Therefore, $\calE\exch_+$ is a sub-CTB of $\calE\exch$, isomorphic to $\cocom\exch_+$.

For $X$ a finite set and $I\subseteq \calE[X]$, we shall say that $I$ is an lower part of $\calE[X]$ if
\begin{align*}
&\forall \sim,\sim'\in \calE[X],\: \sim\subseteq \sim'\mbox{ and }\sim\in I&\Longrightarrow \sim'\in I.
\end{align*}
The set of elements $I\in \calE\exch[X]$ which are lower parts of $\calE[X]$ is denoted by $\calE\exch_-[X]$. This defines a species $\calE\exch_-$, isomorphic to $\sym\exch_-$ by construction.
Therefore, $\calE\exch_-$ is a sub-CTB of $\calE\exch$, isomorphic to $\cocom\exch_-$.

\begin{remark}
We shall give another description of $\sym\exch_+$ with the help of hypergraphs in Corollary \ref{cor7.2.20}.
\end{remark}

It remains the case of $\cocom\exch_{\pp}$. The $\sym$-invariance shows that if $X=\{x_1,\ldots,x_k\}$ is a finite set, then $\cocom\exch_{\pp}[X]$ is reduced to a single element, namely
\[\{(x_{\sigma(1)},\ldots,x_{\sigma(k)})\mid \sigma\in \sym_k\}.\]
Therefore, $\cocom\exch_{\pp}$ is isomorphic to the CTB $\Com$.

\subsection{Reciprocity principle} 

\begin{defi}\label{defi5.3.3}
Let $\calB$ be a CTB. We shall say that $\calB$ satisfies the reciprocity principle if for any finite set $X$, 
\begin{align*}
&\forall b\in \calB[X],&\phi^\calB_+(b)(-T)&=(-1)^{|X|}\phi^\calB_-(b)(T).
\end{align*}\end{defi}

\begin{prop}\label{prop5.3.4}
Let $\calB$ be a CTB. The following conditions are equivalent:
\begin{enumerate}
\item $\calB$ satisfies the reciprocity principle.
\item For any finite set $X$, and any $b\in \calB[X]$, $\phi^\calB_+(b)(-1)=(-1)^{|X|}$.
\end{enumerate}\end{prop}

\begin{proof}
$1.\Longrightarrow 2$. Let $X$ be a finite set, and $b\in \calB[X]$. 
\[\phi^\calB_+(b)(-1)=(-1)^{|X|}\phi^\calB_-(b)(1)=(-1)^{|X|} \epsilon^\calB_-(b)=(-1)^{|X|}.\]

$2.\Longrightarrow 1$. Let us put, for any finite set $X$, and any $b\in \calB[X]$,
\[\psi(b)=(-1)^{|X|} \phi^\calB_+(b)(-T)\in \K[T].\]
Let $X$, $Y$ be disjoint finite sets. For $a\in \calB[X]$, and $b\in \calB[Y]$,
\begin{align*}
\psi(ab)&=(-1)^{|X|+|Y|}\phi^\calB_+(ab)(-T)\\
&=(-1)^{|X|}(-1)^{|Y|}\phi^\calB_+(a)(-T)\phi^\calB_+(b)(-T)\\
&=\psi(a)\psi(b).
\end{align*}
For $b\in \calB[X]$,
\begin{align*}
\psi(b)(T+T')&=(-1)^{|X|}\phi^\calB_+(b)(-T-T')\\
&=(-1)^{|X|}\sum_{X=X_1\sqcup X_2}(\phi^\calB_+\otimes \phi^\calB_+)\circ \Delta_{X_1,X_2}(b)(-T,-T')\\
&=\sum_{X=X_1\sqcup X_2}(-1)^{|X_1|}(-1)^{|X_2|}(\phi^\calB_+\otimes \phi^\calB_+)\circ \Delta_{X_1,X_2}(b)(-T,-T')\\
&=\sum_{X=X_1\sqcup X_2}(\psi\otimes \psi)\circ \Delta_{X_1,X_2}(b)(T,T').
\end{align*}
So $\psi$ is a polynomial invariant. Moreover, by hypothesis on $\phi^\calB_+$,
\[\psi(b)(1)=(-1)^{|X|}\phi^\calB_+(b)(-1)=1=\epsilon^\calB_-(b).\]
By unicity, $\psi=\phi^\calB_-$. \end{proof}

\begin{cor}\label{cor5.3.5}
The CTB of posets satisfies the reciprocity principle.
\end{cor}

\begin{proof}
Let $X$ be a finite set and $P\in \calP[X]$. Let us prove that $\phi_+^\calP(P)(-1)=(-1)^{|X|}$, by induction on $|X|$. It is obvious if $X=\emptyset$, as then $P=1_\calP$ and $\phi_+^\calP(P)=1$.
Let us assume the result at all ranks $<|X|$, with $|X|\geq 1$. We denote by $\min(P)$ the set of minimal elements of the poset $P$. The combinatorial interpretation of $\phi^\calP_+(P)$ gives that, considering the elements of $X$ sent to $1$ by a valid map $\sigma$,
\[\phi^\calP_+(P)(T)=\sum_{Y\subseteq \min(P)}\phi^\calP_+(P_{\mid X\setminus Y})(T-1)
=\phi^\calP_+(T-1)+\sum_{\emptyset\subsetneq Y\subseteq \min(P)}\phi^\calP_+(P_{\mid X\setminus Y})(T-1).\]
Evaluating in $T=0$, the induction hypothesis gives
\begin{align*}
0&=\phi^\calP_+(P)(0)\\
&=\phi^\calP_+(-1)+\sum_{\emptyset\subsetneq Y\subseteq \min(P)}\phi^\calP_+(P_{\mid X\setminus Y})(-1)\\
&=\phi^\calP_+(-1)+\sum_{\emptyset\subsetneq Y\subseteq \min(P)}(-1)^{|X|-|Y|}\\
&=\phi^\calP_+(-1)-(-1)^{|X|}+\sum_{Y\subseteq \min(P)}(-1)^{|X|-|Y|}\\
&=\phi^\calP_+(-1)-(-1)^{|X|}+(-1)^{|X|-|\min(P)|}(1-1)^{|\min(P)|}\\
&=\phi^\calP_+(-1)-(-1)^{|X|},
\end{align*}
as $\min(P)\neq \emptyset$. Therefore, $\phi^\calP_+(P)(-1)=(-1)^{|X|}$. \end{proof}

\begin{prop}\label{prop5.3.6}\begin{enumerate}
\item Let $\calB$ be a CTB. There exists a subs-CTB $\calB^{\mathrm{recip}}$ of $\calB$, maximal for the inclusion, such that for any finite set $X$, 
\begin{align*}
&\forall b\in \calB^{\mathrm{recip}}[X],&\phi^\calB_+(b)(-T)&=(-1)^{|X|}\phi^\calB_-(b).
\end{align*}
Moreover, $\bfD\calB[X]\subseteq \calB^{\mathrm{recip}}$.
\item Let $\calB,\calB'$ be two CTBs, and $\psi:\calB\longrightarrow \calB'$ be a morphism of CTBs, such that $\epsilon^{\calB'}_+\circ \psi=\epsilon^\calB_+$. Then $\psi(\calB^{\mathrm{recip}})\subseteq \calB'^{\mathrm{recip}}$.
\item $\Gamma^\calB_+(\calB^{\mathrm{recip}})\subseteq \exch_+^{\mathrm{recip}}$.
\end{enumerate}\end{prop}

\begin{proof}
1. We consider the subspecies $\calB'$ defined as follows: for any finite set $X$,
\[\calB'[X]=\left\{b\in \calB[X]\mid \phi^\calB_+(b)(-1)=(-1)^{|X|}\right\}.\]
As $\phi^\calB_+$ is compatible with the product, $\calB'$ is a twisted subalgebra of $\calB$.
We put $\calB^{\mathrm{recip}}=\co\calB'$. Then, by Lemma \ref{lem5.3.1}, it is a twisted subbialgebra of $\calB$. Moreover, for any finite set $X$, for any $b\in \calB^{\mathrm{recip}}[X]$, 
\[\phi^\calB_+(b)(-1)=(-1)^{|X|}.\]
By Proposition \ref{prop5.3.4}, $\calB^{\mathrm{recip}}$ satisfied the reciprocity principle.

Let $\calB''$ be a sub-CTDB of $\calB$, satisfying the reciprocity principle.
By Proposition \ref{prop5.3.4}, $\calB''$ is included in $\calB'$. By maximality of $\calB^{\mathrm{recip}}$, $\calB''\subseteq \calB^{\mathrm{recip}}$.
For any finite set $X$, and any $b\in \bfD\calB[X]$, $\phi^\calB_+(b)=T^{|X|}$, so $\phi^\calB_+(b)(-1)=(-1)^{|X|}$. Therefore, by maximality, $\bfD\calB\subseteq \calB^{\mathrm{recip}}$.

2. Let us consider $\calB''=\psi(\calB^{\mathrm{recip}})$ It is a sub-CTDB of $\calB'$. By Proposition \ref{prop3.5.11}, $\phi^{\calB'}_+\circ \psi=\phi^\calB_+$.
If $X$ is a finite set, and $b\in \calB^{\mathrm{recip}}[X]$,
\[\phi^{\calB'}_+(\psi(b))(-1)=\phi^\calB_+(b)(-1)=(-1)^{|X|}.\]
By Proposition \ref{prop5.3.4}, $\calB''$ satisfies the reciprocity principle. By maximality, $\calB''\subseteq \calB'^{\mathrm{recip}}$.

3. By Proposition \ref{prop3.5.13} and \ref{prop3.5.11}, $\epsilon_+\circ \Gamma^\calB_+=\epsilon^\calB_+$. We apply the second item and obtain the result. 
\end{proof}

\begin{remark}
As a consequence, for any CTB with the reciprocity principle, $\Gamma^\calB_+(\calB)\subseteq \exch_+^{\mathrm{recip}}$.
This is the case with $\calB=\calP$, the CTB of posets. It is not difficult to show that if $|X|\leq 2$, then $\exch_+^{\mathrm{recip}}[X]=\Gamma^\calP_+(\calP[X])$.
However, if $|X|\geq 3$, this is a strict inclusion, as shown by the next examples.
\end{remark}

\begin{example}\label{ex5.3.2}
Let un consider $I\in \exch_+[3]$, defined by
\[I=\left\{(1,2,3),(1,3,2),(2,1,3),(2,3,1),(12,3),(1,23),(2,13)\right\}.\]
Consequently, 
\[\Phi_+(I)=4\frac{T(T-1)(T-2)}{6}+3\dfrac{T(T-1)}{2}=\dfrac{(4T+1)(T-1)T}{6},\]
so $\Phi_+(I)(-1)=-1$.
Then 
\[m(I)=I\sqcup\{(13,2),(23,1),(123)\},\]
which gives
\[\Phi_-(I)=4\frac{T(T-1)(T-2)}{6}+5\dfrac{T(T-1)}{2}+T=\dfrac{(4T-1)(T+1)T}{6}.\]
Moreover,
\begin{align*}
\Delta_{\{1\},\{2,3\}}(I)&=\{(1)\}\otimes \{(2,3),(3,2),(23)\},\\
\Delta_{\{2\},\{1,3\}}(I)&=\{(2)\}\otimes \{(1,3),(3,1),(13)\},\\
\Delta_{\{1,2\},\{3\}}(I)&=\{(1,2),(2,1),(12)\}\otimes \{(3)\},\\
\Delta_{\{1,3\},\{2\}}(I)&=\{(1,3)\}\otimes \{(2)\},\\
\Delta_{\{2,3\},\{1\}}(I)&=\{(2,3)\}\otimes \{(1)\}.
\end{align*}
All the elements on the right correspond to posets, so belong to $\exch_+^{\mathrm{recip}}$. 
Therefore, $I\in \exch_+^{\mathrm{recip}}[3]$. However, it does not correspond to any poset, 
as $\Phi_+(I)$ is not a strict Erhahrt polynomial, see Remark \ref{rk3.5.3}.
\end{example}

\begin{example}\label{ex5.3.3}
Another example of element of $\exch_+^{\mathrm{recip}}$, which does not correspond to a poset, is given by
\[J=\left\{(1,2,3),(1,3,2),(2,1,3),(2,3,1),(3,1,2),(12,3),(13,2),(1,23),(2,13)\right\}.\]
We obtain
\[m(J)=J\sqcup\{(3,12),(23,1),(123)\},\]
so
\begin{align*}
\Phi_+(J)&=5\frac{T(T-1)(T-2)}{6}+4\dfrac{T(T-1)}{2}=\dfrac{(5T+2)(T-1)T}{6},\\
\Phi_-(J)&=5\frac{T(T-1)(T-2)}{6}+6\dfrac{T(T-1)}{2}+T=\dfrac{(5T-2)(T+1)T}{6}.
\end{align*}
\end{example}

\begin{remark}
it is possible to prove that any $K\in \exch_+^{\mathrm{recip}}[X]$, with $|X|=3$, is, either the image of a poset, or is isomorphic to $I$ or $J$.
\end{remark}

\begin{prop}
Let $\calB$ be a non-degenerate cocommutative CTB.
Then $\calB^{\mathrm{recip}}=\bfD\calB$.
\end{prop}

\begin{proof}
We already observed that $\bfD\calB\subseteq \calB^{\mathrm{recip}}$. Let $X$ be a finite set, and $b\in \calB^{\mathrm{recip}}[X]$. 
As $\calB$ is non-degenerate, $\phi^\calB_-(b)=T^{|X|}$. By the reciprocity principle on $\calB^{\mathrm{recip}}$,
\[\phi^\calB_+(b)=(-1)^{|X|}\phi^\calB_-(b)(-T)=T^{|X|}.\]
By proposition \ref{prop3.5.10}, $b\in \bfD\calB[X]$.
\end{proof}

\subsection{Alternating signs}

\begin{defi}\label{defi5.3.8}
Let $\calB$ be a CTB. We shall say that it satisfies the alternating signs condition if for any finite set $X$, and any $b\in \calB[X]$, putting
\[\phi^\calB_+(b)=a_0T+\cdots+a_nT^n,\]
then, for $i\in \llbracket\mathrm{val}(\phi^\calB_+(b));\:\deg(\phi^\calB_+(b))\rrbracket$, $a_i\neq 0$ and its sign is $(-1)^{|X|+i}$.
\end{defi}

\begin{example}\begin{enumerate}
\item The CTB $\calG$ of graphs satisfies the alternating signs condition \cite{Rota1964}.
\item The twisted subalgebra of the CTB $\calP$ of posets generated by total orders is a sub-CTDB of $\calP$,
and it satisfies the alternating signs condition: indeed, if $\leq_P$ is a total order on a finite set $X$, denoting by $n$ the cardinality of $X$,
\[\phi^\calP_+(P)=\dfrac{T(T-1)\ldots (T-k+1)}{k!}.\]
As the product of two polynomials with the alternating signs condition also satisfies the alternating signs condition, for any product $P$ of total orders, $\phi^\calP_+(P)$ satisfies the alternating signs condition.
\item The CTB $\calH\calG$ of hypergraphs does not:
if $G\in \calH\calG[4]$ is given by
\[\HE(G)=\{\{1,2,3\},\:\{1,2,4\}\},\]
then $\phi^{\calH\calG}_+(G)=T^4-2T^2-T$. 
The CTBs $\calP$ of posets and $\calO\calG$ of oriented graphs do not neither: if $P\in \calP[5]$ and $G\in \calO\calG[5]$ are given by
\begin{align*}
<_P=E(G)&=\{(1,2),\:(1,3),\:(1,4),\:(1,5)\},
\end{align*}
then $\phi^\calP_+(P)=\phi^{\calO\calG}_+(G)=\dfrac{6T^5-15T^4+10T^3-T}{30}$.
\end{enumerate}\end{example}

\begin{prop}\label{prop5.3.9}\begin{enumerate}
\item Let $\calB$ be a CTB. There exists a sub-CTB $\calB^{\mathrm{alt}}$ of $\calB$, satisfying the alternating signs condition, and maximal for the inclusion. 
\item Let $\calB,\calB'$ be two CTBs, and $\psi:\calB\longrightarrow \calB'$ be a morphism of CTBs, such that $\epsilon^{\calB'}_+\circ \psi=\epsilon^\calB_+$. Then $\psi(\calB^{\mathrm{alt}})\subseteq \calB'^{\mathrm{alt}}$.
\item $\Gamma^\calB_+(\calB^{\mathrm{alt}})\subseteq \exch_+^{\mathrm{alt}}$.
\end{enumerate}
\end{prop}

\begin{proof}
1. For any finite set $X$, we put
\[\calB'[X]=\{b\in \calB[X]\mid \phi^\calB_+(b)\mbox{ has the alternating signs condition}\}.\]
This defines a subspecies of $\calB$, and it is not difficult to prove that it is a twisted subalgebra,
as the product of two polynomials with the alternating sign condition also has the alternating condition. We then put $\calB^{\mathrm{alt}}=\co\calB'$.

2. Then, by Proposition \ref{prop3.5.11}, $\phi^{\calB'}_+\circ \psi=\phi^\calB_+$. 
Therefore, $\psi(\calB^{\mathrm{alt}})$ is a sub-CTB of $\calB'$, with the alternating signs condition. By maximality, $\psi(\calB^{\mathrm{alt}})\subseteq \calB'^{\mathrm{alt}}$.

3. By Proposition \ref{prop3.5.13} and \ref{prop3.5.11}, $\epsilon_+\circ \Gamma^\calB_+=\epsilon^\calB_+$. We apply the second item and obtain the result. 
\end{proof}

We shall give more results on the alternating sign conditions when $\calB$ is cocommutative in Theorem \ref{thoe7.2.30}.

\chapter{An example based on rooted forests}

\label{chap6}

\section*{Introduction}

We study in this chapter a CTB on rooted trees. The most famous Hopf algebra of rooted trees is perhaps the one due to Connes and Kreimer \cite{Connes1998,Brouder2004,Foissy3}, related to Butcher's group of Runge-Kutta's methods, whose construction can be easily led to the level of CTBs. 
For example, the Butcher-Connes-Kreimer coproduct gives
\begin{align*}
\Delta^{\mathrm{BCK}}\left(\hspace{-2mm}\begin{array}{c}\tquatreun\end{array}\hspace{-2mm}\right)&=1\otimes \tquatreun+\tquatreun\otimes 1+ 3\ttroisun\otimes \tun+3\tdeux\otimes\tun\tun+\tun\otimes \tun\tun\tun,\\
\Delta^{\mathrm{BCK}}\left(\hspace{-2mm}\begin{array}{c}\tquatredeux\end{array}\hspace{-2mm}\right)&=1\otimes\tquatredeux+\tquatredeux\otimes 1+\ttroisun\otimes\tun+\ttroisdeux\otimes\tun+\tdeux\otimes\tdeux+\tdeux\otimes\tun\tun+\tun\otimes\tdeux\tun,\\
\Delta^{\mathrm{BCK}}\left(\hspace{-2mm}\begin{array}{c}\tquatrequatre\end{array}\hspace{-2mm}\right)&=1\otimes\tquatrequatre+\tquatrequatre\otimes 1+2\ttroisdeux\otimes\tun+\tdeux\otimes\tun\tun+\tun\otimes\ttroisun,\\
\Delta^{\mathrm{BCK}}\left(\hspace{-2mm}\begin{array}{c}\tquatrecinq\end{array}\hspace{-2mm}\right)&=1\otimes\tquatrecinq+\tquatrecinq\otimes 1+\ttroisdeux\otimes\tun+\tdeux\otimes\tdeux+\tun\otimes\ttroisdeux.
\end{align*}
This can be seen as a sub-CTB of $\calP$, when one considers rooted forests as Hasse graphs of posets, as in Definition \ref{defi6.1.1}. We consider here a variant,
 defined in \cite{Bruned2015} and studied in \cite{Xiao2026} in the context of the use of Hopf algebras for linguistics \cite{Marcolli2025},  where the coproduct is given by
\begin{align*}
\Delta^L\left(\hspace{-2mm}\begin{array}{c}\tquatreun\end{array}\hspace{-2mm}\right)&=1\otimes \tquatreun+\tquatreun\otimes 1+ 3\ttroisun\otimes \tun+3\tdeux\otimes\tun\tun,\\
\Delta^L\left(\hspace{-2mm}\begin{array}{c}\tquatredeux\end{array}\hspace{-2mm}\right)&=1\otimes\tquatredeux+\tquatredeux\otimes 1+\ttroisdeux\otimes\tun+\tdeux\otimes\tdeux,\\
\Delta^L\left(\hspace{-2mm}\begin{array}{c}\tquatrequatre\end{array}\hspace{-2mm}\right)&=1\otimes\tquatrequatre+\tquatrequatre\otimes 1+2\ttroisdeux\otimes\tun,\\
\Delta^L\left(\hspace{-2mm}\begin{array}{c}\tquatrecinq\end{array}\hspace{-2mm}\right)&=1\otimes\tquatrecinq+\tquatrecinq\otimes 1.
\end{align*}
This new coproduct is obtained from $\Delta^{\mathrm{BCK}}$ by keeping only the terms where the number of leaves is conserved. 
An interesting quotient is given on reduced trees, that is to say trees with no vertex with only one child. This gives, after reduction (that is to say, the contraction between a vertex and its unique child, if this exists),
\begin{align*}
\Delta^{r\ell}\left(\hspace{-2mm}\begin{array}{c}\tquatreun\end{array}\hspace{-2mm}\right)&=1\otimes \tquatreun+\tquatreun\otimes 1+ 3\ttroisun\otimes \tun+3\tun\otimes\tun\tun,\\
\Delta^{r\ell}\left(\hspace{-2mm}\begin{array}{c}\ttroisun\end{array}\hspace{-2mm}\right)&=1\otimes\ttroisun+\ttroisun \otimes 1+\tun\otimes\tun+\tun\otimes\tun,\\
\Delta^{r\ell}\left(\hspace{-2mm}\begin{array}{c}\ttroisun\end{array}\hspace{-2mm}\right)&=1\otimes\ttroisun+\ttroisun\otimes 1+2\tun\otimes\tun,\\
\Delta^{r\ell}\left(\hspace{-2mm}\begin{array}{c}\tun\end{array}\hspace{-2mm}\right)&=1\otimes\tun+\tun\otimes 1.
\end{align*}
Note that this coproduct does not conserve the number of vertices anymore, but still conserves the number of leaves.

We shall see that this coproduct $\Delta^L$ can be seen into two different ways at the level of CTBs. The first one, called $\calR\calF$ is based on vertices, see Theorem \ref{theo6.1.3}. For example,
\begin{align*}
\Delta^L\left(\hspace{-2mm}\begin{array}{c}\tdquatreun{$x$}{$t$}{$z$}{$y$}\end{array}\hspace{-2mm}\right)&=1\otimes \tdquatreun{$x$}{$t$}{$z$}{$y$}+\tdquatreun{$x$}{$t$}{$z$}{$y$}\otimes 1+\tdtroisun{$x$}{$t$}{$z$}\otimes \tdun{$y$}+\tdtroisun{$x$}{$t$}{$y$}\otimes \tdun{$z$}+\tdtroisun{$x$}{$z$}{$y$}\otimes \tdun{$t$}+\tddeux{$x$}{$t$}\otimes\tdun{$y$}\tdun{$z$}+\tddeux{$x$}{$z$}\otimes\tdun{$y$}\tdun{$t$}+\tddeux{$x$}{$y$}\otimes\tdun{$z$}\tdun{$t$},\\
\Delta^L\left(\hspace{-2mm}\begin{array}{c}\tdquatredeux{$x$}{$t$}{$y$}{$z$}\end{array}\hspace{-2mm}\right)&=1\otimes\tdquatredeux{$x$}{$t$}{$y$}{$z$}+\tdquatredeux{$x$}{$t$}{$y$}{$z$}\otimes 1+\tdtroisdeux{$x$}{$y$}{$z$}\otimes\tdun{$t$}+\tddeux{$x$}{$t$}\otimes\tddeux{$y$}{$z$},\\
\Delta^L\left(\hspace{-2mm}\begin{array}{c}\tdquatrequatre{$x$}{$y$}{$t$}{$z$}\end{array}\hspace{-2mm}\right)&=1\otimes\tdquatrequatre{$x$}{$y$}{$t$}{$z$}+\tdquatrequatre{$x$}{$y$}{$t$}{$z$}\otimes 1+\tdtroisdeux{$x$}{$y$}{$z$}\otimes\tdun{$t$}+\tdtroisdeux{$x$}{$y$}{$t$}\otimes\tdun{$z$},\\
\Delta^L\left(\hspace{-2mm}\begin{array}{c}\tdquatrecinq{$x$}{$y$}{$z$}{$t$}\end{array}\hspace{-2mm}\right)&=1\otimes\tdquatrecinq{$x$}{$y$}{$z$}{$t$}+\tdquatrecinq{$x$}{$y$}{$z$}{$t$}\otimes 1.
\end{align*}
Its main interest is that the corresponding CTB on rooted forests is quasi-cocommutative, which allows to give a convenient formula for its antipode, using the reciprocity principle for posets (Corollary \ref{cor6.1.7}). 
Its major inconvenient is that it is not compatible with the reduction, which does not conserve the set of vertices, so does not give a CTB version of the coproduct $\Delta^\ell$. 
The second one, called $\calR\calF_\ell$, is based on leaves (Theorem \ref{theo6.2.2}). For example,
\begin{align*}
\Delta^\ell\left(\hspace{-2mm}\begin{array}{c}\tdquatreun{}{$z$}{$y$}{$x$}\end{array}\hspace{-2mm}\right)&=1\otimes \tdquatreun{}{$z$}{$y$}{$x$}+\tdquatreun{}{$z$}{$y$}{$x$}\otimes 1+\tdtroisun{}{$z$}{$y$}\otimes \tdun{$x$}
+\tdtroisun{}{$z$}{$x$}\otimes \tdun{$y$}+\tdtroisun{}{$y$}{$x$}\otimes \tdun{$z$}+\tddeux{}{$z$}\otimes\tdun{$x$}\tdun{$y$}+\tddeux{}{$y$}\otimes\tdun{$x$}\tdun{$z$}+\tddeux{}{$x$}\otimes\tdun{$y$}\tdun{$z$},\\
\Delta^\ell\left(\hspace{-2mm}\begin{array}{c}\tdquatredeux{}{$y$}{}{$x$}\end{array}\hspace{-2mm}\right)&=1\otimes\tdquatredeux{}{$y$}{}{$x$}+\tdquatredeux{}{$y$}{}{$x$}\otimes 1+\tdtroisdeux{}{}{$x$}\otimes\tdun{$y$}+\tddeux{}{$y$}\otimes\tddeux{}{$x$},\\
\Delta^\ell\left(\hspace{-2mm}\begin{array}{c}\tdquatrequatre{}{}{$y$}{$x$}\end{array}\hspace{-2mm}\right)&=1\otimes\tdquatrequatre{}{}{$y$}{$x$}+\tdquatrequatre{}{}{$y$}{$x$}\otimes 1+\tdtroisdeux{}{}{$x$}\otimes\tdun{$y$}+\tdtroisdeux{}{}{$y$}\otimes\tdun{$x$},\\
\Delta^\ell\left(\hspace{-2mm}\begin{array}{c}\tdquatrecinq{}{}{}{$x$}\end{array}\hspace{-2mm}\right)&=1\otimes\tdquatrecinq{}{}{}{$x$}+\tdquatrecinq{}{}{}{$x$}\otimes 1.
\end{align*}
It turns out that the combinatorics of leaves induces a morphism $\Xi$ of CTBs from $\calR\calF_\ell$ to the CTB of hypergraphs $\calH\calG$, with the coproduct $\Delta^{\cap,\subset}$ (Theorem \ref{theo6.2.3}). 
This morphism $\Xi$ is not injective, and its kernel defines a structure of CTB on the species $\calR\calR\calF$ of reduced forests,
which is consequently isomorphic to a sub-object of $\calH\calG$ (Theorem \ref{theo6.2.8}). For example, on reduced rooted forests,
\begin{align*}
\Delta^{r\ell}\left(\hspace{-2mm}\begin{array}{c}\tdquatreun{}{$z$}{$y$}{$x$}\end{array}\hspace{-2mm}\right)&=1\otimes \tdquatreun{}{$z$}{$y$}{$x$}+\tdquatreun{}{$z$}{$y$}{$x$}\otimes 1+\tdtroisun{}{$z$}{$y$}\otimes \tdun{$x$}
+\tdtroisun{}{$z$}{$x$}\otimes \tdun{$y$}+\tdtroisun{}{$x$}{$y$}\otimes \tdun{$z$}+\tdun{$z$}\otimes\tdun{$x$}\tdun{$y$}+\tdun{$y$}\otimes\tdun{$x$}\tdun{$z$}+\tdun{$x$}\otimes\tdun{$y$}\tdun{$z$},\\
\Delta^{r\ell}\left(\hspace{-2mm}\begin{array}{c}\tdtroisun{}{$y$}{$x$}\end{array}\hspace{-2mm}\right)&=1\otimes\tdtroisun{}{$y$}{$x$}+\tdtroisun{}{$y$}{$x$}\otimes 1+\tdun{$x$}\otimes\tdun{$y$}+\tdun{$y$}\otimes\tdun{$x$},\\
\Delta^{r\ell}\left(\hspace{-2mm}\begin{array}{c}\tdtroisun{}{$y$}{$x$}\end{array}\hspace{-2mm}\right)&=1\otimes\tdtroisun{}{$y$}{$x$}+\tdtroisun{}{$y$}{$x$}\otimes 1+\tdun{$x$}\otimes\tdun{$y$}+\tdun{$y$}\otimes\tdun{$x$},\\
\Delta^{r\ell}\left(\hspace{-2mm}\begin{array}{c}\tdun{$x$}\end{array}\hspace{-2mm}\right)&=1\otimes\tdun{$x$}+\tdun{$x$}\otimes 1.
\end{align*}
We show in Theorem \ref{theo6.3.1} that the application of the bosonic Fock functor to both $\calR\calF$ and $\calR\calF_\ell$ gives the same Hopf algebra on rooted trees, with the coproduct $\Delta$ (Theorem \ref{theo6.3.1}), with quotient Xiao's Hopf algebra of reduced rooted forests.
The quasi-cocommutativity of $\calR\calF$ allows to find an easy cancellation-free formula for the antipode of $\calF[\calR\calF]=\calF[\calR\calF_\ell]$.
Grouping the terms according the terms of this formula according to the reduction, 
we obtain a formula for Xiao's Hopf algebra of reduced forests, with a combinatorial interpretation of the coefficients in terms of partitions of the vertices (Proposition \ref{prop6.3.4}), answering a question of \cite{Xiao2026}.

\section{A combinatorial twisted bialgebra on rooted forests based on vertices}

\subsection{Construction}

\begin{defi}\label{defi6.1.1}
For any finite set $X$, we denote by $\calR\calF[X]$ the set of (non planar) rooted forests whose set of vertices is $X$.
For example,
\begin{align*}
\calR\calF[0]&=\{1\},\\
\calR\calF[1]&=\left\{\begin{array}{c}
\tdun{$1$}
\end{array}\right\},\\
\calR\calF[2]&=\left\{\begin{array}{c}
\tdun{$1$}\tdun{$2$},\tddeux{$1$}{$2$},\tddeux{$2$}{$1$}
\end{array}\right\},\\
\calR\calF[3]&=\left\{\begin{array}{c}
\tdun{$1$}\tdun{$2$}\tdun{$3$},\tddeux{$1$}{$2$}\tdun{$3$},\tddeux{$2$}{$1$}\tdun{$3$},\tddeux{$1$}{$3$}\tdun{$2$},\tddeux{$3$}{$1$}\tdun{$2$},\tddeux{$2$}{$3$}\tdun{$1$},\tddeux{$3$}{$2$}\tdun{$1$},\\
\tdtroisun{$1$}{$3$}{$2$},\tdtroisun{$2$}{$3$}{$1$},\tdtroisun{$3$}{$2$}{$1$},
\tdtroisdeux{$1$}{$2$}{$3$},\tdtroisdeux{$1$}{$3$}{$2$},\tdtroisdeux{$2$}{$1$}{$3$},
\tdtroisdeux{$2$}{$3$}{$1$},\tdtroisdeux{$3$}{$1$}{$2$},\tdtroisdeux{$3$}{$2$}{$1$}
\end{array}\right\}.
\end{align*}
This defines a species $\calR\calF$.
If $F\in \calR\calF[X]$, its edges are oriented upwards, from the root to the leaves.
For any $x,y\in X$, we shall say that $x\leq_F y$ if there exists an oriented path in $F$ from $x$ to $y$.
Then $(X,\leq_F)$ is a poset, and $F$ is the Hasse graph of this poset.
In this way, $\calR\calF$ can be seen a subspecies of $\calP$.
\end{defi}

\begin{remark}
For any finite set $X$, of cardinality $n$, $|\calR\calF[X]|=(n+1)^{n-1}$, see \cite[Entry A000272]{Sloane}.
\end{remark}

Note that $\calR\calF$ is in fact a twisted subalgebra of $\calP$. The product of two rooted forests is their disjoint union.

\begin{defi}
Let $X$ be a finite set and $F\in \calR\calF[X]$.
\begin{enumerate}
\item We denote by $L(F)$ the set of leaves of $F$.
\item For $x\in X$, we denote by $L_x(F)$ the set of leaves of $F$ which are descendants of $x$. It is a nonempty subset of $L(F)$.
\item Let $Y\subseteq X$.
\begin{itemize}
\item The restriction $F_{\mid Y}$ is the oriented graph whose vertices are the elements of $Y$ and the edges all the edges of $F$ between these vertices. It belongs to $\calR\calF[Y]$.
\item $Y$ is convex if for any $x,z\in Y$ and $y\in X$,
\[x\leq_Fy\leq_F z\Longrightarrow y\in Y.\]
\item $Y$ is connected if $F_{\mid Y}$ is connected, that is to say is a rooted tree.
\item $Y$ is leaf-compatible if for any $x\in Y$, $L_x(F)\cap Y\neq \emptyset$.
\item $Y$ is an upper part of $F$ if for any $x,y\in X$,
\[(x\in Y\mbox{ and }x\leq_F y\Longrightarrow y\in Y).\]
\item $Y$ is a lower part of $F$ if for any $x,y\in X$,
\[(y\in Y\mbox{ and }x\leq_F y\Longrightarrow x\in Y).\]
\end{itemize}\end{enumerate}
\end{defi}

\begin{remark}
If $F\in \calR\calF[X]$ and $Y\subseteq X$, then:
\begin{enumerate}
\item $Y$ is an upper part of $F$ if, and only if, $X\setminus Y$ is a lower part of $F$.
\item If $Y$ is a lower or an upper part of $F$, then $Y$ is convex.
\end{enumerate}
\end{remark}

\begin{theo}\label{theo6.1.3}
We define a coproduct $\Delta^L$ on $\calR\calF$ as follows. If $(X,Y)$ is a pair of disjoint finite sets, then
\begin{align*}
&\forall F\in \calR\calF[X\sqcup Y],&\Delta^L_{X,Y}(F)&=\begin{cases}
F_{\mid X}\otimes F_{\mid Y}\mbox{ if $X$ is a leaf-compatible lower part of $F$}\\
\hspace{1.8cm}\mbox{and $Y$ is a leaf-compatible upper part of $F$},\\
0\mbox{ otherwise}.
\end{cases}
\end{align*}
Then $(\calR\calF,m,\Delta^L)$ is a quasi-cocommutative CTB.
\end{theo}

\begin{proof}
Let $X,Y,X',Y'$ be finite sets such that $X\sqcup Y=X'\sqcup Y$', $F\in \calR\calF[X]$ and $G\in \calR\calF[Y]$.
Then $X'$ is a lower part of $FG$ (and $Y'$ is an upper part of $FG$) if, and only if, $X'\cap X$ is a lower part of $F$ and $X'\cap Y$ is a lower part of $G$ (and $Y'\cap X$ is an upper part of $F$ and $Y'\cap Y$ is an upper part of $G$).
Moreover, as $L(FG)=L(F)\sqcup L(G)$, $X'$ is leaf-compatible for $FG$ if, and only if, $X'\cap X$ is leaf-compatible for $F$ and $X'\cap Y$ is leaf-compatible for $Y$. Similarly, $Y'$ is leaf-compatible for $FG$ if, and only if, $Y'\cap X$ is leaf-compatible for $F$ and $Y'\cap Y$ is leaf-compatible for $Y$.
Finally, for any $Z\subseteq X\sqcup Y$, $(FG)_{\mid Z}=F_{\mid X\cal Z}G_{\mid Y\cap Z}$.
All this implies that
\[\Delta^L_{X',Y'}(FG)=\Delta^L_{X'\cap X,Y'\cap Y}(F)\Delta^L_{X'\cap Y,Y'\cap Y}(G).\]

Let $X$, $Y$, $Z$ be pairwise disjoint finite sets, and $F\in \calR\calF[X\sqcup Y\sqcup Z]$. We define two conditions:
\begin{description}
\item[C1.] $X\sqcup Y$ is a leaf-compatible lower part of $F$, $Z$ is a leaf-compatible lower part of $F$, $X$ is a leaf-compatible lower part of $F_{\mid X\sqcup Y}$, and $Y$ is a leaf-compatible upper part of $F_{\mid X\sqcup Y}$.
\item[C2.] $X$ is a leaf-compatible lower part of $F$, $Y\sqcup Z$ is a leaf-compatible upper part of $F$, $Y$ is a leaf-compatible lower part of $F_{\mid Y\sqcup Z}$, and $Z$ is a leaf-compatible upper part of $F_{\mid Y\sqcup Z}$.
\end{description} 
Then
\begin{align*}
\left(\Delta^L_{X,Y}\otimes \id_{\calR\calF[Z]}\right)\circ \Delta^L_{X\sqcup Y,Z}(F)&=\begin{cases}
F_{\mid X}\otimes F_{\mid Y}\otimes F_{\mid Z}\mbox{ if C1},\\
0\mbox{ otherwise}.
\end{cases}\\
\left(\id_{\calR\calF[X]}\otimes \Delta^L_{Y,Z}\right)\circ \Delta^L_{X\sqcup Y,Z}(F)&=\begin{cases}
F_{\mid X}\otimes F_{\mid Y}\otimes F_{\mid Z}\mbox{ if C2},\\
0\mbox{ otherwise}.
\end{cases}\end{align*}
Therefore, proving the coassociativity of $\Delta^L$ is equivalent to prove that C1 $\Longleftrightarrow $ C2. We shall prove the following result:
if $(X',Y')$ is a pair of disjoint finite sets, $F\in \calR\calF[X'\sqcup Y]'$ such that $X'$ is a lower part of $F$ and $Y'$ is an upper part of $F$, then $X'$ and $Y'$ are leaf-compatible if, and only if,
\[L(F)=L\left(F_{\mid X'}\right)\sqcup L\left(F_{\mid Y'}\right).\]
Note that we always have $L(F)\subseteq L\left(F_{\mid X'}\right)\sqcup L\left(F_{\mid Y'}\right)$.

$\Longrightarrow$. Let $x$ be a leaf of $F_{\mid X'}$. As $X"$ is leaf-compatible, $L_x(F)\cap X'\neq \emptyset$. Moreover, $L_x(F)\cap X'\subseteq L_x\left(F_{\mid X'}\right)=\{x\}$, so $L_x(F)\cap X'=\{x\}$, and $x\in L(F)$. Similarly, the leaves of $F_{\mid Y'}$ are leaves of $F$.

$\Longleftarrow$. Let $x\in X'$. Then $\emptyset\neq L_x\left(F_{\mid X'}\right)=L_x(F)\cap X'$, so $X$' is leaf-compatible. Similarly, $Y'$ is leaf-compatible.

Consequently, we can rewrite C1 and C2:
\begin{description}
\item[C1.] $X\sqcup Y$ is a lower part of $F$, $Z$ is an upper part of $F$, $X$ is a lower part of $F_{\mid X\sqcup Y}$, and $Y$ is an upper part of $F_{\mid X\sqcup Y}$, and 
$L(F)=L\left(F_{\mid X}\right)\sqcup L\left(F_{\mid Y}\right)\sqcup L\left(F_{\mid Z}\right)$.
\item[C2.] $X$ is a lower part of $F$, $Y\sqcup Z$ is an upper part of $F$, $Y$ is a lower part of $F_{\mid Y\sqcup Z}$, and $Z$ is an upper part of $F_{\mid Y\sqcup Z}$, and 
$L(F)=L\left(F_{\mid X}\right)\sqcup L\left(F_{\mid Y}\right)\sqcup L\left(F_{\mid Z}\right)$.
\end{description} 

C1 $\Longrightarrow$ C2. Let $x\in X$ and $y\in V(F)$, such that $y\leq_F x$. As $X\sqcup Y$ is a lower part of $F$, $y\in X\sqcup Y$. Therefore, $y\leq_{F_{\mid X\sqcup Y}} x$. As $X$ is a lower part of $F_{\mid X\sqcup Y}$, $y\in X$. So $X$ is a lower part of $F$, and consequently $Y\sqcup Z$ is an upper part of $F$. 

Let $y\in Z$ and $x\in Y\sqcup Z$, such that $x\leq_{F_{\mid Y\sqcup Z}}y$. Then $x\leq_F y$.
If $x\in Z$, as $Z$ is an upper part of $F$, $y\in Z$, which contradicts $y\in Y$. So $x\in Y$, which proves
$Y$ is a lower part of $F_{\mid Y\sqcup Z}$, and $Z$ is an upper part of $F_{\mid Y\sqcup Z}$.

C2 $\Longrightarrow$ C1. Let $y\in X\sqcup Y$ and $x\in V(F)$, such that $x\leq_F y$. If $x\in Z$, then $x\leq_{F_{\mid Y\sqcup Z}} y$. As $Z$ is an upper part of $F_{\mid Y\sqcup Z}$, then $y\in Z$, which contradicts $y\in Y$.
So $x\in X\sqcup Y$, which gives that $X\sqcup Y$ is a lower part of $F$ and $Z$ is an upper part of $F$.

Let $x\in X$ and $y\in X\sqcup Y$, such that $y\leq_{F_{\mid X\sqcup Y}} x$. Then $y\leq_F x$. 
As $x$ is a lower part of $F$, $y\in X$, so $X$ is a lower part of $F_{\mid X\sqcup Y}$, and $Y$ is an upper part of $F_{\mid X\sqcup Y}$.

We have finally proved the coassociativity of $\Delta^L$.

Let $X$ be a finite set, and $F\in \calR\calF[X]$. Then
\begin{align*}
\Delta^L_{\emptyset,X}(F)&=1\otimes F,&\Delta^L_{X,\emptyset}(F)&=F\otimes 1.
\end{align*}
So $\Delta^L$ is counitary. Moreover, if $(X_1,\ldots,X_k)$ and $(X_{\sigma(1)},\ldots,X_{\sigma(k)})$ both belong to $\calC_-[F]$, by definition of $\Delta^L$, the for any $i\in [k]$,
\[F_{\mid_i X_i}=F_{\mid_{\sigma^{-1}(i) X_i}}=F_{\mid X_i},\]
so $(\calR\calF,m,\Delta^L)$ is a quasi-cocommutative.
\end{proof}

\begin{example}
We only give the non-zero coproducts.
\begin{align*}
\Delta^L_{\emptyset,[4]}\left(\begin{array}{c}
\tdquatrequatre{$1$}{$2$}{$4$}{$3$}\end{array}\right)&=1\otimes \tdquatrequatre{$1$}{$2$}{$4$}{$3$},&
\Delta^L_{\emptyset,[4]}\left(\begin{array}{c}
\tdquatrequatre{$1$}{$2$}{$4$}{$3$}\end{array}\right)&=\tdquatrequatre{$1$}{$2$}{$4$}{$3$}\otimes 1,\\
\Delta^L_{[3],\{4\}}\left(\begin{array}{c}
\tdquatrequatre{$1$}{$2$}{$4$}{$3$}\end{array}\right)&=\tdtroisdeux{$1$}{$2$}{$3$}\otimes \tdun{$4$},&
\Delta^L_{\{1,2,4\},\{3\}}\left(\begin{array}{c}
\tdquatrequatre{$1$}{$2$}{$4$}{$3$}\end{array}\right)&=\tdtroisdeux{$1$}{$2$}{$4$}\otimes \tdun{$3$}.
\end{align*}\end{example}

\subsection{Relations on rooted forests}
 
Let us describe the set of equivalences associated to any element of the quasi-cocommutative CTB $\calR\calF$, following Definition \ref{defi4.2.5}. We firstly describe $\calC^L_-(F)$ for any rooted forest, where $\calC_-^L(F)$ is the set of compositions $(X_1,\ldots, X_k)\in \comp[X]$ such that $\Delta_{X_1,\ldots,X_k}(F)\neq 0$.

\begin{lemma}\label{lem6.1.4}
Let $X$ be a finite set, and $F\in \calR\calF[X]$. Let $(X_1,\ldots,X_k)\in \comp[X]$.
Then $(X_1,\ldots,X_k)\in \calC^L_-(F)$ if, and only if, the two following conditions hold:
\begin{enumerate}
\item For any $i\in [k]$, $X_i$ is leaf-compatible in $F$.
\item For any $x,y\in X$, $(x\in X_i,\: y\in X_j,\: x\leq_F y)\Longrightarrow i\leq j$. 
\end{enumerate}\end{lemma}

\begin{proof}
$\Longrightarrow$. We proceed by induction on $k$. If $k=1$, then these two conditions are obviously satisfied. Let us assume the result at rank $k-1$. Then
\[\Delta_{X_1,\ldots,X_k}(F)=\left(\Delta_{X_1,\ldots,X_{k-1}}\otimes \id_{\calR\calF[X_k]}\right)\circ \Delta_{X_1\sqcup \cdots \sqcup X_{k-1},X_k}(F)\neq 0,\]
so $X_1\sqcup \cdots \sqcup X_{k-1}$ is a lower part of $F$, $X_k$ is an upper part of $F$, and both are leaf-compatible in $F$. Moreover, $(X_1,\ldots,X_{k-1})\in \calC_-^L\left(F_{\mid X_1\sqcup \cdots \sqcup X_{k-1}}\right)$.
By the induction hypothesis, $X_i$ is leaf-compatible in $F_{\mid X_1\sqcup \cdots \sqcup X_{k-1}}$.
As $X_1\sqcup \cdots \sqcup X_{k-1}$ is leaf-compatible in $F$, $L\left(F_{\mid X_1\sqcup \cdots \sqcup X_{k-1}}\right)\subseteq L(F)$, so $X_i$ is leaf-compatible in $F$. 

Let us assume that $x\in X_i$, $y\in X_j$, and $x\leq_F y$. 
\begin{itemize}
\item If $i=k$, as $X_k$ is a lower part of $F$, $y\in X_j$, so $i=k\leq j=k$.
\item If $j=k$, obviously $i\leq j$.
\item If $i,j<k$, then $x\leq_{F_{\mid X_1\sqcup \cdots \sqcup X_{k-1}}} y$. By the induction hypothesis, $i\leq j$.
\end{itemize}
In all cases, $i\leq j$.

$\Longleftarrow$. The second condition implies that $X_1\sqcup \cdots \sqcup X_{k-1}$ is a lower part of $F$ and $X_k$ is an upper part of $F$. The first condition implies that both are leaf-compatible, so
\[\Delta_{X_1\sqcup \cdots \sqcup X_{k-1},X_k}(F)\neq 0.\]
The induction hypothesis gives that $(X_1,\ldots,X_{k-1})\in \calC_-^L\left(F_{\mid X_1\sqcup \cdots \sqcup X_{k-1}}\right)$, so 
\[\Delta_{X_1,\ldots,X_k}(F)=\left(\Delta_{X_1,\ldots,X_{k-1}}\otimes \id_{\calR\calF[X_k]}\right)\circ \Delta_{X_1\sqcup \cdots \sqcup X_{k-1},X_k}(F)\neq 0,\]
and $(X_1,\ldots,X_k)\in \calC_-^L(F)$.
\end{proof}

\begin{prop}\label{prop6.1.5}
Let $X$ be a finite set, and $F\in \calR\calF[X]$. For any $\sim\in \calE[X]$, $\sim\in \calE^L[F]$ if, and only if, for any $Y\in X/{\sim}$, $Y$ is convex, connected and leaf-compatible in $F$. 
\end{prop}

\begin{proof}
Let $\sim \in \calE^L[F]$. There exists $(X_1,\ldots,X_k)\in \calC^L_-(F)$, such that:
\begin{enumerate}
\item $X/{\sim}=\{X_1,\ldots,X_k\}$.
\item For any $i\in [k]$, $F_{\mid X_i}$ is indecomposable.
\end{enumerate}
By Lemma \ref{lem6.1.4}, for any $i\in [k]$, $X_i$ is leaf-compatible, and is connected as $F_{\mid X_i}$ is indecomposable. Let us show that it is convex. Let $x$, $y$, $z\in X$, such that $x\leq_F y\leq_F z$, and $x$, $z\in X_i$. Let $j\in [k]$, such that $y\in X_j$.
By the second condition in Lemma \ref{lem6.1.4}, $j\leq i$ as $y\leq_F z$, and $i\leq j$ as $x\leq_F y$. So $y\in X_i$.

Let us assume that for any $Y\in X/{\sim}$, $Y$ is convex, connected and leaf-compatible in $F$. 
Let us consider the relation $\leq_F'$ on $X$ defined as the transitive closure of $\leq_F\cup \sim$.
It is a quasi-order, and it induces a quasi-order $\leq_{F/{\sim}}$ on $X/{\sim}$. Let us assume that $x\leq_F' y$ and $y\leq_F' x$.
There exists elements of $X$, such that
\begin{align*}
x&=x_0\leq_F x'_0\sim x_1\leq_F x'_1\cdots \sim x_p\leq_F x'_p=y,\\
y&=y_0\leq_F y'_0\sim y_1\leq_F y'_1\cdots \sim y_q\leq_F y'_q=x.
\end{align*}
As the elements of $X/{\sim}$ are convex, all the $x_i$, $x'_i$, $y_i$ and $y'_i$ are in the same class, so in particular $x\sim y$. Consequently, $\leq_{F/{\sim}}$ is an order on $X/{\sim}$.
Let $\leq_T$ be a total order on $X/{\sim}$ extending $\leq_{F/{\sim}}$, and let us index $X/{\sim}=\{X_1,\ldots,X_k\}$ according to this order. If $x\in X_i$, $y\in X_j$, and $x\leq_F y$, then $X_i\leq_{F/{\sim}} X_j$, 
so $X_i \leq_T X_j$, and finally $i\leq j$. 
Therefore, $(X_1,\ldots,X_k)\in \calC_-^L(F)$, according to Lemma \ref{lem6.1.4}.
Moreover, for any $i\in [k]$, $F_{\mid X_i}$ is indecomposable, as it is connected.
So $\sim \in \calE^L[F]$.\end{proof}

\begin{prop}
Let $F\in \calR\calF[X]$, and $\sim \in \calE^L[F]$. Then $\calC_+^L(F,\sim)=\calC_+(F/{\sim})$,
where $F/{\sim}\in \calP[X/{\sim}]$ is defined as the order on $X/{\sim}$ induced by the transitive closure of $\leq_F\cup \sim$.
\end{prop}

\begin{proof}
We already noticed that $\leq_{F/{\sim}}$ is an order on $X/{\sim}$ in the proof of Proposition \ref{prop6.1.5}. Let $(X_1/{\sim_1},\ldots,X_k/{\sim_k})\in \comp[X/{\sim}]$.
Then it belongs to $\calC_-^L(F,\sim)$ if, and only if, $(X_1,\ldots,X_k)\in \calC_-^L(F)$ and for any $i\in [k]$,
$\sim_i=\sim_{F_{\mid X_i}}^{(\ic)}$. This equivalent to $(X_1/{\sim_1},\ldots,X_k/{\sim_k})\in \calC_-[F/{\sim}]$ and for any $i\in [k]$, $\left(F/{\sim}\right)_{\mid X_i}\in \bfD\calP[X_i/{\sim_i}]$, that is to say it is equivalent to
$(X_1/{\sim_1},\ldots,X_k/{\sim_k})\in \calC_+(F/{\sim})$.
\end{proof}

\begin{cor}\label{cor6.1.7}
Let $X$ be a finite set, and $F\in \calR\calF[X]$. Then 
\[S_{\calR\calF}(F)=\sum_{\sim\in \calE^L[F]} (-1)^{\cl(\sim)} F\mid \sim.\]
\end{cor}

\begin{proof}
We use Corollary \ref{cor4.2.13}. We obtain
\[S_{\calR\calF}(F)=\sum_{\sim\in \calE^L[F]} \Phi_+(\calC_+(F/{\sim}))(-1) F\mid \sim.\]
By the reciprocity principle for posets (Corollary \ref{cor5.3.5}), for any $\sim\in \calC_+(F/{\sim})$,
\[\Phi_+(\calC_+(F/{\sim}))(-1)=\phi_+^\calP(F/{\sim})(-1)=(-1)^{|X/{\sim}|}=(-1)^{\cl(\sim),}\]
which gives the result.
\end{proof}

\begin{example}\begin{align*}
S_{\calR\calF}\left(\hspace{-2mm}\begin{array}{c}
\tdcinqsix{$1$}{$3$}{$2$}{$5$}{$4$}
\end{array}\hspace{-2mm}\right)&=-\tdcinqsix{$1$}{$3$}{$2$}{$5$}{$4$}+\tdun{$4$}\tdquatredeux{$1$}{$3$}{$2$}{$5$}+\tdun{$5$}\tdquatredeux{$1$}{$3$}{$2$}{$4$}+\tddeux{$1$}{$3$}\tdtroisun{$2$}{$5$}{$4$}-\tdun{$3$}\tdun{$4$}\tdtroisdeux{$1$}{$2$}{$5$}-\tdun{$3$}\tdun{$5$}\tdtroisdeux{$1$}{$2$}{$4$}
+\tddeux{$1$}{$3$}\tddeux{$2$}{$5$}\tdun{$4$}+\tddeux{$1$}{$3$}\tddeux{$2$}{$4$}\tdun{$5$}.
\end{align*}\end{example}

Let us now study the universal polynomial invariants of $\calR\calF$.

\begin{prop}\label{prop6.1.8}
Let $X$ be a finite set, and $F\in \calR\calF[X]$.
\begin{enumerate}
\item For any $N\in \N_{>0}$, $\phi_-^{\calR\calF}(F)(N)$ is the number of maps $c:X\longrightarrow [N]$ such that the two following conditions hold:
\begin{enumerate}
\item For any $x$,$y\in X$, $x\leq_F y\Longrightarrow c(x)\leq c(y)$.
\item For any $i\in [N]$, if there exists $x\in X$, such that $c(x)=i$, then there exists $y\in L(F)$, such that $c(y)=i$.
\end{enumerate}
\item $\phi_+^{\calR\calF}(F)=\phi_{\pp}^{\calR\calF}(F)=\begin{cases}
T^{|X|} \mbox{ if $F$ has no edge},\\
0\mbox{ otherwise}.
\end{cases}$
\end{enumerate}
Therefore, $\calR\calF$ is degenerate, and $\nd\calR\calF=\bfD\calR\calF$. \end{prop}

\begin{proof}
We use Lemma \ref{lem6.1.4} for the description of $\calC_-^L(F)$. If $(X_1,\ldots,X_k)\in \calC^-(F)$, let us consider the surjective map $c:X\longrightarrow [k]$, sending $x\in X_i$ to $i$ for $i\in [k]$. 
The second condition of Lemma \ref{lem6.1.4} implies condition (a). As any $X_i$ is leaf-compatible, condition (b) is satisfied. Conversely, if $c:X\longrightarrow [k]$ is a surjective map satisfying conditions (a) and (b), putting $X_i=c^{-1}(i)$ for $i\in [k]$,
condition (b) implies that $X_i$ is leaf-compatible for any $i\in [k]$, and condition (a) implies that $X_i$ is convex for any $i\in [k]$ and condition 2 of Lemma \ref{lem6.1.4}. 
Consequently, elements of $\calC_-^L(F)$ are in bijection with surjective maps satisfying conditions (a) and (b). The formula for $\phi_-^{\cal\calF}(F)$ follows.

Let us assume that $\calC_+^L(F)\neq \emptyset$, and let $(X_1,\ldots,X_k)\in \calC_+^L(F)$. For any $i\in [k]$, $X_i$ is leaf-compatible, and $F_{\mid X_i}$ has no edge. Therefore, $X_i\subseteq L(F)$. Consequently, $L(F)=X$, and $F$ has no edge. 
Conversely, if $F$ has no edge, $\calC_+(F)=\comp[X]$. This implies the formula for $\phi_+^{\calR\calF}(F)$. As $\phi_{\pp}^{\calR\calF}(F)$ is the monomial of degree $|X|$ of $\phi_+^{\calR\calF}(F)$, the formula for $\phi_{\pp}^{\calR\calF}(F)$ follows.
\end{proof}

\begin{example}
Let $X$ be a finite set, and let $F\in \calR\calF[X]$. If $F$ is a ladder (that is to say a tree with only one leaf), then $\phi_-^{\calR\calF}(F)=T$: if $c:X\longrightarrow [N]$ is a valid coloration, in the sense of Proposition \ref{prop6.1.8}, then all the vertices should have the same color, so $\phi_-^{\calR\calF}(F)(N)=N$.
If $F$ is a corolla (that is to say a tree such that all the vertices different from the root are leaves), then $\phi_-^{\calR\calF}(F)=T^{|X|-1}$: 
valid colorations $c:X\longrightarrow [N]$ in the sense of Proposition \ref{prop6.1.8} are obtained by arbitrarily coloring the $|X|-1$ leaves, and the color the root by the minimal chosen color, so $\phi_-^{\calR\calF}(F)(N)=N^{|X|-1}$.
Here are other examples.
\begin{align*}
\phi_-^{\calR\calF}\left(\hspace{-2mm}\begin{array}{c}
\tquatrequatre
\end{array}\hspace{-2mm}\right)=\phi_-^{\calR\calF}\left(\hspace{-2mm}\begin{array}{c}
\tcinqonze
\end{array}\hspace{-2mm}\right)&=\frac{T(3T-1)}{2},&\phi_-^{\calR\calF}\left(\hspace{-2mm}\begin{array}{c}
\tcinqdeux
\end{array}\hspace{-2mm}\right)&=\frac{T(4T^2-3T+1)}{2},\\
\phi_-^{\calR\calF}\left(\hspace{-2mm}\begin{array}{c}
\tcinqcinq
\end{array}\hspace{-2mm}\right)=\phi_-^{\calR\calF}\left(\hspace{-2mm}\begin{array}{c}
\tcinqhuit
\end{array}\hspace{-2mm}\right)&=T(2T-1),&\phi_-^{\calR\calF}\left(\hspace{-2mm}\begin{array}{c}
\tcinqsix
\end{array}\hspace{-2mm}\right)&=\frac{T(8T^2-3T+1)}{6},\\
\phi_-^{\calR\calF}\left(\hspace{-2mm}\begin{array}{c}
\tcinqdix
\end{array}\hspace{-2mm}\right)&=T^3,&\phi_-^{\calR\calF}\left(\hspace{-2mm}\begin{array}{c}
\tcinqtreize
\end{array}\hspace{-2mm}\right)&=T^2.
\end{align*}
We leave to the reader as an exercise to prove that if $F\in \calR\calF[X]$ has two leaves, then 
\[\phi_-^{\calR\calF}(F)=\frac{T(nT-n+2)}{2},\]
where $n$ is the number of vertices of $F$ which are ascendants of only one leaf of $F$, including the two leaves themselves.
\end{example} 

\begin{remark}
It is not difficult to prove that, for any finite set $X$, and any $F\in\calF\calR[X]$, the degree of $\phi^{\calR\calF}_-(F)$ is the number of leaves of $F$. 
\end{remark}

\begin{prop}
Let $X$ be a finite set, $F\in \calR\calF[X]$, $N\geq 1$ and $f:X\longrightarrow [n]$ be a bijection. Then the set of $\calR\calF$-partitions of $(F,f)$ is $[N]^X$ if $F$ is discrete, and $\emptyset$ otherwise.
\end{prop}

\begin{proof}
By Lemma \ref{lem6.1.4}, the parts of any element of $\calC_{\pp}(F)$ are reduced to a single leaf, by the leaf-compatibility condition. Therefore, if $(F,f)$ has a $\calR\calF$-partition, then any of its vertex is a leaf, so $F$ is discrete. Conversely, if $F$ is discrete, then $L(F)=X$, and by Lemma \ref{lem6.1.4}, $\calC_{\pp}(F)=\comp_n[X]$.
As a consequence, any map $g:X\longrightarrow [N]$ is an $\calR\calF$-partition of $F$.
\end{proof}

Therefore, if $\Omega$ is a totally ordered set, $\ell:X\longrightarrow\Omega$ is a map, and $F\in \calR\calF[X]$,
by Example \ref{ex5.2.6}, the formal flow is given by
\[F_\Omega^{\calR\calF}(F,\ell)=\begin{cases}
T^n\mbox{ if $F$ is discrete},\\
0\mbox{ otherwise}.
\end{cases}\]

\section{A combinatorial twisted bialgebra on rooted forests based on leaves}

\subsection{Definition}

\begin{defi}\label{defi6.2.1}
Let $X$ be a finite set. We denote by $\calR\calF_\ell[X]$ the set of isoclasses of rooted forests whose set of leaves is indexed by the elements of $X$. For example,
\begin{align*}
\calR\calF_\ell[0]&=\{1\},&
\calR\calF_\ell[1]&=\left\{\begin{array}{c}
\tdun{$1$},\tddeux{}{$1$},\tdtroisdeux{}{}{$1$},\tdquatrecinq{}{}{}{$1$},\ldots\end{array}\right\},&
\calR\calF_\ell[2]&=\left\{\begin{array}{c}
\tdun{$1$}\tdun{$2$},\tdun{$1$}\tddeux{}{$2$},\tddeux{}{$1$}\tdun{$2$},\tddeux{}{$1$}\tddeux{}{$2$},\ldots\\
\tdtroisun{}{$2$}{$1$},\tdquatredeux{}{$2$}{}{$1$},\tdquatredeux{}{$1$}{}{$2$},\tdquatrequatre{}{}{$2$}{$1$},\ldots
\end{array}\right\}.
\end{align*}
This defines a species $\calR\calF_\ell$. Note that if $X$ is non-empty, then $\calR\calF_\ell[X]$ is infinite.
\end{defi}

\begin{notation}
Let $X$ be a finite set and $F\in \calR\calF_\ell[X]$. 
\begin{enumerate}
\item Let $v$ be a vertex of $F$. We denote by $L_v(F)$ the set of leaves of $F$ which are descendants of $v$: this is a non-empty subset of $X$.
\item Let $Y\subseteq X$. We denote by $F_{\mid_\subset Y}$ the sub-forest of $F$ whose vertices are the vertices $v\in V(F)$ such that $L_v(F)\subseteq Y$, and whose edges are all the edges between these vertices.
\item Let $Y\subseteq X$. We denote by $F_{\mid_\cap Y}$ the sub-forest of $F$ whose vertices are the vertices $v\in V(F)$ such that $L_v(F)\cap Y\neq \emptyset$, and whose edges are all the edges between these vertices.
\end{enumerate}\end{notation}

\begin{remark}
Let $X$ be a finite set, $F\in \calR\calF_\ell[X]$ and $Y\subseteq X$. Then, by construction, $F_{\mid_\cap Y}$ and $F_{\mid_\subset Y}$ both belong to $\calR\calF_\ell[Y]$. Moreover, 
\[V\left(F_{\mid_\subset Y}\right)\sqcup V\left(F_{\mid_\cap X\setminus Y}\right)=V(F).\]
Indeed, if $v\in V(F)$, then, either $L_v(F)\subseteq Y$, or $L_v(F)\cap \cal V(F)\setminus Y\neq \emptyset$, 
and these two conditions are incompatible.
\end{remark}

With the disjoint union product, $\calR\calF_\ell$ is a twisted algebra. Let us give it a coproduct.

\begin{theo}\label{theo6.2.2}
Let us define a coproduct on $\calR\calF_\ell$ as follows. Let $(X,Y)$ be a pair of disjoint finite sets, and $F\in \calR\calF_\ell[X\sqcup Y]$. We put
\[\Delta^\ell_{X,Y}(F)=F_{\mid_\cap X}\otimes F_{\mid_\subset Y}.\]
Then $(\calR\calF_\ell,m,\Delta^\ell)$ is a CTB.
\end{theo} 

\begin{proof}
Let $X$ and $Y$ be disjoint finite sets, $F\in \calR\calF[X]$ and $G\in \calR\calF[Y]$. If $X'\sqcup Y'=X\sqcup Y$, then 
\begin{align*}
(FG)_{\mid_\subset X'}&=F_{\mid_\subset X'\cap X}G_{\mid_\subset X'\cap Y},&
(FG)_{\mid_\cap Y'}&=F_{\mid_\cap Y'\cap X}G_{\mid_\cap Y'\cap Y}.
\end{align*}
Therefore, $\Delta_{X',Y'}(FG)=\Delta_{X'\cap X,Y'\cap X}(F)\Delta_{X'\cap Y,Y'\cap Y}(G)$.

Let $X$, $Y$, $Z$ be pairwise disjoint finite sets, and $F\in \calR\calF[X\sqcup Y\sqcup Z]$.
\begin{align*}
\left(\Delta^\ell_{X,Y}\otimes \id_{\calR\calF_\ell[Z]}\right)\circ \Delta^\ell_{X\sqcup Y,Z}(F)&=
\left(F_{\mid_\cap X\sqcup Y}\right)_{\mid_\cap X}\otimes \left(F_{\mid_\cap X\sqcup Y}\right)_{\mid_\subset Y}\otimes F_{\mid_\cap Z},\\
\left(\id_{\calR\calF_\ell[X]}\otimes \Delta^\ell_{Y,Z}\right)\circ \Delta^\ell(F)&=F_{\mid_\cap X}\otimes 
\left(F_{\mid_\subset Y\sqcup Z}\right)_{\mid_\cap Y}\otimes \left(F_{\mid_\subset Y\sqcup Z}\right)_{\mid_\subset Z}.
\end{align*}
Let $x\in V(F)$.
\begin{align*}
x\in V\left(\left(F_{\mid_\cap X\sqcup Y}\right)_{\mid_\cap X}\right)&\Longleftrightarrow L_x(F)\cap (X\sqcup Y)\neq \emptyset \mbox{ and }L_x\left(F_{\mid_\cap X\sqcup Y}\right)\cap X\neq \emptyset\\
&\Longleftrightarrow L_x(F)\cap (X\sqcup Y)\neq \emptyset \mbox{ and } L_x(F)\cap (X\sqcup Y)\cap X\neq \emptyset\\
&\Longleftrightarrow L_x(F)\cap X\neq \emptyset\\
&\Longleftrightarrow x\in V\left(F_{\mid_\cap X}\right),
\end{align*}
so $\left(F_{\mid_\cap X\sqcup Y}\right)_{\mid_\cap X}=F_{\mid_\cap X}$.
\begin{align*}
x\in V\left(\left(F_{\mid_\subset Y\sqcup Z}\right)_{\mid_\subset Z}\right)&\Longleftrightarrow L_x(F)\subseteq Y\sqcup Z\mbox{ and }L_x\left(F_{\mid_\subset Y\sqcup Z}\right)\subseteq Z\\
&\Longleftrightarrow L_x(F)\subseteq Y\sqcup Z\mbox{ and }L_x(F)\subseteq Z\\
&\Longleftrightarrow L_x(F)\subseteq Z\\
&\Longleftrightarrow x\in V\left(F_{\mid_\subset Z}\right),
\end{align*}
so $\left(F_{\mid_\subset Y\sqcup Z}\right)_{\mid_\subset Z}=F_{\mid_\subset Z}$.
\begin{align*}
x\in V\left(\left(F_{\mid_\cap X\sqcup Y}\right)_{\mid_\subset Y}\right)&\Longleftrightarrow
L_x(F)\cap (X\sqcup Y)\neq \emptyset \mbox{ and }L_x\left(F_{\mid_\cap X\sqcup Y}\right)\subseteq Y\\
&\Longleftrightarrow L_x(F)\cap (X\sqcup Y)\neq \emptyset \mbox{ and }L_x(F)\cap (X\sqcup Y)\subseteq Y\\
&\Longleftrightarrow L_x(F)\subseteq Y\sqcup Z\mbox{ and }L_x(F)\cap Y\neq \emptyset\\
&\Longleftrightarrow L_x(F)\subseteq Y\sqcup Z\mbox{ and }L_x\left(F_{\mid_\subset Y\sqcup Z}\right)\cap Y\neq \emptyset\\
&\Longleftrightarrow x\in V\left(\left(F_{\mid_\subset Y\sqcup Z}\right)_{\mid_\cap Y}\right),
\end{align*}
so $\left(F_{\mid_\cap X\sqcup Y}\right)_{\mid_\subset Y}=\left(F_{\mid_\subset Y\sqcup Z}\right)_{\mid_\cap Y}$.
As a consequence, $\Delta^\ell$ is coassociative.

Let $X$ be a finite set, and $F\in \calR\calF^\ell[X]$.
\begin{align*}
\Delta^\ell_{\emptyset,X}(F)&=1\otimes F,&\Delta^\ell_{X,\emptyset}(F)&=F\otimes 1.
\end{align*}
So $\Delta^\ell$ is counitary. 
\end{proof}

\begin{example}
\begin{align*}
\Delta^\ell_{\emptyset,[2]}\left(\hspace{-1mm}\begin{array}{c}
\tdquatrequatre{}{}{$2$}{$1$}
\end{array}\hspace{-1mm}\right)&=1\otimes \tdquatrequatre{}{}{$2$}{$1$},&
\Delta^\ell_{[2],\emptyset}\left(\hspace{-1mm}\begin{array}{c}
\tdquatrequatre{}{}{$2$}{$1$}
\end{array}\hspace{-1mm}\right)&=\tdquatrequatre{}{}{$2$}{$1$}\otimes 1,\\
\Delta^\ell_{\{1\},\{2\}}\left(\hspace{-1mm}\begin{array}{c}
\tdquatrequatre{}{}{$2$}{$1$}
\end{array}\hspace{-1mm}\right)&=\tdtroisdeux{}{}{$1$}\otimes \tdun{$2$},&
\Delta^\ell_{\{2\},\{1\}}\left(\hspace{-1mm}\begin{array}{c}
\tdquatrequatre{}{}{$2$}{$1$}
\end{array}\hspace{-1mm}\right)&=\tdtroisdeux{}{}{$2$}\otimes \tdun{$1$}.
\end{align*}
This example shows that $(\calR\calF_\ell,m,\Delta^\ell)$ is not quasi-cocommutative.
\end{example}

\begin{prop}
Let $X$ be a finite set, and $F\in \calR\calF_\ell[X]$.
\begin{enumerate}
\item $\phi_-^{\calR\calF_\ell}(F)=\phi_{\pp}^{\calR\calF_\ell}(F)=T^{|X|}$.
\item If $F$ is a tree, $\phi_+^{\cal\calF_\ell}(F)=T(T-1)\cdots (T-|X|+1)$.
\end{enumerate}
\end{prop}

\begin{proof}
As $\calC_-^\ell(F)=\comp[X]$, $\phi_-^\ell(F)=T^{|X|}$. 
Let us assume that $F$ is a tree. If $(X_1,\ldots,X_k) \in \calC^\ell_+(F)$, then for any $i\in [k]$, $F_{\mid_i X_i}$ is a tree, with only one leaf by definition of $\calC_+^\ell(F)$. Therefore, the $X_i$'s are singletons. 
As a consequence, $\calC_+^{\calR\calF_\ell}(F)=\comp_{|X|}[X]$, so $\phi_+^{\cal\calF_\ell}(F)=T(T-1)\cdots (T-|X|+1)$. As $\phi_{\pp}^{\calR\calF_\ell}(F)$ is the monomial of $\phi_+^{\calR\calF_\ell}(F)$ of degree $|X|$, the formula for $\phi_{\pp}^{\calR\calF_\ell}(F)$ follows.
\end{proof}

\subsection{Morphisms to hypergraphs and reduced forests}

\begin{theo}\label{theo6.2.3}
Let $X$ be a finite set, and $F\in \calR\calF_\ell[X]$. The hypergraph $\Xi(F)\in \calH\calG[X]$ is defined by
\[\HE(\Xi(F))=\{L_x(F)\mid x\in V(F)\}.\]
Then $\Xi:(\calR\calF_\ell,m,\Delta^\ell)\longrightarrow (\calH\calG,m,\Delta^{\cap,\subset})$ is a morphism of CTBs.
\end{theo}

\begin{proof}
Let $(X,Y)$ be a pair of disjoint sets, $F\in \calR\calF_\ell[X]$ and $G\in \calR\calF_\ell[Y]$. Then $V(FG)=V(F)\sqcup V(G)$, and
\begin{align*}
&\forall x\in V(FG)=V(F)\sqcup V(G),&L_x(FG)&=\begin{cases}
L_x(F)\mbox{ if }x\in V(F),\\
L_x(G)\mbox{ if }x\in V(G).
\end{cases}
\end{align*}
Consequently,
\[\HE(\Xi(FG))=\{L_x(F)\mid x\in V(F)\}\sqcup \{L_x(G)\mid x\in V(G)\}=\HE(\Xi(F))\sqcup \HE(\Xi(G))=\HE(\Xi(F)\Xi(G)).\]
So $\Xi(FG)=\Xi(F)\Xi(G)$.

Let $(X,Y)$ be a pair of disjoint sets, and $F\in \calR\calF_\ell[X\sqcup Y]$. Then
\begin{align*}
HE\left(\Xi\left(F_{\mid_\cap X}\right)\right)&=\left\{L_x\left(F_{\mid_\subset X}\right)\mid x\in V\left(F_{\mid_\cap X}\right)\right\}\\
&=\left\{L_x(F)\cap X\mid x\in V(F),\: L_x(F)\cap X\neq \emptyset\right\}\\
&=\HE\left(\Xi(F)_{\mid_\cap X}\right).
\end{align*}
So $\Xi\left(F_{\mid_\cap X}\right)=\Xi(F)_{\mid_\cap X}$. Moreover,
\begin{align*}
HE\left(\Xi\left(F_{\mid_\subset Y}\right)\right)&=\left\{L_x\left(F_{\mid_\subset Y}\right)\mid x\in V\left(F_{\mid_\subset Y}\right)\right\}\\
&=\left\{L_x(F)\mid x\in V(F),\: L_x(F)\subseteq Y\right\}\\
&=\HE\left(\Xi(F)_{\mid_\subset Y}\right).
\end{align*}
So $\Xi\left(F_{\mid_\subset Y}\right)=\Xi(F)_{\mid_\subset Y}$. Therefore,
\[\Delta_{X,Y}\circ \Xi(F)=\Xi(F)_{\mid_\cap X}\otimes \Xi(F)_{\mid_\subset Y}
=\Xi\left(F_{\mid_\cap X}\right)\otimes \Xi\left(F_{\mid_\subset Y}\right)=(\Xi\otimes \Xi)\circ \Delta_{X,Y}(F).\]
We proved that $\Xi$ is a morphism of CTBs.
\end{proof}

This morphism is not injective: for example, any rooted tree in $\calR\calF_\ell[1]$ is sent to the hypergraph on $[1]$, which unique hyperedge is $[1]$. 
To study the kernel of $\iota$, we introduce reduced rooted forests.

\begin{defi}
Let $\calF$ be a rooted forest. We shall say that it is reduced if for any vertex $x$ of $F$, the number of children of $x$ is different from 1.
For any finite set $X$, the set of isoclasses of reduced rooted forests whose leaves are indexed by $X$ is denoted by $\calR\calR\calF_\ell[X]$. For example,
\begin{align*}
\calR\calR\calF_\ell[0]&=\{1\},&
\calR\calR\calF_\ell[1]&=\left\{\begin{array}{c}
\tdun{$1$}\end{array}\right\},&
\calR\calR\calF_\ell[2]&=\left\{\begin{array}{c}
\tdun{$1$}\tdun{$2$},\tdtroisun{}{$2$}{$1$}
\end{array}\right\},
\end{align*}
\begin{align*}
\calR\calR\calF_\ell[3]&=\left\{\begin{array}{c}
\tdun{$1$}\tdun{$2$}\tdun{$3$},\tdun{$1$}\tdtroisun{}{$3$}{$2$},\tdun{$2$}\tdtroisun{}{$3$}{$1$},\tdun{$3$}\tdtroisun{}{$2$}{$1$},\tdquatreun{}{$3$}{$2$}{$1$},\tdcinqsix{}{$3$}{}{$2$}{$1$},\tdcinqsix{}{$2$}{}{$3$}{$1$},\tdcinqsix{}{$1$}{}{$3$}{$2$}
\end{array}\right\}.
\end{align*}\end{defi}
 
The cardinality of $\calR\calR\calF_\ell[X]$ according to $|X|$ can be found in \cite[Entry A006351]{Sloane}.

\begin{defi}
Let $F$ be a rooted forest. We denote by $\red(F)$ the forest obtained by contraction of all the edges of $F$ with only one child.
Then $\red(F)$ is a reduced rooted forest, and $L(F)=L(\red(F))$.
\end{defi}

\begin{example}\begin{align*}
\red(\tdun{$1$})=\red\left(\hspace{-1mm}\begin{array}{c}\tddeux{}{$1$}\end{array}\hspace{-2mm}\right)=\red\left(\hspace{-1mm}\begin{array}{c}\tdtroisdeux{}{}{$1$}\end{array}\hspace{-2mm}\right)
=\red\left(\hspace{-1mm}\begin{array}{c}\tdquatrecinq{}{}{}{$1$}\end{array}\hspace{-2mm}\right)&=\cdots=\tdun{$1$},\\
\red\left(\hspace{-2mm}\begin{array}{c}\tdquatredeux{}{$2$}{}{$1$}\end{array}\hspace{-2mm}\right)=\red\left(\hspace{-2mm}\begin{array}{c}\tdquatrequatre{}{}{$2$}{$1$}\end{array}\hspace{-2mm}\right)&=\tdtroisun{}{$2$}{$1$}.
\end{align*}\end{example}

The map $\red$ is a species morphism from $\calR\calF_\ell$ to $\calR\calR\calF_\ell$, as it conserves the set of leaves. Moreover, $\red\circ \red=\red$. 
It cannot be interpreted as a species morphism on $ \calR\calF$, as it does not conserve the set of vertices. Moreover, $\calR\calR\calF_\ell$ is a twisted subalgebra of $\calR\calF_\ell$, as the disjoint union of two reduced rooted forest is also a reduced rooted forest,
but not a twisted subbialgebra, as shown by the following example:
\[\Delta^\ell_{\{1\},\{2\}}\left(\hspace{-2mm}\begin{array}{c}\tdtroisun{}{$2$}{$1$}\end{array}\hspace{-2mm}\right)=\tddeux{}{$1$}\otimes \tdun{$2$}.\]

\begin{lemma}\label{lem6.2.6}
Let $X$ be a finite set and $F\in \calR\calF_\ell[X]$. Then $\Xi(F)=\Xi(\red(F))$.
\end{lemma}

\begin{proof}
Let $x$ be a vertex of $F$ with only one child, denoted by $y$, and let $F'$ be the rooted forest obtained after the contraction of the edge of $F$ from $x$ to $y$ (that is to say, obtained after the identification of $x$ and $y$). Then
\[L_x(F)=L_y(F)=L_x(F'),\]
and if $z\in V(F)$ is different from $x$ and $y$, $L_z(F)=L_z(F')$. Therefore, as $L(F)=L(F')$, $\Xi(F)=\Xi(F')$. Iterating this for all the vertices of $F$ with only one child, we finally obtain that $\Xi(F)=\Xi(\red(F))$.
\end{proof}

\begin{lemma}\label{lem6.2.7}
The restriction of $\red$ to $\cal\calR\calF_\ell$ is injective.
\end{lemma}

\begin{proof}
Let $F\in \calR\calR\calF_\ell[X]$, and let $x$ and $y$ be vertices of $F$. Then:
\begin{itemize}
\item If $x\leq_F y$, then $L_y(F)\subseteq L_x(F)$.
\item If $x\leq_F y$ and $x\neq y$, as $x$ has at least one child which is not on the path from $x$ to $y$ (as $F$ is reduced), $L_y(F)\subsetneq L_x(F)$.
\item If $x$ and $y$ are not comparable for $\leq_F$, then $L_x(F)\cap L_y(F)\neq \emptyset$.
\end{itemize}
Therefore, if $x\neq y$, then $L_x(F)\neq L_y(F)$, and 
\[x\leq_F y\Longleftrightarrow L_y(F)\subseteq L_x(F).\]
Consequently, the poset $(V(F),\leq_F)$ is isomorphic to the poset $(\HE(\Xi(F)),\supseteq)$, and $F$ is the Hasse graph of $(\HE(\Xi(F)),\supseteq)$.
Therefore, if $F,G\in \calR\calR\calF_\ell[X]$ and $\Xi(F)=\Xi(G)$, then $F=G$. 
\end{proof}

\begin{theo}\label{theo6.2.8}
The kernel of the linearization of $\Xi$ is given by
\[\ker(\Xi)[X]=\vect(F-\red(F)\mid F\in \calR\calF_\ell[X]),\]
for any finite set $X$. Consequently, $\calR\calR\calF_\ell$ inherits a structure of CTB, and in particular a coproduct $\Delta^{r\ell}$, such that $\Xi:(\calR\calR\calF_\ell,m,\Delta^{r\ell})\longrightarrow (\calH\calG,m,\Delta^{\cap,\subset})$ is an injective morphism of CTBs.
\end{theo}

\begin{proof}
Lemma \ref{lem6.2.6} implies that for any $F\in \calR\calF_\ell[X]$, $F-\red(F)\in \ker(\Xi)[X]$. Let us consider $b\in \ker(\xi)[X]$. We can write
\[b=\sum_{F\in \calR\calF_\ell[X]} a_FF=\sum_{F\in \calR\calR\calF_\ell[X]}\left(\sum_{G\in \calR\calF_\ell[X],\: \red(G)=F}a_GG\right).\]
By Lemma \ref{lem6.2.6},
\[\Xi(b)=\sum_{F\in \calR\calR\calF_\ell[X]}\left(\sum_{G\in \calR\calF_\ell[X],\: \red(G)=F}a_G\Xi(G)\right)=\sum_{F\in \calR\calR\calF_\ell[X]}\left(\sum_{G\in \calR\calF_\ell[X],\: \red(G)=F}a_G\right)\Xi(F)=0.\]
By Lemma \ref{lem6.2.7}, the $\Xi(F)$'s, when $F$ runs in $\calR\calR\calF_\ell[X]$, are different hypergraphs. Therefore, for any $F\in \calR\calR\calF_\ell[X]$,
\[\sum_{G\in \calR\calF_\ell[X],\: \red(G)=F}a_G=0.\]
Therefore,
\begin{align*}
\sum_{G\in \calR\calF_\ell[X],\: \red(G)=F}a_GG&=a_FF+\sum_{G\in \calR\calF_\ell[X]\setminus \{F\},\: \red(G)=F}a_GG\\
&=-\sum_{G\in \calR\calF_\ell[X]\setminus \{F\},\: \red(G)=F}a_GF+\sum_{G\in \calR\calF_\ell[X]\setminus \{F\},\: \red(G)=F}a_GG\\
&=\sum_{G\in \calR\calF_\ell[X]\setminus \{F\},\: \red(G)=F}a_G(G-\red(G)).
\end{align*}
Summing over all $F$, we obtain that $b\in =\vect(F-\red(F)\mid F\in \calR\calF_\ell[X])$.

Therefore, we can identify the linear species $\ell\calR\calF_\ell/\ker(\Xi)$ with $\ell \calR\calR\calF_\ell$. Therefore, $\calR\calR\calF_\ell$ inherits a structure of CTB, making it isomorphic to $\im(\Xi)$. 
\end{proof}

\begin{remark}
The product of $\calR\calR\calF_\ell$ is again the disjoint union. If $(X,Y)$ is a pair of disjoint finite sets, and $F\in \calR\calR\calF_\ell[X\sqcup Y]$,
\[\Delta^{r\ell}_{X,Y}(F)=\red\left(F_{\mid_\cap X}\right)\otimes \red\left(F_{\mid_\subset Y}\right) =\red\left(F_{\mid_\cap X}\right)\otimes F_{\mid_\subset Y},\]
observing that $F_{\mid_\subset Y}$ is always a reduced forest, whereas $F_{\mid_\cap X}$ is not always reduced.
\end{remark}

\section{Bosonic Fock functor on $\calR\calF$ and $\calR\calF_\ell$}

\subsection{Two identical objects}

As vector spaces, $\calF[\calR\calF]$ and $\calF[\calR\calF_\ell]$ are the same: they have for basis the set of rooted forests
\[\left\{\begin{array}{c}
1;\tun;\tun\tun,\tdeux;\tun\tun\tun,\tdeux\tun,\ttroisun,\ttroisdeux;\tun\tun\tun\tun,\tdeux\tun\tun,\tdeux\tdeux,\ttroisun\tun,\ttroisdeux\tun,\tquatreun,\tquatredeux,\tquatrequatre,\tquatrecinq,\ldots
\end{array}\right\}.\]
They also share the same product, given by the disjoint union. Let us now show that they share the same coproduct:

\begin{theo}\label{theo6.3.1}
As bialgebras, $\calF[(\calR\calF,m,\Delta^L)]=\calF[(\calR\calF_\ell,m,\Delta^\ell)]$.
\end{theo}

\begin{proof}
Let $F$ be a rooted forest. We denote by $\calA(F)$ the set of lower parts $X\subseteq V(F)$, such that both $X$ and $V(F)\setminus X$ are leaf-compatible. Note that $\Delta^L_{X,V(F)\setminus X}\neq 0$ if, and only if, $X\in \calA(F)$. We denote by $\calB(F)$ the set of subsets of $L(F)$. We then define a map
\begin{align*}
\theta&:\left\{\begin{array}{rcl}
\calA(F)&\longrightarrow&\calB(F)\\
X&\longrightarrow& L\left(F_{\mid X}\right).
\end{array}\right.
\end{align*}
As $X$ is leaf-compatible, indeed $L\left(F_{\mid X}\right)\subseteq L(F)$, so $\theta$ is well-defined.
Conversely, we define
\begin{align*}
\theta'&:\left\{\begin{array}{rcl}
\calB(F)&\longrightarrow&\calA(F)\\
Y&\longmapsto&\{x\in V(F)\mid L_x(F)\cap Y\neq \emptyset\}.
\end{array}\right.
\end{align*}
Let us show that $\theta'$ is well-defined. Let $Y\subseteq L(F)$. Let $y\in \theta(Y)$ and $x\in V(F)$, such that $x\leq_F y$. Then $L_x(F)\subseteq L_y(F)$, so $L_x(F)\cap Y\neq \emptyset$. So $\theta'(Y)$ is a lower part of $F$. 
If $x\in \theta'(Y)$, then by definition $L_x(F)\cap Y\neq \emptyset$. Moreover, if $y\in Y$, then $L_y(F)=\{y\}$,
so $Y\subseteq \theta'(Y)$: consequently, $L_x(F)\cap \theta'(Y)\neq \emptyset$ and $\theta'(F)$ is leaf-compatible.
Let $x\in V(F)\setminus \theta'(Y)$. By definition, $L_x(F)\subseteq V(F)\setminus Y$, so $L_x(F)\subseteq V(F)\setminus \theta'(Y)$: we obtain that $V(F)\setminus \theta'(Y)$ is leaf-compatible. Consequently, $\theta'$ is well-defined.

Let $X\in \calA(F)$. We put $Y=L(F_{\mid X})$. If $x\in X$, then $L_x(F)\cap Y=L_x\left(F_{\mid Y}\right)\neq \emptyset$, so $x\in \theta'(Y)$. If $x\in \theta'(Y)$, then $L_x(F)\cap Y\neq \emptyset$. 
Let $y\in L_x(F)\cap Y$. Then $x\leq_F y$ and $y\in Y\subseteq X$. As $X$ is a lower part of $F$, $x\in X$. 
We obtain that $\theta'(Y)=\theta'\circ \theta(X)=X$, so $\theta'\circ\theta=\id_{\calA(F)}$.

Let $Y\in \calB[F]$. We put $\theta'(Y)=X$. Let $x\in \theta(X)$. Then $x$ is a leaf of $F_{\mid X}$.
As $X$ is leaf-compatible, $L\left(F_{\mid X}\right)=L(F)\cap X$, so $x$ is a leaf of $F$. 
By definition of $X$, $L_x(F)\cap Y=\{x\}\cap Y\neq \emptyset$, so $x\in Y$.
Let $y\in Y$. Then obviously $y\in L\left(F_{\mid X}\right)=\theta(X)$. As a consequence, 
$\theta\circ \theta'(Y)=Y$, so $\theta\circ \theta'=\id_{\calB(F)}$.

We proved that $\theta$ and $\theta'$ are bijections, inverse one from the other. Moreover, if $X\in \calA(F)$,
we obtained that $F_{\mid X}=F_{\mid_\cap \theta(X)}$. Taking the complement (on the vertices), 
$F_{\mid V(F)\setminus X}=F_{\mid_\subset L(F)\setminus \theta(X)}$. Therefore,
\begin{align*}
\Delta^L(\overline{F})&=\sum_{X\in \calA(F)} \overline{F_{\mid X}}\otimes \overline{F_{\mid V(F)\setminus X}}=\sum_{Y\in \calB(F)} \overline{F_{\mid_\cap Y}}\otimes \overline{F_{\mid_\subset L(F)\setminus Y}}=\Delta^\ell(\overline{F}),
\end{align*}
where we put $Y=\theta(X)$. Therefore, on $\calF[\calR\calF]=\calF[\calR\calF_\ell]$, $\Delta^L=\Delta^\ell$. 
\end{proof}

\begin{remark}
Of course, they are not identical as graded bialgebras: $\calF[(\calR\calF,m,\Delta^L)]$ is graded by the number of vertices, whereas $\calF[(\calR\calF_\ell,m,\Delta^\ell)]$ is graded by the number of leaves.
\end{remark}

Applying the Fock functor, we obtain also a bialgebra on reduced forests $\calF[\calR\calR\calF_\ell]$. The product is given by disjoint union, and for any reduced rooted forest $F$,
\[\Delta^{r\ell}(F)=\sum_{Y\subseteq L(F)} \overline{\red\left(F_{\mid_\cap Y}\right)}\otimes \overline{F_{\mid_\subset V(F)\setminus Y}}.\]
This Hopf algebra appears in \cite{Xiao2026}, where it is studied and related to the Hopf algebra of trees appearing in \cite{Bruned2015}. The following map is a surjective morphism of bialgebras:
\[\calF[\red]:\left\{\begin{array}{rcl}
\left(\calF[\calR\calF_\ell],m,\Delta^L\right)&\longrightarrow&\left(\calF[\calR\calR\calF_\ell],m,\Delta^{r\ell}\right)\\
\overline{F}&\longmapsto&\overline{\red(F)}.
\end{array}\right.\]
The following map is an injective morphism of bialgebras:
\[\calF[\Xi]:\left\{\begin{array}{rcl}
\left(\calF[\calR\calR\calF_\ell],m,\Delta^{r\ell}\right)&\longrightarrow&\left(\calF[\calH\calG],m,\Delta^{\cap,\subset}\right)\\
\overline{F}&\longmapsto&\overline{\Xi(F)}.
\end{array}\right.\]

\subsection{Antipode for rooted forests and reduced forests}

Both $(\calF[\calR\calR\calF],m,\Delta^L)$ and $(\calF[\calR\calR\calF_\ell],m,\Delta^{r\ell})$ are Hopf algebras.
The antipode of the first one is given by $\calF\left[S_{\calR\calF}\right]$, so, for any rooted forest $F$,
\[S_{\calF[\calR\calF]}(\overline{F})=\sum_{\sim\in \calE^L[\calF]}(-1)^{\cl(\sim)} \overline{F\mid\sim}.\]
As $\calF[\red]$ is a bialgebra morphism from $\calF[\calR\calF]$ to $\calF[\calR\calR\calF_\ell]$,
it is automatically a Hopf algebra morphism, so
\[S_{\calF[\calR\calR\calF_\ell]}\circ \red=\red\circ S_{\calF[\calR\calF]}.\]
This allows to give a formula for this antipode, solving a question in \cite{Xiao2026}. For this, we need the following result:

\begin{prop}
Let $F$ be a reduced rooted forest, and $I$ be a nonempty subset of $L(F)$. We shall say that $I$ if tree-compatible
if all the elements of $I$ belong to the same tree of $F$. If so, there exists a unique $v_I\in V(F)$,
such that for any $v\in V(F)$, any element of $I$ is a descendant of $v$ if, and only if, $v\leq_F v_I$.
We denote by $F_I$ the subtree of $F$ whose vertices are the elements $v\in V(F)$ such that there exists $y\in I$, with $v_I\leq_F v\leq_F y$, and all the edges between these vertices.
\end{prop}

\begin{proof}
Let us denote by $V_I(F)$ the set of vertices $v$ of $F$ such that any element of $I$ is a descendant of $v$. As all the elements of $I$ belong to the same tree of $F$, the root of this tree belongs to $V_I(F)$, which is in consequence not empty. 
Let $x$, $y\in V_I(F)$. Let us choose $v\in I$ (which is not empty). Then $x\leq_F v$ and $y\leq_F v$, by definition of $V_I(F)$. As $F$ is a rooted forest, $x$ and $y$ are comparable for $\leq_F$. Therefore, $(V_I(F),\leq_F)$ is totally ordered, so as a minimal element. 
\end{proof}

\begin{defi}
Let $F$ be a reduced forest.
\begin{enumerate}
\item Let $\sim\in \calE[L(F)]$. We shall say that $\sim$ is tree-compatible if:
\begin{itemize}
\item Any $I\in L(F)/{\sim}$ is tree-compatible.
\item If $I$ and $J$ are two different elements of $X/{\sim}$, then $V(F_I)\cap V(F_J)=\emptyset$.
\end{itemize}
We denote by $\calE^{LC}(F)$ the set of tree-compatible equivalences on $L(F)$.
\item Let $\sim$ be a tree-compatible equivalence on $L(F)$. We put
\[F\mid_\ell \sim=\prod_{I\in L(F)/{\sim}} \red\left(F_I\right).\]
Note that $F\mid_\ell \sim$ is a reduced forest, with the same set of leaves as $F$.
\item Let $\sim$ be a tree-compatible equivalence on $L(F)$. We denote by $n_F(\sim)$ the number of equivalences $\sim'\in \calE^L(F)$ such that $\sim'\cap L(F)^2=\sim$.
\end{enumerate}\end{defi}

\begin{prop} \label{prop6.3.4}
Let $F$ be a reduced forest. Then
\[S_{\calF[\calR\calR\calF_\ell]}(F)=\sum_{\sim\in \calE^{LC}(F)} (-1)^{\cl(\sim)} n_F(\sim) F\mid_\ell \sim.\]
\end{prop}

\begin{proof}
As noticed before, as $\red(F)=F$,
\begin{align}
\label{EQ4.2.1}
S_{\calF[\calR\calR\calF_\ell]}(F)&=\red\circ S_{\calF[\calR\calF]}(F)=
\sum_{\sim\in \calE^L[\calF]}(-1)^{\cl(\sim)} \red(F\mid\sim).
\end{align}
Let $\sim\in \calE^L(F)$. We put $\iota(\sim)=\sim\cap L(F)^2$. By definition of $\calE^L(F)$, any $X\in V(F)/{\sim}$ is connected for $F$, so any $X\in L(F)/{\iota(\sim)}$ is tree-compatible. 
By definition of $F\mid\sim$, the elements of $L(F)/{\iota(\sim)}$ are the sets of leaves of the trees of $F\mid\sim$. Let $I\in L(F)/{\iota(\sim)}$, and let us denote by $T$ the subtree of $F\mid\sim$ such that $I=L(T)$. 
Then all the elements of $I$ belong to the same subtree of $F$, so to the same tree of $F$. Moreover, if $r_I$ is the root of $T$, then all the elements of $I$ are descendant of $r_I$, so $r_I\leq_F v_I$.
As a consequence, $F_I\subseteq T$, and, if $I$ and $J$ are two different elements of $\pi_\sim$, then $F_I$ and $F_J$ are included in two different trees of $F\mid\sim$, so they are disjoint. We proved that $\sim\in \calE^{LC}(F)$.
Moreover, we already noticed that $r_I\leq_F v_I$. As $V(T)$ is a convex subset of $F$, containing $r_I$ and all the elements of $I$, necessarily $V(F_I)\subseteq V(T)$. Let $x\in V(T)$, such that $r_I\leq_T x\leq_T v_I$. 
If $x$ has at least two children, then $I=L_{v_I}(T)\subsetneq L_x(T)$, which contradicts the fact that $L(T)=I$.
Consequently, $x$ has only one child. This implies that $\red(T)=\red(F_I)$. Denoting by $T_1,\ldots,T_k$ the trees of $F\mid \sim$, we finally obtain that
\[\red(F\mid\sim)=\red(T_1)\ldots \red(T_k)=\prod_{I\in \pi_\sim}\red(F_I)=F\mid_\ell \iota(\sim).\]
Note that $k=\cl(\sim)=|\iota(\sim)|$. Rewriting (\ref{EQ4.2.1}), we obtain
 \begin{align*}
S_{\calF[\calR\calR\calF_\ell]}(F)&=\sum_{\sim \in \calE^{LC}(F)} (-1)^{\cl(\sim)}|\{\sim'\in \calE^L(F)\mid \iota(\sim')=\sim\}| F\mid_\ell \sim\\
&=\sum_{\sim \in \calE^{LC}(F)}(-1)^{\cl(\sim)}n_F(\sim) F\mid_\ell \sim. \qedhere 
\end{align*} \end{proof}

\begin{remark} 
The formula of Proposition \ref{prop6.3.4} is multiplicity-free: 
If $\sim\neq \sim'$, then $F\mid_\ell \sim\neq F\mid_\ell \sim'$, as each tree of $F\mid_\ell \sim$ corresponds to an element of $L(F)/{\sim}$, whereas each tree of $F\mid_\ell \sim'$ corresponds to an element of $L(F)/{\sim'}$.
Moreover, it is not difficult to show that if $\sim\in \calE^{LC}(F)$, then $n_F(\sim)>0$.
\end{remark}

\begin{example}
\begin{align*}
S_{\calF[\calR\calR\calF_\ell]}\left(\hspace{-1mm}\begin{array}{c}\tdcinqsix{}{$3$}{}{$2$}{$1$}\end{array}\hspace{-2mm}\right)&=-\tdcinqsix{}{$3$}{}{$2$}{$1$}+\tdun{$1$}\tdtroisun{}{$3$}{$2$}+\tdun{$2$}\tdtroisun{}{$3$}{$1$}+2\tdun{$3$}\tdtroisun{}{$2$}{$1$}-4\tdun{$1$}\tdun{$2$}\tdun{$3$}.
\end{align*}
The coefficients in this expression are explained by the following array, where the corresponding elements of $\calE^L(F)$ are given by colors.
\begin{align*}
\begin{array}{|c|c|}
\hline\tdcinqsix{}{$3$}{}{$2$}{$1$}&\tdcinqsixA{}{$3$}{}{$2$}{$1$}\\
\hline\tdun{$1$}\tdtroisun{}{$3$}{$2$}&\tdcinqsixB{}{$3$}{}{$2$}{$1$}\\
\hline\tdun{$2$}\tdtroisun{}{$3$}{$1$}&\tdcinqsixC{}{$3$}{}{$2$}{$1$}\\
\hline\tdun{$3$}\tdtroisun{}{$2$}{$1$}&\tdcinqsixD{}{$3$}{}{$2$}{$1$},\tdcinqsixE{}{$3$}{}{$2$}{$1$}\\
\hline\tdun{$1$}\tdun{$2$}\tdun{$3$}&\tdcinqsixF{}{$3$}{}{$2$}{$1$},\tdcinqsixG{}{$3$}{}{$2$}{$1$},\tdcinqsixH{}{$3$}{}{$2$}{$1$},\tdcinqsixI{}{$3$}{}{$2$}{$1$}\\
\hline\end{array}\end{align*}\end{example}

\begin{remark}
It is not so difficult to show now that the antipode of $\calR\calR\calF_\ell$ is given by
\[S_{\calR\calR\calF_\ell}(F)=\sum_{\sim\in \calE^{LC}(F)} (-1)^{\cl(\sim)} n_F(\sim) F\mid_\ell \sim.\]
Let us give the steps of a proof:
\begin{enumerate}
\item Let us fix a set $\Omega$. Adapting the proof of the formula for the antipode of $\calF[\calR\calR\calF_\ell]$, we obtain a similar formula for the antipode of $\calF\circ \calF_\Omega[\calR\calR\calF_\ell]$. 
\item Observe that if $f:X\longrightarrow \Omega$ is injective, then the canonical map 
\[\left\{\begin{array}{rcl}
\calR\calR\calF_\ell\boxtimes \calD_\Omega[X]&\longrightarrow&\calF\circ \calF_\Omega[\calR\calR\calF_\ell]\\
(F,f)&\longmapsto&\overline{(F,f)}
\end{array}\right.\]
is injective. Therefore, we deduce the formula for the antipode of any element of $\calR\calR\calF_\ell\boxtimes \calD_\Omega[X]$, and then of any element of $\calR\calR\calF_\ell[X]$ by Proposition \ref{prop5.2.2}.
\end{enumerate}\end{remark}

\chapter{The non-degenerate and cocommutative case}

\label{chap7}

\section*{Introduction}

We now restrict ourselves to non-degenerate cocommutative CTBs, and show how the preceding results simplify.
For example, if $\calB$ is a cocommutative 	CTB, the polynomials invariants $\phi_\pm^\calB$ take their values in $\Z[T]$ and are unitary (and $\phi_-^\calB$ is in fact reduced to a single unitary monomial). We prove the converse when $\calB$ is double (Theorem \ref{theo7.1.2}).
A major simplification is the existence of a (partial) restriction map (Proposition \ref{prop7.1.1}).
Another one is the replacement of $\exch_+$ by several CTBs of hypergraphs, as we shall see.
Firstly, if $\calB$ is a non-degenerate cocommutative CTB, we show how to associate, to any convenient subspecies of $\calB$, a morphism of CTBs from $\calB$ to the CTB $\calH\calG$ of hypergraphs (Theorem \ref{theo7.2.1}), 
and show that any morphism of CTBs from $\calB$ to $\exch_+$ can be written as the composition of such a morphism with $\Gamma_+^{\calH\calG}$ (Corollary \ref{cor7.2.2}). Secondly, we use this result with two particular convenient subspecies of $\calB$, the subspecies of indecomposable elements (giving the morphism $\theta_{\CS}^\calB$),
and the subspecies of minimal indecomposable elements (giving the morphism $\theta_{\AC}^\calB$). 
Let us detail these two morphisms.

If $b\in \calB[X]$, $\theta_{\CS}^\calB(b)$ is the hypergraph on $X$ whose hyperedges are the subsets $Y\subseteq X$ such that $b_{\mid Y}$ indecomposable. Such a hypergraph is a connectedness system (Definition \ref{defi7.2.5}):
\begin{itemize}
\item If $Y\subseteq X$ is of cardinality one, then $Y$ is a hyperedge of $\theta_{\CS}^\calB(b)$.
\item If $Y$ and $Z$ are two intersecting hyperedges of $\theta_{\CS}^\calB(b)$, then $Y\cup Z$ is also a hyperedge of $\theta_{\CS}^\calB(b)$.
\end{itemize}
Connectedness systems define a sub-CTDB of $\calH\calG$ denoted by $\calC\calS$ (Theorem \ref{theo7.2.7}). We prove that $\theta_{\calC\calS}^\calB$ is the unique morphism of CTBs from $\calB$ to $\calC\calS$ that conserves the number of indecomposable components;
if $\calB$ is double, then it is the unique morphism of CTDBs from $\calB$ to $\calC\calS$ (Theorem \ref{theo7.2.8}). 
As a consequence, all the fundamental polynomial invariants of cocommutative CTDBs are in fact chromatic polynomials of connectedness systems. 
Unfortunately, if $\calB$ is a CTDB, the morphism $\Gamma_+^\calB$ from $\calB$ to $\exch_+$ does not always induce a morphism of CTDBs from $\calB$ to $\exch_+^{\double}$, and 
we prove the existence of a sub-CTDB $\calB^{\double}$ of $\calB$, maximal for the inclusion, such that $\Gamma_+^{\calB^{\double}}$ does, and that this construction is functorial (Propositions \ref{prop7.2.12} and \ref{prop7.2.13}).

If $b\in \calB[X]$, $\theta_{\AC}^\calB(b)$ is the hypergraph on $X$ whose hyperedges are the subsets $Y\subseteq X$ with $b_{\mid Y}$ indecomposable, and minimal for this property. Such a hypergraph is called an antichain cover (Definition \ref{defi7.2.14}):
\begin{itemize}
\item Any hyperedge of $\theta_{\AC}^\calB(b)$ is of cardinality $\geq 2$.
\item If $Y\subseteq Z$ are two hyperedges of $\theta_{\AC}^\calB(b)$, then $Y=Z$.
\end{itemize}
This defines a sub-CTB $\calA\calC$ of $\calH\calG$ (Theorem \ref{theo7.2.15}), unfortunately not a sub-CTDB.
We prove that $\Gamma_+^{\calA\calC}$ is injective, and that $\calA\calC$ is a maximal subspecies of $\calH\calG$ for this property (Theorem \ref{theo7.2.19}). 
This implies that $\Gamma_+^{\calA\calC}(\calA\calC)$ is the greatest cocommutative sub-CTB of $\exch_+$ (Corollary \ref{cor7.2.20}), and that any morphism from a cocommutative CTB to $\exch_+$ takes its values in $\Gamma_+^{\calA\calC}(\calA\calC)$.
As we mentioned before, $\calA\calC$ is not a sub-CTDB of $\calH\calG$. We prove the existence of a sub-CTB $\calA\calC^{\double}$, maximal for the inclusion, which can be made double, and that this double structure is unique (Proposition \ref{prop7.2.24}). Moreover, $\Gamma_+^{\calA\calC^{\double}}$ is a morphism of CTDBs from $\calA\calC^{\double}$ to $\exch_+^{\double}$.

We finally show how the fundamental polynomial invariant, composed with $\theta_{\calC\calS}^\calB$, can be used to give a formula for the antipode of $\calB$, even if $\calB$ is not double (Proposition \ref{prop7.2.29}). 
The chapter ends with a result on the alternating sign condition in the non-degenerate and cocommutative condition
(Theorem \ref{thoe7.2.30}). 
This implies that all the fundamental polynomial invariants from a cocommutative CTB with the alternating sign condition are in fact the ones with no "hole" in their coefficients, and are equal to the chromatic polynomial of some graph.

\section{The restriction map}

\subsection{Construction}

Let $\calB$ be a cocommutative CTB.

\begin{prop}\label{prop7.1.1}
Let $\calB$ be a non-degenerate cocommutative CTB. 
\begin{enumerate}
\item Let $X$ be a finite set, and $b\in \calB[X]$. Then $\phi^\calB_+(b)\in \Z[T]$ is unitary, of degree $|X|$, and $\phi^\calB_-(b)=T^{|X|}$, and $\calC_-(b)=\comp[X]$.
\item For any finite sets $Y\subseteq X$, there exists a restriction map
\[\left\{\begin{array}{rcl}
\calB[X]&\longrightarrow&\calB[Y]\\
b&\longmapsto&b_{\mid Y},
\end{array}\right.\]
such that for any finite set $X$, any $b\in \calB[X]$, and any $(X_1,\ldots,X_k)\in \comp[X]$,
\[\Delta_{X_1,\ldots,X_k}(b)=b_{\mid_1 X_1}\otimes\cdots\otimes b_{\mid_k X_k}.\]
\end{enumerate}\end{prop}

\begin{proof}
1. From Proposition \ref{prop3.5.12}, with $f=\epsilon^\calB_+$, $\phi^\calB_+(b)\in \Z[T]$.
By Proposition \ref{prop3.5.8}, $f_{\pp}^\calB(b)\in \Z$. As it belongs to $[0;1]$, it is equal to $0$ or to $1$.
As $\calB$ is non-degenerate, it is equal to 1, so $\phi^\calB_+(b)$ is unitary, of degree $|X|$. 
By lemma \ref{lem3.5.16}, $\phi^\calB_-(b)=T^{|X|}$, and $\calC_-(b)=\comp[X]$.

2. Let $Y\subseteq X$ be two finite sets, and $b\in \calB[X]$. If $Y=\emptyset$, we put $b_{\mid Y}=1_\calB$. If $X=Y$, we put $b_{\mid X}=b$. Otherwise, $(Y,X\setminus Y)\in \calC_-(b)$, and we put
\[\Delta_{Y,X\setminus Y}(b)=b_{\mid Y}\otimes b_{\mid_R X\setminus Y}.\]
As a consequence, in any case, $\Delta_{Y,X\setminus Y}(b)=b_{\mid Y}\otimes b_{\mid_RX\setminus Y}$.

Let $(X_1,\ldots,X_k)\in \comp[X]$, and $b\in \calB[X]$. Let us prove that 
\[\Delta_{X_1,\ldots,X_k}(b)=b_{\mid X_1}\otimes\cdots \otimes b_{\mid X_k},\]
If $k=1$, $\Delta_X(b)=b=b_{\mid X}$. Let us assume that $k\geq 2$. As $\calB$ is cocommutative, if $i\in [k]$,
\begin{align*}
&b_{\mid_i X_i}\otimes b_{\mid_1 X_1}\otimes\cdots\otimes b_{\mid_{i-1}X_{i-1}}\otimes b_{i\mid_{+1}X_{i+1}}\otimes\cdots\otimes b_{\mid_k X_k}\\
&=\Delta_{X_i,X_1,\ldots,X_{i-1},X_{i+1},\ldots,X_k}(b)\\
&=(\id_{\calB[X_i]}\otimes \Delta_{X_1,\ldots,X_{i-1},X_{i+1},\ldots,X_k})\circ \Delta_{X_i,X\setminus X_i}(b)\\
&=b_{\mid X_i}\otimes \Delta_{X_1,\ldots,X_{i-1},X_{i+1},\ldots,X_k}(b_{\mid_R X\setminus X_i}).
\end{align*}
So $b_i=b_{\mid _i X_i}$, for any $i\in [k]$.
\end{proof}

\begin{remark}\label{rk7.1.1}
Let $\calB$ be a non-degenerate cocommutative CTB. The coassociativity of $\Delta$ implies that if $Z\subseteq Y\subseteq X$ are finite sets, and if $b\in \calB[X]$, then
\[\left(b_{\mid Y}\right)_{\mid Z}=b_{\mid Z}.\]
The compatibility of $\Delta$ and $m$ implies that, if $X$, $Y$ are disjoint finite sets, and $Z\subseteq X\sqcup Y$, for any $a\in \calB[X]$, and $b\in \calB[Y]$,
\[(ab)_{\mid Z}=a_{\mid X\cap Z}b_{\mid Y\cap Z}.\]
If $\calB$ is a a non-degenerate cocommutative CTDB, the compatibility between $\delta$ and $\Delta$ implies that, for any finite set $X$, any $b\in \calB[X]$, and any $\sim\in \calE[b]$, 
if $X/{\sim}=(X'/{\sim}')\sqcup (X''/{\sim}'')$, then $\sim' \in \calE[b_{\mid X'}']$, $\sim''\in \calE[b_{\mid X''}]$, and
\begin{align*}
\left(b_{\mid X'}\right)/{\sim}'&=\left(b/{\sim}\right)_{\mid X'/{\sim}'},&\left(b_{\mid X''}\right)/{\sim}''&=\left(b/{\sim}\right)_{\mid X''/{\sim}''},&b\mid\sim&=\left(b_{\mid X'}\mid \sim'\right)\left(b_{\mid X''}\mid \sim''\right).
\end{align*}
\end{remark}

\begin{remark}
We proved that if $\calB$ is a non-degenerate cocommutative CTB, then, for any finite set $X$, and any $b\in \calB[X]$, $\phi^\calB_-(b)=\phi^\calB_{\pp}(b)=T^{|X|}$.
However, there exist non-degenerate twisted bialgebras, with the same $\phi^\calB_-$, and not cocommutative. Take for example $(\calH\calG,m,\Delta^{\subset,\cap})$ (which is not cocommutative).
For any finite set $X$, for any $G\in \calH\calG[X]$, $\phi^{\calH\calG}_-(G)=T^{|X|}$, as $\calC_-(G)=\comp[X]$. 
\end{remark}

\subsection{Polynomial invariants}

\begin{theo} \label{theo7.1.2}
Let $\calB$ be a non-degenerate CTDB. The following assertions are equivalent:
\begin{enumerate}
\item $\calB$ is cocommutative.
\item For any finite set $X$, for any $b\in \calB[X]$, $\phi^\calB_+(b)\in \Z[T]$.
\item For any finite set $X$, for any $b\in \calB[X]$, $\phi^\calB_-(b)\in \Z[T]$.
\item For any finite set $X$, for any $b\in \calB[X]$, $\phi^\calB_+(b)$ is unitary.
\item For any finite set $X$, for any $b\in \calB[X]$, $\phi^\calB_-(b)$ is unitary.
\item For any finite set $X$, for any $b\in \calB[X]$, $f^\calB_{\pp}(b)=1$.
\item For any finite set $X$, for any $b\in \calB[X]$, $\phi^\calB_-(b)=T^{|X|}$.
\end{enumerate}\end{theo}

\begin{proof}
By Lemma \ref{lem3.5.16}, 
\[3.\Longleftrightarrow 4.\Longleftrightarrow 5.\Longleftrightarrow 6.\Longleftrightarrow 7.\]
We already noticed that $1.\Longrightarrow 2.$, see Remark \ref{rk3.5.5}.

$2.\Longrightarrow 6$. Let $X$ be a finite set, and $b\in \calB[X]$. As $\calB$ is non-degenerate, $f^\calB_{\pp}(b)\neq 0$. Moreover, it belongs to $[0;1]$. As $f^\calB_{\pp}(b)$ is the leading coefficient of $\phi^\calB_+(b)$, it is an integer, so it is equal to $1$.

$6.\Longrightarrow 1$. \textit{First step}. Let us prove that for any finite set $X$, for any $b\in \calB[X]$, $\calC_-(b)=\comp[X]$. By hypothesis, $f^\calB_{\pp}(b)=1$, so $\comp_{|X|}[X]\subseteq \calC_-(b)$.
Consequently, $\calC_-(b)$ contains all the maximal elements of $\comp[X]$. As $\calC_-(b)\in \exch_-[X]$, $\calC_-(b)=\comp[X]$. 

\textit{Second step}. Let $X$, $Y$ be disjoint finite sets, and $b\in \calB[X\sqcup Y]$. By the first step, $(X$, $Y)$ and $(Y,X)$ belong to $\calC_-(b)$. We put
\begin{align*}
\Delta_{X,Y}(b)&=b_{\mid_LX}\otimes b_{\mid_RY},&\Delta_{Y,X}(b)&=b_{\mid_LY}\otimes b_{\mid_RX}.
\end{align*}
Let us put $\sim_X=\sim_{b_{\mid_LX}}^{(M)}$, and $\sim_Y=\sim_{b_{\mid_RY}}^{(M)}$, and $\sim=\sim_X\sqcup \sim_Y$. Then
\begin{align*}
(\Delta_{X/{\sim}_X,Y/{\sim}_Y}\otimes \id_{\calB[X\sqcup Y]})\circ \delta_\sim(b)&=m_{1,3,24}\circ (\delta_{\sim_X}\otimes \delta_{\sim_Y})\circ\Delta_{X,Y}(b)\neq 0,
\end{align*}
so $\sim\in \calE[b]$. Moreover, by the first step, $(Y/{\sim}_Y,X/{\sim}_X)\in \calE[b/{\sim}]$, so 
\begin{align*}
m_{1,3,24}\circ (\delta_{\sim_Y}\otimes \delta_{\sim_X})\circ\Delta_{Y,X}(b)&=(\Delta_{Y/{\sim}_Y,X/{\sim}_X}\otimes \id_{\calB[X\sqcup Y]})\circ \delta_\sim(b)\neq 0.
\end{align*}
We conclude that $\sim_X\in \calE[b_{\mid_RX}]$ and $\sim_Y\in \calE[b_{\mid_LY}]$. As a consequence, $\sim_X=\sim_{b_{\mid_LX}}^{(M)}\subseteq \sim_{b_{\mid_RX}}^{(M)}$,
and $\sim_Y=\sim_{b_{\mid_RY}}^{(M)}\subseteq \sim_{b_{\mid_LY}}^{(M)}$. By symmetry in $X$ and $Y$,
\begin{align*}
\sim_X&=\sim_{b_{\mid_LX}}^{(M)}=\sim_{b_{\mid_RX}}^{(M)},&\sim_Y&=\sim_{b_{\mid_RY}}^{(M)}=\sim_{b_{\mid_LY}}^{(M)}.
\end{align*}
Consequently,
\begin{align*}
\delta_{\sim_X}(b_{\mid_LX})&=b'_X\otimes b_{\mid_LX},&\delta_{\sim_X}(b_{\mid_RX})&=b''_X\otimes b_{\mid_RX},\\
\delta_{\sim_Y}(b_{\mid_LY})&=b'_Y\otimes b_{\mid_LY},&\delta_{\sim_Y}(b_{\mid_RY})&=b''_Y\otimes b_{\mid_RY},
\end{align*}
for certain elements $b'_X,b''_X\in \calB[X]$ and $b'_Y,b''_Y\in \calB[Y]$. We obtain
\begin{align*}
(\Delta_{X/{\sim}_X,Y/{\sim}_Y}\otimes \id_{\calB[X\sqcup Y]})\circ \delta_\sim(b)&=(b/{\sim})_{\mid_LX/{\sim}_X}\otimes (b/{\sim})_{\mid_RY/{\sim}_Y}\otimes b\mid\sim\\
=m_{1,3,24}\circ (\delta_{\sim_X}\otimes \delta_{\sim_Y})\circ\Delta_{X,Y}(b)&=b'_X\otimes b''_Y\otimes (b_{\mid_LX})(b_{\mid_RY}),
\end{align*}
so $b\mid\sim= (b_{\mid_LX})(b_{\mid_RY})$. By symmetry in $X$ and $Y$, $b\mid\sim=(b_{\mid_LY})(b_{\mid_RX})$.By unicity of the decomposition into indecomposable elements, $b_{\mid_LX}=b_{\mid_RX}$ 
and $b_{\mid_RY}=b_{\mid_LY}$. Therefore, 
\[\Delta_{Y,X}(b)=\tau_{\calB[X],\calB[Y]}\circ \Delta_{X,Y}(b).\]
So $\calB$ is cocommutative. 
\end{proof}

Here is an example of a non cocommutative CTB, such that for any finite set $X$, and any $b\in \calB[X]$, $\phi^\calB_+(b)\in \Z[T]$. Of course, it is not double.

\begin{example}
For any finite set $X=\{x,y,z\}$ of cardinality 3, we put
\[I(x,y,z)=\comp_3[X]+\{(xy,z)+(xz,y)\}.\]
Note that this depends on the chosen order on $x$, $y$, $z$. If $X$ is a set of cardinality $2$, we put
\[J(X)=\comp_2[X].\]
If $X$ is a set of cardinality 1, we put
\[K(X)=\comp[X].\]
All these elements belong to $\exch_+$. We denote by $\calB$ the twisted subalgebra of $\exch$ generated by these elements. If $X=\{x_1,\ldots,x_k\}$ is a finite set, 
\[\comp[X]=K(\{x_1\})\squplus\cdots\squplus K(\{x_k\})\in \calB[X].\]
Let $X=\{x,y,z\}$ be a finite set of cardinality 3. Then
\begin{align*}
\Delta_{\{x\},\{y,z\}}(I(x,y,z))&=K(\{x\})\otimes J(\{y,z\}),&\Delta_{\{y,z\},\{x\}}(I(x,y,z))&=J[\{y,z\}]\otimes K(\{x\}),\\
\Delta_{\{y\},\{x,z\}}(I(x,y,z))&=K(\{y\})\otimes J(\{x,z\}),&\Delta_{\{x,z\},\{y\}}(I(x,y,z))&=\comp[\{x,z\}]\otimes K(\{y\}),\\
\Delta_{\{z\},\{x,y\}}(I(x,y,z))&=K(\{z\})\otimes J(\{x,y\}),&\Delta_{\{x,y\},\{z\}}(I(x,y,z))&=\comp[\{x,y\}]\otimes K(\{z\}).
\end{align*}
Let $X=\{x,y\}$ is a set of cardinality 2. Then
\begin{align*}
\Delta_{\{x\},\{y\}}(J(\{x,y\}))&=K(\{x\})\otimes K(\{y\}).
\end{align*}
As a consequence, $\calB$ is a non cocommutative sub-CTB of $\exch$. The map $\Gamma^\calB_+$ is obviously the canonical injection. Therefore, for any finite set $X$, 
\begin{align*}
&\forall I\in \calB[X],&\phi_+^\calB(I)&=\sum_{k=0}^\infty |I\cap\comp_k[X]|\binom{T}{k}.
\end{align*}
Therefore, if $X=\{x,y,z\}$ is a set of cardinality 3, 
\[\phi_+^\calB(I(x,y,z))=6\binom{T}{3}+2\binom{T}{2}=T(T-1)^2\in \Z[T].\]
If $X$ is a set of cardinality 2,
\[\phi_+^\calB(J(X))=2\binom{T}{2}=T(T-1)\in \Z[T].\]
And if $X$ is a set of cardinality 1, $\phi^\calB_+(K(X))=T$. By multiplicativity, for any finite set $X$, for any $b\in \calB[X]$, $\phi_+^\calB(b)\in \Z[T]$.
\end{example}

\section{Morphisms to hypergraphs}

\subsection{Construction} 

\begin{theo}\label{theo7.2.1}
Let $\calB$ be a non-degenerate cocommutative CTB. 
\begin{enumerate}
\item Let $\calC$ be a convenient subspecies of $\calB$ (see Definition \ref{defi3.2.4}). For any finite set $X$, any $b\in \calB[X]$, we define the hypergraph $\theta^\calB_\calC(b)$ by
\[\HE\left(\theta^\calB_\calC(b)\right)=\{Y\in \calS_{\geq 1}[X]\mid b_{\mid Y}\in \calC[Y]\}.\]
Then $\theta^\calB_\calC:\calB\longrightarrow \calH\calG$ is a morphism of CTBs. Moreover, $\epsilon_\calC=\epsilon^{\calH\calG}_+\circ \theta^\calB_\calC$ is a combinatorial character of $\calB$, such that for any finite set $X$, and any $b\in \calB[X]$,
\[\epsilon_\calC(b)=\begin{cases}
1\mbox{ if }\forall Y\in \calS_{\geq 2}[X],\: b_{\mid Y}\notin \calC[Y],\\
0\mbox{ otherwise}.
\end{cases}\]
\item Let $\theta:\calB\longrightarrow\calH\calG$ be a morphism of CTBs such that for any finite set $X$, and any $b\in \calB[X]$, $\theta(b)\neq 0$.
There exists a unique convenient subspecies $\calC$ of $\calB$, such that $\theta=\theta^\calB_\calC$.
\end{enumerate}\end{theo}

\begin{proof}
1. Let $(X$, $Y)$ be a pair of disjoint finite sets, $a\in \calB[X]$, $b\in \calB[Y]$. Let us prove that $\theta^\calB_\calC(ab)=\theta^\calB_\calC(a)\theta^\calB_\calC(b)$. Let $Z$ be a hyperedge of $\theta^\calB_\calC(ab)$. 
Then $(ab)_{\mid Z}=a_{\mid X\cap Z}b_{\mid Y\cap Z} \in \calC[Z]$, so is indecomposable. Therefore, $X\cap Z=\emptyset$ or $Y\cap Z=\emptyset$.
\begin{itemize}
\item If $X\cap Z=\emptyset$, then $Z\subseteq Y$ and $(ab)_{\mid Z}=b_{\mid Z}\in \calC[Z]$, so $Z$ is a hyperedge of $\theta^\calB_\calC(b)$, and consequently a hyperedge of $\theta^\calB_\calC(a)\theta^\calB_\calC(b)$.
\item If $Y\cap Z=\emptyset$, similarly $Z$ is a hyperedge of $\theta^\calB_\calC(a)\theta^\calB_\calC(b)$.
\end{itemize}
Let $Z$ be a hyperedge of $\theta^\calB_\calC(a)\theta^\calB_\calC(b)$. Therefore, $Z$ is a hyperedge of $\theta^\calB_\calC(a)$, or of $\theta^\calB_\calC(b)$.
\begin{itemize}
\item If $Z$ is a hyperedge of $\theta^\calB_\calC(a)$, then $Z\subseteq X$ and $(ab)_{\mid Z}=a_{\mid Z}\in \calC[Z]$, so $Z$ is a hyperedge of $\theta^\calB_\calC(ab)$.
\item If $Z$ is a hyperedge of $\theta^\calB_\calC(b)$, similarly $Z$ is a hyperedge of $\theta^\calB_\calC(ab)$.
\end{itemize}
So $\theta^\calB_\calC(ab)=\theta^\calB_\calC(a)\theta^\calB_\calC(b)$.

Let $Y\subseteq X$ be two finite sets, and let $b\in \calB[X]$. Let us prove that $\theta^\calB_\calC(b)_{\mid Y}=\theta^\calB_\calC(b_{\mid Y})$, which will imply that $\theta^\calB_\calC$ is a morphism of twisted coalgebras. Let $Z\subseteq Y$, nonempty.
\begin{align*}
Z\in \HE\left(\theta^\calB_\calC(b)_{\mid Y}\right)&\Longleftrightarrow Z\in \HE\left(\theta^\calB_\calC(b)\right)\\
&\Longleftrightarrow b_{\mid Z}\in \calC[Z]\\
&\Longleftrightarrow \left(b_{\mid Y}\right)_{\mid Z}\in \calC[Z]\\
&\Longleftrightarrow Z\in \HE\left(\theta^\calB_\calC(b_{\mid Y})\right).
\end{align*}
So $\theta^\calB_\calC(b)_{\mid Y}=\theta^\calB_\calC(b_{\mid Y})$.

By composition, $\epsilon^{\calH\calG}_+\circ\theta^\calB_\calC$ is a combinatorial character of $\calB$. For any finite set $X$, and any $b\in \calB[X]$, $\epsilon^{\calH\calG}_+\circ\theta^\calB_\calC(b)=1$ if $\theta^\calB_\calC(b)$ has no nontrivial hyperedge, and $0$ otherwise:
this gives the announced formula for $\epsilon_\calC(b)$.

2. \textit{Existence}. Let $\theta:\calB\longrightarrow \calH\calG$ be a morphism of CTBs such that for any finite set $X$, and any $b\in \calB[X]$, $\theta(b)\neq 0$. For any finite set $X$, 
we consider the set $\calC[X]$ of elements of $\calB[X]$ such that $X$ is a hyperedge of $\theta(b)$.
If $X$ is a finite set and $b\in \calC[X]$, then $\theta(b)$ is a connected hypergraph, as $X$ is a hyperedge of $\theta(b)$. By Proposition \ref{prop2.2.6},
$1=\ic(\theta(b))\geq \ic(b)$, so $\ic(b)=1$: $b$ is indecomposable. Consequently, $\calC$ is a convenient subspecies of $\calB$.

Let $X$ be a finite set, and $b\in \calB[X]$. Let us prove that $\theta(b)=\theta^\calB_\calC(b)$ by induction on $|X|$. If $|X|=0$, then $b=1_\calB$, and $\theta(b)=\theta^\calB_\calC(b)=1_{\calH\calG}$. Let us assume the result at all ranks $<|X|$.
Let $Y$ be a hyperedge of $\theta(b)$. If $Y=X$, then $b\in \calC[X]$, so $b_{\mid X}=b\in \calC[X]$, and $X$ is a hyperedge of $\theta^\calB_\calC(b)$. Otherwise, as both $\theta^\calB_\calC$ and $\theta$ are morphism of twisted coalgebras, and applying the induction hypothesis to $Y$,
$Y$ is a hyperedge of 
\[\theta(b)_{\mid Y}=\theta(b_{\mid Y})=\theta^\calB_\calC(b_{\mid Y})=\theta^\calB_\calC(b)_{\mid Y},\]
so is a hyperedge of $\theta^\calB_\calC(b)$. 
Let $Y$ be a hyperedge of $\theta^\calB_\calC(b)$. If $Y=X$, then $b=b_{\mid X}\in \calC[X]$. By definition, $X$ is a hyperedge of $\theta(b)$. Otherwise, $Y$ is a hyperedge of 
\[\theta^\calB_\calC(b)_{\mid Y}=\theta^\calB_\calC(b_{\mid Y})=\theta(b_{\mid Y})=\theta(b)_{\mid Y},\]
so is a hyperedge of $\theta(b)$. As a conclusion, $\theta=\theta^\calB_\calC$.

2. \textit{Unicity}. Let $\calC,\calD$ be two convenient subspecies of $\calB$, such that $\theta^\calB_\calC=\theta^\calB_\calD$. Let $X$ be a finite set. If $|X|=0$, $\calC[\emptyset]=\calD[\emptyset]=\emptyset$. Let us now assume that $|X|\geq 1$. 
Let $b\in \calC[X]$. Then $b_{\mid X}=b\in \calC[X]$. By definition of $\theta^\calB_\calC$, $X$ is a hyperedge of $\theta^\calB_\calC(b)=\theta^\calB_\calD(b)$. By definition of $\calC_D$, $b_{\mid X}=b\in \calD[X]$. We proved that $\calC[X]\subseteq \calD[X]$. By symmetry, $\calC[X]=\calD[X]$.
\end{proof}

\begin{remark}
There may exist morphisms $\phi:\calB\longrightarrow \calH\calG$ of CTBs, taking the value $0$ on $\calB$. To obtain them, we can apply Proposition \ref{prop2.1.8}:
for a well-chosen biideal $\calI$ of $\calB$, we obtain a morphism $\phi':\calB\setminus\calI\longrightarrow \calH\calG$, which does not take the value $0$ on $\calB\setminus\calI$. 
\end{remark}

\begin{cor}\label{cor7.2.2}
Let $\calB$ be a non-degenerate cocommutative CTB, and $\calC$ a convenient subspecies of $\calB$. Then $\Gamma^{\calH\calG}_+\circ\theta^\calB_\calC=\Gamma^\calB_{\epsilon_\calC}$.
\end{cor}

\begin{proof}
By composition, $\Gamma^{\calH\calG}_+\circ\theta^\calB_\calC$ is a morphism of CTBs from $\calB$ to $\exch$. Moreover,
\[\epsilon_+\circ \Gamma^{\calH\calG}_+\circ\theta^\calB_\calC=\epsilon^{\calH\calG}_+\circ \theta^\calB_\calC=\epsilon_\calC.\]
By Theorem \ref{theo3.2.13} (unicity), $\Gamma^{\calH\calG}_+\circ\theta^\calB_\calC=\Gamma^\calB_{\epsilon_\calC}$.
\end{proof}

\begin{example}
Let $\calC_0$ and $\calC_1$ defined as follows: for any finite set $X$,
\begin{align*}
\calC_1[X]&=\begin{cases}
\calB[X]\mbox{ if }|X|=1,\\
\emptyset\mbox{ otherwise},
\end{cases}&
\calC_0[X]=\emptyset.
\end{align*}
Both $\calC_0$ and $\calC_1$ are convenient. For any finite set $X$, and any $b\in \calB[X]$, 
\begin{align*}
\HE\left(\theta^\calB_{\calC_1}(b)\right)&=\{\{x\}\mid x\in X\},&\HE\left(\theta^\calB_{\calC_0}(b)\right)&=\emptyset.
\end{align*}\end{example}

\begin{example}\label{ex7.2.2}
For any finite set $X$, let $\calC_{\CS}[X]$ be the set of indecomposable elements in $\calB[X]$. Then $\calC_{\CS}$ is convenient. 
For any finite set $X$, and any $b\in \calB[X]$, $\theta^\calB_{\calC_{\CS}}(b)$ is the hypergraph on $X$ whose hyperedges are the subsets $Y$ of $X$ such that $b_{\mid Y}$ is indecomposable.
This defines a morphism $\theta^\calB_{\CS}=\theta^\calB_{\calC_{\CS}}:\calB\longrightarrow\calH\calG$ which will be studied in Theorem \ref{theo7.2.8}.
\end{example}

\begin{example}\label{ex7.2.3}
For any finite set $X$ of cardinality $\geq 2$, let $\calC_{\AC}[X]$ be the set of indecomposable elements of $\calB[X]$, such that if $Y\in \calS_{\geq 2}[X]$, different from $X$, $b_{\mid Y}$ is not indecomposable (minimal indecomposable elements).
If $|X|<2$, we put $\calC_{\AC}[X]=\emptyset$. Then $\calC_{\AC}$ is convenient.
This defines a morphism $\theta^\calB_{\AC}=\theta^\calB_{\calC_{\AC}}:\calB\longrightarrow\calH\calG$ which will be studied in Theorem \ref{theo7.2.16}.
\begin{itemize}
\item For $\calB=\calG$, $\calC_{\AC}[X]$ is empty if $|X|\neq 2$, and reduced to the graph $G$ whose only edge is $X$ otherwise.
\item For $\calB=\calH\calG$, for any finite set $X$ of cardinality $\geq 2$, $\calC_{\AC}[X]$ is the set of hypergraphs $G\in \calH\calG[X]$ whose unique non-trivial hyperedge is $X$.
\item For $\calB=\calP$, $\calC_{\AC}[X]$ is empty if $|X|\geq 2$, and reduced to two elements if $X=\{x,y\}$ is of cardinality 2.
These elements are given by the two total orders
\begin{align*}
x<y,&&y<x.
\end{align*}\end{itemize}\end{example}

\begin{example}
For any finite set $X$, let $\calC[X]$ be the set of hypergraphs $G\in \calH\calG[X]$ such that $X\in \HE(G)$. Such a hypergraph is connected, so is indecomposable: $\calC$ is convenient. Then $\theta^\calB_\calC:\calH\calG\longrightarrow\calH\calG$ is $\id_{\calH\calG}$.
\end{example}

Different convenient subspecies can give the same $\epsilon_\calC$, as shown in the next Proposition:

\begin{prop}\label{prop7.2.3}
Let $\calB$ be a non-degenerate cocommutative CTB, and let $\calC$ be a convenient subspecies. 
\begin{enumerate}
\item $\epsilon_\calC=\epsilon^\calB_+$ if, and only if, $\calC_{\AC}\subseteq \calC$.
\item $\epsilon_\calC=\epsilon^\calB_-$ if, and only if, for any finite set $X$ of cardinality $\geq 2$, $\calC[X]=\emptyset$.
If this holds, $\theta_\calC^\calB$ takes its values in $\bfD\calH\calG$. 
\end{enumerate} 
\end{prop}

\begin{proof}
1. $\Longrightarrow$. Let $X$ be a finite set, and $b\in \calC_{\AC}[X]$. Then $b$ is indecomposable and $|X|\geq 2$, so $\epsilon^\calB(b)=\epsilon_\calC(b)=0$.
Consequently, there exists $Y\in \calS_{\geq 2}[X]$, such that $b_{\mid Y}\in \calC[Y]$. So $b_{\mid Y}$ is indecomposable. By definition of $\calC_{\AC}$, $X=Y$, so $b_{\mid Y}=b\in \calC[X]$.

1. $\Longleftarrow$. Let $X$ be a finite set, and $b\in \calB[X]$.
If $\epsilon_\calC(b)=0$, then there exists $Y\subseteq X$, of cardinality $\geq 2$, such that $b_{\mid Y}\in \calC[Y]$. 
Therefore, $b_{\mid Y}\notin \bfD\calB[Y]$. As $\bfD\calB$ is a twisted subcoalgebra of $\calB$, $b\notin \bfD\calB[X]$,
so $\epsilon^\calB_+(b)=0$.

Let us now assume that $\epsilon^\calB_+(b)=0$. There exists $Z\in \calS_{\geq 2}[X]$, such that $b_{\mid Z}$ is indecomposable (take for example an indecomposable component of $b$). 
Let us choose such a $Z$, of minimal cardinality. Therefore, if $Y\subseteq Z$, such that $b_{\mid Z}$ is indecomposable, then necessarily $|Z|=1$ or $Z=Y$. We conclude that $b_{\mid Z}\in \calC_{\AC}[Z]$, so $\epsilon_\calC(b)=0$.

As both $\epsilon^\calB_+$ and $\epsilon_\calC$ are combinatorial characters of $\calB$, they are equal.

2. $\Longrightarrow$. Let $X$ be a finite set of cardinality $\geq 2$, and $b\in \calB[X]$. Then $\epsilon_+^\calB(b)=\epsilon_\calC(b)=1$. As a consequence, $b_{\mid X}=b\notin \calC[X]$, so $\calC[X]=\emptyset$.

2. $\Longleftarrow$. Let $X$ be a finite set, and $b\in \calB[X]$. If $Y\subseteq X$, of cardinality $\geq 2$, then $b_{\mid Y}\notin \calC[Y]=\emptyset$, so $\epsilon_\calC(b)=1=\epsilon_-^\calB(b)$.

Let us assume that $\epsilon_\calC=\epsilon_-^\calB$. Then $\calC$ is reduced to its components of degree 1. 
Therefore, for any finite set $X$, and any $b\in \calB[X]$, $\HE\left(\theta_\calC^\calB(b)\right)$ only contains singletons.
Consequently, it belongs to $\bfD\calH\calG[X]$. 
\end{proof}

\subsection{More results on convenient subspecies}

\begin{prop}\label{prop7.2.4}\begin{enumerate}
\item Let $\calB$ be a non-degenerate cocommutative CTB.
\begin{enumerate}
\item Let $\calC$ be a convenient subspecies of $\calB$. For any finite set $X$, we denote by $\calB_\calC[X]$ the set of elements $b\in \calB[X]$ such that for any $Y\subseteq X$, $b_{\mid Y}\notin \calC[Y]$.
Then $\calB_\calC$ is a sub-CTDB of $\calB$.
\item Let $\calG$ be a sub-CTDB of $\calB$. There exists a convenient subspecies $\calC$ of $\calB$, such that $\calG=\calB_\calC$.
\end{enumerate}
\item Let $\calB$ be a non-degenerate cocommutative CTDB.
\begin{enumerate}
\item Let $\calC$ be a convenient subspecies of $\calB$. For any finite set $X$, we put
\[\calB_{\delta\calC}[X]=\{b\in \calB[X]\mid \forall \sim\in \calE[b],\:\forall Y\subseteq X/{\sim},\:(b/{\sim})_{\mid Y}\notin \calC[Y]\}.\]
Then $\calB_{\delta\calC}$ is a sub-CTDB of $\calB$. Moreover, $\calB_{\delta\calC}\subseteq \calB_\calC$, with the equality if, and only if, $\calB_\calC$ is a sub-CTDB of $\calB$. 
\item Let $\calG$ be a sub-CTDB of $\calB$. There exists a convenient species $\calC$ of $\calB$, such that $\calG=\calB_{\delta\calC}$.
\end{enumerate}\end{enumerate}\end{prop}

\begin{proof}
1. $(a)$ By definition, $1_\calB\in \calB_\calC[\emptyset]$. Let $X$, $Y$ be disjoint finite sets, $b\in \calB_\calC[X]$, and $b'\in \calB_\calC[Y]$. Let $Z\subseteq X\sqcup Y$. Then $(bb')_{\mid Z}=(b_{\mid X\cap Z})(b'_{\mid Y\cap Z})$. 
If $X\cap Z\neq \emptyset$ and $Y\cap Z\neq \emptyset$, then $(bb')_{\mid Z}$ is decomposable, so do not belong to $\calC[Z]$, as the elements of $\calC[Z]$ are indecomposable. 
If $Z\subseteq X$, then $(bb')_{\mid Z}=b_{\mid Z}\notin \calC[Z]$, as $b\in \calB_\calC[X]$. 
Similarly, if $Z\subseteq Y$, $(bb')_{\mid Z}=b'_{\mid Z}\notin\calC[Z]$. As a conclusion, $bb'\in \calB_\calC[X\sqcup Y]$.

Let $Y\subseteq X$ be finite sets, and $b\in \calB_\calC[X]$. For any $Z\in \calS_{\geq 2}[Y]$, $\geq 2$, $\left(b_{\mid Y}\right)_{\mid Z}=b_{\mid Z}\notin\calC[Z]$, so $b_{\mid Y}\in \calB_\calC[Y]$. 
This finally proves that $\calB_\calC$ is a twisted subbialgebra of $\calB$.


1. $(b)$ Let $\calG$ be a twisted subbialgebra of $\calB$. Let $X$ be a finite set. Let $\calC[X]$ be the set of indecomposable elements of $\calB[X]$ which do not belong to $\calG[X]$.
Obviously, $\calC$ is convenient.

$\calG[X]\subseteq \calB_\calC[X]$. Let $b\in \calG[X]$, such that $b\notin \calB_\calC[X]$. There exists $Y\subseteq X$, such that $b_{\mid Y}\in \calC[Y]$.
As $\calB$ is non-degenerate, $(Y,X\setminus Y)\in \calC_-(b)$, so, as $\calG$ is a twisted subbialgebra of $\calB$, $b_{\mid Y}\in \calG[Y]$. This is a contradiction: by hypothesis, the elements of $\calC[Y]$ do not belong to $\calG[Y]$. This is a contradiction, so $\calG[X]\subseteq \calB_\calC[X]$.

Let us prove that for any finite set $X$, $\calB_\calC[X]\subseteq \calG[X]$, by induction on $|X|$. If $X=\emptyset$, then $\calG[\emptyset]=\calB_\calC[\emptyset]=\{1_\calB\}$. Let us assume the result at all ranks $<|X|$.
Let us decompose $b$ into indecomposable elements: 
\[b=b_1\cdots b_k,\]
with $b_i\in \calB[X_i]$, indecomposable, for any $i\in [k]$. As $\calB_\calC$ is a twisted subbialgebra, for any $i\in [k]$, $b_i=b_{\mid X_i}\in \calB_\calC[X_i]$.
\begin{itemize}
\item If $k\geq 2$, then, by the induction hypothesis applied to $b_i$, for any $i\in [k]$, $b_i\in \calG[X_i]$. As $\calG$ is a twisted subbialgebra, $b=b_1\cdots b_k\in \calG[X]$.
\item If $k=1$, then $b$ is indecomposable. By definition of $\calB_\calC$, $b=b_{\mid X}$ does not belong to $\calC[X]$, so, by definition of $\calC$, it belongs to $\calG[X]$.
\end{itemize}
As a conclusion, $\calB_\calC[X]=\calG[X]$.

2. $(a)$. Let $X$, $Y$ be disjoint finite sets, $a\in \calB_{\delta\calC}[X]$ and $b\in \calB_{\delta_C}[Y]$.
Let $\sim\in \calE[ab]$, and $Z\subseteq (X\sqcup Y)/{\sim}$. By compatibility between $\delta$ and the product of $\calB$, $\sim=\sim_X\sqcup \sim_Y$, with $\sim_X=\sim\cap X^2$ and $\sim_Y=\sim\cap Y^2$, and $ab/{\sim}=(a/{\sim}_X)(b/{\sim}_Y)$, from Remark \ref{rk7.1.1}. 
\begin{itemize}
\item If $Z\cap X/{\sim}_X\neq \emptyset$ and $Z\cap Y/{\sim}_Y\neq \emptyset$, then 
\[(ab/{\sim})_{\mid Z}=(a/{\sim}_X)_{\mid Z\cap X/{\sim}_X}(b/{\sim}_Y)_{\mid Z\cap Y/{\sim}_Y},\]
so is not indecomposable: it does not belong to $\calC[Z]$.
\item If $Z\subseteq X/{\sim}_X$, then 
\[(ab/{\sim})_{\mid Z}=(a/{\sim}_X)_{\mid Z}.\]
As $a\in \calB_{\delta\calC}[X]$, this does not belong to $\calC[Z]$.
\item If $Z\subseteq Y/{\sim}_Y$, similarly $(ab/{\sim})_{\mid Z}\notin \calC[Z]$.
\end{itemize}
Therefore, $ab\in \calB_{\delta\calC}[X\sqcup Y]$: $\calB_{\delta\calC}$ is a twisted subalgebra of $\calB$.

Let $X$, $Y$ be disjoint finite sets, and $b\in \calB_{\delta\calC}[X\sqcup Y]$.
Let $\sim_X\in \calE[b_{\mid X}]$ and $Z\subseteq X/{\sim}_X$. We choose any $\sim_Y\in \calE[Y]$. By the compatibility between $\Delta$ and $\delta$, $\sim=\sim_X\sqcup \sim_Y\in \calE[b]$, and, from remark \ref{rk7.1.1},
\[(b_{\mid X}/{\sim}_X)_{\mid Z}=\left(b/{\sim}_{\mid X/{\sim}_X}\right)_{\mid Z}=(b/{\sim})_{\mid Z}\notin \calC[Z],\]
as $b\in \calB_{\delta_\calC}[X]$. Consequently, $b_{\mid X}\in \calB_{\delta\calC}[X]$. Similarly, $b_{\mid Y}\in \calB_{\delta\calC}[Y]$, so $\calB_{\delta\calC}$ is a twisted subcoalgebra of $\calB$.

Let $X$ be a finite set, $b\in \calB[X]$, and $\sim\in \calE[b]$. The two first items imply that $b\mid\sim\in \calB_{\delta\calC}[X]$. Let $\overline{\sim}'\in \calE[b/{\sim}]$ and $Z\subseteq X/{\sim}$. We denote by $\sim'\in \calE[X]$ the unique equivalence containing $\sim$ and inducing $\overline{\sim}'$ on $X/{\sim}$. 
By the coassociativity of $\delta$, $\sim'\in \calE[b]$, and 
\[((b/{\sim})/{\overline{\sim}'})_{\mid Z}=(b/{\sim}')_{\mid Z}\notin \calC[Z],\]
as $b\in \calB_{\delta\calC}[X]$. Therefore, $b/{\sim}\in \calB_{\delta\calC}[X/{\sim}]$, and $\calB_{\delta\calC}$ is a twisted double subbialgebra of $\calB$.

Let $X$ be a finite set, and $b\in \calB_{\delta\calC}[X]$. Let $Y\subseteq X$. For $\sim=\sim_b^{(m)}$,
\[b_{\mid Y}=(b/{\sim})_{\mid Y}\notin \calC[Y],\]
so $b\in \calB_\calC[X]$. Therefore, $\calB_{\delta\calC}\subseteq \calB_\calC$.

If $\calB_{\calC}=\calB_{\delta\calC}$, then it is a double subbialgebra of $\calB$. Let us assume that $\calB_{\calC}$ is a twisted double subbialgebra of $\calB$. 
Let $X$ be a finite set, and $b\in \calB_\calC[X]$. Let $\sim\in \calE[b]$ and $Y\subseteq X/{\sim}$.
As $\calB_\calC$ is a twisted double subbialgebra of $\calB$, $b/{\sim}\in \calB_\calC[X/{\sim}]$,
so $(b/{\sim})_{\mid Y}\notin \calC[Y]$. We obtain that $\calB_\calC=\calB_{\delta\calC}$.

2. $(b)$. Let $\calG$ be a twisted double subbialgebra of $\calB$. By 1. $(b)$, there exists a convenient subspecies $\calC$ such that $\calG=\calB_\calC$. As this is a twisted double subbialgebra of $\calB$, by 2. $(a)$, 
$\calG=\calB_{\delta\calC}$.
\end{proof}

\begin{remark}\label{rk7.2.2}
It may happen that two different $\calC$ give the same $\calB_\calC$. For example, if for any finite set $X$ of cardinality $1$, $\calC[X]=\calB[X]$, then $\calB_\calC$ is reduced to the unit of $\calB$. 
In order to avoid that two different $\calC$ give the same twisted subbialgebra $\calB_\calC$, we can restrict ourselves to sets $\calC$ with the following minimality condition: if $b\in \calC[X]$ and $Y\subsetneq X$, then $b_{\mid Y}\notin \calC[Y]$.
For example:
\begin{itemize}
\item For the CTB $\calG$ of graphs (Example \ref{ex2.1.2}), $\bfD\calG$ is given by the convenient species defined as follows: for any finite set $X$, 
\[\calC[X]=\begin{cases}
\{G_X\}\mbox{ if }|X|=2,\\
\emptyset\mbox{ otherwise},
\end{cases}\]
where $G_X\in \calG[X]$ is defined by $E(G)=\{X\}$.
 \item For the CTB $\calH\calG$ of hypergraphs (Example \ref{ex2.1.4}), $\bfD\calH\calG$ is given by the convenient species defined as follows: for any finite set $X$, 
\[\calC[X]=\begin{cases}
\{G\in \calH\calG[X]\mid \{X\}\subseteq \HE(G)\subseteq \{X\}\sqcup \calS_1[X]\}\mbox{ if }|X|\geq 2,\\
\emptyset\mbox{ otherwise}.
\end{cases}\]
In other words, the elements of $\calC[X]$ are the hypergraphs on $X$ with $X$ for single non-trivial hyperedge.
\end{itemize}\end{remark}

\begin{example}\begin{enumerate}
\item If $\calC[X]=\emptyset$ for any finite set $X$, then $\calB_\calC=\calB$.
\item Let $\calC[X]$ be the set of indecomposable elements of $\calB[X]$ if $|X|\geq 2$, and $\emptyset$ otherwise. Then $\calB_\calC=\bfD\calB$. 
\item If for any set $X$ if cardinality 1, $\calC[X]=\calB[X]$, then, for any finite set $X$,
\[\calB_\calC[X]=\begin{cases}
\{1_\calB\}\mbox{ if }X=\emptyset,\\
\emptyset\mbox{ otherwise}.
\end{cases}\]
\end{enumerate}\end{example}

\subsection{The map $\theta^\calB_{\CS}$ and connectedness systems}

Let us give the definition of a connectedness system \cite[Entry \href{https://oeis.org/A072446}{A072446}]{Sloane}.

\begin{defi}\label{defi7.2.5}
Let $G$ be a hypergraph. We shall say that $G$ is a connectedness system if: 
\begin{enumerate}
\item For any $x\in V(G)$, $\{x\}\in \HE(G)$.
\item If $A,B\in \HE(G)$ and $A\cap B\neq\emptyset$, then $A\cup B\in \HE(G)$. 
\end{enumerate}
The subspecies of $\calH\calG$ of connectedness systems is denoted by $\calC\calS$. 
\end{defi}

\begin{remark}\label{rk7.2.3} A connected connectedness system is called a clutter \cite{Wiseman}.
The following array gives the cardinality of $\calC\calS[X]$ and $\calH\calG[X]$ according to the cardinality of $|X|$. We denote by $\calH\calG_+[X]$ the set of hypergraphs on $X$ which contains all the singletons $\{x\}$, for $x\in X$.
\[\begin{array}{|c||c|c|c|c|c|c|}
\hline |X|&1&2&3&4&5&6\\
\hline\hline |\calC\calS[X]|&1&2&12&420&254076&18689059680\\
\hline |\calH\calG_+[X]|&1&2&16&2048&67108864&144115188075855872\\
\hline |\calH\calG[X]|&2&8&128&32768&2147483648&9223372036854775808\\
\hline\end{array}\]
Note that if $|X|=n$,
\begin{align*}
 |\calH\calG[X]|&=2^{2^n-1},& |\calH\calG_+[X]|&=2^{2^n-n-1}.
\end{align*}\end{remark}

\begin{lemma}\label{lem7.2.6}
Let $G$ be a hypergraph on a finite set $X$. The following assertions are equivalent:
\begin{enumerate}
\item $G$ is a connectedness system.
\item For any $Y\subseteq X$, $G_{\mid Y}$ is connected if, and only if, $Y\in \HE(G)$.
\end{enumerate}\end{lemma}

\begin{proof}
$2.\Longrightarrow 1$. For any $x\in X$, $G_{\mid\{x\}}$ is connected, so $\{x\}\in \HE(G)$. 
Let $A,B\in \HE(G)$, such that $A\cap B\neq\emptyset$. Then $G_{\mid A}$ and $G_{\mid B}$ are connected, as $A$ and $B$ are hyperedges of $G$. 
As $A\cap B\neq\emptyset$, $G_{\mid A\cup B}$ is connected. Therefore, $A\cup B\in \HE(G)$.

$1.\Longrightarrow 2$. Let $Y\subseteq X$. Obviously, if $Y\in \HE(G)$, then $G_{\mid Y}$ is connected. Let us assume that $G_{\mid Y}$ is connected. As the singletons of $X$ are hyperedges of $X$, $Y$ contains hyperedges of $G$. 
Let us consider a hyperedge $Z\subseteq Y$ of $G$, maximal for the inclusion. Let us assume that $Z\subsetneq Y$. 
As $G_{\mid Y}$ is connected, there exists $Z'\in \HE(G)$, such that $\emptyset\subsetneq Z\cap Z'\subsetneq Z$ and $Z'\subseteq Y$.
Therefore, as $G$ is a connectedness system, $Z\subsetneq Z\cup Z'\subseteq Y$, and $Z\cup Z'\in \HE(G)$. This contradicts the maximality of $Z$. So $Z=Y\in \HE(G)$.\end{proof}

\begin{theo}\label{theo7.2.7}
$\calC\calS$ is a sub-CTDB of $(\calH\calG,m,\Delta^{\subset,\subset},\delta^\subset)$.
\end{theo}

\begin{proof}
Let $X$, $Y$ be finite sets, $G\in \calC\calS[X]$, and $H\in \calC\calS[Y]$. It is not difficult to show that $GH\in \calC\calS[X\sqcup Y]$. 
Let $X\subseteq Y$ be finite sets, and $G\in \calC\calS[X]$. It is not difficult to show that $G_{\mid Y}\in \calC\calS[Y]$. Therefore, $\calC\calS$ is a twisted subbialgebra of $\calH\calG$.
It remains to show that $\calC\calS$ is stable under $\delta$. Let $X$ be a finite set, $G\in \calC\calS[X]$, and $\sim\in \calE[G]$. Observe that
\[G\mid \sim=\prod_{Y\in X/{\sim}} G_{\mid Y}.\]
As we already know that $\calC\calS$ is a sub-CTDB, $G\mid\sim\in \calC\calS[X]$. Let us now prove that $G/{\sim}\in \calC\calS[X/{\sim}]$. We consider the canonical surjection
\[p_\sim:\left\{\begin{array}{rcl}
X&\longrightarrow&X/{\sim}\\
x&\longmapsto&\overline{x}.
\end{array}\right.\]
Let $Y',Z'\in \HE(G/{\sim})$, intersecting. We put $Y'=p_\sim(Y)$ and $Z'=p_\sim(Z)$, with $Y,Z\in \HE(G)$. As $Y'\cap Z'\neq \emptyset$, there exist $y\in Y$, $z\in Z$, such that $y\sim z$. 
As $\sim\in \calE[G]$, the hypergraph $G_{\mid \overline{y}}$ is connected.
As $G$ is a connectedness system, by Lemma \ref{lem7.2.6}, $\overline{y}\in \HE(G)$. As $\overline{y}\cap Y\neq \emptyset$, $\overline{y}\cup Y\in \HE(G)$. As $z\in \overline{y}\cap Z$, $\overline{y}\cup Y\cup Z\in \HE(G)$. 
Therefore, as $\overline{y}\subseteq Y'$,
\[p_\sim(\overline{y}\cup Y\cup Z)=Y'\cup Z'\in \HE(G/{\sim}).\]
So $G/{\sim}\in \calC\calS[X/{\sim}]$. As a consequence, $\calC\calS$ is a sub-CTDB of $\calH\calG$. 
\end{proof}

\begin{remark}
Let $X$ be a finite set, and $G\in \calC\calS[X]$. By Lemma \ref{lem7.2.6}, for any $\sim\in\calE[X]$,
\[\sim\in \calE[G]\Longleftrightarrow X/{\sim} \subseteq \HE(G).\]
\end{remark}

Recall that $\theta^\calB_{\CS}$ is defined in Example \ref{ex7.2.2}.

\begin{theo}\label{theo7.2.8}\begin{enumerate}
\item Let $(\calB,m,\Delta)$ be a non-degenerate cocommutative CTB. 
The morphism $\theta^\calB_{\CS}$ takes its values in $\calC\calS$. Moreover, it is the unique morphism of CTBs from $\calB$ to $\calC\calS$ such that for any finite set $X$, for any $b\in \calB[X]$,
\[\ic\left(\theta^\calB_{\CS}(b)\right)=\ic(b).\]
\item If $(\calB,m,\Delta,\delta)$ is a non-degenerate cocommutative CTDB, then $\theta^\calB_{\CS}$ is the unique morphism of CTDBs from $\calB$ to $\calC\calS$. 
\end{enumerate}\end{theo}

\begin{proof}
1. \textit{First step}. Let $X$ be a finite set, and $b\in \calB[X]$. Let us first prove that $\theta^\calB_{\CS}(b)$ is a connectedness system.
Let $Y,Z\in \HE\left(\theta^\calB_{\CS}(b)\right)$, such that $Y\cap Z\neq \emptyset$. Let us assume that $b_{\mid Y\cup Z}=b_1b_2$, with $b_i\in \calB[X_i]$ for $i\in [2]$, with $X\cup Y=Z_1\sqcup Z_2$. 
The compatibility between $\Delta$ and the product implies that
\[b_{\mid Y}={b_1}_{\mid Y\cap Z_1}{b_2}_{\mid Y\cap Z_2}.\]
As $b_{\mid Y}$ is indecomposable, $Y\cap Z_1=Y$ or $Y\cap Z_2=Y$, so $Y\subseteq Z_1$ or $Y\subseteq Z_2$. 
Similarly, $Z\subseteq Z_1$ or $Z\subseteq Z_2$. As $Y\cap Z\neq \emptyset$ and $Z_1\cap Z_2=\emptyset$, $Y,Z\subseteq Z_1$ or $Y,Z\subseteq Z_2$, which gives $Z_1=Y\cup Z$ or $Z_2=Y\cup Z$. 
We obtain that $b_{\mid Y\cup Z}$ is indecomposable, so $Y\cup Z\in \HE\left(\theta^\calB_{\CS}(b)\right)$: $\theta^\calB_{\CS}(b)\in \calC\calS[X]$.

\textit{Second step}. Let us prove that for any finite set $X$, and any $b\in \calB[X]$, $\ic\left(\theta^\calB_{\CS}(b)\right)=\ic(b)$. Let us decompose $b$ into indecomposable components:
\[b=b_1\cdots b_k,\]
with $b_i\in \calB[X_i]$, indecomposable. Note that $\ic(b)=k$. Then
\[\theta^\calB_{\CS}(b)=\theta^\calB_{\CS}(b_1)\cdots \theta^\calB_{\CS}(b_k).\]
By definition, $X_k\in \HE\left(\theta^\calB_{\CS}(b)\right)$, so $\theta^\calB_{\CS}(b)_{\mid X_k}$ is connected, hence, indecomposable.
We have then found the decomposition of $\theta^\calB_{\CS}(b)$ into indecomposable components, so $\ic\left(\theta^\calB_{\CS}(b)\right)=k=\ic(b)$.

\textit{Third step}. Let $\Psi:\calB\longrightarrow \calC\calS$ be a morphism of CTBs, such that for any finite set $X$, and any $b\in \calB[X]$, $\ic(\Psi(b))=\ic(b)$. Note that this implies that $\Psi(b)\neq 0$ for any $b$.
Let us prove that $\theta^\calB_{\CS}(b)=\Psi(b)$ by induction on $|X|$. If $X=\emptyset$, then $b=1_B$, and $\Psi(b)=1_{\calC\calS}=\theta^\calB_{\CS}(b)$. Let us assume the result at all ranks $<|X|$, with $|X|\geq 1$. We now have to prove that $\HE\left(\theta^\calB_{\CS}(b)\right)=\HE(\Psi(b))$.

Let $Y\in \HE\left(\theta^\calB_{\CS}(b)\right)$. If $Y=X$, then, by definition of $\theta^\calB_{\CS}(b)$, $b$ is indecomposable, so $\ic(b)=1$. Therefore, $\ic(\Psi(b))=1$. By Lemma \ref{lem7.2.6}, $X\in \HE\left(\Psi(b)\right)$. Otherwise, $Y\in \HE\left(\theta^\calB_{\CS}(b_{\mid Y})\right)$.
As $\Psi$ is a morphism of twisted coalgebras, and by the induction hypothesis applied to $b_{\mid Y}$,
\[Y\in \HE\left(\theta^\calB_{\CS}(b_{\mid Y})\right)=\HE(\Psi(b_{\mid Y}))=\HE(\Psi(b)_{\mid Y})\subseteq \HE(\Psi(b)).\]
So $\HE\left(\theta^\calB_{\CS}(b)\right)\subseteq \HE(\Psi(b))$.

Let $Y\in \HE(\Psi(b))$. Then $Y\in \HE(\Psi(b)_{\mid Y})=\HE(\psi(b_{\mid Y}))$. Moreover, $\Psi(b)_{\mid Y}$ is connected, as $Y$ is a hyperedge of $\Psi(b)$,
so $\Psi(b)_{\mid Y}$ is indecomposable. As $\Psi$ conserves $\ic$, $b_{\mid Y}$ is indecomposable. Therefore, $Y\in \HE\left(\theta^\calB_{\CS}(b)\right)$. As a conclusion, $\theta^\calB_{\CS}(b)=\Psi(b)$.

2. \textit{Fourth step}. Let us now assume that $\calB$ is a CTDB. Let $X$ be a finite set, $b\in \calB[X]$, and $\sim \in \calE[X]$. We determine $\calE[b]$ thanks to Lemma \ref{lem4.2.3}:
\begin{align*}
\sim\in \calE[b]&\Longleftrightarrow \forall Y\in X/{\sim},\: b_{\mid y}\mbox{ is indecomposable}\\
&\Longleftrightarrow \forall Y\in X/{\sim}, \: Y\in \HE\left(\theta^\calB_{\CS}(b)\right)\\
&\Longleftrightarrow \forall Y\in X/{\sim},\: \theta^\calB_{\CS}(b)_{\mid Y}\mbox{ is connected}\\
&\Longleftrightarrow \sim \in \calE^\subset\left[\theta^\calB_{\CS}(b)\right].
\end{align*}
We used Lemma \ref{lem7.2.6} for the third equivalence. So $\calE[b]=\calE^\subset\left[\theta^\calB_{\CS}(b)\right]$. Moreover, if $\sim\in \calE[b]$,
\[\theta^\calB_{\CS}(b\mid \sim)=\theta^\calB_{\CS}\left(\prod_{Y\in X/{\sim}} b_{\mid Y}\right)=\prod_{Y\in X/{\sim}} \theta^\calB_{\CS}(b_{\mid Y})=\prod_{Y\in X/{\sim}} \theta^\calB_{\CS}(b)_{\mid Y}=\theta^\calB_{\CS}(b)\mid \sim.\]

Let us now prove that $\theta^\calB_{\CS}(b/{\sim})=\theta^\calB_{\CS}(b)/{\sim}$, for $\sim\in \calE[b]$. Firstly, observe that if $(X'/{\sim}_{X'})\sqcup (X''/{\sim}_{X''})=X/{\sim}$, by the compatibility between $\delta$ and $\Delta$,
$(b/{\sim})_{\mid X'/{\sim}_{X'}}=\left(b_{\mid X'}\right)/{\sim}_{X'}$.
Let $Y'=Y/{\sim}_Y$ be a hyperedge of $\theta^\calB_{\CS}(b/{\sim})$. By the preceding observation,
$(b/{\sim})_{\mid Y/{\sim}_Y}=\left(b_{\mid Y}\right)/{\sim_Y}$ is indecomposable. Moreover, by the compatibility between $\Delta$ and $\delta$, $\sim_Y\in \calE[b_{\mid Y}]$. As $\calB$ is non-degenerate, it is regular (Proposition \ref{prop4.1.16}), so $\ic\left(\left(b_{\mid Y}\right)/{\sim}_Y\right)=\ic\left(b_{\mid Y}\right)=1$, as $Y/{\sim}$ is a hyperedge of $\theta^\calB_{\CS}(b/{\sim})$. Therefore, $Y\in \HE\left(\theta^\calB_{\CS}(b)\right)$, and $Y'=p_\sim(Y)\in \HE\left(\theta^\calB_{\CS}(b)/{\sim}\right)$. 

Let $Y'$ be a hyperedge of $\theta^\calB_{\CS}(b)/{\sim}$. There exists $Y\in \theta^\calB_{\CS}(b)$, such that $Y'=p_\sim(Y)$.
we put $Y=\{y_1,\ldots,y_k\}$. As $\sim\in \calE[b]$, for any $i\in [k]$, $b_{\mid \overline{y_i}}$ is indecomposable,
so $\overline{y_i}$ is a hyperedge of $\theta^\calB_{\CS}(b)$. As $b$ is a connectedness system:
\begin{itemize}
\item $Y\cap \overline{y_1}\neq \emptyset$, so $Y\cup \overline{y_1}$ is a hyperedge of $\theta^\calB_{\CS}(b)$.
\item $Y\cap \overline{y_2}\neq \emptyset$, so $Y\cup \overline{y_1}\cup \overline{y_2}$ is a hyperedge of $\theta^\calB_{\CS}(b)$.
\item $\vdots$
\item $Y\cup \overline{y_1}\cup\cdots \cup \overline{y_k}$ is a hyperedge of $\theta^\calB_{\CS}(b)$. 
\end{itemize}
As a consequence, putting $\tilde{Y}=Y\cup \overline{y_1}\cup\cdots \cup \overline{y_k}$, $\tilde{Y}$ is a hyperedge of $\theta^\calB_{\CS}(b)$, so $b_{\mid \tilde{Y}}$ is indecomposable. As $\calB$ is regular, $\left(b_{\mid \tilde{Y}}\right)/{\sim}_{\tilde{Y}}=(b/{\sim})_{\mid Y'}$ is indecomposable, so $Y'$ is a hyperedge of $\theta^\calB_{\CS}(b/{\sim})$. 

We proved that $\theta^\calB_{\CS}(b)/{\sim}=\theta^\calB_{\CS}(b/{\sim})$. We finally obtain that $\theta^\calB_{\CS}$ is compatible with $\delta$.

\textit{Final step}. It finally remains to prove that $\theta^\calB_{\CS}$ is the unique morphism of CTDBs
from $\calB$ to $\calC\calS$. If $\Psi$ is such a morphism, by Proposition \ref{prop4.3.4}, for any finite set $X$, and any $b\in \calB[X]$, $\ic(\Psi(b))=\Psi(b)$. As a consequence of the fourth step, $\Psi=\theta^\calB_{\CS}$.
\end{proof}

\begin{cor}
Let $(\calB,m,\Delta)$ be a non-degenerate cocommutative CTB. 
Then $\epsilon^{\calC\calS}_+\circ \theta^\calB_{\CS}=\epsilon^\calB_+$, $\Gamma^{\calC\calS}_+\circ \theta^\calB_{\CS}=\Gamma^\calB_+$, and $\phi^{\calC\calS}_+\circ\theta^\calB_{\CS}=\phi^\calB_+$.
\end{cor}

\begin{proof}
By Proposition \ref{prop7.2.3}, as $\calC_{\AC}\subseteq\calC_{\CS}$, $\epsilon^{\calC\calS}_+\circ \theta^\calB_{\CS}=\epsilon_{\calC_{\CS}}=\epsilon^\calB_+$. 
Consequently, by Theorem \ref{theo3.2.9}, $\Gamma^{\calC\calS}_+\circ \theta^\calB_{\CS}=\Gamma^\calB_+$, and finally $\phi^{\calC\calS}_+\circ\theta^\calB_{\CS}=\phi^\calB_+$.
\end{proof}

\begin{remark}
We obtain that for any $n\in \N$, the polynomials of $\calP_+(n)$ coming from a cocommutative CTB are the chromatic polynomials of connectedness systems with $n$ vertices.
\end{remark}

\begin{remark}
Example \ref{ex4.3.2}-2 shows that $\Gamma^{\calC\calS}_+$ is not a morphism of CTDBs taking its values in $\exch_+^{\double}$, as both $G$ and $G'$ are connectedness systems.
\end{remark}

If $\calB$ is degenerate, we can use Proposition \ref{prop4.1.16} to go back to the non-degenerate case:

\begin{cor}
Let $\calB$ be a cocommutative CTB. The following is a morphism of CTBs:
\[\theta^\calB_{\CS}:\left\{\begin{array}{rcl}
\calB&\longrightarrow&\calC\calS\\
b&\longmapsto&\begin{cases}
\theta^\calB_{\CS}(b)\mbox{ if }\calC_+(b)\neq \emptyset,\\
0\mbox{ otherwise}.
\end{cases}
\end{array}\right.\]
It satisfies $\epsilon^{\calC\calS}_+\circ \theta^\calB_{\CS}=\epsilon^\calB_+$. Moreover, if $\calB$ is a regular cocommutative CTDB, then $\theta^\calB_{\CS}$ is the unique morphism of CTDBs from $\calB$ to $\calC\calS$. 
\end{cor}

\begin{proof}
Let us assume that $\calB$ is a CTDB.
By composition of the canonical surjection $\calB\longrightarrow \nd\calB$ and of $\theta^{\nd\calB}_{\CS}:\nd\calB\longrightarrow \calC\calS$, we obtain that $\theta^\calB_{\CS}:\calB\longrightarrow\calC\calS$ is a morphism of CTDBs.
Let $\Psi:\calB\longrightarrow\calC\calS$ be a morphism of CTDBs. Let $X$ be a finite set, and $b\in \calB[X]$, such that $\calC_+(b)=\emptyset$. Let us prove that $\Psi(b)=0$. 
As $\Psi$ is a double bialgebra morphism and $\calB$ and $\calC\calS$ are both regular, 
\[\epsilon^{\calC\calS}_+\circ \Psi=\epsilon_\delta\circ \Psi=\epsilon_\delta=\epsilon^\calB_+.\]
So $\phi^{\calC\calS}_+\circ\Psi=\phi^\calB_+$. In particular, as $\calC_+(b)=\emptyset$, $\phi^\calB_+(b)=\phi^{\calC\calS}_+(\Psi(b))=0$. If $\Psi(b)\neq 0$, as $\calC\calS$ is non-degenerate, $\phi^{\calC\calS}_+(\psi(b))\neq 0$: this is a contradiction. So $\Psi(b)=0$. 

Consequently, $\Psi$ can be factorized as $\Psi=\Psi'\circ \varpi$, where $\varpi:\calB\longrightarrow \nd\calB$ is the canonical surjection.
As $\Psi$ is a morphism of CTDBs, so is $\Psi'$. As $\nd\calB$ is non-degenerate, $\Psi':\nd\calB\longrightarrow \calC\calS$ is $\theta^{\nd\calB}_{\CS}$, so $\Psi$ is $\theta^\calB_{\CS}$. 
\end{proof}

As a consequence:

\begin{prop} $\theta^{\calC\calS}_{\CS}=\id_{\calC\calS}$, and $\theta^{\calH\calG}_{\CS}\circ \theta^{\calH\calG}_{\CS}=\theta^{\calH\calG}_{\CS}$.
\end{prop}

\begin{proof}
Observe that $\id_{\calC\calS}$ is an endomorphism of the CTB $\calC\calS$ which conserves $\ic$. By Theorem \ref{theo7.2.8} (unicity part), it is $\theta^{\calC\calS}_{\CS}$.
As $\theta^{\calH\calG}_{\CS}$ takes its values in $\calC\calS$, $\theta^{\calH\calG}_{\CS}\circ \theta^{\calH\calG}_{\CS}=\theta^{\calC\calS}_{\CS}\circ \theta^{\calH\calG}_{\CS}=\theta^{\calH\calG}_{\CS}$.
\end{proof}

\begin{prop}\label{prop7.2.12}
Let $\calB$ be a non-degenerate cocommutative CTDB. There exists a sub-CTDB $\calB^{\double}$ of $\calB$, maximal for the inclusion, 
such that $\Gamma^\calB_+(\calB^{\double})\subseteq \exch_+^{\double}$ and the restriction of $\Gamma^\calB_+$ to $\calB^{\double}$ is a morphism of CTDBs.
\end{prop}

\begin{proof}
Let $X$ be a finite set. We put
\[\calC[X]=\{b\in \calB[X]\mid b\mbox{ is indecomposable, and } \ic(\calC_+(b))>1\}.\]
Then $\calC$ is a convenient subspecies of $\calB$. We put $\calB^{\double}=\calB_{\delta\calC}$. By Proposition \ref{prop7.2.4}, $\calB^{\double}$ is a sub-CTDB of $\calB$.
Let $X$ be a finite set, and $b\in \calB[X]$, indecomposable. For $\sim=\sim_b^{(m)}\in \calE[b]$ and $Y=X$, by definition of $\calB^{\double}$, $b=(b/{\sim})_{\mid Y}\notin \calC[X]$, so $\ic(\calC_+(b))=1$.
Let now $b\in \calB[X]$, not necessarily indecomposable. Let $b=b_1\cdots b_k$ be its decomposition into indecomposable elements, with $b_i\in \calB[X_i]$. As $\calB^{\double}$ is a twisted subbialgebra of $\calB$, for any $i\in [k]$, $b_i\in\calB^{\double}[X_i]$ and is indecomposable. 
Consequently, $\ic(\calC_+(b_i))=1$. Therefore,
\[\ic(\calC_+(b))=\ic(\calC_+(b_1)\squplus\cdots\squplus \calC_+(b_k))=\calC_+(b_1)+\cdots+\calC_+(b_k)=1+\cdots+1=k=\ic(b).\]
By Theorem \ref{theo4.3.5}, $\Gamma^\calB_+:\calB^{\double}\longrightarrow \exch_+^{\double}$ is a morphism of CTDBs.

Let $\calB'$ be a sub-CTDB of $\calB$, such that $\Gamma^{\calB'}_+:\calB'\longrightarrow \exch_+^{\double}$ is a morphism of CTDBs.
Let us prove that for any finite set $X$, $\calB'[X]\subseteq \calB^{\double}[X]$. Let $b\in \calB'[X]$, $\sim\in \calE[b]$ and $Y\subseteq X/{\sim}$. As $\calB'$ is a sub-CTDB of $\calB$, $(b/{\sim})_{\mid Y}\in \calB'[Y]$.
By Theorem \ref{theo4.3.5}, $\ic(\calC_+((b/{\sim})_{\mid Y}))=\ic((b/{\sim})_{\mid Y})$, so $(b/{\sim})_{\mid Y}\notin \calC[Y]$. Consequently, $b\in \calB_{\delta\calC}[X]$.
\end{proof}

\begin{prop}\label{prop7.2.13}
Let $\calB$ be a be a non-degenerate cocommutative CTDB. Then $\theta^\calB_{\CS}(\calB^{\double})\subseteq \calC\calS^{\double}$.
\end{prop}

\begin{proof}
Let $\calB'=\theta^\calB_{\CS}(\calB^{\double})$. As $\theta^\calB_{\CS}$ is a morphism of CTDBs, $\calB'$ is a sub-CTDB of $\calC\calS$. Let $X$ be a finite set, and $b\in \calB^{\double}[X]$. 
By Proposition \ref{prop4.3.4} for the morphisms $\theta^\calB_{\CS}$ and $\Gamma^\calB_+$,
\[\ic(b)=\ic\left(\theta^\calB_{\CS}(b)\right)=\ic\left(\Gamma^\calB_+(b)\right).\]
Moreover, $\Gamma^{\calC\calS}_+\circ \theta^\calB_{\CS}=\Gamma^\calB_+$. Therefore, for any $G=\theta^\calB_{\CS}(b)\in \theta^\calB_{\CS}(\calB^{\double}[X])$, $\ic\left(\Gamma^{\calC\calS}_+(G)\right)=\ic(G)$. 
By Theorem \ref{theo4.3.5}, $\Gamma^{\calB'}_+:\calB'\longrightarrow \exch_+^{\double}$ is a morphism of CTDBs. By definition of $\calC\calS^{\double}$, $\calB'\subseteq \calC\calS^{\double}$.
\end{proof}

\begin{remark}
As a consequence, the union of the $\Gamma^\calB_+(\calB^{\double})$'s, where $\calB$ is a non-degenerate cocommutative CTB, is also a cocommutative CTB, quotient of $\calC\calS^{\double}$ through $\theta^{\calC\calS^{\double}}_{\CS}$.
\end{remark}

\subsection{The map $\theta^\calB_{\AC}$ and antichain covers}

\begin{defi}\label{defi7.2.14}
Let $G$ be a hypergraph. We shall say that it is an antichain cover if:
\begin{itemize}
\item For any $Y\in \HE(G)$, $|Y|\geq 2$.
\item For any $Y\subseteq Z\in \HE(G)$, $Y=Z$. 
\end{itemize}
The subspecies of $\calH\calG$ of antichain covers is denoted by $\calA\calC$. 
\end{defi}

\begin{remark}
The number of antichain covers of a given order can be found in Entry \href{https://oeis.org/A006126}{A006126} of \cite{Sloane}. Let us complete the array of Remark \ref{rk7.2.3}.
\[\begin{array}{|c||c|c|c|c|c|c|}
\hline |X|&1&2&3&4&5&6\\
\hline\hline |\calA\calC[X]|&1&2&9&114&6894&7785062\\
\hline |\calC\calS[X]|&1&2&12&420&254076&18689059680\\
\hline |\calH\calG_+[X]|&1&2&16&2048&67108864&144115188075855872\\
\hline |\calH\calG[X]|&2&8&128&32768&2147483648&9223372036854775808\\
\hline\end{array}\]
\end{remark}

\begin{theo}\label{theo7.2.15}
$\calA\calC$ is a sub-CTB of $(\calH\calG,m,\Delta^{\subset,\subset})$.
\end{theo}

\begin{proof}
Let $X$, $Y$ be finite sets, $G\in \calA\calC[X]$, and $H\in \calA\calC[Y]$. It is not difficult to show that $GH\in \calA\calC[X\sqcup Y]$. 
Let $X\subseteq Y$ be finite sets, and $G\in \calA\calC[X]$. It is not difficult to show that $G_{\mid Y}\in \calA\calC[Y]$. Therefore, $\calA\calC$ is a twisted subbialgebra of $\calH\calG$.
\end{proof}

\begin{remark}\begin{enumerate}
\item Graphs are antichain covers, so $\calG\subseteq \calA\calG$, ans this inclusion is a morphism of CTBs. 
\item The subspecies $\calA\calC$ is not stable under $\delta$. Let $G\in \calH\calG[4]$, defined by
\[\HE(G)=\{\{1,4\}, \:\{2,4\}, \:\{1,2,3\}\}.\]
It is an antichain cover. Let $\sim\in \calE[4]$ whose classes are $\{1\}$, $\{3\}$ and $\{2,4\}$. Then $\sim\in \calE^\subset[G]$, and
\[\HE(G/{\sim})=\{\{\overline{1},\overline{2}\},\:\{\overline{1},\overline{2},\overline{3}\}\}.\]
So $G/{\sim}$ is not an antichain cover.
\end{enumerate}\end{remark}

\begin{theo}\label{theo7.2.16}
Let $(\calB,m,\Delta)$ be a non-degenerate cocommutative CTB. The morphism $\theta^\calB_{\AC}$ takes its values in $\calA\calC$. 
Moreover, $\epsilon^{\calA\calC}_+\circ \theta^\calB_{\AC}=\epsilon^\calB_+$, $\Gamma^{\calA\calC}_+\circ \theta^\calB_{\AC}=\Gamma^\calB_+$, and $\phi^{\calA\calC}_+\circ\theta^\calB_{\AC}=\phi^\calB_+$.
\end{theo}

\begin{proof}
Let $X$ be a finite set, and $b\in \calB[X]$. Let $Y\in \HE\left(\theta^\calB_{\AC}(b)\right)$.
Then $b_{\mid Y}\in \calC_{\AC}[Y]$, so $|Y|\geq 2$. Let $Y\subseteq Z\in \HE\left(\theta^\calB_{\AC}(b)\right)$. Then $b_{\mid Y}\in \calC_{\AC}[Y]$ and $b_{\mid Z}\in \calC_{\AC}[Z]$. Therefore, $b_{\mid Y}$ and $b_{\mid Z}$ are indecomposable.
Moreover, $(b_{\mid Z})_{\mid Y}=b_{\mid Z}$ so, by definition of $\calC_{\AC}$, $Y=Z$.
Therefore, $\theta^\calB_{\AC}(b)\in \calA\calC[X]$.

By Proposition \ref{prop7.2.3}, $\epsilon^{\calA\calC}_+\circ \theta^\calB_{\AC}=\epsilon_{\calC_{\AC}}=\epsilon^\calB_+$. Consequently, by Theorem \ref{theo3.2.9}, $\Gamma^{\calA\calC}_+\circ \theta^\calB_{\AC}=\Gamma^\calB_+$, and finally $\phi^{\calA\calC}_+\circ\theta^\calB_{\AC}=\phi^\calB_+$.
\end{proof}

\begin{remark}
We obtain that for any $n\in \N$, the polynomials of $\calP_+(n)$ coming from a cocommutative CTB are the chromatic polynomials of antichain covers with $n$ vertices.
\end{remark}

\begin{cor}
Let $(\calB,m,\Delta)$ be a cocommutative CTB. 
The following is a morphism of CTBs:
\[\theta^\calB_{\AC}:\left\{\begin{array}{rcl}
\calB&\longrightarrow&\calA\calC\\
b&\longmapsto&\begin{cases}
\theta^\calB_{\AC}(b)\mbox{ if }\calC_+(b)\neq \emptyset,\\
0\mbox{ otherwise}.
\end{cases}
\end{array}\right.\]
It satisfies $\epsilon^{\calA\calC}_+\circ \theta^\calB_{\AC}=\epsilon^\calB_+$. 
\end{cor}

\begin{prop}
$\theta^{\calA\calC}_{\AC}=\id_{\calA\calC}$, and $\theta^{\calC\calS}_{\AC}\circ\theta^{\calH\calG}_{\CS}=\theta^{\calH\calG}_{\AC}$.
\end{prop}

\begin{proof}
Let $X$ be a finite set, and $G\in \calA\calC[X]$. Let $Y$ be a hyperedge of $\theta^{\calA\calC}_{\AC}(G)$. By definition, $Y$ is a hyperedge of $G$. Let $Y$ be a hyperedge of $G$. As $G\in \calA\calC[X]$, $Y$ is minimal for the inclusion, so $Y$ is a hyperedge of $\theta^{\calA\calC}_{\AC}$. So $\theta^{\calA\calC}_{\AC}(G)=G$.

Let $Y$ be a hyperedge of $\theta^{\calH\calG}_{\AC}(G)$. Then it is a hyperedge of $G$, so $G_{\mid Y}$ is connected: $Y$ is a hyperedge of $\theta^{\calH\calG}(G)$. If $Z\subseteq Y$ is a nontrivial hyperedge of $\theta^{\calH\calG}(G)$, then $G_{\mid Z}$ is connected,
so contains a nontrivial hyperedge $Z'\subseteq Z$: therefore, $Z'\subseteq Y$. By minimality of $Y$, $Z'=Z=Y$. So $Y$ is a hyperedge of $\theta^{\calC\calS}_{\AC}\circ\theta^{\calH\calG}_{\CS}(G)$.

Let $Y$ be a hyperedge of $\theta^{\calC\calS}_{\AC}\circ\theta^{\calH\calG}_{\CS}(G)$. Then $Y$ is a hyperedge of $\theta^{\calH\calG}_{\CS}(G)$, so $G_{\mid Y}$ is connected. If $Y$ is not a hyperedge of $G$, then it contains a nontrivial hyperedge $Z\subsetneq Y$.
Therefore, $G_{\mid Z}$ is connected, so $Z$ is a hyperedge of $\theta^{\calH\calG}_{\CS}(G)$: this contradicts the minimality of $Y$. So $Y$ is a hyperedge of $G$. If it is not minimal for the inclusion, then there exists a non-trivial hyperedge $Z\subsetneq Y$.Then $G_{\mid Z}$ is connected,
so $Z$ is a hyperedge of $\theta^{\calH\calG}_{\CS}(G)$: this contradicts the minimality of $Y$. So $Y$ is a hyperedge of $\theta^{\calH\calG}_{\AC}(G)$. Therefore, $\theta^{\calC\calS}_{\AC}\circ\theta^{\calH\calG}_{\CS}(G)=\theta^{\calH\calG}_{\AC}(G)$.
\end{proof}

\begin{theo}\label{theo7.2.19}\begin{enumerate}
\item $\Gamma^{\calA\calC}_+$ is injective.
\item Let $\calH\calG'$ be a subspecies of $\calH\calG$ which contains $\calA\calC$, such that $\Gamma^{\calH\calG'}_+$ is injective. Then $\calH\calG'=\calA\calC$.
\item Let $\calB$ be a cocommutative CTB. Then $\Gamma^\calB_+(\calB)\subseteq \Gamma^{\calA\calC}_+(\calA\calC)\sqcup \{0\}$.
\item Let $(\calB,m,\Delta)$ be a non-degenerate cocommutative CTB. For any finite set $X$, and any $b\in \calB[X]$, $\ic(\calC_+(b))=\ic\left(\theta^\calB_{\AC}(b)\right)$.
\end{enumerate}\end{theo}

\begin{proof}
1. Let $X$ be a finite set, and $G\in \calA\calC[X]$. Let $Y\in \calS_{\geq 2}[X]$. Let us prove that $Y$ is a hyperedge of $G$ if, and only if, for any $\pi\in \calC_+(G)$, $Y$ is not a block of $\pi$.
This will imply that $\Gamma^{\calA\calC}_+$ is injective.
If $Y$ is a nontrivial hyperedge of $G$, then $G_{\mid Y}\notin \bfD\calH\calG[Y]$, so for any $\pi\in \calC_+(G)$, $Y$ is not a block of $\pi$.
Let us assume that for any $\pi\in \calC_+(G)$, $Y$ is not a block of $\pi$. We put $X\setminus Y=\{x_1,\ldots,x_k\}$. Then, by hypothesis, $(Y,\{x_1\},\ldots,\{x_k\})\notin \calC_+(G)$, so $G_{\mid Y}\notin \bfD\calH\calG[Y]$, as for any $i\in [k]$, obviously $G_{\mid \{x_i\}}\in \bfD\calH\calG[\{x_i\}]$.
So $Y$ contains a nontrivial edge $Z$ of $G$. As $G$ is an antichain cover, $Y=Z\in\HE(G)$.

2. Let $X$ be a finite set, and $G\in \calH\calG'[X]$. Then $\Gamma^{\calA\calC}_+\circ \theta^{\calH\calG'}_{\AC}(G)=\Gamma^{\calH\calG'}_+\circ \theta^{\calH\calG'}_{\AC}(G)=\Gamma^{\calH\calG'}_+(G)$. Both $G$ and $\theta^{\calH\calG'}_{\AC}(G)$ belong to $\calH\calG'[X]$, as $\calA\calC[X]\subseteq \calH\calG'[X]$. 
As $\Gamma^{\calH\calG'}_+$ is injective, $G=\theta^{\calH\calG'}_{\AC}(G)$, so $G\in \calA\calC[X]$.

3. Let $X$ be a finite set, and $b\in \calB[X]$. Let us assume that $\theta^\calB_{\AC}(b)\neq 0$. Then $\theta^\calB_{\AC}(b)\in \calA\calC[X]$, and $\Gamma^\calB_+(b)=\Gamma^{\calA\calC}_+\circ \theta^\calB_{\AC}(b)\in \Gamma^{\calA\calC}_+(\calA\calC[X])$.

4. As $\Gamma^\calB_+:\calA\calC\longrightarrow \exch_+$ is injective, by Proposition \ref{prop2.2.6}, 
\[\ic(\calC_+(b))=\ic\left(\Gamma^\calB_+(b)\right)=\ic\left(\Gamma^{\calA\calC}_+\circ\theta^\calB_{\AC}(b)\right)=\ic\left(\theta^\calB_{\AC}(b)\right). \qedhere\]
\end{proof}

\begin{remark}
It may happen that $\ic(b)<\ic\left(\theta^\calB_{\AC}(b)\right)$. Take for example $\calB=\calH\calG$ and $G\in \calH\calG[3]$, whose hyperedges are $\{1,2\}$ and $\{1,2,3\}$. Then $\ic(G)=1$. The hypergraph $\theta^{\calH\calG}_{\AC}(G)$ has only one nontrivial hyperedge, namely $\{1,2\}$, so $\ic\left(\theta^{\calH\calG}_{\AC}(G)\right)=2$. 
\end{remark}

\begin{remark}
As a consequence, the union of the $\Gamma^\calB_+(\calB)$'s, for all cocommutative CTBs $\calB$, is also a cocommutative CTB, isomorphic to $\calC\calA$ through $\Gamma^{\calA\calC}_+$.
\end{remark}

\begin{cor}\label{cor7.2.20}
The two maps $\Gamma^{\calA\calC}_+:\calA\calC\longrightarrow \cocom\exch_+$ and $\theta_{\AC}^{\cocom\exch_+}:\cocom\exch_+\longrightarrow \calA\calC$ are isomorphisms of CTBs, and are inverse one from the other.
\end{cor}

\begin{proof}
As $\calA\calC$ is a cocommutative CTB, $\Gamma^{\calA\calC}_+(\calA\calC)$ is a cocommutative sub-CTB of $\exch_+$, so, bis included in $\cocom\exch_+$.
Moreover, as $\Gamma^{\calA\calC}_+$ is injective (Theorem \ref{theo7.2.19}), it is isomorphic to $\calA\calC$ through $\Gamma^{\calA\calC}_+$. Let $\calB$ be a cocommutative sub-CDB of $\exch_+$. By restriction, \[\Gamma_+^\calB=\left(\Gamma^{\exch_+}_+\right)_{\mid \calB}=\id_\calB.\]
So $\Gamma^\calB_+(\calB)=\calB$. By Theorem \ref{theo7.2.19}-2, $\calB\subseteq \Gamma^{\calA\calC}_+(\calA\calC)$.
By maximality, $\Gamma^{\calA\calC}_+(\calA\calC)=\cocom\exch_+$. So $\Gamma^{\calA\calC}_+$ is an isomorphism from $\calA\calC$ to $\cocom\exch_+$.
Let us consider $\phi=\Gamma_+^{\calA\calC}\circ\theta_{\AC}^{\cocom\exch_+}$. 
\[\epsilon\circ \phi=\epsilon_+\circ \Gamma_+^{\calA\calC}\circ\theta_{\AC}^{\cocom\exch_+}=\epsilon_+^{\calA\calC}\circ \theta_{\AC}^{\cocom\exch_+}=\epsilon_+^{\cocom\exch_+}.\]
Therefore, by the terminal property of $\exch$, $\phi:\cocom\exch_+\longrightarrow \cocom\exch_+\subseteq\exch$ is the unique morphism of CTBs such that $\epsilon\circ\phi=\epsilon_+^{\cocom\exch_+}$: it is obviously the canonical injection, so $\phi=\id_{\cocom\exch_+}$. As $\Gamma^{\calA\calC}_+$ is an isomorphism, $\theta_{\AC}^{\cocom\exch_+}$ is its inverse.
\end{proof}

\begin{cor}\label{cor7.2.21}
Let $\calB$ be a non-degenerate cocommutative CTDB. For any finite set $X$, we put
\[\calC[X]=\{b\in \calB[X]\mid b\mbox{ indecomposable, and }\ic\left(\theta^\calB_{\AC}(b)\right)>1\}.\]
Then $\calB^{\double}=\calB_{\delta\calC}$.
\end{cor}

\begin{proof}
From Theorem \ref{theo7.2.19}, for any finite set $X$,
\[\calC[X]=\{b\in \calB[X]\mid b\mbox{ indecomposable, and }\ic(\calC_+(b))>1\}.\]
The proof of Proposition \ref{prop7.2.12} shows that $\calB^{\double}=\calB_{\delta\calC}$.
\end{proof}

Unfortunately, $\calA\calC$ is not a sub-CTDB of $\calH\calG$, but has a greatest double sub-object. To describe it, we need to redefine the contraction of an antichain cover.

\begin{defi}
Let $X$ be a finite set, $G\in \calA\calC[X]$, and $\sim\in \calE[X]$.
\begin{enumerate}
\item We shall say that $\sim\in \calE^\subset[G]$ if, for any $Y\in X/{\sim}$, $G_{\mid Y}$ is connected.
\item If so, the hypergraph $G:{\sim}\in \calH\calG[X/{\sim}]$ and the antichain cover $G\mid\sim\in \calA\calC[X]$ are defined by
\begin{align*}
\HE{G/{\sim}})&=\{p_\sim(Y)\mid Y\in \HE(G),\: |p_\sim(Y)|\geq 2\},\\
\HE(G\mid \sim)&=\{Y\in \HE(G)\mid \forall x,y\in Y,\: x\sim y\}.
\end{align*}\end{enumerate}\end{defi}

The unique difference with hypergraphs is the definition of $G/{\sim}$, where we do not allow trivial hyperedges anymore.

\begin{lemma}\label{lem7.2.23}
Let $\calB$ be a sub-CTDB of $\calA\calC$, admitting a coproduct $\delta'$ making it a CTDB. Then, for any finite set $X$, any $G\in \calB[X]$, and any $\sim\in \calE[X]$,
\[\delta'_\sim(G)=\begin{cases}
G/{\sim}\otimes G\mid\sim \mbox{ if }\sim\in \calE^\subset[G],\\
0\mbox{ otherwise}.
\end{cases}\]
\end{lemma}

\begin{proof}
As $\calA\calC$ is non-degenerate, $\calB$ is non-degenerate. Moreover, by Theorem \ref{theo7.2.19}, $\Gamma^\calB_+$ is injective. By Corollary \ref{cor4.3.2}, $\delta$ is unique, and is given by
\[\delta'=\left({\Gamma^\calB_+}^{-1}\otimes {\Gamma^\calB_+}^{-1}\right)\circ \delta\circ \Gamma^\calB_+,\]
where $\delta$ is the coproduct of $\exch_+^{\double}$. Let $X$ be a finite set, and $G\in \calB[X]$. We put
\[\delta'_\sim(G)=\begin{cases}
G/'{\sim}\otimes G\mid'\sim\mbox{ if }\sim\in \calE'[G],\\
0\mbox{ otherwise}.
\end{cases}\]
For any $\sim\in \calE[X]$, by Lemma \ref{lem4.2.3}, as $\calC_-(b)=\comp[X]$,
\[\sim\in\calE'[G]\Longleftrightarrow \forall Y\in X/{\sim},\: G_{\mid Y}\mbox{ is connected}\Longleftrightarrow \sim\in \calE^\subset[G].\]
So $\calE[G]=\calE^\subset[G]$. Still by Lemma \ref{lem4.2.3}, $G\mid'\sim=G\mid \sim$. It remains to consider $G/'\sim$. If is given by
\[\Gamma_+^\calB(G/'\sim)=\Gamma^\calB_+(G)/{\sim}.\]
By definition of an antichain cover, for any $H\in \calA\calC[Y]$,
\begin{align}
\label{EQ7.1}\Gamma^{\calA\calC}_+(H)=\{(Y_1,\ldots,Y_k)\in \comp[Y]\mid \forall i\in [k],\: Y_i\notin \HE(H)\}.
\end{align}
Moreover, as $m(\gamma^\calB_+(G))=\comp[X]$,
\begin{align}
\label{EQ7.2}
(X_1/{\sim_1},\ldots,X_k/{\sim_k})\in \Gamma^\calB_+(G)/{\sim}\Longleftrightarrow \forall i\in [k],\: \sim_i=\sim_{G_{\mid X_i}}^{\ic}.
\end{align}
Combining (\ref{EQ7.1}) and (\ref{EQ7.2}), we obtain that the hyperedges of $G/'\sim$ are the $p_\sim(Y)$ of cardinality $\geq 2$, where $Y\in \HE(G)$. 
\end{proof}

\begin{prop}\label{prop7.2.24}
For any finite set $X$, we put
\[\calA\calC^{\double}[X]=\{G\in \calA\calC[X]\mid \forall \sim\in \calE^\subset[G],\: G/{\sim}\in \calA\calC[X/{\sim}]\}.\]
Then $\calA\calC^{\double}$ is a CTDB with the coproduct $\delta'$ of Lemma \ref{lem7.2.23}. It is the greatest sub-CTB of $\calA\calC$ that can be made double.
Moreover, $\Gamma^{\calA\calC^{\double}}_+$ is a morphism of CTDBs from $\calA\calC^{\double}$ to $\exch_+^{\double}$.
\end{prop} 

\begin{proof}
Let $X$, $Y$ be disjoint finite sets. Let $G\in \calA\calC^{\double}[X]$ and $H\in \calA\calC^{\double}[Y]$. Let $\sim\in \calE^\subset[GH]$. Putting $\sim_X=\sim\cap X^2$ and $\sim_Y=\sim\cap Y^2$, then $\sim=\sim_X\sqcup \sim_Y$ and
\[GH/{\sim}=(G/{\sim_X})(G/{\sim_Y}).\]
As $G/{\sim_X}\in \calA\calC[X/{\sim_X}]$ and $H/{\sim_Y}\in \calA\calC[Y/{\sim_Y}]$, $(GH)/{\sim}\in \calA\calC[(X\sqcup Y)/{\sim}]$. So $GH\in \calA\calC^{\double}[X\sqcup Y]$. So $\calA\calC^{\double}$ is a twisted subalgebra of $\calA\calC$.

Let $G\in \calA\calC^{\double}[X\sqcup Y]$. Let $\sim_X\in \calE^\subset[G_{\mid X}]$. We choose an arbitrary $\sim_Y\in \calE^\subset[G_{\mid Y}]$. Then $\sim=\sim_X\sqcup\sim_Y\in \calE^\subset[G]$, and
\[(G_{\mid X})/{\sim_X}=(G/{\sim})_{\mid X/{\sim_X}}.\]
As $G\in \calA\calC^{\double}[X\sqcup Y]$, $G/{\sim}\in \calA\calC[(X\sqcup Y)/{\sim}]$, so $(G/{\sim})_{\mid X/{\sim_X}}\in \calA\calC[X/{\sim_X}]$. We obtain that $\calA\calC^{\double}$ is a twisted subcoalgebra of $\calA\calC$.

Let $X$ be a finite set, $G\in \calA\calC^{\double}[X]$, and $\sim\in \calE^\subset[G]$. By the two preceding observations, $G\mid\sim\in \calA\calC^{\double}[X]$. By hypothesis, $G/{\sim}\in \calA\calC[X/{\sim}]$.
Let $\overline{\sim}'\in \calE^\subset[G/{\sim}]$. Let $\sim'\in \calE[X]$ be the unique equivalence containing $\sim$, inducing $\overline{\sim}'$ on $X/{\sim}$. Then $\sim'\in \calE^\subset[G]$ and $(G/{\sim})/{\overline{\sim}'}=G/{\sim'}$. By hypothesis on $G$, this belongs to $\calA\calC[X/{\sim'}]$.
So $G/{\sim}\in \calA\calC^{\double}[X/{\sim}]$. We obtain that $\calA\calC^{\double}$ is stable under $\delta'$. It is not difficult to see that it is a CTDB.

Let $\calB\subseteq \calA\calC$ be a sub-CTDB, with a coproduct making it double. By Lemma \ref{lem7.2.23}, this coproduct is necessarily $\delta'$.
Let $X$ be a finite set, and $G\in \calB[X]$. Then, for any $\sim\in \calE[G]$, $G/{\sim}\in \calB[X/{\sim}]\subseteq \calA\calC[X/{\sim}]$, so $G\in \calA\calC^{\double}[X]$.
Consequently, $\calB\subseteq \calA\calC^{\double}$.

By Theorem \ref{theo7.2.19}, $\Gamma^{\calA\calC^{\double}}_+={\Gamma^{\calA\calC}_+}_{\mid \calA\calC^{\double}}$ is injective. By Corollary \ref{cor4.3.6}, $\Gamma^{\calA\calC^{\double}}_+:\calA\calC^{\double}\longrightarrow\exch_+^{\double}$ is a morphism of CTDBs.
\end{proof}

\begin{cor}\label{cor7.2.25}
$\Gamma^{\calA\calC^{\double}}_+(\calA\calC^{\double})$ is the greatest cocommutative subcoalgebra of $\exch_+^{\double}$. It is a sub-CTDB of $\exch_+^{\double}$.
\end{cor}

\begin{proof}
As $\Gamma^{\calA\calC^{\double}}_+$ is a morphism of CTDBs and $\calA\calC^{\double}$ is cocommutative, $\Gamma^{\calA\calC^{\double}}_+(\calA\calC^{\double})$ is a cocommutative sub-CTDB of $\exch_+^{\double}$.
Let us denote by $\calB$ the greatest combinatorial twisted cocommutative subcoalgebra of $\exch_+^{\double}$. By Proposition \ref{prop5.3.2}, it is a sub-CTDB of $\exch_+^{\double}$. 
By Corollary \ref{cor7.2.20}, it is included in $\Gamma^{\calA\calC}_+(\calA\calC)$. As $\Gamma^{\calA\calC}_+$ is injective, $\calB'={\Gamma^{\calA\calC}_+}^{-1}(\calB)$ is a sub-CTDB of $\calA\calC$, and can be given a coproduct $\delta'$ making it double.
By Proposition \ref{prop7.2.24}, $\calB'\subseteq \calA\calC^{\double}$, so $\calB\subseteq \Gamma^{\calA\calC^{\double}}_+(\calA\calC^{\double})$.
\end{proof}

\subsection{Results on $\calC\calS^{\double}$}

\begin{notation}
\begin{enumerate}
\item Let $X$ be a finite set of cardinality $\geq 2$. The hypergraph $G_X\in \calH\calG[X]$ is defined by
\[\HE(G_X)=\{X\}\sqcup\{\{x\}\mid x\in X\}.\]
\item Let $Y\subsetneq X$ be two finite sets, with $|Y|\geq 2$. The hypergraph $G_{X,Y}\in \calH\calG[X]$ is defined by
\[\HE(G_{X,Y})=\{X,Y\}\sqcup\{\{x\}\mid x\in X\}.\]
\end{enumerate}\end{notation}

\begin{remark}
Let $Y\subseteq X$ be finite sets, with $|Y|\geq 2$. Then $\theta_{\AC}(G_{X,Y})$ has a unique hyperedge, namely $Y$, and
\[\ic\left(\theta^{\calH\calG}_{\AC}(G_{X,Y})\right)=|X|-|Y|+1\geq 2,\]
whereas $\ic(G_{X,Y})=1$.
\end{remark}

\begin{lemma}\label{lem7.2.26}
Let $X$ be a finite set of cardinality $2$, and $G\in \calC\calS[X]$, connected. There exists $\sim\in \calE[G]$, such that $|X/{\sim}|\geq 2$ and $G/{\sim}=G_{X/{\sim}}$.
\end{lemma}

\begin{proof}
As $G$ is connected and belongs to $\calC\calS[X]$, by Lemma \ref{lem7.2.6}, $X\in \HE(G)$. We proceed by induction on $|X|$. If $|X|=2$, then $G=G_X$. Let us assume the result at all ranks $<|X|$.
If $X$ is the only non-trivial hyperedge of $G$, then $G=G_X$ and we take $\sim=\sim_G^{(m)}$.
Otherwise, let us choose another non-trivial hyperedge $Y\subseteq X$ of $G$, and let $\sim'\in \calE[X]$, whose classes are $Y$ and singletons. Then $\sim'\in \calE[G]$, and $2\leq |X/{\sim}|<|X|$.
The quotient $G/{\sim}'\in \calC\calS[X/{\sim}']$ is also connected. Therefore, there exists $\overline{\sim}\in \calE[G/{\sim}']$, such that $|(X/{\sim}')/{\overline{\sim}}|\geq 2$ and $(G/{\sim}')/{\overline{\sim}}=G_{(X/{\sim}')/{\overline{\sim}}}$. 
Let $\sim\in \calE[X]$ be the unique equivalence on $X$ containing $\sim'$ and inducing $\overline{\sim'}$ on $X/{\sim}$. By coassociativity of $\delta$, $\sim\in \calE[G]$, and $G/{\sim}=(G/{\sim}')/{\overline{\sim}}=G_{X/{\sim}}$.
\end{proof}

\begin{theo}\label{theo7.2.27}
For any finite set $X$, we put
\[\calC[X]=\{G_{Y,X}\mid Y\subsetneq X,\mid |Y|\geq 2\}.\]
Then 
\[\calC\calS^{\double}=\calC\calS_{\delta \calC}\subsetneq \theta^{\calH\calG}_{\CS}\circ\theta^{\calH\calG}_{\AC}(\calC\calS)=\theta^{\calH\calG}_{\CS}\circ\theta^{\calH\calG}_{\AC}(\calH\calG)\subsetneq \calC\calS_\calC.\]
\end{theo}

\begin{proof}
Obviously, $\theta^{\calH\calG}_{\CS}\circ\theta^{\calH\calG}_{\AC}(\calC\calS)\subseteq\theta^{\calH\calG}_{\CS}\circ\theta^{\calH\calG}_{\AC}(\calH\calG)$, as $\calC\calS\subseteq \calH\calG$. If $X$ is a finite set and $G\in \calH\calG[X]$, then 
\[\theta^{\calH\calG}_{\CS}\circ\theta^{\calH\calG}_{\AC}(G)=\theta^{\calH\calG}_{\CS}\circ \theta^{\calH\calG}_{\AC}\circ \theta^{\calH\calG}_{\CS}(G),\]
so, as $\theta^{\calH\calG}_{\CS}(G)\in \calC\calS[X]$, $\theta^{\calH\calG}_{\CS}\circ\theta^{\calH\calG}_{\AC}(G)\in \theta^{\calH\calG}_{\CS}\circ\theta^{\calH\calG}_{\AC}(\calC\calS)[X]$.

For any finite set $X$, we put
\[\calD[X]=\{G\in \calC\calS[X]\mid G\mbox{ indecomposable},\: \ic\left(\theta^{\calH\calG}_{\AC}(G)\right)>\ic(G)\}.\]
By Corollary \ref{cor7.2.21}, $\calC\calS^{\double}=\calC\calS_{\delta\calD}$. If $Y\subsetneq X$ and $|Y|\geq 2$, then $\ic\left(\theta^{\calH\calG}_{\AC}(G_{X,Y})\right)=|X|-|Y|+1>1=\ic(G_{X,Y})$ so $G_{X,Y}\in \calD[X]$.
Consequently, $\calC\subseteq \calD$, which implies that $\calC\calS^{\double}=\calC\calS_{\delta\calD}\subseteq \calC\calS_{\delta\calC}$. Let us assume that $\calC\calS^{\double}\subsetneq \calC\calS_{\delta\calC}$.
Let $X$ be a finite set, which we choose of minimal cardinality, such that $\calC\calS^{\double}[X]\subsetneq \calC\calS_{\delta\calC}[X]$. Let $G\in \calC\calS_{\delta\calC}[X]\setminus \calC\calS^{\double}[X]$.
Let $G=G_1\cdots G_k$ be the decomposition of $G$ into indecomposable components. At least one of them, which we denote by $G_i\in \calC\calS[X_i]$, does not belong to $\calC\calS^{\double}$: 
otherwise, as $\calC\calS^{\double}$ is a twisted subbialgebra of $\calC\calS$, $G\in \calC\calS^{\double}[X]$. By minimality of $|X|$, $X_i=X$, so $G$ is indecomposable. As $G\in \calC\calS[X]$, $X$ is a hyperedge of $G$.
As $G\notin \calC\calS_{\delta\calD}[X]$, there exists $\sim\in \calE[G]$, and $Y\subseteq X/{\sim}$, such that $\ic(\theta^{\calH\calG}_{\AC}((G/{\sim})_{\mid Y}))>1$, and $(G/{\sim})_{\mid Y}$ indecomposable. 
As $\calC\calS_{\delta\calC}$ is a twisted double subbialgebra of $\calC\calS$, $(G/{\sim})_{\mid Y}\in \calC\calS_{\delta\calC}[Y]$. If $\sim\neq \sim_G^{(m)}$ or $Y\neq X/{\sim}$, then $|Y/{\sim}|<|X|$. By minimality of $X$, $(G/{\sim})_{\mid Y}\in \calC\calS^{\double}[Y]$,
so $\ic(\theta^{\calH\calG}_{\AC}((G/{\sim})_{\mid Y})>1$: this is a contradiction. Consequently, $\sim=\sim_G^{(m)}$ and $Y=X$, and we deduce that $\ic\left(\theta^{\calH\calG}_{\AC}(G)\right)>1$. Let us denote by $G_i\in \calA\calC[X_i]$ the indecomposable components of $\theta^{\calH\calG}_{\AC}(G)$, for $i\in [k]$, with $k\geq 2$. 
If $X$ is the unique non-trivial hyperedge of $G$, then $\theta^{\calH\calG}_{\AC}(G)=G$ is indecomposable: this is a contradiction. Let $Y$ be another non-trivial hyperedge of $G$, which we choose of minimal cardinality.
Then $Y\subsetneq X$ and $|Y|\geq 2$. Moreover, by minimality, $Y$ is a hyperedge of $\theta^{\calH\calG}_{\AC}(G)$. We assume that it is included into $X_1$, up to a permutation of the $X_i$'s. 
We put $G'_1=G_{\mid X_1}$. Then $G_1=\theta^{\calH\calG}_{\AC}(G'_1)$, as $\theta^{\calH\calG}_{\AC}$ is a morphism of twisted coalgebras. Moreover, $G_1$ is connected and has less hyperedges than $G'_1$, so $G'_1$ is also connected.
By Lemma \ref{lem7.2.26}, there exists $\sim_1\in \calE[G'_1]$, such that $|X_1/{\sim}_1|\geq 2$ and $(G'_1)/{\sim}_1=G_{X_1/{\sim}_1}$. Let us put $\sim=\sim_1\sqcup X_2^2\sqcup\cdots\sqcup X_k^2$. Then $\sim\in \calE[G]$
and $(G/{\sim})=G_{X_1/{\sim}_1,X/{\sim}}$. Therefore, $G\notin \calC\calS_{\delta\calC}[X]$: this is a contradiction. So $\calC\calS^{\double}=\calC\calS_{\delta\calC}$.

Let $X$ be a finite set, and $G\in \calC\calS[X]\setminus \theta^{\calH\calG}_{\CS}\circ \theta^{\calH\calG}_{\AC}(\calC\calS)[X]$. Therefore, $\theta^{\calH\calG}_{\CS}\circ \theta^{\calH\calG}_{\AC}(G)\neq G$. 
Let $Y$ be a hyperedge of $\theta^{\calH\calG}_{\CS}\circ \theta^{\calH\calG}_{\AC}(G)$. Then $\theta^{\calH\calG}_{\AC}(G)_{\mid Y}=\theta^{\calH\calG}_{\AC}(G_{\mid Y})$ is connected. As $\theta^{\calH\calG}_{\AC}(G_{\mid Y})$ has less hyperedges than $G_{\mid Y}$, $G_{\mid Y}$ is connected.
As $G\in \calC\calS[X]$, $Y$ is a hyperedge of $G$. We proved that $\HE(\theta^{\calH\calG}_{\CS}\circ \theta^{\calH\calG}_{\AC}(G))\subseteq \HE(G)$. As $\theta^{\calH\calG}_{\CS}\circ\theta^{\calH\calG}_{\AC}(G)$ and $G$ are different, there exists $Y\in \HE(G)$, not a hyperedge of $\theta^{\calH\calG}_{\CS}\circ\theta^{\calH\calG}_{\AC}(G)$.
Consequently, $\theta^{\calH\calG}_{\AC}(G)_{\mid Y}$ is not connected. Therefore, $\ic(\theta^{\calH\calG}_{\AC}(G)_{\mid Y})>1$, whereas $G_{\mid Y}$ is connected. Consequently, $G\notin \calC\calS_{\delta\calC}[X]=\calC\calS^{\double}[X]$.
By duality, $\calC\calS^{\double}\subseteq \theta^{\calH\calG}_{\CS}\circ\theta^{\calH\calG}_{\AC}(\calC\calS)$.

Let $X$ be a finite set, and $G\in \theta^{\calH\calG}_{\CS}\circ \theta^{\calH\calG}_{\AC}(\calC\calS)[X]$. As $(\theta^{\calH\calG}_{\CS}\circ \theta^{\calH\calG}_{\AC})^2=\theta^{\calH\calG}_{\CS}\circ \theta^{\calH\calG}_{\AC}$, $G=\theta^{\calH\calG}_{\CS}\circ \theta^{\calH\calG}_{\AC}(G)$. For any $Y\subseteq X$, as $ \theta^{\calH\calG}_{\CS}\circ \theta^{\calH\calG}_{\AC}(\calC\calS)$ is a twisted subbialgebra of $\calH\calG$,
$\theta^{\calH\calG}_{\CS}\circ \theta^{\calH\calG}_{\AC}(G)_{\mid Y}=\theta^{\calH\calG}_{\CS}\circ \theta^{\calH\calG}_{\AC}(G_{\mid Y})$, so $\theta^{\calH\calG}_{\CS}\circ \theta^{\calH\calG}_{\AC}(G_{\mid Y})=G_{\mid Y}$. Therefore, $G$ is not of the form $G_{Y,Z}$, with $Z\subsetneq Y$ and $|Z|\geq 2$, as $\theta^{\calH\calG}_{\CS}\circ \theta^{\calH\calG}_{\AC}(G_{Y,Z}) \neq G_{Y,Z}$.
So $G\in \calC\calS_\calC[X]$.

Let $G\in \calC\calS[4]$, whose non-trivial hyperedges are $\{1,2\}$, $\{1,3\}$, $\{2,3\}$, $\{1,2,3\}$ and $\{1,2,3,4\}$. It is not difficult to show that $G\in \calC\calS_\calC[4]$.
Moreover, the non-trivial hyperedges of $\theta^{\calH\calG}_{\CS}\circ\theta^{\calH\calG}_{\AC}(G)$ are $\{1,2\}$, $\{1,3\}$, $\{2,3\}$, and $\{1,2,3\}$, so $\theta^{\calH\calG}_{\CS}\circ \theta^{\calH\calG}_{\AC}(G)\neq G$, and $G\notin \theta^{\calH\calG}_{\CS}\circ \theta^{\calH\calG}_{\AC}(\calC\calS)[4]$.
Consequently, $ \theta^{\calH\calG}_{\CS}\circ \theta^{\calH\calG}_{\AC}(\calC\calS)\subsetneq \calC\calS_\calC$.

Let $G\in \calA\calC[4]$ whose non-trivial hyperedges are $\{1,2\}$, $\{3,4\}$, $\{2,4\}$, and $\{1,3,5\}$.Then $\theta^{\calH\calG}_{\AC}(G)=G$. Let $\sim\in \calE[4]$ whose classes are $\{1,2\}$, $\{3,4\}$, and $\{5\}$. Then
$\theta^{\calH\calG}_{\CS}\circ \theta^{\calH\calG}_{\AC}(G)/{\sim}$ has for non-trivial hyperedges $\{\overline{1},\overline{3}\}$, and $\{\overline{1},\overline{3},\overline{5}\}$. In other words, $\theta^{\calH\calG}_{\CS}\circ \theta^{\calH\calG}_{\AC}(G)=G_{X,Y}$, with $X=[4]/{\sim}$ and $Y=\{\overline{1},\overline{3}\}$.
So $\theta^{\calH\calG}_{\CS}\circ\theta^{\calH\calG}_{\AC}(\theta^{\calH\calG}_{\CS}\circ \theta^{\calH\calG}_{\AC}(G)/{\sim})\neq \theta^{\calH\calG}_{\CS}\circ \theta^{\calH\calG}_{\AC}(G)/{\sim}$, which implies that $\theta^{\calH\calG}_{\CS}\circ \theta^{\calH\calG}_{\AC}(G)/{\sim}\notin \theta^{\calH\calG}_{\CS}\circ\theta^{\calH\calG}_{\AC}(\calH\calG)[[4]/{\sim}]$:
$\theta^{\calH\calG}_{\CS}\circ\theta^{\calH\calG}_{\AC}(\calH\calG)$ is not a twisted double subbialgebra of $\calH\calG$, so is different from $\calC\calS^{\double}$. As a consequence, $\calC\calS^{\double}=\subsetneq \theta^{\calH\calG}_{\CS}\circ\theta^{\calH\calG}_{\AC}(\calC\calS)$.
\end{proof}

\begin{prop}
$\theta^{\calH\calG}_{\AC}:\theta^{\calH\calG}_{\CS}\circ\theta^{\calH\calG}_{\AC}(\calH\calG)\longrightarrow \calA\calC$ and $\theta^{\calH\calG}_{\CS}:\calA\calC\longrightarrow\theta^{\calH\calG}_{\CS}\circ\theta^{\calH\calG}_{\AC}(\calH\calG)$ are two isomorphisms of CTBs, inverse one from the other.
\end{prop}

\begin{proof}
Let $X$ be a finite set. Let $G\in \theta^{\calH\calG}_{\CS}\circ \theta^{\calH\calG}_{\AC}(\calH\calG)[X]$. We put $G=\theta^{\calH\calG}_{\CS}\circ\theta^{\calH\calG}_{\AC}(H)$. Then $\theta^{\calH\calG}_{\AC}(G)\in \calA\calC[X]$. Moreover,
\[\theta^{\calH\calG}_{\CS}\circ\theta^{\calH\calG}_{\AC}(G)=\theta^{\calH\calG}_{\CS}\circ\theta^{\calH\calG}_{\AC}\circ \theta^{\calH\calG}_{\CS}\circ\theta^{\calH\calG}_{\AC}(H)=\theta^{\calH\calG}_{\CS}\circ\theta^{\calH\calG}_{\AC}\circ\theta^{\calH\calG}_{\AC}(H)=\theta^{\calH\calG}_{\CS}\circ\theta^{\calH\calG}_{\AC}(H)=G.\]
Let $G\in \calA\calC[X]$. Then $\theta^{\calH\calG}_{\AC}(G)=G$, and $\theta^{\calH\calG}_{\CS}(G)=\theta^{\calH\calG}_{\CS}\circ\theta^{\calH\calG}_{\AC}(G)\in \theta^{\calH\calG}_{\CS}\circ\theta^{\calH\calG}_{\AC}(\calH\calG)[X]$. Moreover, 
$\theta^{\calH\calG}_{\AC}\circ \theta^{\calH\calG}_{\CS}(G)=\theta^{\calH\calG}_{\AC}(G)=G$.
\end{proof}

\subsection{Antipode in the commutative and cocommutative case}

When $\calB$ is commutative and cocommutative, we obtain the following formula for its antipode:

\begin{prop}\label{prop7.2.29}
Let $\calB$ be a non-degenerate, commutative and cocommutative, CTB. Its antipode is denoted by $S_\calB$.
Let $X$ be a finite set, and $b\in \calB[X]$. We put
\[\calE[b]=\{\sim\in \calE[X]\mid \forall Y\in X/{\sim},\: b_{\mid Y}\mbox{ is indecomposable}\}.\]
For any $\sim\in \calE[X]$, we put
\[b\mid \sim=\prod_{Y\in X/{\sim}} b_{\mid Y}.\]
Then
\[S_\calB(b)=\sum_{\sim\in \calE[b]}\left(\phi^{\calC\calS}_+\left(\theta^\calB_{\CS}(b)/{\sim}\right)(-1)\right) b\mid \sim.\]
\end{prop}

\begin{proof}
By Takeuchi's formula,
\begin{align*}
S_\calB(b)&=\sum_{(X_1,\ldots,X_k)\in \comp[X]} (-1)^k m_{X_1,\dots,X_k}\circ \Delta_{X_1,\ldots,X_k}(b)\\
&=\sum_{(X_1,\ldots,X_k)\in \comp[X]} (-1)^k b_{\mid X_1}\cdots b_{\mid X_k}\\
&=\sum_{\sim'\in \calE[X]} \left((-1)^{\cl(\sim')} \cl(\sim')!\right) b\mid \sim',
\end{align*}
where we used the commutativity and the cocommutativity of $\calB$ for the last equality. Let $\sim'\in \calE[X]$. We put $\sim=\sim_{b\mid\sim'}^{(\ic)}$. Then $\sim\subseteq \sim'$, and $(b\mid \sim')\mid \sim=b\mid \sim$.
Moreover, if $Y\in X/{\sim}$, denoting by $Y'$ the unique element of $X/{\sim}'$ that contains $Y$,
\[b_{\mid Y}=\left(b_{\mid Y'}\right)_{\mid Y}=(b\mid \sim')_{\mid Y}.\]
By definition of $\sim=\sim_{b\mid \sim'}^{(\ic)}$, this is indecomposable. So $\sim\in \calE[b]$. Therefore, 
\begin{align}
\label{EQ7.3} S_\calB(b)&=\sum_{\sim \in \calE[b]}\left(\sum_{\sim\subseteq \sim',\: b\mid\sim=b\mid \sim'}(-1)^{\cl(\sim')} \cl(\sim')!\right)b\mid \sim.
\end{align} 
We put $G=\theta^\calB_{\CS}(b)$. For any $\sim\in \calE[X]$, by Lemma \ref{lem7.2.6},
\begin{align*}
\sim\in \calE[G]&\Longleftrightarrow \forall Y\in X/{\sim},\: G_{\mid Y}\mbox{ is connected}\\
&\Longleftrightarrow \forall Y\in X/{\sim},\: Y\in \HE(G)\\
&\Longleftrightarrow\forall Y\in X/{\sim},\: b_{\mid Y}\mbox{ is indecomposable}\\
&\Longleftrightarrow \sim\in \calE[b],
\end{align*}
so $\calE[G]=\calE[b]$. Let $\sim\in \calE[b]$. Then
\[\phi^{\calC\calS}_+(G/{\sim})=\sum_{(X_1/{\sim}_1,\ldots,X_k/{\sim}_k)\in \calC_+(G/{\sim})}\binom{T}{k},\]
which gives
\[\phi^{\calC\calS}_+(G/{\sim})=\sum_{(X_1/{\sim}_1,\ldots,X_k/{\sim}_k)\in \calC_+(G/{\sim})}(-1)^k.\]
Let $(X_1/{\sim}_1,\ldots,X_k/{\sim}_k)\in \comp[X/{\sim}]$, with $(X_1,\ldots,X_k)\in \comp[X]$, and $\sim=\sim_1\sqcup\cdots \sqcup \sim_k$. Let $\sim'\in \calE[X]$, defined by $X/{\sim}'=\{X_1,\ldots,X_k\}$. Then $\sim\subseteq \sim'$. 
\begin{align*}
&(X_1/{\sim}_1,\ldots,X_k/{\sim}_k)\in \exch_+[G/{\sim}]\\
&\Longleftrightarrow \forall i\in[k],\: (G/{\sim})_{\mid X_i/{\sim}_i}\mbox{ is a discrete hypergraph}\\
&\Longleftrightarrow\forall i\in [k],\:\mbox{the hyperedges of $G$ included in $X_i$ are included in an element of $X/{\sim}$}\\
&\Longleftrightarrow\forall i\in [k],\:G_{\mid X_i}=\left(G_{\mid X_i}\right)\mid\sim_i\\
&\Longleftrightarrow G\mid \sim'=\prod_{i=1}^k\left(G_{\mid X_i}\right)\mid\sim_i\\
&\Longleftrightarrow G\mid \sim'=G\mid\sim.
\end{align*}
We obtain
\[\phi^{\calC\calS}_+(G/{\sim})(-1)=\sum_{\sim\subseteq \sim',\: G\mid\sim=G\mid \sim'}(-1)^{\cl(\sim')} \cl(\sim')!.\]
Let $\sim'\in \calE[X]$, such that $\sim\leq \sim'$. By definition of $G=\theta^\calB_{\CS}(b)$,
\begin{align*}
G\mid\sim'=G\mid\sim&\Longleftrightarrow \forall Y\in X/{\sim},\: \sim'\cap Y^2=\sim_{G_{\mid Y}}^{(\ic)}\\
&\Longleftrightarrow \forall Y\in X/{\sim},\: \sim'\cap Y^2=\sim_{b_{\mid Y}}^{(\ic)}\\
&\Longleftrightarrow b\mid\sim'=b\mid\sim.
\end{align*}
Therefore,
\[\phi^{\calC\calS}_+(G/{\sim})(-1)=\sum_{\sim\subseteq \sim',\:b\mid\sim=b\mid\sim'}(-1)^{\cl(\sim')} \cl(\sim')!.\]
Injecting in (\ref{EQ7.3}), we obtain the result. 
\end{proof}

\begin{remark}
\begin{enumerate}
\item If $\sim\neq \sim'\in \calE[b]$, then $b\mid\sim\neq b\mid \sim'$, as the indecomposable components of $b\mid \sim$ are the elements of $X/{\sim}$, whereas the indecomposable components of $b\mid \sim'$ are the elements of $X/{\sim}'$.
Therefore, we obtained a cancellation-free formula for $S_\calB$.
\item For a given connectedness system $G$, we do not have a combinatorial interpretation of $\phi^{\calC\calS}_+(G)(-1)$. We do have such an interpretation if $G$ is a graph (Stanley's theorem \cite{Stanley1973}): 
$(-1)^{|V(G)|}\phi^\calG_+(G)(-1)$ is the number of acyclic orientations of $G$. For connectedness systems, even the sign of $\phi^{\calC\calS}_+(G)(-1)$ is not clear. Let us give three examples in $\calC\calS[4]$:
\begin{align*}
E(H_1)&=\{\{1,2\},\:\{1,2,3\},\:\{1,2,4\},\: \{1,2,3,4\}\},\\
E(H_2)&=\{\{1,2\},\:\{1,3,4\},\:\{2,3,4\},\:\{1,2,3\},\:\{1,2,3,4\}\},\\
E(H_3)&=\{\{1,2,3,4\}\},
\end{align*}
Then
\begin{align*}
\phi^{\calC\calS}_+(H_1)&=T(T-1)^2(T+1),&\phi_+^{\calC\calS}(H_1)(-1)&=0,\\
\phi^{\calC\calS}_+(H_2)&=T(T-1)(T^2-2),&\phi_+^{\calC\calS}(H_1)(-1)&=-2,\\
\phi^{\calC\calS}_+(H_2)&=T(T-1)(T^2+T+1),&\phi_+^{\calC\calS}(H_1)(-1)&=2.
\end{align*}\end{enumerate}\end{remark}

\subsection{Alternating sign condition in the cocommutative case}

\begin{theo}\label{thoe7.2.30}
Let $\calB$ be a non-degenerate cocommutative CTB. 
For any finite set $X$, and any $b\in \calB[X]$, the graph $\theta^\calB_G(b)\in \calG[X]$ is defined by
\[E\left(\theta^\calB_G(b)\right)=\{Y\in \calS_2[X]\mid b_{\mid Y}\mbox{ is indecomposable}\}.\]
Then $\theta^\calB_G:\calB\longrightarrow \calG$ is a morphism of CTBs. 
Moreover, the following conditions are equivalent:
\begin{enumerate}
\item $\calB$ satisfies the alternating signs condition.
\item For any finite set $X$, and any $b\in \calB[X]$, putting
\[\phi^\calB_+(b)=a_0+\cdots+a_nT^n,\]
then, for any $i\in\left\llbracket\mathrm{val}(\phi^\calB_+(b));\deg(\phi^\calB_+(b))\right\rrbracket$, $a_i\neq 0$.
\item $\epsilon^\calG_+\circ \theta^\calB_G=\epsilon^\calB_+$.
\item $\phi^\calB_+\circ \theta^\calB_G=\phi^\calB_+$.
\end{enumerate}
\end{theo}

\begin{proof}
Let $X$, $Y$ be disjoint finite sets, $a\in \calB[X]$ and $b\in \calB[Y]$. For any $Z\subseteq X\sqcup Y$, 
\[(ab)_{\mid Z}=a_{\mid X\cap Z}b_{\mid Y\cap Z}.\]
Consequently, $(ab)_{\mid Z}$ is indecomposable if, and only if, ($Z\subseteq X$ and $a_{\mid Z}$ is indecomposable) or ($Z\subseteq Y$ and $b_{\mid Z}$ is indecomposable). We obtain that
\[E\left(\theta^\calB_G(ab)\right)=E\left(\theta^\calB_G(a)\right)\sqcup E\left(\theta^\calB_G(b)\right),\]
so $\theta^\calB_G(ab)=\theta^\calB_G(a)\theta^\calB_G(b)$: $\theta^\calB_G$ is a morphism of twisted algebras.

Let $Y\subseteq X$, and $b\in \calB[X]$. For any $Z\subseteq Y$, $\left(b_{\mid Y}\right)_{\mid Z}=b_{\mid Z}$. Therefore, 
\[E\left(\theta^\calB_G(b_{\mid Y})\right)=\{Z\in E\left(\theta^\calB_G(b)\right)\mid Z\subseteq Y\}=E\left(\theta^\calB_G(b)\right)_{\mid Y}.\]
This implies that $\theta^\calB_G$ is a morphism of twisted coalgebras.

$1.\Longrightarrow 2$. Obvious.

$2.\Longrightarrow 3$. Let $X$ be a finite set, and $b\in \calB[X]$. Let us assume that $\epsilon^\calB_+(b)=1$. Then $\ic(b)=|X|$. As $|X|\geq \ic\left(\theta^\calB_G(b)\right)\geq \ic(b)$, $\ic\left(\theta^\calB_G(b)\right)=|X|$, so $\epsilon^\calG_+\circ \theta^\calB_G(b)=1$.

Let us now assume that $\epsilon^\calG_+\circ \theta^\calB_G(b)=1$, and let us prove that $\epsilon^\calB_+(b)=1$.
We proceed by induction on $|X|$. If $|X|\leq 1$, then $b\in \bfD\calB[X]=\calB[X]$, so $\epsilon^\calB_+(b)=1$. Let us assume the result at all ranks $<|X|$, and $|X|\geq 2$. Let $(X_1,\ldots,X_k)\in \comp[X]$, with $k\geq 2$. By the induction hypothesis,
\[{\epsilon^\calB_+}^{\otimes k}\circ \Delta_{X_1,\ldots,X_k}(b)={\epsilon^\calG_+}^{\otimes k}\circ {\theta^\calB_G}^{\otimes k}\circ \Delta_{X_1,\ldots,X_k}(b)={\epsilon^\calG_+}^{\otimes k}\circ \Delta_{X_1,\ldots,X_k}\left(\theta^\calB_G(b)\right)=1.\]
For the last equality, $\theta^\calB_G(b)\in \bfD\calG[X]$, and, as $\bfD\calB[X]$ is a twisted subbialgebra of $\calG$,
\[ \Delta_{X_1,\ldots,X_k}\left(\theta^\calB_G(b)\right)\in \bfD\calG[X_1]\otimes\cdots\otimes \bfD\calG[X_k].\]
Therefore, $(X_1,\ldots,X_k)\in \calC_+(b)$. As a conclusion,
\[\calC_+(b)=\comp[X]\mbox{ or }\comp[X]\setminus \{(X)\}.\]
In the second case, $\phi^\calB_+(b)=T^{|X|}-T$. As we assume condition 2, $|X|=2$. As $\phi^\calB_+(b)\neq T^2$, $\epsilon^\calB_+(b)=0$, contradiction. So $\calC_+(b)=\comp[X]$, and finally $\epsilon^\calB_+(b)=1$.

As both $\epsilon^\calB_+$ and $\epsilon^\calG_+\circ \theta^\calB_G$ are both combinatorial characters of $\calB$, they are equal.

$3.\Longrightarrow 4$. See Proposition \ref{prop3.5.11}-2.

$4.\Longrightarrow 1$. For any finite set $X$, and any $b\in \calB[X]$, $\phi^\calB_+(b)=\phi^\calG_+\left(\theta^\calB_G(b)\right)$. As $\calG$ satisfies the alternating signs condition, it is also the case for $\calB$.
\end{proof}

\begin{remark}
As a consequence, if $\calB$ is a cocommutative CTB with the alternating signs condition (or even with the seemingly weaker condition 2 of Theorem \ref{thoe7.2.30}), the non-zero polynomials $\phi^\calB_+(b)$ obtained from $\calB$ are chromatic polynomials of graphs. 
\end{remark}

\chapter{Structure of $\exch$, $\exch_\pm$ and $\exch_{\pp}$}

\section*{Introduction}

This last chapter does not consider anymore the combinatorial structure of $\exch$ and its sub-objects, but aims to understand their bialgebraic underlying structures.
The main tool is the generalization of Loday's rigidity theorem for infinitesimal bialgebras \cite{Loday2006}, which we rise to the level of species. 
We introduce twisted infinitesimal bialgebras (Definition \ref{defi8.1.3}). The main examples are free twisted algebras, with the concatenation product and deconcatenation coproduct. They are in fact the unique examples (under a condition of connectedness),
by the rigidity Theorem \ref{theo8.1.5}. 
We shall use this Theorem to show that $\exch$, $\exch_+$, $\exch_-$ and $\exch_{\pp}$ are all cofree twisted coalgebras, with the introduction of complementary products.

The problem of the product is solved by Theorem \ref{theo8.1.9}, which states that if $\calB$ is a commutative cofree twisted bialgebra, then it is isomorphic to a shuffle twisted bialgebra. The proof of this Theorem is widely based on the eulerian idempotent and its properties. 
The introduction of a new product $\diamond$ on $\exch$, $\exch_+$ and $\exch_{\pp}$, and $\blackdiamond$ on $\exch_-$
then proves that they are isomorphic to shuffle twisted CTBs (Theorems \ref{theo8.1.11} and \ref{theo8.1.14}), in a nontrivial way.

\section{Twisted infinitesimal bialgebras}

\subsection{Definition and rigidity theorem}

Let us recall this definition, due to \cite{Loday2006}, and the main results on these objects.

\begin{defi}\label{defi8.1.1}
An infinitesimal bialgebra is a triple $(A,m,\Delta)$, such that:
\begin{enumerate}
\item $(A,m)$ is a unitary and associative algebra.
\item $(A,\Delta)$ is a counitary and coassociative coalgebra.
\item For any $a,b\in A$, $\Delta(ab)=(a\otimes 1)\Delta(b)+\Delta(a)(1\otimes b)-a\otimes b$.
\end{enumerate}
\end{defi}

\begin{example}
Let $V$ be a vector space. We equip the tensor algebra $T(V)$ with the concatenation product $m_\conc$ and the deconcatenation coproduct $\Delta_\conc$. Then $(A,m_\conc,\Delta_\conc)$ is an infinitesimal bialgebra.
\end{example}

\begin{theo}\label{theo8.1.2} \cite{Loday2006}
Let $(A,m,\Delta)$ be an infinitesimal bialgebra, connected as a coalgebra. Then $(A,m,\Delta)$ is isomorphic to $(T(\prim(A)),m_\conc,\Delta_\conc)$. 
\end{theo}

Let us give an infinitesimal version of these results.

\begin{defi}\label{defi8.1.3}
A twisted infinitesimal bialgebra is a triple $(\calQ,m,\Delta)$ such that:
\begin{enumerate}
\item $(\calQ,m)$ is a unitary and associative twisted algebra.
\item $\calQ[\emptyset]$ is one-dimensional, generated by the unit $1_\calQ$ of $\calQ$.
\item $(\calQ,\Delta)$ is a counitary and coassociative twisted coalgebra.
\item For any finite sets $X$, $Y,X',Y'$, such that $X\sqcup Y=X'\sqcup Y'$,
\begin{align*}
&\forall a\in \calQ[X], \: \forall b\in \calQ[Y],&\Delta_{X',Y'}(ab)&=\begin{cases}
\Delta_{X',X-X'}(a)(1_\calQ\otimes b)\mbox{ if }X'\subseteq X,\\
(a\otimes 1_\calQ)\Delta_{Y-Y',Y'}(b)\mbox{ if }Y''\subseteq Y,\\
a\otimes b\mbox{ if }(X$, $Y)=(X',Y'),\\
0\mbox{ otherwise}.
\end{cases}\end{align*}\end{enumerate}\end{defi}

\begin{remark}
As $\calQ[\emptyset]=\K 1_\calQ$, for any finite set $X$, for any $b\in \calQ[X]$,
\begin{align*}
\Delta_{\emptyset,X}(b)&=1_\calQ\otimes x,&\Delta_{X,\emptyset}(b)&=b\otimes 1_\calQ.
\end{align*}
Consequently, if $(X',Y')=(X$, $Y)$,
\begin{align*}
\Delta_{X',X-X'}(a)(1_\calQ\otimes b)&=\Delta_{X,\emptyset}(a)(1_\calQ\otimes b)=(a\otimes 1_\calQ)(1_\calQ\otimes b)=a\otimes b,\\
(a\otimes 1_\calQ)\Delta_{Y-Y',Y'}(b)&=(a\otimes 1_\calQ)\Delta_{\emptyset,Y}(b)=(a\otimes 1_\calQ)(1_\calQ\otimes b)=a\otimes b.
\end{align*}
Therefore, the definition of twisted infinitesimal bialgebra is consistent.
\end{remark}

\begin{example}
Let $\calR$ be a linear species, such that $\calR[\emptyset]=\{0\}$. We put
\[T(\calR)[X]=\bigoplus_{(X_1,\ldots,X_k)\in \comp[X]}\calR[X_1]\otimes \cdots \otimes \calR[X_k].\]
We give it the concatenation product: if $X$, $Y$ are disjoint finite sets, $(X_1,\ldots,X_k)\in \comp[X]$, $(Y_1,\ldots,Y_l)\in \comp[Y]$, $x_i\in \calP[X_i]$ for $i\in [k]$ and $y_i\in \calP[Y_j]$ for $j\in [l]$,
\[m_\conc(x_1\cdots x_k\cdots y_1\cdots y_l)=x_1\cdots x_ky_1\cdots y_l.\]
We also give it the deconcatenation coproduct: if $X$ is a finite set, $(X_1,\ldots,X_k)\in \comp[X]$, $x_i\in \calP[X_i]$ for $i\in [k]$, and $X'\sqcup X''=X$,
\[\Delta_\conc(x_1\cdots x_k)=\begin{cases}
x_1\cdots x_i\otimes x_{i+1}\cdots x_k\mbox{ if }X'=X_1\sqcup \cdots \sqcup X_i,\\
0\mbox{ otherwise}.
\end{cases}\]
Then $(T(\calP),m_\conc,\Delta_\conc)$ is a twisted infinitesimal bialgebra. 
\end{example}

\begin{remark}
If $(\calQ,m,\Delta)$ is an infinitesimal bialgebra, then $\calF[\calQ]$ is an infinitesimal bialgebra. 
\end{remark}

Let $\calQ$ be an infinitesimal twisted bialgebra. Then it is connected as a coalgebra, and $\id_\calQ\in \Hom(\calQ,\calQ)_1$, with the notations of the first chapter. We can apply Proposition \ref{prop1.2.7}.

\begin{prop}
Let $\calQ$ be an infinitesimal twisted bialgebra. Then it has an antipode, denoted by $S_\calQ$.
For any finite set $X$, we put
\begin{align*}
\prim(\calQ)[X]&=\{x\in \calQ[X]\mid \forall (X',X'')\in \comp_2[X],\: \Delta_{X,'X''}(x)=0\},\\
\calQ^2[X]&=\sum_{(X',X'')\in \comp_2[X]}m_{X',X''}(\calQ[X']\otimes \calQ[X'']).
\end{align*}
This defines two linear subspecies of $\calQ$.
Then, for any nonempty finite set $X$, 
\[\calQ[X]=\prim(\calQ)[X]\oplus \calQ^2[X],\]
and the projection on $\prim(\calQ)$ in this direct sum is $-S_\calQ$. 
\end{prop}

\begin{proof}
\textit{First step}. We shall use the ideas of the proof of Theorem \ref{theo1.2.19}.
By Proposition \ref{prop1.2.7}., $S_\calQ=(1+(\id_\calQ-\iota_\calQ))^{-1}$ . In other word, for any finite set $X$,
\begin{align*}
&\forall x\in \calQ[X],&S_\calQ(x)&=\sum_{k=1}^{|X|} (-1)^k \sum_{(X_1,\ldots,X_k)\in \comp_k[X]} m_{X_1,\ldots,X_k}\circ \Delta_{X_1,\ldots,X_k}(x).
\end{align*}
As in the proof of Theorem \ref{theo1.2.19}, $m^*:\Hom(\calQ,\calQ)\longrightarrow \Hom(\calQ\otimes\calQ,\calQ)$ and $\Delta^*:\Hom(\calQ\otimes\calQ,\calQ)$ are continuous maps, but, 
as neither $\Delta$ is a morphism of twisted algebras nor $m$ is a morphism of twisted coalgebras, they are not algebra morphisms.

Let $(X,Y)$ be a pair of disjoint finite nonempty sets. We put $\rho_\calQ=\id_\calQ-\iota_\calQ$. 
As $\rho_\calQ$ and $\id_\calQ$ coincide on any $\calQ[Z]$ if $Z\neq \emptyset$, 
\[\rho_\calQ\circ m_{X,Y}=m_{X,Y}\circ (\rho_\calQ\otimes \rho_\calQ)=m_{X,Y}.\]
Let $x\in \calQ[X]$ and $y\in \calQ[Y]$.
\begin{align*}
\rho_\calQ^{\star 2}\circ m_{X,Y}(x\otimes y)&=\rho_\calQ^{\star 2}(xy)\\
&=\sum_{X\sqcup Y=Z'\sqcup Z''}m_{Z',Z''}\circ (\rho_\calQ\otimes \rho_\calQ)\circ \Delta_{Z,'Z''}(xy)\\
&=\sum_{X=X'\sqcup X''} m_{X',X''\sqcup Y}\circ (\rho_\calQ\otimes \rho_\calQ)\left(\Delta_{X,'X''}(x)(1_\calQ\otimes y)\right)\\
&+\sum_{Y=Y'\sqcup Y''} m_{X\sqcup Y',Y''}\circ (\rho_\calQ\otimes \rho_\calQ)\left((x\otimes 1_\calQ)\Delta_{Y',Y''}(y)\right)-m_{X,Y}\circ(\rho_\calQ\otimes \rho_\calQ)(x\otimes y)\\
&=\sum_{\substack{X=X'\sqcup X'',\\X'\neq \emptyset}} m_{X',X''\sqcup Y}\left(\Delta_{X,'X''}(x)(1_\calQ\otimes y)\right)\\
&+\sum_{\substack{Y=Y'\sqcup Y'',\\Y''\neq \emptyset}} m_{X\sqcup Y',Y''}\left((x\otimes 1_\calQ)\Delta_{Y',Y''}(y)\right)-m_{X,Y}(x\otimes y)\\
&=\sum_{(X',X'')\in \comp_2[X]} m_{X',X''\sqcup Y}\left(\Delta_{X,'X''}(x)(1_\calQ\otimes y)\right)+m_{X,Y}(x\otimes y)\\
&+\sum_{(Y',Y'')\in \comp_2[Y]} m_{X\sqcup Y',Y''}\left((x\otimes 1_\calQ)\Delta_{Y',Y''}(y)\right)+m_{X,Y}(x\otimes y)-m_{X,Y}(x\otimes y)\\
&=\sum_{(X',X'')\in \comp_2[X]}m_{X',X'',Y}(\Delta_{X',X''}(x)\otimes y)\\
&+\sum_{(Y',Y'')\in \comp_2[Y]} m_{X',X'',Y''}(x\otimes \Delta_{Y',Y''}(y))+m_{X,Y}(x\otimes y)\\
&=m_{X,Y}\circ (\rho_\calQ^{\star 2}\otimes \rho_\calQ+\rho_\calQ\otimes \rho_\calQ^{\otimes 2}+\rho_\calQ\otimes\rho_\calQ).
\end{align*}
We used the infinitesimal compatibility for the third equality, and the fact that $\rho_\calQ(1_\calQ)=0$ and is the identity on any $\calQ[Z]$ with $Z\neq \emptyset$ for the fourth one. Then, a direct induction proves the following: putting
\[S_n=\sum_{k=1}^{n-1} \rho_\calQ^{\star k}\otimes \rho_\calQ^{\star (n-k)},\]
for $n\geq 1$, with the convention $_0=0$, then, for any $n\geq 1$,
\[\rho_\calQ^{\star n}\circ m_{X,Y}=m_{X,Y}\circ (S_{n+1}+S_n).\]
Therefore, by continuity of $m^*$, 
\begin{align*}
S_\calQ\circ m_{X,Y}&=\sum_{n=0}^\infty (-1)^n \rho_\calQ^{\star n}\circ m_{X,Y}\\
&=\sum_{n=1}^\infty (-1)^n \rho_\calQ^{\star n}\circ m_{X,Y}\\
&=m_{X,Y}\circ\left(\sum_{n=1}^\infty (-1)^n (S_{n+1}+S-n)\right)\\
&=m_{X,Y}\circ\left(-\sum_{n=2}^\infty (-1)^n S_n+\sum_{n=1}^\infty (-1)^n S_n\right)\\
&=-m_{X,Y}\circ S_1\\
&=0.
\end{align*}
So $S_\calQ\circ m_{X,Y}=0$, if both $X$ and $Y$ are nonempty.

We also have $\Delta_{X,Y}\circ \rho_\calQ=\rho_\calQ\circ \Delta_{X,Y}=\Delta_{X,Y}$. Let $x\in \calQ[X\sqcup Y]$. 
If $X'\sqcup X''=Y\sqcup Z$, we shall use the Sweedler's notation
\[\Delta_{X,'X''}(x)=x_{\mid_1X'}\otimes x_{\mid_2X''},\]
and so on. Then
\begin{align*}
\rho_\calQ^{\star 2}(x)&=\sum_{Z'\sqcup Z''=X\sqcup Y} m_{Z',Z''}\circ (\rho_\calQ\otimes \rho_\calQ)\circ \Delta_{Z',Z''}(x)\\
&=\sum_{(Z',Z'')\in \comp_2[X\sqcup Y]} x_{\mid_1Z'}x_{\mid_2Z''},
\end{align*}
so
\begin{align*}
\Delta_{X,Y}\circ \rho_\calQ^{\star 2}(x)&=\sum_{(Z',Z'')\in \comp_2[X\sqcup Y]} \Delta_{X,Y}(x_{\mid_1Z'}x_{\mid_2Z''})\\
&=\sum_{(X',X'')\in \comp_2[X]} x_{\mid_1X'}x_{\mid_2X''}\otimes x_{\mid_3Y}+x\otimes y
+\sum_{(Y',Y'')\in \comp_2[Y]}x_{\mid_1X}\otimes y_{\mid_2Y'}y_{\mid_3Y''}+0\\
&=(\rho_\calQ^{\star 2}\otimes \rho_\calQ+\rho_\calQ\otimes \rho_\calQ^{\otimes 2}+\rho_\calQ\otimes\rho_\calQ)\circ \Delta_{X,Y}(x).
\end{align*}
We used the infinitesimal compatibility for the second equality. Therefore, an induction proves that any $n\geq 1$,
\[\Delta_{X,Y}\circ \rho_\calQ^{\star n}=(S_{n+1}+S_n)\circ \Delta_{X,Y}.\]
By continuity of $\Delta^*$,
\begin{align*}
\Delta_{X,Y}\circ S_\calQ(x)&=\sum_{n=0}^\infty (-1)^n\Delta_{X,Y}\circ \rho_\calQ^{\star n}\\
&=\sum_{n=1}^\infty (-1)^n\Delta_{X,Y}\circ \rho_\calQ^{\star n}\\
&=\left(\sum_{n=1}^\infty (-1)^n (S_{n+1}+S-n)\right)\circ \Delta_{X,Y}\\
&=\left(-\sum_{n=2}^\infty (-1)^n S_n+\sum_{n=1}^\infty (-1)^n S_n\right)\circ \Delta_{X,Y}\\
&=-S_1\circ \Delta_{X,Y}\\
&=0.
\end{align*}
Therefore, $S_\calQ(x)$ is primitive.

If $x\in \prim(\calQ)[X]$, then $S_\calQ(x)=-x$, as only the term $k=1$ contributes. As $S_\calQ$ takes its values in $\prim(\calQ)$, $-S_\calQ$ is a projection on $\prim(\calQ)$.
We denote by $\calQ^2$ the linear subspecies of $\calQ$ defined by $\calQ^2[\emptyset]=(0)$ and, if $X\neq \emptyset$,
\[\calQ^2[X]=\sum_{(X',X'')\in \comp[X]} m_{X',X''}(\calQ[X']\otimes \calQ[X'']).\]
We already observed that $\calQ^2\subseteq \ker(-S_\calQ)$. Moreover, if $X\neq \emptyset$ and $x\in \calQ[X]$,
\begin{align*}
x&=-S_\calQ(x)+x+S_\calQ(x)\\
&=\underbrace{S_\calQ(x)}_{\in \prim(\calQ)[X]}+\underbrace{\sum_{k=2}^\infty \sum_{(X_1,\ldots,X_k)\in \comp_k[X]} (-1)^km_{X_1,\ldots,X_k}\circ \Delta_{X_1,\ldots,X_k}(x)}_{\in \calQ^2[X]},
\end{align*}
so $\calQ[X]=\prim(\calQ)[X]+\calQ^2[X]$. As a consequence of all these results, $\calQ[X]=\prim(\calQ)[X]\oplus\calQ^2[X]$, and $-S_\calQ$ is the projection on $\prim(\calQ)$ which vanishes on $\calQ^2[X]$. 
\end{proof}

Let us now state the twisted rigidity theorem:

\begin{theo}\label{theo8.1.5}
Let $(\calQ,m,\Delta)$ be a twisted infinitesimal bialgebra. For any finite set $X$, we put
\[\prim(\calQ)[X]=\{x\in \calQ[X]\mid \forall (X',X'')\in \comp[X],\: \Delta_{X,'X''}(x)=0\}.\]
Then $(\calQ,m,\Delta)$ is isomorphic to $(T(\prim(\calQ)),m_\conc,\Delta_\conc)$.
\end{theo}

\begin{proof} 
Let us firstly prove that $\calQ$ is generated by $\prim(\calQ)$ as a twisted algebra.
Let us denote by $\calQ'$ the twisted subalgebra of $\calQ$ generated by $\prim(\calQ)$.
Let us consider a finite set $X$, and $x\in \calQ[X]$, and let us show that $x\in \calQ'[X]$ by induction on $|X|$.
If $X=\emptyset$, we can restrict ourselves to $x=1_\calQ$, as $\calQ[\emptyset]$ is one-dimensional, and the result is then obvious. Let us assume that $|X|\geq 1$ and the result true for any finite set $Y$ with $|Y|<|X|$.
As $\calQ[X]=\prim(\calQ)[X]\oplus \calQ^2[X]$, we can write $x$ under the form
\[x=y+\sum_{i=1}^n x'_ix''_i,\]
with $y\in \prim(\calQ)[X]$ and $x'_i\in \calQ[X'_i]$, $x''_i\in \calQ[X''_i]$ for a certain $(X'_i,X''_i)\in \comp_2[X]$.
Applying the induction hypothesis, $x'_i \in \calQ'[X'_i]$ and $x''_i\in \calQ'[X''_i]$ for any $i\in [n]$.
As obviously $y\in \calQ'[X]$, $x\in \calQ'[X]$. So $\calQ[X]=\calQ'[X]$.

\textit{Third step}. As a consequence, we obtain a surjective morphism of twisted algebras
\[\Phi:\left\{\begin{array}{rcl}
T(\prim(\calQ))&\longrightarrow&\calQ\\
x_1\cdots x_k\in \prim(\calQ)[X_1]\otimes \cdots \otimes \prim(\calQ)[X_k]&\longmapsto&m_{X_1,\ldots,X_k}(x_1\otimes \cdots \otimes x_k).
\end{array}\right.\]
Let us prove that $\Phi$ is a twisted coalgebra morphism. Let us consider, for any finite set $X$,
\[\calQ'[X]=\{x\in T(\prim(\calQ))[X]\mid \forall (X',X'')\in \comp[X], \:(\Phi\otimes \Phi)\circ \Delta_{X',X''}(x)=\Delta_{X',X''}\circ \Phi(x)\}.\]
Obviously, $1\in \calQ'[\emptyset]$ and, for any finite set $X$, $\prim(\calQ)[X]\subseteq \calQ'[X]$. 
As both $T(\prim(\calQ))$ and $\calQ$ satisfy the infinitesimal compatibility, $\calQ'$ is a twisted subalgebra of $T(\prim(\calQ))$. Therefore, as $\prim(\calQ)\subseteq \calQ'$, $\calQ'=T(\prim(\calQ))$. 
As a conclusion, $\Phi$ is a coalgebra morphism.

It remains to prove that $\Phi$ is injective. Let us assume that it is not the case: then $\ker(\Phi)$ is a nonzero subspecies of $T(\prim(\calQ))$. 
Let us take a finite set $X$, of minimal cardinality, such that $\ker(\Phi)[X]\neq \{0\}$, and let $x\in \calQ[X]$, nonzero, such that $\Phi(x)=0$. As $\Phi$ is a coalgebra morphism, for any $(X',X'')\in \comp[X]$,
\[(\Phi\otimes \Phi)\circ \Delta_{X',X''}(x)=\Delta_{X',X''}\circ \Phi(x)=0.\]
By minimality of $|X|$, $\Phi$ is injective on $T(\prim(\calQ))[X']$ and $T(\prim(\calQ))[X'']$. 
Therefore, $\Delta_{X',X''}(x)=0$, so $x\in \prim(T(\prim(\calQ)))[X]=\prim(\calQ)[X]$. As a consequence, $\Phi(x)=x=0$:
this is a contradiction. So $\Phi$ is injective.
\end{proof}

\subsection{Cocommutative products on a cofree twisted coalgebra}

\begin{theo}[\textbf{Universal property of $T(\calR)$}]
Let $\calQ$ be a twisted coalgebra, such that $\calQ[1]$ is $1$-dimensional, and $\calR$ be a linear species such that $\calR[\emptyset](0)$.
Let $\phi:\calQ\longrightarrow\calR$ be a species morphism. There exists a unique morphism $\Phi:\calQ\longrightarrow T(\calR)$ of twisted coalgebras such that $\pi\circ \Phi=\phi$, where $\pi:T(\calR)\longrightarrow \calR$ is the canonical surjection. For any nonempty finite set $X$, 
\begin{align}
\label{EQ8.1}\Phi_X&=\sum_{k=1}^{|X|} \sum_{X_1,\ldots,X_k\in \comp_k[X]} m_{X_1,\ldots,X_k}\circ \left(\phi_{X_1}\otimes\cdots\otimes \phi_{X_k}\right)\circ \Delta_{X_1,\ldots,X_k},
\end{align}
where $m$ is the concatenation product of $T(\calR)$.
\end{theo}

\begin{proof}
As $\calQ[\emptyset]$ is 1-dimensional, there exists a unique element $1_\calQ\in \calQ[\emptyset]$ such that $\varepsilon_\Delta(1_\calQ)=1$. Then $\Delta_{\emptyset,\emptyset}(1_\calQ)=1_\calQ\otimes 1_\calQ$, 
and any morphism $\phi:\calQ\longrightarrow T(\calR)$ of twisted coalgebras sends necessarily $1_\calQ$ to $1$.

\textit{Existence}. Let us consider $\Phi$ defined by (\ref{EQ8.1}). Then 
\[(\pi\circ \Phi)_X=m_X\circ \phi_X\circ \Delta_X+0=\phi_X,\]
so $\pi\circ\Phi=\phi$. Let $(X,Y)$ be a pair of disjoint finite sets, and $x\in \calQ[X\sqcup Y]$. 
If $X=\emptyset$, then 
\[{\Delta_{\conc}}_{X,Y}\circ \Phi(x)=1\otimes \Phi(x)=\Phi(1_\calQ)\otimes \Phi(x)=(\Phi\otimes \Phi)\circ \Delta_{X,Y}(x).\]
The same can be proved if $Y=\emptyset$. Let us now assume that both $X$ and $Y$ are nonempty.
\begin{align*}
{\Delta_{\conc}}_{X,Y}\circ\Phi(x)&=\sum_{k=1}^{|X|+|Y|} \sum_{\substack{X_1,\ldots,X_k\in \comp_k[X\sqcup Y],\\ X_1\sqcup \ldots \sqcup X_i=X,\\ X_{i+1}\sqcup\ldots\sqcup X_k=Y}}
\left(m_{X_1,\ldots,X_i}\circ\left(\phi_{X_1}\otimes \cdots\otimes \phi_{X_i}\right)\circ \Delta_{X_1,\ldots,X_i}\right) \\
&\hspace{1cm}
\otimes \left(m_{X_{i+1},\ldots,X_k}\circ\left(\phi_{X_{i+1}}\otimes \cdots\otimes \phi_{X_k}\right)\circ \Delta_{X_{i+1},\ldots,X_k}\right)\circ \Delta_{X,Y}(x)\\
&=\sum_{i=1}^{|X|}\sum_{j=1}^{|Y|}\sum_{\substack{(X_1,\ldots,X_i)\in \comp_i[X],\\ (Y_1,\ldots,Y_j)\in \comp_j[Y]}}
\left(m_{X_1,\ldots,X_i}\circ\left(\phi_{X_1}\otimes \cdots\otimes \phi_{X_i}\right)\circ \Delta_{X_1,\ldots,X_i}\right) \\
&\hspace{1cm}
\otimes \left(m_{Y_1,\ldots,Y_j}\circ\left(\phi_{Y_1}\otimes \cdots\otimes \phi_{Y_j}\right)\circ \Delta_{Y_1,\ldots,Y_j}\right)\circ \Delta_{X,Y}(x)\\
&=(\Phi\otimes \Phi)\circ \Delta_{X,Y}(x).
\end{align*}
So $\Phi$ is a morphism of twisted coalgebras.

\textit{Existence}. Let $\phi,\psi:\calQ\longrightarrow T(\calR)$ be two morphism of twisted coalgebras, such that $\pi\circ \phi=\pi\circ\psi$. Let us show that for any finite set $X$, $\phi_X=\psi_X$. We proceed by induction on $|X|$.
We already observed this when $X=\emptyset$. Let us assume the result at all ranks $<|X|$. 
If $(X',X'')\in \comp[X]$, for any $x\in \calQ[X]$,
\[{\Delta_{\conc}}_{X',X''}\circ\phi(x)=(\phi\otimes\phi)\circ\Delta_{X',X''}(x)=(\psi\otimes\psi)\circ\Delta_{X',X''}(x)={\Delta_{\conc}}_{X',X''}\circ\psi(x),\]
so $\phi(x)-\psi(x)\in \prim(T(\calR))[X]=\calR[X]$. As $\pi\circ\phi(x)-\pi\circ\psi(x)=0$, $\phi(x)=\psi(x)$.
\end{proof}

\begin{prop}
Let $\calR$ be a linear species such that $\calR[\emptyset]=(0)$. The following defines an associative and commutative product on $T(\calR)$, making it a bialgebra: for any pairwise disjoint nonempty finite sets $X_i$ and $x_i\in \calR[X_i]$ for i in $[k+l]$,
\begin{align}
\label{EQ8.2}
x_1\cdots x_k \shuffle x_{k+1}\cdots x_{k+l}&=\sum_{\sigma\in \sh(k,l)} x_{\sigma^{-1}(1)}\cdots x_{\sigma^{-1}(k+l)}.
\end{align}
\end{prop}

\begin{proof}
We consider $\calQ=T(\calR)\otimes T(\calR)$, which is a connected twisted coalgebra, and 
\[\phi:\left\{\begin{array}{rcl}
\calQ&\longrightarrow&\calR\\
x_1\cdots x_k\otimes x_{k+1}\cdots x_l&\longmapsto&\begin{cases}
x_1\mbox{ if }k+l=1,\\
0\mbox{ otherwise}.
\end{cases}
\end{array}\right.\]
We denote by $\shuffle:\calQ\longrightarrow T(\calR)$ the unique morphism of twisted coalgebras such that $\pi\circ\shuffle=\phi$. It is not difficult to show that $\shuffle$ satisfies (\ref{EQ8.2}).

Let us show that $\shuffle$ is associative. By composition, both $\shuffle\circ (\shuffle\otimes \id_{T(\calR)})$
and $\shuffle\circ(\id_{T(\calR)}\otimes \shuffle)$ are morphisms of twisted coalgebras from $T(\calR)^{\otimes 3}$ to $T(\calR)$. Moreover, if $x_1,\ldots,x_{k+l+m}$ are elements of different $\calR[X_i]$'s,
\begin{align*}
\pi\circ\shuffle\circ (\shuffle\otimes \id_{T(\calR)})(x_1\cdots x_k\otimes x_{k+1}\cdots x_{k+l}\otimes x_{k+l+1}\cdots x_{k+l+m})\\
=\pi\circ \shuffle\circ(\id_{T(\calR)}\otimes \shuffle)(x_1\cdots x_k\otimes x_{k+1}\cdots x_{k+l}\otimes x_{k+l+1}\cdots x_{k+l+m})&=\begin{cases}
x_1\mbox{ if }k+l+m=1,\\
0\mbox{ otherwise}.
\end{cases}\end{align*}
By the unicity in the universal property, $\shuffle\circ (\shuffle\otimes \id_{T(\calR)})=\shuffle\circ(\id_{T(\calR)}\otimes \shuffle)$.

Let us now show that $\shuffle$ is commutative. Both $\shuffle$ and $\shuffle \circ \tau$ are morphisms of twisted coalgebras from $T(\calR)^{\otimes 2}$ to $T(\calR)$. Moreover, if $x_1,\ldots,x_{k+l}$ are elements of different $\calR[X_i]$'s,
\begin{align*}
\pi\circ \shuffle(x_1\cdots x_k\otimes x_{k+1}\cdots x_{k+l})\\
=\pi\circ \shuffle\circ \tau(x_1\cdots x_k\otimes x_{k+1}\cdots x_{k+l})
&=\begin{cases}
x_1\mbox{ if }k+l=1,\\
0\mbox{ otherwise}.
\end{cases}\end{align*}
By the unicity in the universal property, $\shuffle=\shuffle\circ \tau$, so $\shuffle$ is commutative. 
\end{proof}

\begin{lemma}\label{lem8.1.8}
Let $\calR$ be a linear species such that $\calR[\emptyset]=(0)$. Let $\varpi:T(\calR)\longrightarrow \calR$ be any projection. We denote by $\phi$ the endomorphism of the twisted coalgebra $T(\calR)$ such that $\pi\circ\phi=\varpi$.
\begin{enumerate}
\item $\phi$ is an automorphism.
\item Let us assume that $m$ is a product on $T(\calR)$, making it a twisted bialgebra with the coproduct $\Delta_{\conc}$.
Then $\phi$ is a morphism of twisted bialgebras from $(T(\calR),m, \Delta_{\conc})$ to $(T(\calR),\shuffle,\Delta_{\conc})$ if, and only if, for any pair $(X,Y)$ of disjoint finite nonempty sets, 
$\varpi\circ m_{X,Y}=0$.
\end{enumerate}\end{lemma}

\begin{proof}
1. Let $X_1,\ldots,X_k$ be pairwise disjoint finite nonempty sets, and $x_i\\in \calR[X_i]$ for $i\in [k]$. Then
\begin{align*}
\phi(x_1\ldots x_k)&=\sum_{p=1}^l\sum_{1\leq i_1<\ldots<i_p=k}
\varpi(x_1\cdots x_{i_k})\cdots \varpi(x_{i_{p-1}+1}\cdots x_k)
\end{align*}
As $\varpi$ is a projection on $\calR$, the unique term in $\phi(x_1\ldots x_k)$ of length $k$ is $x_1\cdots x_k$,
and the other terms are of length $<k$. A triangularity argument shows that $\phi$ is bijective.

2. $\Longrightarrow$. If $\phi$ is a morphism of twisted bialgebras, then if $(X,Y)$ is a pair of disjoint finite nonempty sets,
\[\varpi\circ m_{X,Y}=\pi\circ \phi\circ m_{X,Y}=\pi\circ \shuffle_{X,Y}\circ (\phi\otimes \phi)=0,\]
as $\pi\circ \shuffle_{X,Y}=0$.

$\Longleftarrow$. Let us assume that for any pair $(X,Y)$ of of disjoint finite nonempty sets, $\varpi\circ m_{X,Y}=0$.
Let us consider $\psi=\phi\circ m$ and $\psi'=\shuffle \circ (\phi\otimes \phi)$. Both are morphism of twisted coalgebras from $T(\calR)^{\otimes 2}$ to $T(\calR)$, by composition. By the unicity in the universal property of $T(\calR)$, it is enough to prove that $\pi\circ \psi=\pi\circ \psi'$. Let $(X,Y)$ be a pair of disjoint finite sets. 
If $X=\emptyset$, for any $y\in T(\calR)[Y]$,
\[\pi\circ\psi(1\otimes y)=\pi\circ\phi(y)=\varpi(y),\]
whereas
\[\pi\circ\psi'(1\otimes y)=\pi\circ \shuffle(1\otimes \phi(y))=\pi\circ\phi(y)=\psi(y),\]
The same result holds if $Y=\emptyset$. If both $X$ and $Y$ are nonempty, then for $x\in T(\calR)[X]$ and $y\in T(\calR)[Y]$,
\[\pi\circ \psi(x\otimes y)=\psi\circ m_{X,Y}(x\otimes y)=0,\]
whereas
\[\pi\circ \psi'(x\otimes y)=\pi\circ \shuffle_{X,Y}\circ (\phi_X\otimes \phi_Y)=0,\]
as $\pi\circ \shuffle_{X,Y}=0$. Consequently, $\pi\circ \psi=\pi\circ \psi'$, so $\psi=\psi'$, and $\phi$ is a morphism of twisted bialgebras. 
\end{proof}

\begin{theo}\label{theo8.1.9}
Let $\calR$ be a linear species such that $\calR[\emptyset]=(0)$, and $m$ a product on $T(\calR)$ such that $(T(\calR),m,\Delta_{\conc})$ is a commutative twisted bialgebra. Then $(T(\calR),m,\Delta_{\conc})$ is isomorphic to $(T(\calR),\shuffle,\Delta_{\conc})$.
\end{theo}

\begin{proof}
We shall use the eulerian idempotent $\ln(\id_{T(\calR)})$. Let us fix a integer $n$. By Theorem \ref{theo1.2.19}:
\begin{itemize}
\item If $x\in \calR[n]$, then $x$ is primitive, so $\ln(\id_\calQ)(x)=x$.
\item If $(X',X'')\in \comp[n]$, then $\ln(\id_\calQ)\circ m_{X',X''}=0$.
\end{itemize}
Therefore, 
\[\calR[n]\cap \left(\sum_{(X,'X'')\in \comp[X]} m_{X',X''}(T(\calR)[X']\otimes T(\calR)[X'])\right)=(0).\]
By semi-simplicity of the representations of the symmetric group $\sym_n$, there exists a subspace $\calS[n]$ of $T(\calR)[n]$, stable under all the bijections $T(\calR)[\sigma]$ for $\sigma\in \sym_n$, such that 
\begin{align*}
T(\calR)[n]&=\calR[n]\oplus \calS[n],&\sum_{(X,'X'')\in \comp[n]} m_{X',X''}(T(\calR)[X']\otimes T(\calR)[X'])\subseteq \calS[n].
\end{align*}
Let now $X$ be any finite set. We denote its cardinality by $n$. Choosing a bijection $\sigma:[n]\longrightarrow X$,
we put $\calS[X]=T(\calR)[\sigma](\calS[n])$. As $\calS[n]$ is stable under the action of the symmetric group $\sym_n$, $\calS[X]$ does not depend on the choice of $\sigma$, and this defines a subspecies $\calS$ of $T(\calR)$. Moreover, for any finite set $X$, as $T(\calR)[\sigma](\calR[n])=\calR[X]$,
\begin{align*}
T(\calR)[X]&=\calR[X]\oplus \calS[X],&\sum_{(X,'X'')\in \comp[X]} m_{X',X''}(T(\calR)[X']\otimes T(\calR)[X'])\subseteq \calS[X].
\end{align*}
Let $\varpi$ the projection on $\calR$ which vanishes on $\calS$, and let $\phi$ the associated endomorphism of the twisted coalgebra $T(\calR)$, as in Lemma \ref{lem8.1.8}. This Lemma gives that $\phi$ is an isomorphism of twisted bialgebras from $(T(\calR),m,\Delta_{\conc})$ to $(T(\calR),\shuffle,\Delta_{\conc})$, as by construction, 
if $(X,Y)$ is a pair of finite nonempty sets, then $\varpi\circ m_{X,Y}=0$. 
\end{proof}

\subsection{Infinitesimal product on $\exch$}

\begin{defi}
Let $X$, $Y$ be disjoint finite sets, $I\in \exch[X]$, $J\in \exch[Y]$. We put
\[I\diamond J=\{(\pi,\pi')\mid \pi\in I,\: \pi'\in J\}\subseteq \comp[X\sqcup Y].\]
\end{defi}

\begin{theo}\label{theo8.1.11}
$(\exch,\diamond,\Delta)$ is a twisted infinitesimal bialgebra. Moreover, $\exch_+$ and $\exch_{\pp}$ are twisted subbialgebras of $(\exch,\diamond,\Delta)$.
\end{theo}

\begin{proof}
\textit{First step}. Let $X$, $Y$ be disjoint finite sets, $I\in \exch[X]$, $J\in \exch[Y]$.
Let us prove that $I\diamond J\in \exch[X\sqcup Y]$. Let $(\pi_1,\pi_2),(\pi'_1,\pi'_2)\in I\diamond J$, with $\pi_i,\pi'_i\in \comp[X_i]$ for $i\in [2]$, and $X_1\sqcup X_2=X\sqcup Y$. According to the form of the elements of the elements of $I\diamond J$, two cases are possible.
\begin{enumerate}
\item $X\subseteq X_1$. Then, we can put $\pi_1=(\pi_{11},\pi_{12})$ and $\pi'_1=(\pi_{11}',\pi_{12}')$, with $\pi_{11},\pi'_{11} \in \comp[X]$. By definition of $I\diamond J$, $\pi_{11}$, $\pi'_{11}\in I$, and $(\pi_{12},\pi_2)$, $(\pi'_{12},\pi'_2)\in J$. As $J\in \exch[Y]$, $(\pi'_{12},\pi_2)\in J$. 
By definition of $I\diamond J$, $(\pi'_{11}, \pi'_{12},\pi_2)=(\pi'_1,\pi_2)\in I\diamond J$.
\item $X_1\subseteq X$. Then, $Y\subseteq Y_1$, and we can put $\pi_2=(\pi_{21},\pi_{22})$ and $\pi'_2=(\pi_{21}',\pi_{22}')$, with $\pi_{22},\pi'_{22} \in \comp[Y]$. 
By definition of $I\diamond J$, $\pi_{22}$, $\pi'_{22}\in J$, and $(\pi_1,\pi_ {21})$, $(\pi'_1,\pi'_{21})\in I$. As $I\in \exch[X]$, $(\pi'_1,\pi_{21})\in I$. 
By definition of $I\diamond J$, $(\pi'_1,\pi_{21},\pi_{22})=(\pi'_1,\pi_2)\in I\diamond J$.
\end{enumerate}
Therefore, $I\diamond J\in \exch[X\sqcup Y]$.

Let us assume that $I\in \exch_{\pp}[X]$ and $J\in \exch_{\pp}[Y]$. 
If $\pi\in I$ and $\pi'\in J$, then 
\[\ell(\pi,\pi')=\ell(\pi)+\ell(\pi')=|X|+|Y|=|X\sqcup Y|,\]
so $I\diamond J\in \exch_{\pp}[X\sqcup Y]$. 
Let us now assume that $I\in \exch_+[X]$ and $J\in \exch_+[Y]$. Let us prove that $I\diamond J\in \exch_+[X\sqcup Y]$. Let $\pi\in I\diamond J$ and $\pi'\in \comp[X\sqcup Y]$, such that $\pi\leq \pi'$.
By definition of $I\diamond J$, $\pi=(\pi_1,\pi_2)$, with $\pi_1\in I$ and $\pi_2\in J$. By definition of the order on $\comp$, $\pi'=(\pi'_1,\pi_2)$, with $\pi_i\leq \pi'_i$ for $i\in [2]$.
As $I\in \exch_+[X]$ and $J\in \exch_+[Y]$, $\pi'_1\in I$ and $\pi'_2\in J$, so $(\pi'_1,\pi'_2)=\pi'\in I\diamond J$.

\textit{Second step}. Let us prove that $\diamond$ is associative. Let $X$, $Y,Z$ be pairwise disjoint finite sets, $I\in \exch[X]$, $J\in\exch[Y]$, $K\in \exch[Z]$. Then
\[(I\diamond J)\diamond K=I\diamond (J\diamond K)=\{(\pi,\pi',\pi'')\mid \pi\in I,\:\pi'\in J,\:\pi''\in K\}.\]

\textit{Last step}. It remains to show the infinitesimal compatibility. Let $X$, $Y,X',Y'$ be finite sets, such that $X\sqcup Y=X'\sqcup Y'$, $I\in \exch[X]$ and $J\in \exch[Y]$. Let us assume that $\Delta_{X',Y'}(I\diamond J)\neq 0$.
Then, by definition of $\Delta$, there exists $\pi\in \comp[X]$, $\pi'\in \comp[Y]$, such that $(\pi,\pi')\in I\diamond J$.
By definition of $I\diamond J$, $(\pi,\pi')=(\varpi,\varpi')$, with $\varpi\in \comp[X]$ and $\varpi'\in \comp[Y]$.
Consequently, $X\subseteq X'$ or $X'\subseteq X$. Dually, if $X$ is not included in $X'$, nor $X'$ is included in $X$,
then $\Delta_{X',Y'}(I\diamond J)=0$.

Let us assume that $X\subseteq X'$. Let $\pi\in (I\diamond J)_{\mid_LX'}$, and $\pi'\in (I\diamond J)_{\mid_RY'}$. Then $(\pi,\pi')\in I\diamond J$, by the exchange condition. By definition of $I\diamond J$, $\pi=(\varpi,\varpi')$, with $\varpi\in I$, and $(\varpi',\pi')\in J$.
So $\pi' \in J_{\mid_R Y'}$, and $\pi\in I\diamond J_{\mid_LY\setminus Y'}$. 
We proved that $(I\diamond J)_{\mid_LX'}\subseteq I\diamond J_{\mid_LY\setminus Y'}$ and $(I\diamond J)_{\mid_RY'}\subseteq J_{\mid_R Y'}$. 

Conversely, let $\pi\in I\diamond J_{\mid_LY\setminus Y'}$ and $\pi'\in J_{\mid_R Y'}$.
Then $\pi=(\pi_1,\pi_2)$, with $\pi_1\in I$ and $\pi_2\in J_{\mid_LY\setminus Y'}$. By the exchange condition, $(\pi_2,\pi')\in J$. So $(\pi_1,\pi_2,\pi')\in I\diamond J$, which gives $\pi\in (I\diamond J)_{\mid_LX'}$ and $\pi'\in (I\diamond J)_{\mid_RY'}$. We proved that $(I\diamond J)_{\mid_LX'}= I\diamond J_{\mid_LY\setminus Y'}$ and $(I\diamond J)_{\mid_RY'}= J_{\mid_R Y'}$. Therefore,
\[\Delta_{X',Y'}(I\diamond J)= I\diamond J_{\mid_LY\setminus Y'}\otimes J_{\mid_R Y'}=(I\otimes 1)\diamond \Delta_{Y\setminus Y',Y'}(J).\]

Similarly, if $X'\subseteq X$, we obtain
\[\Delta_{X',Y'}(I\diamond J)=\Delta_{X',X\setminus X'}(I)\diamond(1\otimes J).\]
So $(\exch,\diamond,\Delta)$ is a twisted infinitesimal bialgebra.
\end{proof}

Combining Theorems \ref{theo8.1.5} and \ref{theo8.1.9}: 

\begin{cor}\begin{enumerate}
\item As a twisted bialgebra, $\exch$ is isomorphic to $(T(\prim(\exch)),\shuffle,\Delta_\conc)$.
\item As a twisted bialgebra, $\exch_+$ is isomorphic to $(T(\prim(\exch_+)),\shuffle,\Delta_\conc)$.
\item As a twisted bialgebra, $\exch_{\pp}$ is isomorphic to $(T(\prim(\exch_{\pp})),\shuffle,\Delta_\conc)$.
\end{enumerate}\end{cor}

\begin{remark}
This does not work for $\exch_-$. For example, if $X$, $Y$ are disjoint sets, $\{(X)\}\in \exch_-[X]$, $\{(Y)\}\in \exch_-[Y]$,
but
\[\{(X)\}\diamond \{(Y)\}=\{(X,Y)\}\notin \exch_-[X\sqcup Y].\]
\end{remark}

\subsection{Infinitesimal structure on $\exch_-$}

\begin{defi}
Let $X$, $Y$ be disjoint finite sets, $I\in \exch_-[X]$ and $J\in \exch_-[Y]$. We put
\[I\blackdiamond J=m(I\diamond J)\subseteq\comp[X\sqcup Y].\] 
\end{defi}

\begin{theo}\label{theo8.1.14}
$(\exch_-,\blackdiamond,\Delta)$ is a twisted infinitesimal bialgebra.
\end{theo}

\begin{proof}
\textit{First step}. 
Let $X$, $Y$ be disjoint finite sets. We already proved that $I\diamond J\in \exch[X\sqcup Y]$, so $m(I\diamond J)\in \exch_-[X\sqcup Y]$. If $X$ is a finite set and $I\in \exch_-[X]$,
\[I\blackdiamond 1=m(I\diamond 1)=m(I)=I,\]
as $I\in \exch_-[I]$. Similarly, $1\blackdiamond I=I$, so $1$ is a unit for $\blackdiamond$. 
Let $X$, $Y,Z$ be pairwise disjoint finite sets, $I\in \exch_-[X]$, $J\in \exch_-[Y]$, $K\in \exch_-[Z]$. Let $\pi\in \comp[X\sqcup Y\sqcup Z]$.
\begin{align*}
\pi\in (I\blackdiamond J)\blackdiamond K&\Longleftrightarrow \exists \varpi_1\in I\blackdiamond J,\: \exists \pi_3\in K,\: \pi\leq (\varpi_1,\pi_3)\\
&\Longleftrightarrow \exists (\pi_1,\pi_2,\pi_3)\in I\times J\times K,\:\varpi_1\in \comp[X\sqcup Y],\: \varpi_1\leq (\pi_1,\pi_2),\: \pi\leq (\varpi_1,\pi_3)\\
&\Longleftrightarrow \exists (\pi_1,\pi_2,\pi_3)\in I\times J\times K,\: \pi\leq (\pi_1,\pi_2,\pi_3)\\
&\Longleftrightarrow \exists (\pi_1,\pi_2,\pi_3)\in I\times J\times K,\:\varpi_2\in \comp[Y\sqcup Z],\: \varpi_2\leq (\pi_2,\pi_3),\:\pi\leq (\pi_1,\varpi_2)\\
&\Longleftrightarrow \exists \pi_1\in I,\: \exists \varpi_2\in J\blackdiamond J,\: \pi\leq (\pi_1,\varpi_2)\\
&\Longleftrightarrow \pi\in I\blackdiamond(J\blackdiamond K).
\end{align*}
So $\blackdiamond$ is associative.

\textit{Second step}. Let us prove the infinitesimal compatibility. Let $X$, $Y,X',Y'$ be finite sets, such that $X\sqcup Y=X'\sqcup Y'$, $I\in \exch_-[X]$ and $J\in \exch_-[Y]$. Let us assume that $\Delta_{X',Y'}(I\blackdiamond J)\neq 0$.
Let $\pi'\in \comp[X']$, $\pi''\in \comp[X'']$, such that $(\pi',\pi'')\in I\blackdiamond J$. Therefore, there exists $\varpi\in I\diamond J$, such that $(\pi',\pi'')\leq \varpi$. By definition of $I\diamond J$, there exists $\varpi'\in I$, $\varpi''\in J$, such that $(\pi',\pi'')\leq \varpi=(\varpi',\varpi'')$. By definition of the order, $X'\subseteq X$ or $X\subseteq X'$. 
Dually, if $X$ is not included in $X'$, nor $X'$ is included in $X$,
then $\Delta_{X',Y'}(I\blackdiamond J)=0$.

Let us assume that $X\subseteq X'$. Let $\pi\in (I\blackdiamond J)_{\mid_LX'}$, and $\pi'\in (I\blackdiamond J)_{\mid_RY'}$. Then $(\pi,\pi')\in I\blackdiamond J$, by the exchange condition. By definition of $I\blackdiamond J$, $\pi\leq(\varpi,\varpi')$, with $\varpi\in I$, and $(\varpi',\pi')\in J$.
Therefore, $\varpi'\in J_{\mid_L Y\setminus Y'}$, and $\pi\in I\blackdiamond J_{\mid_LY\setminus Y'}$, and $\pi' \in J_{\mid_R Y'}$. 
We proved that $(I\blackdiamond J)_{\mid_LX'}\subseteq I\blackdiamond J_{\mid_LY\setminus Y'}$ and $(I\blackdiamond J)_{\mid_RY'}\subseteq J_{\mid_R Y'}$. 

Conversely, let $\pi\in I\blackdiamond J_{\mid_LY\setminus Y'}$ and $\pi'\in J_{\mid_R Y'}$.
Then $\pi\leq (\pi_1,\pi_2)$, with $\pi_1\in I$ and $\pi_2\in J_{\mid_LY\setminus Y'}$. By the exchange condition, $(\pi_2,\pi')\in J$. So $(\pi,\pi')\leq (\pi_1,\pi_2,\pi')\in I\diamond J$, so
$(\pi,\pi')\in I\blackdiamond J$, which finally gives $\pi\in (I\blackdiamond J)_{\mid_LX'}$ and $\pi'\in (I\blackdiamond J)_{\mid_RY'}$. We proved that $(I\blackdiamond J)_{\mid_LX'}= I\blackdiamond J_{\mid_LY\setminus Y'}$ and $(I\blackdiamond J)_{\mid_RY'}= J_{\mid_R Y'}$. Therefore,
\[\Delta_{X',Y'}(I\blackdiamond J)= I\blackdiamond J_{\mid_LY\setminus Y'}\otimes J_{\mid_R Y'}=(I\otimes 1)\blackdiamond \Delta_{Y\setminus Y',Y'}(J).\]

Similarly, if $X'\subseteq X$, we obtain
\[\Delta_{X',Y'}(I\blackdiamond J)=\Delta_{X',X\setminus X'}(I)\blackdiamond(1\otimes J).\]
So $(\exch_-,\blackdiamond,\Delta)$ is a twisted infinitesimal bialgebra.. \end{proof}

Combining Theorems \ref{theo8.1.5} and \ref{theo8.1.9}: 

\begin{cor}
As a twisted coalgebra, $\exch_-$ is isomorphic to $(T(\prim(\exch_-)),\shuffle,\Delta_\conc)$.
\end{cor}

\bibliographystyle{amsplain}
\bibliography{biblio}

 \end{document}

%% file: arbres.tex
\newcommand{\tun}{\begin{tikzpicture}[line cap=round,line join=round,>=triangle 45,x=0.5cm,y=0.5cm]
\clip(-0.2,-0.1) rectangle (0.2,0.2);
\begin{scriptsize}
\draw [fill=black] (0.,0.) circle (1pt);
\end{scriptsize}
\end{tikzpicture}}

\newcommand{\tdeux}{\begin{tikzpicture}[line cap=round,line join=round,>=triangle 45,x=0.5cm,y=0.5cm]
\clip(-.2,-.1) rectangle (0.2,0.7);
\draw [line width=.5pt] (0.,0.5)-- (0.,0.);
\begin{scriptsize}
\draw [fill=black] (0.,0.) circle (1pt);
\draw [fill=black] (0.,0.5) circle (1pt);
\end{scriptsize}
\end{tikzpicture}}

\newcommand{\ttroisun}{\begin{tikzpicture}[line cap=round,line join=round,>=triangle 45,x=0.5cm,y=0.5cm]
\clip(-0.5,-0.1) rectangle (0.5,0.7);
\draw [line width=0.5pt] (0.,0.)-- (-0.3,0.5);
\draw [line width=0.5pt] (0.,0.)-- (0.3,0.5);
\begin{scriptsize}
\draw [fill=black] (-0.3,0.5) circle (1pt);
\draw [fill=black] (0.,0.) circle (1pt);
\draw [fill=black] (0.3,0.5) circle (1pt);
\end{scriptsize}
\end{tikzpicture}}
\newcommand{\ttroisdeux}{\begin{tikzpicture}[line cap=round,line join=round,>=triangle 45,x=0.5cm,y=0.5cm]
\clip(-.2,-.1) rectangle (0.2,1.2);
\draw [line width=0.5pt] (0.,0.5)-- (0.,0.);
\draw [line width=0.5pt] (0.,0.5)-- (0.,1.);
\begin{scriptsize}
\draw [fill=black] (0.,0.) circle (1pt);
\draw [fill=black] (0.,0.5) circle (1pt);
\draw [fill=black] (0.,1.) circle (1pt);
\end{scriptsize}
\end{tikzpicture}}

\newcommand{\tquatreun}{\begin{tikzpicture}[line cap=round,line join=round,>=triangle 45,x=0.5cm,y=0.5cm]
\clip(-0.5,-0.1) rectangle (0.5,0.7);
\draw [line width=0.5pt] (0.,0.)-- (-0.3,0.5);
\draw [line width=0.5pt] (0.,0.)-- (0.3,0.5);
\draw [line width=0.5pt] (0.,0.)-- (0.,0.5);
\begin{scriptsize}
\draw [fill=black] (-0.3,0.5) circle (1.0pt);
\draw [fill=black] (0.,0.) circle (1.0pt);
\draw [fill=black] (0.3,0.5) circle (1.0pt);
\draw [fill=black] (0.,0.5) circle (1.0pt);
\end{scriptsize}
\end{tikzpicture}}
\newcommand{\tquatredeux}{\begin{tikzpicture}[line cap=round,line join=round,>=triangle 45,x=0.5cm,y=0.5cm]
\clip(-0.5,-0.1) rectangle (0.5,1.2);
\draw [line width=0.5pt] (0.,0.)-- (-0.3,0.5);
\draw [line width=0.5pt] (0.,0.)-- (0.3,0.5);
\draw [line width=0.5pt] (-0.3,0.5)-- (-0.3,1.);
\begin{scriptsize}
\draw [fill=black] (-0.3,0.5) circle (1.0pt);
\draw [fill=black] (0.,0.) circle (1.0pt);
\draw [fill=black] (0.3,0.5) circle (1.0pt);
\draw [fill=black] (-0.3,1.) circle (1.0pt);
\end{scriptsize}
\end{tikzpicture}}

\newcommand{\tquatrequatre}{\begin{tikzpicture}[line cap=round,line join=round,>=triangle 45,x=0.5cm,y=0.5cm]
\clip(-0.5,-0.1) rectangle (0.5,1.2);
\draw [line width=0.5pt] (0.,0.)-- (0.,0.5);
\draw [line width=0.5pt] (0.,0.5)-- (0.3,1.);
\draw [line width=0.5pt] (0.,0.5)-- (-0.3,1.);
\begin{scriptsize}
\draw [fill=black] (0.,0.) circle (1.0pt);
\draw [fill=black] (0.,0.5) circle (1.0pt);
\draw [fill=black] (-0.3,1.) circle (1.0pt);
\draw [fill=black] (0.3,1.) circle (1.0pt);
\end{scriptsize}
\end{tikzpicture}}
\newcommand{\tquatrecinq}{\begin{tikzpicture}[line cap=round,line join=round,>=triangle 45,x=0.5cm,y=0.5cm]
\clip(-.2,-.1) rectangle (0.2,1.7);
\draw [line width=0.5pt] (0.,0.)-- (0.,0.5);
\draw [line width=0.5pt] (0.,0.5)-- (0.,1.);
\draw [line width=0.5pt] (0.,1.)-- (0.,1.5);
\begin{scriptsize}
\draw [fill=black] (0.,0.) circle (1.0pt);
\draw [fill=black] (0.,0.5) circle (1.0pt);
\draw [fill=black] (0.,1.) circle (1.0pt);
\draw [fill=black] (0.,1.5) circle (1.0pt);
\end{scriptsize}
\end{tikzpicture}}

\newcommand{\tcinqdeux}{\begin{tikzpicture}[line cap=round,line join=round,>=triangle 45,x=0.5cm,y=0.5cm]
\clip(-0.5,-0.1) rectangle (0.5,1.2);
\draw [line width=0.5pt] (0.,0.)-- (-0.3,0.5);
\draw [line width=0.5pt] (0.,0.)-- (0.3,0.5);
\draw [line width=0.5pt] (0.,0.)-- (0.,0.5);
\draw [line width=0.5pt] (-0.3,0.5)-- (-0.3,1.);
\begin{scriptsize}
\draw [fill=black] (-0.3,0.5) circle (1.0pt);
\draw [fill=black] (0.,0.) circle (1.0pt);
\draw [fill=black] (0.,0.5) circle (1.0pt);
\draw [fill=black] (0.3,0.5) circle (1.0pt);
\draw [fill=black] (-0.3,1.) circle (1.0pt);
\end{scriptsize}
\end{tikzpicture}}

\newcommand{\tcinqcinq}{\begin{tikzpicture}[line cap=round,line join=round,>=triangle 45,x=0.5cm,y=0.5cm]
\clip(-0.5,-0.1) rectangle (0.5,1.2);
\draw [line width=0.5pt] (0.,0.)-- (-0.3,0.5);
\draw [line width=0.5pt] (-0.3,0.5)-- (-0.3,1.);
\draw [line width=0.5pt] (0.,0.)-- (0.3,0.5);
\draw [line width=0.5pt] (0.3,0.5)-- (0.3,1.);
\begin{scriptsize}
\draw [fill=black] (-0.3,0.5) circle (1.0pt);
\draw [fill=black] (0.,0.) circle (1.0pt);
\draw [fill=black] (-0.3,1.) circle (1.0pt);
\draw [fill=black] (0.3,0.5) circle (1.0pt);
\draw [fill=black] (0.3,1.) circle (1.0pt);
\end{scriptsize}
\end{tikzpicture}}
\newcommand{\tcinqsix}{\begin{tikzpicture}[line cap=round,line join=round,>=triangle 45,x=0.5cm,y=0.5cm]
\clip(-0.7,-0.1) rectangle (0.5,1.2);
\draw [line width=0.5pt] (0.,0.)-- (-0.3,0.5);
\draw [line width=0.5pt] (0.,0.)-- (0.3,0.5);
\draw [line width=0.5pt] (-0.3,0.5)-- (-0.6,1.);
\draw [line width=0.5pt] (-0.3,0.5)-- (0.,1.);
\begin{scriptsize}
\draw [fill=black] (-0.3,0.5) circle (1.0pt);
\draw [fill=black] (0.,0.) circle (1.0pt);
\draw [fill=black] (0.3,0.5) circle (1.0pt);
\draw [fill=black] (-0.6,1.) circle (1.0pt);
\draw [fill=black] (0.,1.) circle (1.0pt);
\end{scriptsize}
\end{tikzpicture}}

\newcommand{\tcinqhuit}{\begin{tikzpicture}[line cap=round,line join=round,>=triangle 45,x=0.5cm,y=0.5cm]
\clip(-0.5,-0.1) rectangle (0.5,1.7);
\draw [line width=0.5pt] (0.,0.)-- (-0.3,0.5);
\draw [line width=0.5pt] (-0.3,0.5)-- (-0.3,1.);
\draw [line width=0.5pt] (0.,0.)-- (0.3,0.5);
\draw [line width=0.5pt] (-0.3,1.)-- (-0.3,1.5);
\begin{scriptsize}
\draw [fill=black] (-0.3,0.5) circle (1.0pt);
\draw [fill=black] (0.,0.) circle (1.0pt);
\draw [fill=black] (-0.3,1.) circle (1.0pt);
\draw [fill=black] (0.3,0.5) circle (1.0pt);
\draw [fill=black] (-0.3,1.5) circle (1.0pt);
\end{scriptsize}
\end{tikzpicture}}

\newcommand{\tcinqdix}{\begin{tikzpicture}[line cap=round,line join=round,>=triangle 45,x=0.5cm,y=0.5cm]
\clip(-0.5,-0.1) rectangle (0.5,1.7);
\draw [line width=0.5pt] (0.,0.)-- (0.,0.5);
\draw [line width=0.5pt] (0.,0.5)-- (0.3,1.);
\draw [line width=0.5pt] (0.,0.5)-- (-0.3,1.);
\draw [line width=0.5pt] (0.,0.5)-- (0.,1.);
\begin{scriptsize}
\draw [fill=black] (0.,0.) circle (1.0pt);
\draw [fill=black] (0.,0.5) circle (1.0pt);
\draw [fill=black] (-0.3,1.) circle (1.0pt);
\draw [fill=black] (0.3,1.) circle (1.0pt);
\draw [fill=black] (0.,1.) circle (1.0pt);
\end{scriptsize}
\end{tikzpicture}}
\newcommand{\tcinqonze}{\begin{tikzpicture}[line cap=round,line join=round,>=triangle 45,x=0.5cm,y=0.5cm]
\clip(-0.5,-0.1) rectangle (0.5,1.7);
\draw [line width=0.5pt] (0.,0.)-- (0.,0.5);
\draw [line width=0.5pt] (0.,0.5)-- (0.3,1.);
\draw [line width=0.5pt] (0.,0.5)-- (-0.3,1.);
\draw [line width=0.5pt] (-0.3,1.)-- (-0.3,1.5);
\begin{scriptsize}
\draw [fill=black] (0.,0.) circle (1.0pt);
\draw [fill=black] (0.,0.5) circle (1.0pt);
\draw [fill=black] (-0.3,1.) circle (1.0pt);
\draw [fill=black] (0.3,1.) circle (1.0pt);
\draw [fill=black] (-0.3,1.5) circle (1.0pt);
\end{scriptsize}
\end{tikzpicture}}

\newcommand{\tcinqtreize}{\begin{tikzpicture}[line cap=round,line join=round,>=triangle 45,x=0.5cm,y=0.5cm]
\clip(-0.5,-0.1) rectangle (0.5,1.7);
\draw [line width=0.5pt] (0.,0.)-- (0.,0.5);
\draw [line width=0.5pt] (0.,0.5)-- (0.,1.);
\draw [line width=0.5pt] (0.,1.)-- (-0.3,1.5);
\draw [line width=0.5pt] (0.,1.)-- (0.3,1.5);
\begin{scriptsize}
\draw [fill=black] (0.,0.) circle (1.0pt);
\draw [fill=black] (0.,0.5) circle (1.0pt);
\draw [fill=black] (0.,1.) circle (1.0pt);
\draw [fill=black] (-0.3,1.5) circle (1.0pt);
\draw [fill=black] (0.3,1.5) circle (1.0pt);
\end{scriptsize}
\end{tikzpicture}}

\newcommand{\tdun}[1]{\begin{tikzpicture}[line cap=round,line join=round,>=triangle 45,x=0.5cm,y=0.5cm]
\clip(-.2,-.1) rectangle (0.5,0.5);
\begin{scriptsize}
\draw [fill=black] (0.,0.) circle (1pt);
\end{scriptsize}
\draw(0.3,0.1) node {\tiny #1};
\end{tikzpicture}}

\newcommand{\tddeux}[2]{\begin{tikzpicture}[line cap=round,line join=round,>=triangle 45,x=0.5cm,y=0.5cm]
\clip(-.2,-.1) rectangle (0.5,1.);
\draw [line width=.5pt] (0.,0.5)-- (0.,0.);
\begin{scriptsize}
\draw [fill=black] (0.,0.) circle (1pt);
\draw [fill=black] (0.,0.5) circle (1pt);
\end{scriptsize}
\draw(0.3,0.1) node {\tiny #1};
\draw(0.3,0.6) node {\tiny #2};
\end{tikzpicture}}

\newcommand{\tdtroisun}[3]{\begin{tikzpicture}[line cap=round,line join=round,>=triangle 45,x=0.5cm,y=0.5cm]
\clip(-0.5,-0.1) rectangle (0.5,1.2);
\draw [line width=0.5pt] (0.,0.)-- (-0.3,0.5);
\draw [line width=0.5pt] (0.,0.)-- (0.3,0.5);
\begin{scriptsize}
\draw [fill=black] (-0.3,0.5) circle (1pt);
\draw [fill=black] (0.,0.) circle (1pt);
\draw [fill=black] (0.3,0.5) circle (1pt);
\end{scriptsize}
\draw(0.35,0.1) node {\tiny #1};
\draw(0.3,0.8) node {\tiny #2};
\draw(-0.3,0.8) node {\tiny #3};
\end{tikzpicture}}
\newcommand{\tdtroisdeux}[3]{\begin{tikzpicture}[line cap=round,line join=round,>=triangle 45,x=0.5cm,y=0.5cm]
\clip(-.2,-.1) rectangle (0.5,1.5);
\draw [line width=0.5pt] (0.,0.5)-- (0.,0.);
\draw [line width=0.5pt] (0.,0.5)-- (0.,1.);
\begin{scriptsize}
\draw [fill=black] (0.,0.) circle (1pt);
\draw [fill=black] (0.,0.5) circle (1pt);
\draw [fill=black] (0.,1.) circle (1pt);
\end{scriptsize}
\draw(0.3,0.1) node {\tiny #1};
\draw(0.3,0.6) node {\tiny #2};
\draw(0.3,1.1) node {\tiny #3};
\end{tikzpicture}}

\newcommand{\tdquatreun}[4]{\begin{tikzpicture}[line cap=round,line join=round,>=triangle 45,x=0.5cm,y=0.5cm]
\clip(-0.5,-0.1) rectangle (0.5,1.2);
\draw [line width=0.5pt] (0.,0.)-- (-0.3,0.5);
\draw [line width=0.5pt] (0.,0.)-- (0.3,0.5);
\draw [line width=0.5pt] (0.,0.)-- (0.,0.5);
\begin{scriptsize}
\draw [fill=black] (-0.3,0.5) circle (1.0pt);
\draw [fill=black] (0.,0.) circle (1.0pt);
\draw [fill=black] (0.3,0.5) circle (1.0pt);
\draw [fill=black] (0.,0.5) circle (1.0pt);
\end{scriptsize}
\draw(0.35,0.1) node {\tiny #1};
\draw(0.3,0.8) node {\tiny #2};
\draw(0.,0.8) node {\tiny #3};
\draw(-0.3,0.8) node {\tiny #4};
\end{tikzpicture}}
\newcommand{\tdquatredeux}[4]{\begin{tikzpicture}[line cap=round,line join=round,>=triangle 45,x=0.5cm,y=0.5cm]
\clip(-0.7,-0.1) rectangle (0.5,1.7);
\draw [line width=0.5pt] (0.,0.)-- (-0.3,0.5);
\draw [line width=0.5pt] (0.,0.)-- (0.3,0.5);
\draw [line width=0.5pt] (-0.3,0.5)-- (-0.3,1.);
\begin{scriptsize}
\draw [fill=black] (-0.3,0.5) circle (1.0pt);
\draw [fill=black] (0.,0.) circle (1.0pt);
\draw [fill=black] (0.3,0.5) circle (1.0pt);
\draw [fill=black] (-0.3,1.) circle (1.0pt);
\end{scriptsize}
\draw(0.35,0.1) node {\tiny #1};
\draw(0.3,0.8) node {\tiny #2};
\draw(-0.6,0.6) node {\tiny #3};
\draw(-0.3,1.3) node {\tiny #4};
\end{tikzpicture}}

\newcommand{\tdquatrequatre}[4]{\begin{tikzpicture}[line cap=round,line join=round,>=triangle 45,x=0.5cm,y=0.5cm]
\clip(-0.5,-0.1) rectangle (0.5,1.7);
\draw [line width=0.5pt] (0.,0.)-- (0.,0.5);
\draw [line width=0.5pt] (0.,0.5)-- (0.3,1.);
\draw [line width=0.5pt] (0.,0.5)-- (-0.3,1.);
\begin{scriptsize}
\draw [fill=black] (0.,0.) circle (1.0pt);
\draw [fill=black] (0.,0.5) circle (1.0pt);
\draw [fill=black] (-0.3,1.) circle (1.0pt);
\draw [fill=black] (0.3,1.) circle (1.0pt);
\end{scriptsize}
\draw(0.3,0.1) node {\tiny #1};
\draw(0.3,0.5) node {\tiny #2};
\draw(0.3,1.3) node {\tiny #3};
\draw(-0.3,1.3) node {\tiny #4};
\end{tikzpicture}}
\newcommand{\tdquatrecinq}[4]{\begin{tikzpicture}[line cap=round,line join=round,>=triangle 45,x=0.5cm,y=0.5cm]
\clip(-.2,-.1) rectangle (0.5,2.);
\draw [line width=0.5pt] (0.,0.)-- (0.,0.5);
\draw [line width=0.5pt] (0.,0.5)-- (0.,1.);
\draw [line width=0.5pt] (0.,1.)-- (0.,1.5);
\begin{scriptsize}
\draw [fill=black] (0.,0.) circle (1.0pt);
\draw [fill=black] (0.,0.5) circle (1.0pt);
\draw [fill=black] (0.,1.) circle (1.0pt);
\draw [fill=black] (0.,1.5) circle (1.0pt);
\end{scriptsize}
\draw(0.3,0.1) node {\tiny #1};
\draw(0.3,0.6) node {\tiny #2};
\draw(0.3,1.1) node {\tiny #3};
\draw(0.3,1.6) node {\tiny #4};
\end{tikzpicture}}
\newcommand{\tdcinqsix}[5]{\begin{tikzpicture}[line cap=round,line join=round,>=triangle 45,x=0.5cm,y=0.5cm]
\clip(-0.7,-0.1) rectangle (0.5,1.7);
\draw [line width=0.5pt] (0.,0.)-- (-0.3,0.5);
\draw [line width=0.5pt] (0.,0.)-- (0.3,0.5);
\draw [line width=0.5pt] (-0.3,0.5)-- (-0.6,1.);
\draw [line width=0.5pt] (-0.3,0.5)-- (0.,1.);
\begin{scriptsize}
\draw [fill=black] (-0.3,0.5) circle (1.0pt);
\draw [fill=black] (0.,0.) circle (1.0pt);
\draw [fill=black] (0.3,0.5) circle (1.0pt);
\draw [fill=black] (-0.6,1.) circle (1.0pt);
\draw [fill=black] (0.,1.) circle (1.0pt);
\end{scriptsize}
\draw(0.35,0.) node {\tiny #1};
\draw(0.3,0.8) node {\tiny #2};
\draw(-0.6,0.5) node {\tiny #3};
\draw(0.,1.3) node {\tiny #4};
\draw(-0.6,1.3) node {\tiny #5};
\end{tikzpicture}}

\newcommand{\tdcinqsixA}[5]{\begin{tikzpicture}[line cap=round,line join=round,>=triangle 45,x=0.5cm,y=0.5cm]
\clip(-0.7,-0.1) rectangle (0.5,1.7);
\draw [line width=0.5pt,color=red] (0.,0.)-- (-0.3,0.5);
\draw [line width=0.5pt,color=red] (0.,0.)-- (0.3,0.5);
\draw [line width=0.5pt,color=red] (-0.3,0.5)-- (-0.6,1.);
\draw [line width=0.5pt,color=red] (-0.3,0.5)-- (0.,1.);
\begin{scriptsize}
\draw [fill=red,color=red] (-0.3,0.5) circle (1.0pt);
\draw [fill=red,color=red] (0.,0.) circle (1.0pt);
\draw [fill=red,color=red] (0.3,0.5) circle (1.0pt);
\draw [fill=red,color=red] (-0.6,1.) circle (1.0pt);
\draw [fill=red,color=red] (0.,1.) circle (1.0pt);
\end{scriptsize}
\draw(0.35,0.) node {\tiny #1};
\draw(0.3,0.8) node {\tiny #2};
\draw(-0.6,0.5) node {\tiny #3};
\draw(0.,1.3) node {\tiny #4};
\draw(-0.6,1.3) node {\tiny #5};
\end{tikzpicture}}

\newcommand{\tdcinqsixB}[5]{\begin{tikzpicture}[line cap=round,line join=round,>=triangle 45,x=0.5cm,y=0.5cm]
\clip(-0.7,-0.1) rectangle (0.5,1.7);
\draw [line width=0.5pt,color=red] (0.,0.)-- (-0.3,0.5);
\draw [line width=0.5pt,color=red] (0.,0.)-- (0.3,0.5);
\draw [line width=0.5pt,color=black] (-0.3,0.5)-- (-0.6,1.);
\draw [line width=0.5pt,color=red] (-0.3,0.5)-- (0.,1.);
\begin{scriptsize}
\draw [fill=red,color=red] (-0.3,0.5) circle (1.0pt);
\draw [fill=red,color=red] (0.,0.) circle (1.0pt);
\draw [fill=red,color=red] (0.3,0.5) circle (1.0pt);
\draw [fill=red,color=blue] (-0.6,1.) circle (1.0pt);
\draw [fill=red,color=red] (0.,1.) circle (1.0pt);
\end{scriptsize}
\draw(0.35,0.) node {\tiny #1};
\draw(0.3,0.8) node {\tiny #2};
\draw(-0.6,0.5) node {\tiny #3};
\draw(0.,1.3) node {\tiny #4};
\draw(-0.6,1.3) node {\tiny #5};
\end{tikzpicture}}

\newcommand{\tdcinqsixC}[5]{\begin{tikzpicture}[line cap=round,line join=round,>=triangle 45,x=0.5cm,y=0.5cm]
\clip(-0.7,-0.1) rectangle (0.5,1.7);
\draw [line width=0.5pt,color=red] (0.,0.)-- (-0.3,0.5);
\draw [line width=0.5pt,color=red] (0.,0.)-- (0.3,0.5);
\draw [line width=0.5pt,color=red] (-0.3,0.5)-- (-0.6,1.);
\draw [line width=0.5pt,color=black] (-0.3,0.5)-- (0.,1.);
\begin{scriptsize}
\draw [fill=red,color=red] (-0.3,0.5) circle (1.0pt);
\draw [fill=red,color=red] (0.,0.) circle (1.0pt);
\draw [fill=red,color=red] (0.3,0.5) circle (1.0pt);
\draw [fill=red,color=red] (-0.6,1.) circle (1.0pt);
\draw [fill=red,color=blue] (0.,1.) circle (1.0pt);
\end{scriptsize}
\draw(0.35,0.) node {\tiny #1};
\draw(0.3,0.8) node {\tiny #2};
\draw(-0.6,0.5) node {\tiny #3};
\draw(0.,1.3) node {\tiny #4};
\draw(-0.6,1.3) node {\tiny #5};
\end{tikzpicture}}

\newcommand{\tdcinqsixD}[5]{\begin{tikzpicture}[line cap=round,line join=round,>=triangle 45,x=0.5cm,y=0.5cm]
\clip(-0.7,-0.1) rectangle (0.5,1.7);
\draw [line width=0.5pt,color=red] (0.,0.)-- (-0.3,0.5);
\draw [line width=0.5pt,color=black] (0.,0.)-- (0.3,0.5);
\draw [line width=0.5pt,color=red] (-0.3,0.5)-- (-0.6,1.);
\draw [line width=0.5pt,color=red] (-0.3,0.5)-- (0.,1.);
\begin{scriptsize}
\draw [fill=red,color=red] (-0.3,0.5) circle (1.0pt);
\draw [fill=red,color=red] (0.,0.) circle (1.0pt);
\draw [fill=red,color=blue] (0.3,0.5) circle (1.0pt);
\draw [fill=red,color=red] (-0.6,1.) circle (1.0pt);
\draw [fill=red,color=red] (0.,1.) circle (1.0pt);
\end{scriptsize}
\draw(0.35,0.) node {\tiny #1};
\draw(0.3,0.8) node {\tiny #2};
\draw(-0.6,0.5) node {\tiny #3};
\draw(0.,1.3) node {\tiny #4};
\draw(-0.6,1.3) node {\tiny #5};
\end{tikzpicture}}

\newcommand{\tdcinqsixE}[5]{\begin{tikzpicture}[line cap=round,line join=round,>=triangle 45,x=0.5cm,y=0.5cm]
\clip(-0.7,-0.1) rectangle (0.5,1.7);
\draw [line width=0.5pt,color=black] (0.,0.)-- (-0.3,0.5);
\draw [line width=0.5pt,color=blue] (0.,0.)-- (0.3,0.5);
\draw [line width=0.5pt,color=red] (-0.3,0.5)-- (-0.6,1.);
\draw [line width=0.5pt,color=red] (-0.3,0.5)-- (0.,1.);
\begin{scriptsize}
\draw [fill=red,color=red] (-0.3,0.5) circle (1.0pt);
\draw [fill=red,color=blue] (0.,0.) circle (1.0pt);
\draw [fill=red,color=blue] (0.3,0.5) circle (1.0pt);
\draw [fill=red,color=red] (-0.6,1.) circle (1.0pt);
\draw [fill=red,color=red] (0.,1.) circle (1.0pt);
\end{scriptsize}
\draw(0.35,0.) node {\tiny #1};
\draw(0.3,0.8) node {\tiny #2};
\draw(-0.6,0.5) node {\tiny #3};
\draw(0.,1.3) node {\tiny #4};
\draw(-0.6,1.3) node {\tiny #5};
\end{tikzpicture}}

\newcommand{\tdcinqsixF}[5]{\begin{tikzpicture}[line cap=round,line join=round,>=triangle 45,x=0.5cm,y=0.5cm]
\clip(-0.7,-0.1) rectangle (0.5,1.7);
\draw [line width=0.5pt,color=red] (0.,0.)-- (-0.3,0.5);
\draw [line width=0.5pt,color=black] (0.,0.)-- (0.3,0.5);
\draw [line width=0.5pt,color=red] (-0.3,0.5)-- (-0.6,1.);
\draw [line width=0.5pt,color=black] (-0.3,0.5)-- (0.,1.);
\begin{scriptsize}
\draw [fill=red,color=red] (-0.3,0.5) circle (1.0pt);
\draw [fill=red,color=red] (0.,0.) circle (1.0pt);
\draw [fill=red,color=blue] (0.3,0.5) circle (1.0pt);
\draw [fill=red,color=red] (-0.6,1.) circle (1.0pt);
\draw [fill=red,color=green] (0.,1.) circle (1.0pt);
\end{scriptsize}
\draw(0.35,0.) node {\tiny #1};
\draw(0.3,0.8) node {\tiny #2};
\draw(-0.6,0.5) node {\tiny #3};
\draw(0.,1.3) node {\tiny #4};
\draw(-0.6,1.3) node {\tiny #5};
\end{tikzpicture}}

\newcommand{\tdcinqsixG}[5]{\begin{tikzpicture}[line cap=round,line join=round,>=triangle 45,x=0.5cm,y=0.5cm]
\clip(-0.7,-0.1) rectangle (0.5,1.7);
\draw [line width=0.5pt,color=green] (0.,0.)-- (-0.3,0.5);
\draw [line width=0.5pt,color=black] (0.,0.)-- (0.3,0.5);
\draw [line width=0.5pt,color=black] (-0.3,0.5)-- (-0.6,1.);
\draw [line width=0.5pt,color=green] (-0.3,0.5)-- (0.,1.);
\begin{scriptsize}
\draw [fill=red,color=green] (-0.3,0.5) circle (1.0pt);
\draw [fill=red,color=green] (0.,0.) circle (1.0pt);
\draw [fill=red,color=blue] (0.3,0.5) circle (1.0pt);
\draw [fill=red,color=red] (-0.6,1.) circle (1.0pt);
\draw [fill=red,color=green] (0.,1.) circle (1.0pt);
\end{scriptsize}
\draw(0.35,0.) node {\tiny #1};
\draw(0.3,0.8) node {\tiny #2};
\draw(-0.6,0.5) node {\tiny #3};
\draw(0.,1.3) node {\tiny #4};
\draw(-0.6,1.3) node {\tiny #5};
\end{tikzpicture}}

\newcommand{\tdcinqsixH}[5]{\begin{tikzpicture}[line cap=round,line join=round,>=triangle 45,x=0.5cm,y=0.5cm]
\clip(-0.7,-0.1) rectangle (0.5,1.7);
\draw [line width=0.5pt,color=black] (0.,0.)-- (-0.3,0.5);
\draw [line width=0.5pt,color=blue] (0.,0.)-- (0.3,0.5);
\draw [line width=0.5pt,color=red] (-0.3,0.5)-- (-0.6,1.);
\draw [line width=0.5pt,color=black] (-0.3,0.5)-- (0.,1.);
\begin{scriptsize}
\draw [fill=red,color=red] (-0.3,0.5) circle (1.0pt);
\draw [fill=red,color=blue] (0.,0.) circle (1.0pt);
\draw [fill=red,color=blue] (0.3,0.5) circle (1.0pt);
\draw [fill=red,color=red] (-0.6,1.) circle (1.0pt);
\draw [fill=red,color=green] (0.,1.) circle (1.0pt);
\end{scriptsize}
\draw(0.35,0.) node {\tiny #1};
\draw(0.3,0.8) node {\tiny #2};
\draw(-0.6,0.5) node {\tiny #3};
\draw(0.,1.3) node {\tiny #4};
\draw(-0.6,1.3) node {\tiny #5};
\end{tikzpicture}}

\newcommand{\tdcinqsixI}[5]{\begin{tikzpicture}[line cap=round,line join=round,>=triangle 45,x=0.5cm,y=0.5cm]
\clip(-0.7,-0.1) rectangle (0.5,1.7);
\draw [line width=0.5pt,color=black] (0.,0.)-- (-0.3,0.5);
\draw [line width=0.5pt,color=blue] (0.,0.)-- (0.3,0.5);
\draw [line width=0.5pt,color=black] (-0.3,0.5)-- (-0.6,1.);
\draw [line width=0.5pt,color=green] (-0.3,0.5)-- (0.,1.);
\begin{scriptsize}
\draw [fill=red,color=green] (-0.3,0.5) circle (1.0pt);
\draw [fill=red,color=blue] (0.,0.) circle (1.0pt);
\draw [fill=red,color=blue] (0.3,0.5) circle (1.0pt);
\draw [fill=red,color=red] (-0.6,1.) circle (1.0pt);
\draw [fill=red,color=green] (0.,1.) circle (1.0pt);
\end{scriptsize}
\draw(0.35,0.) node {\tiny #1};
\draw(0.3,0.8) node {\tiny #2};
\draw(-0.6,0.5) node {\tiny #3};
\draw(0.,1.3) node {\tiny #4};
\draw(-0.6,1.3) node {\tiny #5};
\end{tikzpicture}}

\newcommand{\tddun}[1]{\begin{tikzpicture}[line cap=round,line join=round,>=triangle 45,x=0.5cm,y=0.5cm]
\clip(-.2,-.1) rectangle (0.7,0.7);
\begin{scriptsize}
\draw [fill=black] (0.,0.) circle (1pt);
\end{scriptsize}
\draw(0.4,0.1) node {\tiny #1};
\end{tikzpicture}}

\newcommand{\tdddun}[1]{\begin{tikzpicture}[line cap=round,line join=round,>=triangle 45,x=0.5cm,y=0.5cm]
\clip(-.2,-.1) rectangle (1.0,0.7);
\begin{scriptsize}
\draw [fill=black] (0.,0.) circle (1pt);
\end{scriptsize}
\draw(0.6,0.1) node {\tiny #1};
\end{tikzpicture}}

\newcommand{\tdddeux}[2]{\begin{tikzpicture}[line cap=round,line join=round,>=triangle 45,x=0.5cm,y=0.5cm]
\clip(-.2,-.1) rectangle (0.7,1.);
\draw [line width=.5pt] (0.,0.5)-- (0.,0.);
\begin{scriptsize}
\draw [fill=black] (0.,0.) circle (1pt);
\draw [fill=black] (0.,0.5) circle (1pt);
\end{scriptsize}
\draw(0.4,0.1) node {\tiny #1};
\draw(0.3,0.6) node {\tiny #2};
\end{tikzpicture}}

\newcommand{\tddddeux}[2]{\begin{tikzpicture}[line cap=round,line join=round,>=triangle 45,x=0.5cm,y=0.5cm]
\clip(-.2,-.1) rectangle (0.7,1.);
\draw [line width=.5pt] (0.,0.5)-- (0.,0.);
\begin{scriptsize}
\draw [fill=black] (0.,0.) circle (1pt);
\draw [fill=black] (0.,0.5) circle (1pt);
\end{scriptsize}
\draw(0.3,0.1) node {\tiny #1};
\draw(0.4,0.6) node {\tiny #2};
\end{tikzpicture}}

%% file: posets.tex
\newcommand{\ptroisun}{\begin{tikzpicture}[line cap=round,line join=round,>=triangle 45,x=0.5cm,y=0.5cm]
\clip(-0.5,-0.1) rectangle (0.5,0.7);
\draw [line width=0.5pt] (0.,0.5)-- (-0.3,0.);
\draw [line width=0.5pt] (0.,0.5)-- (0.3,0.);
\begin{scriptsize}
\draw [fill=black] (-0.3,0.) circle (1pt);
\draw [fill=black] (0.,0.5) circle (1pt);
\draw [fill=black] (0.3,0.) circle (1pt);
\end{scriptsize}
\end{tikzpicture}}

\newcommand{\pdtroisun}[3]{\begin{tikzpicture}[line cap=round,line join=round,>=triangle 45,x=0.5cm,y=0.5cm]
\clip(-0.7,-0.1) rectangle (0.7,1.2);
\draw [line width=0.5pt] (0.,0.5)-- (-0.3,0.);
\draw [line width=0.5pt] (0.,0.5)-- (0.3,0.);
\begin{scriptsize}
\draw [fill=black] (-0.3,0.) circle (1pt);
\draw [fill=black] (0.,0.5) circle (1pt);
\draw [fill=black] (0.3,0.) circle (1pt);
\end{scriptsize}
\draw(0.,0.8) node {\tiny #1};
\draw(-0.6,0.1) node {\tiny #2};
\draw(0.6,0.1) node {\tiny #3};
\end{tikzpicture}}

%% file: graphes.tex
\newcommand{\grtroisun}{\begin{tikzpicture}[line cap=round,line join=round,>=triangle 45,x=0.5cm,y=0.5cm]
\clip(-0.5,-0.1) rectangle (0.5,0.7);
\draw [line width=0.5pt] (0.,0.)-- (-0.3,0.5);
\draw [line width=0.5pt] (0.,0.)-- (0.3,0.5);
\draw [line width=0.5pt] (-0.3,0.5)-- (0.3,0.5);
\begin{scriptsize}
\draw [fill=black] (-0.3,0.5) circle (1pt);
\draw [fill=black] (0.,0.) circle (1pt);
\draw [fill=black] (0.3,0.5) circle (1pt);
\end{scriptsize}
\end{tikzpicture}}

\newcommand{\grdun}[1]{\begin{tikzpicture}[line cap=round,line join=round,>=triangle 45,x=0.5cm,y=0.5cm]
\clip(-.2,-.1) rectangle (0.5,0.5);
\begin{scriptsize}
\draw [fill=black] (0.,0.) circle (1pt);
\end{scriptsize}
\draw(0.3,0.1) node {\tiny #1};
\end{tikzpicture}}

\newcommand{\grddeux}[2]{\begin{tikzpicture}[line cap=round,line join=round,>=triangle 45,x=0.5cm,y=0.5cm]
\clip(-.2,-.1) rectangle (0.5,1.);
\draw [line width=.5pt] (0.,0.5)-- (0.,0.);
\begin{scriptsize}
\draw [fill=black] (0.,0.) circle (1pt);
\draw [fill=black] (0.,0.5) circle (1pt);
\end{scriptsize}
\draw(0.3,0.1) node {\tiny #1};
\draw(0.3,0.6) node {\tiny #2};
\end{tikzpicture}}

\newcommand{\grdtroisun}[3]{\begin{tikzpicture}[line cap=round,line join=round,>=triangle 45,x=0.5cm,y=0.5cm]
\clip(-0.5,-0.1) rectangle (0.5,1.2);
\draw [line width=0.5pt] (0.,0.)-- (-0.3,0.5);
\draw [line width=0.5pt] (0.,0.)-- (0.3,0.5);
\draw [line width=0.5pt] (-0.3,0.5)-- (0.3,0.5);
\begin{scriptsize}
\draw [fill=black] (-0.3,0.5) circle (1pt);
\draw [fill=black] (0.,0.) circle (1pt);
\draw [fill=black] (0.3,0.5) circle (1pt);
\end{scriptsize}
\draw(0.35,0.1) node {\tiny #1};
\draw(0.3,0.8) node {\tiny #2};
\draw(-0.3,0.8) node {\tiny #3};
\end{tikzpicture}}
\newcommand{\grdtroisdeux}[3]{\begin{tikzpicture}[line cap=round,line join=round,>=triangle 45,x=0.5cm,y=0.5cm]
\clip(-0.5,-0.1) rectangle (0.5,1.2);
\draw [line width=0.5pt] (0.,0.)-- (-0.3,0.5);
\draw [line width=0.5pt] (0.,0.)-- (0.3,0.5);
\begin{scriptsize}
\draw [fill=black] (-0.3,0.5) circle (1pt);
\draw [fill=black] (0.,0.) circle (1pt);
\draw [fill=black] (0.3,0.5) circle (1pt);
\end{scriptsize}
\draw(0.35,0.1) node {\tiny #1};
\draw(0.3,0.8) node {\tiny #2};
\draw(-0.3,0.8) node {\tiny #3};
\end{tikzpicture}}